%% file: arxiv.tex
\documentclass[10pt]{article}
\usepackage[margin = 0.75in]{geometry}

\usepackage[numbers, compress, sort, square]{natbib}

\usepackage[utf8]{inputenc} % allow utf-8 input
\usepackage[T1]{fontenc}    % use 8-bit T1 fonts
\usepackage{hyperref}       % hyperlinks
\usepackage{url}            % simple URL typesetting
\usepackage{booktabs}       % professional-quality tables
\usepackage{amsfonts}       % blackboard math symbols
\usepackage{nicefrac}       % compact symbols for 1/2, etc.
\usepackage{microtype}      % microtypography
\usepackage{xcolor}         % colors

\usepackage{xparse}

\usepackage{amsmath, amsthm, amssymb, amsfonts}

\usepackage[mathscr]{euscript}

\DeclareSymbolFont{rsfs}{U}{rsfs}{m}{n}
\DeclareSymbolFontAlphabet{\mathrsfs}{rsfs}

\usepackage[mathscr]{euscript}

\usepackage{upgreek}

\usepackage{calligra}
\usepackage{enumitem}

\usepackage{nccmath}
\usepackage{booktabs}
\usepackage{empheq}
\usepackage{pifont}
\usepackage{enumitem}
\usepackage{graphicx} % more modern
\usepackage{wrapfig}
\usepackage{subfigure}
\usepackage{empheq}

\usepackage{multirow}

\usepackage{algorithm}
\usepackage{algpseudocode}

\usepackage{thmtools}
\usepackage{thm-restate}

\usepackage{nicefrac}
\usepackage{booktabs}
\usepackage{multirow}
\usepackage{rotating}
\usepackage{bm}
\usepackage{comment}
\usepackage{dsfont}
\usepackage[framemethod=tikz]{mdframed}
\usetikzlibrary{shapes}
\usetikzlibrary{decorations.markings}
\usetikzlibrary{plotmarks}
\usetikzlibrary{patterns}
\usepackage{ctable}

\definecolor{myblue}{HTML}{D2E4FC}
\definecolor{Gray}{gray}{0.92}

\global\mdfdefinestyle{exampledefault}{%
outerlinewidth=0pt,innerlinewidth=0pt,
roundcorner=5pt,backgroundcolor=myblue,
leftmargin=0.2cm,rightmargin=0.2cm
}

\usepackage{color, colortbl}
\definecolor{mydarkblue}{rgb}{0,0.08,0.45}

\usepackage{empheq}

\definecolor{myteal}{RGB}{27,158,119}
\definecolor{myorange}{RGB}{217,95,2}
\definecolor{myred}{RGB}{231,41,138}
\definecolor{mypurple}{RGB}{152,78,163}
\definecolor{mygreen}{RGB}{0,100,0}

\newtheorem{definition}{Definition}[section]

\newtheorem{proposition}{Proposition}[section]
\newtheorem{lemma}{Lemma}[section]
\newtheorem{theorem}{Theorem}[section]
\newtheorem{remark}{Remark}[section]
\newtheorem{corollary}{Corollary}[section]
\newtheorem{assumption}{Assumption}

\def \EE {\mathbb{E}}
\DeclareMathOperator{\Tr}{Tr}

\title{Homogenization of $\ell_2$-Adversarial Training in High-Dimensions: Exact Dynamics under Stochastic Gradient Descent}

\author{%
  Fabrizzio Sabelli\\
  Université de Montréal\\
  \texttt{fabrizzio.sabelli@umontreal.ca} \\
}
\include{macros}

\begin{document}

\maketitle

\maketitle
\begin{abstract}
    \noindent We develop a framework for analyzing the learning dynamics of $\ell_2$-adversarial training of single-index models on Gaussian mixtures in the high-dimensional limit under streaming stochastic gradient descent (SGD). We derive deterministic equivalents for a broad class of statistics of the SGD iterates, including the adversarial risk and distance to adversarial optimality, in terms of the solution to a system of ODEs. We use them to study two idealized learning rate schedules: the Polyak stepsize and exact line search. In the case of $\ell_2$-adversarial least squares with a single class, we show that, unlike noiseless standard least squares, no constant learning rate guarantees monotone descent of SGD towards a minimizer of the adversarial risk. We identify anisotropic covariance and a mismatch in ridge parameters as the main sources of suboptimality of exact line search relative to the Polyak stepsize. We also introduce a stochastic differential equation (SDE), called adversarial homogenized SGD, that captures the evolution of statistics of the iterates of SGD. For $\ell_2$-adversarial least squares, using this SDE, we show the evolution of the risk is equivalent, up to dimension-free constants, to that of SGD on standard least squares with an adaptive learning rate and adaptive $\ell_2$-regularization. When the dynamics converge, the limiting adversarial risk and SGD iterate are determined by a fixed-point equation, with the limiting iterate being equivalent to the solution of a ridge regression problem whose regularization parameter is the limiting effective regularization of SGD.
\end{abstract}
\section{Introduction}
State-of-the-art machine learning models have been shown to be vulnerable to small input perturbations, better known as adversarial examples \citep{goodfellow2014Intrig, goodfellow2015harnessing, fawzi2018analysis, ilays2019adversarial, yuan2019adversarial}. These perturbations, often imperceptible to the human eye, are designed to maximize the model's error and result in large drops in model performance. 

Adversarial training has emerged as one of the most effective and persistent defenses against such perturbations \citep{goodfellow2015harnessing, madry2018comp, Bai2021training}. It consists of training on adversarial examples in order to increase the robustness of the model to unseen adversarial examples. It replaces the standard risk with a robust risk in which each sample is evaluated against its worst-case perturbation. When minimizing this objective with stochastic gradient descent (SGD), the update depends not only on the current iterate and the newly sampled data point, but also on the adversarial perturbation of this data point, thereby modifying the training dynamics of SGD. This dependence manifests empirically in the higher sensitivity of adversarial training to optimization hyperparameters, including the learning rate schedule, training time, weight decay, and early stopping \citep{rice2020overfittingadversariallyrobustdeep,pang2021bagtricksadversarialtraining, dabouei2022revisitingouteroptimizationadversarial}.

In this work, we analyze the dynamics of adversarial training of generalized linear models under $\ell_2$-adversarial perturbations and SGD. We extend the frameworks developed in \citep{collinswoodfin2023hitting, collinswoodfin2025Exact} to account for adversarial training of single-index models. In these frameworks, the objective function is defined as follows
\begin{equation} \label{eq:intro_function}
  \min_{X \in \mathbb{R}^{d}} \Big \{ \mathcal{R}_{\lambda}(X) \defas \EE_{a,I,\epsilon} [f_I( X^\top a  , \epsilon) ] + \frac{\lambda}{2} \|X\|_2^2
  \Big \},
\end{equation}
where $\lambda \geq 0$ is a regularization parameter, $(a, I, \epsilon) \sim \mathcal{D} \subset \mathbb{R}^d \times \mathbb{R}\times \mathbb{R}$, $X, X^{\star} \in \R^{d}$ and $\{f_i \, : \, \mathbb{R}^{2}\to \mathbb{R}\}_{i=1}^2$ satisfy certain regularity and growth assumptions which we formulate below. The $\ell_2$-adversarial counterpart to \eqref{eq:intro_function} takes the form of the following min-max optimization problem
\begin{equation} \label{eq:intro_function:adv}
  \min_{X \in \mathbb{R}^{d}} \Big \{ \mathcal{R}_{\lambda}^{\operatorname{adv}}(X) \defas \EE_{a,I,\epsilon} \left[\max_{\|\Delta a \| \leq \delta}f_I( X^\top (a + \Delta a), \epsilon) \right]+ \frac{\lambda}{2}\|X\|_2^2  \Big \},
\end{equation}
where $\delta>0$ corresponds to the maximum $\ell_2$-perturbation parameter which limits the perturbation of $a$ within a euclidean ball of radius $\delta$. Given that the ball $\{\Delta a ; \|\Delta a\|\leq \delta\}$ is convex, centered, the inner product is linear and the $\ell_2$-norm is self-dual then $\{\langle \Delta a, X\rangle : \|\Delta a\|\leq \delta\}=[-\delta \| X\|, \delta \|X\|]$ from which we obtain the following simplification
\footnote{
    If $f_i$ is convex and lower semicontinuous in its first argument, the maximization over $s \in [-1,1]$ is attained at an endpoint; under the differentiability assumptions imposed below, lower semicontinuity is automatic.
}
\begin{align}\label{eq:main:adversarial:risk}
    \mathcal{R}_{\lambda}^{\operatorname{adv}}(X) \defas \EE_{a,i,\epsilon} [\max_{|s| \leq 1}f_i( X^\top a  + \delta s \|X\|, \epsilon) ]+ \frac{\lambda}{2}\|X\|_2^2, \quad   \mathcal{R}^{\operatorname{adv}}(X) \defas  \mathcal{R}_{0}^{\operatorname{adv}}(X) .
\end{align}
 We suppose $a \in \R^d$ is distributed according to a Gaussian mixture with $2$ classes. We minimize \eqref{eq:main:adversarial:risk} using streaming SGD and explicitly solve the maximization problem at each step. At each SGD update, a new data triplet $(a,I, \epsilon)$ is sampled such that $a$ belongs to class $I$ and $\epsilon \sim \mathcal{N}(0, \eta^2)$ is independent label noise. Our main goal is not to study the optimal method for computing adversarial examples, but rather to provide a framework for analyzing how the adversarial objective affects the dynamics of SGD.
 
 \begin{figure}[t!]
    \centering
    \includegraphics[width=0.45\textwidth]{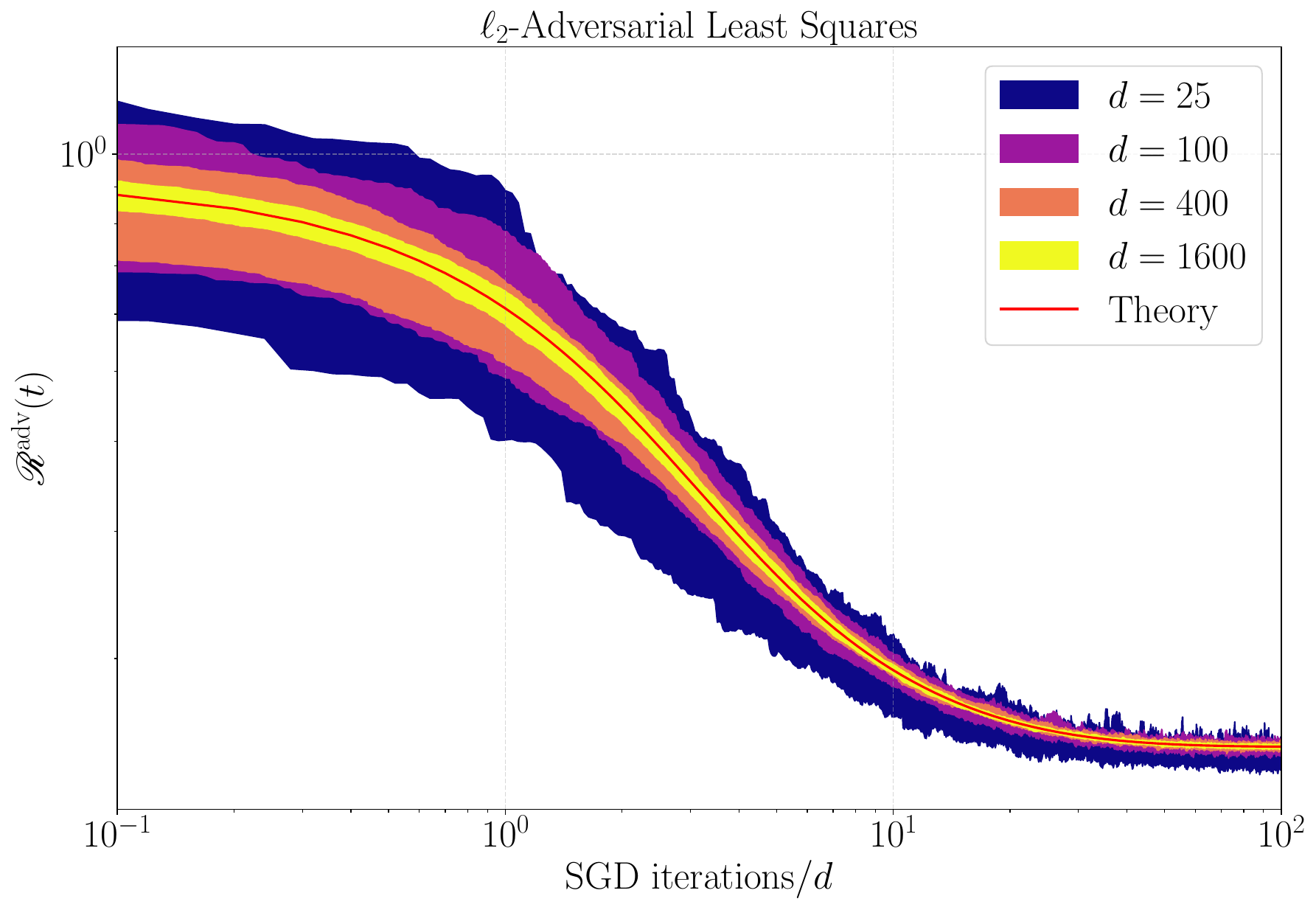}
    \includegraphics[width=0.45\textwidth]{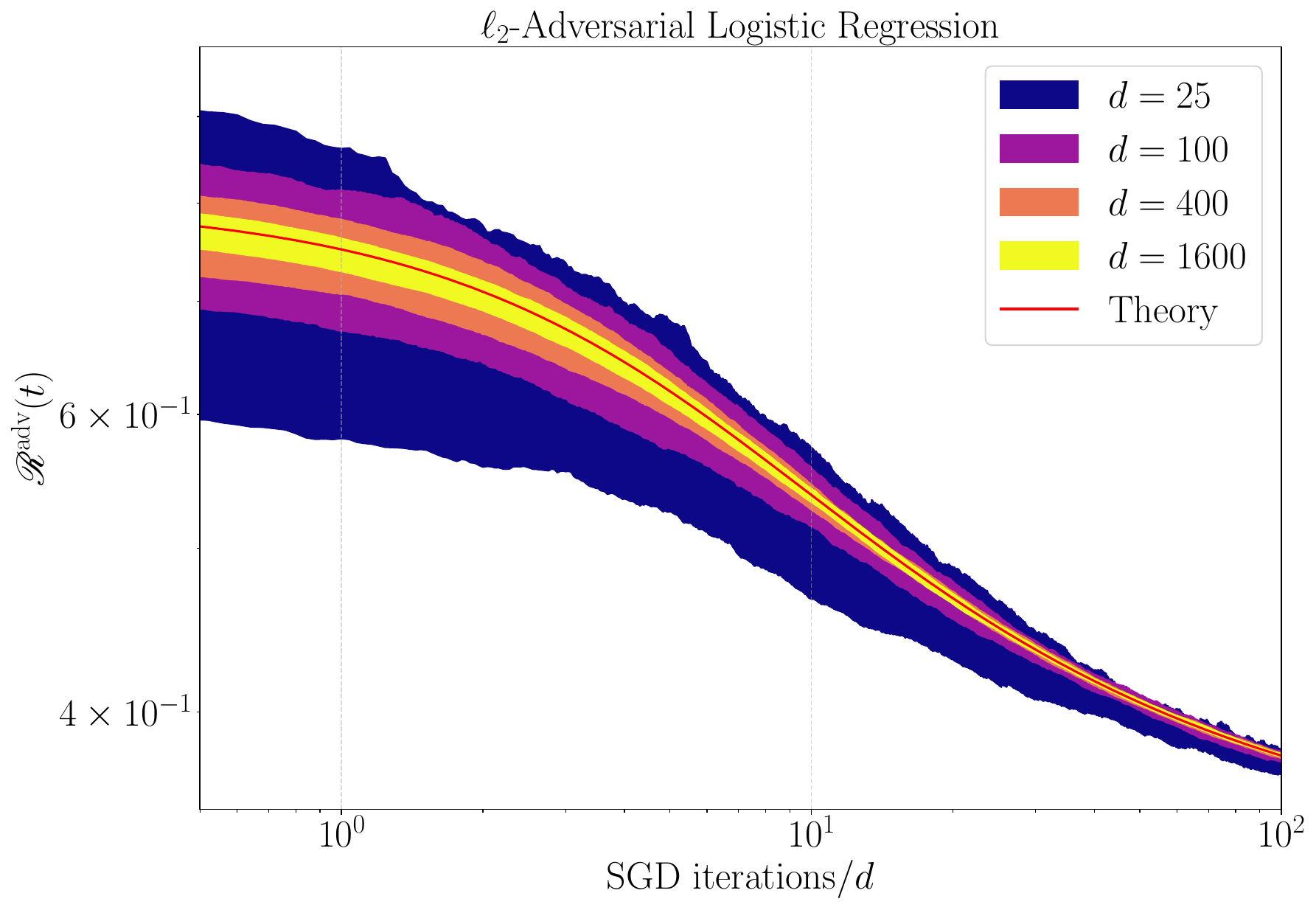}
    \caption{\textbf{Concentration of $\ell_2$-adversarial risk} on noiseless $\ell_2$-adversarial least squares with a single class $a\sim N(0,K)$ (left) and noiseless binary logistic regression with hard labels on a mixture of Gaussians (right) with different means and same covariance. As dimension $d$ increases, in both plots the adversarial risk concentrates around the deterministic limit \textcolor{red}{(red)} described by the system of ODEs \eqref{eq:ODE:V_i(t):def} as predicted by Theorem~\ref{thm:main:conc:St:general:with:stop:adv}. See Appendix~\ref{app:captions} for simulation details.}
    \label{fig:concentration}
\end{figure}

\paragraph{Main contributions.} 
\begin{itemize}
    \item We study $\ell_2$-adversarial training of single-index models when input dimension $d$ is proportional to the number of samples under data generated from a Gaussian mixture with two classes. We derive a deterministic equivalent which consists of a system of $d$-dimensional coupled ODEs which allows us to predict the learning curve of SGD for general covariances and means. See Figure~\ref{fig:concentration}. As the dimension $d$ increases, the adversarial risk $\mathcal{R}^\Adv(X_k)$ concentrates around $\mathrsfs{R}^\Adv(t)$ with overwhelming probability (see Theorem~\ref{thm:main:conc:St:general:with:stop:adv}).
    \item We introduce a SDE called \emph{Adversarial Homogenized SGD} which captures the dynamics of SGD as $d \to \infty$, even for large learning rate or above the convergence threshold. For data sampled from a single class $a \sim \mathcal{N}(0, K)$, using this SDE we show that the behaviour of the adversarial risk of $\ell_2$-adversarial least squares is equivalent up to constants (independent of $d$) to the adversarial risk of SGD on standard least squares with an adaptive learning rate and adaptive $\ell_2$-regularization parameter.
    \item We analyze the deterministic equivalents to characterize the optimal learning rate (\emph{Polyak stepsize}) which guarantees descent towards $X^{\star, \Adv} \in \operatorname{argmin}_X\; \mathcal{R}^{\operatorname{adv}}(X)$ and the optimal learning rate (\emph{exact line search}) which guarantees decrease in $\mathcal{R}^{\Adv}$ at the next SGD update. For $\ell_2$-adversarial least squares with data sampled from a single class $a \sim \mathcal{N}(0, K)$, we show that no constant learning rate guarantees monotone descent of SGD towards $X^{\star,\Adv}$, in contrast to noiseless standard least squares. For a given state of the deterministic equivalents, we identify the two main factors driving suboptimality of exact line search in comparison to the Polyak stepsize: anisotropic covariance and a mismatch in ridge parameters (see  \eqref{eq:ratio:linesearch:Polyak}). 
    \item For $\ell_2$-adversarial least squares with data sampled from a single class $a \sim \mathcal{N}(0, K)$, we show that the deterministic equivalents of $\mathcal{R}(X_k)$ and $\|X_k\|^2$ solve a system of two coupled nonlinear Volterra equations. We identify the ratio $\|X_k\| / \sqrt{2 \mathcal{R}(X_k)}$ as the key statistic governing the dynamics of SGD. Whenever this quantity converges as $k \to \infty$, for each fixed $d \geq 1$ we characterize its limit as the solution of a fixed-point equation, which also yields the limiting values of the deterministic equivalents of $\mathcal{R}(X_k)$, $\|X_k\|^2$ and $\mathcal{R}^\Adv(X_k)$. We further characterize the limiting iterate $X_{\infty}$ of SGD as the solution to a ridge regression problem whose regularization parameter is the limiting effective regularization of SGD. Taking the learning rate $\gamma \downarrow 0$ recovers the results of \citep{xing2021same_minimizer}. Under power law behaviour of the spectrum of $K$ and $X^\star$, we characterize the fixed-point equation as $d\to \infty$.   
\end{itemize}

\paragraph{Notation.} Define $\mathbb{R}_+ = [0, \infty)$. We say an event holds \emph{with overwhelming probability, w.o.p.,} 
if there exists a function $\omega: \N \to \R$ with $\omega(d)/\log d \to
\infty $ such that the event holds with probability at least
$1-e^{-\omega(d)}$. We let $\mathds{1}_{A}(x)$ be the indicator function of
the set $A$ where it is $1$ if $x \in A$ and $0$ otherwise. For a matrix $A \in \mathbb{R}^{m \times d}$, we use $\|A\|$ to denote the
Frobenius norm, $\|A\|_{\opt}$ to denote the operator-2
norm. For normed vector spaces $\mathcal{A}$, $\mathcal{B}$  with norms $\|\cdot\|_{\mathcal{A}}$ and $\|\cdot\|_{\mathcal{B}}$, respectively, and for $\alpha \geq 0$, we say a function $F \, : \, \mathcal{A} \to \mathcal{B}$ is \textit{$\alpha$-pseudo-Lipschitz} with constant $L$ if for any $A, \widehat{A} \in \mathcal{A}$, we have 
\[
\|F(A) - F(\widehat{A})\|_{\mathcal{B}} \le L\|A-\widehat{A}\|_{\mathcal{A}} (1 + \|A\|_{\mathcal{A}}^{\alpha} + \|\widehat{A}\|_{\mathcal{A}}^{\alpha} ).
\]
We write $f(t) \asymp g(t)$ if there exist \textit{absolute} constants $C, c > 0$ such that $c \cdot g(t) \le f(t) \le C \cdot g(t)$ for all $t$. If the constants depend on parameters, e.g., $\alpha$, then we write $\asymp_{\alpha}$. We write $f(t) \sim g(t)$ as $t\to a$ if $\lim_{t\to a}\frac{f(t)}{g(t)}=1$. Finally, the operator $\oplus$ denotes vector/matrix concatenation.
\section{Related Work}\label{sect:related:work}
\paragraph{Adversarial training.} Adversarial robustness and examples have been extensively studied in machine learning
\citep{goodfellow2014Intrig, goodfellow2015harnessing, madry2018comp, 2017adversarialcomp, yuan2019adversarial, Bai2021training}. Adversarial training has also been studied extensively in the case of deep neural networks \citep{goodfellow2015harnessing, madry2018comp, shafahi2019free, rice2020overfittingadversariallyrobustdeep, zhai2019adversariallyrobustgeneralizationjust, chen2021robust}. From a theoretical perspective, the generalization and regularization properties of adversarial training as well as the robust versus standard accuracy trade-off have been studied in linear models \citep{raghunathan2020understandingmitigatingtradeoffrobustness, dan2020sharpstatisticalguaranteesadversarially, clarysse2022adversarialtraininghurtrobust, ribeiro2023regularizationpropertiesadversariallytrainedlinear, pmlr-v258-ribeiro25a, xing2021same_minimizer, Ribeiro2023over}. Under Gaussian data, exact asymptotics for these properties are derived in the proportional regime  \citep{tanner2024highdimensionalstatisticalmodel, vilucchio2024geometryregularizationadversarialtraining, javanmard2020precisetradeoffsadversarialtraining, hassani2024curseoverparametrizationadversarialtraining, pmlr-v235-dohmatob24c}. These properties were also investigated for Gaussian mixture models (GMM) \citep{dobriban2022provabletradeoffsadversariallyrobust, dan2020sharpstatisticalguaranteesadversarially, bhagoji2019lowerboundsadversarialrobustness, javanmard2022precisestatisticalanalysisclassification, taheri2021asymptoticbehavioradversarialtraining}. Using the theory of stochastic modified equations \citep{li2018stochasticmodifiedequationsdynamics}, \citep{gu2023adversarialtraininggradientdescent} derive a SDE which approximates the min-max optimization dynamics under SGD in the vanishing learning rate limit. Our diffusion approximation of SGD differs from \citep{gu2023adversarialtraininggradientdescent} as it does not rely on a vanishing stepsize but rather the input dimension $d\to \infty$. 

\paragraph{Deterministic dynamics of stochastic algorithms in high-dimensions.}
As mentioned in \citep{collinswoodfin2024highline,collinswoodfin2025Exact}, the theory of deterministic dynamics of stochastic algorithms on isotropic Gaussian data has a long history \citep{biehl1994line,biehl1995learning,saad1995dynamics,saad1995exact}. These rigorous results were also extended to other models \citep{goldt2019dynamics, wang2019solvable, arnaboldi2023highdimensional, damian2023smoothinglandscapeboostssignal, arnaboldi2024escapingmediocritytwolayernetworks, dandi2024benefits} under isotropic Gaussian data. Deterministic equivalents were also proven for multi-pass SGD with small mini-batches \citep{PPAP01} and with momentum \citep{LeeChengPaquettePaquette}. Other high-dimensional deterministic equivalents yielding different types of dynamics have also been studied \citep{Mignacco_2021, gerbelot2022rigorous, celentano2026highdimensionalasymptoticsordermethods, fan2026high, chandrasekher2021sharp, bordelon2022learningcurvessgdstructured, arous2023highdimensionallimittheoremssgd}. The works \citep{ben2022high, arous2026localgeometryhighdimensionalmixture, refinetti2021classifyinghighdimensionalgaussianmixtures} study the high-dimensional training dynamics of streaming SGD under GMM with isotropic covariance and the work \citep{Mignacco_2021} studies the analogous setting under multi-pass SGD.

In recent years, a line of work has made significant contributions to understanding the effect of non-identity covariance on the training dynamics of SGD under data from a single class \citep{CollinsWoodfinPaquette01, balasubramanian2024highdimensionalscalinglimitsfluctuations, goldt2022gaussian, Yoshida, Goldt, collinswoodfin2023hitting, xiao2026exactriskcurvessignsgd, marshall2024clipclipdynamicssgd}. As mentioned in \citep{collinswoodfin2023hitting}, the non-identity covariance modifies the optimization landscape and affects convergence properties.
The work \citep{collinswoodfin2025Exact} extends \citep{collinswoodfin2023hitting} to Gaussian mixture data with non-identity covariance and where the number of classes grows logarithmically with the input dimension $d$. Our work extends \citep{collinswoodfin2023hitting, collinswoodfin2025Exact} by accounting for the $\ell_2$-adversarial objective under Gaussian mixture data with two classes and non-identity covariance.
\subsection{Model Set-up} \label{sec:main_setup}
We study the setting of streaming SGD characterized by a sequence of independent samples $\{a_k, y_k\}$. Here $y_k$ is the target which is either a function of data triplet $(a_k, I_k, \epsilon_k)$ or a hard label. We formalize this in the following assumptions analogously to \citep{collinswoodfin2023hitting, collinswoodfin2025Exact}.
\begin{assumption}[Class and data structure]
\label{assumption:data}
The triplet $(a, I, \epsilon) \sim \mathcal{D}$ is distributed as follows: $a$ is generated from a mixture of Gaussians with 2 classes such that $a$ belongs to class $1$ with probability $p_1$ and to class $2$ with probability $p_2=1-p_1$. Conditioned on belonging to class $i$, $a\mid i\sim \mathcal{N}(\mu_i, K_i)$ with covariance $K_i \in \mathbb{R}^{d \times d}$ such that $\| K_i \|_{\opt} \leq C$ for some constant $C$ (independent of $d$). Furthermore, $K_1$ and $K_2$ commute, that is $[K_1, K_2] = 0$. We also scale the means $\|\mu_i\|^2 \leq C$. Finally, let $\mu \defas \mu_1 \oplus \mu_2 \in \R^{d\times 2}$. 
\end{assumption}
\begin{assumption}[Target model and label noise]\label{assumption:target}
    The targets satisfy one the following models:
   %\vspace{-5pt} 
    \begin{enumerate}[label=\alph*)]
        \item (Hard Label) $y \mid i = i$ with $i \in \{-1,1\}$
        \item (Soft Label) For a fixed vector $X^{\star}$ and a link function $g: \R \to \R$ such that $y\mid i = g(X^{\star\top}a; \epsilon)$ for $a \mid i \sim \mathcal{N}(\mu_i, K_i)$. The label noise $\epsilon \sim \mathcal{N}(0, \eta^2)$, independent of $(a,i)$.
    \end{enumerate}
\end{assumption}
Given that the data $a$ is sampled from a Gaussian mixture, the adversarial risk \eqref{eq:main:adversarial:risk} can be decomposed as follows
\begin{equation}\label{eq:nlgc}
    \mathcal{R}_{\lambda}^{\operatorname{adv}}(X) = p_1 \cdot \EE_{v,\epsilon} [\Psi_1(X; \sqrt{K_1}v + \mu_1, \epsilon) ] + p_2 \EE_{v,\epsilon} [\Psi_2(X; \sqrt{K_2}v + \mu_2 , \epsilon)  ]+ \frac{\lambda}{2}\|X\|_2^2,
\end{equation}
where we define $\Psi_i(X; \sqrt{K_i}v+ \mu_i, \epsilon) \defas\max_{|s| \leq  1}f_i( X^\top (\sqrt{K_i}v+ \mu_i) + \delta s \|X\|,   \epsilon)$ for $i=1,2$ and $v \sim \mathcal{N}(0, \Id_d)$. Due to the maximization in \eqref{eq:main:adversarial:risk}, we must first formalize when the objective is differentiable. 
\begin{assumption}[Pseudo-lipschitz $f$] 
\label{assumption:pseudo_lipschitz}
The function $f_i :\R^3 \to \R$ is $\alpha$-pseudo-Lipschitz with constant $L(f)>0$ for $i=1,2$. 
%In other words, for all $r_1, r_2 \in \R^2$ and $\epsilon \in \R$
%\[
%|f_i(r; \epsilon) - f_i(\widehat{r}; \epsilon)| \le L(f) \| r - \widehat{r} \| (1 + \|r\|^\alpha + \|\widehat{r} \|^\alpha + \|\epsilon\|^{\alpha}).
%\]    
\end{assumption}
As we are minimizing \eqref{eq:main:adversarial:risk}, we consider $f_i(X^\top a, \epsilon)$ as a function of its first variable $X^\top a$ with parameter $\epsilon$ and write $f_i(X^\top a, \epsilon)=f_i(X^\top a;\epsilon)$. We denote the (almost everywhere) derivative with respect to the first variable as $f_i'$. Assumption~\ref{assumption:pseudo_lipschitz} also implies that $\Psi_i(X; a, \epsilon)$ is differentiable (almost everywhere) with respect to $X$. We formalize the conditions for this statement to hold.
\begin{theorem}[Danskin's theorem \citep{Danskin1967TheTO}]\label{thm:danskin}
    Let $f_i^{\operatorname{adv}}(X, s) = f_i( X^\top a + \delta s \|X\| , \epsilon)$ then $f^{\operatorname{adv}}:\mathbb{R}^d  \times [-1, 1] \to \mathbb{R}$ is differentiable (almost everywhere) in its two arguments. Denote $S_0(X)$ as the set of maximizing points
    \begin{align}
        S_0(X) = \left\{\bar{s}: \bar{s} = \operatorname{argmax}_{|s| \leq  1}f_i^{\operatorname{adv}}(X, s)\right\}.
    \end{align}
    If $S_0$ contains a single unique element $\bar{s}$ and $X\not=0$, then the gradient of $\Psi_i(X; a,\epsilon)$ is given by
    \begin{equation}
        \nabla_X \Psi_i(X; a,\epsilon) =  f_i'( X^\top a + \delta \bar{s} \|X\|   , \epsilon)\left(a + \delta \bar{s} \frac{X}{\|X\|} \right).
    \end{equation}
    Moreover, define the subgradient of $f_i^{\operatorname{adv}}$ at $X$
    \begin{equation}
    \partial  \Psi_i(X; a,\epsilon)= \operatorname{conv}\left\{ \nabla_X\left[f_i( X^\top a + \delta s \|X\| ,  \epsilon)\right]: s \in S_0(X)\right\}.
    \end{equation}
    If $ \partial  \Psi_i(X; a,\epsilon) = \{v\} $ (i.e. $\nabla_X\left[f_i( X^\top a + \delta s \|X\| ,  \epsilon)\right] = v$ for all $s \in S_0(X)$) then $ \Psi_i(X; a,\epsilon)$ is differentiable at $X$ with $\nabla_X \Psi_i(X; a,\epsilon) = v$. %Furthermore, if $f_i^{\operatorname{adv}}$ is convex then $\Psi_i$ is convex.
\end{theorem}
Theorem~\ref{thm:danskin} implies that the function $\Psi_i$ is not differentiable if there exists $s_1 < s_2 \in [-1,1]$ such that 
\begin{equation}
    \begin{aligned}   
    &f_i( X^\top a + \delta s_1 \|X\|) = f_i( X^\top a + \delta s_2 \|X\|)= \max_{s \in [-1,1]} f_i( X^\top a + \delta s \|X\|) 
    \\
    \text{and} \quad &f_i'( X^\top a + \delta s_1 \|X\|)(a + \delta s_1 \tfrac{X}{\|X\|}) \not = f_i'( X^\top a + \delta s_2 \|X\|)(a + \delta s_2 \tfrac{X}{\|X\|}).
\end{aligned}
\end{equation}
In our theoretical analysis, we will impose a stopping time to study the behaviour of SGD until the first time at which Danskin's theorem does not hold. It is easy to see that if $X = 0$ then Danskin's theorem does not hold. Hence, we require the initialization of SGD to be non-zero. 

% Finally, given $X \not = 0$ non-decreasing losses also satisfy the conditions of Danskin's theorem as these losses ensure differentiability of $\Psi$ across the whole path of SGD. Indeed, if $f$ is non-decreasing then for any $s_1 < s_2 \in [-1,1]$, we have
% \begin{equation}
%     f( X^\top a + \delta s_1 \|X\|) \leq f( X^\top a + \delta s_2 \|X\|) \leq f( X^\top a + \delta \|X\|).
% \end{equation}
% Now, since $f$ is non-decreasing and differentiable in its first argument, if $s_1 < s_2$ and 
% \begin{equation}
%     f( X^\top a + \delta s_1 \|X\|) = f( X^\top a + \delta s_2 \|X\|),
% \end{equation}
% then, by the mean-value theorem for any $c \in (s_1, s_2)$ with $X \not = 0 $ 
% \begin{equation}
%     0=f( X^\top a + \delta s_1 \|X\|) - f( X^\top a + \delta s_2 \|X\|) = \delta( s_1 -s_2) \|X\|f'( X^\top a + \delta c \|X\|) .
% \end{equation}
% Hence, it follows that $f'( X^\top a + \delta c \|X\|) = 0$ for all $c \in (s_1, s_2)$. 

In order to keep a unified notation for hard and soft labels, we build upon the notation introduced in \citep{collinswoodfin2025Exact}:
\begin{equation}
    \hat{X} \defas \begin{cases}
        X \oplus X^\star\in \mathbb{R}^{d \times 2} \quad \text{(soft labels)} \\
        X \oplus 0 \in  \mathbb{R}^{d \times 2} \quad \text{(hard labels)}, \\
        \end{cases}  \text{and}\quad W = \hat{X} \oplus \mu \in \R^{d\times 4}.
\end{equation}
For $i=1,2$, let $\{(\lambda^{(i)}_j, \omega_i)\}_{1\leq j \leq d}$ be the orthonormal eigenvalue-eigenvector pairs of $K_i$ where the eigenvectors of $K_1$ and $K_2$ are the same since they commute. For $i=1,2$, we also define the following set of block matrices
\begin{equation}\label{eq:stats:BBcheck}
     \widehat{B}_i(W) \defas \begin{bmatrix}
            U_i & 0\\
            0 & \|X\|^2
            \end{bmatrix},\quad  U_i(W) \defas \begin{bmatrix}
                B_i(W)  & m_i(W) \\
                m_i(W)^\top & \|\mu_i\|^2
            \end{bmatrix},\quad B_i(W)\defas \hat{X}^\top K_i \hat{X},
    \end{equation}
and $m_i(W) \defas \hat{X}^\top \mu_i$.\footnote{We write $\widehat{B}_i = \widehat{B}_i(W)$, $U_i = U_i(W)$ and $B_i = B_i(W)$ when clear from context.} From \eqref{eq:nlgc}, it is clear that  $\mathcal{R}^{{\operatorname{adv}}}(X)$ is either an expectation over the correlated Gaussians $X^\top (\sqrt{K_i}v+ \mu_i)$ and $X^{\star\top}(\sqrt{K_i}v+ \mu_i)$ in the soft label setting or an expectation over the Gaussian $X^\top ((\sqrt{K_i}v+ \mu_i))$ in the hard label setting. We thus write $\mathcal{R}^{{\operatorname{adv}}}(X) \defas p_1h_1(\widehat{B}_1) + p_2h_2(\widehat{B}_2)$ for some well-behaved functions $h_i \, : \, \mathbb{R}^{4 \times 4} \to \mathbb{R}$.
\begin{assumption}[Risk representation]
  \label{assumption:risk} 
There exists an open convex set $\mathcal{U}\subset \mathbb{R}^{4\times 4}$ such that $\widehat{B}_1(W_0), \widehat{B}_2(W_0) \in \mathcal{U}$. Provided $\widehat{B}_i(W) \in \mathcal{U}$ for $i=1,2$, there exists functions $h_i: \mathbb{R}^{4\times 4}  \rightarrow \mathbb R$ such that $\mathcal{R}^{{\operatorname{adv}}}(X) \defas p_1 h_1(\widehat{B}_1) + p_2h_2(\widehat{B}_2)$ is differentiable and satisfies
\[
\nabla_X \mathcal{R}^{{\operatorname{adv}}}(X) = p_1 \mathbb{E}_{ v, \epsilon} [\nabla_X \Psi_1(X;  (\sqrt{K_1}v+ \mu_1), \epsilon) ]+ p_2\mathbb{E}_{v, \epsilon} [\nabla_X \Psi_2(X; (\sqrt{K_2}v+ \mu_2), \epsilon)]. 
\]
Furthermore, $h_1$ and $h_2$ are continuously differentiable on $\mathcal{U}$ and their derivatives $\nabla h_1$ and $\nabla h_2$ are $\alpha$-pseudo-Lipschitz with constant $L( h)>0$.
%such that for all $\widehat{B}_1, \widehat{B}_2 \in \mathcal{U}$ and $i=1,2$
 %\[
 %\|\nabla h_i(\widehat{B}_1)-\nabla h_i(\widehat{B}_2)\| \le L(h) \|\widehat{B}_1-\widehat{B}_2\| (1 + \|\widehat{B}_1\|^{\alpha} + \|\widehat{B}_2\|^{\alpha}).
 %\]
\end{assumption}
Given the definition of $\mathcal{R}^{{\operatorname{adv}}}$, it follows that the set $\mathcal{U}$ will not include points where $X = 0$. Moreover, we formulate $\mathcal{U}$ such that for all $(\widehat{B}_1(W), \widehat{B}_2(W)) \in \mathcal{U} \times \mathcal{U}$, Danskin's theorem holds (see \eqref{eq:example:U} for an example). It will be useful to view the gradients in block form
\begin{equation}
    \nabla h_i 
=
\left[ \begin{array}{c|c|c|c} 
   \partial_{11} h_i &\partial_{12} h_i & \partial_{13} h_i  & 0 \\
    \hline
   \partial_{21} h_i &\partial_{22} h_i &  \partial_{23} h_i & 0\\
    \hline
    \partial_{31}h_i & \partial_{32} h_i  &\partial_{33} h_i & 0 \\
    \hline
    0&0 & 0 &\partial_{44} h_i
    \end{array} \right] \in \R^{4\times 4}.
 \end{equation}
 From this representation, we derive an explicit formula in Lemma~\ref{lem:derivative_adv_risk} for $\nabla_X \mathcal{R}^{{\operatorname{adv}}}$ in terms of $\nabla h_i$. We place an assumption requiring some regularity of the Fisher information matrices of the gradients. 
\begin{assumption}[Fisher matrix] 
\label{assumption:fisher} 
For $i=1,2$, define $\mathcal{I}_i(\widehat{B}_i) \defas \mathbb{E}_{v,\epsilon}[(f'_i( X^\top ((\sqrt{K_i}v+ \mu_i)) + \bar{s}_i \delta  \| X\|;  \epsilon))^2]$ where the function $\mathcal{I}_i \, : \, \mathbb{R}^{4 \times 4} \to \mathbb{R}$ and $\bar{s}_i= \operatorname{argmax}_{|s| \leq 1}f_i( X^\top ((\sqrt{K_i}v+ \mu_i)) + \delta s \| X\|,  \epsilon)$. Furthermore, $\mathcal{I}_i$ is $\alpha$-pseudo-Lipschitz with constant $L(\mathcal{I})>0$. 
% such that for $\widehat{B}_1, \widehat{B}_2 \in \mathcal{U}$ 
 %\[
 %\|\mathcal{I}_i(\widehat{B}_1)-\mathcal{I}(\widehat{B}_2)\| \le L(\mathcal{I}) \|\widehat{B}_1-\widehat{B}_2\| (1 + \|\widehat{B}_1\|^{\alpha} + \|\widehat{B}_2\|^{\alpha}).
 %\]
\end{assumption}
A large class of regression problems fit this set of assumptions, notably generalized linear models  such as such as least squares and binary logistic regression (see Appendix~\ref{sec:lsq_analysis}, \ref{ex:binary:logreg} and \citep[Appendix D]{collinswoodfin2023hitting} for more examples). 
 \subsection{Algorithmic Set-up}
 \label{sec:algorithmic_setup_main}
We consider \textit{online/streaming} SGD with a deterministic learning rate $\gamma_k$ given by the update
    \begin{align} \label{eq:sgd}
    X_{k+1} &= X_{k} - \frac{\gamma_k}{d}\left( f_{I_{k+1}}'(g_{k, I_{k+1}} ;  \epsilon_{k+1})(a_{k+1, I_{k+1}} + \delta s_{k, {I_{k+1}}} \frac{X_k}{\|X_k\|} ) +  \lambda X_k\right)
    \end{align}
    where $a_{k} = \sqrt{K_{I_{k}}}v_k + \mu_{I_{k}}$ for $v_k \sim \mathcal{N}(0, \Id_d)$, $I_{k} \in \{1,2\}$ is randomly sampled and we define $r_{k, I_{k+1}} = \hat{X}_k^\top a_{k+1, I_{k+1}}$ such that
    \begin{equation}
        \begin{gathered}
        g_{k, I_{k+1}} = r_{k, I_{k+1}} + \begin{bmatrix}
            \delta s_{k, I_{k+1}}\|X_k\| & 0 
        \end{bmatrix}^\top,\\
        s_{k, {I_{k+1}}} = \operatorname{argmax}_{|s| \leq  1}f_{I_{k+1}}( X_k^\top a_{k+1, I_{k+1}}+ \delta s \|X_k\|; \epsilon_{k+1}).
    \end{gathered}
    \end{equation}
To perform our analysis, we require an assumption on the initialization $X_0$ and $X^{\star}$.
 \begin{assumption}[Initialization and signal]
 \label{assumption:scaling}
 For some constant $C$ independent of $d$, the initialization point $X_0 \in \R^{d} \setminus \{0\}$ and $X^\star\in \R^{d}$ satisfy $\max \{ \| X_0 \|, \| X^\star \| \} \le C$.
 \end{assumption}
 Note that the initialization $X_0$ and $X^{\star}$ may be random. In this case, the following results are stated while conditioning on their values. Finally, we place an assumption on the learning rate schedule.
 \begin{assumption}[Bounded learning rate]\label{assumption:lr}
     There is a $\bar{\gamma} < \infty $ and a continuous function $\gamma :\mathbb{R}^+\to \mathbb{R}^+$ which is bounded by $\bar{\gamma}$ such that $\gamma_k = \gamma(k/d)$.
 \end{assumption}
 \section{Deterministic Dynamics for SGD under $\ell_2$-Adversarial Training}\label{sect:mainresult:SGD}
From Assumption~\ref{assumption:risk} and Assumption~\ref{assumption:fisher}, the adversarial risk $\mathcal{R}^{{\operatorname{adv}}}(X) $ and functions $\mathcal{I}_i$ are entirely determined by the matrices $\widehat{B}_i(W)$ which involve the matrix products $B_i= \hat{X}^\top K_i \hat{X}$, $\hat{X}^\top \mu_i$ and the scalar $X^\top X$. Note that these matrices are embedded into the following polynomials in $K_i$: $W^\top K_i W$ and $W^\top W$. Hence, if we are able to compute a deterministic equivalent for $W^\top q(K_i)W^\top$ for an arbitrary polynomial $q$, then we can study the deterministic dynamics of these statistics.

Through the use of Cauchy's integral formula, we can recover these statistics using the resolvents $R(z_i;K_i) \defas (K_i-z_i\Id_d)^{-1}$ of $K_i$ with $z_i \in \mathbb{C} \setminus \text{spec}(K_i)$. For any polynomial $q$ we have
\begin{equation*}
    W^\top q(K_i)W = -\frac{1}{2\pi \operatorname{i}} \oint_{\Gamma} q(z_i) W^\top R(z_i;K_i)W \dif z_i,
\end{equation*}
where $\Gamma$ is a contour enclosing the eigenvalues of $K_i$. Since we are interested in tracking simultaneously $ W^\top q(K_1)W$ and $ W^\top q(K_2)W$, we introduce $\mathcal{K} \defas \{K_1, K_2\}$ and $\mathscr{R}(z; \mathcal{K}) \defas R(z_1;K_1) R(z_2;K_2)$. We also introduce a fixed contour $\Gamma = \{ w \in \mathbb{C} \, : \, |w| = \max\{1, 2 \|\mathcal{K}\|_{\operatorname{op}} \} \}$ with $\|\mathcal{K}\|_{\operatorname{op}} \defas \max\{ \|K_1\|_{\opt}, \|K_2\|_{\opt}\}$. Note that by construction the distance between the spectra of $K_1$ and $K_2$ and the contour is always greater or equal than $\frac{1}{2}$. 

For $z \in \Gamma^2$, we write $z =(z_1 ,z_2)$ and when integrating over $z_1$ and $z_2$ simultaneously, we write for any function $f: \mathbb{C}^2 \to \mathbb{C}$
\[
\oint f(z) \text{D}z \defas \frac{-1}{4\pi^2 }\oint_{\Gamma^2} f(z)  \dif z_1 \dif z_2.
\]
We introduce a norm over the contour $\Gamma^2$.
    \begin{definition}\label{eq:main:gamma:normdef}
       For a continuous function $H:\mathbb{C}^2 \to \mathbb{C}^{4\times 4}$, define 
        \[
            \|H\|_{\Gamma} \defas \sup_{z\in \Gamma^2}\|H(z)\|.
        \]
    \end{definition}
Similarly to \citep{collinswoodfin2025Exact}, we will be interested in computing deterministic equivalents for the following matrices 
\begin{equation}\label{eq:main:def:S(W,z)}
    S(W,z) \defas \hat{X}^\top \mathscr{R}(z; \mathcal{K})\hat{X}\quad \text{and}\quad M(W,z) \defas \hat{X}^\top \mathscr{R}(z; \mathcal{K})\mu \in \mathbb{C}^{2\times 2},
\end{equation}
which regroup under
\begin{equation}\label{eq:def:og:main:Z}
    Z(W,z) \defas W^\top \mathscr{R}(z; \mathcal{K})W = \begin{bmatrix}
        S(W,z) & M(W,z) \\
        M(W,z)^\top & \mu^\top \mathscr{R}(z; \mathcal{K}) \mu
    \end{bmatrix} \in \mathbb{C}^{4\times 4}.
\end{equation}
Let $M_i(W,z)$ denote the $i$th column of $M(W,z)$ then from Cauchy's integral formula, we can recover $\widehat{B}_i(W)$ as follows
\begin{equation}\label{eq:Bcheck:contourint}
    \widehat{B}_i(W) = \begin{bmatrix}
         \oint z_i S(W,z) \Dif z & \oint M_i(W,z) \Dif z  & 0 \\
         \oint M_i(W,z)^\top \Dif z & \|\mu_i\|^2 & 0 \\
         0& 0 &   \oint S_{11}(W,z) \Dif z .
    \end{bmatrix}
\end{equation}
Given that we are studying the dynamics of SGD, the deterministic equivalent of $Z(W_k, z)$ will evolve with time. This deterministic equivalent which we denote $\mathcal{Z}(t,z)$ corresponds to the solution of the following integral-differential equation \eqref{eq:ODE_resolvent_2}.
\begin{mdframed}[style=exampledefault]
    \textbf{Integro-Differential Equation for $\mathcal{Z}(t, z)$.} For any contour $\Gamma \subset \mathbb{C}$ enclosing the eigenvalues of $K_1$ and $K_2$, we have an expression for the derivative of $\mathcal{Z}$ for $z\in \Gamma^2$:
    \begin{equation}\label{eq:ODE_resolvent_2}
        \dif \mathcal{Z}(t, \cdot) = \mathscr{F}(z, \mathcal{Z}(t, \cdot)) \, \dif t,\quad \text{with initial condition}\quad  \mathcal{Z}(0, z) = Z(W_0, z),
        \end{equation}
        where
        \begin{equation}
            \begin{gathered}
            \mathcal{Z}(t, z) \defas  \begin{bmatrix}
                \mathcal{S}(t, z) &  \mathcal{M}(t, z) \\
                \mathcal{M}(t, z)^\top & \mu^\top \mathscr{R}(z;\mathcal{K})\mu
            \end{bmatrix},
            \\
            \text{and}\quad
            \mathscr{F}(z, \mathcal{Z}(t, \cdot)) \defas  \begin{bmatrix}
                \mathscr{F}_\mathcal{S}(z, \mathcal{Z}(t, \cdot)) & \mathscr{F}_\mathcal{M}(z,\mathcal{Z}(t, \cdot)) \\
                \mathscr{F}_\mathcal{M}(z, \mathcal{Z}(t, \cdot))^\top & 0_{2\times 2}
            \end{bmatrix},
        \end{gathered}
        \end{equation}
The maps $\mathscr{F}_\mathcal{S}$ and $\mathscr{F}_\mathcal{M}$ are defined as follows
\begin{equation}\label{eq:ODE_resolvent_2:M}
    \begin{aligned}
        &\dif \mathcal{M}(t,z) 
        =
        -2\gamma(t) \sum_{i=1}^2 p_i\bigg(H_{1,i}^\top( \widehat{\mathrsfs{B}}_i(t))  \left(z_i\mathcal{M}(t,z)- \frac{1}{2\pi\operatorname{i}} \oint_{\Gamma} \mathcal{M}(t, z) \dif z_i\right)
        \\
        &
        + H_{2,i} (\widehat{\mathrsfs{B}}_i(t)) \otimes \mu^\top \mathscr{R}(z; \mathcal{K})\mu_i
        + \left(\frac{\lambda}{2} + \partial_{44} h_i (\widehat{\mathrsfs{B}}_i(t) )\right)   \operatorname{D}  \mathcal{M}(t,z)\bigg)\dif t,
    \end{aligned}
    \end{equation}
    and
    \begin{equation}\label{eq:ODE_resolvent_2:S}
        \begin{aligned}
    &\dif \mathcal{S}(t,z) 
    =
    -2\gamma(t) \sum_{i=1}^2 p_i\bigg(\left(z_i\mathcal{S}(t,z)- \frac{1}{2\pi\operatorname{i}} \oint_{\Gamma} \mathcal{S}(t, z) \dif z_i\right) H_{1,i}(\widehat{\mathrsfs{B}}_i(t)) 
    \\
    &
     + H_{1,i}^\top( \widehat{\mathrsfs{B}}_i(t))  \left(z_i\mathcal{S}(t,z)- \frac{1}{2\pi\operatorname{i}} \oint_{\Gamma} \mathcal{S}(t, z) \dif z_i\right)
    \\
    &
    + \mathcal{M}_i(t, z) \otimes H_{2,i} (\widehat{\mathrsfs{B}}_i(t)) 
    + H_{2,i} (\widehat{\mathrsfs{B}}_i(t)) \otimes \mathcal{M}_i(t, z)
    \\
    &
    + \left(\frac{\lambda}{2} + \partial_{44} h_i (\widehat{\mathrsfs{B}}_i(t) )\right) (\mathcal{S}(t,z) \operatorname{D}  +  \operatorname{D}  \mathcal{S}(t,z))\bigg)\dif t
    \\
    &
    + \frac{\gamma(t)^2}{d} \sum_{i=1}^2 p_i
      \mathcal{I}_i(\widehat{\mathrsfs{B}}_i(t))\Tr(( K_i + \mu_i \mu_i^\top ) \mathscr{R}(z;\mathcal{K})) \operatorname{D}  \, \dif t,
\end{aligned}
\end{equation}
\begin{gather}
\text{where} \quad
H_{1,i}(\widehat{\mathrsfs{B}}_i)=\left [ \begin{array}{c|c} 
   \partial_{11} h_i(\widehat{\mathrsfs{B}}_i)  & 0 \\
    \hline
   \partial_{21} h_i(\widehat{\mathrsfs{B}}_i) & 0 
    \end{array} \right ],\quad
    H_{2,i}(\widehat{\mathrsfs{B}}_i)= \left [ \begin{array}{c} 
        \partial_{13} h_i(\widehat{\mathrsfs{B}}_i)   \\
         \hline
        0 
         \end{array} \right ],
         \nonumber
         \\
         \operatorname{D} = \left [ \begin{array}{c|c} 
            1 & 0\\
            \hline
            0 & 0
            \end{array} \right ]
           \quad \text{and}\quad 
           \widehat{\mathrsfs{B}}_i(t) = \begin{bmatrix}
               \mathrsfs{B}_i(t) & \mathfrak{m}_i(t) & 0  \\
               \mathfrak{m}_i(t)^\top & \|\mu_i\|^2 & 0 \\
               0 & 0 & \oint \mathcal{S}_{11}(t,z) \Dif z
            \end{bmatrix}.\label{eq:def:B(t):main}
\end{gather}
Here $\mathrsfs{B}_i(t) = \oint z_i \mathcal{S}(t,z) \Dif z \in \R^{2\times 2}$, $\mathfrak{m}_i(t) = \oint \mathcal{M}_i(t) \Dif z$ where $\mathcal{M}_i(t) \in \R^{2}$ denotes the $i$th column of $\mathcal{M}$, associated with $\mu_i$. For $i=1,2$, the functions $h_i:\R^{4\times 4} \to \R$ and $\mathcal{I}_i:\R^{4\times 4} \to \R$ are defined in Assumptions~\ref{assumption:risk}~and~\ref{assumption:fisher} respectively. Finally, here $W_0 = \hat{X}_0 \oplus \mu$ with $\hat{X}_0 = X_0 \oplus X^{\star}$ in the soft label setting and $\hat{X}_0 = X_0 \oplus 0$ in the hard label setting.
\end{mdframed}
\begin{remark}
    We note that, in the hard label setting, since $X^{\star} = 0$ and from the definition of $\mathcal{Z}(t,z)$, $\mathcal{S}(t,z)$ and $\mathcal{M}(t,z)$, the only non-zero terms in the integro-differential equation \eqref{eq:ODE_resolvent_2} are $\mathcal{S}_{11}(t,z)$ and the first row in $M(t,z)$. 
\end{remark}

Note that the integral-differential equation in \eqref{eq:ODE_resolvent_2} corresponds to the adversarial version of the equations presented in \citep{collinswoodfin2023hitting} and \citep{collinswoodfin2025Exact}. Hence, the two equations follow a very similar structure which differs due to the need to track the statistic $\|X_k\|^2$ as presented in \eqref{eq:stats:BBcheck} and due to different definitions for the functions $h_i$ and $\mathcal{I}_i$ defined in Assumptions~\ref{assumption:risk}~and~\ref{assumption:fisher} respectively.

Since SGD is a discrete time process, we introduce the following time parametrization to compare it to its continuous time limit:
    \begin{equation*}
        k\text{ iterations of SGD} = \lfloor td \rfloor,\quad \text{for } t\in \R^+ \text{ the continuous time parameter}.
    \end{equation*}
Here $\lfloor\cdot \rfloor$ denotes the floor function. We will show that $\mathcal{Z}(t,z)$ approximates $Z(W_{\lfloor td \rfloor},z)$ as $d\to \infty$ and formalize this in Theorem~\ref{thm:main:conc:St:SGD}. We thus denote $\mathcal{Z}(t,z)$ as the deterministic equivalent of $Z(W_{\lfloor td \rfloor},z)$. Similarly, it follows from Theorem~\ref{thm:main:conc:St:general:with:stop:adv} which we will introduce later that $\widehat{\mathrsfs{B}}_i(t)$, $\widehat{\mathrsfs{B}}_{i,44}(t)$, $\mathrsfs{B}_i(t)$ and $\mathfrak{m}_i(t)$ as defined in \eqref{eq:def:B(t):main} respectively correspond to the deterministic equivalents of $\widehat{B}_i(W_{\lfloor td \rfloor})$, $\|X_{\lfloor td \rfloor}\|^2$, $B_i(W_{\lfloor td \rfloor})$ and $\hat{X}_{\lfloor td \rfloor}^\top \mu_i$.

Another statistic we can extract from $Z(\widehat{W}_{\lfloor td \rfloor}, z)$ and $\mathcal{Z}(t,z)$ through contour integration is the norm of the iterates $X_{\lfloor td \rfloor}$. We extend this to $W_{\lfloor td \rfloor}$ by introducing the following pair of statistics:
\begin{equation}\label{eq:def:N(t):main}
    \begin{aligned}
        N(W_k) \defas \|W_k\|^2 = \|X_k\|^2 + \|X^{\star}\|^2 + \|\mu_1\|^2 + \|\mu_2\|^2\quad\text{and} \quad \mathrsfs{N}(t)  \defas  \oint \Tr(\mathcal{Z}(t,z)) \Dif z.
    \end{aligned}
\end{equation}
Here $\mathrsfs{N}(t) $ corresponds to the deterministic equivalent of $N(W_{\lfloor td \rfloor})$. We are now ready to present Theorem~\ref{thm:main:conc:St:SGD}.
\begin{theorem}\label{thm:main:conc:St:SGD}
    Suppose that the previous assumptions holds. Let $\vartheta_M$ be the first time either $(\widehat{\mathrsfs{B}}_1(t), \widehat{\mathrsfs{B}}_2(t))$ or $(\widehat{B}_1(W_{\lfloor td \rfloor}), \widehat{B}_2(W_{\lfloor td \rfloor}))$ exits $\mathcal{U} \times \mathcal{U}$ or that $\mathrsfs{N}(t)\geq M$ or $N(W_{\lfloor td \rfloor}) \geq M$. Then for any $T,M>0$ and any $\varepsilon  \in (0,\frac{1}{2})$, with overwhelming probability
    \[
        \sup_{0\leq t\leq T \wedge \vartheta_M}\|Z(W_{\lfloor td \rfloor}, z) - \mathcal{Z}(t,z)\|_\Gamma \leq Cd^{-\varepsilon},
    \]
    for some constant $C$ (independent of $d$) and where $\mathcal{Z}(t,z)$ is the solution to the integro-differential equation \eqref{eq:ODE_resolvent_2}.
\end{theorem}
Theorem~\ref{thm:main:conc:St:SGD} is stated in terms of the stopping time $\vartheta_M$. In Proposition~\ref{prop:main:nonexplosiveness}, we formalize a condition under which we may remove this stopping time which we prove in Section~\ref{sect:remove:stopping:time}. 

 When $K_1$ and $K_2$ are finite-dimensional matrices, to facilitate the numerical simulations of the integral-differential equation in \eqref{eq:ODE_resolvent_2} we can decompose the equation in terms of the eigenspaces of the covariance matrices $K_1$ and $K_2$ since they commute. Let $((\lambda_j^{(\ell)},\omega_j): 1 \leq j \leq d)$ be the eigenvalues and orthonormal eigenvectors of $K_\ell$ for $\ell=1,2$ and define the pairs of statistics
\begin{equation}
    \begin{gathered}
    V_j(W_k) \defas \hat{X}_k^\top \omega_j \omega_j^\top \hat{X}_k = \frac{-1}{4\pi^2}\oint_{\Gamma_j^2} S(W_{k}, z )\dif z_1 \dif z_2 \quad \text{and}\quad\mathrsfs{V}_j(t) \defas \begin{bmatrix}
        \mathrsfs{V}_{j,11}(t) & \mathrsfs{V}_{j,12}(t)  \\
        \mathrsfs{V}_{j,21}(t) & \mathrsfs{V}_{j,22}(t) 
    \end{bmatrix},
    \\
    m_{j,i}(W_k) \defas \hat{X}_k^\top \omega_j \omega_j^\top \mu_i = \frac{-1}{4\pi^2}\oint_{\Gamma_j^2} M_i(W_{k}, z )\dif z_1 \dif z_2 \quad \text{and}\quad\mathfrak{m}_{j,i}(t) \defas \begin{bmatrix}
        \mathfrak{m}_{j,i,1}(t)
        \\
        \mathfrak{m}_{j,i,2}(t)
    \end{bmatrix},
\end{gathered}
\end{equation}
where $\Gamma_j$ is a contour enclosing solely the eigenvalues $\lambda_j^{(1)}$ and $\lambda_j^{(2)}$ defined analogously to $\Gamma$. Similarly to \citep{collinswoodfin2023hitting,collinswoodfin2025Exact}, from the partial fractions decomposition of meromorphic functions, we have the representations for $i=1,2$
\begin{equation}\label{eq:S(t):def:average}
    \begin{gathered}
    \mathcal{S}(t,z) = \sum_{j=1}^d \frac{\mathrsfs{V}_j(t)}{\prod_{\ell=1}^2(\lambda^{(\ell)}_j - z_\ell)}\quad \text{and}\quad\mathcal{M}_i(t,z) = \sum_{j=1}^d \frac{\mathfrak{m}_{j,i}(t)}{\prod_{\ell=1}^2(\lambda^{(\ell)}_j - z_\ell)}.
    \end{gathered}
\end{equation}
Summing over the eigenspaces of $K_1$ and $K_2$, we also have the representations for $\ell=1,2$
\begin{equation}
    \begin{gathered}
   \mathrsfs{B}_\ell(t) = \sum_{j=1}^d \lambda_{j}^{(\ell)} \mathrsfs{V}_j(t),\quad \mathfrak{m}_{\ell} = \sum_{j=1}^d \mathfrak{m}_{j,\ell} (t) \quad\text{and}\quad \widehat{\mathrsfs{B}}_{\ell,44}(t)  =  \sum_{j=1}^d \mathrsfs{V}_{j,11}(t).
    \end{gathered}
\end{equation}
Given the relation \eqref{eq:S(t):def:average}, we obtain the $d$-dimensional set of coupled ODEs
\begin{equation}\label{eq:ODE:V_i(t):def}
    \begin{aligned}
        &\frac{\dif \mathrsfs{V}_j(t)}{\dif t} = -2\gamma(t) \sum_{i=1}^2 p_i \bigg(\lambda_j^{(i)}(\mathrsfs{V}_j(t) H_{1,i}(\widehat{\mathrsfs{B}}_i(t)) + H_{1,i}(\widehat{\mathrsfs{B}}_i(t))^\top \mathrsfs{V}_j(t) ) 
        + \mathfrak{m}_{j,i}(t) \otimes H_{2,i}(\widehat{\mathrsfs{B}}_i(t)) 
        \\
        &
        \quad
        + H_{2,i}(\widehat{\mathrsfs{B}}_i(t)) \otimes \mathfrak{m}_{j,i}(t) 
        +\left(\frac{\lambda}{2} + \partial_{44} h_i (\widehat{\mathrsfs{B}}_i(t) )\right) (\mathrsfs{V}_j(t)\operatorname{D} + \operatorname{D} \mathrsfs{V}_j(t))\bigg)
        \\
        &
        \quad
        +\frac{\gamma(t)^2}{d} \left(\sum_{i=1}^2 p_i \mathcal{I}_i(\widehat{\mathrsfs{B}}_i(t))(\lambda_j^{(i)} + (\omega_j^\top \mu_i)^2)\right) \operatorname{D},
        \\
        &
        \frac{\dif \mathfrak{m}_{j,\ell}(t)}{\dif t}  = -2\gamma(t) \sum_{i=1}^2 p_i \bigg( \lambda_j^{(i)}H_{1,i}(\widehat{\mathrsfs{B}}_i(t))^\top \mathfrak{m}_{j,\ell}(t)  + \langle \omega_j, \mu_\ell\rangle \langle \omega_j, \mu_i\rangle \cdot H_{2,i}(\widehat{\mathrsfs{B}}_i(t))
        \\
         &\quad+ \left(\frac{\lambda}{2} + \partial_{44} h_i (\widehat{\mathrsfs{B}}_i(t)) \right)  \operatorname{D}\mathfrak{m}_{j,\ell}(t)\bigg),
    \end{aligned}
\end{equation}
with initial condition $\mathrsfs{V}_j(0) = \hat{X}_0^\top \omega_j \omega_j^\top \hat{X}_0$ and $\mathfrak{m}_{j,\ell}(0) = \hat{X}_0^\top \omega_j \omega_j^\top\mu_{\ell}$. See Appendix~\ref{sec:IntegroDiffAnalysis} for more details on the derivation of the system of ODEs. Note that we show in Appendix~\ref{sec:IntegroDiffAnalysis} that $\mathrsfs{V}_j(t)$ are positive semidefinite matrices for $1\leq j \leq d$ and $t\geq 0$. It follows that $\mathrsfs{B}_\ell(t)$ is positive semidefinite, $\widehat{\mathrsfs{B}}_{\ell, 44}(t) \geq 0$ and $\mathrsfs{N}(t) \geq 0$ for all $t\geq 0$. 

Assumption~\ref{assumption:risk} entails that the ODEs \eqref{eq:ODE:V_i(t):def} have locally unique solutions whilst $(\widehat{\mathrsfs{B}}_1(t), \widehat{\mathrsfs{B}}_2(t)) \in \mathcal{U} \times \mathcal{U}$.
\begin{remark}\label{remark:equiv:sols:of:odes}
     By Lemma 2 in \citep{collinswoodfin2025Exact}, Theorem~\ref{thm:main:conc:St:SGD} can be reformulated with $\mathcal{Z}(t,z)$ constructed with $\mathcal{S}(t,z)$ and $\mathcal{M}(t,z)$ defined as in \eqref{eq:S(t):def:average}. We omit the proof as the argument is almost identical. Hence, we can recover a solution to \eqref{eq:ODE_resolvent_2} using a solution to the system of coupled ODEs \eqref{eq:ODE:V_i(t):def}.
\end{remark}
    \subsection{Adversarial Homogenized SGD}\label{sect:mainresult:AdvHSGD}
    An alternative way to compute the deterministic equivalent of $Z(W, z)$ as defined in \eqref{eq:main:def:S(W,z)} is through the use of Itô calculus. We introduce \textit{adversarial homogenized SGD} (AdvHSGD)    
\begin{equation}\label{eq:main:AdvHSGD:def}
        \dif \WHSGD_t = -\gamma(t) \nabla \mathcal{R}_{\lambda}^{\operatorname{adv}}(\WHSGD_t) \dif t + \gamma(t) \sqrt{\frac{1}{d} \sum_{i=1}^2 p_i \EE_{v, \epsilon} [  f'_i( \rho_{t, i}; \epsilon)^2  ]  (K_i + \mu_i\mu_i^\top)} \dif B_t,
        \end{equation}
        where the initial condition is given by $\mathscr{X}_0 = X_0$, $(B_t, t \ge 0)$ is a $\R^{d}$ standard Brownian motion and for $I \in \{1,2\}$ and $v\sim \mathcal{N}(0, \Id_d)$ we have
        \begin{equation}
            \begin{gathered}
        \widehat{\mathscr{X}}_t \defas \mathscr{X}_t \oplus X^{\star}\in \R^{d \times 2},\quad\mathscr{W}_t \defas \widehat{\mathscr{X}}_t \oplus \mu\in \R^{d \times 4},
        \\
        \rho_{t,I} \defas \widehat{\mathscr{X}}_t^\top(\sqrt{K_{I}}v + \mu_{I})+ \begin{bmatrix}\delta u_{t,I} \|\WHSGD_t\| & 0 \end{bmatrix}^\top,
        \\
        \text{and}\quad u_{t,I} = \operatorname{argmax}_{|u| \leq  1}f_I(\widehat{\mathscr{X}}_t^\top (\sqrt{K_{I}}v + \mu_{I}) + \begin{bmatrix} \delta u\|\mathscr{X}_t \|& 0 \end{bmatrix}^\top, \epsilon).
            \end{gathered}
        \end{equation}
In Appendix~\ref{sect:DoobDecomp:HSGD}, we perform a Doob decomposition elementwise on the function $Z(\HSGD_t, z)$ and show \eqref{eq:ODE_resolvent_2} corresponds to the expression obtained by dropping the martingale terms which arise from Itô's lemma. Through the following theorem, we see that AdvHSGD concentrates around \eqref{eq:ODE_resolvent_2} with overwhelming probability. 
\begin{theorem}\label{thm:main:conc:St:HSGD}
    Suppose that the previous assumptions holds. Let $\Theta_M$ be the first time $(\widehat{\mathrsfs{B}}_1(t), \widehat{\mathrsfs{B}}_2(t))$ exits $\mathcal{U} \times \mathcal{U}$ or that $\mathrsfs{N}(t)\geq M$. Then for any $T,M>0$ and any $\varepsilon  \in (0,\frac{1}{2})$, with overwhelming probability
    \[
        \sup_{0\leq t\leq T \wedge \Theta_M}\|Z(\HSGD_{t}, z) - \mathcal{Z}(t,z)\|_\Gamma \leq Cd^{-\varepsilon},
    \]
    for some constant $C$ (independent of $d$) and where $\mathcal{Z}(t,z)$ is the solution to the integro-differential equation \eqref{eq:ODE_resolvent_2} or the system of coupled ODEs~\eqref{eq:ODE:V_i(t):def}.
\end{theorem}
Notice that the stopping time in Theorem~\ref{thm:main:conc:St:SGD} is formulated in terms of both $(\widehat{B}_1(W_{\lfloor td \rfloor}), \widehat{B}_2(W_{\lfloor td \rfloor}))$ and $(\widehat{\mathrsfs{B}}_1(t), \widehat{\mathrsfs{B}}_2(t))$ but the stopping time in Theorem~\ref{thm:main:conc:St:HSGD} in terms of $(\widehat{\mathrsfs{B}}_1(t), \widehat{\mathrsfs{B}}_2(t))$. 

 It turns out that $\mathcal{Z}(t,z)$ remaining bounded (i.e. $\mathrsfs{N}(t) < M$) and having sufficiently regular paths (i.e. $(\widehat{\mathrsfs{B}}_1(t), \widehat{\mathrsfs{B}}_2(t)) \in \mathcal{U} \times \mathcal{U}$) are sufficient conditions for Theorem~\ref{thm:main:conc:St:SGD} and ~\ref{thm:main:conc:St:HSGD} to hold. It is also possible to reverse the roles in Theorem~\ref{thm:main:conc:St:SGD}~and~\ref{thm:main:conc:St:HSGD} using the conditions $(\widehat{B}_1(W_{\lfloor td \rfloor}), \widehat{B}_2(W_{\lfloor td \rfloor})) \in \mathcal{U} \times \mathcal{U} $ and $N(W_{\lfloor td \rfloor}) < M$ or $(\widehat{B}_1(\HSGD_{\lfloor td \rfloor}), \widehat{B}_2(\HSGD_{\lfloor td \rfloor})) \in \mathcal{U}\times \mathcal{U}$ and $N(\HSGD_{t}) < M$ respectively. For more details and a stronger result, see Section~\ref{sect:remove:stopping:time}, specifically Theorem~\ref{thm:main:S(W):oneST}, Lemma~\ref{lem:prop:normGamma:normW} and Remark~\ref{remark:equiv:stopping:time:N(t)}.
 \\

 The minimizer of the adversarial population risk is defined as follows.
\begin{definition}\label{def:adversarial:ground:truth}
    Let $X^{\star,\operatorname{adv}}$ denote the adversarial ground truth. In other words
    \begin{equation}
        X^{\star,\operatorname{adv}} \in \underset{X}{\operatorname{argmin}}  \left\{\mathcal{R}^{\operatorname{adv}}(X) =\EE_{a,I,\epsilon} [\max_{|s| \leq 1}f_I( X^\top a  + \delta s \|X\|, \epsilon) ]\right\}.
    \end{equation}
\end{definition}
In some cases such as least squares \citep{xing2021same_minimizer}, there exists a threshold $\delta_0$ such that for any $\delta < \delta_0$, $X^{\star,\operatorname{adv}} = X^\star$. We will assume that the norm of $X^{\star,\Adv}$ is bounded by a constant independent of dimension.
\begin{assumption} \label{assumption:norm:Xstar:adv}
    The adversarial ground truth $X^{\star,\Adv}$ satisfies $\| X^{\star, \Adv} \|  \le C$ for some constant $C$ independent of $d$. 
\end{assumption}

 Given that we are also interested in statistics involving the adversarial ground truth $X^{\star,\operatorname{adv}}$ such as the descent to adversarial optimality $\|X- X^{\star,\operatorname{adv}}\|^2$, we introduce $\widehat{W} = W \oplus X^{\star,\Adv}$. To recover statistics of the form $\widehat{W}^\top q(K_1,K_2)\widehat{W}$, we must also keep track of the statistic $A_{\mathcal{K}}(\widehat{W},z) \defas X^\top \mathscr{R}(z;\mathcal{K})X^{\star,\operatorname{adv}}$. 
Hence, analogously to $Z$ defined in \eqref{eq:def:og:main:Z}, we introduce $Z^\Adv$ as follows
\begin{equation}\label{eq:def:Zadv:main}
    Z^\Adv(\widehat{W},z) \defas \widehat{W}^\top \mathscr{R}(z; \mathcal{K})\widehat{W} = \begin{bmatrix}
        Z(W,z) & (\star)_0\\
        (\star)_0^\top & (X^{\star, \Adv})^\top\mathscr{R}(z;\mathcal{K})(X^{\star, \Adv})
    \end{bmatrix},
\end{equation}
with $A_{\mathscr{K}}(\widehat{W},z) \defas X^\top\mathscr{R}(z;\mathcal{K}) X^{\star, \Adv}$ and where $(\star)_0$ corresponds to 
\[ 
    (\star)_0 \defas  \begin{bmatrix}
       A_{\mathscr{K}}(\widehat{W},z) & (X^{\star})^\top \mathscr{R}(z;\mathcal{K}) X^{\star, \Adv} & \mu_1^\top \mathscr{R}(z;\mathcal{K}) X^{\star, \Adv} & \mu_2^\top \mathscr{R}(z;\mathcal{K}) X^{\star, \Adv} 
    \end{bmatrix}^\top.
\]
We also require $ \widehat{\mathrsfs{B}}_i^{\Adv}(t)$ which we define similarly to $ \widehat{\mathrsfs{B}}_i(t)$:
    \begin{equation}\label{eq:def:B:hat:adv:detequiv}
        \widehat{\mathrsfs{B}}_i^{\Adv}(t) \defas \begin{bmatrix}
            \widehat{\mathrsfs{B}}_i(t) &(\star) \\
            (\star)^\top & \|X^{\star,\Adv}\|^2
         \end{bmatrix} \in \R^{5 \times 5},  
        \end{equation}
    with $(\star)$ given by 
    \[
        (\star) \defas  \begin{bmatrix}
            \oint z_i  \mathcal{A}_{\mathscr{K}}(t,z) \Dif z & (X^{\star})^\top  (X^{\star, \Adv}) & \mu_i^\top X^{\star, \Adv} &0
        \end{bmatrix}^\top.
    \]
Similarly to $\mathcal{U}$, we also introduce the convex set $\mathcal{U}^\Adv \subset \R^{5\times 5}$ defined below
    \begin{equation}\label{eq:Uadv:set}
        \mathcal{U}^\Adv \subseteq \left\{\begin{bmatrix}
        B & u_1 \\
        u_2^\top & w
        \end{bmatrix}:\quad B\in \mathcal{U} ,\quad u_1,u_2 \in \R^4\quad\text{and}\quad  w\in \R\right\}.
    \end{equation}
Let $\mathcal{A}_{\mathcal{K}}(t,z)$ denote the deterministic equivalent of $ A_{\mathscr{K}}(\widehat{W}_{td},z) = X_{td}^\top \mathscr{R}(z;\mathcal{K}) X^{\star, \Adv}$. Similarly to \eqref{eq:S(t):def:average}, when $K_1$ and $K_2$ are finite-dimensional matrices, we can decompose \eqref{eq:det:equiv:align:Xstar:Adv} in terms of the eigenspaces of the covariance matrices $K_1$ and $K_2$. Define the pair of statistics
\begin{equation}
    \begin{gathered}
    U_{j}(\widehat{W}_k) \defas X_k^\top \omega_j \omega_j^\top X^{\star, \Adv} = \frac{-1}{4\pi^2}\oint_{\Gamma_j^2} A_{\mathcal{K}}(\widehat{W}_k, z)\dif z_1 \dif z_2 \quad \text{and its deterministic equiv. }   \mathrsfs{U}_{j}(t),
\end{gathered}
\end{equation}
where $\Gamma_i$ is a contour enclosing solely the eigenvalues $\lambda_j^{(1)}$ and $\lambda_j^{(2)}$. From the partial fractions decomposition of meromorphic functions, we have the representation 
\begin{equation}\label{eq:align_K(t):def:average}
    \begin{gathered}
       \mathcal{A}_{\mathcal{K}}(t,z) = \sum_{j=1}^d \frac{\mathrsfs{U}_{j}(t)}{\prod_{\ell=1}^2(\lambda^{(\ell)}_j - z_\ell)}.
    \end{gathered}
\end{equation}
From Cauchy's integral formula, we have $\mathrsfs{A}_{K_i}(t) = \oint z_i  \mathcal{A}_{\mathcal{K}}(t,z) \Dif z$ and summing over the eigenspaces of $K_1$ and $K_2$, for $i=1,2$ we obtain $\mathrsfs{A}_{K_i}(t) = \sum_{j=1} \lambda_j^{(i)} \mathrsfs{U}_{j}(t)$ where $\mathrsfs{A}_{K_i}(t)$ denotes the deterministic equivalent of $X_{td}^\top K_i X^{\star, \Adv}$. Thus, we obtain the coupled ODEs
\begin{equation}\label{eq:ODE:align_K(t):def:main}
    \begin{aligned}
        &
        \frac{\dif \mathrsfs{U}_{j}(t)}{\dif t}  = -2\gamma(t) \sum_{i=1}^2 p_i \bigg( \left(\lambda_j^{(i)} \partial_{11} h_i(\widehat{\mathrsfs{B}}_i(t)) + \left(\frac{\lambda}{2} + \partial_{44} h_i (\widehat{\mathrsfs{B}}_i(t))\right)\right) \mathrsfs{U}_{j}(t)  
        \\
         &\quad+\left(\lambda_j \partial_{12}h_i (\widehat{\mathrsfs{B}}_i(t)) \langle X^{\star }, \omega_j\rangle+ \partial_{13} h_i(\widehat{\mathrsfs{B}}_i(t)) \langle \mu_i ,\omega_j \rangle \right) \langle X^{\star, \Adv},  \omega_j\rangle 
         \bigg),
    \end{aligned}
\end{equation}
with initial condition $ \mathrsfs{U}_{j}(0) = X_0^\top \omega_j \omega_j^\top X^{\star, \Adv}$. As in Remark~\ref{remark:equiv:sols:of:odes}, using an almost identical argument to Lemma 2 in \citep{collinswoodfin2025Exact}, we can recover $\mathcal{Z}^\Adv(t,z)$ from a solution to the system of coupled ODEs \eqref{eq:ODE:V_i(t):def} and \eqref{eq:ODE:align_K(t):def:main}. See \eqref{eq:ODE_resolvent_2:align} for the integro-differential equation of $\mathcal{Z}^\Adv(t,z)$ analogous to \eqref{eq:ODE_resolvent_2} in Section~\ref{sect:general:adv:concentration}. We then obtain a parallel result to Theorems~\ref{thm:main:conc:St:SGD} and \ref{thm:main:conc:St:HSGD}.
\begin{theorem}\label{thm:main:conc:St:adv}
    Suppose that the previous assumptions holds. Let $\Theta_M^\Adv$ be the first time $(\widehat{\mathrsfs{B}}^\Adv_1(t), \widehat{\mathrsfs{B}}^\Adv_2(t))$ exits $\mathcal{U}^\Adv \times \mathcal{U}^\Adv$ or that $\mathrsfs{N}(t) + \|X^{\star,\Adv}\|^2\geq M$. Then for any $T,M>0$ and any $\varepsilon  \in (0,\frac{1}{2})$, with overwhelming probability
    \[
        \sup_{0\leq t\leq T \wedge \Theta_M}\|Z^\Adv(\widehat{\HSGD}_{t}, z) - \mathcal{Z}^\Adv(t,z)\|_\Gamma + \|Z^\Adv(\widehat{W}_{\lfloor td \rfloor}, z) - \mathcal{Z}^\Adv(t,z)\|_\Gamma \leq Cd^{-\varepsilon},
    \]
    for some constant $C>0$ (independent of $d$) and where $\mathcal{Z}^\Adv(t,z)$ is the solution to the integro-differential equation \eqref{eq:ODE_resolvent_2:align} or the system of coupled ODEs~\eqref{eq:ODE:V_i(t):def} and \eqref{eq:ODE:align_K(t):def:main}.
\end{theorem}
\subsection{Deterministic Dynamics for General Statistics }\label{sect:mainresult:SGD:generalstats}
As mentioned previously, we can recover statistics of the form $\widehat{W}^\top q(K_1,K_2)\widehat{W}$ for a polynomial $q$ from Cauchy's integral formula. It is then natural to expect Theorems~\ref{thm:main:conc:St:SGD},~\ref{thm:main:conc:St:HSGD}~and~\ref{thm:main:conc:St:adv} extend to a more general set of polynomial functions of $K_1$ and $K_2$.
\begin{assumption}\label{assumption:smooth:stats:Adv}
    Denote $\widetilde{X} = X \oplus X^{\star} \oplus X^{\star, \Adv}$ and let $q$ be a polynomial. Define 
    \[
        Q_q(\widetilde X,\mu)
        \defas 
        \begin{bmatrix}
            \widetilde X^\top q(K_1,K_2)\widetilde X 
            & \widetilde X^\top \mu \\
            \mu^\top \widetilde X 
            & \mu^\top \mu
        \end{bmatrix}
        \in \mathbb R^{5\times 5}.
    \]
    Suppose there exists an open convex set
    $\overline{\mathcal U}\subset \mathbb R^{5\times 5} \times \R$ such that $ (\widehat{B}_1(W), \widehat{B}_2(W))\in \mathcal{U} \times \mathcal{U}$ implies $(Q_q(\widetilde X,\mu), \|X\|^2)\in \overline{\mathcal U}$ where $\mathcal U$ is the admissible set from
    Assumption~\ref{assumption:risk}. The statistic $\varphi$ satisfies the composite structure,
   \[
        \varphi(X) = g\left(\begin{bmatrix}
            Q_q(\tilde{X})& 0 \\
             0 &  \|X\|^2
         \end{bmatrix}\right),
    \]
    where $g: \R^{6\times 6} \to \R$ is $\alpha$-pseudo-Lipschitz on $\overline{ \mathcal{U}}^\Adv $. Here $\overline{\mathcal{U}}^\Adv \subset \R^{6\times 6}$ is defined as follows
    \begin{equation}
        \overline{\mathcal{U}}^\Adv \subseteq \left\{\begin{bmatrix}
        A & u_1 \\
        u_2^\top & w
        \end{bmatrix}:\quad (A, w) \in \overline{\mathcal{U}} \quad \text{and}\quad u_1,u_2 \in \R^5\quad\right\}.
    \end{equation}
\end{assumption}
It is easy to see that $\widehat{B}_i(W)$ and $h_i$ hence $\mathcal{R}_{\lambda}^{{\operatorname{adv}}}$ satisfy Assumption~\ref{assumption:smooth:stats}. Other examples of important statistics which satisfy this assumption are $\|X\|^2$, $\|X-X^{\star}\|^2$ and $\|X-X^{\star, \Adv}\|^2$. For statistics satisfying Assumption~\ref{assumption:smooth:stats:Adv}, the deterministic equivalent is defined as follows
\begin{equation}\label{eq:main:def:phi(t):adv}
    \phi^\Adv(t) =g\left( \begin{bmatrix}
        G_{1}  & G_{2}^\top & 0  \\
        G_{2} & \mu^\top \mu & 0 \\
        0 & 0 & G_{3}
    \end{bmatrix}\right),
\end{equation}
where $ G_{2} = \oint \mathcal{M}(t,z)\Dif z \oplus (\mu^\top X^{\star,\Adv})\in \R^{2\times 3}$, $G_{3}= \oint\mathcal{S}_{11}(t,z)\Dif z$ and 
\begin{equation*}
    \begin{aligned}
        G_{1}  =  \begin{bmatrix}
            \oint q(z)\mathcal{S}_{11}(t,z)\Dif z & \oint q(z)\mathcal{S}_{12}(t,z)\Dif z & \oint q(z)\mathcal{A}_{\mathscr{K}}(t,z)\Dif z \\
            \oint q(z)\mathcal{S}_{21}(t,z)\Dif z & X^{\star \top}q(K_1,K_2)X^{\star} & X^{\star\top}q(K_1,K_2)X^{\star, \Adv}  \\
            \oint q(z)\mathcal{A}_{\mathscr{K}}(t,z)\Dif z &  X^{\star\top}q(K_1,K_2)X^{\star, \Adv} & X^{\star, \Adv \top}q(K_1,K_2)X^{\star, \Adv}
        \end{bmatrix},
    \end{aligned}
\end{equation*}
or where $G_{3}= \sum_{j=1}^d \mathrsfs{V}_{j,11}(t)$, 
\begin{equation*}
    G_{2}= \begin{bmatrix} \sum_{j=1}^d \mathfrak{m}_{j,1}(t) & \sum_{j=1}^d \mathfrak{m}_{j,2}(t) \\
                            (X^{\star, \Adv})^\top  \mu_1 & (X^{\star, \Adv})^\top  \mu_2
    \end{bmatrix},
\end{equation*} 
and
\begin{equation*}
    \begin{aligned}
        G_{i,1}  =  \begin{bmatrix}
            \sum_{j=1}^d \mathrsfs{V}_{j,11}(t) q(\lambda_j^{(1)},\lambda_j^{(2)} ) &  \sum_{j=1}^d \mathrsfs{V}_{j,12}(t) q(\lambda_j^{(1)},\lambda_j^{(2)} ) &   \sum_{j=1} \mathrsfs{U}_{j}q(\lambda_j^{(1)},\lambda_j^{(2)} )(t) \\
            \sum_{j=1}^d \mathrsfs{V}_{j,21}(t) q(\lambda_j^{(1)},\lambda_j^{(2)} ) & X^{\star \top}q(K_1,K_2)X^{\star} & X^{\star\top}q(K_1,K_2)X^{\star, \Adv}  \\
            \sum_{j=1} \mathrsfs{U}_{j}(t) q(\lambda_j^{(1)},\lambda_j^{(2)} )&  X^{\star\top}q(K_1,K_2)X^{\star, \Adv} & X^{\star, \Adv \top}q(K_1,K_2)X^{\star, \Adv}
        \end{bmatrix}.
    \end{aligned}
\end{equation*}
    We may finally present our main result. 
\begin{theorem}\label{thm:main:conc:St:general:with:stop:adv}
    Suppose that the previous assumptions holds. Let $\Theta_M^\Adv$ be the first time $(\widehat{\mathrsfs{B}}^\Adv_1(t), \widehat{\mathrsfs{B}}^\Adv_2(t))$ exits $\mathcal{U}^\Adv \times \mathcal{U}^\Adv$ or that $\mathrsfs{N}(t) + \|X^{\star, \Adv}\|^2\geq M$. For any function satisfying Assumption~\ref{assumption:smooth:stats:Adv}, for any $T,M>0$ and any $\varepsilon  \in (0,\frac{1}{2})$, with overwhelming probability
    \[
        \sup_{0\leq t\leq T \wedge \Theta_M}|\varphi(\WHSGD_{t}) - \varphi(X_{\lfloor td \rfloor})|+ |\varphi(\WHSGD_{t})- \phi(t)| \leq Cd^{-\varepsilon},
    \]
    for some constant $C$ (independent of $d$) and where $\phi(t)$ is defined in \eqref{eq:main:def:phi(t):adv}.
\end{theorem}
Note that Section~\ref{sect:conc:proof:any:stat} contains an analogous result for statistics not involving $X^{\star, \Adv}$ in the case that Assumption~\ref{assumption:norm:Xstar:adv} fails. Provided $(\widehat{\mathrsfs{B}}^\Adv_1(t), \widehat{\mathrsfs{B}}^\Adv_2(t))$, $(\widehat{B}_1^\Adv(\widehat{W}_{\lfloor td \rfloor}), \widehat{B}^\Adv_2(\widehat{W}_{\lfloor td \rfloor}))$ or $(\widehat{B}^\Adv_1(\widehat{\HSGD}_{t}), \widehat{B}^\Adv_2(\widehat{\HSGD}_{t}))$ stays in $\mathcal{U}^\Adv \times \mathcal{U}^\Adv$, we formalize a condition under which it is possible to remove the stopping time $\Theta_M^\Adv$. 
\begin{proposition}[Non-explosiveness]\label{prop:main:nonexplosiveness}
    Suppose that the assumptions stated in this framework hold and suppose that the objective function $f$ is $\alpha$-pseudo-Lipschitz with $\alpha\leq 1$. Then there exists a constant $C$ such that 
    \[
    \mathrsfs{N}(t) \leq (1+ \mathrsfs{N}(0))e^{C t},
    \]
    for $C=C (L(f), \eta, \delta, \|\mathcal{K}\|_{\operatorname{op}}, \max_{i=1,2}\|\mu_i\|, \lambda, \bar{\gamma})$ and all time $t$ such that $(\widehat{\mathrsfs{B}}_1^\Adv(t), \widehat{\mathrsfs{B}}_2^\Adv(t))$ is in $\mathcal{U}^\Adv\times \mathcal{U}^\Adv$.
    \end{proposition}
For more details, see Corollary~\ref{cor:bounded_iterates:adv} and Proposition~\ref{prop:nonexplosiveness}. Hence, from Proposition~\ref{prop:main:nonexplosiveness}, we obtain a simplication of Theorem~\ref{thm:main:conc:St:general:with:stop:adv}.
\begin{corollary}
    Suppose that the previous assumptions holds. Suppose $\mathcal{U}^\Adv = \R^{5\times 5}$ and that $f$ is $\alpha$-pseudo-Lipschitz with $\alpha \leq 1$. For any function satisfying Assumption~\ref{assumption:smooth:stats:Adv}, for any $T>0$ and any $\varepsilon  \in (0,\frac{1}{2})$, with overwhelming probability
    \[
        \sup_{0\leq t\leq T }|\varphi(\WHSGD_{t}) - \varphi(X_{\lfloor td \rfloor})|+ |\varphi(\WHSGD_{t})- \phi(t)| \leq Cd^{-\varepsilon},
    \]
    for some constant $C>0$ (independent of $d$) and where $\phi(t)$ is defined in \eqref{eq:main:def:phi(t):adv}.
\end{corollary}

\section{$\ell_2$-Adversarial Training v.s. $\ell_2$-Regularization}\label{sect:AdvHSGD:vs:HSGD}
From \eqref{sect:mainresult:AdvHSGD}, AdvHSGD for $\ell_2$-adversarial least squares with $\lambda= 0$ on a single data class $a \sim N(0, K)$ (i.e. $p_1 = 1$, $K_1 = K$, $\mu_1 = 0$) with soft labels is the solution to the SDE 
\begin{equation}\label{eq:advhsgd}
    \dif \WHSGD_t = -\gamma^{\operatorname{eff}}(t)\left(\nabla \mathcal{R}(\WHSGD_t) + \lambda^{\operatorname{eff}}(t)\WHSGD_t\right) \dif t + \gamma(t)\sqrt{\frac{2}{d} \mathcal{R^\Adv}(\WHSGD_t) K} \dif B_t,
\end{equation}
where we define 
\[
    \gamma^{\operatorname{eff}}(t) = \gamma(t)\bigg(1 + \delta \sqrt{\frac{2}{\pi}}\frac{\|\WHSGD_t\|}{\sqrt{2\mathcal{R}(\WHSGD_t)}}\bigg)\quad \text{and}\quad  \lambda^{\operatorname{eff}}(t) = \frac{\delta^2 + \delta \sqrt{\frac{2}{\pi}}\frac{\sqrt{2\mathcal{R}(\WHSGD_t)}}{\|\WHSGD_t\|}}{1 + \delta \sqrt{\frac{2}{\pi}}\frac{\|\WHSGD_t\|}{\sqrt{2\mathcal{R}(\WHSGD_t)}}}.
\]
 Here $\gamma^{\operatorname{eff}}(t)$ and $\lambda^{\operatorname{eff}}(t)$ correspond to the effective learning rate and regularization induced by adversarial training. \paragraph{Effective learning rate $\gamma^{\operatorname{eff}}$ and regularization $\lambda^{\operatorname{eff}}$}
The schedules of $\gamma^{\operatorname{eff}}$ and $\lambda^{\operatorname{eff}}$ are entirely determined by the path of the ratio $q(\WHSGD_t)\defas \left(\|\WHSGD_t\|^2/ (2\mathcal{R}(\WHSGD_t))\right)^{1/2}$. The effective learning rate is a strictly increasing function of $q(\WHSGD_t)$, whereas the effective regularization parameter is a strictly decreasing function of $q(\WHSGD_t)$. We may reexpress AdvHSGD \eqref{eq:AdvHSGD:app} as follows
\begin{equation}\label{eq:AdvHSGD:app:reformulation}
    \dif \WHSGD_t = -\gamma^{\operatorname{eff}}(t)\left(\nabla \mathcal{R}(\WHSGD_t) + \lambda^{\operatorname{eff}}(t)\WHSGD_t\right) \dif t + \gamma(t) \sqrt{\frac{2}{d} G(q(\WHSGD_t))\mathcal{R}(\WHSGD_t) K} \dif B_t,
\end{equation}
where $G(q) \defas 1+ 2\delta \sqrt{\frac{2}{\pi}} q + \delta^2 q^2$. Consider $\gamma(t) \equiv \gamma$ for some $\gamma >0$. When $q(\WHSGD_t)$ is large (i.e. $\mathcal{R}(\WHSGD_t) \ll \|\WHSGD_t\|^2 $), we have 
\[
    \gamma^{\operatorname{eff}}(t)\lambda^{\operatorname{eff}}(t) = \gamma\left(\delta^2 + \delta \sqrt{\frac{2}{\pi}} \frac{1}{q(\WHSGD_t)}\right) \approx \gamma \delta^2,
\]
and $\lambda^{\operatorname{eff}}(t) \to 0 $ as $q(\WHSGD_t) \to \infty$. Hence, when $q(\WHSGD_t)$ is large, $\gamma^{\operatorname{eff}}(t)$ grows linearly in $q(\WHSGD_t)$ which implies the descent term of \eqref{eq:AdvHSGD:app:reformulation} favors the least squares component relative to the $\ell_2$-regularization component and acts similarly to least squares with learning rate $\gamma^{\operatorname{eff}}(t)$. When $q(\WHSGD_t) \ll 1$ (i.e. $\|\WHSGD_t\|^2 \ll \mathcal{R}(\WHSGD_t)$), the effective learning rate $\gamma^{\operatorname{eff}}(t) \approx \gamma$ and
\[
     \gamma^{\operatorname{eff}}(t)\lambda^{\operatorname{eff}}(t) = \gamma\left(\delta^2 + \delta \sqrt{\frac{2}{\pi}} \frac{1}{q(\WHSGD_t)}\right) \approx \gamma \delta \sqrt{\frac{2}{\pi}} \frac{1}{q(\WHSGD_t)}.
\]
Thus, for small $q(\WHSGD_t) \ll 1$, adversarial training favors the $\ell_2$-regularization component and acts similarly to least squares with constant learning rate $\gamma$ and large regularization parameter. This is further emphasized by the following decomposition of $\nabla \mathcal{R}^{\Adv}$:
\begin{equation}
    \nabla \mathcal{R}^{\operatorname{adv}}(\WHSGD_t) =   \nabla \mathcal{R}(\WHSGD_t) + \delta^2 \WHSGD_t + \delta \sqrt{\frac{2}{\pi }} \left( q(\WHSGD_t)  \nabla \mathcal{R}(\WHSGD_t) + \frac{1}{ q(\WHSGD_t) }   \WHSGD_t\right).
\end{equation}
Furthermore, $q(\WHSGD_t) \ll 1$ implies $G(q(\WHSGD_t))\approx 1$ from which it follows that the diffusion coefficient of \eqref{eq:AdvHSGD:app:reformulation} acts similarly to the diffusion coefficient of Homogenized SGD \citep{collinswoodfin2023hitting, CollinsWoodfinPaquette01, PPAP01}:
\begin{equation}
    \gamma(t) \sqrt{\frac{2}{d} G(q(\WHSGD_t))\mathcal{R}(\WHSGD_t) K} \dif B_t \approx  \gamma(t) \sqrt{\frac{2}{d} \mathcal{R}(\WHSGD_t) K} \dif B_t.
\end{equation}
\paragraph{Relation to ridge regularization}
It has been shown that, under suitable assumptions on the data and adversarial radius $\delta$, the minimizer of $\ell_2$-adversarial least squares can coincide with the ridge regression solution \citep{ribeiro2023regularizationpropertiesadversariallytrainedlinear}. We emphasize this connection by reexpressing the drift term of AdvHSGD \eqref{eq:AdvHSGD:app} as follows 
\begin{equation}
   -\gamma(t) \nabla \mathcal{R}^{\operatorname{adv}}(\WHSGD_t) = -\gamma^{\operatorname{eff}}(t)\left(\left(K+\lambda^{\operatorname{eff}}(t) \Id_d\right) \WHSGD_t - K X^\star\right).
\end{equation}
Conditional on the parameters $\gamma^{\operatorname{eff}}(t)$ and $\lambda^{\operatorname{eff}}(t)$, the drift is that of ridge regression with adaptive learning rate and regularization. This connection becomes explicit by projecting onto the eigenbasis of $K$. Recall $(\lambda_j, \omega_j)$ the eigenvalue-eigenvectors pairs of $K$. Then by Itô's formula
\begin{equation}
    \dif\langle \WHSGD_t, \omega_j\rangle= -\gamma^{\operatorname{eff}}(t)\left(\left(\lambda_j +\lambda^{\operatorname{eff}}(t)\right)\langle \WHSGD_t, \omega_j\rangle- \lambda_j \langle X^\star, \omega_j\rangle \right)\dif t+ \dif  \mathcal{M}_j(\WHSGD_t),
\end{equation}
where $\dif \mathcal{M}_j(\WHSGD_t)$ is a martingale term associated with eigenvector $\omega_j$ which vanishes as $d\to \infty$ (see Theorem~\ref{thm:main:conc:St:general:with:stop:adv}). Let $\rho_j(t)$, $\tilde{\lambda}^{\operatorname{eff}}(t)$ $q(t)$ denote the deterministic equivalents of $\varphi(\WHSGD_t)$, $\lambda^{\operatorname{eff}}(t)$ and $q(\WHSGD_t)$ in the sense of Theorem~\ref{thm:main:conc:St:general:with:stop:adv}. In Proposition~\ref{prop:fixed:point:main}, assuming $q(t) \to q_\infty$, we show that as $t\to\infty$
\begin{equation}
    \rho_j(t) \to \frac{\lambda_j}{\lambda_j + \tilde{\lambda}^{\operatorname{eff}}(q_\infty)} \langle X^\star, \omega_j\rangle.
\end{equation}
Hence, whenever $q(t) \underset{t \to \infty}{\longrightarrow} q_\infty$ and as $d\to \infty$, SGD  converges to the ridge optimum associated with $ \tilde{\lambda}^{\operatorname{eff}}(q_\infty)$:
\[
    X_{\infty} = \left(K +  \tilde{\lambda}^{\operatorname{eff}}(q_\infty)\right)^{-1}KX^\star.
\] 
In particular, directions with $\lambda_j\gg  \tilde{\lambda}^{\operatorname{eff}}(q_\infty)$ are nearly identical to $X^\star$, whereas directions with $\lambda_j\ll \tilde{\lambda}^{\operatorname{eff}}(q_\infty)$ are strongly attenuated. \\

\paragraph{Relation to Homogenized SGD}
AdvHSGD has a similar structure to the SDE called Homogenized SGD (HSGD) \citep{collinswoodfin2023hitting, CollinsWoodfinPaquette01, PPAP01} which approximates the dynamics of statistics of the iterates of SGD on $\ell_2$-regularized least squares with fixed parameter $\lambda^{\operatorname{Reg}} >0$. Let $\lambda^{\operatorname{Reg}}$ vary with time such that we define $\lambda^{\operatorname{Reg}}(t)$ and $\gamma^{\operatorname{Reg}}(t)$ analogously to $\lambda^{\operatorname{eff}}(t)$ and $\gamma^{\operatorname{eff}}(t)$. Then the SDE for $\ell_2$-regularized least squares with adaptive learning rate $\gamma^{\operatorname{Reg}}(t)$ amounts to
\begin{equation}\label{eq:HSGD}
    \begin{aligned}
    \dif \WHSGD_t^{\operatorname{Reg}} &= -\gamma^{\operatorname{Reg}}(t) \left(\nabla \mathcal{R}(\WHSGD_t^{\operatorname{Reg}}) + \lambda^{\operatorname{Reg}}(t) \WHSGD_t^{\operatorname{Reg}} \right)\dif t 
    \\
    &
    \qquad + \gamma(t) \sqrt{\frac{2}{d}\bigg( \mathcal{R}(\WHSGD_t^{\operatorname{Reg}})+ \delta \sqrt{\frac{2}{\pi}} \|\WHSGD_t^{\operatorname{Reg}}\| \sqrt{2 \mathcal{R}(\WHSGD_t^{\operatorname{Reg}})}  + \frac{\delta^2}{\pi}\|\WHSGD_t^{\operatorname{Reg}}\|^2\bigg)K} \dif B_t.
    \end{aligned}
\end{equation}
Using similar methods to \citep{collinswoodfin2024highline, CollinsWoodfinPaquette01,collinswoodfin2023hitting}, it would be possible to show that statistics of the SDE \eqref{eq:LRHSGD:reg:adaptive} approximate the dynamics of SGD as $d\to \infty$ on $\ell_2$-regularized least squares with adaptive learning rate $ \gamma^{\operatorname{Reg}}(\lfloor td \rfloor)$ and regularization $\lambda^{\operatorname{Reg}}(\lfloor td \rfloor) $. See Figure~\ref{fig:HSGD} for numerical evidence supporting this claim. See Appendix~\ref{sect:compareADVHSGD:to:HSGD} for a justification behind the derivation of the deterministic equivalents \eqref{eq:R:reg:det:equiv:volterra} and \eqref{eq:l2norm2:reg:det:equiv:volterra} for HSGD. 

\begin{figure}[t!]
    \centering
    \includegraphics[width=0.47\textwidth]{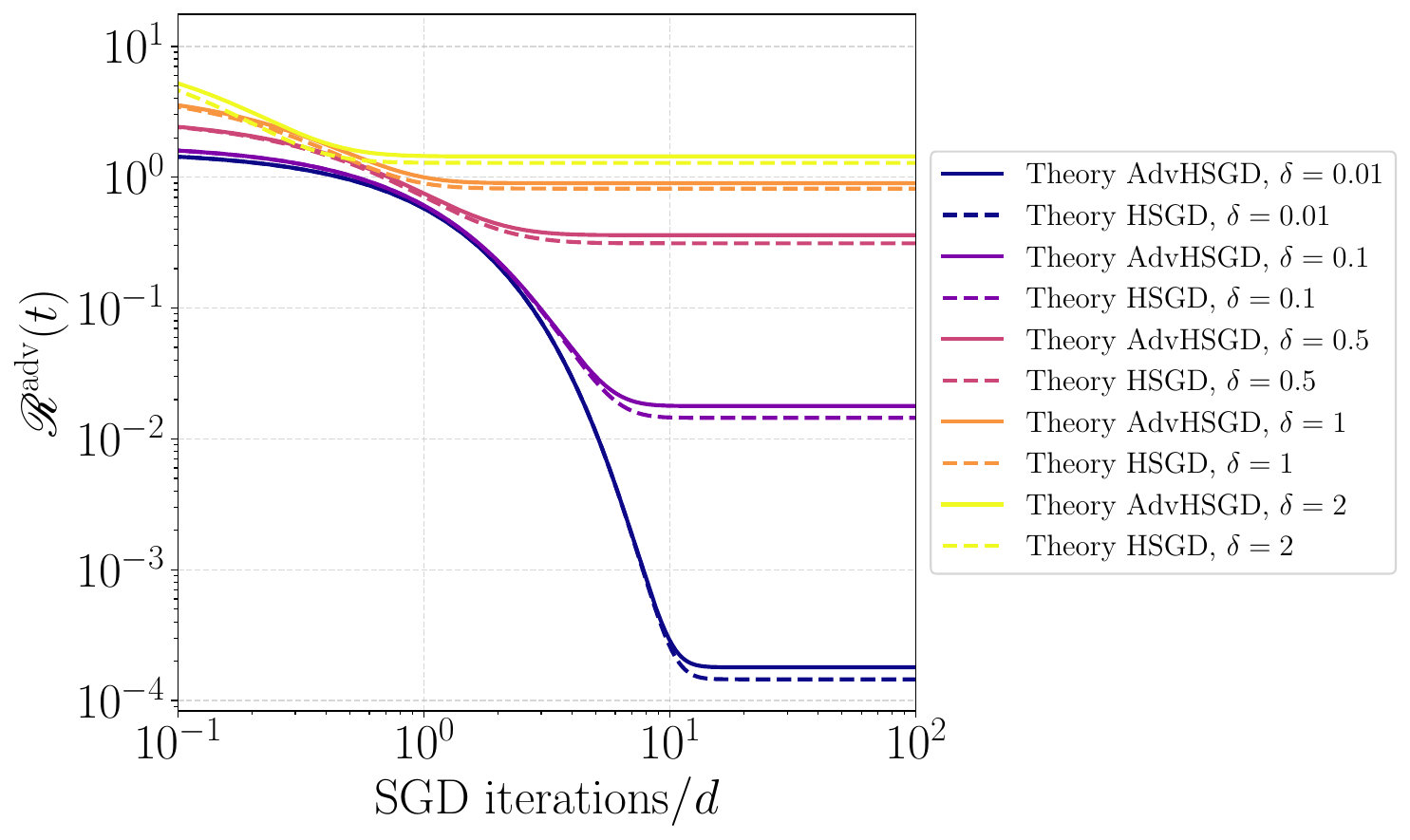}
    \includegraphics[width=0.47\textwidth]{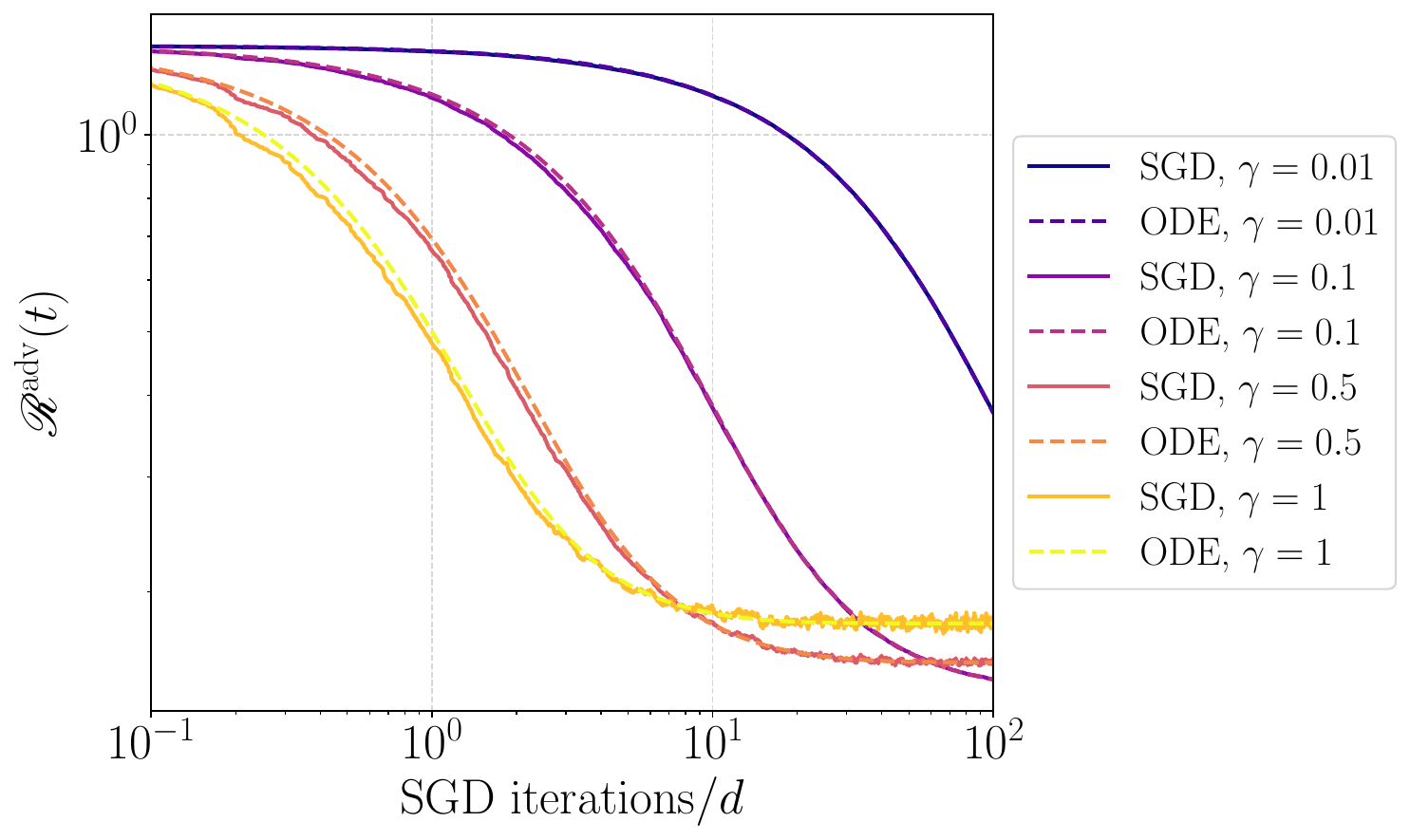}
    \caption{\textbf{Comparison of $\ell_2$-adversarial least squares and $\ell_2$-regularized least squares with adaptive learning rate and regularization.} The left plot compares the paths of the deterministic equivalents of $\mathcal{R}^{\Adv}$ for AdvHSGD and HSGD and confirms Proposition~\ref{prop:adv:reg:volterra:stability}. The right plot compares the path of $\mathcal{R}^{\Adv}(X_k)$ for SGD with adaptive learning rate $\gamma^{\operatorname{Reg}}(t)$ and regularization $\lambda^{\operatorname{Reg}}(t)$ versus the deterministic equivalent computed from \eqref{eq:R:reg:det:equiv:volterra} and \eqref{eq:l2norm2:reg:det:equiv:volterra} and shows close agreement between them for a variety of constant $\gamma(t) \equiv \gamma$. The $\ell_2$-adversarial risk is computed using $\mathcal{R}^{\Adv}(X_k) = \mathcal{R}(X_k) + \delta \sqrt{\frac{2}{\pi}} \sqrt{2\mathcal{R}(X_k)}\|X_k\| + \frac{\delta^2}{2}\|X_k\|^2$. See Appendix~\ref{app:captions} for simulation details.}
    \label{fig:HSGD}
\end{figure}

The drift terms of \eqref{eq:advhsgd} and \eqref{eq:HSGD} are identical and their diffusion terms differ up to $\frac{\delta^2}{2}$ versus $\frac{\delta^2}{\pi}$. Given the same initialization, it turns out that the paths of the deterministic equivalents of AdvHSGD and HSGD are very similar. See Figure~\ref{fig:HSGD} for numerical experiments comparing the paths of their deterministic equivalents for a variety of $\delta$.
\begin{proposition}
    \label{prop:adv:reg:volterra:stability}
    Let $\mathscr{R}(t)$ and $\Phi(t)$ be the analogous deterministic equivalents of $\mathcal{R}(X_t^{\operatorname{Reg}})$ and $\|X_t^{\operatorname{Reg}}\|^2$ for HSGD with adaptive learning rate $ \gamma^{\operatorname{Reg}}(t)$ and regularization $\lambda^{\operatorname{Reg}}(t) $. Assume $\mathrsfs{R}(0)=\mathscr{R}(0)$ and $\widehat{\mathrsfs{B}}_{44}(0)=\Phi(0)$. Assume further that for every $t\in[0,T]$, $ 0<r_- \leq \mathrsfs{R}(t),\mathscr{R}(t)\leq r_+<\infty$ and $0<n_- \leq \widehat{\mathrsfs{B}}_{44}(t),\Phi(t)\leq n_+<\infty$. Then there exists a constant $C$, independent of $d$, such that
    \[
        \sup_{0\leq t\leq T}\left(\left|\mathrsfs{R}(t)-\mathscr{R}(t)\right|
        +
        \left|
        \widehat{\mathrsfs{B}}_{44}(t)-\Phi(t)
        \right|\right)
        \leq
        C\left(\|K\|_{\opt}, \bar{\gamma}, T, r_-, r_+, n_-, n_+\right)
                \delta^2
                \left(
                \frac{1}{2}-\frac{1}{\pi}
                \right).            
    \]
    \end{proposition}
    We defer the proof to Appendix~\ref{sect:compareADVHSGD:to:HSGD}. Note that it is possible to relax the boundedness assumption on the deterministic equivalents of the least squares risks and squared $\ell_2$-norms by introducing an analogous theory of the set $\mathcal{U}$ (see Assumption~\ref{assumption:risk}) for HSGD. 
    
    Under the assumptions presented in Propositions \ref{prop:q(t):bounded}, \ref{prop:isotropic_Rasymp} and \ref{prop:powerlaw_Rasymp} and Corollaries \ref{cor:powerlaw_R_B_44_bounded:isotropic} and \ref{cor:powerlaw_R_B_44_bounded:plaw}, for isotropic covariance and power law covariance when $d\to\infty$ (See Section~\ref{sect:plaw}), the deterministic equivalents of the risk and squared $\ell_2$-norm remain bounded from above and below by constants independent of dimension. Using an almost identical argument, we could show this holds in the case of $\mathscr{R}(t)$ and $\Phi(t)$. It follows that, under the assumptions presented in these Corollaries and up to constants independent of $d$, Proposition~\ref{prop:adv:reg:volterra:stability} holds. From Theorem~\ref{thm:concentration_statistic:adv} and similar arguments to Theorem 4 in \citep{collinswoodfin2023hitting} and Theorem 2.1 in \citep{collinswoodfin2024highline}, we can translate this back to SGD. 

\section{Idealized Exact Line Search and Polyak Stepsize}
In this section, we consider the idealized algorithms: \emph{exact line search} and \emph{Polyak stepsize}. See \citep{collinswoodfin2024highline} for a discussion of stochastic line search and Polyak-type stepsizes. Let $\mathrsfs{D}^{2,\Adv}(t)$ be the deterministic equivalent of $\|X_{\lfloor td \rfloor } - X^{\star, \Adv}\|^2$ where we recall $X^{\star,\Adv} = \operatorname{argmin}_X \mathcal{R}^\Adv(X)$. In deterministic optimization, the Polyak stepsize is the learning rate schedule chosen for maximum decrease of the distance to optimality at the following iteration. Exact line search is the learning rate strategy chosen for maximum decrease of the objective function at the next iteration. 

Analogously to \citep{collinswoodfin2024highline}, given that $\mathcal{R}^\Adv(X_{\lfloor td \rfloor})$ and $\|X_{\lfloor td \rfloor } - X^{\star, \Adv}\|^2$ are close to $\mathrsfs{R}^\Adv(t)$ and $\mathrsfs{D}^{2,\Adv}(t)$ when $d\to \infty$ (See Theorem~\ref{thm:main:conc:St:general:with:stop:adv}), we can ask this to hold for the deterministic equivalents. We restrict the results in this section to $\ell_2$-adversarial least squares (see start of Section~\ref{sect:AdvHSGD:vs:HSGD}) as the general expressions derived are too cumbersome. See Appendix~\ref{sect:Polyakstepsize} and \ref{sect:linesearch} for general results and the proofs. Given that we do not have access to $X^\star$ nor $X^{\star, \Adv}$ nor distributional knowledge of $a$ in practice, these are idealized algorithms. Nevertheless, they serve as a foundation for developing more practical algorithms. 

\paragraph{Polyak Stepsize.} The stability threshold for descent towards $X^{\star, \Adv}$ corresponds to the largest learning rate such that $\dif \mathrsfs{D}^{2,\Adv}(t) < 0$, which we call $\gamma_{k}^{\operatorname{Stable, adv}}$. For $\ell_2$-adversarial least squares, translating back deterministic equivalents to SGD, this corresponds to 
\begin{equation}\label{eq:polyak:sgd:leastsquares:main}
        \gamma_{k}^{\operatorname{Stable, adv}}
         = 
      \frac{2\mathcal{R}(X_k) - \eta^2 - (X_k- X^\star)^\top K(X^{\star, \Adv}- X^\star)+ \tilde{\lambda}_k^{\operatorname{eff}}(\|X_k\|^2  -\langle X_k , X^{\star, \Adv}\rangle) }{\frac{\Tr(K) }{2d}\left(2\mathcal{R}(X_k) + \tilde{\lambda}_k^{\operatorname{eff}} \|X_k\|^2\right)}.
\end{equation}
By Lemma~\ref{lem:strong:con:Radv}, for any $\delta>0$, $ \mathcal{R}^{\operatorname{Adv}}$ is $\upsilon$-strongly convex with unique minimizer $X^{\star, \Adv}$ which implies that for any $X$ and combination of $\delta$ and $K$, there exists a learning rate sequence for which descent towards $X^{\star,\Adv}$ is possible. However, at the optimum $X = X^{\star,\Adv}$, we have $\gamma_{k}^{\operatorname{Polyak, adv}}=0$ from which it follows that no constant learning rate schedule guarantees global monotone descent towards $X^{\star,\Adv}$. This differs from standard noiseless least squares, where the stability threshold for descent towards $X^\star$ is the constant $\gamma_{k}^{\operatorname{Stable}}=\frac{2d}{\Tr(K)}$ \citep{loizou2021stochastic, paquetteSGD2021}.

The numerator of stability threshold above admits a natural decomposition into three contributions:
\begin{enumerate}
\item $2\mathcal{R}(X_k) - \eta^2$ is the standard descent contribution towards $X^\star$
\item $(X_k- X^\star)^\top K(X^{\star, \Adv}- X^\star)$ is the correction for the shift between $X^\star$ and $X^{\star, \Adv}$
\item $\tilde{\lambda}_k^{\operatorname{eff}}(\|X_k\|^2  -\langle X_k , X^{\star, \Adv}\rangle) $ arises from the adversarial perturbation of the input
\end{enumerate} 
This also holds for our general results and hints that a practical learning rate schedule for $\ell_2$-adversarial training must take into consideration the difference in objective between the standard loss and adversarial loss and must take into consideration the adversarial perturbation of the input. 

 The Polyak stepsize corresponds to $ \gamma^{\operatorname{Polyak, adv}}(t) \in  \operatorname{argmin}_{\gamma} \dif \mathrsfs{D}^{2,\Adv}(t) $. In the case of $\ell_2$-adversarial least squares, this corresponds to $  \gamma_{k}^{\operatorname{Polyak, adv}}=  \tfrac{1}{2} \gamma_{k}^{\operatorname{Stable, adv}}$. 

 \paragraph{Exact Line Search.} Inspired by \citep{collinswoodfin2024highline}, we denote the greedy learning rate $\gamma^{\operatorname{line}}(t) \in \operatorname{argmin}_{\gamma}  \dif \mathrsfs{R}^\Adv(t)$ which produces the largest decreases in adversarial risk at the next iteartion as $\emph{exact line search}$. For $\ell_2$-adversarial least squares, translating back deterministic equivalents to SGD, this corresponds to 
\begin{equation}\label{eq:linesearch:ls:sgd}
    \begin{aligned}
    &\gamma^{\operatorname{line}}_k
    = \frac{ \|\nabla \mathcal{R}^{\Adv}(X_k)\|^2}
    {
        \frac{2}{d}\left(1 + \delta \sqrt{\frac{2}{\pi}}\frac{\|X_k\|}{\sqrt{2 \mathcal{R}(X_k)}}\right)
        \mathcal{R}^\Adv(X_k)\left(
    \operatorname{Tr}(K^2)
    +
    \tilde{\lambda}^{\mathrm{eff}}_k\operatorname{Tr}(K)
    \right)}.
    \end{aligned}
\end{equation}
 Here $ \tilde{\lambda}^{\operatorname{eff}}_k$ is the effective learning rate of SGD as defined in Section~\eqref{sect:AdvHSGD:vs:HSGD}. It turns out $\gamma^{\operatorname{line}}(t)$, the deterministic equivalent of \eqref{eq:linesearch:ls:sgd}, can be expressed in terms of a ridge regression objective with fixed parameter $\tilde{\lambda}^{\operatorname{eff}}(t)$. See Appendix~\ref{sect:linesearch} for the proof of Proposition~\ref{prop:linesearch:explicit}.
\begin{proposition}\label{prop:linesearch:explicit}
    Consider $\ell_2$-adversarial least squares and suppose $K$ has $0<n \leq d$ distinct eigenvalues $\lambda_n < \cdots < \lambda_1$ each with algebraic multiplicity $\#\lambda_j$. Let $J_+ \defas \{j : \lambda_j >0\}$. Then the learning rate schedule of exact line search for all $t\geq 0$ where $\mathrsfs{E}(t) >0$ is given by
    \begin{equation}\label{eq:linesearch:ridge:general}
        \gamma^{\operatorname{line}}(t) = \frac{\sum_{j=1}^n \alpha_j(t) (\lambda_j + \tilde{\lambda}^{\operatorname{eff}}(t))}{\sum_{j=1}^n \frac{\#\lambda_j}{d}\lambda_j (\lambda_j + \tilde{\lambda}^{\operatorname{eff}}(t)) } \left( 1- \frac{
            \Phi_{\min}(\tilde{\lambda}^{\operatorname{eff}}(t))
            }{
             \Phi_{\tilde{\lambda}^{\operatorname{eff}}(t)}(t)
            }\right).
    \end{equation}
    Here $\Phi_{\tilde{\lambda}^{\operatorname{eff}}(t)}(t)
    \defas
    \mathrsfs R(t)
    +
    \frac{\tilde{\lambda}^{\operatorname{eff}}(t)}{2}
    \widehat{\mathrsfs B}_{44}(t)$ is the ridge regression objective with fixed parameter $\tilde{\lambda}^{\operatorname{eff}}(t)$ and $\Phi_{\operatorname{min}}(\tilde{\lambda}^{\operatorname{eff}}(t))
    =
    \frac{1}{2}
    \sum_{j \in J_+} \lambda_j
    \frac{\tilde{\lambda}^{\operatorname{eff}}(t)}
    {\lambda_j+\tilde{\lambda}^{\operatorname{eff}}(t)}
    \mathrsfs V_{j,22}(t) + \frac{\eta^2}{2}$ is its minimum. Here
    where 
    \[
\mathcal{E}_j(t)
\defas
\frac{1}{2}\left(
\left(\lambda_j+\tilde{\lambda}^{\operatorname{eff}}(t)\right)\mathrsfs V_{j,11}(t)
-
2\lambda_j\mathrsfs V_{j,12}(t)
+
\frac{\lambda_j^2}{\lambda_j+\tilde{\lambda}^{\operatorname{eff}}(t)}
\mathrsfs V_{j,22}(t)\right),
\]
denotes the $\lambda_j$ eigenspace's contribution to the excess risk $\mathcal{E}(t)\defas\Phi_{\tilde{\lambda}^{\operatorname{eff}}(t)}(t) - \Phi_{\operatorname{min}}(\tilde{\lambda}^{\operatorname{eff}}(t))$ such that $\mathcal{E}(t) = \sum_{j=1}^n \mathcal{E}_j(t)$ and, $0\leq \alpha_{j} \defas \frac{ \mathrsfs{E}_j}{\mathrsfs{E}(t)}\leq 1$ for $1\leq j\leq n$. In the isotropic setting, the learning rate schedule further simplifies to
    \[
        \gamma^{\operatorname{line}}(t)
        =
        1-
        \frac{
        \Phi_{\min}(\tilde{\lambda}^{\operatorname{eff}}(t))
        }{
         \Phi_{\tilde{\lambda}^{\operatorname{eff}}(t)}(t)
        }.
    \]
    Note that when $\mathcal{E}(t) =0$, we have reached the optimum.
\end{proposition}

\begin{figure}[t!]
    \centering
    \includegraphics[width=0.47\textwidth]{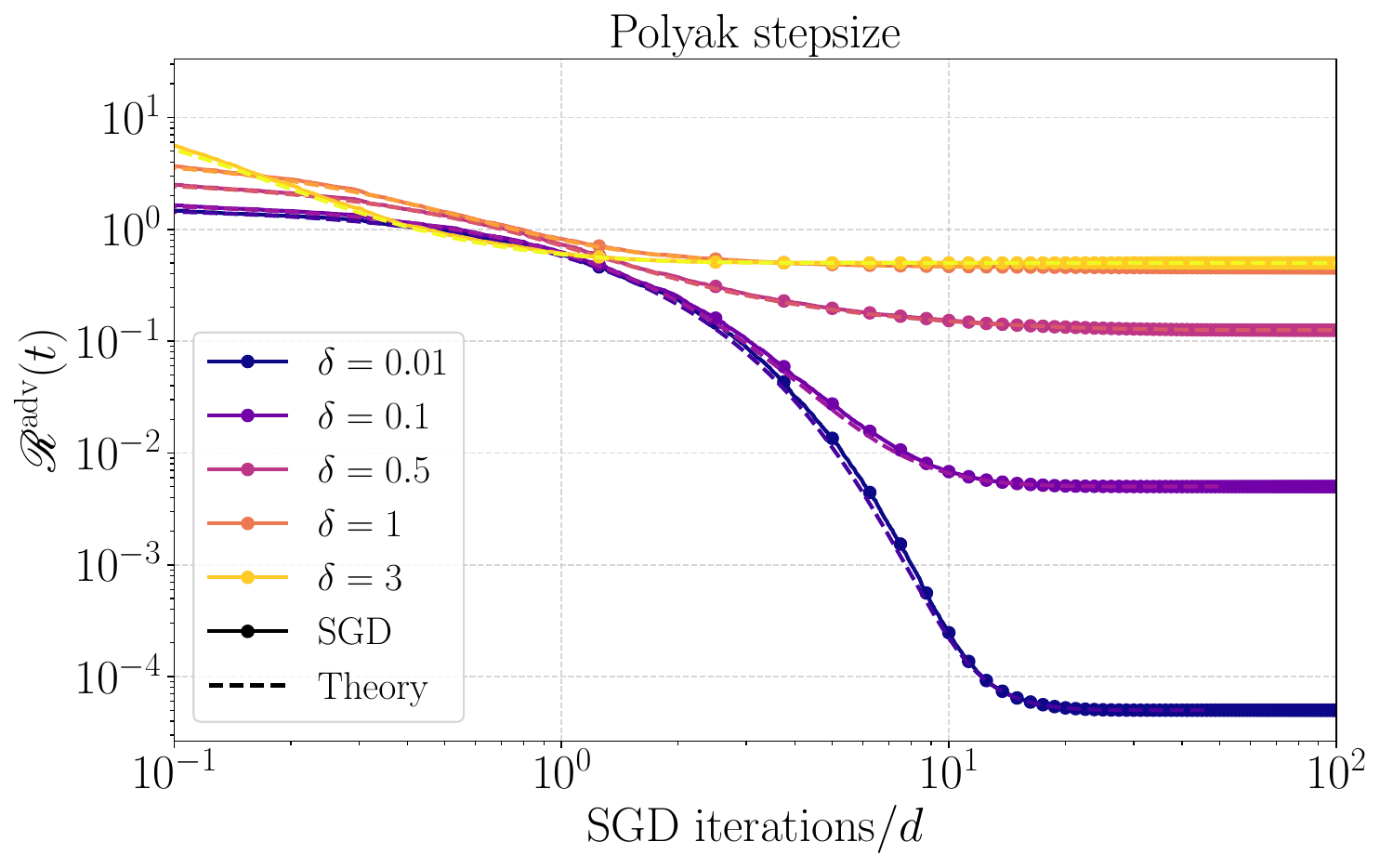}
    \includegraphics[width=0.47\textwidth]{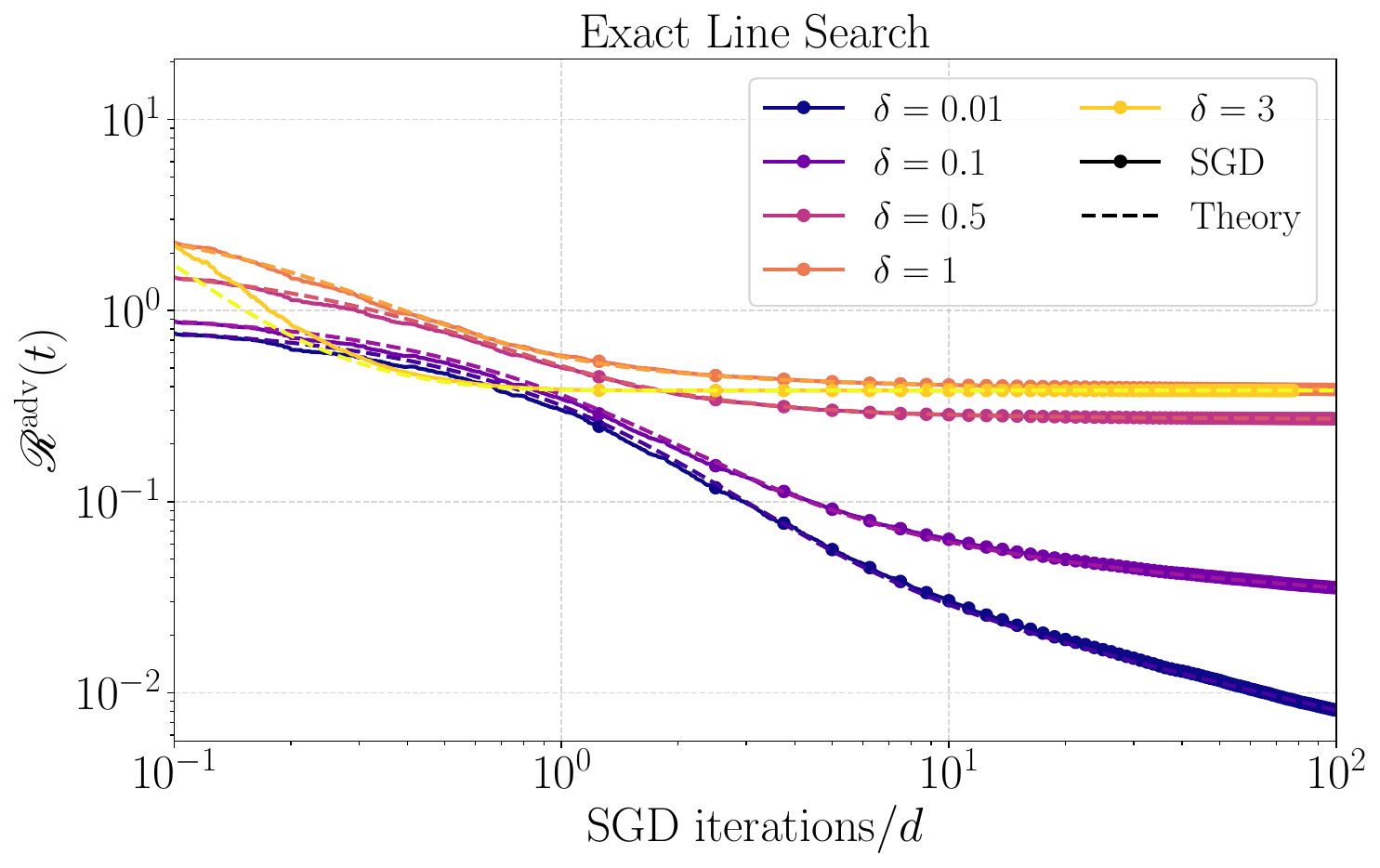}
    \caption{\textbf{SGD with exact line search $\gamma_{k}^{\operatorname{line}}$ or Polyak stepsize $\gamma_{k}^{\operatorname{Polyak, adv}}$ } matches closely the path of our system of ODEs \eqref{eq:ODE:V_i(t):def} with deterministic learning rates schedules $\gamma(t)=\gamma^{\operatorname{line}}(t)$ and $\gamma(t)=\gamma^{\operatorname{Polyak, adv}}(t)$ for $\mathcal{R}^\Adv(X_k)$ on noiseless $\ell_2$-adversarial least squares. See Appendix~\ref{app:captions} for simulation details.}
    \label{fig:polyak:linesearch:evidence}
\end{figure}
 It turns out the deterministic equivalent of the Polyak stepsize \eqref{eq:linesearch:ls:sgd} can also be expressed in terms of a ridge regression objective similarly to Proposition~\ref{prop:linesearch:explicit}.
\begin{proposition}\label{lem:Polyak:ridge:general}
    Consider $\ell_2$-adversarial least squares and assume $K$ has an arbitrary number $0<n \leq d$ distinct eigenvalues. For $1\leq j \leq n$, define
    \[
        s_j^\star
        \defas
        \begin{cases}
        \displaystyle
        \frac{\lambda_j}{\lambda_j+\tilde{\lambda}_\star},
        &  0<\lambda_j+\tilde{\lambda}_\star<\infty\\
        0
        & \lambda_j=\tilde{\lambda}_\star=0 \text{ or } \tilde{\lambda}_\star=\infty,
        \end{cases}
        \quad
    \vartheta^\star_j(t) \defas \begin{cases}
\left(
\tilde{\lambda}_\star
-
\tilde{\lambda}^{\mathrm{eff}}(t)
\right)
s_j^\star & \text{if}\quad \tilde{\lambda}_\star<\infty \\
\lambda_j & \text{if} \quad \tilde{\lambda}_\star=\infty
        \end{cases}.
    \]
    Then the deterministic equivalent of the Polyak stepsize \eqref{eq:polyak:sgd:leastsquares:main} for all $t\geq 0$ where $\widehat{\mathrsfs{B}}_{44}(t) > 0$ is given by 
    \begin{equation}\label{eq:Polyak:ridge:general}
        \begin{gathered}
        \gamma^{\operatorname{Polyak},\Adv}(t) = \frac{
            \Phi_{\tilde{\lambda}^{\operatorname{Polyak}}_\star(t)}(t) - \Phi_{\min}(\tilde{\lambda}_{\star}) 
        }{
        \frac{\operatorname{Tr}(K)}{d}
        \left(
        \Phi_{\tilde{\lambda}^{\mathrm{eff}}(t)}(t)
        \right)},\quad
        \text{for}\quad  
            \tilde{\lambda}^{\operatorname{Polyak}}(t) \defas \widetilde{\lambda}^{\mathrm{eff}}(t) + \sum_{j=1}^d\vartheta^\star_j(t) \frac{ \mathrsfs V_{j,12}(t) } { \widehat{\mathrsfs B}_{44}(t)}.
        \end{gathered}
    \end{equation}
    Here $\Phi_{\tilde{\lambda}^{\operatorname{Polyak}}_\star(t)}(t) \defas
    \mathrsfs R(t)
    +
    \frac{\tilde{\lambda}^{\operatorname{Polyak}}_\star(t)}{2}\widehat{\mathrsfs B}_{44}(t)$ is the ridge regression objective with fixed parameter $\tilde{\lambda}^{\operatorname{Polyak}}_\star(t)$ and $\Phi_{\operatorname{min}}(\tilde{\lambda}^{\operatorname{Polyak}}_\star(t))$ is its minimum defined in Proposition~\ref{prop:linesearch:explicit}. When $\tilde{\lambda}_\star = \infty$, we note that 
    \begin{equation}\label{eq:convention:limit:lambdastar}
    \Phi_{\min}(\tilde{\lambda}_\star)
    \defas
    \lim_{\lambda\to\infty}\Phi_{\min}(\lambda)
    =
    \frac{\eta^2}{2}
    +\frac{1}{2}
    (X^\star)^\top KX^\star.
    \end{equation}
\end{proposition}
See Appendix~\ref{sect:linesearch} for the proof of Proposition~\ref{lem:Polyak:ridge:general}. Note that Theorems~\ref{thm:main:conc:St:general:with:stop:adv}, ~\ref{thm:concentration_statistic} and \ref{thm:concentration_statistic:adv} are formulated in terms of the set $\mathcal{U}$ in which $\widehat{\mathrsfs{B}}_{44}(t)>0$. Hence, one can view Lemma~\ref{lem:Polyak:ridge:general} as remaining valid for all time these theorems hold.
\begin{remark}
    If we set $\gamma(t)=\gamma^{\operatorname{line}}(t)$ or $\gamma(t)=\gamma^{\operatorname{Polyak, adv}}(t)$ in our system of ODEs \eqref{eq:ODE:V_i(t):def} specialized to $\ell_2$-adversarial least squares (see \eqref{eq:odes:linrreg} in Appendix~\ref{sec:lsq_analysis}), then Theorem~\ref{thm:main:conc:St:general:with:stop:adv} still seems to hold in the case running SGD with the respective learning rate schedules $\gamma_{k}^{\operatorname{Polyak, adv}}$ or $\gamma_{k}^{\operatorname{line}}$. See Figure~\ref{fig:polyak:linesearch:evidence} for numerical evidence. This is to be expected as it would be possible to extend the framework of \citep{collinswoodfin2024highline} to $\ell_2$-adversarial training.
\end{remark}

\paragraph{Comparison of Exact Line Search and Polyak stepsize.}
We now compare the deterministic equivalents \eqref{eq:polyak:sgd:leastsquares:main} and \eqref{eq:linesearch:ls:sgd}. The representations \eqref{eq:linesearch:ridge:general} and \eqref{eq:Polyak:ridge:general} allow us to compare exact line search and the Polyak stepsize locally, at a fixed state of the limiting ODE. We use this comparison to identify the factors that can cause exact line search to be suboptimal relative to the Polyak stepsize. Given a fixed state of the coupled ODEs \eqref{eq:odes:linrreg}, the ratio of $\gamma^{\operatorname{line}}(t)/\gamma^{\operatorname{Polyak},\Adv}(t)$ decomposes into a factor driven by anisotropy and a factor driven by the mismatch in $\tilde{\lambda}^{\mathrm{eff}}(t)$ and $\tilde{\lambda}^{\operatorname{Polyak}}_\star(t)$:
\begin{equation}\label{eq:ratio:linesearch:Polyak}
\boxed{
\frac{
    \gamma^{\operatorname{line}}(t)
}{
\gamma^{\operatorname{Polyak},\Adv}(t)
}
=
\underbrace{
\frac{
\frac{\operatorname{Tr}(K)}{d}
\sum_{r=1}^n
\alpha_r(t)
\left(
\lambda_r+\tilde{\lambda}^{\mathrm{eff}}(t)
\right)
}{
\sum_{r=1}^n
\frac{\#\lambda_r}{d}
\lambda_r
\left(
\lambda_r+\tilde{\lambda}^{\mathrm{eff}}(t)
\right)
}}_{Spectral factor}
\cdot
\underbrace{
\frac{
\Phi_{\tilde{\lambda}^{\mathrm{eff}}(t)}(t)
-
\Phi_{\min}\!\left(\tilde{\lambda}^{\mathrm{eff}}(t)\right)
}{
\Phi_{\tilde{\lambda}_{\star}^{\operatorname{Polyak}}(t)}(t)
-
\Phi_{\min}(\tilde{\lambda}_{\star})
}}_{Ridge factor}.
}
\end{equation}
From a simple computation, the first term in \eqref{eq:ratio:linesearch:Polyak} is equal to one if and only if
\begin{equation}\label{eq:condition:equality:polyak}
    \sum_{r=1}^n
    \alpha_r(t)
    \lambda_r = \frac{\Tr(K^2)}{\Tr(K)}.
\end{equation}
One natural way for equality in \eqref{eq:condition:equality:polyak} to hold is that the excess risk weights $\alpha_j(t)$ match each eigenspace's contribution to the trace: 
\[
\alpha_j(t) = \frac{\# \lambda_j \cdot \lambda_r}{\Tr(K)}.
\] 

\begin{figure}[t!]
    \centering
    \includegraphics[width=0.30\textwidth]{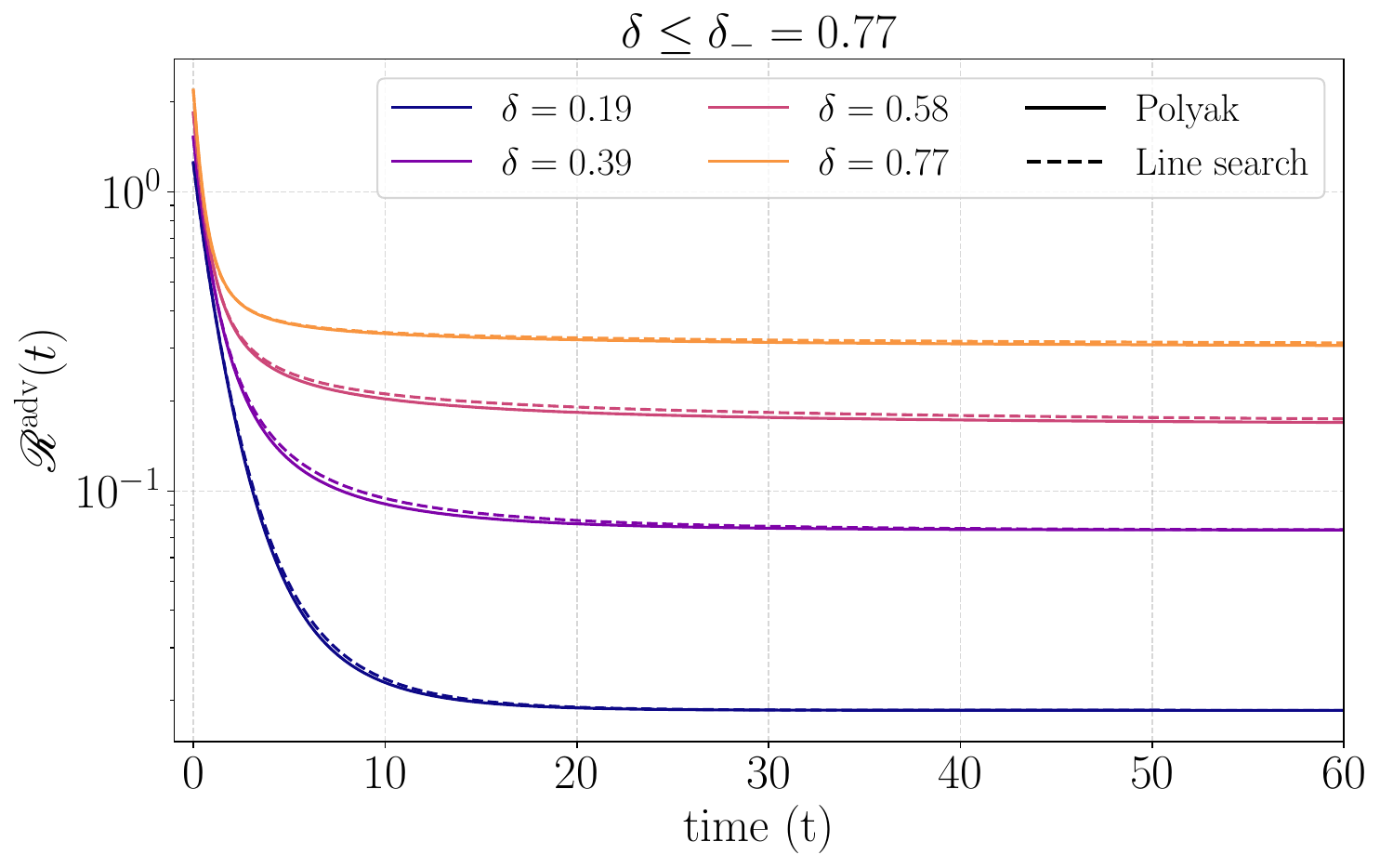}
    \includegraphics[width=0.30\textwidth]{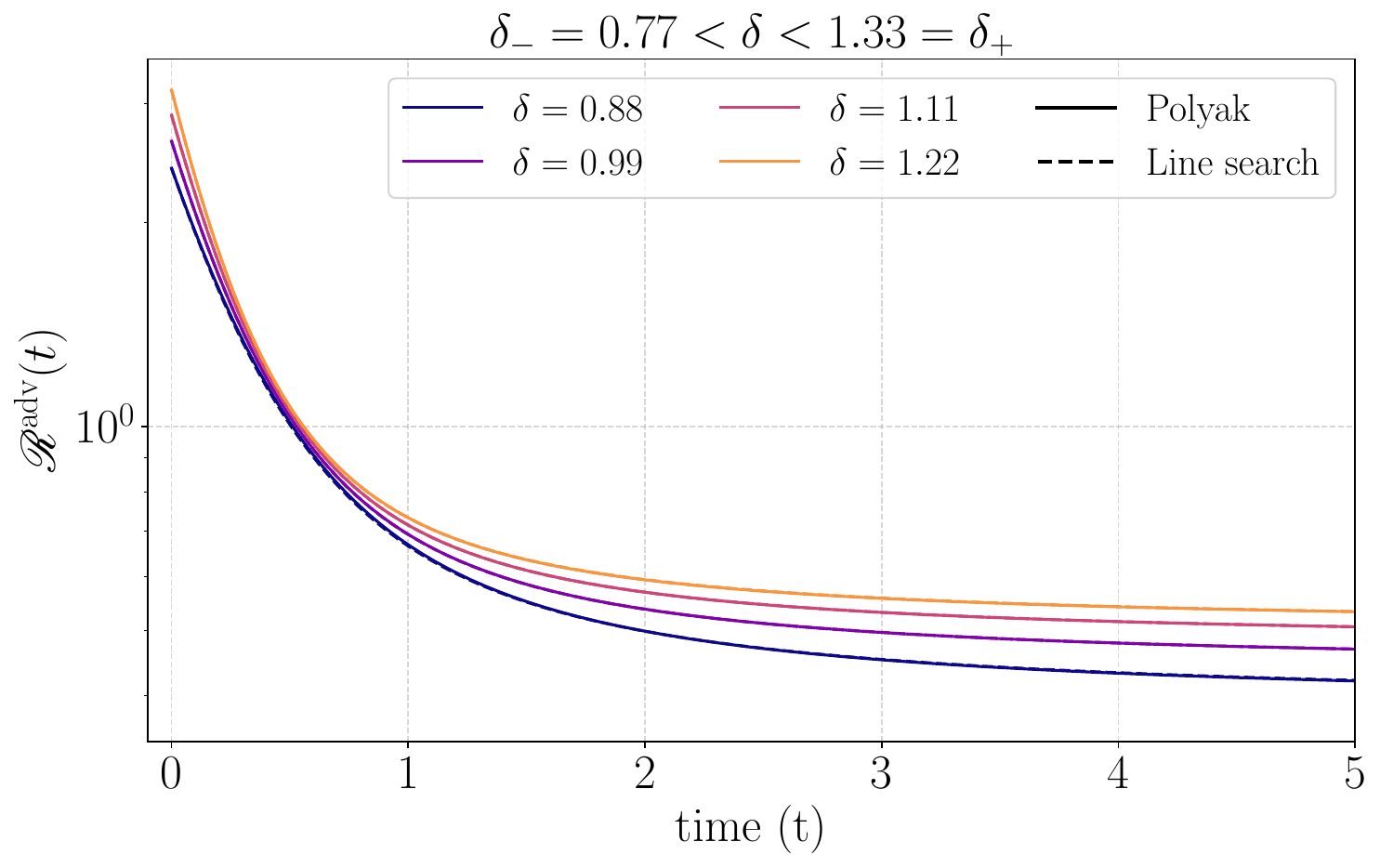}
    \includegraphics[width=0.30\textwidth]{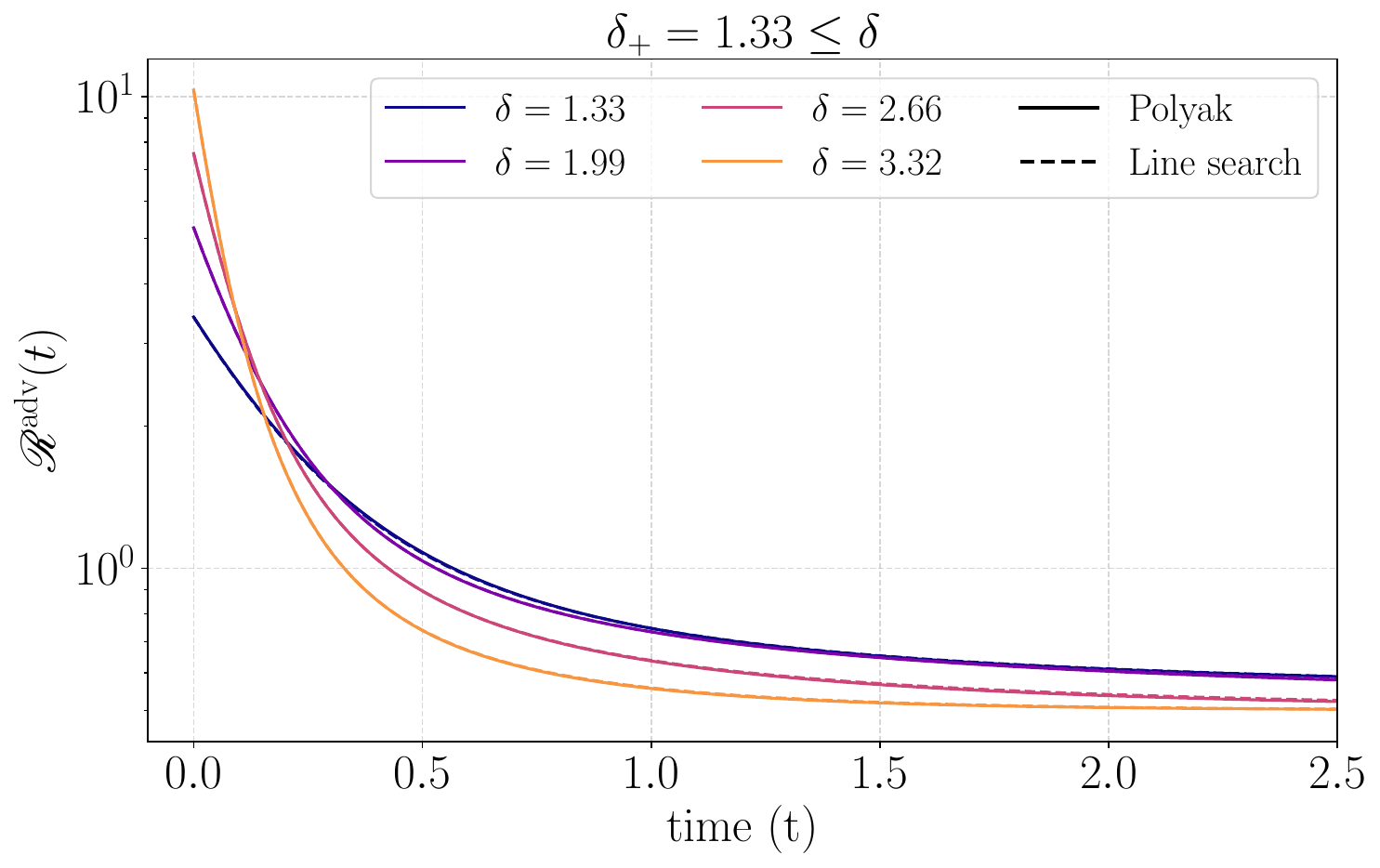}
    \caption{\textbf{Comparison between Exact Line Search and Polyak Stepsize under weak anisotropy} on noiseless $\ell_2$-adversarial least squares for different values of $\delta$ in the three regimes of $X^{\star,\Adv}$ (See Proposition~\ref{prop:fixed:point:gamma:to:0:main}). The three plots illustrate the convergence of the $\ell_2$-adversarial risk and that, under weak anisotropy, exact line search and the Polyak stepsize perform similarly. See Figure~\ref{fig:polyak:delta} for results in the case of strong anisotropy. See~\eqref{eq:ratio:linesearch:Polyak} for more details and Section~\ref{app:captions} for simulation details.}
    \label{fig:polyak:anisotropy}
\end{figure}

 If $K$ is nearly isotropic then this condition holds approximately so the spectral factor is nearly one and the mismatch is mainly driven by the ridge factor. Conversely, in the presence of strong anisotropy, the spectral factor can be far from one. When the remaining ridge excess $\mathcal{E}(t)=\Phi_{\tilde{\lambda}^{\operatorname{eff}}(t)}(t) - \Phi_{\operatorname{min}}(\tilde{\lambda}^{\operatorname{eff}}(t))$ is concentrated in low-curvature eigenspaces, we have $ \sum_{r=1}^n\alpha_r(t)\lambda_r  <  \frac{\Tr(K^2)}{\Tr(K)}$ and, for a fixed ODE state, exact line search chooses a smaller stepsize than Polyak. When the excess is sufficiently concentrated in high-curvature eigenspaces, the opposite can occur. Thus anisotropy can cause greedy line search to deviate substantially from the Polyak stepsize. Under the normalization $\frac{1}{d}\Tr(K)=1$, the spectral factor in
 \eqref{eq:ratio:linesearch:Polyak} takes the form
 \[
 \frac{
 \sum_{r=1}^n\alpha_r(t)\lambda_r+\tilde{\lambda}^{\mathrm{eff}}(t)
 }{
 \frac{1}{d}\Tr(K^2)+\tilde{\lambda}^{\mathrm{eff}}(t).
 }
 \]
 Thus the discrepancy in the spectral factor is governed by how much
 $\sum_r\alpha_r(t)\lambda_r$ differs from 
 $\frac{1}{d}\Tr(K^2)$ and by $\tilde{\lambda}^{\mathrm{eff}}(t)$. Strong anisotropy increases
 $\frac{1}{d}\Tr(K^2)$, while concentration of the remaining ridge excess in
 low-curvature eigenspaces decreases $\sum_r\alpha_r(t)\lambda_r$. Both effects
 make the line search stepsize smaller relative to the Polyak stepsize. In contrast, a large effective regularization parameter dampens this discrepancy. We emphasize this for $K$ with the eigenvalues generated by \citep{collinswoodfin2024highline}  
 \begin{equation}\label{eq:s_values:eigs}
 \lambda_i
 =
 \frac{
 d
 }{
 \sum_{j=1}^d
 \left(\frac{j}{d+1}\right)^{-2/s}
 }\left(\frac{i}{d+1}\right)^{-2/s},
 \qquad
 i=1,\ldots,d,
 \qquad
 s>2.
 \end{equation}
 By construction, $\frac{1}{d}\Tr(K)=1$. As $s\downarrow2$, the spectrum becomes
 more spread out and $\frac{1}{d}\Tr(K^2)$ increases.
 Consequently, when $\tilde{\lambda}^{\mathrm{eff}}(t)$ is small and the
 remaining ridge excess is not concentrated in high-curvature eigenspaces, the
 spectral factor can be far below one.
 
 For a fixed ODE state, the three adversarial regimes (see Proposition~\ref{prop:fixed:point:gamma:to:0:main} and \citep[Proposition 1]{xing2021same_minimizer}) determine how strongly this anisotropy can
 affect the same-state comparison. If $\delta\le\delta_-$, then
 $\tilde{\lambda}^{\mathrm{eff}}(t)\to\tilde{\lambda}_\star=0$, so the damping
 from the ridge term disappears near the optimum. This is the regime most
 analogous to standard least squares as strong anisotropy can make exact line
 search substantially more conservative than Polyak \citep{collinswoodfin2024highline}. If
 $\delta_-<\delta<\delta_+$, then
 $\tilde{\lambda}^{\mathrm{eff}}(t)\to\tilde{\lambda}_\star\in(0,\infty)$, so
 the ridge term partially damps the spectral discrepancy. 
 
\begin{figure}[t!]
    \centering
    \includegraphics[width=0.30\textwidth]{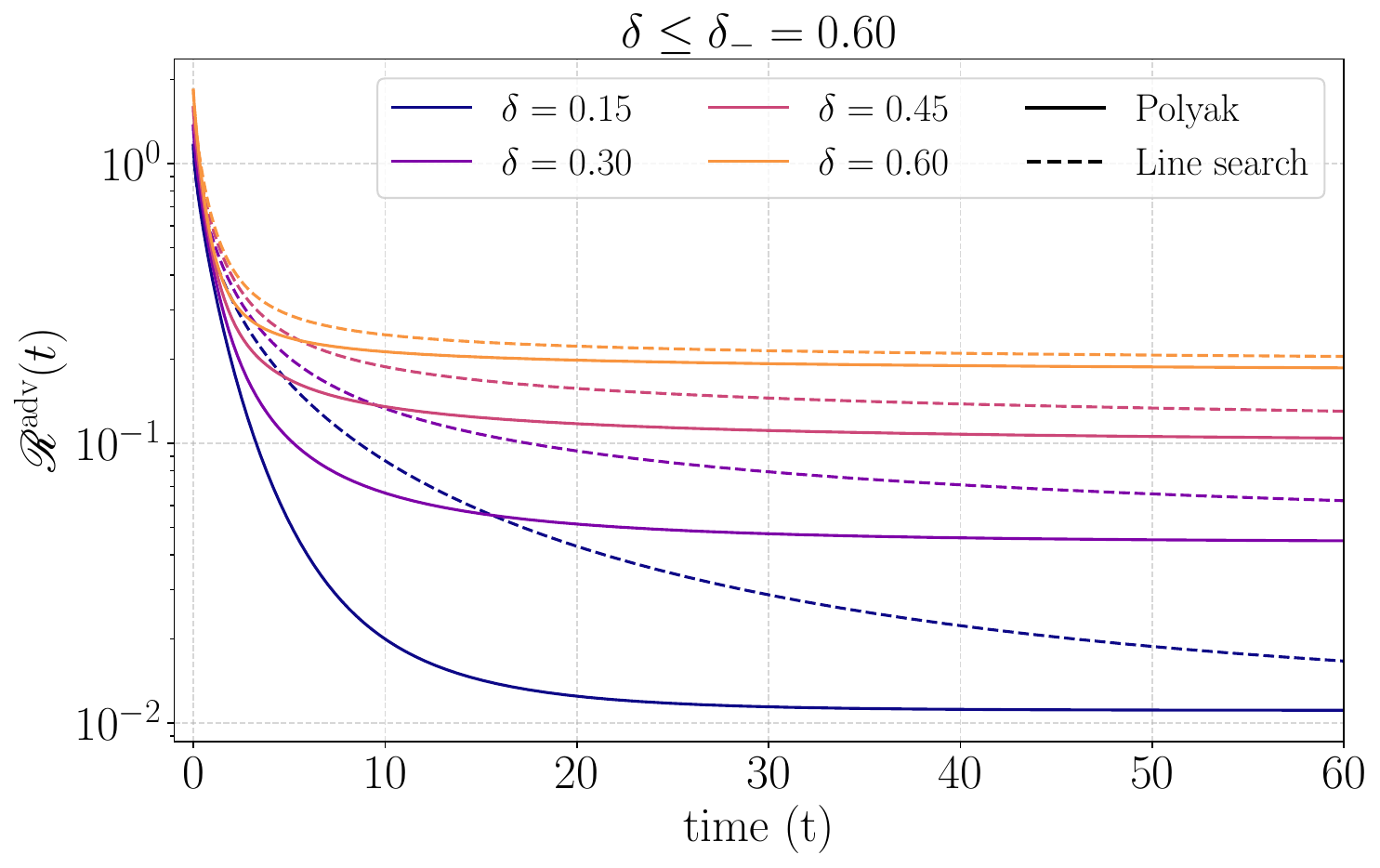}
    \includegraphics[width=0.30\textwidth]{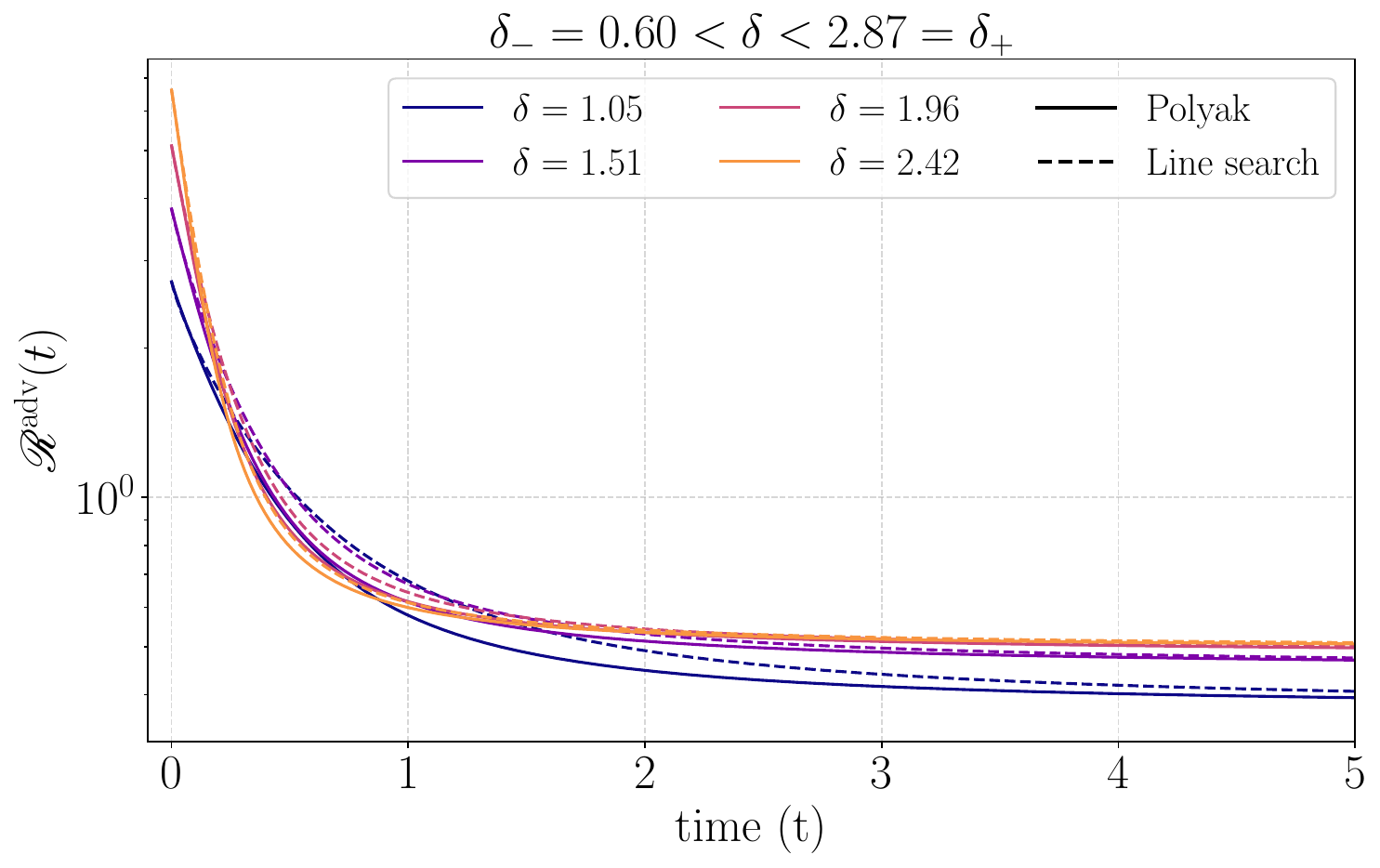}
    \includegraphics[width=0.30\textwidth]{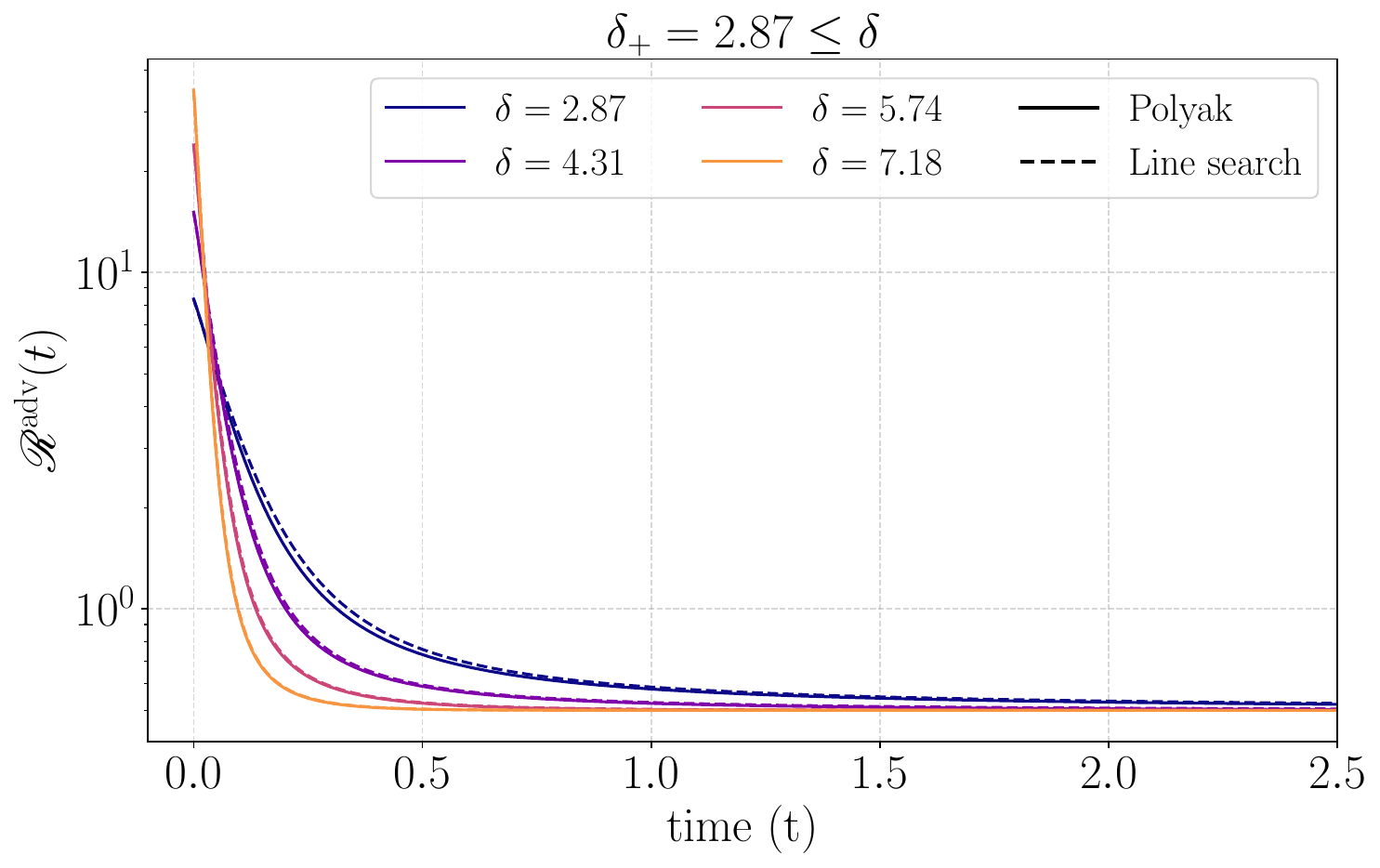}
    \caption{\textbf{Comparison for Exact Line Search and Polyak Stepsize under strong anisotropy} on noiseless $\ell_2$-adversarial least squares for different values of $\delta$ in the three regimes of $X^{\star,\Adv}$ (See Proposition~\ref{prop:fixed:point:gamma:to:0:main}). The three plots illustrate the convergence of the $\ell_2$-adversarial risk and how $\delta$ and $\tilde{\lambda}^{\operatorname{eff}}(t)$ mitigate the influence of strong anisotropy on the discrepancy between the Polyak stepsize and exact line search. See Figure~\ref{fig:polyak:anisotropy} for results in the case of weak anisotropy. See~\eqref{eq:ratio:linesearch:Polyak} for more details and Section~\ref{app:captions} for simulation details.}
    \label{fig:polyak:delta}
\end{figure}

 If $\delta\ge\delta_+$, then $\tilde{\lambda}^{\mathrm{eff}}(t) \to \infty$ as $t\to \infty$ so the spectral factor approaches one. In this
 large $\delta$ regime, anisotropy no longer creates a large mismatch
 between line search and Polyak and any remaining discrepancy is driven by the
 ridge factor. See Figure~\ref{fig:polyak:anisotropy} and Figure~\ref{fig:polyak:delta} for numerical evidence illustrating the impact of $\delta$ and anisotropy on the convergence of the $\ell_2$-adversarial risk.

 \section{$\ell_2$-Adversarial Linear Regression Analysis}\label{sect:lsadv:analysis:main}
 In this section, we will analyze the training dynamics of $\ell_2$-adversarial least squares with regularization parameter $\lambda=0$, with $\gamma_k =\gamma(t) \equiv \gamma$ and under the assumption that the data $a\sim \mathcal{N}(0, K)$. Given that this is equivalent to $p_1=1$ and $p_2=0$ with $\mu_1=0$, we drop the subscript $i$ in this section. In this setting, the adversarial risk simplifies to
 \begin{equation}
    \begin{aligned}
    &\mathcal{R}^{{\operatorname{adv}}}(X) 
    =
    \mathcal{R}(X)+ \delta\sqrt{\frac{2}{\pi}}\|X\| \sqrt{2\mathcal{R}(X)} +\frac{\delta^2}{2}  \|X\|^2
    \\
    \end{aligned}
    \end{equation}
    From \eqref{eq:Def:Radv:lr} and $\mathrsfs{R}^{{\operatorname{adv}}}(t) = h \circ (\widehat{\mathrsfs{B}}(t))$, the deterministic equivalent of $\mathcal{R}^{{\operatorname{adv}}}$ is given by
    \begin{equation}\label{eq:Radv(t):detequiv:main}
        \mathrsfs{R}^{{\operatorname{adv}}}(t) =\mathrsfs{R}(t) + \delta\sqrt{\frac{2}{\pi}}\sqrt{2\widehat{\mathrsfs{B}}_{44}(t)\mathrsfs{R}(t)} +\frac{\delta^2}{2} \widehat{ \mathrsfs{B}}_{44}(t).
    \end{equation}
    It is clear that we lose differentiability at $\mathrsfs{R}(t)= \frac{1}{2} \left(\mathrsfs{B}_{11}(t) - \mathrsfs{B}_{12} - \mathrsfs{B}_{21} + \mathrsfs{B}_{22}\right)+ \frac{\eta^2}{2} = 0$ and $\widehat{\mathrsfs{B}}_{44}(t) = 0$. To apply the theory in this paper, we work on the set $\mathcal{U}$ defined as follows
    \begin{equation}\label{eq:example:U}
        \mathcal{U} \defas \left\{ \widehat{B}: \widehat{B}_{11} - \widehat{B}_{12} - \widehat{B}_{21} + \widehat{B}_{22} + \eta > 0, \quad \widehat{B}_{44} >0\right\}
    \end{equation}
    As in Section~\ref{sect:AdvHSGD:vs:HSGD}, we introduced the effective learning rate $\tilde{\gamma}^{\operatorname{eff}}_k$ and regularization parameter $\tilde{\lambda}^{\operatorname{eff}}_k$ whose dynamics are determined from $q(X_k)\defas \left(\|X_k\|^2/ (2\mathcal{R}(X_k))\right)^{1/2}$. Their respective deterministic equivalents are given by 
    \begin{equation}\label{eq:adapt:gamma:lambda:main}
     q(t) \defas \sqrt{\frac{\widehat{\mathrsfs{B}}_{44}(t)}{2\mathrsfs{R}(t)}},\quad  \tilde{\gamma}^{\operatorname{eff}}(t) = \gamma \left(1 + \delta \sqrt{\frac{2}{\pi}} q(t)\right)
        \quad \text{and}\quad 
        \tilde{\lambda}^{\operatorname{eff}}(t) = \frac{ \delta^2 + \delta \sqrt{\frac{2}{\pi}} \frac{1}{q(t)}}{1 + \delta \sqrt{\frac{2}{\pi}}q(t)}.
    \end{equation}
    We also introduce the following integral terms
    \begin{equation}\label{eq:adapt:gamma:lambda:int:main}
        \tilde{\Gamma}(t) = \int_0^t \tilde{\gamma}^{\operatorname{eff}}(s) \dif s\quad \text{and}\quad  \Lambda(t) = \int_0^t \tilde{\gamma}^{\operatorname{eff}}(s)\tilde{\lambda}^{\operatorname{eff}}(s) \dif s.
    \end{equation}
    As derived in detail in Appendix~\ref{sec:lsq_analysis}, $ \mathrsfs{R}(t)$ and $\widehat{\mathrsfs{B}}_{44}(t) $ solve a system of coupled nonlinear Volterra equations
    \begin{equation}\label{eq:def:R(t):det:equiv:lr:reform:main}
        \begin{aligned}
            \mathrsfs{R}(t)&= \frac{1}{2}\langle K , (e^{-K\tilde{\Gamma}(t)  - \Lambda(t)}\left(X_0 + \left(\int_0^t  K\tilde{\gamma}^{\operatorname{eff}}(s) e^{K \tilde{\Gamma}(s) + \Lambda(s)}\dif s\right) X^\star \right)- X^{\star})^{\otimes 2}\rangle + \frac{\eta^2}{2}
            \\
            &
            \qquad + \frac{\gamma}{d} \int_0^t  \tilde{\gamma}^{\operatorname{eff}}(s)  \left(\mathrsfs{R}(s) + \frac{\tilde{\lambda}^{\operatorname{eff}}(s)}{2}\widehat{\mathrsfs{B}}_{44}(s)\right) \tr \big ( K^2 e^{-2K(\tilde{\Gamma}(t) - \tilde{\Gamma}(s)) -2(\Lambda(t) -\Lambda(s))} \big ) \, \dif s.
        \end{aligned}
    \end{equation}
    and 
    \begin{equation}\label{eq:def:l2norm2:det:equiv:lr:reform:main}
        \begin{aligned}
            \widehat{\mathrsfs{B}}_{44}(t) 
            &= \left\| e^{-K\tilde{\Gamma}(t)  -\Lambda(t)}\left(X_0 +\left(\int_0^t  K\tilde{\gamma}^{\operatorname{eff}}(s) e^{K \tilde{\Gamma}(s) + \Lambda(s)}\dif s\right) X^\star\right)\right\|^2
            \\
            &
            \qquad+ \frac{2\gamma}{d} \int_0^t\tilde{\gamma}^{\operatorname{eff}}(s)  \left(\mathrsfs{R}(s) + \frac{\tilde{\lambda}^{\operatorname{eff}}(s)}{2}\widehat{\mathrsfs{B}}_{44}(s)\right) \tr \big ( K e^{-2K(\tilde{\Gamma}(t) - \tilde{\Gamma}(s)) -2(\Lambda(t) -\Lambda(s))} \big ) \, \dif s.
        \end{aligned}
        \end{equation}
        See Appendix~\ref{sec:lsq_analysis} for more details. Setting $\delta=0$, then $\mathrsfs{R}(t)$ satisfies the linear convolution Volterra equation derived in \citep{collinswoodfin2023hitting, CollinsWoodfinPaquette01, PPAP01}. The central object of the dynamics is $q(t)$ which satisfies the following nonlinear ODE:
 \begin{equation}\label{eq:ode:q(t):main}
     \begin{aligned}
         \frac{\dif}{\dif t} (q(t))^2 
         &
         = 
         \frac{\gamma}{\mathrsfs{R}(t)} (q(t))^2\left( 1 + \delta\sqrt{\frac{2}{\pi}}q(t) \right) \sum_{j=1}^d \lambda_j^2 \left(\mathrsfs{V}_{j,11}(t) - 2 \mathrsfs{V}_{j,12}(t) + \mathrsfs{V}_{j,22}(t)\right)
         \\
         &
         +
         \frac{\gamma}{\mathrsfs{R}(t)} (\delta^2 (q(t))^2 - 1)  \sum_{j=1}^d \lambda_j\left(\mathrsfs{V}_{j,11}(t) - \mathrsfs{V}_{j,12}(t) \right)-2 \gamma q(t)\left(\delta^2 q(t) + \delta \sqrt{\frac{2}{\pi}}\right)
         \\
         &
         +
         \frac{\gamma^2}{d} \left(1+ 2 \delta \sqrt{\frac{2}{\pi}} q(t) + \delta^2 (q(t))^2\right)(\Tr(K) - \Tr(K^2) (q(t))^2).
     \end{aligned}
 \end{equation}
Assuming that $X_0 \not = X^\star$ then we show that $q(t)$ is bounded above and below for all time $t>0$.
\begin{proposition}\label{prop:q(t):bounded}
    Suppose $\mathrsfs{R}(0), \widehat{\mathrsfs{B}}_{44}(0) >0$ then for all $t\geq 0$, $\mathrsfs{R}(t) >0$ from which it follows $q(t)$ is uniformly bounded by constants for all time $t\geq 0$
    \begin{equation}
        \min(q(0), q_L) \leq q(t) \leq \max(q(0), q_U),
    \end{equation}
    where $q_L$ and $q_U$ are constants which depend on $\Tr(K)$, $\Tr(K^2)$, $\gamma$ and $\delta$ and are defined in Proposition~\ref{prop:q(t):bound:cond}. Suppose further that there exists a $c >0$ such that $\frac{\Tr(K^2)}{d} \geq c$ uniformly in $d$ and $0 < c_0 \leq q(0) \leq C_0 < \infty$ where $c_0, C_0 >0$ are independent of $d$ then $q_L$ and $q_U$ are uniform bounds with respect to $d$.
\end{proposition}
See Appendix~\ref{proof:prop:q(t):bounded} for the proof of Proposition~\ref{prop:q(t):bounded}. Given the assumptions presented in Proposition~\ref{prop:fixed:point:main} and $q(t) \to q_\infty$, we characterize $q_\infty$ as the unique admissible fixed point of \eqref{eq:fixed:point:main}. 

\begin{figure}[t!]
    \centering
    \includegraphics[width=0.30\textwidth]{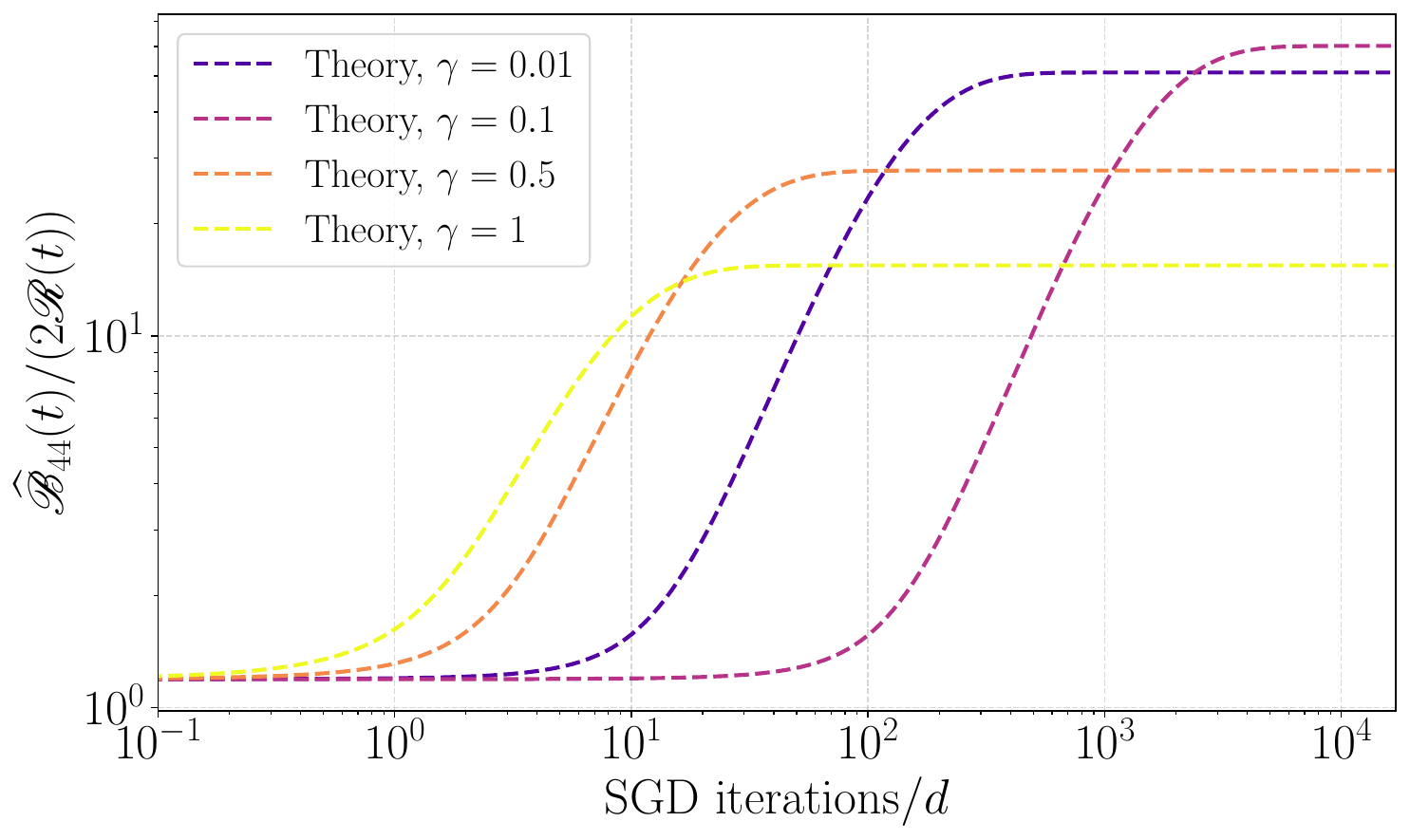}
    \includegraphics[width=0.30\textwidth]{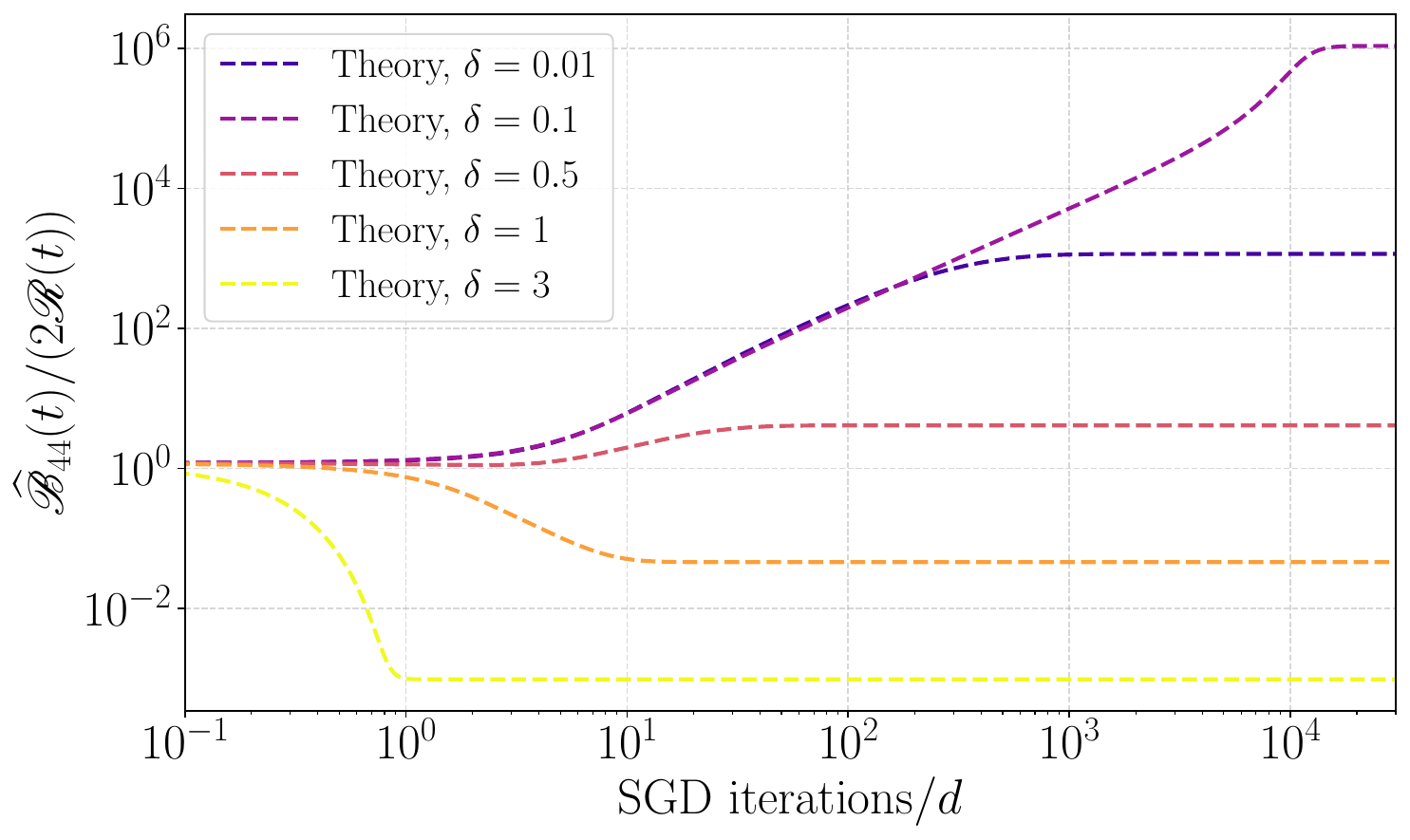}\includegraphics[width=0.30\textwidth]{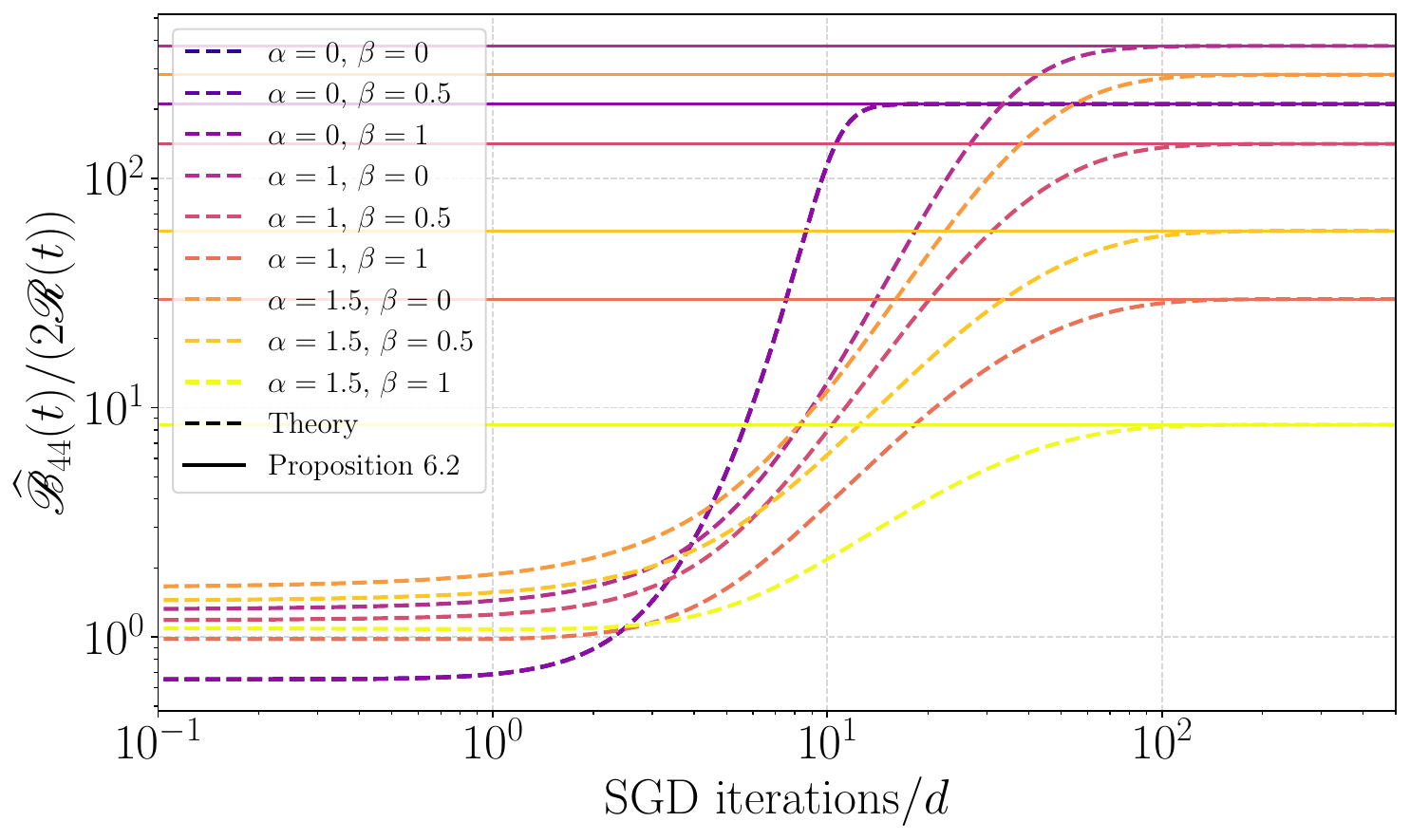} 
    \caption{\textbf{Numerical evidence that $q(t) \defas \sqrt{\frac{\widehat{\mathrsfs{B}}_{44}(t)}{2\mathrsfs{R}(t)}}$ converges} for noiseless $\ell_2$-adversarial least squares. The first (left) and second (middle) plots provide evidence that $q(t)$ converges for a variety of constant learning rates $\gamma$ and $\delta$. Here we fix either parameter and vary the other according to the values presented in the legends. The third plot provides evidence that $q(t)$ converges for a variety of $X^\star$ and covariances $K$ under different power law setups (See Assumption~\ref{ass:powerlaw:main}) and validates that the fixed point equation in Proposition~\ref{prop:fixed:point:main} predicts the limiting value $q_\infty$. See Appendix~\ref{app:captions} for simulation details.}
    \label{fig:q}
\end{figure}

     \begin{proposition}\label{prop:fixed:point:main}
             Suppose $\mathrsfs{R}(0), \widehat{\mathrsfs{B}}_{44}(0) >0$, $(X^\star)^\top K X^\star> 0$, $\sup_{t\geq 0} \mathrsfs{R}^\Adv(t)< \infty$, $\delta >0$ and $q(t) \to q_{\infty} \in (0,\infty)$. Let $r_\star\approx 8.35\times 10^{-5}$. If $\delta \leq \sqrt{\frac{\lambda_{\operatorname{min}}^+}{ r_\star}}$ then $q_\infty$ is the unique fixed point satisfying the fixed point equation
         \begin{equation}\label{eq:fixed:point:main}
             q^2_\infty = (1 - G(q_\infty) C_{\mathrsfs{R}}(q_\infty)) \frac{L_\infty(q_\infty)}{2 F_\infty(q_\infty)} + \frac{1}{2}G(q_\infty)C_\mathrsfs{B}(q_\infty)
         \end{equation}
         and the admissibility condition $G(q_\infty)C_{\mathrsfs{R}}(q_\infty) < 1$ for $G(q) = 1+ 2 \delta \sqrt{\frac{2}{\pi}} q + \delta^2 q^2$ and $C_\mathrsfs{R}(q) = \frac{\gamma}{2d}\sum_{j=1}^d \frac{\lambda_j^2}{D_j(q)}$ with $D_j(q) = \lambda_j(1+ \delta \sqrt{\frac{2}{\pi}}q)  + \delta^2 + \delta \sqrt{\frac{2}{\pi}}\frac{1}{q}$. The limiting values of $\mathrsfs{R}$, $\widehat{\mathrsfs{B}}_{44}$ and $\mathrsfs{R}^\Adv$ are given by
         \begin{equation*}
             \begin{gathered}
                 \mathrsfs{R}(\infty) = \frac{F_{\infty}(q_\infty)}{1- G(q_\infty)C_{\mathrsfs{R}}(q_\infty)},\quad 
                 \widehat{\mathrsfs{B}}_{44}(\infty) = 2 q^2_\infty  \mathrsfs{R}(\infty),\quad
                 \text{and}\quad \mathrsfs{R}^\Adv(\infty) = G(q_\infty)\mathrsfs{R}(\infty),
             \end{gathered}
         \end{equation*}
     where $C_{\mathrsfs{B}}(q) = \frac{\gamma}{d}\sum_{j=1}^d \frac{\lambda_j}{D_j(q)}$, $F_\infty(q) = \frac{1}{2}\sum_{j=1}^d \lambda_j \left(\frac{\tilde{\lambda}^{\operatorname{eff}}(q)}{\lambda_j +\tilde{\lambda}^{\operatorname{eff}}(q)}\right)^2(x_j^\star)^2$, $ L(\infty) = \sum_{j=1}^d  \left(\frac{\lambda_j}{\lambda_j + \tilde{\lambda}^{\operatorname{eff}}(\infty)}\right)^2( x_j^\star)^2$, and $\tilde{\lambda}^{\operatorname{eff}}(q) = \frac{\delta^2 + \delta \sqrt{\frac{2}{\pi}} \frac{1}{q}}{1 + \delta \sqrt{\frac{2}{\pi}}q }$.
     \end{proposition}
     \begin{remark}
        Note that a general condition for the fixed point to be unique and admissible can be found in Proposition~\ref{prop:fixed:point} and that the condition provided in Proposition~\ref{prop:fixed:point:main} is a simplication. We omit this condition from the main text given its lengthiness. We also note that we could predict the limiting values of a larger set of statistics of $q(t)$ (e.g. distance to optimality) than just $\mathrsfs{R}$, $\widehat{\mathrsfs{B}}_{44}$ and $\mathrsfs{R}^\Adv$. See Figure~\ref{fig:q} for numerical evidence that $q(t)$ converges and the fixed point equation \eqref{eq:fixed:point} predicts its limit.
     \end{remark}
     The fixed point equation \eqref{eq:fixed:point:main} is hard to interpret. Even for the simpler setting of $K=\Id_d$, \eqref{eq:fixed:point:main} remains a polynomial of degree 6 in $q$. To provide insights on the limiting values $q_\infty$ and $\tilde{\lambda}(\infty)$, we derive $q_\infty$ as $\gamma \to 0$. Note that we recover the results of \citep[Proposition 1]{xing2021same_minimizer} in the limit $\gamma \to 0$.
     \begin{proposition}\label{prop:fixed:point:gamma:to:0:main}
         Suppose the assumptions of Proposition~\ref{prop:fixed:point:main} hold. Define $J_+ \defas \{j: \lambda_j >0 \}$ and define the two thresholds
         \begin{equation}
             \delta_- = \sqrt{\frac{2}{\pi}}\sqrt{\frac{\sum_{j \in J_+} (x_j^\star)^2}{\sum_{j \in J_+}\frac{ (x_j^\star)^2}{\lambda_j}}}
             \quad \text{and}\quad  \delta_+ = \sqrt{\frac{\pi}{2}}  \sqrt{\frac{\sum_{j \in J_+} \lambda_j^2 (x_j^\star)^2}{ \sum_{j \in J_+}\lambda_j (x_j^\star)^2}}.
         \end{equation}
         For $\gamma >0$, let $q_{\gamma}(\infty)$ be an admissible solution of the fixed point equation of Proposition~\ref{prop:fixed:point:main}. Then as $\gamma \to 0$
         \begin{equation}\label{eq:result:q:limits}
             \begin{gathered}
             \text{if $\delta \leq \delta_-$ then}\quad q_{\gamma}(\infty) \to \infty\quad\text{and}\quad \tilde{\lambda}(q_{\gamma}(\infty)) \to \tilde{\lambda}_\star=0,\\
             \text{if $\delta_-< \delta < \delta_+$ then}\quad q_{\gamma}(\infty) \to q_{\star} \in (0,\infty)\quad\text{and}\quad \tilde{\lambda}(q_{\gamma}(\infty)) \to \tilde{\lambda}_\star \in (0,\infty),\\
             \text{and if $\delta \geq \delta_+$ then}\quad q_{\gamma}(\infty) \to 0\quad\text{and}\quad \tilde{\lambda}(q_{\gamma}(\infty)) \to \tilde{\lambda}_\star=\infty.
             \end{gathered}
         \end{equation}
         If $\delta_-< \delta < \delta_+$, then $(q_\star, \tilde{\lambda}_{\star})$ is the unique solution of 
         \begin{equation}\label{eq:fixed:point:gamma:0:main}
             q_\star = \sqrt{ \frac{\sum_{j \in J_+} \frac{\lambda_j^2}{\left(\lambda_j +  \tilde{\lambda}_\star\right)^2} (x_j^\star)^2}{ \sum_{j \in J_+}  \lambda_j \frac{ \tilde{\lambda}_\star^2}{
             \left(\lambda_j +  \tilde{\lambda}_\star\right)^2 
                } (x_j^\star)^2}}
             \quad\text{and}\quad \tilde{\lambda}_{\star}  = \frac{\delta^2 +\delta \sqrt{\frac{2}{\pi}} \frac{1}{q_\star}}{1+ \delta \sqrt{\frac{2}{\pi}} q_\star}.
         \end{equation}
         If $K= \Id_d$ then $\delta_- = \sqrt{\frac{2}{\pi}}$, $\delta_+ = \sqrt{\frac{\pi}{2}}$ and for $\delta_- < \delta < \delta_+$
         \begin{equation}
             q_{\star} = \frac{1- \delta \sqrt{\frac{2}{\pi}}}{\delta \left(\delta -\sqrt{\frac{2}{\pi}}\right)} \quad \text{and}\quad \tilde{\lambda}_{\star} = \frac{\delta^2 - \delta \sqrt{\frac{2}{\pi}}}{1- \delta \sqrt{\frac{2}{\pi}}}.
         \end{equation}
     \end{proposition}
     See Appendix~\ref{proof:prop:q(t):bounded} for the proofs of Propositions~\ref{prop:fixed:point:main} and~\ref{prop:fixed:point:gamma:to:0:main}.

    \subsection{Power law behaviour at $d \to \infty$}\label{sect:plaw}
    We will now analyze the dynamics of the deterministic equivalents $\mathrsfs{R}(t)$, $\mathrsfs{R}^\Adv(t)$ and $\widehat{\mathrsfs{B}}_{44}(t)$ in the limit $d\to \infty$ under the following power law assumption. See Appendix~\ref{app:power:law} for the proofs of the results presented in this section. Recall $(\lambda_j, \omega_j)$ are the eigenvalue-eigenvector pairs of $K$.
    \begin{assumption}\label{ass:powerlaw:main}
        Suppose $\alpha \geq 0$, $\beta \geq  0$. Define $x_{0,j} \defas \langle X_0, \omega_j\rangle$ and $x^\star_j\defas \langle X^\star, \omega_j\rangle$ such that
        \[
            \lambda_j = \left(\frac{j}{d}\right)^\alpha \quad \text{and} \quad (x^\star_j)^2 = \frac{1}{d} \lambda_j^{-\beta}.
        \]
    \end{assumption}
    \begin{remark}\label{rem:equiv:powerlaw}
        This power law parametrization is inspired by the parametrizations presented in \citep{collinswoodfin2024highline,collinswoodfin2025Exact}. The parameter $\beta_{\operatorname{HL}}$ of \citep{collinswoodfin2024highline} satisfies the following relation with our parameter $\alpha$:
        \[
            \alpha = \frac{1}{1-\beta_{\operatorname{HL}}}.
        \]
        Hence for large $d$, the eigenvalues $\{\lambda_j\}_{j=1}^d$ follow a power law distribution with density function $f(\lambda) = \frac{1}{\alpha} \lambda^{\frac{1}{\alpha}-1}$ and support $(0,1]$. $\alpha = 0$ corresponds to isotropic covariance and $\alpha > 1$ corresponds to power law setup of \citep{collinswoodfin2024highline} with $0 < \beta_{\operatorname{HL}}< 1$.
        See Remark 4 in \citep{collinswoodfin2025Exact} for more details on the role of $\alpha$ or Section 6 in \citep{pmlr-v162-wei22a} for a discussion on power laws. 
    \end{remark}
    Let $J_+ \defas \{j : \lambda_j >0\}$. From \citep[Proposition 1]{xing2021same_minimizer} and Proposition~\ref{prop:fixed:point:gamma:to:0:main}, we see that $X^{\star, \Adv}$ admits the following representation depending on $\delta$
    \begin{equation}
        X^{\star, \Adv} = \begin{cases}
            \sum_{j \in J_+} \frac{\lambda_j}{\lambda_j + \tilde{\lambda}_\star} \langle X^\star, \omega_j\rangle \omega_j & \text{if} \quad \delta < \delta_+  \\
            0  & \text{if}\quad  \delta \geq \delta_+ 
        \end{cases}
    \end{equation}
    Thus, depending on the value of $\delta$, we can view $X^{\star, \Adv}$ as a compromise between minimizing the least squares risk $\mathcal{R}$ and $\|\cdot\|^2$. For $\delta \leq \delta_- $ we have $\mathcal{R}(X^{\star, \Adv}) = 0$ and for $\delta \geq \delta_+$ we have $\|X^{\star, \Adv}\|^2 = 0$. A natural question is then 
    \begin{center}
        \begin{minipage}{0.85\textwidth}
        \begin{mdframed}[
          linewidth=0.8pt,
          innerleftmargin=10pt,
          innerrightmargin=10pt,
          innertopmargin=8pt,
          innerbottommargin=8pt
        ]
        \centering
         Can SGD lead to $\mathcal{R}^\Adv(X_k) >0$ and either $\mathcal{R}(X_k) = 0$ or $\|X_k\|^2 =0$ when $k\to\infty$?\label{question}
        \end{mdframed}
        \end{minipage}
    \end{center}
    For finite $d\geq 1$, whenever $q(t) \to q_\infty$, Proposition~\ref{prop:fixed:point:main} states that $\mathrsfs{R}(\infty) > 0 $ and $\widehat{\mathrsfs{B}}_{44}(\infty) >0$. The following propositions and corollaries lead a formulation of Proposition~\ref{prop:fixed:point:main} as $d \to \infty$. We introduce some notation first. Define $\ell(u) \defas \Lambda(t)$ then $\ell'(u) = \tilde{\lambda}^{\operatorname{eff}}(\tilde{\Gamma}^{-1}(u))$ and define $F^\Adv$ as the sum of the two first terms in \eqref{eq:def:R(t):det:equiv:lr:reform:main} and \eqref{eq:def:l2norm2:det:equiv:lr:reform:main}:
    \begin{equation*}
        \begin{aligned}
            F^\Adv(u) &
            =
            \frac{1}{2}\sum_{j=1}^d (\lambda_j + \delta^2 )\bigg(e^{-\lambda_j u - \ell(u) } \left( x_{0,j} - x^\star_j\right) 
            \\
            &
            + x^\star_j \left(\frac{\delta^2}{\lambda_j + \delta^2 } - \left(\int_0^{u} \ell'(s) e^{-\lambda_j (u-s) - (\ell(u)-\ell(s))} \dif s\right)\right) \bigg)^2
            + \frac{1}{2} \sum_{j=1}^d \frac{\delta^2 \lambda_j}{ \lambda_j + \delta^2} (x^\star_j)^2.
        \end{aligned}
    \end{equation*} 
    It turns out that under the power law Assumption~\ref{ass:powerlaw:main}, the adversarial risk $\mathrsfs{R}^\Adv(t)$ is bounded above and below by constants independent of time $t$ when $d\to \infty$. Note that we can obtain a similar result for $\alpha = 0$ for finite $d$.
        \begin{proposition}\label{prop:isotropic_Rasymp}
            Suppose $\mathrsfs{R}(0), \widehat{\mathrsfs{B}}_{44}(0) >0$, Assumption~\ref{ass:powerlaw:main} is satisfied with $\alpha= 0$, $\gamma < 1$ then there exists constants $C_1, C_2>0$ such that for all $t\geq 0$
            \[
            \begin{aligned}
                C_1 \leq F^\Adv(\tilde{\Gamma}(t)) \leq \mathrsfs{R}^\Adv(t) \leq &2F^\Adv(\tilde{\Gamma}(t))
                \\
                &
                +  4\gamma^2(1+\delta^2)\int_0^t e^{-2\gamma(1+\delta^2)(1-\gamma)(t-s)} F^\Adv(\tilde{\Gamma}(s)) \dif s \leq C_2,
            \end{aligned}
            \]
            where $C_1  = \frac{\delta^2}{2(1+\delta^2)}$, $ C_2 = 2 C_F \frac{1+\gamma}{1-\gamma}$ with $C_{F} = (1+\delta^2) (2 \|X_0\|^2 +3)  + \frac{\delta^2}{2(1+\delta^2)}$.
            \end{proposition}
        \begin{proposition}\label{prop:powerlaw_Rasymp}
            Suppose $\mathrsfs{R}(0), \widehat{\mathrsfs{B}}_{44}(0) >0$, Assumption \ref{ass:powerlaw:main} is satisfied with $\alpha >0$, $\alpha \beta<1$ and $ \gamma < \tfrac{1}{2\left(\delta^2 + \frac{1}{\alpha + 1}\right)}$ with $\delta >0$. Then there exists constants $C_1, C_2>0$ such that for $t\geq 1$
            \[
                C_1 \leq F^\Adv(\tilde{\Gamma}(t)) \leq \mathrsfs{R}^\Adv(t) \leq  2F^\Adv(\tilde{\Gamma}(t)) + 2C\int_0^t\overline{\mathscr{K}}(t-s) F^\Adv(\tilde{\Gamma}(s)) \dif s \leq C_2,
            \]
            where $C$ and $\overline{\mathscr{K}}$ are defined in Lemma~\ref{lem:R_up_low_bound:2} and where $C_1  = \frac{\delta^2}{2(1+\delta^2)}\frac{1}{1+ \alpha(1-\beta)}$ and $ C_2 =2(1+\delta^2) \|X_0\|^2 + \left(3(1+\delta^2) + \frac{1}{2}\right)\frac{1}{1-\alpha \beta}$.
            \end{proposition}
            Using Propositions~\ref{prop:isotropic_Rasymp}~and~\ref{prop:powerlaw_Rasymp}, we can derive bounds for $\mathrsfs{R}(t)$ and $\widehat{\mathrsfs{B}}_{44}(t)$ in terms of $L(\tilde{\Gamma}(t))$, $F(\tilde{\Gamma}(t))$ and $F^\Adv(\tilde{\Gamma}(t))$.
            \begin{corollary}\label{cor:powerlaw_R_B_44_bounded:isotropic}
                Suppose $\mathrsfs{R}(0), \widehat{\mathrsfs{B}}_{44}(0) >0$, Assumption \ref{ass:powerlaw:main} is satisfied and $\gamma < 1$ with $\delta >0$. Then there exists constants $C_1, C_2, C_3, C_4>0$ such that for $t\geq 0$, $C_1 \leq \mathrsfs{R}(t) \leq C_2$ and $C_3 \leq \widehat{\mathrsfs{B}}_{44}(t) \leq C_4$. Furthermore, for $t\geq 0$, $\mathrsfs{R}(t)$ and $\widehat{\mathrsfs{B}}_{44}(t)$ satisfy the following inequalities $F(\tilde{\Gamma}(t)) \leq \mathrsfs{R}(t)$, $L(\tilde{\Gamma}(t)) \leq \widehat{\mathrsfs{B}}_{44}(t) $ and
                \[
                \begin{aligned}
                     \mathrsfs{R}(t) &\leq F(\tilde{\Gamma}(t)) 
                    + 2\gamma^2 \int_0^t e^{-2\gamma(1 + \delta^2)(t-s)} F^\Adv(\tilde{\Gamma}(s)) \dif s 
                    \\
                    &\qquad + 4\gamma^4(1+\delta^2)\int_0^t  e^{-2\gamma(1+\delta^2)(2-\gamma)(t-s)}  F^\Adv(\tilde{\Gamma}(s)) \dif s,
                    \\
                    \widehat{\mathrsfs{B}}_{44}(t)  &\leq L(\tilde{\Gamma}(t)) 
                    + 4\gamma^2 \int_0^t e^{-2\gamma(1 + \delta^2)(t-s)} F^\Adv(\tilde{\Gamma}(s)) \dif s 
                    \\
                    &\qquad + 8\gamma^4(1+\delta^2)\int_0^t  e^{-2\gamma(1+\delta^2)(2-\gamma)(t-s)}  F^\Adv(\tilde{\Gamma}(s)) \dif s.
                \end{aligned}
                \]
                \end{corollary}
        
                            \begin{corollary}\label{cor:powerlaw_R_B_44_bounded:plaw}
                                Suppose $\mathrsfs{R}(0), \widehat{\mathrsfs{B}}_{44}(0) >0$, Assumption \ref{ass:powerlaw:main} is satisfied and $\gamma$ satisfies the threshold of Proposition~\ref{prop:powerlaw_Rasymp} with $\delta >0$. Then there exists constants $C_1, C_2, C_3, C_4>0$ such that for $t\geq 0$, $C_1 \leq \mathrsfs{R}(t)$, $C_3 \leq \widehat{\mathrsfs{B}}_{44}(t)$ and for $t\geq 1$, $\mathrsfs{R}(t) \leq C_2$ and $\widehat{\mathrsfs{B}}_{44}(t) \leq C_4$. Furthermore, for $t\geq 0$, $\mathrsfs{R}(t)$ and $\widehat{\mathrsfs{B}}_{44}(t)$ satisfy the following inequalities $F(\tilde{\Gamma}(t)) \leq \mathrsfs{R}(t)$, $L(\tilde{\Gamma}(t)) \leq \widehat{\mathrsfs{B}}_{44}(t) $ and for $t\geq 1$
                                \[
                                \begin{gathered}
                                     \mathrsfs{R}(t) \leq F(\tilde{\Gamma}(t)) 
                                    + \int_0^t \overline{\mathscr{K}}_{\mathrsfs{R}}(t-s ) F^\Adv(\tilde{\Gamma}(s)) \dif s +C \int_0^t  (\overline{\mathscr{K}}_{\mathrsfs{R}} *\overline{\mathscr{K}})(t-s) F^\Adv(\tilde{\Gamma}(s)) \dif s,
                                    \\
                                    \widehat{\mathrsfs{B}}_{44}(t)  \leq L(\tilde{\Gamma}(t)) 
                                    + \int_0^t \overline{\mathscr{K}}_{\mathrsfs{B}}(t-s ) F^\Adv(\tilde{\Gamma}(s)) \dif s +C \int_0^t  (\overline{\mathscr{K}}_{\mathrsfs{B}} *\overline{\mathscr{K}})(t-s) F^\Adv(\tilde{\Gamma}(s)) \dif s,
                                \end{gathered}
                                \]
                                where $C$ and $\overline{\mathscr{K}}$ are defined in Lemmas~\ref{lem:R_up_low_bound:2}~and~\ref{lem:powerlaw:Linearreg}, $\overline{\mathscr{K}}_{\mathrsfs{R}}(x) \defas 2\gamma^2  \int_0^1 y^{2\alpha} e^{-2\gamma(y^{\alpha} + \delta^2)x} \dif y$ and $\overline{\mathscr{K}}_{\mathrsfs{B}}(x) \defas 4\gamma^2 \int_0^1   y^{\alpha} e^{-2\gamma(y^{\alpha} + \delta^2)x} \dif y$.
                                \end{corollary}
                            From Corollaries~\ref{cor:powerlaw_R_B_44_bounded:isotropic}~and~\ref{cor:powerlaw_R_B_44_bounded:plaw} and given that SGD is close to $\mathrsfs{R}(t)$ and $\widehat{\mathrsfs{R}}_{44}(t)$ when $d\to \infty$, we see that $\mathcal{R}(X_k) > 0$ and $\|X_k\|^2 >0$ for all $k >0$. In Proposition~\ref{prop:fixed:point:main}, we provided a general sufficient condition under which, if $q(t) \to q_\infty$ as $t\to \infty$, $q_\infty$ is the unique solution to the fixed point equation \eqref{eq:fixed:point:main} for each fixed $d\geq 1$. In the following Proposition, we provide sufficient conditions for different regimes of $\alpha \geq 0$ when $d\to \infty$.
                            
                            \begin{proposition}\label{prop:fixed:point:plaw}
                                Suppose the assumptions of Proposition~\ref{prop:fixed:point:main} hold. Suppose $\alpha \geq 0$ and $\alpha \beta < 1$. Let $q_d(t)$ denote the term $q(t)$ associated with dimension $d$ and $0 < c_0 \leq q_d(0) \leq C_0 < \infty$ uniformly in $d$ and where $c_0, C_0 >0$ are independent of $d$. For each $d \geq 1$, assume  as $t\to \infty$ that $q_d(t)$ converges to $q_{d,\infty}$ which solves the fixed point equation presented in Proposition~\ref{prop:fixed:point:main}. Assume moreover that one of the following conditions holds for $r_\star \approx 8.35\times 10^{-5}$ and $s_\star \approx1.77 \times 10^{-3}$:
                                \[
                                \begin{gathered}
                                   \text{for}\quad\alpha = 0: \delta \leq \sqrt{\frac{1}{ r_\star}}, \quad \text{for}\quad
                                   0< \alpha < 1: \delta < \sqrt{\frac{1-\alpha}{ r_\star}} ,\\
                                   \text{or for} \quad \alpha \geq 1:  \delta < \sqrt{\frac{\kappa_{\alpha} - s_\star}{r_\star}} \quad \text{with}\quad \kappa_{\alpha} = \frac{\alpha+ 1}{\alpha+2} \frac{1}{2\alpha + 1}\left(\frac{2\alpha}{2\alpha + 1}\right)^{2\alpha} > s_\star.
                                \end{gathered}
                                \]
                            Note that for $\alpha \geq 1$, we have $\kappa_{\alpha} > s_{\star}$ if $\alpha < \alpha_0 \approx 102.87$. Then $q_{d,\infty}$ converges as $d\to \infty$ to $q_{\infty} \defas \lim_{d\to \infty}q_{d, \infty}$ which is the unique fixed point satisfying 
                                \begin{equation}
                                    q^2_\infty = (1 - G(q_\infty) C_{\mathrsfs{R}}(q_\infty)) \frac{L_\infty(q_\infty)}{2 F_\infty(q_\infty)} + \frac{1}{2}G(q_\infty)C_\mathrsfs{B}(q_\infty)
                                \end{equation}
                                and the admissibility condition $G(q_\infty)C_{\mathrsfs{R}}(q_\infty) < 1$ where
                                \[
                                \begin{gathered}
                                    C_\mathrsfs{R}(q) = \frac{\gamma}{2} \int_0^1 \frac{y^{2\alpha}}{y^{\alpha}(1+ \delta \sqrt{\frac{2}{\pi}}q)  + \delta^2 + \delta \sqrt{\frac{2}{\pi}}\frac{1}{q}} \dif y, \quad C_{\mathrsfs{B}}(q) = \gamma\int_0^1 \frac{y^{\alpha}}{y^{\alpha}(1+ \delta \sqrt{\frac{2}{\pi}}q)  + \delta^2 + \delta \sqrt{\frac{2}{\pi}}\frac{1}{q}} \dif y,
                                    \\
                                    _{\infty}(q) = \int_0^1  y^{-\alpha \beta}\left(\frac{y^{\alpha}}{y^{\alpha} + \tilde{\lambda}^{\operatorname{eff}}(q)}\right)^2 \quad \text{and}\quad F_\infty(q) = \frac{1}{2}\int_0^1  y^{\alpha(1-\beta)} \left(\frac{\tilde{\lambda}^{\operatorname{eff}}(q)}{y^{\alpha} +\tilde{\lambda}^{\operatorname{eff}}(q)}\right)^2 \dif y.
                                \end{gathered}
                                \] 
                            Here $\mathrsfs{R}$, $\widehat{\mathrsfs{B}}_{44}$, $\mathrsfs{R}^\Adv$ and $\tilde{\lambda}^{\operatorname{eff}}(q)$ are identical to Proposition~\ref{prop:fixed:point:main}.
                            \end{proposition} 
                            The three conditions derived in Proposition~\ref{prop:fixed:point:plaw} are sufficient conditions under which power law covariance with $\alpha \geq0$ converges to the unique admissible fixed point. The isotropic condition is sharp. For $0 < \alpha < 1$, the condition gives a dimension-free bound. The $\kappa_{\alpha}$ condition is more conservative near $\alpha = 0$, but remains useful beyond $\alpha = 1$.
    
\section{SGD and AdvHSGD are Approximate Solutions}\label{sect:SGD:AdvHSGD:approx:sols}
\paragraph{Outline of the section.} In this section, we will prove the main result Theorem~\ref{thm:main:conc:St:general:with:stop:adv}. The section is divided as follows. Section~\ref{sect:doob:decomp} introduces the Doob decompositions of SGD and AdvHSGD and provides an overview of the necessary steps to prove SGD and AdvHSGD are approximate solutions to the integro-differential equation \eqref{eq:ODE_resolvent_2}. Section~\ref{sect:approx:sols:theory} introduces the theory of approximate solutions to the integro-differential equation \eqref{eq:ODE_resolvent_2}. It also includes proofs on the stability of approximate solutions and the proofs that AdvHSGD and SGD are approximate solutions to \eqref{eq:ODE_resolvent_2}. Section~\ref{sect:remove:stopping:time} contains a rigorous justification for the removal of the stopping time in Theorems~\ref{thm:main:conc:St:SGD} and \ref{thm:main:conc:St:HSGD}. Section~\ref{sect:conc:proof:any:stat} contains analogous results to Theorems \ref{thm:main:conc:St:adv} and \ref{thm:main:conc:St:general:with:stop:adv} whenever Assumption~\ref{assumption:norm:Xstar:adv} fails. Section~\ref{sect:general:adv:concentration} introduces the theory of approximate solutions to the integro-differential equation \eqref{eq:ODE_resolvent_2:align} and the rigorous justification for the removal of the stopping time in Theorems~\ref{thm:main:conc:St:adv} and \ref{thm:main:conc:St:general:with:stop:adv}.

Finally, Appendix~\ref{sect:prelim:proofs} contains the details of the derivations of the Doob decompositions for AdvHSGD and SGD and supporting results for the proofs that AdvHSGD and SGD are approximate solutions to \eqref{eq:ODE_resolvent_2}. We choose to delegate these results to the Appendix as they are quite computational.

\subsection{The Doob Decomposition of SGD and AdvHSGD}\label{sect:doob:decomp}
In order to effectively compare the discrete-time iterates of SGD to its continuous time limit, recall we introduce the following time parametrization
\begin{equation*}
    k\text{ iterations of SGD} = \lfloor td \rfloor,\quad \text{for } t\in \R^+ \text{ the continuous time parameter}.
\end{equation*}
For simplicity of notation, we denote $X_{td} \defas  X_{\lfloor td \rfloor}$, $\hat{X}_{td} \defas \hat{X}_{\lfloor td \rfloor}$ and $W_{td}  \defas W_{\lfloor td \rfloor}$. The goal of Section~\ref{sect:SGD:AdvHSGD:approx:sols} is to show AdvHSGD and SGD approximately solve \eqref{eq:ODE_resolvent_2}. Before diving into the theory of approximate solutions to integral-dfferential equations, we provide intuition behind the existence of the deterministic limits for statistics of the iterates of SGD and AdvHSGD. 

Since the entries of  $Z(W_k, z)$ are quadratic statistics of the entries of $W_k$, a natural starting point is to consider the Taylor expansion of a quadratic function of $W_{k+1}$ centered at $W_k$. Formally, let a quadratic function denote a $C^{\infty}$ function $\varphi \, : \, \R^{d} \to \mathbb{R}$ such that $(\nabla^{(j)} \varphi)(X) \equiv 0$ for all $X \in  \R^{d}$ and $j \ge 3$. Here we define $\nabla^{(j)} \varphi \defas \nabla^{(j)}_X \varphi$. 

Recall we defined $\hat{X}_{k} \defas X_k \oplus X^{\star}$, $W_{k} \defas \hat{X}_{k} \oplus \mu$, $a_{k+1, I_{k+1}} = \sqrt{K_{I_{k+1}}}v_{k+1} + \mu_{I_{{k+1}}}$ for $v_{k+1} \sim \mathcal{N}(0, \Id_d)$, $I_{k+1} \in \{1,2\}$, $r_{k, I_{k+1}} = \hat{X}_k^\top a_{k+1, I_{k+1}}$ and
    \begin{equation}
        \begin{gathered}
        g_{k, I_{k+1}} = r_{k, I_{k+1}} + \begin{bmatrix}
            \delta s_{k, I_{k+1}}\|X_k\| & 0 
        \end{bmatrix}^\top,\\
        s_{k, {I_{k+1}}} = \operatorname{argmax}_{|s| \leq  1}f_{I_{k+1}}( X_k^\top a_{k+1, I_{k+1}}+ \delta s \|X_k\|; \epsilon_{k+1}).
    \end{gathered}
    \end{equation}
By Taylor expansion, setting $\Delta_k \defas f'_{I_{k+1}}( g_{k, I_{k+1}}, \epsilon_{k+1}) (a_{k+1, I_{k+1}} + \delta s_{k, I_{k+1}} \frac{X_k}{\|X_k\|}) + \lambda X_k$, we have
\begin{equation} \label{eq:taylor_theorem}
    \begin{aligned}
    \varphi(X_{k+1}) 
    =
    \varphi(X_{k}) - \frac{\gamma_k}{d} \ip{\nabla \varphi(X_k), \Delta_k}
    + 
     \frac{\gamma_k^2}{2d^2} \cdot \ip{\nabla^2\varphi(X_k), \Delta_k^{\otimes 2}}
    \end{aligned}
\end{equation}
To derive the Doob decomposition, the idea is to iteratively condition on $r_{k, I_{k+1}}$, $\hat{X}_k$ and $I_{k+1}$. For this, we define the $\sigma$-algebras, 
\[
\begin{gathered}
\mathcal{G}_{k,i} \defas \sigma ( \{\hat{X}_j\}_{j=0}^k, \{r_i\}_{i=0}^k, \{I_j\}_{j=1}^k, \epsilon_{k+1}, I_{k+1} = i)\quad
 \text{and}\quad \mathcal{F}_{k} \defas \sigma ( \{\hat{X}_j\}_{j=0}^k ).
\end{gathered}
\]
It is clear that $s_{k, I_{k+1}}$ is measurable with respect to  $\mathcal{G}_{k,i}$ such that $s_{k,I_{k+1}}\mid \mathcal{G}_{k,i}$ corresponds to a deterministic maximization problem. This leads to the decomposition for SGD 
\begin{equation} \label{eq:main:taylor_expansion_result}
    \begin{aligned}
        \varphi(X_{k+1}) 
        & 
        =
        \varphi(X_k) - \frac{\gamma_k}{d} \ip{\nabla \varphi(X_k), \nabla \mathcal{R}^{{\operatorname{adv}}}_\lambda(X_k)}
        \\
        &
        + 
        \frac{\gamma_k^2}{2d^2} \sum_{i=1}^2 p_i \ip{\nabla^2 \varphi(X_k), K_i + \mu_i \mu_i^\top }\EE[(f_{i}'(g_{k, i} , \epsilon_{k+1}))^2\, | \, \mathcal{F}_k ]
        \\ 
        & + \Delta \mathcal{M}_k^{\text{Grad}}(\varphi) + \Delta \mathcal{M}_k^{\text{Hess}}(\varphi) + \EE[ \mathcal{E}_k^{\text{Hess}}(\varphi) \, |\, \mathcal{F}_k] 
    \end{aligned}
\end{equation}
where $\mathcal{M}_k^{\grad}(\varphi) $ and $\mathcal{M}_k^{\hess}(\varphi)$ are martingale increments and $ \mathcal{E}_k^{\text{Hess}}(\varphi)$ is an error term:
\begin{equation}
    \begin{aligned}
    \Delta \mathcal{M}_k^{\grad}(\varphi) 
    &
    = \frac{\gamma_k}{d} f_{I_{k+1}}'(g_{k, I_{k+1}} )\ip{\nabla \varphi (X_k),  a_{k+1, I_{k+1}} + \delta s_{k, {I_{k+1}}} \tfrac{X_k}{\|X_k\|} }  
    \\&\qquad- 
    \frac{\gamma_k}{d}\sum_{i=1}^2 p_i \EE \big [ f_{i}'(g_{k, i}, \epsilon )\ip{ \nabla \varphi (X_k), a_{k+1, i} + \delta s_{k, i} \tfrac{X_k}{\|X_k\|}  } \, | \, \mathcal{F}_k \big ],
    \\
   \Delta \mathcal{M}_k^{\text{Hess}}(\varphi)
    &
    = 
    \frac{\gamma_k^2}{2d^2} \bigg ( \ip{\nabla^2 \varphi(W_k), \Delta_k^{\otimes 2} } - \EE[ \ip{\nabla^2 \varphi(X_k), \Delta_k^{\otimes 2}} \, | \, \mathcal{F}_k] \bigg ) ,
    \\
    \text{and} \quad \mathcal{E}_k^{\text{Hess}}(\varphi)
    &
    = 
    \frac{\gamma_k^2}{2d^2}\sum_{i=1}^2 p_i \langle \nabla^2 \varphi(W_k), \Delta^{\otimes 2} - (K_i +\mu_i\mu_i^\top)f_i'(g_{k, i}, \epsilon_{k+1})^2\rangle.
\end{aligned}
\end{equation}
From Assumption~\ref{assumption:lr}, recall we defined $\gamma_t = \gamma(t/d)$ from which it follows $\gamma_{\lfloor td \rfloor} = \gamma(\lfloor td \rfloor/d)$. Similarly, we define $\gamma_{td}   = \gamma_{\lfloor td \rfloor}$ for simplicity of notation. With the following forward difference
\[
(\Delta \varphi)(X_{j}) \defas \varphi(X_{j+1}) - \varphi(X_{j}),
\]
we can use our continuous time equivalents to sum up (integrate) and obtain
\begin{equation}
\varphi(X_{td}) = \varphi(X_0) + \sum_{j=0}^{\lfloor td \rfloor - 1} (\Delta \varphi)(W_j) \defas \varphi(X_0) + \int_0^{t} d \cdot (\Delta \varphi)(X_{sd}) \, \dif s + \xi_{td}(\varphi),
\end{equation}
where $ \displaystyle |\xi_{td}(\varphi)|= \bigg | \int_{(\lfloor td\rfloor-1)/d}^t d \cdot \Delta \varphi(X_{sd})  \,\dif s \bigg | \le 2\max_{0 \le j \le \lfloor td \rfloor} \{ | \Delta \varphi(X_j) | \}$. With this, we obtain our Doob decomposition for SGD
\begin{align}
\varphi(X_{td}) 
&
= \varphi(X_0) - \int_0^{t} \gamma_{sd} \ip{ \nabla \varphi(X_{sd}), \nabla \mathcal{R}_{\lambda}^{{\operatorname{adv}}}(X_{sd})} \, \dif s \label{eq:main:integral_1} \\
&
        + 
        \frac{1}{2d} \int_0^{t} \gamma_{sd}^2 \sum_{i=1}^2 p_i \EE[ f_i'(g_{sd,i}, \epsilon)^{2} \, | \, \mathcal{F}_{sd} ] \ip{\nabla^2 \varphi(X_{sd}), K_i + \mu_i\mu_i^\top} \, \dif s \label{eq:main:integral_2}
        \\
        & +  \sum_{j=0}^{\lfloor td \rfloor-1} \left(\Delta \mathcal{M}_{j}^{\text{Grad}}(\varphi) + \Delta \mathcal{M}_{j}^{\text{Hess}}(\varphi) + \EE [\mathcal{E}_{j}^{\text{Hess}}(\varphi) \, | \, \mathcal{F}_{j} ]\right)
        + \xi_{td}(\varphi). \label{eq:main:error_terms_integrated}
\end{align}
We refer to Section~\ref{sect:doob:decomp:SGD:derivation} for more details on the derivation \eqref{eq:main:taylor_expansion_result}. Recall from \eqref{eq:main:AdvHSGD:def} the definition of AdvHSGD 
\begin{equation}
\dif \WHSGD_t = -\gamma(t) \nabla \mathcal{R}_{\lambda}^{\operatorname{adv}}(\WHSGD_t) \dif t + \gamma(t) \sqrt{\frac{1}{d} \sum_{i=1}^2 p_i \EE_{v, \epsilon} [  f'_i( \rho_{t, i}, \epsilon)^2  ]  (K_i + \mu_i\mu_i^\top)} \dif B_t,
\end{equation}
where the initial condition is given by $\mathscr{X}_0 = X_0$, $(B_t, t \ge 0)$ is a $\R^{d}$ standard Brownian motion and for $I \in \{1,2\}$ and $v\sim \mathcal{N}(0, \Id_d)$ we have
\begin{equation}
    \begin{gathered}
\widehat{\mathscr{X}}_t \defas \mathscr{X}_t \oplus X^{\star}\in \R^{d \times 2},\quad\mathscr{W}_t \defas \widehat{\mathscr{X}}_t \oplus \mu\in \R^{d \times 4},
\\
\rho_{t,I} \defas \widehat{\mathscr{X}}_t^\top(\sqrt{K_{I}}v + \mu_{I})+ \begin{bmatrix}\delta u_{t,I} \|\WHSGD_t\| & 0 \end{bmatrix}^\top,
\\
\text{and}\quad u_{t,I} = \operatorname{argmax}_{|u| \leq  1}f(\widehat{\mathscr{X}}_t^\top (\sqrt{K_{I}}v + \mu_{I}) + \begin{bmatrix} \delta u\|\mathscr{X}_t \|& 0 \end{bmatrix}^\top, \epsilon).
    \end{gathered}
\end{equation}
By It\^{o}'s lemma \cite[Thm. 33, Chapt. 2]{protter2005stochastic} applied to a quadratic statistic, we obtain
\begin{equation} \label{eq:main:HSGD_statistic}
    \begin{aligned}
        \dif \varphi(\WHSGD_t) 
        &= 
        \ip{\nabla \varphi(\WHSGD_t), \dif \WHSGD_t} + \frac{1}{2}\ip{\nabla^2 \varphi(\WHSGD_t), (\dif \WHSGD_t)^{\otimes 2} } 
        \\
        &
        =
        -\gamma(t) \ip{ \nabla \varphi(\WHSGD_t), \nabla \mathcal{R}^{{\operatorname{adv}}}_{\lambda}(\WHSGD_t) } \, \dif t
        \\
        &
        \qquad + \frac{\gamma^2(t)}{2d}\sum_{i=1}^2 p_i \EE_{v, \epsilon} [f_i'(\rho_{t,i}, \epsilon)^{2} ]\Tr (\nabla^2 \varphi(\WHSGD_t) (K_i + \mu_i\mu_i^\top ) )\dif t
        \\
        &
        \qquad 
        + \dif \mathcal{M}_{t}^{\operatorname{AdvHSGD}}(\varphi),
         \\
         \text{where} &\quad\dif \mathcal{M}_{t}^{\operatorname{AdvHSGD}}(\varphi) =  \frac{\gamma(t)}{\sqrt{d}}   \left\langle \sqrt{ \sum_{i=1}^2 p_i \EE_{v, \epsilon} [  f'_i( \rho_{t, i}, \epsilon)^2  ](K_i + \mu_i\mu_i^\top)}, \nabla\varphi(\WHSGD_t) \otimes \dif B_t \right\rangle
    \end{aligned}
\end{equation}
corresponds to the martingale increment for AdvHSGD. By integrating, we obtain the Doob decomposition for $\varphi(\WHSGD_t)$
\begin{align}
\varphi(\WHSGD_t) 
&
= \varphi(W_0) - \int_0^t \gamma(s) \ip{\nabla \varphi(\WHSGD_s), \nabla \mathcal{R}^{{\operatorname{adv}}}_{\lambda}(\WHSGD_s)} \, \dif s + \int_0^t \dif \mathcal{M}_s^{\operatorname{AdvHSGD}}(\varphi)\label{eq:main:doob:decomp:AdvHSGD}
\\
&
\qquad
+
\frac{1}{2 d} \int_0^t \gamma^2(s) \sum_{i=1}^2p_i\EE_{v, \epsilon} [f'_i(\rho_{s,i}, \epsilon)^2]\Tr (\nabla^2 \varphi(\WHSGD_s) (K_i + \mu_i\mu_i^\top ) )   \, \dif s. \label{eq:main:doob:decomp:AdvHSGD:vanish}
\end{align}
Equipped with these Doob decompositions, we must now prove that the term \eqref{eq:main:error_terms_integrated} and the martingale increment \eqref{eq:main:doob:decomp:AdvHSGD:vanish} become negligible as $d \to \infty$. The other terms \eqref{eq:main:integral_1}, \eqref{eq:main:integral_2} and \eqref{eq:main:doob:decomp:AdvHSGD} survive the limit while also being well-approximated by our deterministic system of coupled ODEs. 

In Section~\ref{sec:error_bounds}, we show that the terms \eqref{eq:main:error_terms_integrated} and \eqref{eq:main:doob:decomp:AdvHSGD:vanish} vanish as $d\to\infty$. Having this in mind, we finally show in Section~\ref{sect:SGD:AdvHSGD:approx:sols} that SGD and AdvHSGD on $S$ are $(\varepsilon, M,T)$-approximate solutions.

In the following section, we introduce the theory of approximate solutions to the integro-differential equation and an equivalent formulation of our deterministic system of ODEs in terms of statistics of the resolvents of $K_1$ and $K_2$.
\subsection{Theory of Approximate Solutions to the Integro-Differential Equation}\label{sect:approx:sols:theory}
In this section, we build upon the theory of approximate solutions to integro-differential equations presented in Section 4 of \citep{collinswoodfin2023hitting} and Section 5 in \citep{collinswoodfin2025Exact}. Recall by Assumptions~\ref{assumption:risk} and Assumption~\ref{assumption:fisher}, we have for $i=1,2$
\[
\begin{gathered}
     \mathcal{R}^{{\operatorname{adv}}}(X) \defas p_1h_1(\widehat{B}_1(W)) + p_2h_2(\widehat{B}_2(W)),\\
     \text{and}\quad 
 \mathcal{I}_i(\widehat{B}_i(W)) \defas \mathbb{E}_{v,\epsilon}[(f_i'( X^\top ((\sqrt{K_i}v+ \mu_i)) + \bar{s}_i \delta  \| X\|;  \epsilon))^2],
\end{gathered}
\]
where $h_i:\R^{4\times 4} \to \R$ and $I_i:\R^{4\times 4} \to \R$ are $\alpha$-pseudo-Lipschitz functions with constants $L(h), L(\mathcal{I})$ respectively and $h_i$ are continuously differentiable provided $(\widehat{B}_1(W), \widehat{B}_2(W) )\in \mathcal{U} \times \mathcal{U}$. We also recall the notation we previously introduced
\begin{equation}
    \widehat{B}_i(W) \defas \begin{bmatrix}
           U_i & 0\\
           0 & \|X\|^2
           \end{bmatrix},\quad  U_i(W) \defas \begin{bmatrix}
               B_i(W)  & m_i(W) \\
               m_i(W)^\top & \|\mu_i\|^2
           \end{bmatrix},\quad B_i(W)\defas \hat{X}^\top K_i \hat{X}
   \end{equation}
and $m_i(W) \defas \hat{X}^\top \mu_i$. Here $W$ is a placeholder value for AdvHSGD or SGD. 
In this section, we are interested in approximate solutions to the integro-differential equation \eqref{eq:ODE_resolvent_2}.

Fundamentally, the adversarial framework requires rederiving the results presented in \citep{collinswoodfin2023hitting, collinswoodfin2025Exact} to account for the new structure of the deterministic equivalents involving $\partial_{44} h$ and $\mathrsfs{B}_{44}(t)$. Hence, the proofs will build upon the two papers while adding the necessary technicalities. The general idea behind the framework presented in \citep{collinswoodfin2023hitting,collinswoodfin2025Exact} is that both $Z(\HSGD_t, z)$ and $Z(W_{\lfloor td \rfloor}, z)$ which are respectively functions of AdvHSGD and SGD are approximate solutions to \eqref{eq:ODE_resolvent_2}. 

Recall we previously defined the contour $\Gamma = \{ w \in \mathbb{C} \, : \, |w| = \max\{1, 2 \|\mathcal{K}\|_{\operatorname{op}} \} \}$ with $\|\mathcal{K}\|_{\operatorname{op}} \defas \max\{ \|K_1\|_{\opt}, \|K_2\|_{\opt}\}$ such that the distance of any element of $\Gamma$ to the spectrum of $K_1$ and $K_2$ is at least $\frac{1}{2}$. Recall for $z \in \Gamma^2$, we write $z =(z_1 ,z_2)$ and when integrating over $z_1$ and $z_2$ simultaneously, we write for any function $f: \mathbb{C}^2 \to \mathbb{C}$
\[
\oint f(z) \text{D}z \defas \frac{-1}{4\pi^2 }\oint_{\Gamma^2} f(z)  \dif z_1 \dif z_2.
\]
Recall we introduced a norm over the contour $\Gamma^2$ for a continuous function $H:\mathbb{C}^2 \to \mathbb{C}^{4\times 4}$ in Definition~\ref{eq:main:gamma:normdef}:
    \[
        \|H\|_{\Gamma} \defas \sup_{z\in \Gamma^2}\|H(z)\|.
    \]
\begin{remark}
    The results and proofs in this section will be presented for the soft label setting. Note that proofs may be adapted to the hard label setting by considering $X^{\star} = 0$ and very similar arguments.
\end{remark}

Lemma 3 in \citep{collinswoodfin2025Exact} shows that up to constants depending on $\|\mathcal{K}\|_{\operatorname{op}}$, the norm $\|\cdot\|_{\Gamma}$ applied to the statistics $Z(W_{\lfloor td\rfloor}, \cdot)$ or $\mathcal{Z}(t, \cdot)$ is proportional to the norm of the iterates. We adapt this result slightly as our contour $\Gamma$ differs from their contour $\Omega$. We omit this proof however as it is almost identical. This result is analogous to Lemma 5 in \citep{collinswoodfin2023hitting}.
\begin{lemma}[Lemma 3, \citep{collinswoodfin2025Exact}]\label{lem:prop:normGamma:normW}
    Recall $\mathrsfs{N}(t) \defas \oint\Tr(\mathcal{Z}(t,z))\Dif z$. Then there exists a constant $C>0$ depending on $\|\mathcal{K}\|_{\operatorname{op}}$ and such that
    \[
        C \leq \frac{\|Z(\HSGD_t, \cdot)\|_{\Gamma}}{\|\HSGD_{t}\|^2},\frac{\|Z(W_{td}, \cdot)\|_{\Gamma}}{\|W_{td}\|^2}, \frac{\|\mathcal{Z}(t, \cdot)\|_{\Gamma}}{\mathrsfs{N}(t)} \leq 4.
    \]
    
\end{lemma}
 Using the norm $\|\cdot\|_{\Gamma}$, we introduce a stopping time $\tau_M$ analogous to $\vartheta_M$ in Theorem~\ref{thm:main:conc:St:SGD}.
\begin{definition}[Stopping time $\tau_M$]\label{def:stoppingtime}
    For any $M>0$ and a continuous function $\mathrsfs{Z}: [0, \infty] \times \mathbb{C}^2 \to \R^{4\times 4}$, define
    \[
        \tau_M(\mathrsfs{Z}) \defas \inf\bigg\{t \geq 0 : \|\mathrsfs{Z}(t,\cdot)\|_{\Gamma}> M \quad \text{or}\quad \widehat{\mathcal{B}}_1(t) \not \in \mathcal{U} \quad \text{or}\quad \widehat{\mathcal{B}}_2(t) \not \in \mathcal{U}\bigg\},
    \]
    where we define
    \begin{equation}\label{eq:def:Z(t,z):stability:single}
        \mathrsfs{Z}(t, \cdot) \defas \begin{bmatrix}
            \mathrsfs{S}(t,\cdot) & \mathrsfs{M}(t,\cdot)\\
            \mathrsfs{M}(t,\cdot)^\top & \mu^\top \mathscr{R}(z;\mathcal{K})\mu
        \end{bmatrix}.
    \end{equation}
    with
\[
    \widehat{\mathcal{B}}_i(t) \defas \begin{bmatrix}
        \mathcal{B}_{i}(t)&  \mathscr{M}_{i}(t)&0 \\
                (\mathscr{M}_{i}(t))^\top&  \|\mu_i\|^2 &0 \\
                0 & 0 &  \mathcal{V}_{11}(t)
     \end{bmatrix},
\]
and 
\[
    \mathcal{B}_{i}(t) \defas
    \oint z_i \mathrsfs{S}(t,z) \, \Dif z,
  \quad 
  \mathscr{M}_{i}(t)\defas\oint \mathrsfs{M}_{i}(t,z) \, \Dif z,\quad
\mathcal{V}(t) \defas \oint \mathrsfs{S}(t,z) \, \Dif z,
\]
where $\mathrsfs{M}_{i} $ denotes the $i$th column of $\mathrsfs{M}$ associated with $\mu_i$ and $\mathcal{V}_{11}$ denotes the element in the position $(1,1)$ of $\mathcal{V}$
\end{definition}
When clear from context, we omit the dependence $  \tau_M=\tau_M(\mathrsfs{Z}) $. The stopping time $\tau_M$ differs from $\vartheta_M$ introduced in Theorem~\ref{thm:main:conc:St:SGD} as they respectively bound $\|\mathcal{Z}(t,\cdot)\|_\Gamma$ and $\mathrsfs{N}(t)$. Note that Lemma~\ref{lem:prop:normGamma:normW} implies there exists constants $c,C>0$ such that $\vartheta_{c \cdot M} \leq \tau_M \leq \vartheta_{C\cdot M}$ so the stopping times are interchangeable. We are now ready to introduce the notion of approximate solutions to integro-differential equations. 
\begin{definition}[$(\varepsilon, M, T )$-approximate solution to the integro-differential equation]\label{def:main:approx:sols}
    For $M, T, \varepsilon > 0$ and a continuous function $\mathrsfs{Z}: [0, \infty] \times \mathbb{C}^2 \to \R^{4\times 4}$, we denote $\mathrsfs{Z}$ an $(\varepsilon, M, T)$-approximate solution of \eqref{eq:ODE_resolvent_2} if
    \[
    \sup_{0\leq t\leq T\wedge \tau_M} \left\| \mathrsfs{Z}(t, \cdot) -\mathrsfs{Z}(0, \cdot) -\int_0^t \mathscr{F}(\cdot, \mathrsfs{Z}(s, \cdot)) \dif s\right\|_{\Gamma} \leq \varepsilon,
    \]
    where $\mathrsfs{Z}(0, \cdot) = W_0^\top \mathscr{R}(\cdot; \mathcal{K})W_0$ and $W_0 = X_0 \oplus X^{\star} \oplus \mu$.
\end{definition}
It is clear from the definition~\ref{def:main:approx:sols} and the definition of $\mathcal{Z}$ in \eqref{eq:ODE_resolvent_2} that $\mathcal{Z}$ is an approximate solution with $\varepsilon = 0$. In the next subsection, we will show $Z(W_{td})$ and $Z(\HSGD_{t})$ are $(\varepsilon, M, T)$-approximate solutions to \eqref{eq:ODE_resolvent_2}. Before proving these results, we study the stability of the approximate solutions of \eqref{eq:ODE_resolvent_2} analogous to Proposition 10 in \citep{collinswoodfin2023hitting} and Corollary 1 in \citep{collinswoodfin2025Exact}.
\begin{proposition}[Stability]\label{prop:stability}
    For all $(\varepsilon, M, T)$-approximate solutions $\mathrsfs{Z}_1$ and $\mathrsfs{Z}_2$, there exists a constant $C= C(\bar{\gamma},\lambda, M, \alpha, L(h), L(\mathcal{I}), \|\mathcal{K}\|_{\operatorname{op}},\max_{i=1,2} \|\mu_i\|, T)>0$ such that
    \[
        \sup_{0\leq t \leq T\wedge \widehat{\tau}_M} \|\mathrsfs{Z}_1(t,\cdot) - \mathrsfs{Z}_2(t,\cdot)\|_\Gamma \leq C\varepsilon,
    \]
    where $\widehat{\tau}_M =\min\{\tau_M(\mathrsfs{Z}_1), \tau_M(\mathrsfs{Z}_2)\}$.
\end{proposition}
\begin{proof}
    Since $\mathrsfs{Z}_1$ and $\mathrsfs{Z}_2$ are $(\varepsilon, M, T)$-approximate solutions then for $i=1,2$ we can write
    \begin{equation}
        \mathrsfs{Z}_i(t,\cdot) = \mathrsfs{Z}_i(0,\cdot) + \int_0^t \mathscr{F}(\cdot, \mathrsfs{Z}_i(s,\cdot)) \dif s + \varepsilon_t( \mathrsfs{Z}_i),
    \end{equation}
    where $\varepsilon_t( \mathrsfs{Z}_i)$ are errors terms satisfying $\sup_{0\leq t\leq T \wedge \widehat{\tau}_M}\| \varepsilon_t( \mathrsfs{Z}_i)\|_{\Gamma} \leq \varepsilon$. By construction, we have $\widehat{\tau}_M  \leq \tau_M(\mathrsfs{Z}_1)$, $\tau_M(\mathrsfs{Z}_2)$ so we will work with $\widehat{\tau}_M$ throughout the proof. Since $\mathrsfs{Z}_1(0,\cdot) = \mathrsfs{Z}_2(0,\cdot)$, proving the result amounts to showing the function $\mathscr{F}$ is Lipschitz with respect to its second input:
    \[
        \| \mathscr{F}(\cdot, \mathrsfs{Z}_1(t\wedge \widehat{\tau}_M,\cdot)) - \mathscr{F}(\cdot, \mathrsfs{Z}_2(t\wedge \widehat{\tau}_M,\cdot)) \|_{\Gamma} \leq L(\mathcal{F})\|\mathrsfs{Z}_1(t\wedge \widehat{\tau}_M , \cdot) - \mathrsfs{Z}_2(t\wedge \widehat{\tau}_M, \cdot)\|_{\Gamma}.
    \]
    From the definition of of $\mathscr{F}$, we have that
    \begin{equation}
        \begin{aligned}
            \| \mathscr{F}(\cdot, \mathrsfs{Z}_1(t\wedge \widehat{\tau}_M,\cdot)) &- \mathscr{F}(\cdot, \mathrsfs{Z}_2(t\wedge \widehat{\tau}_M,\cdot)) \|_{\Gamma}
            \\
            &\leq \| \mathscr{F}_{\mathcal{S}}(\cdot, \mathrsfs{Z}_1(t\wedge \widehat{\tau}_M,\cdot)) - \mathscr{F}_{\mathcal{S}}(\cdot, \mathrsfs{Z}_2(t\wedge \widehat{\tau}_M,\cdot)) \|_{\Gamma} 
            \\&\qquad
            + 2 \| \mathscr{F}_{\mathcal{M}}(\cdot, \mathrsfs{Z}_1(t\wedge \widehat{\tau}_M,\cdot)) - \mathscr{F}_{\mathcal{M}}(\cdot, \mathrsfs{Z}_2(t\wedge \widehat{\tau}_M,\cdot)) \|_{\Gamma}.
        \end{aligned}
    \end{equation}
    Hence, showing $\mathscr{F}$ is Lipschitz amounts to showing $\mathscr{F}_{\mathcal{S}}$ and $\mathscr{F}_{\mathcal{M}}$ are Lipschitz. To show this, we build upon the proofs of Proposition 10 in \citep{collinswoodfin2023hitting} and account for the contributions of the adversarial terms. Similarly to the notation introduced in Definition~\ref{def:stoppingtime}, we will represent the structure of $\mathrsfs{Z}_j$ for $j=1,2$ as follows
    \begin{equation}\label{eq:def:Z(t,z):stability}
        \mathrsfs{Z}_{j}(t, \cdot) \defas \begin{bmatrix}
            \mathrsfs{S}_{j}(t,\cdot) & \mathrsfs{M}_{j}(t,\cdot)\\
            \mathrsfs{M}_{j}(t,\cdot)^\top & \mu^\top \mathscr{R}(z;\mathcal{K})\mu
        \end{bmatrix}.
    \end{equation}
    For simplicity of notation, we define for $j=1,2$
    \begin{equation}\label{eq:notation:stability:proof}
        \begin{gathered}
            \mathcal{B}_{i,j}(t) \defas
                 \oint z_i \mathrsfs{S}_j(t,z) \, \Dif z,
               \quad 
               \mathscr{M}_{i,j}(t)\defas\oint \mathrsfs{M}_{i,j}(t,z) \, \Dif z,\quad
            \mathcal{V}_{j}(t) \defas \oint \mathrsfs{S}_{j}(t,z) \, \Dif z,\\
            \text{and}\quad\widehat{\mathcal{B}}_{i,j}(t) \defas \left[ \begin{array}{ccc} 
                \mathcal{B}_{i,j}(t)&  \mathscr{M}_{i,j}(t)&0 \\
                \mathscr{M}_{i,j}(t)& \|\mu_i\|^2   &0 \\
                0 & 0 &  \mathcal{V}_{j, 11}(t)
            \end{array} 
            \right] .
        \end{gathered}
    \end{equation}
    where $\mathrsfs{M}_{i,j} $ denotes the $i$th column of $\mathrsfs{M}_{j}$ associated with $\mu_i$ and $\mathcal{V}_{j, 11}$ denotes the element in the position $(1,1)$ of $\mathcal{V}_j$. We will work with $s\leq \widehat{\tau}_M$ and omit the stopping time in the next steps where we show $\mathscr{F}$ is Lipschitz. First, since $\oint |q(z_i)| \Dif |z|$ is uniformly bounded by a constant depending on $\|\mathcal{K}\|_{\opt}$, it follows from the definition of the Frobenius norm
    \begin{equation}\label{eq:Fstuff:bounds}
        \begin{gathered}
            \|\mathrsfs{S}_{j}(t,\cdot)\|_{\Gamma} \leq M, \quad \|z_i\mathrsfs{S}_{j}(t,\cdot)\|_{\Gamma}\leq C(M, \|\mathcal{K}\|_{\opt}),
            \\
            \|\mathrsfs{M}_{i,j}(t,\cdot)\|_{\Gamma}\leq \|\mathrsfs{M}_{j}(t,\cdot)\|_{\Gamma} \leq M, \quad \|z_i\mathrsfs{M}_{j}(t,\cdot)\|_{\Gamma}\leq C(M, \|\mathcal{K}\|_{\opt}),
            \\ \left\|\frac{-1}{2\pi \operatorname{i}}\oint_{\Gamma} \mathrsfs{M}_{j}(t,z) \dif z_i\right\|_{\Gamma} \leq C(M, \|\mathcal{K}\|_{\opt})\quad \text{and}\quad \left\|\frac{-1}{2\pi \operatorname{i}}\oint_{\Gamma} \mathrsfs{S}_{j}(t,z) \dif z_i\right\|_{\Gamma} \leq C(M, \|\mathcal{K}\|_{\opt}).
        \end{gathered}
    \end{equation}
    Similarly, we have the following inequalities which also follow from linearity and the definition of the Frobenius norm
    \begin{equation}\label{eq:Fstuff:Lipschitz}
        \begin{gathered}
            \|\mathrsfs{S}_{1}(t,\cdot) - \mathrsfs{S}_{2}(t,\cdot)\|_{\Gamma} \leq \|\mathrsfs{Z}_{1}(t,\cdot) - \mathrsfs{Z}_{2}(t,\cdot)\|_{\Gamma},
            \\
            \|\mathrsfs{M}_{i,1}(t,\cdot) - \mathrsfs{M}_{i, 2}(t,\cdot)\|_{\Gamma}\leq \|\mathrsfs{M}_{1}(t,\cdot) - \mathrsfs{M}_{2}(t,\cdot)\|_{\Gamma} \leq \|\mathrsfs{Z}_{1}(t,\cdot) - \mathrsfs{Z}_{2}(t,\cdot)\|_{\Gamma},
            \\
            \left\|\frac{-1}{2\pi \operatorname{i}}\oint_{\Gamma}\mathrsfs{M}_{1}(t,z) -\mathrsfs{M}_{2}(t,z) \dif z_i\right\|_{\Gamma}\leq C(\|\mathcal{K}\|_{\opt})  \|\mathrsfs{Z}_{1}(t,\cdot) - \mathrsfs{Z}_{2}(t,\cdot)\|_{\Gamma},
            \\
            \left\|\frac{-1}{2\pi \operatorname{i}}\oint_{\Gamma}\mathrsfs{S}_{1}(t,z) -\mathrsfs{S}_{2}(t,z) \dif z_i\right\|_{\Gamma}\leq C(\|\mathcal{K}\|_{\opt})  \|\mathrsfs{Z}_{1}(t,\cdot) - \mathrsfs{Z}_{2}(t,\cdot)\|_{\Gamma}.
        \end{gathered}
    \end{equation}
    Hence, we have
    \begin{equation}
        \| \mathcal{B}_{i,j}(s)\| \leq C(\|\mathcal{K}\|_{\opt})\|\mathrsfs{S}_j(\cdot,z)\|_{\Gamma} \leq C(\|\mathcal{K}\|_{\opt})\|\mathrsfs{Z}_j(\cdot,z)\|_{\Gamma} \leq C(\|\mathcal{K}\|_{\opt}, M).
    \end{equation}
    Similarly, we obtain for $C = C(\|\mathcal{K}\|_{\opt}, M)$
    \begin{equation}\label{eq:bounds:BBV:Lipschitzproof}
        \begin{gathered}
        \|\mathscr{M}_{i,j}(\cdot) \| \leq C,
        \quad \|\mu_i\| \leq C,
        \quad
            \text{and}\quad
            \|\mathcal{V}_{j,11}(s)\|\leq \|\mathcal{V}_{j}(s)\|\leq
            C
            \\
            \text{from which it follows that}\quad 
            \|\widehat{\mathcal{B}}_{i,j}(s)\| \leq C.
        \end{gathered}
    \end{equation}
    From the definitions of $H_{1,i}$ and $H_{2,i}$ in \eqref{eq:ODE_resolvent_2} and Assumption~\ref{assumption:risk}, $\nabla h_i$ is $\alpha$-pseudo-Lipschitz so using \eqref{eq:bounds:BBV:Lipschitzproof} it follows for $\ell=1,2$
    \begin{equation}\label{eq:H:Lipschitz:working}
        \begin{aligned}
            \|H_{\ell, i}( \widehat{\mathcal{B}}_{i,1}(s)) - H_{\ell, i}( \widehat{\mathcal{B}}_{i,2}(s))\|
            &
            \leq
             L(h)(1+ \|\widehat{\mathcal{B}}_{i,1}(s)\|^{\alpha} + \|\widehat{\mathcal{B}}_{i,2}(s)\|^{\alpha})\|\widehat{\mathcal{B}}_{i,1}(s)-\widehat{\mathcal{B}}_{i,2}(s)\|
             \\
            &
            \leq
            C(M, \alpha, L(h), \|\mathcal{K}\|_{\operatorname{op}}) \|\widehat{\mathcal{B}}_{i,1}(s)-\widehat{\mathcal{B}}_{i,2}(s)\|.
        \end{aligned}
    \end{equation}
    Once again since $\oint |q(z_i)| \Dif |z|$ is uniformly bounded and from \eqref{eq:Fstuff:Lipschitz}, then we have
    \begin{equation}
    \begin{gathered}\label{eq:bound:V:Lipschitz}
            \|\mathcal{B}_{i,1}(s)-\mathcal{B}_{i,2}(s)\| \leq   C(\|\mathcal{K}\|_{\operatorname{op}}) \|\mathrsfs{Z}_1(s,\cdot) - \mathrsfs{Z}_2(s,\cdot)\|_{\Gamma},
            \\
            \|\mathscr{M}_{i,1}(s)-\mathscr{M}_{i,2}(s)\| \leq   C(\|\mathcal{K}\|_{\operatorname{op}}) \|\mathrsfs{Z}_1(s,\cdot) - \mathrsfs{Z}_2(s,\cdot)\|_{\Gamma},
            \\
            \|\mathcal{V}_{1, 11}(s)-\mathcal{V}_{2, 11}(s)\|  \leq\|\mathcal{V}_{1}(s)-\mathcal{V}_{2}(s)\| \leq   C(\|\mathcal{K}\|_{\operatorname{op}}) \|\mathrsfs{Z}_1(s,\cdot) - \mathrsfs{Z}_2(s,\cdot)\|_{\Gamma}.
    \end{gathered}   
\end{equation}
Hence, we obtain
\begin{equation}\label{eq:bound:B(s):lipschitz}
    \|\widehat{\mathcal{B}}_{i,1}(s) - \widehat{\mathcal{B}}_{i,2}(s)\|  \leq C(\|\mathcal{K}\|_{\operatorname{op}}) \|\mathrsfs{Z}_1(s,\cdot) - \mathrsfs{Z}_2(s,\cdot)\|_{\Gamma}.
\end{equation}
Hence, inserting \eqref{eq:bound:B(s):lipschitz} into \eqref{eq:H:Lipschitz:working}, we obtain
\begin{equation}\label{eq:H:Lipschitz:final}
    \|H_{\ell,i}( \widehat{\mathcal{B}}_{i,1}(s)) - H_{\ell, i}( \widehat{\mathcal{B}}_{i,2}(s))\| \leq C(M, \alpha, L(h), \|\mathcal{K}\|_{\operatorname{op}}) \|\mathrsfs{Z}_1(s,\cdot) - \mathrsfs{Z}_2(s,\cdot)\|_{\Gamma}.
\end{equation}
From Assumption~\ref{assumption:risk} and \eqref{eq:bounds:BBV:Lipschitzproof}, we obtain the bound
\begin{equation}\label{eq:main:bound:stability:H}
    \|H_{\ell,i}( \widehat{\mathcal{B}}_{i,j}(s))\| \leq C L(h)(1+ \|\widehat{\mathcal{B}}_{i,j}(s)\|)^{1+\alpha} \leq C(\|\mathcal{K}\|_{\operatorname{op}}, M, \alpha, L(h)).
\end{equation} 
Similarly to \eqref{eq:H:Lipschitz:working} and \eqref{eq:main:bound:stability:H}, from Assumption~\ref{assumption:risk} and \eqref{eq:bounds:BBV:Lipschitzproof}, we have
\begin{equation}\label{eq:H:Lipschitz:partial44}
    \begin{aligned}
        |\partial_{44} h_i (\widehat{\mathcal{B}}_{i,1}(s)) - \partial_{44} h_i (\widehat{\mathcal{B}}_{i,2}(s))|
        &
        \leq
        C( M, \alpha, L(h), \|\mathcal{K}\|_{\operatorname{op}}) \|\mathrsfs{Z}_1(s,\cdot) - \mathrsfs{Z}_2(s,\cdot)\|_{\Gamma},\\
        |\frac{\lambda}{2} +\partial_{44} h_i (\widehat{\mathcal{B}}_{i,j}(s))| &\leq C(\lambda, \|\mathcal{K}\|_{\operatorname{op}}, M, \alpha, L(h)).
    \end{aligned}
\end{equation}
It is also clear from \eqref{eq:def:Z(t,z):stability} that $\|\mu^\top \mathscr{R}(z;\mathcal{K})\mu_i\| \leq \|\mu^\top \mathscr{R}(z;\mathcal{K})\mu\| \leq \|\mathrsfs{Z}_j\|_{\Gamma} \leq M$. Thus, it follows from \eqref{eq:Fstuff:bounds}, \eqref{eq:Fstuff:Lipschitz}, \eqref{eq:H:Lipschitz:final}, \eqref{eq:main:bound:stability:H} and \eqref{eq:H:Lipschitz:partial44} and summing over $i$ that
\begin{equation}\label{eq:main:stability:bound:M}
    \begin{aligned}
        \| \mathscr{F}_{\mathcal{M}}(\cdot,& \mathrsfs{Z}_1(t\wedge \widehat{\tau}_M,\cdot)) - \mathscr{F}_{\mathcal{M}}(\cdot, \mathrsfs{Z}_2(t\wedge \widehat{\tau}_M,\cdot)) \|_{\Gamma}
        \\
        &
        \leq C(\lambda, \bar{\gamma}, M, \alpha, L(h), \|\mathcal{K}\|_{\operatorname{op}}) \|\mathrsfs{Z}_1(s)-\mathrsfs{Z}_2(s)\|.
    \end{aligned}
\end{equation}
From Assumption~\ref{assumption:fisher}, \eqref{eq:bounds:BBV:Lipschitzproof} and \eqref{eq:bound:B(s):lipschitz}, we obtain
\begin{equation}
    |\mathcal{I}_i(\widehat{\mathcal{B}}_{i,1}(s)) - \mathcal{I}_i(\widehat{\mathcal{B}}_{i,2}(s))| \leq C(L(\mathcal{I}), \|\mathcal{K}\|_{\operatorname{op}}, M, \alpha)\|\mathrsfs{Z}_1(s,\cdot) - \mathrsfs{Z}_2(s,\cdot)\|_{\Gamma}.
\end{equation}
As shown in Lemma A.3 of \citep{collinswoodfin2023hitting}, from the Neumann series expansion of the resolvent $\sup_{z_i\in \Gamma} \|R(z_i;K_i)\|_{\operatorname{op}} \leq 2$ from which it follows that $\sup_{z\in \Gamma^2}\|\mathscr{R}(z; \mathcal{K})\|_{\opt} \leq 4$. Given the previous two bounds and using the identity $K_iR(z_i;K_i) = \Id_d + z_i R(z_i;K_i)$, it follows
\begin{equation}\label{eq:main:stability:bound:4}
    \begin{aligned}
    &\sup_{z\in\Gamma^2}\left|\frac{\gamma_s^2}{d} \Tr((K_1 + \mu_1\mu_1^\top) \mathscr{R}(z;\mathcal{K}))\left(\mathcal{I}_i(\widehat{\mathcal{B}}_{1,1}(s)) - \mathcal{I}_i(\widehat{\mathcal{B}}_{1,2}(s))\right)\right|
    \\
    &
    \leq
    \frac{C(\bar{\gamma})}{d} \sup_{z \in \Gamma^2}(| \Tr(R(z_2;K_2))| + |z_1||\Tr(\mathscr{R}(z;\mathcal{K}))| + |\mu_1^\top\mathscr{R}(z;\mathcal{K})\mu_1|)|\mathcal{I}_i(\widehat{\mathcal{B}}_{1,1}(s)) - \mathcal{I}_i(\widehat{\mathcal{B}}_{1,2}(s))|
    \\
    &
    \leq
    C(\bar{\gamma},L(\mathcal{I}), \|\mathcal{K}\|_{\operatorname{op}}, \max_{i=1,2} \|\mu_i\|, M, \alpha)\|\mathrsfs{Z}_1(s,\cdot) - \mathrsfs{Z}_2(s,\cdot)\|_{\Gamma}.
    \end{aligned}
\end{equation}
The same result holds with $K_2$, $\mu_2$, $\widehat{\mathcal{B}}_{2,1}(s)$ and $\widehat{\mathcal{B}}_{2,1}(s)$.  Thus, it follows from \eqref{eq:Fstuff:bounds}, \eqref{eq:Fstuff:Lipschitz}, \eqref{eq:H:Lipschitz:final}, \eqref{eq:main:bound:stability:H}, \eqref{eq:H:Lipschitz:partial44} and \eqref{eq:main:stability:bound:4} and summing over $i$ that
\begin{equation}\label{eq:main:stability:bound:S}
    \begin{aligned}
        \| \mathscr{F}_{\mathcal{S}}(\cdot,& \mathrsfs{Z}_1(t\wedge \widehat{\tau}_M,\cdot)) - \mathscr{F}_{\mathcal{S}}(\cdot, \mathrsfs{Z}_2(t\wedge \widehat{\tau}_M,\cdot)) \|_{\Gamma}
        \\
        &
        \leq C(\lambda, \bar{\gamma}, M, \alpha, L(h), L(\mathcal{I}), \|\mathcal{K}\|_{\operatorname{op}}) \|\mathrsfs{Z}_1(s)-\mathrsfs{Z}_2(s)\|.
    \end{aligned}
\end{equation}
    Hence, combining \eqref{eq:main:stability:bound:M} and \eqref{eq:main:stability:bound:S}, we've shown $\mathscr{F}$ is $L(\mathscr{F})$-Lipschitz with respect to its second input where $L(\mathscr{F}) =  C(\bar{\gamma},\lambda, M, \alpha, L(h), L(\mathcal{I}), \|\mathcal{K}\|_{\operatorname{op}}, \max_{i=1,2} \|\mu_i\|)$. It thus follows
    \begin{equation}
        \begin{aligned}
        \sup_{0\leq t\leq T \wedge \widehat{\tau}_M} \| \mathrsfs{Z}_1(t, \cdot) &- \mathrsfs{Z}_2(t, \cdot) \|_{\Gamma} 
        \leq 
        2 \varepsilon +  \sup_{0\leq t\leq T \wedge \widehat{\tau}_M} \int_0^t \|\mathscr{F}(\cdot, \mathrsfs{Z}_1(s,\cdot))  - \mathscr{F}(\cdot, \mathrsfs{Z}_2(s,\cdot)) \|_{\Gamma}\dif s
        \\
        &
        \leq 
        2 \varepsilon +  \sup_{0\leq t\leq T } \int_0^t \|\mathscr{F}(\cdot, \mathrsfs{Z}_1(s \wedge \widehat{\tau}_M,\cdot))  - \mathscr{F}(\cdot, \mathrsfs{Z}_2(s \wedge \widehat{\tau}_M,\cdot)) \|_{\Gamma}\dif s
        \\
        &
        \leq 
        2 \varepsilon +  L(\mathscr{F})\int_0^t \| \mathrsfs{Z}_1(s \wedge \widehat{\tau}_M,\cdot)  -  \mathrsfs{Z}_2(s \wedge \widehat{\tau}_M,\cdot) \|_{\Gamma}\dif s
        \\
        &
        \leq 
        2 \varepsilon +  L(\mathscr{F})\int_0^t\sup_{0\leq u\leq s\wedge \widehat{\tau}_M}\| \mathrsfs{Z}_1(u,\cdot)  -   \mathrsfs{Z}_2(u,\cdot) \|_{\Gamma}\dif s.
        \\
        \end{aligned}
    \end{equation}
Applying Grönwall's inequality, we obtain
\begin{equation}
    \sup_{0\leq t\leq T \wedge \widehat{\tau}_M} \| \mathrsfs{Z}_1(t, \cdot) - \mathrsfs{Z}_2(t, \cdot) \|_{\Gamma}  \leq 2\varepsilon e^{L(\mathscr{F})T},
\end{equation}
from which the result follows.
\end{proof}
Proposition~\ref{prop:stability} can be extended to any function $\varphi(X)$ satisfying Assumption~\ref{assumption:smooth:stats}. This result is analogous to Proposition 11 in \citep{collinswoodfin2023hitting} and Proposition 7 in \citep{collinswoodfin2025Exact}.
\begin{proposition}\label{prop:stability:general}
    Suppose $\varphi: \R^{d}\to \R$ satisfies Assumption~\ref{assumption:smooth:stats} such that 
    \[
        \varphi(X) = g\left(\begin{bmatrix}
            \hat{X}^\top q(K_1, K_2) \hat{X} & \hat{X}^\top \mu &0 \\
            \mu^\top \hat{X} & \mu^\top \mu & 0 \\
            0 & 0 &  \|X\|^2
         \end{bmatrix}\right),
    \]
    with $q$ a polynomial. Let $\mathrsfs{Z}_1$ and $\mathrsfs{Z}_2$ be $(\varepsilon, M, T)$-approximate solutions and for $i=1,2$ define
    \begin{equation}
        \mathrsfs{G}(\mathrsfs{Z}_i(t,\cdot))  \defas  \begin{bmatrix}
            \oint q(z)\mathrsfs{S}_i(t,z)\Dif z &  \oint \mathrsfs{M}_i(t,z)\Dif z & 0  \\
            \oint \mathrsfs{M}_i(t,z)^\top \Dif z &  \mu^\top  \mu  & 0 \\
            0 & 0 & \oint\mathrsfs{S}_{i,11}(t,z)\Dif z
        \end{bmatrix} \in \mathbb{C}^{5\times 5}.
    \end{equation}
    where $\mathrsfs{S}_{i,11}$ denotes the $(1,1)$ element of $\mathrsfs{S}_{i}$. Then there exists a positive constant $C$ such that
    \[
        \sup_{0\leq t \leq T\wedge \widehat{\tau}_M} \left| g(\mathrsfs{G}(\mathrsfs{Z}_1(t,\cdot)))  - g(\mathrsfs{G}(\mathrsfs{Z}_2(t,\cdot)))\right| \leq C\varepsilon,
    \]
    where $C= C(\bar{\gamma},T, q,\lambda, M, \alpha, L(h), L(g), \|K\|_{\operatorname{op}})$ and $\widehat{\tau}_M =\min\{\tau_M(\mathrsfs{Z}_1), \tau_M(\mathrsfs{Z}_2)\}$.
\end{proposition}
\begin{proof}
    Similarly to the proof of Proposition~\ref{prop:stability}, we have $\widehat{\tau}_M  \leq \tau_M(\mathrsfs{Z}_1)$, $\tau_M(\mathrsfs{Z}_2)$ so we can work with $\widehat{\tau}_M$ throughout the proof. For $i=1,2$ we define 
    \begin{equation}
        \mathrsfs{Q}_i(t) =\oint q(z) \mathrsfs{S}_i(t,z) \Dif z,
    \end{equation}
    where $q$ is an arbitrary polynomial. We will work with $s \leq \widehat{\tau}_M$ in order to omit the stopping time. For $i=1,2$, we have 
    \begin{equation}\label{eq:bound:Qi}
    \|\mathrsfs{Q}_i(s)\|  \leq C \oint |q(z)| \|\mathrsfs{S}_i(s,\cdot)\|_{\Gamma} \Dif |z| \leq C(\|\mathcal{K}\|_{\operatorname{op}}, \sup_{z\in \Gamma^2}|q(z)|, M) .
    \end{equation}
    Given that for any $w\in \Gamma$ we have $|w| = \max(1, 2\|\mathcal{K}\|_{\operatorname{op}})\leq 1 + 2 \|\mathcal{K}\|_{\operatorname{op}}$ then $\sup_{z\in \Gamma^2} |q(z)| \leq C(\|\mathcal{K}\|_{\operatorname{op}})$. It follows similarly that $\mathrsfs{Q}$ is Lipschitz
    \begin{equation}\label{eq:lipschitz:Qi}
        \begin{aligned}
        \|\mathrsfs{Q}_1(s) - \mathrsfs{Q}_2(s)\| 
        \leq
         C(\|\mathcal{K}\|_{\operatorname{op}}) \|\mathrsfs{S}_1(s,\cdot) - \mathrsfs{S}_2(s,\cdot)\|_{\Gamma} \leq C(\|\mathcal{K}\|_{\operatorname{op}}) \|\mathrsfs{Z}_1(s,\cdot) - \mathrsfs{Z}_2(s,\cdot)\|_{\Gamma}.
        \end{aligned}
    \end{equation}
    Define for $i=1,2$
    \begin{equation}
        \mathrsfs{U}_i(t) =\oint  \mathrsfs{M}_i(t,z) \Dif z.
    \end{equation}
    Then it follows using similar arguments as above that
    \begin{equation}\label{eq:bound:Mi}
        \|\mathrsfs{U}_j(s)\| \leq  C(\|\mathcal{K}\|_{\operatorname{op}}, M) \quad \text{and}\quad
        \|\mathrsfs{U}_1(s) - \mathrsfs{U}_2(s)\| 
        \leq C(\|\mathcal{K}\|_{\operatorname{op}}) \|\mathrsfs{Z}_1(s,\cdot) - \mathrsfs{Z}_2(s,\cdot)\|_{\Gamma}.
    \end{equation}
    We also have $\left\|\mu^\top  \mu \right\| \leq  C(\max_{j=1,2}\|\mu_j\|)$. For notational simplicity, we define 
    \begin{equation}
        \mathrsfs{G}_i(t)  \defas  \begin{bmatrix}
            \mathrsfs{Q}_i(t) &  \mathrsfs{U}_i(t) & 0 \\
            \mathrsfs{U}_i(t) ^\top &  \mu^\top \mu & 0\\
            0 &
            0 & \mathrsfs{Q}^{0}_{i,11}(t)
        \end{bmatrix},
    \end{equation}
    where $\mathrsfs{Q}^{0}$ is associated with the polynomial of degree $0$. It is clear from \eqref{eq:bound:Qi}, \eqref{eq:lipschitz:Qi} and \eqref{eq:bound:Mi} that 
    \begin{equation}\label{eq:bound:lipschitz:G}
        \begin{gathered}
            \|\mathrsfs{G}_i(s)\| \leq C(\|\mathcal{K}\|_{\operatorname{op}}, \max_{j=1,2}\|\mu_i\|, M)
            \\
            \text{and}\quad \|\mathrsfs{G}_1(s) - \mathrsfs{G}_2(s) \| \leq C(\|\mathcal{K}\|_{\operatorname{op}})\|\mathrsfs{Z}_1(s,\cdot) - \mathrsfs{Z}_2(s,\cdot)\|_{\Gamma}.
        \end{gathered}
    \end{equation}
    Given that $g$ is $\alpha$-pseudo Lipschitz by Assumption~\ref{assumption:smooth:stats} and from \eqref{eq:bound:lipschitz:G}, we obtain
    \begin{equation}
        \begin{aligned}
        |g(\mathrsfs{G}_1(s) ) - g(\mathrsfs{G}_2(s) )|
        &
        \leq L(g) \|\mathrsfs{G}_1(s)  - \mathrsfs{G}_2(s) \| (1+ \|\mathrsfs{G}_1(s)\|^{\alpha} +\|\mathrsfs{G}_2(s)\|^{\alpha} )
        \\
        &
        \leq 
        C(L(g), \max_{j=1,2}\|\mu_j\|, \|\mathcal{K}\|_{\operatorname{op}}, M, \alpha) \|\mathrsfs{Z}_1(s,\cdot) - \mathrsfs{Z}_2(s,\cdot)\|_{\Gamma}.
        \end{aligned}
    \end{equation}
    Taking the supremum from $0 \leq t \leq T$ and applying Proposition~\ref{prop:stability} completes the proof.
\end{proof}
 As mentioned previously, in Section~\ref{sec:error_bounds}, we show for each fixed $z\in \Gamma^2$ that the terms \eqref{eq:main:error_terms_integrated} and \eqref{eq:main:doob:decomp:AdvHSGD:vanish} vanish as $d\to\infty$ with overwhelming probability. From the definition~\ref{def:main:approx:sols}, we must show these results hold simultaneously for each $z\in \Gamma^2$. Since union bounds only hold over countable sets, we use a net argument similar to \citep{collinswoodfin2023hitting, collinswoodfin2025Exact}.  
 
 Recall we defined the contour $\Gamma = \left\{w\in \mathbb{C}: |w| =\max(1, 2\|\mathcal{K}\|_{\operatorname{op}})\right\}$. For a fixed $\xi>0$, we build a $d^{-\xi}$-mesh of $\Gamma_{\xi}\subset\Gamma$ such that for every $w\in \Gamma$, there exists a $\bar{w}\in \Gamma_{\xi}$ with $|w - \bar{w}| < d^{-\xi}$, with $\Gamma_{\xi}$ having cardinality $|\Gamma_{\xi}| = C(|\Gamma|)d^{\xi}$. Here $|\Gamma|$ denotes the contour length and $|\Gamma_{\xi}|$ denotes the cardinality. We aso introduce the following resolvent Lemma from \citep{collinswoodfin2025Exact}.
\begin{lemma}[Lemma 4 \citep{collinswoodfin2025Exact}] \label{lem:resolvent:bound}
    For $z, \bar{z} \in \mathbb{C}^2 \setminus (\operatorname{spec}(K_1) \times \operatorname{spec}(K_2))$, the following formula holds
    \[
    \| \mathscr{R}(z;\mathcal{K}) - \mathscr{R}(\bar{z};\mathcal{K}) \|_{\opt}  \leq  \sum_{i=1}^2 \left(\prod_{j=1}^i \|R(z_i; K_i)\|_{\opt } \right) |z_i - \bar{z}_i| \left(\prod_{j=i}^2 \|R(\bar{z}_i; K_i)\|_{\opt}\right).
    \] 
\end{lemma}
Note that the normal resolvent $R(z_i;K_i) = (K_i-z_i\Id_d)^{-1}$ satisfies a similar identity for $z_i, \bar{z}_{i} \in \mathbb{C} \setminus \operatorname{spec}(K_i)$
\begin{equation}\label{eq:lem:resolvent:bound:classic}
    \|R(z_i;K_i) - R(\bar{z}_{i};K_i)\|_{\opt} \leq |z_i - \bar{z}_{i}| \|R(z_i;K_i)\|_{\opt}\|R(\bar{z}_{i};K_i)\|_{\opt}.
\end{equation}
 We build upon Lemma A.1 in \citep{collinswoodfin2023hitting} and Lemma 5 in \citep{collinswoodfin2025Exact}.
\begin{lemma}[Net argument] \label{lem:net_argument} Fix $T, M > 0, \xi > 0$. Suppose $\Gamma_{\xi}$ is a $d^{-\xi}$ mesh of $\;\Gamma$ with $|\Gamma_{\xi}| = C d^{\xi}$ for some $C > 0$. Let the function $Z(t,z) = Z(W_{td}, z)$ or $Z(t,z)=Z(\HSGD_t, z)$ satisfy
    \begin{equation}
    \label{eq:net_argument_gamma_xi}
    \sup_{0 \le t \le (T\wedge \tau_M )} \| Z(t, \cdot) - Z(0, \cdot) - \int_0^t \mathscr{F}(\cdot, Z(s, \cdot) ) \, \dif s \|_{\Gamma_\xi} \le \varepsilon
    \end{equation}
    with $\tau_M$ as in Definition~\ref{def:stoppingtime}. Then $Z$ is an $(\varepsilon + C  d^{-\xi}, M, T)$-approximate solution to the integro-differential equation, that is, 
    \[
    \sup_{0 \le t \le (T\wedge \tau_M )} \|Z(t, \cdot) - Z(0, \cdot) - \int_0^t \mathscr{F}(\cdot, Z(s, \cdot) ) \, \dif s \|_{\Gamma} \le \varepsilon + C  d^{-\xi},
    \]
    where $C = C(\bar{\gamma}, \alpha, L(\mathcal{I}), M, \|\mathcal{K}\|_{\operatorname{op}}, L(h), \lambda, \max_{i=1,2}\|\mu_i\|, T)>0$.
    \end{lemma}
\begin{proof}
    We will prove the lemma for $Z(t,z) \defas Z(\HSGD_t, z)$ as the result follows for SGD with a similar argument. We will also work with the stopped process $Z(t,z)\defas Z(t\wedge \tau_{M},z)$ and omit the dependence on the stopping time for notational convenience. For a fixed $z\in \Gamma^2$ with $\bar{z} \in \Gamma_{\xi}^2$ such that $|z_i- \bar{z}_i | < d^{-\xi}$ for $i=1,2$, we have
\begin{equation}\label{eq:S_1}
\begin{aligned} 
\big \| Z(t,z) - 
&
Z(0,z) - \int_0^t \mathscr{F}(z, Z(s, \cdot)) \, \dif s  \big \| \\
&
\le
\| Z(t,z)- Z(t, \bar{z}) \| + \| Z(0, z) - Z(0, \bar{z}) \| 
\\
&
+ \int_0^t \| \mathscr{F}(z, Z(s,\cdot)) - \mathscr{F}(\bar{z}, Z(s, \cdot)) \| \, \dif s 
\\
& 
\quad + \big \| Z(t,\bar{z}) - 
Z(0,\bar{z}) - \int_0^t \mathscr{F}(\bar{z}, Z(s, \cdot)) \, \dif s  \big \|.
\end{aligned}
\end{equation}   
We will bound each of these terms separately. The last term is bounded by $\varepsilon$ by the assumption \eqref{eq:net_argument_gamma_xi} stated in the lemma. To bound the other terms, we first introduce some useful identities. From Proposition~\ref{prop:stability}, we have
\begin{equation}\label{eq:netarg:res:bound}
    \sup_{z_i\in \Gamma}\|R(z_i,K_i)\|_{\operatorname{op}} \leq 2 \quad \text{and}\quad  \sup_{z\in \Gamma^2}\|\mathscr{R}(z, \mathcal{K})\|_{\operatorname{op}} \leq 4.
\end{equation}
For the first term, using Lemma~\ref{lem:resolvent:bound}, \eqref{eq:netarg:res:bound} and recalling that we are working with the stopped process we obtain 
\begin{equation}\label{eq:bound:net:1}
    \begin{aligned}
        \| Z(t,z)- Z(t, \bar{z}) \| 
        &
        \leq \|\HSGD_t\|^2 \| \mathscr{R}(z; \mathcal{K}) - \mathscr{R}(\bar{z}; \mathcal{K})\|_{\operatorname{op}}
        \\
        &
        \leq C \left\| \oint Z(t,\tilde{z}) \Dif \tilde{z}\right\|\sum_{i=1}^2 |z_i - \bar{z}_i|
        \\
        &
        \leq Cd^{-\xi}  \oint \|Z(t,\tilde{z}) \|\Dif |\tilde{z}|
        \\
        &\leq C(\|\mathcal{K}\|_{\operatorname{op}}, M )d^{-\xi} ,
    \end{aligned}
\end{equation}
where we used Cauchy's integral formula and the arc length $|\Gamma|  \leq C(\|\mathcal{K}\|_{\operatorname{op}} )$. The same argument can be used to bound the second term in \eqref{eq:S_1}.  It remains to bound the the third term. We have for any $s\leq t \wedge \tau_M$
\begin{equation}\label{eq:bound:F:netarg:0}
    \begin{aligned}
         \|\mathscr{F}(z, Z(s,\cdot))  - \mathscr{F}(\bar{z}, Z(s,\cdot))\| 
        &
        \leq
        \|\mathscr{F}_{\mathcal{S}}(z, Z(s,\cdot))  - \mathscr{F}_{\mathcal{S}}(\bar{z}, Z(s,\cdot))\| 
        \\
        &\qquad + 2\|\mathscr{F}_{\mathcal{M}}(z, Z(s,\cdot))  - \mathscr{F}_{\mathcal{M}}(\bar{z}, Z(s,\cdot))\|.
    \end{aligned}
    \end{equation}
We will use an analogous notation for $\widehat{\mathcal{B}}_i$ as defined in Definition~\ref{def:stoppingtime} in terms of $Z(t,z)$, $S(t,z)$, $M(t,z)$. Before bounding the terms on the right-hand side of \eqref{eq:bound:F:netarg:0}, we derive some preliminary bounds. By Assumption~\ref{assumption:fisher} and the bound derived in \eqref{eq:bounds:BBV:Lipschitzproof}, we have
\begin{equation}\label{eq:net:growth:b:I}
    |\mathcal{I}_i(\widehat{\mathcal{B}}_i(s))| \leq 
    CL(\mathcal{I})(1+ \|\widehat{\mathcal{B}}_i(s)\|)^{1+\alpha} \leq C(\alpha, L(\mathcal{I}), M, \|\mathcal{K}\|_{\operatorname{op}}).
\end{equation}
From \eqref{eq:bound:net:1}, it follows that
\begin{equation}\label{eq:net:bound:zS(t,z)}
    \begin{gathered}
    \|z_i S(s,z) -\bar{z}_i S(s,\bar{z})\| \leq |z_i- \bar{z}_i|\|S(s,z)\| + |\bar{z}_i| \|S(s,z) - S(s,\bar{z})\| \leq C(\|\mathcal{K}\|_{\opt}, M) d^{-\xi},\\
    \quad \text{and similarly }\quad  \|z_i M(s,z) -\bar{z}_i M(s,\bar{z})\| \leq C(\|\mathcal{K}\|_{\opt}, M) d^{-\xi}.
    \end{gathered}
\end{equation}
From \eqref{eq:bound:net:1}, we also obtain
\begin{equation}\label{eq:oint:inequality:lipschitz}
    \left\| \frac{-1}{2\pi \operatorname{i}}\oint_{\Gamma} S(t, (w, z_2)) - S(t, (w, \bar{z}_2)) \dif w \right\| \leq C(\|\mathcal{K}\|_{\opt}) \sup_{w \in \Gamma} \|S(t, (w, z_2)) - S(t, (w, \bar{z}_2))\| \leq C(\|\mathcal{K}\|_{\opt}, M) d^{-\xi}.
\end{equation}
It is clear that this hold if we integrate over the second position in $z = (z_1, z_2)$. We also bound the trace terms in \eqref{eq:bound:F:netarg:0}. From the construction of $\Gamma_{\xi}$, we have from Lemma~\ref{lem:resolvent:bound} and the resolvent bounds in \eqref{eq:netarg:res:bound}
\begin{equation}\label{eq:bound:net:trace}
    \begin{aligned}
    \frac{1}{d}|\Tr((K_i+\mu_i\mu_i^\top)(\mathscr{R}(z;\mathcal{K})- \mathscr{R}(\bar{z};\mathcal{K}))|  &\leq (\|\mathcal{K}\|_{\operatorname{op}} + \max_{j=1,2} \|\mu_j\|^2)\|\mathscr{R}(z;\mathcal{K}) -\mathscr{R}(\bar{z};\mathcal{K})\|_{\operatorname{op}} 
    \\
    &
    \leq C(\|\mathcal{K}\|_{\operatorname{op}}, \max_{j=1,2} \|\mu_j\|^2) d^{-\xi}.
    \end{aligned}
\end{equation} 
Given the definition of $\mathscr{F}_{\mathcal{S}}$ in \eqref{eq:ODE_resolvent_2}, we have
\begin{equation}
    \begin{aligned}
        \|\mathscr{F}_{\mathcal{S}}(z, Z(s,\cdot))  - \mathscr{F}_{\mathcal{S}}(\bar{z}, Z(s,\cdot))\| 
        &
        \leq  
        \frac{\bar{\gamma}^2}{d}\sum_{i=1}^2 p_ i |\mathcal{I}_i(\widehat{\mathcal{B}}_i(s))||\Tr((K_i+\mu_i\mu_i^\top)(\mathscr{R}(z;\mathcal{K}) -  \mathscr{R}(\bar{z};\mathcal{K}))|
        \\
        &
        + 4\bar{\gamma} \sum_{i=1}^2p_i\|H_{2,i}(\widehat{\mathcal{B}}_i(s))\| \| M_i(s,z) - M_i(s,\bar{z})\| 
        \\
        &
        + 4\bar{\gamma} \;p_1\|H_{1,1}(\widehat{\mathcal{B}}_1(s))\| \left\|  \frac{-1}{2\pi \operatorname{i}}\oint_{\Gamma} \left(S(s,(w, z_2)) -S(s,(w, \bar{z_2}))\right) \dif w\right\| 
        \\
        &
        + 4\bar{\gamma} \;p_2\|H_{1,2}(\widehat{\mathcal{B}}_2(s))\| \left\|  \frac{-1}{2\pi \operatorname{i}}\oint_{\Gamma} \left(S(s,(z_1, w)) -S(s,(\bar{z}_1, w))\right) \dif w\right\| 
        \\
        &
        + 4\bar{\gamma} \sum_{i=1}^2p_i\|H_{1,i}(\widehat{\mathcal{B}}_i(s))\| \|z_iS(s,z) -\bar{z}_iS(s,\bar{z})\| 
        \\
        &
        + 8\bar{\gamma}\sum_{i=1}^2 p_i |\lambda +  \partial_{44}h_i(\widehat{\mathcal{B}}_i(s))|\| S(s,z) -S(s,\bar{z})\|,
    \end{aligned}
    \end{equation}
Hence, from \eqref{eq:main:bound:stability:H}, \eqref{eq:H:Lipschitz:partial44}, \eqref{eq:bound:net:1}, \eqref{eq:net:growth:b:I}, \eqref{eq:net:bound:zS(t,z)}, \eqref{eq:oint:inequality:lipschitz} and \eqref{eq:bound:net:trace} we have
\begin{equation}\label{eq:bound:F:netarg:1}
    \begin{aligned}
        \|\mathscr{F}_{\mathcal{S}}(z, Z(s,\cdot))  - \mathscr{F}_{\mathcal{S}}(\bar{z}, Z(s,\cdot))\| 
        &
        \leq
        C(\bar{\gamma},\lambda, \alpha, M, \|\mathcal{K}\|_{\opt}, L(h), L(\mathcal{I}), \max_{j=1,2}\|\mu_i\|)d^{-\xi}.
    \end{aligned}
    \end{equation}
Using a similar argument, we obtain for $\mathscr{F}_{\mathcal{M}}$
\begin{equation}\label{eq:bound:F:netarg:2}
    \begin{aligned}
        \|\mathscr{F}_{\mathcal{M}}(z, Z(s,\cdot))  - \mathscr{F}_{\mathcal{M}}(\bar{z}, Z(s,\cdot))\| 
        &
        \leq
        C(\bar{\gamma},\lambda, \alpha, M, \|\mathcal{K}\|_{\opt}, L(h), \max_{j=1,2}\|\mu_i\|)d^{-\xi}.
    \end{aligned}
    \end{equation}
Combining \eqref{eq:bound:F:netarg:1}, \eqref{eq:bound:F:netarg:2} and integrating from 0 to $t$, we obtain
\begin{equation}\label{eq:bound:net:2}
    \begin{aligned}
    \int_0^t\|\mathscr{F}(z, S(s,\cdot))  - \mathscr{F}(\bar{z}, S(s,\cdot))\|  \dif s \leq C(\bar{\gamma},\lambda, \alpha, M, \|\mathcal{K}\|_{\opt}, L(h), L(\mathcal{I}), \max_{j=1,2}\|\mu_i\|) d^{-\xi} t.
    \end{aligned}
\end{equation}
Taking the supremum over $z\in \Gamma^2$ and the supremum over $0\leq s\leq t\wedge \tau_M$ in \eqref{eq:S_1} and from the bounds \eqref{eq:bound:net:1}, \eqref{eq:bound:net:2} and the assumption stated in \eqref{eq:net_argument_gamma_xi}, the result follows.
\end{proof}
As mentioned previously, in order to show SGD and AdvHSGD are approximates solutions to \eqref{eq:ODE_resolvent_2}, we must control the error terms which arise in the Doob decompositions of AdvHSGD and SGD respectively. In Definition~\ref{def:stoppingtime} and Theorem~\ref{thm:main:conc:St:SGD}, we introduced the stopping times $\tau_M$ and $\vartheta_M$. We formally introduce $\vartheta_M$ for fixed $M>0$: 
\begin{equation}
    \begin{aligned} \label{eq:main:stopping_time_1}
        \vartheta_{M} 
        &
        \defas \inf \{ t \ge 0 \, : \, \|W_{td} \|^2 > M \quad \text{or} \quad \widehat{B}_1(W_{td}) \not \in \mathcal{U}\quad \text{or}\quad \widehat{B}_2(W_{td}) \not \in \mathcal{U} \}\\
        \text{or} \quad  \vartheta_{M} 
        &
        \defas \inf \{t \ge 0 \, : \, 
        \|\HSGD_t\|^2  > M  \quad \text{or}\quad \widehat{B}_1(\HSGD_t) \not \in \mathcal{U}\quad \text{or}\quad \widehat{B}_2(\HSGD_t) \not \in \mathcal{U} \},     
    \end{aligned}
    \end{equation}
Here we overload the notation so that $\vartheta_{M} $ can be applied to either the iterates of SGD or AdvHSGD. Which stopping time is used will be clear from context. By Lemma~\ref{lem:prop:normGamma:normW}, there exists constants $c,C>0$ such that $\vartheta_{c \cdot M} \le \tau_{M} \le  \vartheta_{C \cdot M}$ implying the stopping times $\tau_M$ and $\vartheta_M$ are interchangeable. When clear from context, we will omit the subscript $\vartheta = \vartheta_M$. Hence, in the following proofs of Propositions \ref{prop:approx:sols:SGD} and \ref{prop:approx:sols:HSGD}, we will instead work with the stopped processes $W_{td}^{\vartheta}=W_{(t\wedge \vartheta)d}$ and $\HSGD^{\vartheta}_{t} = \HSGD_{t\wedge \vartheta}$. The following proof builds upon Proposition 12 in \citep{collinswoodfin2023hitting} and Proposition 5 in \citep{collinswoodfin2025Exact}. 
\begin{proposition}[SGD is an approximate solution]\label{prop:approx:sols:SGD}
Fix $T, M >0$ and $0 < \varepsilon < \frac{1}{2}$. Then $Z(W_{td}, z)$ is an $(Cd^{-\varepsilon}, M, T)$-approximate solution w.o.p, that is
\begin{gather*}
            \sup_{0 \le t \le (T\wedge \tau_M )} \| Z(W_{td}, \cdot) - Z(W_0, \cdot) - \int_0^t \mathscr{F}(\cdot, Z(W_{sd}, \cdot) ) \, \dif s \|_{\Gamma} \le C d^{-\varepsilon}\quad w.o.p.
\end{gather*}
for some constant $C>0$ independent of $d$.
\end{proposition}
\begin{proof}
    Given that the elements of $Z(W_{td}, \cdot)$ are at most quadratic functions of $W_{td}$, applying \eqref{eq:main:integral_1}, \eqref{eq:main:integral_2} and \eqref{eq:main:error_terms_integrated} position-wise, we obtain the following Doob decomposition for $Z(W_{td}, \cdot)$ for a fixed $z\in \Gamma^2$
    \begin{equation}\label{eq:doob:decomp:S(W,SGD)}
        \begin{aligned}
            Z(W_{td}, z) &= W_0^\top \mathscr{R}(z;\mathcal{K})W_0 + \int_0^t \mathscr{F}(z, Z(W_{sd}, \cdot)) \dif s
            \\
        & +  \sum_{j=0}^{\lfloor td \rfloor-1} \left(\Delta \mathcal{M}_{j}^{\text{Grad}}(Z) + \Delta \mathcal{M}_{j}^{\text{Hess}}(Z) + \EE [\mathcal{E}_{j}^{\text{Hess}}(Z) \, | \, \mathcal{F}_{j} ]\right)
        + \xi_{td}(Z). 
        \end{aligned}
    \end{equation}
    See Appendix~\ref{sect:doob:decomp:SGD:derivation} for a detailed derivation on the Doob decomposition for an individual quadratic statistic of $W_{td}$ and Appendix~\ref{sect:DoobDecomp:HSGD} for a detailed derivation of \eqref{eq:doob:decomp:S(W,SGD)} but replacing $\HSGD_t$ by $W_{td}$ and $\rho_{t}$ by $g_{td}$. Taking the supremum from $0 \leq t \leq T \wedge \tau_M$, we then have
    \begin{equation}
        \begin{aligned}
        \sup_{0 \le t \le T \wedge \tau_M}
         &
        \| Z(W_{td}, \cdot) - Z(W_0, \cdot) - \int_0^t \mathscr{F}(\cdot, Z(W_{sd}, \cdot) ) \, \dif s \|
        \\
        &
        \leq
        \sup_{0 \le t \le T \wedge  \vartheta_{CM} } \|\mathcal{M}_{td}^{\text{Grad}}(Z(\cdot, z)) \|+\sup_{0 \le t \le T \wedge  \vartheta_{CM} } \|\mathcal{M}_{td}^{\text{Hess}}(Z(\cdot, z)) \|
\\
&
+
\sup_{0 \le t \le T \wedge \vartheta_{CM} } \|\sum_{j=0}^{\lfloor td \rfloor-1}\EE [\mathcal{E}_{j}^{\text{Hess}}(Z) \, | \, \mathcal{F}_{j} ]  \|+\sup_{0 \le t \le T \wedge  \tau_M } \|\xi_{td}(Z) \|.
        \end{aligned}
    \end{equation}
    where we used the fact that $\tau_M \leq \vartheta_{CM}$ for some constant $C>0$. We defer the proof of the bounds on $\mathcal{M}_{td}^{\text{Grad}}$ and $\mathcal{M}_{td}^{\text{Hess}}$ to Proposition~\ref{prop:gradmartingale} and Proposition~\ref{prop:hessmartingale} in Section~\ref{sec:HSGD_martingale}. We also defer the proof of the bound on the Hessian error term to Proposition~\ref{prop:Hessian-error} in Section~\ref{sec:lower_term_Hessian}.
    
    Fix a constant $\beta >0$ and $\Gamma_\beta\subset \Gamma$ such that there exists a $\bar{z} \in \Gamma_{\xi}^2$ with $|z_i-\bar{z}_i| < d^{-\beta}$ for $i=1,2$ and the cardinality $|\Gamma_\beta| = C(\|\mathcal{K}\|_{\operatorname{op}})d^{\beta}$. From Proposition~\ref{prop:gradmartingale} and Proposition~\ref{prop:hessmartingale} and union bounds since the cardinality of $\Gamma_{\beta}$ is polynomial in $d$, we have that for any $\zeta >0$
    \begin{equation}
        \begin{aligned}
            \sup_{z\in \Gamma_{\beta}}\sup_{0\leq t\leq T} \|\mathcal{M}_{(t\wedge \vartheta_{CM})d}^{\text{Grad}}(Z(\cdot, z))\| < d^{-\frac{1}{2}+\zeta}, \quad w.o.p.
        \end{aligned}
    \end{equation}
    and
    \begin{equation}
        \begin{aligned}
            \sup_{z\in \Gamma_{\beta}}\sup_{0\leq t\leq T} \|\mathcal{M}_{(t\wedge \vartheta_{CM})d}^{\text{Hess}}(Z(\cdot, z))\| < d^{-\frac{1}{2}+\zeta}, \quad w.o.p.
        \end{aligned}
    \end{equation}
    From Proposition~\ref{prop:Hessian-error}, the triangle inequality and a union bound over $\Gamma_{\beta}$, we also have
    \begin{equation}
        \begin{aligned}
            \sup_{z\in \Gamma_{\beta}}\sup_{0\leq t\leq T} \sum_{j=0}^{\lfloor (t\wedge \vartheta_{CM})d \rfloor-1}\|\EE [\mathcal{E}_{j}^{\text{Hess}}(Z) \, | \, \mathcal{F}_{j} ]\| \leq C(L(f))^2 d^{-1/2 + \zeta}.
        \end{aligned}
    \end{equation}
    It remains to bound the last error term $\xi_{td}$. We have
    \begin{equation}
        \begin{aligned}
            \sup_{0 \le t \le T \wedge   \tau_M  } \|\xi_{td}(Z) \|_{\Gamma} 
            &
            \leq  \frac{2}{d}\sup_{0 \le t \le T \wedge   \tau_M } \| \Delta(Z_{td}) \|_{\Gamma}  
            = \frac{2}{d}\sup_{0 \le t \le T \wedge   \tau_M  } \| \mathscr{F}(\cdot, Z(W_{td}, \cdot))\|_{\Gamma} 
            \\
            &
            \leq 
            \frac{2}{d}\sup_{0 \le t \le T \wedge   \tau_M  } (\| \mathscr{F}_{\mathcal{S}}(\cdot, Z(W_{td}, \cdot))\|_{\Gamma} + 2\| \mathscr{F}_{\mathcal{M}}(\cdot, Z(W_{td}, \cdot))\|_{\Gamma} ).
        \end{aligned}
    \end{equation}
We will use an analogous notation for $\widehat{\mathcal{B}}_i$ as defined in \eqref{eq:notation:stability:proof} in terms of $Z(t,z)$, $S(t,z)$, $M(t,z)$ and drop the subscript $j$ as we work with a single approximate solution. For $0 \leq t \leq T \wedge  \tau_M$, we have
\begin{equation}
    \begin{aligned}
        \| \mathscr{F}_{\mathcal{S}}(\cdot, Z(W_{td}, \cdot))\|_{\Gamma}
        &
         \leq 
         C(\bar{\gamma}, \|\mathcal{K}\|_{\operatorname{op}}, M)\sum_{i=1}^2 p_i\|H_{1,i}(\widehat{\mathcal{B}}_i(W_{td}))\|  (1+ \sup_{z \in \Gamma^2} \|z_i S(W_{td}, z)\|)
         \\
         &
         +C(\bar{\gamma})\sum_{i=1}^2 p_i\|H_{2,i}(\widehat{\mathcal{B}}_i(W_{td}))\|  \|M_i(W_{td}, \cdot)\|_{\Gamma}
         \\
         &
         +\frac{C(\bar{\gamma})}{d} \sum_{i=1}^2 p_i \sup_{z\in \Gamma}|\Tr((K_i + \mu_i\mu_i^\top  )\mathscr{R}(z;\mathcal{K}))| \cdot |\mathcal{I}_i(\widehat{\mathcal{B}}_i(W_{td}))|
         \\
         &
         +C(\bar{\gamma}) \sum_{i=1}^2 p_i\|H_{1,i}(\widehat{\mathcal{B}}_i(W_{td}))\| \sup_{z \in \Gamma^2}\|z_i Z(W_{td}, z )\|
         \\
         &
         + C(\bar{\gamma}) \sum_{i=1}^2 p_i|\lambda + \partial_{44} h_i(\widehat{\mathcal{B}}_i(W_{td}))| \|S(W_{td}, \cdot)\|_{\Gamma}.
    \end{aligned}
\end{equation}
Using a similar argument as \eqref{eq:bound:net:trace} from \eqref{eq:netarg:res:bound}, we have
\begin{equation}\label{eq:bound:trace:bound:approxsol:SGD}
    \begin{aligned}
    \frac{1}{d}|\Tr((K_i+\mu_i\mu_i^\top)\mathscr{R}(z;\mathcal{K}))|  &\leq (\|\mathcal{K}\|_{\operatorname{op}} + \max_{j=1,2} \|\mu_j\|^2)\|\mathscr{R}(z;\mathcal{K}) \|_{\opt}\leq C(\|\mathcal{K}\|_{\operatorname{op}}, \max_{j=1,2} \|\mu_j\|^2) .
    \end{aligned}
\end{equation} 
From Assumption~\ref{assumption:risk}, Assumption~\ref{assumption:fisher}, since $\sup_{z\in \Gamma^2}\|z_i Z(W_{td},\cdot)\|_{\Gamma} \leq C(\|\mathcal{K}\|_{\operatorname{op}}, M)$ and from \eqref{eq:main:bound:stability:H}, \eqref{eq:H:Lipschitz:partial44}, \eqref{eq:net:growth:b:I} and \eqref{eq:bound:trace:bound:approxsol:SGD}, it follows that
\begin{equation}
    \begin{aligned}
        \| \mathscr{F}_{\mathcal{S}}(\cdot, Z(W_{td}, \cdot))\|_{\Gamma}
        &
         \leq 
         C(\bar{\gamma},\lambda, \alpha, M, \|\mathcal{K}\|_{\opt}, L(h), L(\mathcal{I}), \max_{j=1,2}\|\mu_i\|).
    \end{aligned}
\end{equation}
Using a similar argument, we obtain
\begin{equation}
    \begin{aligned}
        \| \mathscr{F}_{\mathcal{M}}(\cdot, Z(W_{td}, \cdot))\|_{\Gamma}
        &
         \leq 
         C(\bar{\gamma},\lambda, \alpha, M, \|\mathcal{K}\|_{\opt}, L(h), \max_{j=1,2}\|\mu_i\|).
    \end{aligned}
\end{equation}
Hence, there exist a positive constant $C = C(\bar{\gamma},\lambda, \alpha, M, \|\mathcal{K}\|_{\opt}, L(h), L(\mathcal{I}), \max_{j=1,2}\|\mu_i\|)$ such that 
\[
    \| \mathscr{F}(\cdot, Z(W_{td}, \cdot))\|_{\Gamma} \leq C
\]
from which it follows
\begin{equation}
    \sup_{0 \le t \le T \wedge   \tau_M  } \|\xi_{td}(S) \|_{\Gamma} \leq C d^{-1}\quad a.s.
\end{equation}
Combining all these bounds, we then obtain from a union over $\Gamma_{\beta}^2$
\begin{equation}
    \sup_{0 \le t \le T \wedge \tau_M}
        \| Z(W_{td}, \cdot) - Z(W_0, \cdot) - \int_0^t \mathscr{F}(\cdot, Z(W_{sd}, \cdot) ) \, \dif s \|_{\Gamma_\beta} \leq C d^{\zeta - \frac{1}{2}} \quad w.o.p.
\end{equation}
The result follows from the net argument presented in Lemma~\ref{lem:net_argument} with $\zeta = \frac{1}{2}-\beta$ for any $\beta \in (0, \frac{1}{2})$.
\end{proof}
The proof of Proposition~\ref{prop:approx:sols:HSGD} follows a similar argument to the proof of Proposition~\ref{prop:approx:sols:SGD} and builds upon Proposition 13 in \citep{collinswoodfin2023hitting}.
\begin{proposition}[AdvHSGD is an approximate solution]\label{prop:approx:sols:HSGD}
    Fix $T, M >0$ and $0 < \varepsilon < \frac{1}{2}$. Then $Z(\HSGD_t, z)$ is an $(Cd^{-\varepsilon}, M, T)$-approximate solution w.o.p, that is
    \begin{gather*}
                \sup_{0 \le t \le (T\wedge \tau_M )} \| Z(\HSGD_t, \cdot) - Z(\HSGD_0, \cdot) - \int_0^t \mathscr{F}(\cdot, Z(\HSGD_s, \cdot) ) \, \dif s \|_{\Gamma} \le Cd^{-\varepsilon}\quad w.o.p.
    \end{gather*}
    for some constant $C>0$.
    \end{proposition}
    \begin{proof}
        Given that the elements of $Z(\HSGD_t, \cdot)$ are quadratics, from Itô's lemma and \eqref{eq:main:doob:decomp:AdvHSGD}, we obtain the Doob decomposition for $Z(\HSGD_t, \cdot)$
        \begin{equation}
            \begin{aligned}
                Z(\HSGD_t, z) = W_0^\top \mathscr{R}(z;\mathcal{K})W_0 + \int_0^t \mathscr{F}(z, Z(\HSGD_s, \cdot)) \dif s + \int_0^t \dif \mathcal{M}_s^{\operatorname{AdvHSGD}}(Z(\HSGD_s, \cdot)).
            \end{aligned}
        \end{equation}
        See Appendix~\ref{sect:DoobDecomp:HSGD} for a detailed derivation. Showing $Z(\HSGD_t, \cdot)$ is an approximate solution is equivalent to bounding the martingale term. It then follows that for any $z\in \Gamma^2$
        \begin{equation}
            \begin{aligned}
                \sup_{0 \le t \le (T\wedge \tau_M )}\| &Z(\HSGD_{t},z) - Z(W_0,z) - \int_0^t \mathscr{F}(\cdot, Z(\HSGD_{s}, \cdot) ) \, \dif s \| 
                \\
                &\leq \sup_{0 \le t \le \vartheta_{CM}   \wedge T}  \|\mathcal{M}_s^{\operatorname{AdvHSGD}}(Z(\HSGD_s, z))\|
            \end{aligned}
        \end{equation}
        where we used the fact that $ \tau_{M} \leq \vartheta_{CM}$ for some constant $C>0$. We defer the proof of the bound on $\mathcal{M}_s^{\operatorname{AdvHSGD}}$ to Proposition~\ref{prop:HSGD_Martingale_Bound} in Section~\ref{sec:HSGD_martingale}. Fix a constant $\xi >0$ and $\Gamma_\xi\subset \Gamma$ such that there exists a $\bar{z} \in \Gamma_{\xi}^2$ such that $|z_i-\bar{z}_i| < d^{-\xi}$ for $i=1,2$ and the cardinality $|\Gamma_\xi| = C(\|\mathcal{K}\|_{\operatorname{op}})d^{\xi}$. From Proposition~\ref{prop:HSGD_Martingale_Bound} and a union bound since the cardinality of $\Gamma_{\xi}$ is polynomial in $d$, we have that for any $\zeta>0$
        \begin{equation}
            \sup_{z\in \Gamma_{\xi}}\sup_{0 \le t \le T} \| \mathcal{M}^{\operatorname{AdvHSGD}}_{t \wedge \vartheta} (Z(\cdot,z)) \| \le C L(f) \, \, d^{\zeta / 2 - 1/2}\quad w.o.p.
        \end{equation}
        It then follows that 
        \begin{equation}
            \begin{aligned}
                \sup_{0 \le t \le (T\wedge \tau_M )}\| Z(\HSGD_{t}, \cdot) - Z(W_0, \cdot) - \int_0^t \mathscr{F}(\cdot, Z(\HSGD_{s}, \cdot) ) \, \dif s \|_{\Gamma_{\xi}} \leq  C L(f) \, \, d^{\zeta / 2 - 1/2}\quad w.o.p.
            \end{aligned}
        \end{equation}
        The result follows from the net argument presented in Lemma~\ref{lem:net_argument} with $\zeta = 1-2\xi$ for any $\xi \in (0, \frac{1}{2})$.
    \end{proof}
Through Proposition~\ref{prop:stability} on the stability of approximate solutions, we finally formalize the relationship between SGD, AdvHSGD and $\mathcal{Z}(t,\cdot)$ the solution to \eqref{eq:ODE_resolvent_2}.
\begin{theorem}\label{thm:main:S(W)}
    For all $\varepsilon \in (0,\frac{1}{2})$ such that for any $T, M >0$ and $d$ sufficiently large, with overwhelming probability
    \begin{equation}
        \begin{gathered}
          \sup_{0 \leq t \leq T \wedge \widehat{\tau}_{M}(Z(W,\cdot), \mathcal{Z})} \! \! \! \! \!\! \! \! \! \| 
          Z(W_{td }, \cdot)
          - 
          \mathcal{Z}(t,\cdot)
          \|_{\Gamma} \leq Cd^{-\varepsilon}, 
          \\
           \sup_{0 \leq t \leq T \wedge \widehat{\tau}_{M}(Z(\HSGD,\cdot), \mathcal{Z})} \| 
          Z(\HSGD_t, \cdot)
          - 
          \mathcal{Z}(t,\cdot)
          \|_{\Gamma} \leq Cd^{-\varepsilon},
          \\
             \text{and} \qquad \sup_{0 \leq t \leq T \wedge \widehat{\tau}_M(Z(W,\cdot), Z(\HSGD,\cdot) )} \! \! \! \! \!\! \! \! \!  \| 
          Z(W_{ td }, \cdot)
          - 
          Z(\HSGD_t,\cdot)
          \|_{\Gamma} \leq Cd^{-\varepsilon},
          \end{gathered}
        \end{equation}
      where $C>0$ is a constant independent of $d$, where $\mathcal{Z}(t,z)$ is the deterministic solution to the integro-differential equation \eqref{eq:ODE_resolvent_2} and
      \[
      \widehat{\tau}_{M}(\mathrsfs{Z}_1, \mathrsfs{Z}_2) = \min\{\tau_M(\mathrsfs{Z}_1), \tau_M(\mathrsfs{Z}_2) \}.
      \]
\end{theorem}
\begin{proof}
    The proof follows exactly the same steps as Theorem 14 in \citep{collinswoodfin2023hitting}. Similarly, we provide the details for the case $\mathrsfs{Z}_1(t,\cdot)= Z(\HSGD_t,\cdot)$ and  $\mathrsfs{Z}_2(t,\cdot)= Z(W_{td},\cdot)$ as the argument is analogous in both cases $\mathrsfs{Z}_1(t,\cdot)= Z(\HSGD_t,\cdot)$ and  $\mathrsfs{Z}_2(t,\cdot)= \mathcal{Z}(t,\cdot)$ and $\mathrsfs{Z}_1(t,\cdot)= Z(W_{td},\cdot)$ and  $\mathrsfs{Z}_2(t,\cdot)= \mathcal{Z}(t,\cdot)$. For notational convenience, we write $ \widehat{\tau}_{M}=\widehat{\tau}_{M}(\mathrsfs{Z}_1, \mathrsfs{Z}_2)$.

    By Proposition~\ref{prop:approx:sols:SGD} and Proposition~\ref{prop:approx:sols:HSGD}, for some $\varepsilon \in (0,\frac{1}{2})$ we've shown $Z(W_{td},\cdot)$ and $Z(\HSGD_t,\cdot)$ are respectively $(Cd^{-\varepsilon}, M, T)$-approximate solutions with overwhelming probability. (Recall by definition~\ref{def:main:approx:sols} the deterministic solution $\mathcal{Z}(t,\cdot)$ is an $(0, M, T)$.) Hence, by the stability result on solutions of \eqref{eq:ODE_resolvent_2} proven in Proposition~\ref{prop:stability}, it follows that
    \begin{equation}
        \sup_{0 \leq t \leq T \wedge \widehat{\tau}_M} \; \; \! \! \!\! \! \! \!  \| 
      Z(W_{ td }, \cdot)
      - 
      Z(\HSGD_t,\cdot)
      \|_{\Gamma} \leq Cd^{-\varepsilon}, w.o.p.
    \end{equation}
    which concludes the proof.
  \end{proof}
\subsection{A Sufficient Condition to Remove the Stopping-time}\label{sect:remove:stopping:time}
In the previous section, we showed in theorem~\ref{thm:main:S(W)} that $Z(\HSGD_t,\cdot)$ and $Z(W_{td},\cdot)$ both concentrate around the deterministic solution $\mathcal{Z}(t,z)$ of the integro-differential equation \eqref{eq:ODE_resolvent_2} under the conditions that both processes behave sufficiently regularly. We now introduce a two stopping times which we will use to show that only one of SGD or AdvHSGD must behave sufficiently regularly is enough for theorem~\ref{thm:main:S(W)} to hold. Let $\mathrsfs{Z}(t,\cdot)$ be an $(\varepsilon, M, T)$-approximate solution and define the stopping time for $M, \omega > 0$
\begin{equation}\label{eq:def:tau_M,omega:stoppingtime}
    \begin{aligned}
    \tau_{M,\omega} \defas \inf\bigg\{t\geq 0: \|\mathrsfs{Z}(t,\cdot)\|_{\Gamma} &> M \quad \text{or}\quad \inf_{V \in \mathcal{U}^c}\|\widehat{\mathcal{B}}_1(t) - V\| \leq \omega 
    \\
    &\text{or}\quad \inf_{V \in \mathcal{U}^c}\|\widehat{\mathcal{B}}_2(t) - V\| \leq \omega\bigg\},
    \end{aligned}
\end{equation}
where $\widehat{\mathcal{B}}_i(t)$ is defined in Definition~\ref{def:stoppingtime} and $\mathcal{U}^c$ is the complement of $\mathcal{U}$ as defined in Assumption~\ref{assumption:risk}. It is clear that $ \tau_{M,0} =  \tau_{M}$ as defined in Definition~\ref{def:stoppingtime}. The next theorem will replace $ \widehat{\tau}_{M}(\mathrsfs{Z}_1,\mathrsfs{Z}_2)$ in Theorem~\ref{thm:main:S(W)} by the following stopping time which requires only one of the approximate solutions to stay bounded
\begin{equation}\label{eq:def:nu:stoppingtime}
    \begin{aligned}
        \Theta_{M,\omega}(\mathrsfs{Z}_1, &\mathrsfs{Z}_2) \defas \max\left\{ \inf_\{t\geq 0: \|\mathrsfs{Z}_j(t,\cdot)\|_{\Gamma} > M \}\quad \text{for}\quad j=1,2 \right\}\\
        &\quad\wedge \max\left\{\inf\{t\geq 0:\inf_{V \in \mathcal{U}^c}\|\widehat{\mathcal{B}}_i(\mathrsfs{X}_j(t,\cdot)) - V\| \leq \omega\} \quad \text{for}\quad i=1,2\quad j=1,2\right\}
    \end{aligned}
\end{equation}
where we write $\widehat{\mathcal{B}}_i(\mathrsfs{X}_j(t,\cdot))$ to emphasize the dependence on $\mathrsfs{X}_j$. Note that $\Theta = \Theta_{M, 0}$ in Section~\ref{sect:mainresult:SGD}.
\begin{theorem}\label{thm:main:S(W):oneST}
    There exists an $\varepsilon \in (0, \frac{1}{2})$ such that for any $T, M, \omega >0$ and $d$ sufficiently large, with overwhelming probability
    \begin{equation}
        \begin{gathered}
          \sup_{0 \leq t \leq T \wedge \Theta_{M,\omega}(Z(W, \cdot), \mathcal{Z})} \! \! \! \! \!\! \! \! \! \| 
          Z(W_{td }, \cdot)
          - 
          \mathcal{Z}(t,\cdot)
          \|_{\Gamma} \leq Cd^{-\varepsilon}, 
          \\
           \sup_{0 \leq t \leq T \wedge \Theta_{M,\omega}(Z(\HSGD,\cdot), \mathcal{Z})} \| 
          Z(\HSGD_t, \cdot)
          - 
          \mathcal{Z}(t,\cdot)
          \|_{\Gamma} \leq Cd^{-\varepsilon},
          \\
             \text{and} \qquad \sup_{0 \leq t \leq T \wedge \Theta_{M,\omega}(Z(W,\cdot), Z(\HSGD,\cdot) )} \! \! \! \! \!\! \! \! \!  \| 
          Z(W_{ td }, \cdot)
          - 
          Z(\HSGD_t,\cdot)
          \|_{\Gamma} \leq Cd^{-\varepsilon},
          \end{gathered}
        \end{equation}
        where $C>0$ is a constant independent of $d$, where $\mathcal{Z}(t,z)$ is the deterministic solution to the integro-differential equation \eqref{eq:ODE_resolvent_2}.
\end{theorem}
\begin{proof}
    Fix a $\omega >0$ and for two approximate solutions $\mathrsfs{Z}_1$ and $\mathrsfs{Z}_2$, define the stopping time
    \begin{equation}
        \widehat{\tau}_{M+1, 0}^{\mathrsfs{Z}_1, \mathrsfs{Z}_2} \defas \min\{\tau_{M+1, 0}(\mathrsfs{Z}_1), \tau_{M+1, 0}(\mathrsfs{Z}_2)\}.
    \end{equation}
    Similarly to the proof of theorem~\ref{thm:main:S(W)}, the proof follows similar steps in all three cases $\mathrsfs{Z}_1(t,\cdot)= Z(\HSGD_t,\cdot)$ and  $\mathrsfs{Z}_2(t,\cdot)= Z(W_{td},\cdot)$, $\mathrsfs{Z}_1(t,\cdot)= Z(\HSGD_t,\cdot)$ and  $\mathrsfs{Z}_2(t,\cdot)= \mathcal{Z}(t,\cdot)$, $\mathrsfs{Z}_1(t,\cdot)= Z(W_{td},\cdot)$ and  $\mathrsfs{Z}_2(t,\cdot)= \mathcal{Z}(t,\cdot)$. We omit the details of the proof as it follows almost identical steps as the proof of Theorem 15 in \citep{collinswoodfin2023hitting} after noting that for any $0\leq t\leq T\wedge \widehat{\tau}_{M+1, 0}$ and $i=1,2$
    \begin{equation}
        \|\widehat{\mathcal{B}}_i(\mathrsfs{Z}_1(t, \cdot))- \widehat{\mathcal{B}}_i(\mathrsfs{Z}_2(t, \cdot))\| \leq C(\|\mathcal{K}\|_{\operatorname{op}}) \|\mathrsfs{Z}_1(t, \cdot) - \mathrsfs{Z}_2(t, \cdot)\|_{\Gamma} \leq C(\|\mathcal{K}\|_{\operatorname{op}})d^{-\varepsilon}\quad w.o.p.
    \end{equation}
    which follows from Theorem~\ref{thm:main:S(W)} and \eqref{eq:bound:B(s):lipschitz}.
\end{proof}
It immediately follows from Theorem~\ref{thm:main:S(W):oneST} that SGD and AdvHSGD concentrate around $\mathcal{Z}(t,z)$ provided that AdvHSGD or $\mathcal{Z}(t,z)$ stay bounded. Analogously to Corollary 4 in \citep{collinswoodfin2023hitting} and , we formalize this in the following Corollary. We omit the proof as it is identical.
\begin{corollary}
    [Bounded $\mathrsfs{N}$ and concentration] \label{cor:bounded_iterates} Recall we defined the statistic
    $\mathrsfs{N}(t) = \oint\Tr(\mathcal{Z}(t,z)) \Dif z$. Suppose for a fixed $T,  \omega> 0$ that
\begin{equation} \label{eq:concentration_condition}
 \sup_{0 \le t \le T} \mathrsfs{N}(t) \le M  \quad \text{and for $i=1,2$} \quad \inf_{0 \le t \le T} \inf_{V \in \mathcal{U}^c} \|\widehat{\mathrsfs{B}}_i(t)-V\| > \eta \quad \text{hold w.o.p.}
\end{equation}
where $M>0$ is independent of $d$. Then there exists $\varepsilon \in (0, \frac{1}{2})$ so that for $d$ sufficiently large, with overwhelming probability, 
  \begin{equation}
  \begin{gathered} 
    \sup_{0 \leq t \leq T } \| 
    Z(\HSGD_t, \cdot)
    - 
    \mathcal{Z}(t,\cdot)
    \|_{\Gamma} \leq Cd^{-\varepsilon} \quad 
       \text{and} \quad \sup_{0 \leq t \leq T} \| 
    Z(W_{td}, \cdot)
    - 
    \mathcal{Z}(t,\cdot)
    \|_{\Gamma} \leq Cd^{-\varepsilon},
    \end{gathered}
  \end{equation}
  where $C>0$ is a constant independent of $d$. Hence, by the triangle inequality, it follows that
  \begin{equation} 
  \sup_{0 \le t \le T} \| Z(W_{td}, \cdot) - Z(\HSGD_t,\cdot)\|_{\Gamma} \le 2Cd^{-\varepsilon}. 
  \end{equation}
\end{corollary}
\begin{remark}\label{remark:equiv:stopping:time:N(t)}
Intuitively, Corollary~\ref{cor:bounded_iterates} is written in terms of a condition on the size of $\mathrsfs{N}(t)$ instead of $\|\mathcal{Z}(t,\cdot)\|_{\Gamma}$ as per Lemma~\ref{lem:prop:normGamma:normW}, there exists constants $c,C>0$ (independent of $d$) such that $c \mathrsfs{N}(t) \leq \|\mathcal{Z}(t, \cdot)\|_{\Gamma} \leq C  \mathrsfs{N}(t)$. Hence, we obtain the inequality for a fixed $M>0$
\begin{equation}
    \inf\{t\geq 0: \|\mathcal{Z}(t,\cdot)\|_{\Gamma} > CM \} \geq  \inf\{t\geq 0: \mathrsfs{N}(t) > M \},
\end{equation}
from which it is not hard to see that Theorem~\ref{thm:main:S(W):oneST} holds. Corollary~\ref{cor:bounded_iterates} also holds if we replace $(\mathrsfs{N}(t), \widehat{\mathrsfs{B}}_1(t), \widehat{\mathrsfs{B}}_2(t))$ in \eqref{eq:concentration_condition} with the equivalents for either SGD or AdvHSGD by Lemma~\ref{lem:prop:normGamma:normW}. 
\end{remark}
Theorem~\ref{thm:main:conc:St:SGD} and Theorem~\ref{thm:main:conc:St:HSGD} follow from direct applications of Theorem~\ref{thm:main:S(W):oneST} and Lemma~\ref{lem:prop:normGamma:normW} as explained in Remark~\ref{remark:equiv:stopping:time:N(t)}. We also provide a sufficient condition for Corollary~\ref{cor:bounded_iterates} to hold.
\begin{proposition}[Non-explosiveness]\label{prop:nonexplosiveness}
    Suppose that the assumptions stated in this framework hold then the function $\mathrsfs{N}(t)$ is the solution to the following ordinary differential equation 
    \begin{equation}\label{eq:N(t):ode}
        \begin{aligned}
        \frac{\dif}{\dif t}\mathrsfs{N}(t)
        &=
        -2\gamma(t)\bigg(\sum_{i=1}^2p_i \E_{v,\epsilon}[f_i'(q_{t,i}, \epsilon) x_{t,i}] +\lambda(\mathrsfs{N}(t) - \|X^\star\|^2 -\Tr(\mu^\top\mu))\bigg)
        \\
        &
        + \frac{\gamma(t)^2}{d}\sum_{i=1}^2 p_i\Tr(K_i + \mu_i\mu_i^\top )  \E_{v,\epsilon}[f_i'(q_{t,i}, \epsilon)^2] ,
        \end{aligned}
    \end{equation}
    where 
    \begin{equation}
        \begin{aligned}
        q_{t,i} &= x_{t,i}\oplus x_i^{\star} =\sqrt{\mathrsfs{B}_i(t)} v  +\mathfrak{m}_i(t) + \begin{bmatrix}
            \delta u_{t,i} \sqrt{\widehat{\mathrsfs{B}}_{i,44}(t)} & 0
        \end{bmatrix}^\top,
        \\
         u_{t,i} &= \operatorname{argmax}_{|u|\leq 1} f_i(\sqrt{\mathrsfs{B}_i(t)} v  +\mathfrak{m}_i(t) + \begin{bmatrix}
            \delta u \sqrt{\widehat{\mathrsfs{B}}_{i,44}(t)} & 0
        \end{bmatrix}^\top, \epsilon),
    \end{aligned}
    \end{equation}
    for $v\sim \mathcal{N}(0, \Id_2)$ and where $\mathrsfs{B}_i(t)$ is defined in \eqref{eq:ODE_resolvent_2}. Suppose further that the objective function $f$ is $\alpha$-pseudo-Lipschitz with $\alpha\leq 1$. Then there exists $C=C(\bar{\gamma}, L(f), \eta, \delta, \|\mathcal{K}\|_{\operatorname{op}}, \lambda, \max_{i=1,2}\|\mu_i\|^2)$ such that 
    \[
    \mathrsfs{N}(t) \leq (1+ \mathrsfs{N}(0))e^{C t},
    \]
    for all time $t$ such that $(\widehat{\mathrsfs{B}}_1(t), \widehat{\mathrsfs{B}}_2(t))$ is in $\mathcal{U} \times \mathcal{U}$.
    \end{proposition}
    \begin{proof}
        To derive \eqref{eq:N(t):ode}, we first apply Itô's lemma to the statistic $\varphi(\WHSGD_t) \defas \|\WHSGD_t\|^2 + \|X^{\star}\|^2 + \|\mu_1\|^2 + \|\mu_2\|^2$ and obtain
        \begin{equation}
            \begin{aligned}
            \dif \varphi(\WHSGD_t) &
            = -2\gamma(t)\langle \WHSGD_t, \nabla \mathcal{R}^{{\operatorname{adv}}}_{\lambda}(\HSGD_t)\rangle \dif t + \frac{\gamma(t)^2}{d} \sum_{i=1}^2 p_i\Tr(K_i + \mu_i\mu_i^\top ) \E_{a,\epsilon}[f_i'(\rho_{t,i}, \epsilon)^2] \dif t 
            \\
            &
            \qquad + \dif \mathcal{M}_t^{\operatorname{AdvHSGD}}(\varphi),
            \end{aligned}
        \end{equation}
        where  $\mathcal{M}_t^{\operatorname{AdvHSGD}}(\varphi)$ is a martingale. By Assumption~\ref{assumption:fisher}, we know there exist functions $\mathcal{I}_i$ for $i=1,2$ such that $\mathcal{I}_i(\widehat{B}_i(\HSGD_t)) = \E_{a,\epsilon}[f_i'(\rho_{t,i}, \epsilon)^2]$. We also have that
        \begin{equation}
            \begin{aligned}
                \langle \WHSGD_t, \nabla \mathcal{R}^{{\operatorname{adv}}}_{\lambda}(\HSGD_t)\rangle &= \sum_{i=1}^2 p_i (\E_{v,\epsilon}[f_i'(\rho_{t,i}, \epsilon)\langle \WHSGD_t, \sqrt{K_i}v + \mu_i\rangle +  \delta u_{t,i}  \|\WHSGD_t\|]   +\lambda \|\WHSGD_t\|^2.
            \end{aligned}
        \end{equation}
        This is a Gaussian expectation hence there exists some functions $A$ such that
        \begin{equation}
            A(\widehat{B}_i(\HSGD_t)) \defas\E_{v,\epsilon}[f_i'(\rho_{t,i}, \epsilon)\langle \WHSGD_t, \sqrt{K_i}v + \mu_i\rangle +  \delta u_{t,i}  \|\WHSGD_t\|].
        \end{equation}
        We then obtain the representation
        \begin{equation}
            \begin{aligned}
            &\dif \varphi(\WHSGD_t)  
            = -2\gamma(t)\bigg(\sum_{i=1}^2 p_i A(\widehat{B}_i(\HSGD_t)) 
            +\lambda (\varphi(\WHSGD_t) - \|X^{\star}\|^2 - \|\mu_1\|^2 - \|\mu_2\|^2)\bigg) \dif t \\
            &+ \frac{\gamma(t)^2}{d} \sum_{i=1}^2 p_i \Tr(K_i + \mu_i\mu_i^\top) \E_{v,\epsilon}[f_i'(q_{t,i}, \epsilon)^2] \dif t  + \dif \mathcal{M}_t^{\operatorname{AdvHSGD}}(\varphi).
            \end{aligned}
        \end{equation}
        Using Cauchy's integral formula, we also have the representation
        \begin{equation}
            \varphi(\WHSGD_t) = 
              \oint \left({S}_{11}(\HSGD_t,z)+ {S}_{22}(\HSGD_t,z) \right) \Dif z + \|\mu_1\|^2 + \|\mu_2\|^2.
        \end{equation}
        From \eqref{eq:HSGD_exact}, we also have
        \begin{equation}
            \dif \varphi(\WHSGD_t) = 
             \frac{-1}{2 \pi i} \oint \left(\mathscr{F}_{\mathcal{S}, 11}(z, Z(\HSGD_t,z))+  \mathscr{F}_{\mathcal{S}, 22}(z, Z(\HSGD_t,z)) \, \right) \Dif z +\dif \mathcal{M}_t^{\operatorname{AdvHSGD}} (\varphi)
        \end{equation}
        for some martingale $ \mathcal{M}_t^{\operatorname{AdvHSGD}} $. Given that the semi-martingale decomposition of an Itô process is unique and that 
        \begin{equation}
            \mathrsfs{N}(t) = \oint\Tr(\mathcal{Z}(t,z))\Dif z = \oint \left(\mathcal{S}_{11}(t,z)+ \mathcal{S}_{22}(t,z)\right)\Dif z + \|\mu_1\|^2 + \|\mu_2\|^2
        \end{equation} 
        then the terms with finite quadratic variation are equal from which it follows that
        \begin{equation}
            \begin{aligned}
            \frac{\dif}{\dif t}\mathrsfs{N}(t)
            &=
             -2\gamma(t)\bigg(\sum_{i=1}^2 p_i A(\widehat{\mathrsfs{B}}_i(t)) + \lambda (\mathrsfs{N}(t) - \|X^\star\|^2 - \|\mu_1\|^2 - \|\mu_2\|^2)\bigg)
             \\
             &
             + \frac{\gamma(t)^2}{d} \sum_{i=1}^2 p_i \Tr(K_i + \mu_i\mu_i^\top) \E_{v,\epsilon}[f_i'(q_{t,i}, \epsilon)^2] \dif t .
            \end{aligned}
        \end{equation}
        This completes the derivation of \eqref{eq:N(t):ode}. From Assumption~\ref{assumption:pseudo_lipschitz}, the function $f$ is $\alpha$-pseudo-Lipschitz with $\alpha \leq 1$. Hence, from Lemma~\ref{lem:growth_grad_f} we obtain the almost sure bound
        \begin{equation}
            |f_i'(q_{t,i}, \epsilon)| \le C (L(f)) (1 + \|q_{t,i}\| + |\epsilon|)\quad a.s.
        \end{equation}
        Squaring and taking the expectation, we obtain
        \begin{equation}
            \begin{aligned}
            \E_{v,\epsilon}[  |f_i'(q_{t,i}, \epsilon)|^2] &
            % \le 
            % C (L(f)) \E_{v,\epsilon}[ (1 + \|q_{t,i}\| + |\epsilon|)^2] 
            % \\
            % &
            % \leq
            %  C (L(f))\E_{v,\epsilon}[ (1 + \|q_{t,i}\|^2 + |\epsilon|^2)] 
            %  \\
            %  &
            %  \leq C (L(f), \eta, \delta)(1 + \E_{v}[\|\sqrt{\mathrsfs{B}_i(t)} v +\mathfrak{m}_i(t)\|^2] +   \widehat{\mathrsfs{B}}_{i,44}(t))
            %  \\
            %  &
            %  \leq C (L(f), \eta, \delta)(1 + \Tr(\mathrsfs{B}_i(t)) + \|\mathfrak{m}_i(t)\|^2 +   \mathrsfs{N}(t))
            %  \\
            %  &
             \leq C (L(f), \eta, \delta, \max_{i=1,2}\|\mu_i\|^2, \|\mathcal{K}\|_{\operatorname{op}})(1 +  \mathrsfs{N}(t))
            \end{aligned}
        \end{equation}
        since $\widehat{\mathrsfs{B}}_{i,44}(t) = \mathrsfs{N}(t) - \|X^{\star}\|^2- \Tr(\mu^\top\mu )$, since by Cauchy--Schwarz and Lemma~\ref{lem:prop:normGamma:normW}
        \begin{equation}
            \begin{gathered}
            \Tr(\mathrsfs{B}_i(t))  \leq C\|\mathrsfs{B}_i(t)\| \leq C(\|\mathcal{K}\|_{\operatorname{op}}) \|\mathcal{S}(t,\cdot)\|_{\Gamma}  \leq C(\|\mathcal{K}\|_{\operatorname{op}}) \mathrsfs{N}(t),\\
            \text{and from \eqref{eq:upper bound:norm:m(t)}}\quad\|\mathfrak{m}_i\|^2  \leq \|\mu_i\|^2 \mathrsfs{N}(t).
            \end{gathered}
        \end{equation}
        Similarly, we obtain from Cauchy--Schwarz 
        \begin{equation}
            \begin{aligned}
               &\left| \E_{v,\epsilon}[f_i'(q_{t,i}, \epsilon) x_{t,i}]  \right|^2
                % \leq
                % \E_{v,\epsilon}[|f_i'(q_{t,i}, \epsilon)|^2] \E_{v,\epsilon}[|x_{t,i}|^2] 
                % \\ 
                % &
                % \leq
                % C (L(f), \eta, \delta,  \|\mathcal{K}\|_{\operatorname{op}})(1 +  \mathrsfs{N}(t))\left(\mathrsfs{B}_{i,11}(t) + \left([\mathfrak{m}_{i}(t) ]_1+ \delta u_{t,i} \sqrt{\widehat{\mathrsfs{B}}_{i,44}(t)}\right)^2\right)
                % \\
                % &
                % \leq
                % C (L(f), \eta, \delta,  \|\mathcal{K}\|_{\operatorname{op}})(1 +  \mathrsfs{N}(t))(\|\mathrsfs{B}_{i}(t)\| + \|\mathfrak{m}_i(t)\|^2 + \widehat{\mathrsfs{B}}_{i,44}(t))
                % \\
                % &
                \leq
                C (L(f), \eta, \delta,\max_{i=1,2}\|\mu_i\|^2, \|\mathcal{K}\|_{\operatorname{op}})(1 +  \mathrsfs{N}(t))^2.
            \end{aligned}
        \end{equation}
        Combining these inequalities, we obtain
        \begin{equation}
            \begin{aligned}
            \frac{\dif}{\dif t}\mathrsfs{N}(t)
            &\leq C( 1+ \mathrsfs{N}(t)).
            \end{aligned}
        \end{equation}
        where $C = C(\bar{\gamma}, L(f), \eta, \delta, \|\mathcal{K}\|_{\operatorname{op}}, \lambda, \max_{i=1,2}\|\mu_i\|^2)$. The result follows from an application of Grönwall's inequality $(1 + \mathrsfs{N}(t)) \leq (1 + \mathrsfs{N}(0))e^{Ct}$.
    \end{proof}

\subsection{Concentration Result for General Statistics}\label{sect:conc:proof:any:stat}
In this section, we prove a result for general statistics satisfying the following assumption.
\begin{assumption}[Smoothness of statistics, $\varphi$]\label{assumption:smooth:stats}
    Denote $\hat{X} = X \oplus X^{\star}$ and $W = \hat{X} \oplus \mu$. The statistic $\varphi$ satisfies the composite structure,
   \[
        \varphi(X) = g\left(\begin{bmatrix}
            \hat{X}^\top q(K_1, K_2) \hat{X} & \hat{X}^\top \mu &0 \\
            \mu^\top \hat{X} & \mu^\top \mu & 0 \\
            0 & 0 &  \|X\|^2
         \end{bmatrix}\right).
    \]
    where $g: \R^{5\times 5} \to \R$ is $\alpha$-pseudo-Lipschitz on $\overline{\mathcal{U}}$ as defined in Assumption~\ref{assumption:smooth:stats:Adv} and $q$ is a polynomial.
\end{assumption}
 As presented in the proof of Proposition~\ref{prop:stability:general}, recall the deterministic equivalents of $\varphi(\WHSGD_t)$ and $\varphi(W_{td})$ are defined as
\begin{equation}\label{eq:det:equiv:General}
    \phi(t) \defas g\left(\begin{bmatrix}
        \oint q(z)\mathcal{S}(t,z)\Dif z &  \oint \mathcal{M}(t,z)\Dif z & 0  \\
        \oint \mathcal{M}(t,z)^\top \Dif z & \mu^\top  \mu  & 0 \\
        0 & 0 & \oint\mathcal{S}_{11}(t,z)\Dif z
    \end{bmatrix}\right),
\end{equation}
where $\mathcal{S}(t,z)$ and $\mathcal{M}(t,z)$ are block submatrices of $\mathcal{Z}(t,z)$ which solves \eqref{eq:ODE_resolvent_2}. We are now ready to state our main result for the concentration of $\varphi(\WHSGD_t)$ and $\varphi(W_{td})$ around $ \phi(t)$.
\begin{theorem}[Concentration of any statistic] \label{thm:concentration_statistic} Suppose that the assumptions stated in this framework hold. Then there exists an $\varepsilon \in (0, \frac{1}{2})$ such that for any $T, M >0$ and $d$ sufficiently large, with overwhelming probability 
  \begin{equation}
  \begin{gathered}
    \sup_{0 \leq t \leq T \wedge \Theta_{M}(Z(W,\cdot), \mathcal{Z})} | 
    \varphi(W_{ td })
    - 
    \phi(t)
    |\leq Cd^{-\varepsilon}, \\
    \sup_{0 \leq t \leq T \wedge  \Theta_{M}(Z(\HSGD,\cdot), \mathcal{Z}) }| 
    \varphi(\HSGD_t)
    - 
    \phi(t)
    | \leq Cd^{-\varepsilon},
    \\
       \text{and} \qquad \sup_{0 \leq t \leq T \wedge  \Theta_{M}(Z(W,\cdot), Z(\HSGD,\cdot))} | 
    \varphi(W_{td})
    - 
    \varphi(\HSGD_t)
    |\leq Cd^{-\varepsilon},
    \end{gathered}
  \end{equation}
where $C>0$ is a constant independent of $d$, where $\phi$ is defined in \eqref{eq:det:equiv:General} and where the stopping time $\Theta_{M}(\mathrsfs{Z}_1, \mathrsfs{Z}_2)\equiv\Theta_{M, 0}(\mathrsfs{Z}_1, \mathrsfs{Z}_2)$ is defined in \eqref{eq:def:nu:stoppingtime}. 
\end{theorem}
\begin{proof}
    This proof builds upon the proof of Theorem 15 and Theorem 16 in \citep{collinswoodfin2023hitting}. Similarly to the proof of Theorem~\ref{thm:main:S(W):oneST}, we define the stopping time 
    \begin{equation}
        \widehat{\tau}_{M+1, \eta}^{\mathrsfs{Z}_1, \mathrsfs{Z}_2} \defas \min\{\tau_{M+1, \eta}(\mathrsfs{Z}_1), \tau_{M+1, \eta}(\mathrsfs{Z}_2)\},
    \end{equation}
    and consider the case $\mathrsfs{Z}_1(t,\cdot)= Z(\HSGD_t,\cdot)$ and  $\mathrsfs{Z}_2(t,\cdot)= Z(W_{td},\cdot)$ as the other cases follow analogously. Fix $T, M>0$. By Proposition~\ref{prop:approx:sols:HSGD} and Proposition~\ref{prop:approx:sols:SGD}, there exists a $\varepsilon \in (0,\frac{1}{2})$ and constant $C>0$ such that $Z(\HSGD_t, z)$ and $Z(W_{td},z)$ are $(Cd^{-\varepsilon}, M+1, T)$-approximate solutions. (Recall by definition~\ref{def:main:approx:sols} the deterministic solution $\mathcal{Z}(t,\cdot)$ is an $(0, M, T)$.) By Cauchy's integral formula, we have
    \begin{equation}
        \begin{gathered}
         \widehat{\WHSGD}_t^\top q(K_1,K_2)  \widehat{\WHSGD}_t = \oint q(z) S(\HSGD_t,z)\Dif z, \quad \hat{X}_{td}^\top q(K_1, K_2)  \hat{X}_{td}  =\oint q(z) S(W_{td},z)\Dif z,
         \\
         \widehat{\WHSGD}_t^\top  \mu  = \oint M(\HSGD_t,z)\Dif z,\quad  \hat{X}_{td}^\top \mu  =\oint  M(W_{td},z)\Dif z,
         \\
         \|\WHSGD_t\|^2 = \oint S_{11}(\HSGD_t,z)\Dif z\quad\text{and}\quad \|X_{td}\|^2 = \oint S_{11}(W_{td},z)\Dif z.
        \end{gathered}
    \end{equation}
    Hence from Proposition~\ref{prop:stability:general}, it follows that
    \begin{equation}
        \sup_{0\leq t \leq T\wedge \widehat{\tau}_{M+1, 0}}|\varphi(\HSGD_t) - \varphi(W_{td})|\leq C d^{-\varepsilon}, \quad w.o.p.
    \end{equation}
    Using the same argument as the proof of Theorem 15 in \citep{collinswoodfin2023hitting} noting that remark in the proof of Theorem~\ref{thm:main:S(W):oneST}, the stopping time $\widehat{\tau}_{M+1, 0}$ can be replaced with $\Theta_{M,0}$ for sufficiently large $d$.
\end{proof}
\begin{corollary}\label{cor:statistic_bounded}
    Suppose that the assumptions stated in this framework hold and the assumptions in Corollary~\ref{cor:bounded_iterates} hold. Then there exists an $\varepsilon \in (0, \frac{1}{2})$ such that for any $T>0$ and sufficiently large $d$, with overwhelming probability, 
  \begin{equation}
  \begin{gathered} 
    \sup_{0 \leq t \leq T } | 
    \varphi(\WHSGD_t)
    - 
    \phi(t)
    | \leq Cd^{-\varepsilon} \quad 
       \text{and} \qquad \sup_{0 \leq t \leq T} | 
    \varphi(X_{td})
    - 
    \phi(t)
    | \leq Cd^{-\varepsilon}.
    \end{gathered}
  \end{equation}
  Hence, by the triangle inequality, it follows that
  \begin{equation} 
  \sup_{0 \le t \le T} | \varphi(X_{td}) - \varphi(\WHSGD_t) | \le 2Cd^{-\varepsilon}. 
  \end{equation}
\end{corollary}
\subsection{Concentration of General Statistics Involving the Minimizer of the $\ell_2$-Adversarial Risk}\label{sect:general:adv:concentration}
Recall from Definition~\ref{def:adversarial:ground:truth} the minimizer of the adversarial population risk is defined as follows
\begin{definition}
    Let $X^{\star,\operatorname{adv}}$ denote the adversarial ground truth. In other words, we have
    \begin{equation}
        X^{\star,\operatorname{adv}} \in \underset{X}{\operatorname{argmin}}  \left\{\mathcal{R}^{\operatorname{adv}}(X) =\EE_{a,i,\epsilon} [\max_{|s| \leq 1}f_i( X^\top a  + \delta s \|X\|, \epsilon) ]\right\}.
    \end{equation}
\end{definition}
To study the distance to adversarial optimality, we must compute the deterministic equivalent of the statistic $\|X- X^{\star,\operatorname{adv}}\|^2$. We cannot use a similar argument as presented in the proof of Proposition~\ref{prop:nonexplosiveness} given that the first expectation in the drift contribution from Itô's lemma applied to $\|\WHSGD_t - X^{\star,\operatorname{adv}}\|^2$, involves an expectation over the jointly Gaussian scalars $\WHSGD_t^\top a_i$, $(X^\star)^\top a$ and $(X^{\star,\operatorname{adv}})^\top a_i$. Consequently, we must compute the deterministic equivalent for the mean and covariance matrix of the gaussian vector $\widetilde{X}_t^\top a_i$ for $\widetilde{X}_t \defas X_t  \oplus X^\star\oplus X^{\star, \Adv}$ which we define below for $i=1,2$:
\begin{equation}\label{def:det:equiv:covariance:adv}
    \mathrsfs{B}^{\Adv}_i(t) \defas \begin{bmatrix}
        \mathrsfs{B}_{i,11}(t) &  \mathrsfs{B}_{i,12}(t) &\mathrsfs{A}_{K_i}(t) \\
        \mathrsfs{B}_{i,12}(t)  &  \mathrsfs{B}_{i,22}(t) & ( X^{\star})^\top K_i  X^{\star,\operatorname{adv}} \\
        \mathrsfs{A}_{K_i}(t) & ( X^{\star})^\top K  X^{\star,\operatorname{adv}}& ( X^{\star,\operatorname{adv}})^\top K_i  X^{\star,\operatorname{adv}}
    \end{bmatrix} \in \R^{3\times 3},
    \end{equation}
and $\mathfrak{m}_i^\Adv(t) = \mathfrak{m}_i(t) \oplus( \mu_i^\top  X^{\star, \Adv} ) \in \R^3$. We define $\widetilde{X}_k =X_k  \oplus X^\star\oplus X^{\star, \Adv}$ analogously for SGD. Note that the results in this section hold for the soft label setting and hard label setting with $X^\star = 0$.
    
    Here $\mathrsfs{B}_i(t)$ is defined in \eqref{eq:ODE_resolvent_2} and is the deterministic equivalent of the covariance matrix of the Gaussian vector $(\WHSGD_t^\top a_i, (X^\star)^\top a_i)$ and $\mathrsfs{A}_{K_i}(t)$ is the deterministic equivalent of $X_{td}^\top K_i  X^{\star,\operatorname{adv}}$. Since $ \mu_1^\top  X^{\star, \Adv}$ and $ \mu_2^\top  X^{\star, \Adv}$ are constants, computing $\mathfrak{m}_i^\Adv(t) $ and $ \mathrsfs{B}^{\Adv}_i(t)$ amounts to computing the deterministic equivalent of $X_{td}^\top \mathscr{R}(z;\mathcal{K}) X^{\star,\operatorname{adv}}$. 

Given that the argument is extremely similar to the steps presented in Section~\ref{sect:doob:decomp} to Section~\ref{sect:conc:proof:any:stat}, we will only provide a brief sketch of the proof, highlighting the differences whenever they arise. We omit the steps of the proof for AdvHSGD and focus on SGD as the argument is extremely similar.  Define $\widehat{W}_k = W_k \oplus X^{\star, \Adv}$ for SGD and $\widehat{\HSGD}_t$ analogously for AdvHSGD. From the Doob decompositions presented in Section~\ref{sect:doob:decomp}, the doob decomposition of $ \varphi(\widehat{W}_{td}) =X_{td}^\top \mathscr{R}(z;\mathcal{K})  X^{\star, \Adv}$ is given by
\begin{align}
    \varphi(\widehat{W}_{td}) 
    &
    = \varphi(\widehat{W}_0) - \int_0^{t} \gamma_{sd} \ip{ \nabla \varphi(\widehat{W}_{sd}), \nabla \mathcal{R}_{\lambda}^{{\operatorname{adv}}}(X_{sd})} \, \dif s 
 +  \sum_{j=0}^{\lfloor td \rfloor-1} \Delta \mathcal{M}_{j}^{\text{Grad}}(\varphi)
            + \xi_{td}(\varphi). \label{eq:main:error_terms_integrated:align:mean}
    \end{align}
We redefine the $\sigma$-algebras as follows
\[
\begin{gathered}
\mathcal{G}_{k,i} \defas \sigma ( \{\hat{X}_j\}_{j=0}^k, \{r_i\}_{i=0}^k, \epsilon_{k+1}, \{I_j\}_{j=1}^k, I_{k+1} = i, X^{\star, \Adv})\quad
 \text{and}\quad \mathcal{F}_{k} \defas \sigma ( \{\hat{X}_j\}_{j=0}^k , X^{\star, \Adv}).
\end{gathered}
\]
Analogously to $Z$ defined in \eqref{eq:def:og:main:Z}, we introduce $Z^\Adv$ as follows
\begin{equation}\label{eq:def:Zadv}
    Z^\Adv(\widehat{W},z) \defas \widehat{W}^\top \mathscr{R}(z; \mathcal{K})\widehat{W} = \begin{bmatrix}
        Z(W,z) & (\star)_0\\
        (\star)_0^\top & (X^{\star, \Adv})^\top\mathscr{R}(z;\mathcal{K})(X^{\star, \Adv})
    \end{bmatrix}.
\end{equation}
with $A_{\mathscr{K}}(\widehat{W},z) \defas X^\top\mathscr{R}(z;\mathcal{K}) X^{\star, \Adv}$ such that $(\star)_0$ corresponds to the following term 
\[ 
    (\star)_0 \defas  \begin{bmatrix}
       A_{\mathscr{K}}(\widehat{W},z) & (X^{\star})^\top \mathscr{R}(z;\mathcal{K}) X^{\star, \Adv} & \mu_1^\top \mathscr{R}(z;\mathcal{K}) X^{\star, \Adv} & \mu_2^\top \mathscr{R}(z;\mathcal{K}) X^{\star, \Adv} 
    \end{bmatrix}^\top.
\]
We also define $\widehat{B}^\Adv_i(W)$ similarly to $\widehat{B}_i(W)$
\begin{equation}\label{eq:def:B:hat:adv}
    \widehat{\mathcal{B}}_i^{\Adv}(\widehat{W}) \defas \begin{bmatrix}
        \widehat{B}_i(W) &(\star)_1 \\
        (\star)_1^\top & \|X^{\star,\Adv}\|^2
     \end{bmatrix} \in \R^{5 \times 5},  
    \end{equation}
with $(\star)_1$ corresponding to the following term
\[
    (\star)_1 \defas  \begin{bmatrix}
        \oint z_i  A_{\mathscr{K}}(\widehat{W},z) \Dif z & (X^{\star})^\top  (X^{\star, \Adv}) & \mu_i^\top X^{\star, \Adv}  &0
    \end{bmatrix}^\top.
\]
Hence, using a similar argument as presented in Appendix~\ref{sect:DoobDecomp:HSGD}, the deterministic equivalent of $ A_{\mathscr{K}}(\widehat{W}_{td},z) = X_{td}^\top \mathscr{R}(z;\mathcal{K}) X^{\star, \Adv}$ which we denote $\mathcal{A}_{\mathcal{K}}(t,z)$ is defined in terms of the solution to the following integro-differential equation.
 \begin{mdframed}[style=exampledefault]
    \textbf{Integro-Differential Equation for $\mathcal{Z}^{\Adv}(t, z)$.} For any contour $\Gamma \subset \mathbb{C}$ enclosing the eigenvalues of $K_1$ and $K_2$, we have an expression for the derivative of $\mathcal{Z}^\Adv$ for $z\in \Gamma^2$:
    \begin{equation}\label{eq:ODE_resolvent_2:align}
        \dif \mathcal{Z}^\Adv(t, \cdot) = \mathscr{F}^\Adv(z, \mathcal{Z}^\Adv(t, \cdot)) \, \dif t,\quad \text{with initial condition }\mathcal{Z}^\Adv(0, z) = Z^\Adv(W_0, z),
        \end{equation}
        where
        \begin{equation}
            \begin{gathered}
                \mathcal{Z}^\Adv(t, z) \defas  \begin{bmatrix}
                    \mathcal{Z}(t, z) & (\star)_2 \\
                    (\star)_2 ^\top &  (X^{\star, \Adv})^\top\mathscr{R}(z;\mathcal{K})(X^{\star, \Adv})
            \end{bmatrix},
            \\
            (\star)_2 \defas  \begin{bmatrix}
                \mathcal{A}_{\mathscr{K}}(t,z) & (X^{\star})^\top \mathscr{R}(z;\mathcal{K}) X^{\star, \Adv} & \mu_1^\top \mathscr{R}(z;\mathcal{K}) X^{\star, \Adv} & \mu_2^\top \mathscr{R}(z;\mathcal{K}) X^{\star, \Adv} 
             \end{bmatrix}^\top,
            \\
            \text{and}\quad
            \mathscr{F}^\Adv(z, \mathcal{Z}^\Adv(t, \cdot)) \defas  \begin{bmatrix}
                \mathscr{F}(z, \mathcal{Z}(t, \cdot)) & \mathscr{F}_\star(z,\mathcal{Z}^\Adv(t, \cdot)) \\
                \mathscr{F}_\star(z, \mathcal{Z}^\Adv(t, \cdot))^\top & 0
            \end{bmatrix},
            \\
            \text{with}\quad \mathscr{F}_\star(z, \mathcal{Z}^\Adv(t, \cdot)) \defas \begin{bmatrix}
                \mathscr{F}_{\mathcal{A}_{\mathcal{K}}}(z, \mathcal{Z}^\Adv(t, \cdot)) & 0
             \end{bmatrix}^\top \in \R^5.
        \end{gathered}
        \end{equation}
The map $\mathscr{F}$ is defined in \eqref{eq:ODE_resolvent_2} and 
\begin{equation}\label{eq:det:equiv:align:Xstar:Adv}
    \dif \mathcal{A}_{\mathcal{K}}(t, \cdot) = \mathscr{F}_{\mathcal{A}_{\mathcal{K}}}(z, \mathcal{Z}^\Adv(t, \cdot)) \, \dif t,\quad \text{with initial condition}\quad  \mathcal{A}_{\mathcal{K}}(0, z) = A_{\mathscr{K}}(\widehat{W}_0, z),
        \end{equation}
where the map $\mathscr{F}_{\mathcal{A}_{\mathcal{K}}}$ is defined as follows
 \begin{equation}
    \begin{aligned}
        &\mathscr{F}_{\mathcal{A}_{\mathcal{K}}}(z, \mathcal{Z}^\Adv(t, \cdot))
        \defas 
        -2\gamma(t) \sum_{i=1}^2 p_i\bigg(\partial_{11}h_i( \widehat{\mathrsfs{B}}_i(t))  \left(z_i \mathcal{A}_{\mathcal{K}}(t,z)- \frac{1}{2\pi\operatorname{i}} \oint_{\Gamma}  \mathcal{A}_{\mathcal{K}}(t,z) \dif z_i\right)
        \\
        &
        + \partial_{12}h_i (\widehat{\mathrsfs{B}}_i(t)) \left(z_i \langle X^{\star} \otimes X^{\star, \Adv}, \mathscr{R}(z;\mathcal{K})\rangle  - \frac{1}{2\pi\operatorname{i}} \oint_{\Gamma}  \langle X^{\star} \otimes X^{\star, \Adv}, \mathscr{R}(z;\mathcal{K})\rangle  \dif z_i\right)
        \\
        &
        +\partial_{13} h_i(\widehat{\mathrsfs{B}}_i(t)) \langle \mu_i \otimes X^{\star, \Adv},  \mathscr{R}(z;\mathcal{K})\rangle
        + \left(\frac{\lambda}{2} + \partial_{44} h_i (\widehat{\mathrsfs{B}}_i(t))\right)\mathcal{A}_{\mathcal{K}}(t,z)\bigg)\dif t .
    \end{aligned}
 \end{equation}
Here $\widehat{\mathrsfs{B}}_i(t)$ is defined in \eqref{eq:ODE_resolvent_2} and the functions $h_i:\R^{4\times 4} \to \R$ and $\mathcal{I}_i:\R^{4\times 4} \to \R$ are defined in Assumptions~\ref{assumption:risk}~and~\ref{assumption:fisher} respectively. Finally, here $\widehat{W}_0 = \hat{X}_0 \oplus \mu \oplus X^{\star,\Adv} $ where $\hat{X}_0 = X_0 \oplus X^{\star}$ in the soft label setting and $\hat{X}_0 = X_0 \oplus 0$ in the hard label setting.
\end{mdframed}
 Similarly to \eqref{eq:S(t):def:average}, when $K_1$ and $K_2$ are finite-dimensional matrices, we can decompose \eqref{eq:det:equiv:align:Xstar:Adv} in terms of the eigenspaces of the covariance matrices $K_1$ and $K_2$. Define the pair of statistics
 \begin{equation}
     \begin{gathered}
     U_{j}(\widehat{W}_k) \defas X_k^\top \omega_j \omega_j^\top X^{\star, \Adv} = \frac{-1}{4\pi^2}\oint_{\Gamma_j^2} A_{\mathcal{K}}(\widehat{W}_k, z)\dif z_1 \dif z_2 \quad \text{and its deterministic equiv. }   \mathrsfs{U}_{j}(t),
 \end{gathered}
 \end{equation}
 where $\Gamma_i$ is a contour enclosing solely the eigenvalues $\lambda_j^{(1)}$ and $\lambda_j^{(2)}$. From the partial fractions decomposition of meromorphic functions, we have the representation 
 \begin{equation}
     \begin{gathered}
        \mathcal{A}_{\mathcal{K}}(t,z) = \sum_{j=1}^d \frac{\mathrsfs{U}_{j}(t)}{\prod_{\ell=1}^2(\lambda^{(\ell)}_j - z_\ell)}.
     \end{gathered}
 \end{equation}
From the Cauchy integral formula, we have $\mathrsfs{A}_{K_i}(t) = \oint z_i  \mathcal{A}_{\mathcal{K}}(t,z) \Dif z$ and summing over the eigenspaces of $K_1$ and $K_2$, for $i=1,2$ we obtain 
\begin{equation}\label{eq:def:A:K_i:adv}
\mathrsfs{A}_{K_i}(t) = \sum_{j=1} \lambda_j^{(i)} \mathrsfs{U}_{j}(t).
\end{equation}
Thus, we obtain the coupled ODEs
 \begin{equation}\label{eq:ODE:align_K(t):def}
     \begin{aligned}
         &
         \frac{\dif \mathrsfs{U}_{j}(t)}{\dif t}  = -2\gamma(t) \sum_{i=1}^2 p_i \bigg( \left(\lambda_j^{(i)} \partial_{11} h_i(\widehat{\mathrsfs{B}}_i(t)) + \left(\frac{\lambda}{2} + \partial_{44} h_i (\widehat{\mathrsfs{B}}_i(t))\right)\right) \mathrsfs{U}_{j}(t)  
         \\
          &\quad+\left(\lambda_j \partial_{12}h_i (\widehat{\mathrsfs{B}}_i(t)) \langle X^{\star }, \omega_j\rangle+ \partial_{13} h_i(\widehat{\mathrsfs{B}}_i(t)) \langle \mu_i ,\omega_j \rangle \right) \langle X^{\star, \Adv},  \omega_j\rangle 
          \bigg),
     \end{aligned}
 \end{equation}
 with initial condition $ \mathrsfs{U}_{j}(0) = X_0^\top \omega_j \omega_j^\top X^{\star, \Adv}$. Similarly to Remark~\ref{remark:equiv:sols:of:odes}, an almost identical argument to Lemma 2 in \citep{collinswoodfin2025Exact} shows that we can recover a solution to \eqref{eq:det:equiv:align:Xstar:Adv} using a solution to the system of coupled ODEs \eqref{eq:ODE:align_K(t):def}. The following Lemma is analogous to Lemma~\ref{lem:prop:normGamma:normW} and follows an argument almost identical to Lemma 3 in \citep{collinswoodfin2025Exact} and Lemma 5 in \citep{collinswoodfin2023hitting}. 

\begin{lemma}\label{lem:prop:normGamma:normW:align_K}
    Recall $\mathrsfs{N}(t) \defas \oint\Tr(\mathcal{Z}(t,z))\Dif z$. Then there exists a constant $C>0$ depending on $\|\mathcal{K}\|_{\operatorname{op}}$ and such that
    \[
        C \leq \frac{\|Z^\Adv(\widehat{W}_{td},\cdot)\|_{\Gamma}}{\|\widehat{W}_{td}\|^2},\frac{\| Z^\Adv(\widehat{\HSGD}_{t},\cdot)\|_{\Gamma}}{\|\widehat{\HSGD}_{t}\|^2}, \frac{\| \mathcal{Z}^\Adv(t,z\cdot\|_{\Gamma}}{\mathrsfs{N}(t) + \|X^{\star, \Adv}\|^2} \leq4 .
    \]
\end{lemma}
We introduce a stopping time $\psi_M$ analogous to $\tau_M$ in Definition~\ref{def:stoppingtime}.
\begin{definition}[Stopping time $\psi_M$]\label{def:stoppingtime:adv}
    For any $M>0$ and a continuous function $\mathrsfs{Z}^{\Adv}: [0, \infty] \times \mathbb{C}^2 \to \R^{5\times 5}$, define
    \[
        \psi_M(\mathrsfs{Z}^\Adv) \defas \inf\bigg\{t \geq 0 : \|\mathrsfs{Z}^\Adv(t,\cdot)\|_{\Gamma}> M \quad \text{or}\quad  \widehat{\mathcal{B}}_1^{\Adv}(t)  \not \in \mathcal{U}^\Adv \quad \text{or}  \quad \widehat{\mathcal{B}}_2^{\Adv}(t)  \not \in \mathcal{U}^\Adv\bigg\},
    \]
    where we define
    \begin{equation}\label{eq:def:Z(t,z):stability:single:adv}
        \mathrsfs{Z}^\Adv(t, \cdot) \defas \begin{bmatrix}
            \mathrsfs{Z}(t, \cdot)&(\star)_3\\
            (\star)_3^\top & (X^{\star, \Adv})^\top\mathscr{R}(z;\mathcal{K})(X^{\star, \Adv})
        \end{bmatrix}.
    \end{equation}
    with $(\star)_3$ corresponding to the following term
    \[ 
        (\star)_3 \defas  \begin{bmatrix}
            \mathrsfs{O}(t,\cdot) & (X^{\star})^\top \mathscr{R}(z;\mathcal{K}) X^{\star, \Adv} & \mu_1^\top \mathscr{R}(z;\mathcal{K}) X^{\star, \Adv} & \mu_2^\top \mathscr{R}(z;\mathcal{K}) X^{\star, \Adv} 
        \end{bmatrix}^\top
    \]  
    and
\begin{equation}
    \widehat{\mathcal{B}}_i^{\Adv}(t) \defas \begin{bmatrix}
        \widehat{\mathcal{B}}_i(t) &(\star)_4 \\
        (\star)_4^\top & \|X^{\star,\Adv}\|^2
     \end{bmatrix} \in \R^{5 \times 5},  
    \end{equation}
with 
\[
    (\star)_4 \defas  \begin{bmatrix}
        \oint z_i \mathrsfs{O}(t,z) \Dif z & (X^{\star})^\top  (X^{\star, \Adv}) & \mu_i^\top X^{\star, \Adv} & 0
    \end{bmatrix}^\top
\]
where $\mathrsfs{Z}(t,\cdot)$ is defined in Definition~\ref{def:stoppingtime}. Here $\mathrsfs{O}(t,\cdot)$ is either $A_{\mathscr{K}}(\widehat{W}_{td},\cdot)$, $A_{\mathscr{K}}(\widehat{\HSGD}_{t},\cdot)$ or $\mathcal{A}_{\mathcal{K}}(t,\cdot)$. Here $ \mathcal{U}^\Adv$ is defined in \eqref{eq:Uadv:set}.
\end{definition}
When clear from context, we omit the dependence $  \psi_M=\psi_M(\mathrsfs{Z}^\Adv) $. We are now ready to introduce the notion of approximate solutions to integro-differential equations to \eqref{eq:ODE_resolvent_2:align}.
\begin{definition}[$(\varepsilon, M, T )$-approximate solution to the integro-differential equation]\label{def:main:approx:sols:align}
    For $M, T, \varepsilon > 0$ and a continuous function $\mathrsfs{Z}^\Adv: [0, \infty] \times \mathbb{C}^2 \to \R^{5\times 5}$, we denote $\mathrsfs{Z}^\Adv$ an $(\varepsilon, M, T)$-approximate solution of \eqref{eq:ODE_resolvent_2:align} if
    \begin{equation}\label{eq:align:def:approxsols}
    \sup_{0\leq t\leq T\wedge \psi_M} \left\| \mathrsfs{Z}^\Adv(t, \cdot) -\mathrsfs{Z}^\Adv(0, \cdot) -\int_0^t\mathscr{F}^\Adv(z,  \mathrsfs{Z}^\Adv(s, \cdot)) \dif s\right\|_{\Gamma} \leq \varepsilon,
    \end{equation}
    where $\mathrsfs{Z}^\Adv(0, \cdot) = \widehat{W}_0^\top \mathscr{R}(\cdot; \mathcal{K}) \widehat{W}_0$.
\end{definition}
It is clear from the definition~\ref{def:main:approx:sols} and the definition of $ \mathcal{Z}^\Adv$ in \eqref{eq:ODE_resolvent_2:align} that $ \mathcal{Z}^\Adv$ is an approximate solution with $\varepsilon = 0$. In the next subsection, we will show show $\mathcal{Z}^\Adv(\widehat{W}_{td},z)$ and $Z^\Adv(\widehat{\HSGD}_{t},z)$ are $(\varepsilon, M, T)$-approximate solutions to \eqref{eq:align:def:approxsols}

As we've already shown in Propositions~\ref{prop:stability}~and~\ref{prop:stability:general} that approximate solutions to \eqref{eq:ODE_resolvent_2} satisfy stability properties, showing stability of approximate solutions as defined above in Definition~\ref{def:main:approx:sols:align} amounts to showing stability for the term $ \mathrsfs{O}(t,\cdot)$ in $\mathrsfs{Z}^\Adv(t,\cdot)$. As the proofs follow an almost identical argument to the proofs of Propositions~\ref{prop:stability}~and~\ref{prop:stability:general}, we state the results and omit them.
\begin{proposition}[Stability]\label{prop:stability:align}
    For all $(\varepsilon, M, T)$-approximate solutions $\mathrsfs{Z}^\Adv_1$ and $\mathrsfs{Z}^\Adv_2$, there exists a constant $C= C(\bar{\gamma},\lambda, M, \alpha, L(h), L(\mathcal{I}), \|\mathcal{K}\|_{\operatorname{op}},\max_{i=1,2} \|\mu_i\|, T)>0$ such that
    \[
        \sup_{0\leq t \leq T\wedge \widehat{\psi}_M} \|\mathrsfs{Z}^\Adv_1(t,\cdot) - \mathrsfs{Z}^\Adv_2(t,\cdot)\|_\Gamma \leq C\varepsilon,
    \]
    where $\widehat{\psi}_M =\min\{\psi_M(\mathrsfs{Z}^\Adv_1), \psi_M(\mathrsfs{Z}^\Adv_2)\}$.
\end{proposition}
Proposition~\ref{prop:stability} can be extended to any function $\varphi(X)$ satisfying Assumption~\ref{assumption:smooth:stats:Adv}.
\begin{proposition}\label{prop:stability:general:Adv}
    Suppose $\varphi: \R^{d}\to \R$ satisfies Assumption~\ref{assumption:smooth:stats:Adv}. Let $\mathrsfs{Z}^\Adv_1$ and $\mathrsfs{Z}^\Adv_2$ be $(\varepsilon, M, T)$-approximate solutions and for $i=1,2$ define
    \begin{equation}
        \mathrsfs{G}(\mathrsfs{Z}^\Adv_i(t,\cdot))  \defas  \begin{bmatrix}
            G_{i,1}  & G_{i,2}^\top  & 0  \\
            G_{i,2}  & \mu^\top \mu & 0 \\
            0 & 0 & G_{i,3}
        \end{bmatrix},
    \end{equation}
    where $ G_{i,2} = \oint \mathrsfs{M}_i(t,z)\Dif z \oplus (\mu^\top X^{\star,\Adv})\in \R^{2\times 3}$, $G_{i,3}= \oint\mathrsfs{S}_{i,11}(t,z)\Dif z$ and 
    \begin{equation}
        \begin{aligned}
            G_{i,1}  =  \begin{bmatrix}
                \oint q(z)\mathrsfs{S}_{i,11}(t,z)\Dif z & \oint q(z)\mathrsfs{S}_{i,12}(t,z)\Dif z & \oint q(z)\mathrsfs{O}_{i}(t,z)\Dif z \\
                \oint q(z)\mathrsfs{S}_{i,21}(t,z)\Dif z & X^{\star \top}q(K_1,K_2)X^{\star} & X^{\star\top}q(K_1,K_2)X^{\star, \Adv}  \\
                \oint q(z)\mathrsfs{O}_{i}(t,z)\Dif z &  X^{\star\top}q(K_1,K_2)X^{\star, \Adv} & X^{\star, \Adv \top}q(K_1,K_2)X^{\star, \Adv}
            \end{bmatrix}.
        \end{aligned}
    \end{equation}
     Then there exists a positive constant $C= C(\bar{\gamma},\lambda, M, \alpha, L(h), L(g), \|K\|_{\operatorname{op}})$ such that
    \[
        \sup_{0\leq t \leq T\wedge \widehat{\psi}_M} \left| g(\mathrsfs{G}(\mathrsfs{Z}^\Adv_1(t,\cdot)))  - g(\mathrsfs{G}(\mathrsfs{Z}^\Adv_2(t,\cdot)))\right| \leq C\varepsilon,
    \]
    where $\widehat{\psi}_M =\min\{\psi_M(\mathrsfs{Z}^\Adv_1), \psi_M(\mathrsfs{Z}^\Adv_2)\}$.
\end{proposition}
We require a net argument as in Lemma~\ref{lem:net_argument}. The proof follows an almost identical argument once again.
\begin{lemma}[Net argument] \label{lem:net_argument:align} Fix $T, M > 0, \xi > 0$. Suppose $\Gamma_{\xi}$ is a $d^{-\xi}$ mesh of $\;\Gamma$ with $|\Gamma_{\xi}| = C d^{\xi}$ for some $C > 0$. Let the function $Z^\Adv(t,z) =Z^\Adv(\widehat{W}_{td}, z)$ or $Z^\Adv(t,z)=Z^\Adv(\widehat{\HSGD}_t, z)$ satisfy
    \begin{equation}
    \label{eq:net_argument_gamma_xi:align}
    \sup_{0 \le t \le (T\wedge \tau_M )} \| Z^\Adv(t, \cdot) - Z^\Adv(0, \cdot) - \int_0^t \mathscr{F}^\Adv(\cdot, Z^\Adv(s, \cdot))\, \dif s \|_{\Gamma_\xi} \le \varepsilon
    \end{equation}
    with $\psi_M$ as in Definition~\ref{def:stoppingtime:adv}. Then $Z^\Adv$ is an $(\varepsilon + C  d^{-\xi}, M, T)$-approximate solution to the integro-differential equation \eqref{eq:ODE_resolvent_2:align}, that is, 
    \[
    \sup_{0 \le t \le (T\wedge \tau_M )} \|Z^\Adv(t, \cdot) - Z^\Adv(0, \cdot) - \int_0^t \mathscr{F}^\Adv(\cdot, Z^\Adv(s, \cdot) ) \, \dif s \|_{\Gamma} \le \varepsilon + C  d^{-\xi},
    \]
    where $C = C(\bar{\gamma}, \alpha, L(\mathcal{I}), M, \|\mathcal{K}\|_{\operatorname{op}}, L(h), \lambda, \max_{i=1,2}\|\mu_i\|, T)>0$ and 
    \end{lemma}
    Similarly to $\vartheta_M$, We introduce the stopping times for fixed $M>0$: 
\begin{equation}
    \begin{aligned} \label{eq:main:stopping_time_1:adv}
        \vartheta_{M}^\Adv 
        &
        \defas \inf \{ t \ge 0 \, : \, \|\widehat{W}_{td} \|^2 > M \quad \text{or} \quad \widehat{B}^{\Adv}_1(\widehat{W}_{td}) \not \in \mathcal{U}^\Adv\quad \text{or}\quad \widehat{B}^\Adv_2(\widehat{W}_{td}) \not \in \mathcal{U}^\Adv \}\\
        \text{or} \quad  \vartheta_{M}^\Adv 
        &
        \defas \inf \{t \ge 0 \, : \, 
        \|\widehat{\HSGD}_t\|^2  > M  \quad \text{or}\quad \widehat{B}_1^\Adv(\widehat{\HSGD}_t) \not \in \mathcal{U}^\Adv\quad \text{or}\quad \widehat{B}_2^\Adv(\widehat{\HSGD}_t) \not \in \mathcal{U}^\Adv \}. 
    \end{aligned}
    \end{equation}
     By Lemma~\ref{lem:prop:normGamma:normW:align_K}, there exists constants $c,C>0$ such that $\vartheta^\Adv_{c \cdot M} \le \psi_{M} \le  \vartheta^\Adv_{C \cdot M}$ implying the stopping times $\psi_M$ and $\vartheta^\Adv_M$ are interchangeable. 
    
    The proof of the following Theorem~\ref{thm:main:S(W):adv} follows basically the same argument as the proof of Theorem~\ref{thm:main:S(W):oneST}. Indeed, Propositions~\ref{prop:gradmartingale} and ~\ref{prop:HSGD_Martingale_Bound} were proven for an arbitrary quadratic statistic of SGD and AdvHSGD respectively. Hence, using the same arguments as the proof of Propositions~\ref{prop:approx:sols:SGD}~and~\ref{prop:approx:sols:HSGD}, it follows that for fixed $T, M> 0$ and $0 < \varepsilon < \frac{1}{2}$, $Z^\Adv(\widehat{W}_{td }, \cdot)$ and $ Z^\Adv(\widehat{\HSGD}_t,\cdot)$ are $(Cd^{-\varepsilon}, M, T)$ approximate solutions. Finally, the theorem follows from the stability proven in Proposition~\ref{prop:stability:align}.
    \begin{theorem}\label{thm:main:S(W):adv}
        There exists an $\varepsilon \in (0,\frac{1}{2})$ such that for any $T, M >0$ and $d$ sufficiently large, with overwhelming probability
        \begin{equation}
            \begin{gathered}
              \sup_{0 \leq t \leq T \wedge \widehat{\psi}_{M}(Z^\Adv(W),\mathcal{Z}^{\Adv})} \! \! \! \! \!\! \! \! \! \| 
              Z^\Adv(\widehat{W}_{td }, \cdot)
              - 
              \mathcal{Z}^\Adv(t,\cdot)
              \|_{\Gamma} \leq Cd^{-\varepsilon}, 
              \\
               \sup_{0 \leq t \leq T \wedge \widehat{\psi}_{M}(Z^\Adv(\HSGD), \mathcal{Z}^{\Adv})} \| 
               Z^\Adv(\HSGD_t, \cdot)
              - 
              \mathcal{Z}^\Adv(t,\cdot)
              \|_{\Gamma} \leq Cd^{-\varepsilon},
              \\
                 \text{and} \qquad \sup_{0 \leq t \leq T \wedge \widehat{\psi}_M(Z^\Adv(W), Z^\Adv(\HSGD) )} \! \! \! \! \!\! \! \! \!  \| 
                 Z^\Adv(\widehat{W}_{ td }, \cdot)
              - 
              Z^\Adv(\widehat{\HSGD}_t,\cdot)
              \|_{\Gamma} \leq Cd^{-\varepsilon},
              \end{gathered}
            \end{equation}
          where $C>0$ is a constant independent of $d$, where $\mathcal{Z}^\Adv(t,z)$ is the deterministic solution to the integro-differential equation \eqref{eq:det:equiv:align:Xstar:Adv} and
          \[
          \widehat{\psi}_{M}(\mathrsfs{Z}^\Adv_1, \mathrsfs{Z}^\Adv_2) = \min\{\psi_M(\mathrsfs{Z}^\Adv_1), \psi_M(\mathrsfs{Z}^\Adv_2) \}.
          \]
    \end{theorem}
    As in Section~\ref{sect:remove:stopping:time}, we now proceed with the removal of the stopping time $\widehat{\psi}_M$. We now introduce two stopping times analogous to \eqref{eq:def:tau_M,omega:stoppingtime} and\eqref{eq:def:nu:stoppingtime}. Let $\mathrsfs{Z}^\Adv(t,\cdot)$ be an $(\varepsilon, M, T)$-approximate solution then we define the stopping time for $M, \omega > 0$
    \begin{equation}
        \begin{aligned}
        \psi_{M,\omega} \defas \inf\bigg\{t\geq 0: \|\mathrsfs{Z}^\Adv(t,\cdot)\|_{\Gamma} &> M \quad \text{or}\quad \inf_{V \in (\mathcal{U}^{\Adv})^c}\|\widehat{\mathcal{B}}^\Adv_1(t) - V\| \leq \omega 
        \\
        &\text{or}\quad \inf_{V \in (\mathcal{U}^{\Adv})^c}\|\widehat{\mathcal{B}}^\Adv_2(t) - V\| \leq \omega\bigg\},
        \end{aligned}
    \end{equation}
    where $\widehat{\mathcal{B}}^\Adv_i(t)$ is defined in Definition~\ref{def:stoppingtime:adv} and $(\mathcal{U}^\Adv)^c$ is the complement of $\mathcal{U}^\Adv$ as defined in Assumption~\ref{def:stoppingtime:adv}. It is clear that $ \psi_{M,0} =  \psi_{M}$ as defined in Definition~\ref{def:stoppingtime:adv}. We also introduce the stopping time
    \begin{equation}\label{eq:def:nu:stoppingtime:adv}
        \begin{aligned}
            &\Theta_{M,\omega}^\Adv (\mathrsfs{Z}^\Adv_1, \mathrsfs{Z}^\Adv_2) \defas \max\left\{ \inf\{t\geq 0: \|\mathrsfs{Z}^\Adv_j(t,\cdot)\|_{\Gamma} > M\} \quad \text{for}\quad j=1,2 \right\}\\
            &\quad\wedge \max\left\{\inf\{t\geq 0:\inf_{V \in (\mathcal{U}^\Adv)^c}\|\widehat{\mathcal{B}}^\Adv_i(\mathrsfs{X}^\Adv_j(t,\cdot)) - V\| \leq \omega \}\quad \text{for}\quad i=1,2\quad j=1,2\right\}
        \end{aligned}
    \end{equation}
    where we write $\widehat{\mathcal{B}}^\Adv_i(\mathrsfs{X}^\Adv_j(t,\cdot))$ to emphasize the dependence on $\mathrsfs{X}^\Adv_j$. Note that $\Theta^\Adv_M = \Theta_{M, 0}^\Adv$ in Section~\ref{sect:mainresult:SGD}. We omit the proof of the following theorem as it uses the same argument as the proof of Theorem~\ref{thm:main:S(W):oneST} but applied to $\widehat{\mathcal{B}}^\Adv_i(\mathrsfs{X}^\Adv_j)$ instead.
    \begin{theorem}\label{thm:main:S(W):oneST:Adv}
        There exists an $\varepsilon \in (0, \frac{1}{2})$ such that for any $T, M, \omega >0$ and $d$ sufficiently large, with overwhelming probability
        \begin{equation}
            \begin{gathered}
              \sup_{0 \leq t \leq T \wedge \Theta^\Adv_{M,\omega}(Z^\Adv(W, \cdot), \mathcal{Z}^\Adv)} \! \! \! \! \!\! \! \! \! \| 
              Z^\Adv(\widehat{W}_{td }, \cdot)
              - 
              \mathcal{Z}^\Adv(t,\cdot)
              \|_{\Gamma} \leq Cd^{-\varepsilon}, 
              \\
               \sup_{0 \leq t \leq T \wedge \Theta^\Adv_{M,\omega}(Z^\Adv(\HSGD,\cdot), \mathcal{Z}^\Adv)} \| 
              Z^\Adv(\widehat{\HSGD}_t, \cdot)
              - 
              \mathcal{Z}^\Adv(t,\cdot)
              \|_{\Gamma} \leq Cd^{-\varepsilon},
              \\
                 \text{and} \qquad \sup_{0 \leq t \leq T \wedge \Theta_{M,\omega}^{\Adv}(Z^\Adv(W,\cdot), Z^\Adv(\HSGD,\cdot) )} \! \! \! \! \!\! \! \! \!  \| 
              Z^\Adv(\widehat{W}_{ td }, \cdot)
              - 
              Z^\Adv(\widehat{\HSGD}_t,\cdot)
              \|_{\Gamma} \leq Cd^{-\varepsilon},
              \end{gathered}
            \end{equation}
            where $C>0$ is a constant independent of $d$, where $\mathcal{Z}^\Adv(t,z)$ is the deterministic solution to the integro-differential equation \eqref{eq:ODE_resolvent_2:align}.
    \end{theorem}
    Hence, it follows from Theorem~\ref{thm:main:S(W):oneST:Adv} that SGD and AdvHSGD concentrate around $\mathcal{Z}^\Adv(t,z)$ provided that AdvHSGD or $\mathcal{Z}^\Adv(t,z)$ stay bounded. Analogously to Corollary~\ref{cor:bounded_iterates}, we formalize this in the following corollary.
    \begin{corollary}
        [Bounded $\mathrsfs{N}(t)$ and concentration] \label{cor:bounded_iterates:adv} Define $ \widehat{\mathrsfs{B}}_i^{\Adv}(t)$ similarly to $ \widehat{\mathrsfs{B}}_i(t)$:
\begin{equation}
    \widehat{\mathrsfs{B}}_i^{\Adv}(t) \defas \begin{bmatrix}
        \widehat{\mathrsfs{B}}_i(t) &(\star) \\
        (\star)^\top & \|X^{\star,\Adv}\|^2
     \end{bmatrix} \in \R^{6 \times 6},  
    \end{equation}
with $(\star)$ corresponding to the following term
\[
    (\star) \defas  \begin{bmatrix}
        \oint z_i  \mathcal{A}_{\mathscr{K}}(t,z) \Dif z & (X^{\star})^\top  (X^{\star, \Adv}) & \mu_1^\top X^{\star, \Adv} & \mu_2^\top X^{\star, \Adv} &0
    \end{bmatrix}^\top.
\]
Suppose for a fixed $T,  \omega> 0$ that
    \begin{equation} 
     \sup_{0 \le t \le T} \mathrsfs{N}(t) + \|X^{\star,\Adv}\|^2\le M  \quad \text{and for $i=1,2$} \quad \inf_{0 \le t \le T} \inf_{V \in \mathcal{U}^c} \|\widehat{\mathrsfs{B}}^\Adv_i(t)-V\| > \eta \quad \text{hold w.o.p.}
    \end{equation}
    where $M>0$ is independent of $d$. Then there exists $\varepsilon \in (0, \frac{1}{2})$ so that for $d$ sufficiently large, with overwhelming probability, 
      \begin{equation}
      \begin{gathered} 
        \sup_{0 \leq t \leq T } \| 
        Z^\Adv(\widehat{\HSGD}_t, \cdot)
        - 
        \mathcal{Z}^\Adv(t,\cdot)
        \|_{\Gamma} \leq Cd^{-\varepsilon} \quad 
           \text{and} \quad \sup_{0 \leq t \leq T} \| 
        Z^\Adv(\widehat{W}_{td}, \cdot)
        - 
        \mathcal{Z}^\Adv(t,\cdot)
        \|_{\Gamma} \leq Cd^{-\varepsilon},
        \end{gathered}
      \end{equation}
      where $C>0$ is a constant independent of $d$. Hence, by the triangle inequality, it follows that
      \begin{equation} 
      \sup_{0 \le t \le T} \| Z^\Adv(\widehat{W}_{td}, \cdot) - Z^\Adv(\widehat{\HSGD}_t,\cdot)\|_{\Gamma} \le 2Cd^{-\varepsilon}. 
      \end{equation}
    \end{corollary}
    Finally, using the same argument as the proof of Theorem~\ref{thm:concentration_statistic} which we also omit, we obtain a general concentration result for statistics satisfying Assumption~\ref{assumption:smooth:stats:Adv}, As presented in the proof of Proposition~\ref{prop:stability:general:Adv}, recall the deterministic equivalent of $\varphi(\WHSGD_t)$ and $\varphi(X_{td})$ for $\varphi$ satisfying Assumption~\ref{assumption:smooth:stats:Adv} is defined as
    \begin{equation}\label{eq:det:equiv:General:adv}
        \phi^\Adv(t) =g\left( \begin{bmatrix}
            G_{1}  & G_{2}^\top & 0  \\
            G_{2} & \mu^\top \mu & 0 \\
            0 & 0 & G_{3}
        \end{bmatrix}\right),
    \end{equation}
    where $ G_{2} = \oint \mathcal{M}(t,z)\Dif z \oplus (\mu^\top X^{\star,\Adv})\in \R^{2\times 3}$, $G_{3}= \oint\mathrsfs{S}_{11}(t,z)\Dif z$ and 
    \begin{equation}
        \begin{aligned}
            G_{i,1}  =  \begin{bmatrix}
                \oint q(z)\mathcal{S}_{11}(t,z)\Dif z & \oint q(z)\mathcal{S}_{12}(t,z)\Dif z & \oint q(z)\mathcal{A}_{\mathscr{K}}(t,z)\Dif z \\
                \oint q(z)\mathcal{S}_{21}(t,z)\Dif z & X^{\star \top}q(K_1,K_2)X^{\star} & X^{\star\top}q(K_1,K_2)X^{\star, \Adv}  \\
                \oint q(z)\mathcal{A}_{\mathscr{K}}(t,z)\Dif z &  X^{\star\top}q(K_1,K_2)X^{\star, \Adv} & X^{\star, \Adv \top}q(K_1,K_2)X^{\star, \Adv}
            \end{bmatrix}.
        \end{aligned}
    \end{equation}
    \begin{theorem}[Concentration of any statistic] \label{thm:concentration_statistic:adv} Suppose that the assumptions stated in this framework hold. Then there exists an $\varepsilon \in (0, \frac{1}{2})$ such that for any $T, M >0$ and $d$ sufficiently large, with overwhelming probability 
      \begin{equation}
      \begin{gathered}
        \sup_{0 \leq t \leq T \wedge \Theta^\Adv_{M}(Z^\Adv(\widehat{W},\cdot), \mathcal{Z}^\Adv)} | 
        \varphi(W_{ td })
        - 
        \phi^\Adv(t)
        |\leq Cd^{-\varepsilon}, \\
        \sup_{0 \leq t \leq T \wedge  \Theta^\Adv_{M}(Z^\Adv(\widehat{\HSGD},\cdot), \mathcal{Z}^\Adv) }| 
        \varphi(\HSGD_t)
        - 
        \phi^\Adv(t)
        | \leq Cd^{-\varepsilon},
        \\
           \text{and} \qquad \sup_{0 \leq t \leq T \wedge  \Theta^\Adv_{M}(Z^\Adv(\widehat{W},\cdot), Z^\Adv(\widehat{\HSGD},\cdot))} | 
        \varphi(W_{td})
        - 
        \varphi(\HSGD_t)
        | \leq Cd^{-\varepsilon},
        \end{gathered}
      \end{equation}
    where $C>0$ is a constant independent of $d$, where $\phi^\Adv$ is defined in \eqref{eq:det:equiv:General:adv} and where the stopping time $\Theta^\Adv_{M}(\mathrsfs{Z}^\Adv_1, \mathrsfs{Z}^\Adv_2)\equiv\Theta^\Adv_{M, 0}(\mathrsfs{Z}^\Adv_1, \mathrsfs{Z}^\Adv_2)$ is defined in \eqref{eq:def:nu:stoppingtime}. 
    \end{theorem}
    \begin{corollary}\label{cor:statistic_bounded:adv}
        Suppose that the assumptions stated in this framework hold and the assumptions in Corollary~\ref{cor:bounded_iterates:adv} hold. Then there exists an $\varepsilon \in (0, \frac{1}{2})$ such that for any $T>0$ and sufficiently large $d$, with overwhelming probability, 
      \begin{equation}
      \begin{gathered} 
        \sup_{0 \leq t \leq T } | 
        \varphi(\WHSGD_t)
        - 
        \phi^\Adv(t)
        | \leq Cd^{-\varepsilon} \quad 
           \text{and} \qquad \sup_{0 \leq t \leq T} | 
        \varphi(X_{td})
        - 
        \phi^\Adv(t)
        | \leq Cd^{-\varepsilon}.
        \end{gathered}
      \end{equation}
      Hence, by the triangle inequality, it follows that
      \begin{equation} 
      \sup_{0 \le t \le T} | \varphi(X_{td}) - \varphi(\WHSGD_t) | \le 2Cd^{-\varepsilon}. 
      \end{equation}
    \end{corollary}

    \section*{Acknowledgments}
We would like to thank Lucas Benigni for his careful proofreading and helpful feedback throughout the writing of this paper. We would also like to thank Courtney Paquette, Elliot Paquette and Gauthier Gidel for fruitful discussions.
\bibliographystyle{plainnat}
\bibliography{references}
\clearpage
\appendix
\section{Preliminaries for the Proofs}\label{sect:prelim:proofs}
The results in this section build upon \citep{collinswoodfin2023hitting,collinswoodfin2025Exact} in order to fit the current framework. We refer the reader to Section 3 of \citep{collinswoodfin2023hitting} for a thorough overview to contractions over abtract vector spaces and tensor differentiation. To facilitate understanding of the proofs and streamline notation, we will include results from these papers and use similar notation. Recall we defined
$r_{k, i} = \hat{X}_k^\top a_{k+1, i}$ such that 
\begin{equation}
    \begin{gathered}
        g_{k, i} = r_{k,i} + \begin{bmatrix}
            \delta s_{k, i}\|X_k\| & 0 
        \end{bmatrix}^\top
    ,\\
    s_{k, i} = \operatorname{argmax}_{|s| \leq  1}f_{i}( X_k^\top a_{k+1, i}+ \delta s \|X_k\|; \epsilon_{k+1}),
\end{gathered}
\end{equation}
where $a_{k, I_{k}} = \sqrt{K_{I_{k}}}v_k + \mu_{I_{k}}$ for $v_k \sim \mathcal{N}(0, \Id_d)$.

\subsection{Derivative of Special Statistics}\label{sect:deriv:stats}

In this section, we compute the derivative the adversarial risk function $\mathcal{R}^{\operatorname{adv}}$. From Assumption~\ref{assumption:pseudo_lipschitz}, the functions $f_i$ $\alpha$-pseudo-Lipschitz for $i=1,2$. Hence, on any bounded set, the functions $f_i$ and  and $\Psi_i\, :\R^d \to \mathbb{R}$ defined in \eqref{eq:nlgc} are Lipschitz so differentiable a.e. by Radamacher's theorem. Thus, the derivatives $f_i'$ and $\nabla_X \Psi_i$, exist a.e. In section~\ref{sect:conc:pseudoLispchitz}, we derive growth bounds on $f'_i$. We will now to compute the derivative of the risk $\mathcal{R}^{\operatorname{adv}}(X)$.  
\begin{lemma}[Derivatives of the statistic, $\mathcal{R}^{\operatorname{adv}}$] \label{lem:derivative_adv_risk}
Suppose the risk is $\mathcal{R}^{\operatorname{adv}}(X)$ satisfies Assumption~\ref{assumption:risk}. Then, one has
\begin{equation*}
    \begin{aligned}
        \nabla  \mathcal{R}^{\operatorname{adv}}(X) &= 2\sum_{i=1}^2 p_i \bigg(\partial_{11} h_i(\widehat{B}_i(W)) \cdot K_iX + \partial_{12} h_i(\widehat{B}_i(W)) \cdot K_iX^{\star}
        \\
        &
        \qquad
        +
         \partial_{13} h_i(\widehat{B}_i(W)) \cdot  \mu_i  +  \partial_{44} h_i(\widehat{B}_i(W))\cdot X\bigg)
    \end{aligned}
\end{equation*}
where $h:  \R^{4\times 4}  \to \R$ with gradient
\begin{equation*}
    \nabla h_i 
=
\left[ \begin{array}{c|c|c|c} 
   \partial_{11} h_i &\partial_{12} h_i & \partial_{13} h_i  & 0 \\
    \hline
   \partial_{21} h_i &\partial_{22} h_i &  \partial_{23} h_i & 0\\
    \hline
    \partial_{31}h_i & \partial_{32} h_i  &\partial_{33} h_i & 0 \\
    \hline
    0&0 & 0 &\partial_{44} h_i
    \end{array} \right] \in \R^{4\times 4}.
 \end{equation*}
\end{lemma}
\begin{proof}
We first introduce some notation on tensor contractions and tensor differentiation in order to prove a more general result as the proof of this more general result is easier to follow. See Section 3 of \citep{collinswoodfin2023hitting} for more details.
 
 Let $\mathcal{A} \cong \R^{d}$, $\mathcal{O} \cong \R^{\ell}$, $\mathcal{T} \cong \R^{\ell^\star}$ be the finite dimensional Hilbert spaces such that $\mathcal{A} \otimes \mathcal{O} \cong \R^d \otimes \R^\ell \cong \R^{d\ell}$. Let $\oplus$ denotes the concatenation operation of all vectors from the vector space $\mathcal{O} $ with the vectors from the vector space $ \mathcal{T}$. Then define the space $O^{+} = \mathcal{O} \oplus \mathcal{T}$ such that $\mathcal{A} \otimes O^{+} \cong \R^d\otimes (\R^\ell \oplus \R^{\ell^\star})$. We will define the inner product $\langle \cdot, \cdot \rangle$ on each of these spaces as follows. Each space $\mathcal{A}$, $\mathcal{O}$ and $\mathcal{T}$ carries their own inner product. Consequently, the product spaces such as $\mathcal{A}\otimes \mathcal{O}$ have a natural inner product defined as follows
 \begin{equation}
    \langle a\otimes c, b\otimes d\rangle_{\mathcal{A} \otimes \mathcal{O}} = \langle a, b\rangle_{\mathcal{A}} \langle c, d\rangle_{\mathcal{O}},
 \end{equation}
 for simple tensors $a,b\in \mathcal{A}$ and $c,d\in \mathcal{O}$. This is extended to the full space through bilinearity. It is important to note we will always sort the simple tensors in the following order: $\mathcal{A}$ first and then $\mathcal{O}$, if $\mathcal{O}$ is involved in the contraction. We also note if it is crucial to preserve ordering within each space. For example, let $a_i \in \mathcal{A}$ for $i=1,2$ and $o_i \in \mathcal{O}$ for $i=1,2,3$ then we preserve ordering within each space as follows
 \[
    o_1 \otimes a_1 \otimes o_2\otimes a_2 \otimes o_3 \cong a_1 \otimes a_2 \otimes o_1 \otimes o_2 \otimes o_3.
 \]
 However, the following is invalid as it doesn't preserve ordering within $\mathcal{O}$
 \[
    o_1 \otimes a_1 \otimes o_2\otimes a_2 \otimes o_3 \cong a_1 \otimes a_2 \otimes o_2 \otimes o_1 \otimes o_3.
 \]
Before computing the derivatives using tensor notation, we provide a warm up from \citep{collinswoodfin2023hitting} to ease into the notation. Let $X\in \mathcal{A} \otimes \mathcal{O} $, $X^\star \in \mathcal{A} \otimes \mathcal{T}$, $a\in \mathcal{A}$ such that $K \in \mathcal{A}^{\otimes 2}$. We introduce the identity mapping 
\[
 \Id_{\mathcal{A} \otimes \mathcal{O}}: \mathcal{A} \otimes \mathcal{O} \to \mathcal{A} \otimes \mathcal{O} \quad \text{such that} \quad \Id(X) = X.
\]
Now, the derivative of the map $X \mapsto X$ is the linear map $D(X\mapsto X): \mathcal{A} \otimes \mathcal{O} \to \mathcal{L}(\mathcal{A} \otimes \mathcal{O} , \mathcal{A} \otimes \mathcal{O} )$ and is given by the identity mapping
\[
 DX \cong \Id_{\mathcal{A} \otimes \mathcal{O}}.
\]
Similarly, we obtain for $D(\langle X, a\rangle_{\mathcal{A}} ) \in \mathcal{L}(\mathcal{A} \otimes \mathcal{O}, \mathcal{O})$ by linearity
\[
 D(\langle X, a\rangle_{\mathcal{A}}) = \langle DX, a\rangle_{\mathcal{A}} \cong \langle  \Id_{\mathcal{A} \otimes \mathcal{O}}, a\rangle_{\mathcal{A}} \in \mathcal{L}(\mathcal{A} \otimes \mathcal{O}, \mathcal{O}).
\]
Choosing an orthonormal basis $\{e_{\alpha} \otimes f_{o}\}$ of $\mathcal{A} \otimes \mathcal{O}$ then 
\[
\begin{aligned}
    D(\langle X, a\rangle_{\mathcal{A}})&\cong \langle  \Id_{\mathcal{A} \otimes \mathcal{O}}, a\rangle_{\mathcal{A}}\\
    &\cong \sum_{\alpha, o} \langle e_{\alpha} \otimes f_{o}, a\rangle_{\mathcal{A}} \otimes e_{\alpha} \otimes f_{o}\\
    &= \sum_{\alpha, o} \langle e_{\alpha} , a\rangle_{\mathcal{A}}  f_{o}\otimes e_{\alpha} \otimes f_{o}\\
    &\cong \sum_{\alpha, o} \langle e_{\alpha} , a\rangle_{\mathcal{A}}  e_{\alpha}\otimes f_{o}\otimes  f_{o}\\
    &= \sum_{o} a\otimes f_{o}\otimes  f_{o}\\
    &\cong a \otimes \Id_{\mathcal{O}}.
\end{aligned}
\]
This representation corresponds to the tensor definition of the derivative. We can explicitly relate this definition to the operator definition as follows. For $H \in \mathcal{A} \otimes \mathcal{O}$ then the directional derivative
\[
    D(\langle X, a\rangle _{\mathcal{A}})[H] = \lim_{t \downarrow 0} \frac{\langle a, X+ tH\rangle_{\mathcal{A}} - \langle a, X\rangle_{\mathcal{A}}}{t} = \langle a, H\rangle_{\mathcal{A}} = \langle a \otimes \Id_{O}, H\rangle_{\mathcal{A} \otimes \mathcal{O}}. 
\]
 We are now ready to derive the general result. Recall $\hat{X} = X \oplus X^\star \in  \mathcal{A} \otimes \mathcal{O}^+$ and $W = \hat{X} \oplus \mu \in \mathcal{A} \otimes (\mathcal{O}^+ \oplus \mathcal{O} \oplus \mathcal{O})$.
From the definitions \eqref{eq:stats:BBcheck}, we have 
 \begin{equation}
    \begin{aligned}
        \check{B}_i(W)
        &=\begin{bmatrix}
        \ip{\hat{X}\otimes \hat{X}, K_i }_{\mathcal{A}^{\otimes 2}} & \langle \hat{X}, \mu_i \rangle_{\mathcal{A}} & 0_{\mathcal{O}^{+} \otimes  \mathcal{O}}\\
        \langle \mu_i, \hat{X} \rangle_{\mathcal{A}}  & \langle \mu_i, \mu_i \rangle_{\mathcal{A}}  & 0_{\mathcal{O}^{\otimes 2}} \\
        0_{ \mathcal{O} \otimes  \mathcal{O}^{+} } &  0_{\mathcal{O}^{\otimes 2}}& \ip{X, X}_{\mathcal{A}}
    \end{bmatrix}
    \cong  
    \left [\begin{array}{c|c|c} 
       ( \mathcal{O}^+ )^{ \otimes 2} &  \mathcal{O}^+  \otimes  \mathcal{O} &  \mathcal{O}^+  \otimes  \mathcal{O}\\
        \hline
        \mathcal{O} \otimes  \mathcal{O}^+ &  \mathcal{O}^{\otimes 2} & \mathcal{O}^{\otimes 2} \\
        \hline
        \mathcal{O} \otimes  \mathcal{O}^+ & \mathcal{O}^{\otimes 2}  & \mathcal{O}^{\otimes 2} 
        \end{array} \right ]\nonumber.
    \end{aligned}
 \end{equation} 
 Let $h_i:  (\mathcal{O}^{+} \oplus \mathcal{O}\oplus \mathcal{O})^{\otimes 2}  \to \R$ such that  
 \begin{equation*}
    \nabla h_i 
\cong
\left[ \begin{array}{c|c|c|c} 
   \nabla_{11}   h_i & \nabla_{12}  h_i  & \nabla_{13}h_i& 0_{\mathcal{O} \otimes \mathcal{O}} \\
    \hline
   \nabla_{21} h_i &\nabla_{22} h_i& \nabla_{23} h_i & 0_{\mathcal{T} \otimes \mathcal{O}}  \\
   \hline
   \nabla_{31} h_i & \nabla_{32} h_i & \nabla_{33} h_i & 0_{\mathcal{O} \otimes \mathcal{O}} \\
    \hline
    0_{\mathcal{O} \otimes \mathcal{O}}  & 0_{\mathcal{O} \otimes \mathcal{T}}  & 0_{\mathcal{O} \otimes \mathcal{O}} &\nabla_{44} h_i
    \end{array} \right ] \cong  \left[ \begin{array}{c|c|c|c} 
        \mathcal{O} \otimes \mathcal{O} & \mathcal{O} \otimes \mathcal{T} &  \mathcal{O} \otimes \mathcal{O} & \mathcal{O} \otimes \mathcal{O}
        \\
        \hline
        \mathcal{T} \otimes  \mathcal{O} & \mathcal{T} \otimes \mathcal{T} &  \mathcal{T} \otimes \mathcal{O} & \mathcal{T} \otimes \mathcal{O} \\
        \hline
        \mathcal{O} \otimes \mathcal{O}  &  \mathcal{O} \otimes \mathcal{T}  &  \mathcal{O} \otimes \mathcal{O} & \mathcal{O} \otimes \mathcal{O} \\
        \hline
        \mathcal{O} \otimes \mathcal{O} & \mathcal{O} \otimes \mathcal{T} & \mathcal{O} \otimes \mathcal{O} & \mathcal{O} \otimes \mathcal{O}
        \end{array} \right ].
 \end{equation*}
 Here $0_{\mathcal{O} \otimes \mathcal{O}}$ represents the zero tensor of the space $\mathcal{O} \otimes \mathcal{O}$ for example. Recall from \eqref{eq:nlgc} the risk $\mathcal{R}^{\operatorname{adv}}$ decomposes as follows
 \begin{equation}
    \mathcal{R}^{\operatorname{adv}}(X) = p_1 \cdot \EE_{v,\epsilon} [\Psi_1(X; \sqrt{K_1}v +\mu_1, \epsilon) ] + p_2 \EE_{v,\epsilon} [\Psi_2(X; \sqrt{K_2}v +\mu_2, \epsilon)  ].
\end{equation}
From Assumption~\ref{assumption:risk}, we have $\EE_{v,\epsilon} [\Psi_i(X; \sqrt{K_i}v +\mu_i, \epsilon) ]= h_i(\widehat{B}_i(W))$. Hence, computing the derivative of $\mathcal{R}^{\operatorname{adv}}$ amount to applying the chain rule to the functions $h_i$. First, the derivative of the mapping $X \mapsto \hat{X} = X \oplus X^{\star}$ is
\[
\Dif(\hat{X}) = \Dif(X \oplus X^{\star} ) \cong \text{Id}_{\mathcal{A} \otimes \mathcal{O}} \oplus 0_{\mathcal{A} \otimes \mathcal{T}}.
\]

Before applying the chain rule, we need to compute the derivative of the inside function $\Dif (X \mapsto \widehat{B}_i(W)\in \mathcal{L}(\mathcal{A} \otimes \mathcal{O}, (\mathcal{O}^+ \oplus \mathcal{O})^{\otimes 2})$. First, this requires the computation of the derivative of the functions 
\[
\begin{aligned}
    \Dif (X \mapsto \ip{\hat{X} \otimes \hat{X}, K_i}_{\mathcal{A}^{\otimes 2}}) \in \mathcal{L}(\mathcal{A} \otimes \mathcal{O}, (\mathcal{O}^+)^{\otimes 2}),\\
    \Dif (X \mapsto \ip{\mu_i, \hat{X}}_{\mathcal{A}}) \in \mathcal{L}(\mathcal{A} \otimes \mathcal{O}, \mathcal{O} \otimes \mathcal{O}^+),\\
    \Dif (X \mapsto \ip{\hat{X}, \mu_i}_{\mathcal{A}}) \in \mathcal{L}(\mathcal{A} \otimes \mathcal{O}, \mathcal{O}^+ \otimes \mathcal{O}),\\
    \Dif (X \mapsto \ip{X, X}_{\mathcal{A}}) \in \mathcal{L}(\mathcal{A} \otimes \mathcal{O}, \mathcal{O}^{\otimes 2}).
\end{aligned}
\]
 From Section 3.6 in \citep{collinswoodfin2023hitting}, we have
\begin{equation} \label{eq:statistic_inner_function}
\begin{aligned}
    \Dif ( \ip{\hat{X}\otimes \hat{X}, K_i}_{\mathcal{A}^{\otimes 2}}) 
    &
    =
    \ip{\Dif \hat{X} \otimes \hat{X}, K_i}_{\mathcal{A}^{\otimes 2}} + \ip{\hat{X} \otimes \Dif \hat{X}, K_i}_{\mathcal{A}^{\otimes 2}}
    \\
    &
    \cong
    (\text{Id}_{\mathcal{O}} \oplus 0_{\mathcal{T}}) \otimes \ip{\hat{X}, K_i}_{\mathcal{A}} + \ip{\hat{X}, K_i}_{\mathcal{A}} \otimes ( \text{Id}_{\mathcal{O}} \oplus 0_{\mathcal{T}}). 
\end{aligned}
\end{equation}
Once again, from the product rule and choosing an orthonormal basis $\{e_\alpha \otimes f_o\}$ for $\mathcal{A} \otimes \mathcal{O}$, 
\begin{equation}
    \begin{aligned}
    \Dif ( \ip{ X, X}_{\mathcal{A}}) &= \ip{\text{Id}_{\mathcal{A} \otimes \mathcal{O}}, X}_{\mathcal{A}} + \ip{X, \text{Id}_{\mathcal{A} \otimes \mathcal{O}}}_{\mathcal{A}}\\
    &\cong\sum_{o, \alpha} \ip{e_\alpha \otimes f_o, X}_{\mathcal{A}} \otimes e_\alpha \otimes f_o + \ip{X, e_\alpha \otimes f_o}_{\mathcal{A}}\otimes e_\alpha \otimes f_o\\
    &=\sum_{o, \alpha}  f_o \otimes \ip{e_\alpha , X}_{\mathcal{A}} \otimes e_\alpha \otimes f_o + \ip{X, e_\alpha }_{\mathcal{A}} \otimes f_o \otimes e_\alpha \otimes f_o\\
    &=\sum_{o}  f_o \otimes X \otimes f_o + X \otimes f_o \otimes f_o\\
    &\cong I_{\mathcal{O}} \otimes X + X \otimes I_{\mathcal{O}}.
    \end{aligned}
\end{equation}
Using a similar argument, we obtain
\begin{equation}
    \begin{aligned}
    \Dif(\langle \mu_i , \hat{X})_{\mathcal{A}}) 
    &
    \cong 
    \langle \mu_i, \text{Id}_{\mathcal{A} \otimes \mathcal{O}} \oplus 0_{\mathcal{A} \otimes \mathcal{T}}\rangle_{\mathcal{A}} 
    \\
    &
    \cong
    \sum_{o,\alpha}\langle \mu_i, e_{\alpha} \otimes (f_o \oplus 0_{\mathcal{T}})\rangle_{\mathcal{A}}  \otimes  e_{\alpha} \otimes f_o
    \\
    &
    =
    \sum_{o,\alpha}\langle \mu_i, e_{\alpha} \rangle_{\mathcal{A}}   (f_o \oplus 0_{\mathcal{T}}) \otimes e_{\alpha} \otimes f_o
    \\
    &
    =
    \sum_{o,\alpha}\langle \mu_i, e_{\alpha} \rangle_{\mathcal{A}}   e_{\alpha} \otimes (f_o \oplus 0_{\mathcal{T}}) \otimes f_o
    \\
    &
    =
    \sum_{o}\mu_i\otimes (f_o \oplus 0_{\mathcal{T}})   \otimes f_o
    \\
    &
    \cong 
    \mu_i\otimes (\Id_{\mathcal{O}} \oplus 0_{\mathcal{T}})  .
    \end{aligned}
\end{equation}
Similarly, $ \Dif(\langle \hat{X} , \mu_i\rangle _{\mathcal{A}}) =  (\Id_{\mathcal{O}} \oplus 0_{\mathcal{T}}) \otimes \mu_i$. Using these results, we see from the chain rule and product rule
\begin{equation}
    \begin{aligned}
        &\Dif ( \widehat{B}_i(W_i)) 
        \\
        &\cong \begin{bmatrix}
            (\text{Id}_{\mathcal{O}} \oplus 0_{\mathcal{T}} ) \otimes \ip{\hat{X}, K_i}_{\mathcal{A}} + \ip{\hat{X}, K_i}_{\mathcal{A}} \otimes ( \text{Id}_{\mathcal{O}} \oplus 0_{\mathcal{T}} ) &  \mu_i\otimes (\Id_{\mathcal{O}} \oplus 0_{\mathcal{T}}) &  0_{ \mathcal{O}^+ \otimes \mathcal{O} } \\
            (\Id_{\mathcal{O}} \oplus 0_{\mathcal{T}}) \otimes \mu_i & 0_{ \mathcal{O}^{\otimes 2} } & 0_{ \mathcal{O}^{\otimes 2} } \\
            0_{ \mathcal{O} \otimes \mathcal{O}^+ } &0_{ \mathcal{O}^{\otimes 2} } & I_{\mathcal{O}} \otimes X + X \otimes I_{\mathcal{O}}
        \end{bmatrix}
    \end{aligned}
\end{equation}
From the chain rule \citep[See Section 3.5.1]{collinswoodfin2023hitting}, we have $\nabla \mathcal{R}^{\operatorname{adv}}(X)\in \mathcal{L}(\mathcal{A} \otimes \mathcal{O}, \R)$ and
\begin{equation}\label{eq:R:grad:todo}
    \begin{aligned}
   \nabla \mathcal{R}^{\operatorname{adv}}(X) &= p_1\ip{\nabla h_1, D (\widehat{B}_1(W))}_{(\mathcal{O}^{+} \oplus \mathcal{O}\oplus \mathcal{O})^{\otimes 2}} 
   \\
   &\qquad+ p_2\ip{\nabla h_2, D (\widehat{B}_2(W))}_{(\mathcal{O}^{+} \oplus \mathcal{O}\oplus \mathcal{O})^{\otimes 2}},
    \end{aligned}
\end{equation}
where $\nabla h_i$ is evaluated at $ \widehat{B}_i(W)$. 
Denote 
 \begin{equation}
    \nabla g_{1,i} = \left[\begin{array}{c|c} 
       \nabla_{11} h_i &\nabla_{12} h_i  \\
        \hline
       \nabla_{21} h_i &\nabla_{22} h_i \\
    \end{array}\right]\quad \text{and}\quad \nabla g_{2,i} = \left[\begin{array}{c} 
       \nabla_{13} h_i  \\
         \hline
       \nabla_{23} h_i 
     \end{array}\right].
 \end{equation}
 Since $\widehat{B}_i = \widehat{B}_i^\top$ then $\nabla_{j\ell} h_{i} = \nabla_{\ell j} h_{i}$ and through a simple computation, we obtain explicit expression
\begin{equation}\label{eq:simp:Radv:grad}
\begin{aligned}
    \nabla & \mathcal{R}^{\operatorname{adv}}(X)  \\
    &
    \cong 
    \sum_{i=1}^2 p_i\bigg(\ip{ \nabla g_{1,i}, (\text{Id}_{\mathcal{O}} \oplus 0_{\mathcal{T}} ) \otimes \ip{\hat{X}, K_i}_{\mathcal{A}} + \ip{\hat{X}, K_i}_{\mathcal{A}} \otimes ( \text{Id}_{\mathcal{O}} \oplus 0_{\mathcal{T}})}_{(\mathcal{O}^+ )^{\otimes 2}} 
    \\
    &
    \qquad
    +
    \langle \nabla g_{2,i}, \mu_i\otimes (\Id_{\mathcal{O}} \oplus 0_{\mathcal{T}})  \rangle_{\mathcal{O} \otimes \mathcal{O}^+} + \langle(\nabla g_{2,i})^\top , (\Id_{\mathcal{O}} \oplus 0_{\mathcal{T}})   \otimes \mu_i  \rangle_{\mathcal{O}^+ \otimes \mathcal{O}} 
    \bigg)
    \\
    &
    \qquad
    + \ip{\nabla_{44} h_i,I_{\mathcal{O}} \otimes X + X \otimes I_{\mathcal{O}}}_{\mathcal{O}^{\otimes 2}}.
\end{aligned}
\end{equation}
 We will simplify each term separately. Taking the first two terms in \eqref{eq:simp:Radv:grad}, we obtain
\begin{equation}
    \begin{aligned}
        &\ip{ \nabla g_{1,i}, (\text{Id}_{\mathcal{O}} \oplus 0_{\mathcal{T}}) \otimes \ip{\hat{X}, K_i}_{\mathcal{A}} + \ip{\hat{X}, K_i}_{\mathcal{A}} \otimes ( \text{Id}_{\mathcal{O}} \oplus 0_{\mathcal{T}})}_{(\mathcal{O}^+ )^{\otimes 2}}\\
        &
        =\ip{ \nabla g_{1,i}, (\text{Id}_{\mathcal{O}} \oplus 0_{\mathcal{T}}) \otimes \ip{\hat{X}, K_i}_{\mathcal{A}}}_{(\mathcal{O}^+ )^{\otimes 2}} 
        \\
        &
        \qquad
        + \ip{ \nabla g_{1,i}, \ip{\hat{X}, K_i}_{\mathcal{A}} \otimes ( \text{Id}_{\mathcal{O}} \oplus 0_{\mathcal{O}}\oplus 0_{\mathcal{T}})}_{(\mathcal{O}^+ )^{\otimes 2}}\\
        &
        \cong \ip{\ip{ \nabla g_{1,i}, (\text{Id}_{\mathcal{O}} \oplus 0_{\mathcal{T}}) \otimes \hat{X} }_{(\mathcal{O}^+ )^{\otimes 2}}, K_i}_{\mathcal{A}} 
        \\
        &
        \qquad+ \ip{\ip{ \nabla g_{1,i}, \hat{X} \otimes ( \text{Id}_{\mathcal{O}} \oplus 0_{\mathcal{T}})}_{(\mathcal{O}^+ )^{\otimes 2}}, K_i}_{\mathcal{A}}\\
        &= 2\ip{\ip{\nabla_{11} h_i, X }_{\mathcal{O}}, K_i}_{\mathcal{A}} + 2\ip{\ip{\nabla_{12} h_i, X^{\star} }_{\mathcal{T}}, K_i}_{\mathcal{A}}.
    \end{aligned}
\end{equation}
Similarly, for the second term in \eqref{eq:simp:Radv:grad}, we have
\begin{equation}
    \begin{aligned}
        \langle \nabla g_{2,i}, \mu_i\otimes (\Id_{\mathcal{O}} \oplus 0_{\mathcal{T}})  \rangle_{\mathcal{O} \otimes \mathcal{O}^+} 
        &= \langle \nabla_{13}h_i, \mu_i\rangle_{\mathcal{O}} ,
        \\
        \langle(\nabla g_{2,i})^\top , (\Id_{\mathcal{O}} \oplus 0_{\mathcal{T}})   \otimes \mu_i  \rangle_{\mathcal{O}^+ \otimes \mathcal{O}} 
        &= \langle \nabla_{13}h_i, \mu_i\rangle_{\mathcal{O}}. 
    \end{aligned}
\end{equation}
Continuing with the last term in \eqref{eq:simp:Radv:grad}, from a simple computation we have
\begin{equation}
    \begin{aligned}
       \ip{\nabla_{44} h_i,I_{\mathcal{O}} \otimes X + X \otimes I_{\mathcal{O}}}_{\mathcal{O}^{\otimes 2}} =2 \ip{\nabla_{44} h_i, X}_{\mathcal{O}}.
    \end{aligned}
\end{equation}
Putting it all together, we obtain
\begin{equation}
    \begin{aligned}
        \nabla  \mathcal{R}^{\operatorname{adv}}(X) 
        &
        = 
        2\sum_{i=1}^2 p_i\bigg(\ip{\ip{\nabla_{11} h_i, X }_{\mathcal{O}}, K_i}_{\mathcal{A}} + \ip{\ip{\nabla_{12} h_i, X^{\star} }_{\mathcal{T}}, K_i}_{\mathcal{A}}
        \\
        &
        \qquad
        + \langle \nabla_{13}h_i, \mu_i\rangle_{\mathcal{O}}  +  \ip{\nabla_{44} h_i, X}_{\mathcal{O}}\bigg).
    \end{aligned}
\end{equation}
Setting $\mathcal{O} = \mathcal{T} = 1$ and $\mathcal{A} = \R^d$ completes the proof.
\end{proof}
\subsection{Concentration and Pseudo-Lipschitz}\label{sect:conc:pseudoLispchitz}
In order to simplify the proofs in section~\ref{sec:error_bounds}, we will make use of the subgaussian norm $\|\cdot\|_{\psi_2}$ (see e.g., \citep{vershynin2018high} for more details):
\begin{equation}\label{eq:subgaussian_norm}
    \| X \|_{\psi_2} 
    \asymp
    \inf \{ V  > 0 : \forall~t > 0~\Pr( |X| > t) \leq 2 e^{-t^2/V^2}\}.
\end{equation}
We will also adapt some concentration results for pseudo-Lipschitz functions from \citep{collinswoodfin2023hitting} to our adversarial framework.

The proofs in section~\ref{sec:error_bounds} will involve statistics whose gradients are $\alpha$-pseudo-Lipschitz. Before introducing bounds on the norm of the gradient of $\alpha$-pseudo-Lipschitz functions, we define the following
projection operator onto the ball of radius $\beta$ in the vector space $V$, $\text{\rm Proj}_{\beta} \, : \, V \to V$:
\begin{equation}
    \begin{aligned}
        \text{Proj}_{\beta}(x) 
        &
        \defas
         \begin{cases}
        x & \text{ if $\|x\| \le \beta$ } \\
        \beta \frac{x}{\|x\|} & \text{otherwise.}
        \end{cases}
    \end{aligned}
\end{equation}
It can be shown that taking compositions of projections with $\alpha$-pseudo-Lipschitz functions creates Lipschitz functions.
\begin{lemma}[Lemma 3, \citep{collinswoodfin2023hitting}]\label{lem:lipschitz:proj}
    Suppose $f$ is $\alpha$-pseudo-Lipschitz with constant $L$. Then the composition $f \circ \text{\rm Proj}_\beta$ is Lipschitz with constant $L(1 + 2 \beta^{\alpha})$. 
\end{lemma}
By Assumption~\ref{assumption:pseudo_lipschitz} $f$ $\alpha$-pseudo-Lipschitz. We build upon adapt Lemma 3.4 from \citep{collinswoodfin2023hitting} to obtain some bounds on the growth of moments of $ f'$ in terms of $\hat{X} = X \oplus X^{\star}$. We omit the proof as it follows an almost identical argument.
\begin{lemma}[Growth of $f'$] \label{lem:growth_grad_f} Suppose a function $f \, : \, \R^{3} \to \mathbb{R}$ is $\alpha$-pseudo-Lipschitz with constant $L(f)$. Suppose the noise $\epsilon \sim \mathcal{N}(0,\eta^2)$ is independent of $v\sim \mathcal{N}(0, K)$. Then for $p > 0$ and any $g\in \R^2$,
\begin{equation} \label{eq:growth_grad_f}
|  f'(g)|^p \le C(\alpha, p) (L(f))^p (1 + \|g\| + |\epsilon|)^{\max \{1, \alpha p\}},
\end{equation}
Moreover, for $\mu  \in \R^d$, if $g = \hat{X}^\top v +  \hat{X}^\top\mu + \begin{bmatrix}  \delta s\|X\| & 0 \end{bmatrix}^\top$ with $ \hat{X}= X \oplus X^{\star}$ and $|s| \leq 1$, there is a growth rate on $f'(g)$ and sub-Gaussian norm on $g$ in terms of $\hat{X}$, $X$ and $\eta$, 
    \begin{equation} \label{eq:expectation_f_growth}
    \begin{gathered}
        \EE_{v,\epsilon}[|f'(g) |^p] \le C(\alpha, p)  (L(f))^p \big  (1+(\|K\|_{\opt}^{1/2} +\|\mu\|)\|\hat{X} \| + \delta \|X\|+ \eta )^{\max \{1, \alpha p\} }
        \\
            \text{and} \qquad \| (1 + \|g\|+ |\epsilon| )\|_{\psi_2} \le C ( 1+(\|K\|_{\opt}^{1/2} + \|\mu\| )\|\hat{X} \| + \delta \|X\|+ \eta  ).
    \end{gathered}
    \end{equation}
\end{lemma}

\subsection{Doob Decomposition of Adversarial Homogenized SGD} \label{sect:DoobDecomp:HSGD}
The arguments in this section are built upon the results presented in Appendix A of \citep{collinswoodfin2023hitting} and Section 5 in \citep{collinswoodfin2025Exact}. We present the steps for the soft label setting as the hard label setting follows similarly. We are interested in the behaviour of AdvHSGD under the following statistics:
\begin{equation}\label{eq:def:statistics:S(X,z):M(X,z)}
\begin{aligned}
S(W, z) \defas\hat{X}^\top \mathscr{R}(z; \mathcal{K}) \hat{X} \in \mathbb{C}^{2\times 2},\quad M(W, z) \defas \hat{X}^\top \mathscr{R}(z; \mathcal{K}) \mu\in \mathbb{C}^{2\times 2}, \quad \text{for} \quad z \in \Gamma^{2},
\end{aligned}
\end{equation} 
where we recall $\mathcal{K} \defas \{K_1, K_2\}$, $\Gamma \defas \{w\in \mathbb{C}: |w| = \max(1, 2\|\mathcal{K}\|_{\opt})\}$ with $\|\mathcal{K}\|_{\opt}  \defas \max_{i=1,2} \|K_i\|_{\opt}$ and $\mathscr{R}(z; \mathcal{K}) \defas R(z_1;K_1) R(z_2;K_2)$ with $z_i \in \Gamma$ for $i=1,2$. Recall from \eqref{eq:main:AdvHSGD:def} the definition of AdvHSGD 
\begin{equation}
\dif \WHSGD_t = -\gamma(t) \nabla \mathcal{R}_{\lambda}^{\operatorname{adv}}(\WHSGD_t) \dif t + \gamma(t) \sqrt{\frac{1}{d} \sum_{i=1}^2 p_i \EE_{v, \epsilon} [  f'_i( \rho_{t, i}, \epsilon)^2  ]  (K_i + \mu_i\mu_i^\top)} \dif B_t,
\end{equation}
where the initial condition is given by $\mathscr{X}_0 = X_0$, $(B_t, t \ge 0)$ is a $\R^{d}$ standard Brownian motion and for $I \in \{1,2\}$ and $v\sim \mathcal{N}(0, \Id_d)$ we have
\begin{equation}
    \begin{gathered}
\widehat{\mathscr{X}}_t \defas \mathscr{X}_t \oplus X^{\star}\in \R^{d \times 2},\quad\mathscr{W}_t \defas \widehat{\mathscr{X}}_t \oplus \mu\in \R^{d \times 4},
\\
\rho_{t,I} \defas \widehat{\mathscr{X}}_t^\top(\sqrt{K_{I}}v + \mu_{I})+ \begin{bmatrix}\delta u_{t,I} \|\WHSGD_t\| & 0 \end{bmatrix}^\top,
\\
\text{and}\quad u_{t,I} = \operatorname{argmax}_{|u| \leq  1}f(\widehat{\mathscr{X}}_t^\top (\sqrt{K_{I}}v + \mu_{I}) + \begin{bmatrix} \delta u\|\mathscr{X}_t \|& 0 \end{bmatrix}^\top, \epsilon).
    \end{gathered}
\end{equation}
From \eqref{eq:def:statistics:S(X,z):M(X,z)}, it is immediate that $S(\HSGD_t, z)$ is a bilinear form. Hence, we decompose it into block form
\begin{gather*}
S(\HSGD_t,z) 
 =
 \left[ \begin{array}{c|c} 
\WHSGD_t^\top \mathscr{R}(z;\mathcal{K}) \WHSGD_t & \WHSGD_t^\top \mathscr{R}(z;\mathcal{K}) X^{\star}\\
\hline
(X^{\star\top}) \mathscr{R}(z;\mathcal{K}) \WHSGD_t & (X^{\star})^{\top} \mathscr{R}(z;\mathcal{K}) X^{\star}
\end{array}
 \right]
 =
 \left[ \begin{array}{c|c} 
    S_{11}(\HSGD_t,z) & S_{12}(\HSGD_t,z)
    \\
    \hline
    S_{21}(\HSGD_t,z) & S_{22}(\HSGD_t,z)
    \end{array} \right].
\end{gather*}
With this notation established, we derive an explicit equation for $S(\HSGD_t, z)$ in terms of its martingale increment to facilitate the proof Proposition~\ref{prop:approx:sols:HSGD} using Itô calculus.

\paragraph{It\^o calculus applied to $S(\HSGD_t, z)$.}
From Assumption~\ref{assumption:risk}, recall that the adversarial expected risk $\mathcal{R}^{\operatorname{adv}}$ can be decomposed as the sum of two function compositions
\[
    \mathcal{R}^{\operatorname{adv}}(\WHSGD_t) = \sum_{i=1}^2 p_i \cdot (h_i \circ \widehat{B}_i(\widehat{\HSGD}_t)),
\]
for some functions $h_i \, : \, \R^{4 \times4 } \to \mathbb{R}$ for $i=1,2$ and 
\begin{equation}
    \widehat{B}_i(\HSGD_t) \defas \begin{bmatrix}
           U_i(\HSGD_t) & 0\\
           0 & \|\WHSGD_t\|^2
           \end{bmatrix},\quad  U_i(\HSGD_t) \defas \begin{bmatrix}
               B_i(\HSGD_t)  & m_i(\HSGD_t) \\
               m_i(\HSGD_t)^\top & \|\mu_i\|^2
           \end{bmatrix},\quad B_i(\HSGD_t)\defas \widehat{\mathscr{X}}_t^\top K_i \widehat{\mathscr{X}}_t
   \end{equation}
and $m_i(\HSGD_t) \defas \widehat{\mathscr{X}}_t^\top\mu_i$
Recall from Lemma~\ref{lem:derivative_adv_risk} the gradient of the adversarial risk
\begin{equation*}
    \begin{aligned}
        \nabla  \mathcal{R}^{\operatorname{adv}}(X) &= 2\sum_{i=1}^2 p_i \bigg(\partial_{11} h_i(\widehat{B}_i(W)) \cdot K_iX + \partial_{12} h_i (\widehat{B}_i(W))\cdot K_iX^{\star}
        \\
        &
        \qquad +  \partial_{13} h_i(\widehat{B}_i(W)) \cdot  \mu_i +  \partial_{44}  h_i (\widehat{B}_i(W))\cdot X\bigg)
    \end{aligned}
\end{equation*}
 For simplicity of notation, we will omit the composition of the partial derivatives with $ \widehat{B}_i(\HSGD_t)$ in the next steps of the proof. Since $\dif \widehat{\mathscr{X}}_t = \dif \mathscr{X}_t \oplus 0_{d }$ and using the  product rule for It\^o derivatives,
\begin{equation} \label{eq:concentration_S_1}
    \begin{aligned}
        \dif S 
        &
        = 
        \ip{ \dif \widehat{\mathscr{X}}_t \otimes \widehat{\mathscr{X}}_t, \mathscr{R}(z;\mathcal{K})}_{\R^{d\times d}}+ \ip{\widehat{\mathscr{X}}_t \otimes \dif \widehat{\mathscr{X}}_t,  \mathscr{R}(z;\mathcal{K})}_{\R^{d\times d}}
        + \ip{ \dif \widehat{\mathscr{X}}_t \otimes \dif \widehat{\mathscr{X}}_t,  \mathscr{R}(z;\mathcal{K})}_{\R^{d\times d}}
        \\
        &
        = 
        \ip{ (\dif \WHSGD_t \oplus 0_{d}) \otimes \widehat{\mathscr{X}}_t,  \mathscr{R}(z;\mathcal{K})}_{\R^{d\times d}} + \ip{\widehat{\mathscr{X}}_t \otimes (\dif \WHSGD_t \oplus 0_{d}), \mathscr{R}(z;\mathcal{K})} _{\R^{d\times d}}
        \\
        & 
        \quad + 
        \ip{ (\dif \WHSGD_t \oplus 0_{d}) \otimes (\dif \WHSGD_t \oplus 0_{d}),  \mathscr{R}(z;\mathcal{K})}_{\R^{d\times d}}
        \\
        &
        =
        -\gamma(t) \cdot \ip{ ( (\nabla \mathcal{R}^{{\operatorname{adv}}}(\WHSGD_t) + \lambda \WHSGD_t) \oplus 0_{d}) \otimes \widehat{\mathscr{X}}_t,  \mathscr{R}(z;\mathcal{K})}_{\R^{d\times d}} \dif t
        \\
        &
        \quad
        - 
        \gamma(t) \cdot \ip{\widehat{\mathscr{X}}_t \otimes ( (\nabla \mathcal{R}^{{\operatorname{adv}}}(\WHSGD_t) + \lambda \WHSGD_t)  \oplus 0_{d}),  \mathscr{R}(z;\mathcal{K})}_{\R^{d\times d}} \dif t
        \\
        & 
        \quad + \small\begin{bmatrix}  \tfrac{\gamma^2_t}{d}\sum_{i=1}^2 p_i\EE_{v, \epsilon}[ f'_i(\rho_{t,i}, \epsilon)^2] \Tr(( K_i + \mu_i \mu_i^\top ) \mathscr{R}(z;\mathcal{K}))\dif t& 0 \\
            0 & 0 \end{bmatrix} 
        + \dif \mathcal{M}_t^{\operatorname{AdvHSGD}}(S(\HSGD_t, z)), 
    \end{aligned}
\end{equation}
where we define the martingale increment of $S(\HSGD_t,z) $ elementwise such that for $1\leq j,\ell \leq 2$
    \[
    \begin{aligned}
        \big ( \dif &\mathcal{M}_t^{\operatorname{AdvHSGD}}(S(\HSGD_t,z)) \big )_{j\ell}\defas\dif \mathcal{M}^{\operatorname{AdvHSGD}}_t(S_{j\ell}(\HSGD_t,z))
        \\
         &\defas 
         \gamma(t)   \left\langle \sqrt{\frac{1}{d} \sum_{i=1}^2 p_i \EE_{v, \epsilon} [  f'_i( \rho_{t, i}, \epsilon)^2  ](K_i + \mu_i\mu_i^\top)}, \nabla_X (S_{j\ell}(\HSGD_t,z))  \otimes \dif B_t \right\rangle.
    \end{aligned}
\]
Similarly, we define the martingale term elementwise
\[
    \big ( \mathcal{M}_t^{\operatorname{AdvHSGD}}(S(\cdot, z)))_{j\ell}  = \int_0^t \dif \mathcal{M}_s^{\operatorname{AdvHSGD}}(S_{j\ell}(\cdot,z)).
\]
 Plugging in $\nabla \mathcal{R}^{{\operatorname{adv}}}$ into the first term of \eqref{eq:concentration_S_1} and Developing, we have
\begin{align*}
&\ip{ ( (\nabla \mathcal{R}^{{\operatorname{adv}}}(\WHSGD_t) + \lambda \WHSGD_t) \oplus 0_{d}) \otimes \widehat{\mathscr{X}}_t, \mathscr{R}(z;\mathcal{K})}_{\R^{d\times d}} \, \dif t \\
&= 2\sum_{i=1}^2 p_i\bigg(\ip{ (  \partial_{11} h_i \cdot K_i \WHSGD_t  \oplus 0_{d}) \otimes \widehat{\mathscr{X}}_t, \mathscr{R}(z;\mathcal{K})}_{\R^{d\times d}} \, 
\\
&
\qquad 
+\ip{ (\partial_{12} h_i\cdot  K_i X^{\star}  \oplus 0_{d}) \otimes \widehat{\mathscr{X}}_t, \mathscr{R}(z;\mathcal{K})}_{\R^{d\times d}} \, 
+\ip{ (\partial_{13} h_i\cdot   \mu_i  \oplus 0_{d}) \otimes \widehat{\mathscr{X}}_t, \mathscr{R}(z;\mathcal{K})}_{\R^{d\times d}} \, 
\\
&
\qquad
+\ip{ (  \partial_{44} h_i \cdot  \WHSGD_t  \oplus 0_{d}) \otimes \widehat{\mathscr{X}}_t, \mathscr{R}(z;\mathcal{K})}_{\R^{d\times d}} \, \bigg)\dif t
+ \lambda\ip{   (\WHSGD_t \oplus 0_{d}) \otimes \widehat{\mathscr{X}}_t, \mathscr{R}(z;\mathcal{K})}_{\R^{d\times d}} \, \dif t 
\end{align*}
Expanding the terms with $\widehat{\mathscr{X}}_t = \WHSGD_t \oplus X^{\star} $ then we obtain
\begin{equation}\label{eq:S:terms1}
\begin{aligned}
    -\gamma(t) \cdot & \ip{ ( (\nabla \mathcal{R}^{{\operatorname{adv}}}(\WHSGD_t) + \lambda \WHSGD_t) \oplus 0_{d}) \otimes \widehat{\mathscr{X}}_t,\mathscr{R}(z;\mathcal{K})}_{\R^{d\times d}} \, \dif t 
    =
    \left [ \begin{array}{c|c}
        A_1 + \tilde{A}_1 & E + \tilde{E}\\
        \hline
        0 & 0
    \end{array} \right ], 
    \\
    \text{where} 
    \quad A_1 
    &
    = -2\gamma(t) \sum_{i=1}^2 p_i\bigg(
     \partial_{11} h_i\cdot \ip{ \WHSGD_t \otimes \WHSGD_t, K_i \mathscr{R}(z;\mathcal{K})}
    +\partial_{12} h_i \cdot \ip{ X^{\star} \otimes \WHSGD_t, K_i \mathscr{R}(z;\mathcal{K})} 
    \\
    &
    +
    \partial_{13} h_i \cdot \ip{ \mu_i\otimes \WHSGD_t,  \mathscr{R}(z;\mathcal{K})} 
+\partial_{44} h_i \cdot \ip{ \WHSGD_t\otimes \WHSGD_t,  \mathscr{R}(z;\mathcal{K})}  \, \bigg)\dif t,
    \\
    \quad \tilde{A}_1 
    &
    = -\gamma(t) \cdot \lambda \cdot \ip{\WHSGD_t \otimes \WHSGD_t,  \mathscr{R}(z;\mathcal{K})} \, \dif t,
    \\
    \quad
    E 
    &
    = -2\gamma(t) \sum_{i=1}^2 p_i\bigg(
     \partial_{11} h_i\cdot \ip{ \WHSGD_t \otimes X^{\star}, K_i \mathscr{R}(z;\mathcal{K})}
    +\partial_{12} h_i \cdot \ip{ X^{\star} \otimes X^{\star}, K_i \mathscr{R}(z;\mathcal{K})} 
    \\
    &
    +
    \partial_{13} h_i \cdot \ip{ \mu_i\otimes X^{\star},  \mathscr{R}(z;\mathcal{K})} 
+\partial_{44} h_i \cdot \ip{ \WHSGD_t\otimes X^{\star},  \mathscr{R}(z;\mathcal{K})}  \, \bigg)\dif t,
    \\
    \text{and} 
    \quad
    \tilde{E} 
    &
    = -\gamma(t) \cdot \lambda \cdot \ip{ \WHSGD_t \otimes X^{\star}, \mathscr{R}(z;\mathcal{K})} \, \dif t.
\end{aligned}
\end{equation}
Similarly for the second term in \eqref{eq:concentration_S_1}, 
\begin{equation}\label{eq:S:terms2}
    \begin{aligned}
    - \gamma(t) \cdot & \ip{\widehat{\mathscr{X}}_t \otimes ( (\nabla \mathcal{R}^{{\operatorname{adv}}}(\mathscr{X}_t) + \lambda \WHSGD_t)  \oplus 0_d), \mathscr{R}(z;\mathcal{K})}_{\R^{d\times d}}\dif t = \left [
    \begin{array}{c|c}
    A_2 + \tilde{A}_2 & 0\\
    \hline
    C + \tilde{C} & 0
    \end{array}
    \right ],
    \\
    \text{where} \quad
    A_2  
        &
        =
        -2\gamma(t) \sum_{i=1}^2 p_i\bigg(
     \partial_{11} h_i\cdot \ip{ \WHSGD_t \otimes \WHSGD_t, K_i \mathscr{R}(z;\mathcal{K})}
    +\partial_{12} h_i \cdot \ip{ \WHSGD_t \otimes X^{\star}  , K_i \mathscr{R}(z;\mathcal{K})} 
    \\
    &
    +
    \partial_{13} h_i \cdot \ip{  \WHSGD_t \otimes \mu_i, \mathscr{R}(z;\mathcal{K})} 
+\partial_{44} h_i \cdot \ip{ \WHSGD_t\otimes \WHSGD_t,  \mathscr{R}(z;\mathcal{K})}  \, \bigg)\dif t,
    \\
    \quad
    \tilde{A}_2 
    &
    =
    -\gamma(t) \cdot \lambda \cdot \ip{ \WHSGD_t \otimes \WHSGD_t, \mathscr{R}(z;\mathcal{K})} \, \dif t,
    \\
    \quad
    C
    &
    =
    -2\gamma(t) \sum_{i=1}^2 p_i\bigg(
     \partial_{11} h_i\cdot \ip{ X^{\star} \otimes \WHSGD_t, K_i \mathscr{R}(z;\mathcal{K})}
    +\partial_{12} h_i \cdot \ip{ X^{\star}  \otimes X^{\star}  , K_i \mathscr{R}(z;\mathcal{K})} 
    \\
    &
    +
    \partial_{13} h_i \cdot \ip{  X^{\star}  \otimes \mu_i, \mathscr{R}(z;\mathcal{K})} 
+\partial_{44} h_i \cdot \ip{ X^{\star} \otimes \WHSGD_t,  \mathscr{R}(z;\mathcal{K})}  \, \bigg)\dif t, 
    \\
    \text{and}
    \quad 
    \tilde{C} 
    &
    =
    -\gamma(t) \cdot \lambda \cdot \ip{X^{\star} \otimes \WHSGD_t, \mathscr{R}(z;\mathcal{K})} \, \dif t.
\end{aligned}
\end{equation}
Finally, denote the last term in \eqref{eq:concentration_S_1} as
\begin{align*}
    A_3=\frac{\gamma^2_t}{d}\sum_{i=1}^2 p_i\EE_{v, \epsilon}[ f'_i(\rho_{t,i}, \epsilon)^2] \Tr(( K_i + \mu_i \mu_i^\top ) \mathscr{R}(z;\mathcal{K}))\dif t.
\end{align*}
It follows then that 
\[
\begin{aligned}
(\dif S)(\HSGD_t,z) 
&
= 
\left [ 
\begin{array}{c|c}
\dif S_{11} & \dif S_{12} \\
\hline
\dif S_{21} & \dif S_{22}
\end{array}
\right ] 
\\
&
= 
\left [ 
\begin{array}{c|c}
A_1 + \tilde{A}_1 + A_2 + \tilde{A}_2 + A_3  & E + \tilde{E} \\
\hline
C + \tilde{C} & 0
\end{array}
\right ] + \dif \mathcal{M}_t^{\operatorname{AdvHSGD}}(S(\HSGD_t, z)).
\end{aligned}
\]
Through a simple computation, we observe that 
\[
\begin{aligned}
K_1 \mathscr{R}(z;\mathcal{K}) 
&= R(z_2;K_2) + z_1 \mathscr{R}(z;\mathcal{K}).
\end{aligned}
\]
Similarly, we have $K_2 \mathscr{R}(z;\mathcal{K}) =R(z_1;K_1) + z_2 \mathscr{R}(z;\mathcal{K})$. From Cauchy's integral formula, we also note that 
\[
    \widehat{\mathscr{X}}_t^\top R(z_1; K_1) \widehat{\mathscr{X}}_t = \frac{-1}{2\pi\operatorname{i}} \oint_{\Gamma} S(\HSGD_t, z) \dif z_2.
\]
 Similar results hold for $z_2$. Define $M_i(\HSGD_t, z) \defas \widehat{\mathscr{X}}_t^\top \mathscr{R}(z;\mathcal{K})\mu_i$ for $i=1,2$. Using the above identity, the quantities $A_1, A_2, E$, $C$, $\tilde{A}_1, \tilde{A}_2, \tilde{E}$, and $\tilde{C}$ from \eqref{eq:S:terms1} and \eqref{eq:S:terms2} can be reexpressed in terms of $S(\HSGD_t, z) = \ip{\widehat{\mathscr{X}}_t \otimes \widehat{\mathscr{X}}_t, \mathscr{R}(z;\mathcal{K})}$. It follows that
\begin{equation} \label{eq:coupled:AdvHSGD:S}
\begin{aligned}
    &\dif S(\HSGD_t,z) 
    =
    -2\gamma(t) \sum_{i=1}^2 p_i\bigg(\left(z_iS(\HSGD_t,z)- \frac{1}{2\pi\operatorname{i}} \oint_{\Gamma} S(\HSGD_t, z) \dif z_i\right) H_{1,i}(\widehat{B}_i(\HSGD_t)) 
    \\
    &
     + H_{1,i}^\top( \widehat{B}_i(\HSGD_t))  \left(z_iS(\HSGD_t,z)- \frac{1}{2\pi\operatorname{i}} \oint_{\Gamma} S(\HSGD_t, z) \dif z_i\right)
    \\
    &
    + M_i(\HSGD_t, z) \otimes H_{2,i} (\widehat{B}_i(\HSGD_t)) 
    + H_{2,i} (\widehat{B}_i(\HSGD_t)) \otimes M_i(\HSGD_t, z)
    \\
    &
    + \left(\frac{\lambda}{2} + \partial_{44} h_i (\widehat{B}_i(\HSGD_t) )\right) (S(\HSGD_t,z) \operatorname{D}  +  \operatorname{D}  S(\HSGD_t,z))\bigg)\dif t
    \\
    &
    + \frac{\gamma^2_t}{d} \sum_{i=1}^2 p_i\EE_{v, \epsilon}[ f'_i(\rho_{t,i}, \epsilon)^2] \Tr(( K_i + \mu_i \mu_i^\top ) \mathscr{R}(z;\mathcal{K}))\dif t  \operatorname{D}  \, \dif t
    \\
    &
    + 
    \dif \mathcal{M}_t^{\operatorname{AdvHSGD}}(S(\HSGD_t,z)),
\end{aligned}
\end{equation}
\begin{gather*}
\text{where} \quad
H_{1,i}(\widehat{B}_i)=\left [ \begin{array}{c|c} 
   \partial_{11} h_i(\widehat{B}_i)  & 0 \\
    \hline
   \partial_{21} h_i(\widehat{B}_i) & 0 
    \end{array} \right ],\quad
    H_{2,i}(\widehat{B}_i)= \left [ \begin{array}{c} 
        \partial_{13} h_i(\widehat{B}_i)   \\
         \hline
        0 
         \end{array} \right ],\quad
         \operatorname{D}  = \left [ \begin{array}{c|c} 
 1 & 0\\
 \hline
 0 & 0
 \end{array} \right ],
 \\
\text{and initialized with} \; S(0,z) = \hat{X}_0^\top \mathscr{R}(z; \mathcal{K}) \hat{X}_0. 
\end{gather*}
Now, recall by Assumption~\ref{assumption:fisher}, there exists a function $\mathcal{I}_i: \R^{4\times 4}\to \R$ such that
\[
    \EE_{v, \epsilon}[ f'_i(\rho_{t,i}, \epsilon)^2] = \mathcal{I}_i( \widehat{B}_i(\HSGD_t)). 
\]
From the Cauchy integral formula, we have
\begin{equation}
    \widehat{B}_i(W) = \left[ \begin{array}{ccc} 
        \displaystyle\oint z_i S(W,z) \,\Dif z&\displaystyle \oint M_i(W,z) \,\Dif z &0 \\
        \displaystyle\oint M_i(W,z)^\top \,\Dif z&  \displaystyle\oint \mu_i^\top \mathscr{R}(z; \mathcal{K})\mu_i\,\Dif z& 0\\
            0 & 0 & \displaystyle\oint S_{11}(W, z) \, \Dif z
        \end{array} 
        \right] .
\end{equation}
where $M_i(W,z)$ denotes the $i$th columns of $M(W,z)$. We must now derive the Doob decomposition of $M(\HSGD_t, z)$. We first decompose it into block form
\begin{gather*}
    M(\HSGD_t,z) 
     =
     \left[ \begin{array}{c|c} 
    \WHSGD_t^\top \mathscr{R}(z;\mathcal{K})\mu_1 &  \WHSGD_t^\top\mathscr{R}(z;\mathcal{K})\mu_2 \\
    \hline
    X^{\star \top} \mathscr{R}(z;\mathcal{K})\mu_1 &  X^{\star \top}\mathscr{R}(z;\mathcal{K})\mu_2 \\
    \end{array}
     \right ]
     =
     \left[ \begin{array}{c|c} 
        M_{11}(\HSGD_t,z) &  M_{12}(\HSGD_t,z)\\
        \hline
        M_{21}(\HSGD_t,z) & M_{22}(\HSGD_t,z)
        \end{array} \right ] \in \R^{2\times 2}.
    \end{gather*}
 From Itô's lemma, we have
\begin{equation}
    \begin{aligned}
        \dif &M(\HSGD_t,z) =  (\dif \widehat{\WHSGD}_t)^\top \mathscr{R}(z;\mathcal{K})\mu
        \\
        &= -\gamma(t) \langle ((\nabla \mathcal{R}_{\lambda}^{\operatorname{adv}}(\WHSGD_t) + \lambda \mathscr{X}_t) \oplus 0_d) \otimes \mu, \mathscr{R}(z;\mathcal{K})\rangle_{\R^{d\times d}} \dif t + \dif \mathcal{M}_t^{\operatorname{AdvHSGD}}(M(\HSGD_t, z)),
    \end{aligned}
\end{equation}
where for $1\leq j,\ell\leq 2$
\begin{equation}
    \begin{aligned}
    \dif &\mathcal{M}_t^{\operatorname{AdvHSGD}}(M_{j\ell}(\HSGD_t, z)) 
    \\
    &\defas \gamma(t) \left\langle \sqrt{\frac{1}{d} \sum_{i=1}^2 p_i \EE_{v, \epsilon} [  f'_i( \rho_{t, i}, \epsilon)^2  ]  (K_i + \mu_i\mu_i^\top)}, \nabla_X (M_{j\ell}(\HSGD_t, z)) \otimes \dif B_t\right\rangle.
    \end{aligned}
\end{equation}
Using a similar argument as we did in \eqref{eq:S:terms1}, we obtain
\begin{equation}\label{eq:M:terms1}
    \begin{aligned}
        -\gamma(t) \cdot & \ip{ ( (\nabla \mathcal{R}^{{\operatorname{adv}}}(\WHSGD_t) + \lambda \WHSGD_t) \oplus 0_{d}) \otimes \mu,\mathscr{R}(z;\mathcal{K})}_{\R^{d\times d}} \, \dif t 
        =
        \left [ \begin{array}{c|c}
            O_1 + \tilde{O}_1 & P + \tilde{P}\\
            \hline
            0 & 0
        \end{array} \right ], 
        \\
        \text{where} 
        \quad O_1 
        &
        = -2\gamma(t) \sum_{i=1}^2 p_i\bigg(
         \partial_{11} h_i\cdot \ip{ \WHSGD_t \otimes \mu_1, K_i \mathscr{R}(z;\mathcal{K})}
        +\partial_{12} h_i \cdot \ip{ X^{\star} \otimes \mu_1, K_i \mathscr{R}(z;\mathcal{K})} 
        \\
        &
        +
        \partial_{13} h_i \cdot \ip{ \mu_i\otimes \mu_1,  \mathscr{R}(z;\mathcal{K})} 
    +\partial_{44} h_i \cdot \ip{ \WHSGD_t\otimes \mu_1,  \mathscr{R}(z;\mathcal{K})}  \, \bigg)\dif t,
        \\
        \quad \tilde{O}_1 
        &
        = -\gamma(t) \cdot \lambda \cdot \ip{\WHSGD_t \otimes \WHSGD_t,  \mathscr{R}(z;\mathcal{K})} \, \dif t,
        \\
        \quad
        P 
        &
        = -2\gamma(t) \sum_{i=1}^2 p_i\bigg(
         \partial_{11} h_i\cdot \ip{ \WHSGD_t \otimes \mu_2, K_i \mathscr{R}(z;\mathcal{K})}
        +\partial_{12} h_i \cdot \ip{ X^{\star} \otimes \mu_2, K_i \mathscr{R}(z;\mathcal{K})} 
        \\
        &
        +
        \partial_{13} h_i \cdot \ip{ \mu_i\otimes \mu_2,  \mathscr{R}(z;\mathcal{K})} 
    +\partial_{44} h_i \cdot \ip{ \WHSGD_t\otimes \mu_2,  \mathscr{R}(z;\mathcal{K})}  \, \bigg)\dif t,
        \\
        \text{and} 
        \quad
        \tilde{P} 
        &
        = -\gamma(t) \cdot \lambda \cdot \ip{ \WHSGD_t \otimes X^{\star}, \mathscr{R}(z;\mathcal{K})} \, \dif t.
    \end{aligned}
    \end{equation}
It then follows from the resolvent identities presented above that
\begin{equation} \label{eq:coupled:AdvHSGD:M}
    \begin{aligned}
        &\dif M(\HSGD_t,z) 
        =
        -2\gamma(t) \sum_{i=1}^2 p_i\bigg(H_{1,i}^\top( \widehat{B}_i(\HSGD_t))  \left(z_iM(\HSGD_t,z)- \frac{1}{2\pi\operatorname{i}} \oint_{\Gamma} M(\HSGD_t, z) \dif z_i\right)
        \\
        &
        + H_{2,i} (\widehat{B}_i(\HSGD_t)) \otimes \mu^\top \mathscr{R}(z; \mathcal{K})\mu_i
        + \left(\frac{\lambda}{2} + \partial_{44} h_i (\widehat{B}_i(\HSGD_t) )\right)   \operatorname{D}  M(\HSGD_t,z)\bigg)\dif t 
        \\
        &
        +
        \dif \mathcal{M}_t^{\operatorname{AdvHSGD}}(M(\HSGD_t,z)).
    \end{aligned}
    \end{equation}
From the Doob decompositions \eqref{eq:coupled:AdvHSGD:S} and \eqref{eq:coupled:AdvHSGD:M}, we obtain the representation
\begin{equation}\label{eq:HSGD_exact}
\dif Z(\HSGD_t, \cdot) = \mathscr{F}(z, Z(\HSGD_t, \cdot)) \, \dif t + \dif \mathcal{M}^{\operatorname{AdvHSGD}}_t(Z(\HSGD_t, \cdot)), 
\end{equation}
where
\begin{equation}
    \mathscr{F}(z, Z(\HSGD_t, \cdot)) \defas  \begin{bmatrix}
        \mathscr{F}_\mathcal{S}(z, Z(\HSGD_t, \cdot)) & \mathscr{F}_\mathcal{M}(z, Z(\HSGD_t, \cdot)) \\
        \mathscr{F}_\mathcal{M}(z, Z(\HSGD_t, \cdot))^\top & 0_{2\times 2}
    \end{bmatrix},
\end{equation}
and 
\begin{equation}
    Z(\HSGD_0, \cdot) \defas  \begin{bmatrix}
        S(\HSGD_0, \cdot) &  M(\HSGD_0, \cdot) \\
        M(\HSGD_0, \cdot)^\top & \mu^\top \mathscr{R}(z;\mathcal{K})\mu
    \end{bmatrix} = \begin{bmatrix}
        \hat{X}_0^\top \mathscr{R}(z;\mathcal{K})\hat{X}_0&    \hat{X}_0^\top \mathscr{R}(z;\mathcal{K})\mu \\
        \mu^\top \mathscr{R}(z;\mathcal{K})\hat{X}_0 & \mu^\top \mathscr{R}(z;\mathcal{K})\mu
    \end{bmatrix}.
\end{equation}
 The Doob decomposition in the hard label setting can be recovered using similar arguments or simply be setting $X^{\star} = 0$.

\subsection{Doob Decomposition for SGD}\label{sect:doob:decomp:SGD:derivation}
We provide explicit details behind the derivation of the Doob decomposition of SGD for an arbitrary quadratic statistic presented in Section~\ref{sect:doob:decomp}. In order to streamline notation, we write 
\begin{equation}
    \begin{aligned}
        f_{i}(g_{k, i} ) 
        &
        \defas 
        f_{i}(g_{k, i} ;  \epsilon_{k+1})\quad\text{and}\quad
        f_{i}'(g_{k, i} ) 
        &
        \defas 
        f_{i}'(g_{k, i} ;  \epsilon_{k+1}).
    \end{aligned}
\end{equation}
The arguments in this section build upon Appendix~A in \citep{collinswoodfin2023hitting} and Section~5.2.2 in \citep{collinswoodfin2025Exact}.

Recall we defined $a_{k} = \sqrt{K_{I_{k}}}v_k + \mu_{I_{k}}$ for $v_k \sim \mathcal{N}(0, \Id_d)$, $I_{k} \in \{1,2\}$ is randomly sampled, $r_{k, I_{k+1}} = \hat{X}_k^\top a_{k+1, I_{k+1}}$ and
    \begin{equation}
        \begin{gathered}
        g_{k, I_{k+1}} = r_{k, I_{k+1}} + \begin{bmatrix}
            \delta s_{k, I_{k+1}}\|X_k\| & 0 
        \end{bmatrix}^\top,\\
        s_{k, {I_{k+1}}} = \operatorname{argmax}_{|s| \leq  1}f_{I_{k+1}}(r_{k, I_{k+1}}+ \delta s \|X_k\|; \epsilon_{k+1}).
    \end{gathered}
    \end{equation}
By Taylor expansion, setting $\Delta_k \defas f_{I_{k+1}}'(g_{k, I_{k+1}} )(a_{k+1, I_{k+1}} + \delta s_{k, {I_{k+1}}} \frac{X_k}{\|X_k\|} ) +  \lambda X_k$, 
\begin{equation} 
    \begin{aligned}
    \varphi(X_{k+1}) 
    = 
    \varphi(X_k - \frac{\gamma_k}{d} \Delta_k \big ) 
    &
    =
    \varphi(X_{k}) - \frac{\gamma_k}{d} \ip{\nabla \varphi(X_k), \Delta_k}
    + 
     \frac{\gamma_k^2}{2d^2} \cdot \ip{\nabla^2\varphi(X_k), \Delta_k^{\otimes 2}}
    \end{aligned}
\end{equation}
To derive the Doob decomposition, the idea is to iteratively condition on $r_k$, $I_{k+1}$ and $\hat{X}_k$. For this, we define the $\sigma$-algebras, 
\[
\begin{gathered}
    \mathcal{G}_{k,i} \defas \sigma ( \{\hat{X}_j\}_{j=0}^k, \{r_i\}_{i=0}^k, \{I_j\}_{j=1}^k, , \epsilon_{k+1}, I_{k+1} = i)\quad
    \text{and}\quad \mathcal{F}_{k} \defas \sigma ( \{\hat{X}_j\}_{j=0}^k ).
\end{gathered}
\]
It is clear that $s_{k, I_{k+1}}$ is measurable with respect to  $\mathcal{G}_{k,i}$ such that $s_{k, I_{k+1}}\mid \mathcal{G}_{k,i}$ corresponds to a deterministic maximization problem.

\paragraph{Gradient term in Taylor expansion.}
 First, we start by decomposing the gradient term in \eqref{eq:taylor_theorem},
 \begin{equation} \label{eq:grad}
     \frac{\gamma_k}{d} \ip{\nabla \varphi(X_k),  f_{I_{k+1}}'(g_{k, I_{k+1}} )(a_{k+1, I_{k+1}} + \delta s_{k, {I_{k+1}}} \tfrac{X_k}{\|X_k\|} ) +  \lambda X_k}.
 \end{equation} 
We define a martingale increment associated with it as
\begin{equation}
    \begin{aligned}
      \Delta \mathcal{M}^{\text{Grad}}_k(\varphi) 
      &
      \defas 
      \frac{\gamma_k}{d} f_{I_{k+1}}'(g_{k, I_{k+1}} )\ip{\nabla \varphi (\hat{X}_k),  a_{k+1, I_{k+1}} + \delta s_{k, {I_{k+1}}} \tfrac{X_k}{\|X_k\|} }  
      \\&\qquad- 
      \frac{\gamma_k}{d} \EE_{a,I, \epsilon} \big [ f_{I_{k+1}}'(g_{k, I_{k+1}} )\ip{ \nabla \varphi (\hat{X}_k), a_{k+1, I_{k+1}} + \delta s_{k, {I_{k+1}}} \tfrac{X_k}{\|X_k\|}  } \, | \, \mathcal{F}_k \big ].
    \end{aligned}
\end{equation}
From Danskin's theorem~\ref{thm:danskin}, Assumption~\ref{assumption:risk} and interchanging the derivative and expectation, it follows that 
\[
    \EE_{a,I, \epsilon} \big [ f_{I_{k+1}}'(g_{k, I_{k+1}} )\ip{ \nabla \varphi (\hat{X}_k), a_{k+1, I_{k+1}} + \delta s_{k, {I_{k+1}}} \tfrac{X_k}{\|X_k\|}  } \, | \, \mathcal{F}_k \big ] =  \ip{ \nabla \varphi (\hat{X}_k), \nabla \mathcal{R}^{{\operatorname{adv}}}(X_k) }. 
\]
Hence, the doob decomposition of the gradient term in \eqref{eq:grad} is given by
\begin{equation} \label{eq:grad_term}
\begin{aligned}
    \frac{\gamma_k}{d}& \ip{\nabla \varphi(X_k),  f_{I_{k+1}}'(g_{k, I_{k+1}} )(a_{k+1, I_{k+1}} + \delta s_{k, {I_{k+1}}} \tfrac{X_k}{\|X_k\|} ) +  \lambda X_k}
    \\
    &
    =  
    \frac{\gamma_k}{d} \ip{\nabla \varphi (\hat{X}_k), \nabla \mathcal{R}^{{\operatorname{adv}}}(X_k) + \lambda X_k}
    + \Delta \mathcal{M}_k^{\grad},
\end{aligned}
\end{equation}
where 
\begin{equation}
    \begin{aligned}
        \Delta \mathcal{M}^{\text{Grad}}_k(\varphi) 
        &
        \defas 
        \frac{\gamma_k}{d} f_{I_{k+1}}'(g_{k, I_{k+1}} )\ip{\nabla \varphi (\hat{X}_k),  a_{k+1, I_{k+1}} + \delta s_{k, {I_{k+1}}} \tfrac{X_k}{\|X_k\|} }  
        \\&\qquad- 
        \frac{\gamma_k}{d} \EE_{a,I, \epsilon} \big [ f_{I_{k+1}}'(g_{k, I_{k+1}} )\ip{ \nabla \varphi (\hat{X}_k), a_{k+1, I_{k+1}} + \delta s_{k, {I_{k+1}}} \tfrac{X_k}{\|X_k\|}  } \, | \, \mathcal{F}_k \big ].
    \end{aligned}
\end{equation}
\paragraph{Hessian term in the Taylor expansion.}  Next, we continue by computing the Doob decomposition of the Hessian term in the Taylor expansion~\eqref{eq:taylor_theorem}. Developing the outer product, we have
\begin{equation} \label{eq:Hessian_taylor_term}
\begin{aligned}
    \frac{\gamma_k^2}{2d^2} &\ip{ \nabla^2 \varphi(X_k), \big ( f_{I_{k+1}}'(g_{k, I_{k+1}} )(a_{k+1, I_{k+1}} + \delta s_{k, {I_{k+1}}} \tfrac{X_k}{\|X_k\|} ) +  \lambda X_k\big )^{\otimes 2} } 
    \\
    &=\frac{\gamma_k^2}{2d^2}( f_{I_{k+1}}'(g_{k, I_{k+1}} ))^{ 2} \ip{ \nabla^2 \varphi (\hat{X}_k),  a_{k+1, I_{k+1}}^{\otimes 2} }  
    \\
    &
    \qquad
    + \frac{\gamma_k^2}{d^2} \delta s_{k, I_{k+1}} ( f_{I_{k+1}}'(g_{k, I_{k+1}} ))^{2}\ip{ \nabla^2 \varphi (\hat{X}_k),     a_{k+1, I_{k+1}} \otimes \tfrac{X_k}{\|X_k\|}    }
    \\
    &
    \qquad
    +\frac{\gamma_k^2}{2d^2}\delta^2 s_{k, I_{k+1}}^2 ( f_{I_{k+1}}'(g_{k, I_{k+1}} ))^{2} \ip{ \nabla^2 \varphi (\hat{X}_k),   \tfrac{X_k^{\otimes 2}}{\|X_k\|^2}}
    +
    \frac{\gamma_k^2}{2d^2} \ip{ \nabla^2 \varphi(X_k),  ( \lambda X_k )^{\otimes 2} } 
    \\
    &
    \qquad
    +
    \frac{\gamma_k^2}{d^2}   f_{I_{k+1}}'(g_{k, I_{k+1}} )\ip{ \nabla^2 \varphi(X_k), (a_{k+1, I_{k+1}} + \delta s_{k, I_{k+1}} \tfrac{X_k}{\|X_k\|}) \otimes \lambda X_k   }.
\end{aligned}
\end{equation}
The martingale increment associated with the Hessian is given by
\begin{equation} \label{eq:Hessian_Martingale}
    \begin{aligned}
        \Delta \mathcal{M}_k^{\text{Hess}}(\varphi) \defas \frac{\gamma_k^2}{2d^2} \bigg ( \ip{\nabla^2 \varphi(X_k), \Delta_k^{\otimes 2} } - \EE[ \ip{\nabla^2 \varphi(X_k), \Delta_k^{\otimes 2}} \, | \, \mathcal{F}_k] \bigg ) .
    \end{aligned}
\end{equation}
In order to explicitly characterize the martingale increment, we will compute the conditional expectation of \eqref{eq:Hessian_taylor_term} given $\mathcal{F}_k$. We first condition on $\mathcal{G}_k$ and use the following Lemma~\ref{lem:conditioning} to simplify the computation of the conditional expectation of the first term in \eqref{eq:Hessian_taylor_term}. 
\begin{lemma}[Conditioning, Lemma A2 \citep{collinswoodfin2023hitting}] \label{lem:conditioning} Suppose $v \in \R^d$ is distributed $N(0, I_d)$ and $U \in \R^{d\times 2}$ has orthonormal columns. Then 
\begin{equation}
 v \, | \, \ip{U, v}_{\R^d} \sim v- UU^\top v +UU^\top v ,
 \end{equation}
where $v- UU^\top v  \sim \mathcal{N}(0, \Id_d - UU^\top)$ and $UU^\top v \sim \mathcal{N}(0, UU^\top)$ with $v - UU^\top v$ independent of $UU^\top v$.
\end{lemma}
A simple computation yields
\begin{equation} \label{eq:hessian_1}
    \begin{aligned}
        \EE[ &( f_{i}'(g_{k, i} ))^{ 2} \ip{ \nabla^2 \varphi (\hat{X}_k),  a_{k+1, i}^{\otimes 2} }  \, | \, \mathcal{G}_{k,i} ] 
        \\
        &
        =
        \EE [ ( f_{i}'(g_{k, i} ))^{ 2} \ip{ \nabla^2 \varphi(X_k), (a_{k+1, i} -\EE[a_{k+1, i} \, | \, \mathcal{G}_{k,i}])^{\otimes 2}  } \, | \, \mathcal{G}_{k,i}]
        \\
        &
        \qquad +
        \EE_{\epsilon_{k+1}}[( f_{i}'(g_{k, i} ))^{ 2} ]\cdot \ip{ \nabla^2 \varphi(X_k), \EE[a_{k+1, i}  \, | \, \mathcal{G}_{k,i}]^{\otimes 2}}  .
    \end{aligned}
\end{equation}
We will compute these conditional expectations using Lemma~\ref{lem:conditioning}. By Assumption~\ref{assumption:data}, we write $a_{{k+1}, i} = \mu_i + \sqrt{K_i} v_{k+1}$ for $v_{k+1} \sim N(0, I_d)$. It then follows
\[
    v_{k+1} \mid r_{k,i}, \hat{X}_k \overset{d}{=}v_{k+1} \mid \hat{X}_k^\top(\mu_i + \sqrt{K_i} v_{k+1}), \hat{X}_k\overset{d}{=} v_{k+1}\mid\hat{X}_k^\top \sqrt{K_i} v_{k+1}.
\]
In order to apply Lemma~\ref{lem:conditioning}, we apply the QR decomposition to $ \sqrt{K}\hat{X}_k =Q_{k,i}R_{k,i}$ such that $R_{k,i}\in \R^{\ell_i\times \ell_i}$ is upper-triangular and invertible, $Q_{k,i} \in \R^{d\times \ell_i}$ has orthonormal columns and where $\ell_i = \operatorname{rank}( \sqrt{K_i}W_k ) \leq 2$.

Set $\Pi_{k,i} = Q_{k,i}Q_{k,i}^\top$. By Lemma~\ref{lem:conditioning} and since $R_{k,i}$ is invertible, in distribution we have
\[ 
\begin{aligned}
v_{k+1} \, | \,\hat{X}_k^\top \sqrt{K_i} v_{k+1} 
&
\overset{\text{d}}{=} v_{k+1} \, | \,R_{k,i}^\top Q_{k,i}^\top v_{k+1}\overset{\text{d}}{=}  v_{k+1} \, | \, Q_{k,i}^\top v_{k+1}
\\
&\overset{\text{d}}{=}  \big ( v_{k+1} - \Pi_{k,i} v_{k+1} \big ) + \Pi_{k,i} v_{k+1}.
\end{aligned}
\]  
Hence, it follows that 
\begin{equation}
    a_{k+1, i} \mid r_{k,i}, \hat{X}_k \overset{d}{=} \sqrt{K_i}(  v_{k+1} - \Pi_{k,i} v_{k+1} ) +  \sqrt{K_i}\Pi_{k,i} v_{k+1} + \mu_i 
\end{equation}
By Lemma~\ref{lem:conditioning}, we have $(I_d - \Pi_{k,i}) v_{k+1}\sim N(0, I_{d}  - \Pi_{k,i} )$ and $\Pi_{k,i} v_{k+1} \sim N(0, \Pi_{k,i})$ such that $(I_d - \Pi_{k,i})v_{k+1}$ independent of $\Pi_{k,i} v_{k+1}$ from which we obtain
\begin{equation} \label{eq:mean_a}
\EE[ a_{k+1, i} \, | \, \mathcal{G}_{k,i}]=\sqrt{K_i} \Pi_{k,i} v_{k+1} + \mu_i, \quad \text{where $v_{k+1} \sim N(0, I_d)$.}
\end{equation}
 Furthermore, it follows from a simple computation that
\begin{equation}
    \begin{gathered}
        \EE[ ( a_{k+1, i} - \EE[a_{k+1, i} \, | \, \mathcal{G}_{k,i} ] )^{\otimes 2} | \, \mathcal{G}_k]= \sqrt{K_i} (I_d - \Pi_{k,i}) \sqrt{K_i}.
    \end{gathered}
\end{equation}
Putting this all together, we obtain the decomposition of the leading order term of \eqref{eq:hessian_1} 
\begin{equation} \label{eq:hessian_2}
    \begin{aligned}
    \EE[ &( f_{i}'(g_{k, i} ))^{ 2} \ip{ \nabla^2 \varphi (\hat{X}_k),  a_{k+1, i}^{\otimes 2} }  \, | \, \mathcal{G}_{k,i} ] 
    = 
    \EE_{\epsilon_{k+1}}[(f_{i}'(g_{k, i} ))^{ 2}]\cdot\bigg(  \ip{ \nabla^2 \varphi(X_k), K_i  } 
    \\
    &
    \qquad
    -
    \ip{\nabla^2 \varphi(X_k), \sqrt{K_i} \Pi_{k,i} \sqrt{K_i} } 
    + 
   \ip{ \nabla^2 \varphi (\hat{X}_k), \big ( \sqrt{K_i}\Pi_{k,i} v_{k+1} + \mu_i \big )^{\otimes 2} }\bigg).
\end{aligned}
\end{equation}
 With this in mind, we obtain
\begin{equation}
    \begin{aligned}
        \frac{\gamma_k^2}{2d^2} \EE [ (f_{i}'(g_{k, i} ))^{ 2}  &\ip{ \nabla^2 \varphi(X_k), a_{k+1} ^{\otimes 2}  }  \, | \, \mathcal{F}_k ] 
        \\
        &
        = 
        \frac{\gamma_k^2}{2d^2} \sum_{i=1}^2 p_i \EE[(f_{i}'(g_{k, i} ))^{ 2}   \, | \, \mathcal{F}_k ]\ip{\nabla^2 \varphi(X_k), K_i + \mu_i\mu_i^\top}  + \EE[ \mathcal{E}_{k,1}^{\text{Hess}} \, | \, \mathcal{F}_{k}].
    \end{aligned}
\end{equation}
where
\begin{align*}
    \mathcal{E}_{k,1}^{\text{Hess}}(\varphi) 
    &
    \defas
    \frac{\gamma_k^2}{2d^2}\sum_{i=1}^2 p_i\EE_{\epsilon_{k+1}} [(f_{i}'(g_{k, i} ))^{ 2} ]\cdot \bigg(-\ip{\nabla^2 \varphi(X_k), \sqrt{K_i} \Pi_{k,i}\sqrt{K_i}  } 
\\
&
\qquad
+
   \ip{ \nabla^2 \varphi (\hat{X}_k), \big ( \sqrt{K_i} \Pi_{k,i} v_{k+1} \big )^{\otimes 2}  }
   +
   2\ip{ \nabla^2 \varphi (\hat{X}_k), \sqrt{K_i} \Pi_{k,i} v_{k+1}\mu_i^\top }
   \bigg).
\end{align*}
We will show in Section~\ref{sec:lower_term_Hessian} that $ \mathcal{E}_{k,1}^{\text{Hess}}(\varphi)$ is of lower order and will disappear as $d \to \infty$. For the second term term in \eqref{eq:Hessian_taylor_term}, by \eqref{eq:mean_a} and a similar argument as above, we obtain
\begin{equation}
\begin{aligned}
\EE[ & \delta s_{k, I_{k+1}} ( f_{I_{k+1}}'(g_{k, I_{k+1}} ))^{2}\ip{ \nabla^2 \varphi (\hat{X}_k),     a_{k+1, I_{k+1}} \otimes \tfrac{X_k}{\|X_k\|}    } \, | \, \mathcal{F}_k ] \\
&
= 
\sum_{i=1}^2 p_i
\EE [\delta s_{k, i} ( f_{i}'(g_{k, i} ))^{2}\ip{ \nabla^2 \varphi (\hat{X}_k),     (\sqrt{K_i} \Pi_{k,i} v_{k+1} + \mu_i) \otimes \tfrac{X_k}{\|X_k\|}    }  \, | \, \mathcal{F}_k],
\end{aligned}
\end{equation}
with $v_{k+1} \sim N(0,I_d)$. Due to the factor $\tfrac{1}{d^2}$, this lower order term will disappear as $d \to \infty$ (see Section~\ref{sec:lower_term_Hessian}). Hence, we define
\begin{equation}
    \mathcal{E}_{k,2}^{\text{Hess}}(\varphi) \defas \frac{\gamma_k^2}{d^2} \sum_{i=1}^2 p_i\delta s_{k, i} ( f_{i}'(g_{k, i} ))^{2}\ip{ \nabla^2 \varphi (\hat{X}_k),     (\sqrt{K_i} \Pi_{k,i} v_{k+1} + \mu_i) \otimes \tfrac{X_k}{\|X_k\|}    } 
\end{equation}
For the third term and fourth terms in \eqref{eq:Hessian_taylor_term}, we leave them as is and define
\[
\begin{aligned}
    \mathcal{E}_{k,3}^{\text{Hess}}(\varphi) &
    \defas 
    \frac{\gamma_k^2}{2d^2}\sum_{i=1}^2 p_i \delta^2 s_{k, i}^2 ( f_{i}'(g_{k, i} ))^{2} \ip{ \nabla^2 \varphi (\hat{X}_k),   \tfrac{X_k^{\otimes 2}}{\|X_k\|^2}} ,
    \\
    \mathcal{E}_{k,4}^{\text{Hess}}(\varphi)
    & \defas 
    \frac{\gamma_k^2}{2d^2} \ip{ \nabla^2 \varphi(X_k), (\lambda X_k)^{\otimes 2}}.
\end{aligned}
\]
For the last term in \eqref{eq:Hessian_taylor_term}, by \eqref{eq:mean_a} and a similar argument as above, we obtain
\begin{equation}
\begin{aligned}
\EE[ &  f_{I_{k+1}}'(g_{k, I_{k+1}} )\ip{ \nabla^2 \varphi(X_k), (a_{k+1, I_{k+1}} + \delta s_{k, I_{k+1}} \tfrac{X_k}{\|X_k\|}) \otimes \lambda X_k   } \, | \, \mathcal{F}_k ] \\
&
= \sum_{i=1}^2 p_i
\EE [ f_{i}'(g_{k, i} )\ip{ \nabla^2 \varphi(X_k), (\sqrt{K_i} \Pi_{k,i} v_{k+1} + \mu_i + \delta s_{k, i} \tfrac{X_k}{\|X_k\|}) \otimes \lambda X_k   }\, | \, \mathcal{F}_k].
\end{aligned}
\end{equation}
 Finally, due to the $\tfrac{1}{d^2}$, this term will be of lower order (see Section~\ref{sec:lower_term_Hessian}), and thus, we define
\begin{equation}
    \mathcal{E}_{k,5}^{\text{Hess}}(\varphi) \defas \frac{\gamma_k^2}{d^2} \sum_{i=1}^2 p_i f_{i}'(g_{k, i} )\ip{ \nabla^2 \varphi(X_k), (\sqrt{K_i} \Pi_{k,i} v_{k+1} + \mu_i + \delta s_{k, i} \tfrac{X_k}{\|X_k\|}) \otimes \lambda X_k   }. 
\end{equation}

We have thus identified the martingale increments of a \textit{single} update of SGD, that is
\begin{equation} \label{eq: taylor_expansion_result}
    \begin{aligned}
        \varphi(\hat{X}_{k+1}) 
        & 
        =
        \varphi(X_k) - \frac{\gamma_k}{d} \ip{\nabla \varphi(X_k), \nabla \mathcal{R}^{{\operatorname{adv}}}_\lambda(X_k)}
        \\
        &
        + 
        \frac{\gamma_k^2}{2d^2} \sum_{i=1}^2 p_i \ip{\nabla^2 \varphi(X_k), K_i + \mu_i \mu_i^\top }\EE[(f_{i}'(g_{k, i} ))^2\, | \, \mathcal{F}_k ]
        \\ 
        & + \Delta \mathcal{M}_k^{\text{Grad}}(\varphi) + \Delta \mathcal{M}_k^{\text{Hess}}(\varphi) + \EE[ \mathcal{E}_k^{\text{Hess}}(\varphi) \, |\, \mathcal{F}_k] 
        % + \mathcal{E}_k^{\text{Higher}},
    \end{aligned}
\end{equation}
where the error terms look like
\begin{equation} \label{eq:error_terms}
    \begin{aligned}
    \Delta \mathcal{M}_k^{\grad}(\varphi) 
    &
    = \frac{\gamma_k}{d} f_{I_{k+1}}'(g_{k, I_{k+1}} )\ip{\nabla \varphi (\hat{X}_k),  a_{k+1, I_{k+1}} + \delta s_{k, {I_{k+1}}} \tfrac{X_k}{\|X_k\|} }  
    \\&\qquad- 
    \frac{\gamma_k}{d} \EE_{a,I, \epsilon} \big [ f_{I_{k+1}}'(g_{k, I_{k+1}} )\ip{ \nabla \varphi (\hat{X}_k), a_{k+1, I_{k+1}} + \delta s_{k, {I_{k+1}}} \tfrac{X_k}{\|X_k\|}  } \, | \, \mathcal{F}_k \big ],
    \\
   \Delta \mathcal{M}_k^{\text{Hess}}(\varphi)
    &
    = 
    \frac{\gamma_k^2}{2d^2} \bigg ( \ip{\nabla^2 \varphi(X_k), \Delta_k^{\otimes 2} } - \EE[ \ip{\nabla^2 \varphi(X_k), \Delta_k^{\otimes 2}} \, | \, \mathcal{F}_k] \bigg ) ,
    \\
    \text{and} \quad \mathcal{E}_k^{\text{Hess}}(\varphi)
    &
    = 
    \mathcal{E}_{k,1}^{\text{Hess}}(\varphi) +\mathcal{E}_{k,2}^{\text{Hess}}(\varphi) + \mathcal{E}_{k,3}^{\text{Hess}}(\varphi)+ \mathcal{E}_{k,4}^{\text{Hess}}(\varphi) + \mathcal{E}_{k,5}^{\text{Hess}}(\varphi)
    \\
    &
    = \frac{\gamma_k^2}{d^2}\sum_{i=1}^2 p_i\EE_{\epsilon_{k+1}} [(f_{i}'(g_{k, i} ))^{ 2} ]\cdot \bigg(-\ip{\nabla^2 \varphi(X_k), \sqrt{K_i} \Pi_{k,i}\sqrt{K_i}  } 
    \\
    &
    +
       \ip{ \nabla^2 \varphi (\hat{X}_k), \big ( \sqrt{K_i} \Pi_{k,i} v_{k+1} \big )^{\otimes 2}  }
       +
       2\ip{ \nabla^2 \varphi (\hat{X}_k), \sqrt{K_i} \Pi_{k,i} v_{k+1}\mu_i^\top }
       \bigg)
       \\
       &
       +
       \frac{\gamma_k^2}{d^2} \sum_{i=1}^2 p_i\delta s_{k, i} ( f_{i}'(g_{k, i} ))^{2}\ip{ \nabla^2 \varphi (\hat{X}_k),     (\sqrt{K_i} \Pi_{k,i} v_{k+1} + \mu_i) \otimes \tfrac{X_k}{\|X_k\|}    } 
       \\
        &
        +
    \frac{\gamma_k^2}{2d^2}\sum_{i=1}^2 p_i \delta^2 s_{k, i}^2 ( f_{i}'(g_{k, i} ))^{2} \ip{ \nabla^2 \varphi (\hat{X}_k),   \tfrac{X_k^{\otimes 2}}{\|X_k\|^2}} 
    \\
    & 
    +
    \frac{\gamma_k^2}{2d^2} \ip{ \nabla^2 \varphi(X_k), (\lambda X_k)^{\otimes 2}}
    \\
    &
    +
    \frac{\gamma_k^2}{d^2} \sum_{i=1}^2 p_i f_{i}'(g_{k, i} )\ip{ \nabla^2 \varphi(X_k), (\sqrt{K_i} \Pi_{k,i} v_{k+1} + \mu_i + \delta s_{k, i} \tfrac{X_k}{\|X_k\|}) \otimes \lambda X_k   }
    .
\end{aligned}
\end{equation}
The Doob decomposition in the hard label setting can be recovered using similar arguments or simply be setting $X^{\star} = 0$.
\subsection{Error Bounds} \label{sec:error_bounds}
The results in this section build upon Appendix A in \citep{collinswoodfin2023hitting} and Section 5.5 in \citep{collinswoodfin2025Exact}. Recall, we defined $W = \hat{X} \oplus \mu$. We are studying the statistics
\begin{equation}\label{eq:def:stats:S:M:Z}
\begin{gathered}
    S(W,\cdot) = \hat{X}^\top \mathscr{R}(\cdot;\mathcal{K})\hat{X},\quad  M(W,\cdot) = \hat{X}^\top\mathscr{R}(\cdot;\mathcal{K})\mu,
    \\
   \text{and}\quad Z(W, \cdot) \defas  W^\top \mathscr{R}(\cdot;\mathcal{K})W = \begin{bmatrix}
        S(W, \cdot) &  M(W, \cdot) \\
        M(W, \cdot)^\top & \mu^\top \mathscr{R}(\cdot;\mathcal{K})\mu
    \end{bmatrix}.
\end{gathered}
\end{equation}
Recall we defined the fixed contour
\[
\Gamma = \{ w \in \mathbb{C} \, : \, |w| = \max\{1, 2 \|\mathcal{K}\|_{\operatorname{op}} \} \}\quad \text{and}\quad \|\mathcal{K}\|_{\opt} \defas \max\{\|K_1\|_{\opt}, \|K_1\|_{\opt}\}.
\]
Note that by construction the distance between the spectra of $K_1$ and $K_2$ and the contour is always greater or equal than $\frac{1}{2}$. Recall for $z \in \Gamma^2$, we write $z =(z_1 ,z_2)$ and when integrating over all $z_1$ simultaneously, we write for any function $f: \mathbb{C}^2 \to \mathbb{C}$
\[
\oint f(z) \text{D}z \defas \frac{-1}{4\pi^2 }\oint_{\Gamma^2} f(z)  \dif z_1 \dif z_2.
\]
Recall we also defined the norm $\|\cdot \|_{\Gamma}$ for a continuous function $H: \mathbb{C}^2 \to \R^{4\times 4}$:
\[
    \|H\|_{\Gamma} = \max_{z\in \Gamma^2}\|H(z)\|.
\]
In the next subsections, we will control the error terms which arise in the Doob decompositions of adversarial homogenized SGD and SGD respectively. This results will then be used in the proofs of SGD and AdvHSGD being approximate solutions to the integro-differential equation \eqref{eq:ODE_resolvent_2}. Before diving into the proofs, we introduce Lemma 3 from \citep{collinswoodfin2025Exact} on the norm of the derivatives of $Z(W, \cdot)$.
\begin{lemma}\label{lem:S_derivative_bounds} 
There exist constants $c, C > 0$ such that
\[
\begin{gathered}
c \|W\|^2 \le \|Z(W,\cdot)\|_{\Gamma} \le C \|W\|^2, \quad \|\nabla_X Z (W,\cdot)\|_{\Gamma} \le C \|W\|,
\\
\sup_{z\in \Gamma^2}\|\nabla^2_X Z(W,z)\|_{\opt} \le C\quad
\text{and} \quad \|\nabla^2_X Z(W,\cdot)\|_{\Gamma} \le C\sqrt{d}. 
\end{gathered}
\]
\end{lemma}
Lemma~\ref{lem:S_derivative_bounds} differs from Lemma 3 from \citep{collinswoodfin2025Exact} and is analogous to Lemma A3 in \citep{collinswoodfin2023hitting} as we are studying single-index models (i.e. $\ell =1, \ell^{\star}=1, \ell^+ = 4$ in \citep{collinswoodfin2025Exact}).

To control the martingale terms, we will work with the stopped processes $W_{td }^{\vartheta} \defas W_{ t \wedge \vartheta d }$ and $\HSGD_t^{\vartheta} \defas \HSGD_{t \wedge \vartheta}$ where we recall the stopping time $\vartheta$ for fixed $M>0$ which we defined in \eqref{eq:main:stopping_time_1}: 
\begin{equation}
    \begin{aligned} \label{eq:stopping_time_1}
        \vartheta_{M} 
        &
        \defas \inf \{ t \ge 0 \, : \, \|W_{td} \|^2 > M \quad \text{or}\quad \widehat{\mathcal{B}}_1(t) \not \in \mathcal{U} \quad \text{or}\quad \widehat{\mathcal{B}}_2(t) \not \in \mathcal{U}\}\\
        \text{or} \quad  \vartheta_{M} 
        &
        \defas \inf \{t \ge 0 \, : \, 
        \|\HSGD_t\|^2  > M \, \, \quad \text{or}\quad \widehat{\mathcal{B}}_1(t) \not \in \mathcal{U} \quad \text{or}\quad \widehat{\mathcal{B}}_2(t) \not \in \mathcal{U}\},     
    \end{aligned}
    \end{equation}
    Here we overload the notation so that $\vartheta_{M} $ can be applied to either the iterates of SGD or AdvHSGD. When clear from context, we drop the superscript $\vartheta$ for notational simplicity. We will first control the errors with the stopped processes and later derive a sufficient condition under which we can remove the stopping time.

 Throughout the proofs in the following subsections, the constant $C$ may depend on the following various bounded quantities: $\bar{\gamma}$, $T$, $\lambda$, $\delta$, $\alpha$, $\|K\|_{\opt}$, and $M$. However, $C$ does not depend on $d$. The value of $C$ may also change from line to line in these proofs and may differ from the one presented in Lemma~\ref{lem:S_derivative_bounds}.

\subsubsection{Adversarial Homogenized SGD Martingale Error} \label{sec:HSGD_martingale}
Given that the elements of $Z(\HSGD_t, z)$ are at most quadratics in $\WHSGD_t$, it suffices to control the martingale increment that arises in the Doob decomposition of AdvHSGD for an arbitrary quadratic statistic $\varphi \, : \, \R^d \to \mathbb{R}$:
\begin{equation}
\begin{aligned}
    \mathcal{M}^{\operatorname{AdvHSGD}}_{t}(\varphi)
    &
    \defas \int_0^t \dif \mathcal{M}_s^{\operatorname{AdvHSGD}}(\varphi)  
    \\
    &
    = \frac{1}{\sqrt{d}}\int_0^t
    \gamma(s)   \left\langle \sqrt{\sum_{i=1}^2 p_i \EE_{v, \epsilon} [  f'_i( \rho_{s, i}, \epsilon)^2  ](K_i + \mu_i\mu_i^\top)}, \nabla \varphi(\WHSGD_s)\otimes \dif B_s \right\rangle.
\end{aligned}
\end{equation}

To show this martingale increment vanishes as $d \to \infty$, we will control its fluctuations by showing elementwise, its quadratic variation as defined below vanishes. Given a partition $\pi =\{t_0, ...,t_n\}$ of $[0,t]$ such that $0=t_0< 
 t_1<\cdots <t_n=t$ then the $\pi$-quadratic variation of $\mathcal{M}_t^{\operatorname{AdvHSGD}}(Z_{j\ell}(\cdot, z))$ is defined as 
 \[
    Q_{\pi}( \mathcal{M}_t^{\operatorname{AdvHSGD}}(Z_{j\ell}(\cdot, z))) = \sum_{i=0}^{n-1}(\mathcal{M}_{t_{i+1}}^{\operatorname{AdvHSGD}}(Z_{j\ell}(\cdot, z))-\mathcal{M}_{t_i}^{\operatorname{AdvHSGD}}(Z_{j\ell}(\cdot, z)))^2.
 \]
 If there exists process denoted $\langle \mathcal{M}^{\operatorname{AdvHSGD}}(Z_{j\ell}(\cdot, z))\rangle_t$ such that for all sequence of partitions $\pi_\ell$ with $\max_{0\leq i\leq n_\ell-1} |t_{i+1}^{(\ell)} -t_{i}^{(\ell)} | \underset{k\to \infty}{\to} 0$, we have
\[
    Q_{\pi_\ell}(\mathcal{M}_t^{\operatorname{AdvHSGD}}(Z_{j\ell}(\cdot, z))) \underset{\ell\to \infty}{\to} \langle \mathcal{M}_t^{\operatorname{AdvHSGD}}(Z_{j\ell}(\cdot, z))\rangle_t,
\]
then $\langle \mathcal{M}^{\operatorname{AdvHSGD}}(Z_{j\ell}(\cdot, z))\rangle$ is called the \emph{quadratic variation} of $\langle \mathcal{M}^{\operatorname{AdvHSGD}}(Z_{j\ell}(\cdot, z))\rangle$ (see \cite[Chapter 1]{protter2005stochastic} for details). We build upon Proposition A.1 of \citep{collinswoodfin2023hitting}. We provide the proof for the soft label setting as the proof for hard labels follows the same argument with $X^{\star} = 0$.
\begin{proposition}[AdvHSGD small.] \label{prop:HSGD_Martingale_Bound} For $i=1,2$, suppose $f_i \, : \,\R^{3} \to \mathbb{R}$ is $\alpha$-pseudo-Lipschitz function with constant $L(f)$ (see Assumption~\ref{assumption:pseudo_lipschitz}). Let the statistic $Z \, : \, \R^d \to \R^{4 \times 4}$ be defined as in \eqref{eq:def:stats:S:M:Z}. For any $T > 0$,  $\zeta > 0$ and fix $z \in \Gamma^2$, there is some constant $C$ such that, with overwhelming probability, 
\begin{equation}
    \sup_{0 \le t \le T} \| \mathcal{M}^{\operatorname{AdvHSGD}}_{t \wedge \vartheta} (S(\cdot,z)) \| \le C L(f) \, \, d^{\zeta / 2 - 1/2}.
\end{equation}
\end{proposition}
\begin{proof} Let $\varphi(\WHSGD_t) \defas Z_{j\ell}(\HSGD_t,z)$ be the $j\ell$-coordinate of $Z$ for a fixed $z \in \Gamma^2$. Recall the martingale increment elementwise is given by,
    \begin{equation}
        \dif \mathcal{M}_{t}^{\operatorname{AdvHSGD}}(\varphi) =  \frac{\gamma(t)}{\sqrt{d}}   \left\langle \sqrt{ \sum_{i=1}^2 p_i \EE_{v, \epsilon} [  f'_i( \rho_{t, i}, \epsilon)^2  ](K_i + \mu_i\mu_i^\top)}, \nabla_X \varphi(\WHSGD_t) \otimes \dif B_t \right\rangle
    \end{equation}
    The quadratic variation of $\mathcal{M}_t^{\operatorname{AdvHSGD}}$ is 
    \begin{equation} \label{eq:quad_variation_HSGD}
        \begin{aligned}
       \langle &\mathcal{M}^{\operatorname{AdvHSGD}}(\varphi) \rangle_t
       \\
       &
        = \frac{1}{d} \int_0^t \gamma^2(s)  \left\|\left \langle \sqrt{ \sum_{i=1}^2 p_i \EE_{v, \epsilon} [  f'_i( \rho_{s, i}, \epsilon)^2  ](K_i + \mu_i\mu_i^\top)} , \nabla_X \varphi(\WHSGD_s) \right\rangle_{\R^{d}}\right\|^2 \, \dif s. 
        \end{aligned}
    \end{equation}
    Our goal is to show $\displaystyle \sup_{0 \le t \le T} \langle\mathcal{M}^{\operatorname{AdvHSGD}}(\varphi)\rangle_{t \wedge \vartheta}$ is small. Since $\gamma:\R^+ \to \R^+$ is a bounded function by Assumption~\ref{assumption:lr}, we must then show the $\|\cdot\|^2$ is also small. We have 
    \begin{equation}
        \begin{aligned}
            & 
            \left\|\left \langle \sqrt{ \sum_{i=1}^2 p_i \EE_{v, \epsilon} [  f'_i( \rho_{s, i}^{\vartheta}, \epsilon)^2  ](K_i + \mu_i\mu_i^\top)} , \nabla_X \varphi(\WHSGD_s^{\vartheta}) \right\rangle_{\R^{d}}\right\|^2
            \\
            &
          \qquad \leq  
           \left\|\sqrt{\sum_{i=1}^2 p_i \EE_{v, \epsilon} [  f'_i( \rho_{s, i}^{\vartheta}, \epsilon)^2  ](K_i + \mu_i\mu_i^\top)}\right\|_{\opt}^2 \|\nabla_X \varphi(\WHSGD_s^{\vartheta}) \|^2
           \\
           \qquad 
           &
           \qquad
           \leq  
           \sum_{i=1}^2 p_i  \EE_{v, \epsilon} [  f'_i( \rho_{s, i}^{\vartheta}, \epsilon)^2  ]\left\|(K_i + \mu_i\mu_i^\top)\right\|_{\opt} \|\nabla_X \varphi(\WHSGD_s^{\vartheta}) \|^2
           \\
           &
           \qquad
           \leq  
           \sum_{i=1}^2 p_i  \EE_{v, \epsilon} [  f'_i( \rho_{s, i}^{\vartheta}, \epsilon)^2  ] (\|K_i\|_{\opt} +\|\mu_i\|^2 )\|\nabla_X \varphi(\WHSGD_s^{\vartheta}) \|^2.
        \end{aligned}
    \end{equation} 
    By Lemma~\ref{lem:S_derivative_bounds}, we have a bound on 
    \[
        \|\nabla_X \varphi(\WHSGD_s^{\vartheta}) \| = \|\nabla_X Z_{j\ell}(\HSGD_s^{\vartheta},z)\| \le \|\nabla_X Z_{j\ell}(\HSGD_s^{\vartheta},\cdot)\|_{\Gamma} \le C\|\HSGD_s^{\vartheta}\|.
    \] 
    From Assumption~\ref{assumption:data} for $i=1,2$, we have $\|K_i\|_{\operatorname{op}} \leq C$, $p_i\leq 1$ and $\|\mu_i\|^2 \leq C$. From Lemma~\ref{lem:growth_grad_f}, we have a bound on $\EE_{v, \epsilon} [  f'_i( \rho_{s, i}^{\vartheta}, \epsilon)^2  ]$, and since $\|\WHSGD_t^{\vartheta}\| \leq \|\widehat{\WHSGD}_t^{\vartheta}\| \leq \|\widehat{\HSGD}_t^{\vartheta}\|$ then
    \begin{equation} \label{eq:quad_variation_HSGD_1}
        \begin{aligned}
            & \left\|\left \langle \sqrt{ \sum_{i=1}^2 p_i \EE_{v, \epsilon} [  f'_i( \rho_{s, i}^{\vartheta}, \epsilon)^2  ](K_i + \mu_i\mu_i^\top)} , \nabla_X \varphi(\WHSGD_s^{\vartheta}) \right\rangle_{\R^{d}}\right\|^2
            \leq 
            C \|\HSGD_s^{\vartheta}\|^2\sum_{i=1}^2 \EE_{v, \epsilon} [  f'_i( \rho_{s, i}^{\vartheta}, \epsilon)^2  ] 
            \\
            &
            \qquad
            \leq
            C  (L(f))^2  \|\HSGD_s^{\vartheta}\|^2 \sum_{i=1}^2 \big  (1+(\|K_i\|_{\opt}^{1/2} +\|\mu_i\|)\|\widehat{\HSGD}_t^\nu \| + \delta \|\mathscr{X}_t^\nu\|+ \eta )^{\max \{1, 2\alpha \} }
            \\
            &
            \qquad
            \leq
            C  (L(f))^2  M\big  (1+(1+ \delta) \sqrt{M}+ \eta   )^{\max \{1, 2\alpha \} }
        \end{aligned}
    \end{equation}
    Thus, from \eqref{eq:quad_variation_HSGD} and \eqref{eq:quad_variation_HSGD_1}, it follows
    \begin{equation}
        \sup_{0 \le t \le T} \langle \mathcal{M}^{\operatorname{AdvHSGD}}(\varphi)\rangle_{t \wedge \vartheta} \le C  (L(f))^2  d^{- 1}.
    \end{equation}
    From Doob's maximal inequality for non-negative super-martingales and the exponential super-martingale 
    \[
    Z^a_t
    =
    e^{a \mathcal{M}_{t \wedge \vartheta}^{\operatorname{AdvHSGD}}(\varphi) - \frac{a^2}{2} \langle \mathcal{M}^{\operatorname{AdvHSGD}}(\varphi)\rangle_{t \wedge \vartheta} }, \quad \text{for}\quad a,b >0,
    \]
    if $\displaystyle \sup_{0 \le t \le T}  \langle \mathcal{M}^{\operatorname{AdvHSGD}}(\varphi)\rangle_{t \wedge \vartheta}  \le b$ a.s, then it is easy to show
     \[
     \Pr (\sup_{0 \le t \le T} | \mathcal{M}_{t \wedge \vartheta}^{\operatorname{AdvHSGD}}(\varphi) | > p ) \le 2\exp(-p^2 / 2b).
     \]
    Absorbing $\bar{\gamma}^2$ into $C$ and setting $p = \sqrt{2C} L(f) d^{\zeta/2 - 1/2}$ for any $\zeta > 0$, then
    \[
    \Pr( \sup_{0 \le t \le T} | \mathcal{M}_{t \wedge \vartheta}^{\operatorname{AdvHSGD}}(\varphi) | > \sqrt{C} L(f) d^{\zeta/2 - 1/2} ) \le C \exp( - d^{\zeta} ).
    \]
    The result immediately follows after noting that the number of $j\ell$ coordinates is $16$ which is independent of $d$. Indeed, from a union bound we have
    \begin{equation}
        \begin{aligned}
        \Pr(\sup_{1 \leq i,j \leq 4}\sup_{0 \le t \le T} | \mathcal{M}_{t \wedge \vartheta}^{\operatorname{AdvHSGD}}(\varphi) | < p) &= 1 -  \Pr(\sup_{1 \leq i,j \leq 4}\sup_{0 \le t \le T} | \mathcal{M}_{t \wedge \vartheta}^{\operatorname{AdvHSGD}}(\varphi) | > p)\\
        &= 1 - \Pr(\cup_{1\leq i,j \leq 4}\sup_{0 \le t \le T} | \mathcal{M}_{t \wedge \vartheta}^{\operatorname{AdvHSGD}}(\varphi) | > p)\\
        &\geq  1 - C\exp( - d^{\zeta} ).
        \end{aligned}
    \end{equation}
Noting that
\begin{equation}
\begin{aligned}
    \sup_{0 \le t \le T} \| \mathcal{M}^{\operatorname{AdvHSGD}}_{t \wedge \vartheta} (S(\cdot,z)) \| \leq  4\sup_{1 \leq i,j \leq 2} \sup_{0 \le t \le T} | \mathcal{M}^{\operatorname{AdvHSGD}}_{t \wedge \vartheta} (\varphi) |,
\end{aligned}
\end{equation}
 concludes the proof.
\end{proof}

\subsubsection{Bounds on the Martingales $\cM_k^{\grad}$ and $\cM_k^{\hess}$ \label{subsec:SGD_martingales}}
Given that the elements of $Z(W_k, z)$ are at most quadratics in $X_k$, it suffices to control the martingale increment that arise in the Doob decomposition of SGD for an arbitrary quadratic statistic $\varphi \, : \, \R^d \to \mathbb{R}$:
\begin{equation*}
    \begin{aligned}
        \Delta \mathcal{M}_k^{\grad}(\varphi) 
        &
        = \frac{\gamma_k}{d} f_{I_{k+1}}'(g_{k, I_{k+1}} )\ip{\nabla \varphi (\hat{X}_k),  a_{k+1, I_{k+1}} + \delta s_{k, {I_{k+1}}} \tfrac{X_k}{\|X_k\|} }  
        \\&\qquad- 
        \frac{\gamma_k}{d}\sum_{i=1}^2 p_i \EE_{a,\epsilon} \big [ f_{i}'(g_{k, i} )\ip{ \nabla \varphi (\hat{X}_k), a_{k+1, i} + \delta s_{k, i} \tfrac{X_k}{\|X_k\|}  } \, | \, \mathcal{F}_k \big ],
        \\
       \Delta \mathcal{M}_k^{\text{Hess}}(\varphi)
        &
        = 
        \frac{\gamma_k^2}{2d^2} \bigg ( \ip{\nabla^2 \varphi(X_k), \Delta_k^{\otimes 2} } - \EE[ \ip{\nabla^2 \varphi(X_k), \Delta_k^{\otimes 2}} \, | \, \mathcal{F}_k] \bigg ),
    \end{aligned}
\end{equation*}
where we recall $a_{k, I_{k}} = \sqrt{K_{I_{k}}}v_k + \mu_{I_{k}}$ for $v_k \sim \mathcal{N}(0, \Id_d)$, $I_{k} \in \{1,2\}$, $r_{k, I_{k+1}} = \hat{X}_k^\top a_{k+1, I_{k+1}}$ and
    \begin{equation}
        \begin{gathered}
        g_{k, I_{k+1}} = r_{k, I_{k+1}} + \begin{bmatrix}
            \delta s_{k, I_{k+1}}\|X_k\| & 0 
        \end{bmatrix}^\top,\\
        s_{k, {I_{k+1}}} = \operatorname{argmax}_{|s| \leq  1}f_{I_{k+1}}( X_k^\top a_{k+1, I_{k+1}}+ \delta s \|X_k\|; \epsilon_{k+1}).
    \end{gathered}
    \end{equation}
We also define 
\[
 \mathcal{M}^{\grad}_{k}(\varphi) = \sum_{j=0}^{k-1} \Delta \mathcal{M}^{\grad}_j(\varphi)\quad \text{and}\quad   \mathcal{M}^{\text{Hess}}_{k}(\varphi) = \sum_{j=0}^{k-1} \Delta \mathcal{M}^{\hess}_j(\varphi)
\]
Once again, we provide the proofs for the soft label setting as the proofs for hard labels follow the same arguments with $X^{\star} = 0$. In the following proof of Proposition~\ref{prop:gradmartingale}, we build upon Proposition A.2 of \citep{collinswoodfin2023hitting} and Lemma 9 in \citep{collinswoodfin2025Exact}.
\begin{proposition}[Gradient martingale]\label{prop:gradmartingale}
    For $i=1,2$, suppose $f_i \, : \,\R^{3} \to \mathbb{R}$ is $\alpha$-pseudo-Lipschitz function with constant $L(f)$ (see Assumption~\ref{assumption:pseudo_lipschitz}). Let the statistic $Z \, : \, \R^d \to \R^{4 \times 4}$ be defined as in \eqref{eq:def:stats:S:M:Z}. Then, for any $\zeta>0$ and $T > 0$, and with overwhelming probability,
    \beq
    \sup_{0 \leq t\leq T} \|\cM_{\lfloor(t \wedge \vartheta)d\rfloor}^{\grad}(Z(\cdot, z))\|<d^{-\frac12+\zeta}.
    \eeq
    \end{proposition}
\begin{proof}
    Let $\varphi \defas Z_{j\ell}(W,z)$ be the $j\ell$-coordinate of $Z$ for a fixed $z \in \Gamma^2$. Similarly as in the proof of Proposition~\ref{prop:HSGD_Martingale_Bound}, we will be working on the stopped process $\cM_{(t \wedge \vartheta)d}^{\grad}(\varphi )$. To simplify notation, we omit the dependence on the stopping time $\vartheta$ and on $\varphi $ such that $\cM_{k}^{\grad} \defas \cM_{(k/d \wedge \vartheta)d}^{\grad}(\varphi )$. We also omit the dependence on $X$ in $\nabla \varphi = \nabla_X \varphi$.
    
    Define the projected martingale increment
    \begin{equation}
        \begin{aligned}
        &\Delta \mathcal{M}_k^{\grad, \beta} 
    = \frac{\gamma_k}{d}  f_{I_{k+1}}'(\text{Proj}_{\beta}(g_{k, I_{k+1}} )) \cdot\big( \text{Proj}_{\beta} (\ip{\nabla \varphi (\hat{X}_k), \sqrt{K_i} v_{k+1}}) + \ip{\nabla \varphi (\hat{X}_k), \mu_i }
    \\
    &
    + \delta s_{k, {I_{k+1}}} \ip{\nabla \varphi (\hat{X}_k), \tfrac{X_k}{\|X_k\|} }  \big )
    - 
    \frac{\gamma_k}{d} \sum_{i=1}^2 p_i \EE \big [  f_{i}'(\text{Proj}_{\beta}(g_{k, i} )) \cdot  (\text{Proj}_{\beta} (\ip{\nabla \varphi (\hat{X}_k),  \sqrt{K_i}v_{k+1}}) 
    \\
    &
    +\ip{\nabla \varphi (\hat{X}_k), \mu_i } + \delta s_{k, {I_{k+1}}} \ip{\nabla \varphi (\hat{X}_k), \tfrac{X_k}{\|X_k\|} } )  \, | \, \mathcal{F}_k \big ]
        \end{aligned}
    \end{equation}
For some $\zeta >0$, set $\beta = d^{\zeta}$. We will later optimize the value of $\zeta$. Our goal will be to bound $\Delta \mathcal{M}_k^{\grad, \beta} $ and then the difference $|\Delta \mathcal{M}_k^{\grad}  - \Delta \mathcal{M}_k^{\grad, \beta} |$. 

To bound the difference, we will make use of subgaussian tail bounds of the projected quantities taking advantage of the fact that $a_{k,i} = \sqrt{K_i} v_{k+1} + \mu_i$ for $v_{k+1}$ and $v_{k+1}$ is subgaussian. Denoting $\widehat{r}_{k,i} \defas \hat{X}_k^\top\sqrt{K_i}v_{k+1}$ then $\widehat{r}_{k,i}$ is $\|\sqrt{K}\|_{\operatorname{op}}\|\hat{X}_k\|$-subgaussian and $\ip{\sqrt{K}\nabla \varphi (X_k),  v_{k+1} }$ is $\|\sqrt{K}\|_{\operatorname{op}}\|\nabla \varphi (X_k)\|$-subgaussian.

Given that $|s_{k,i}| \leq 1$ and noting that $(g_{k,i}, \epsilon_{k+1}) = (\widehat{r}_{k,i}, \epsilon_{k+1}) + (\hat{X}_k^\top \mu_i, 0) + (\delta s_{k,i} \|X_k\|, 0_{2}) \in \R^3$,  we have by the triangle inequality and Cauchy--Schwarz
\begin{equation}\label{eq:bounds:for:tails:grad}
    \begin{aligned}
   \|(g_{k,i}, \epsilon_{k+1}) \| &\leq  \| (\hat{X}_k^\top\sqrt{K_i}v_{k+1}, \epsilon_{k+1})\| + \| \hat{X}^\top\mu_i\| + \delta \|X_k\|.
   %\\|\langle \nabla \varphi (X_k),  a_{k+1, i} + \delta s_{k,i} \tfrac{X_k}{\|X_k\|} \rangle| &\leq  |\ip{\sqrt{K_i}\nabla \varphi (X_k),  v_{k+1} }| +|\ip{\nabla \varphi (X_k),  \mu_i}|+\delta \| \nabla \varphi (X_k)\| .
    \end{aligned}
\end{equation}

From Lemma~\ref{lem:S_derivative_bounds}, we have
\begin{equation}\label{eq:bound:grad:varphi}
\|\nabla \varphi(X_k)\| = \|\nabla  Z_{j\ell}(W_k,z)\| \leq \sup_{z \in \Gamma}\|\nabla_X  Z_{j\ell}(W_k,z)\|  \leq \|\nabla_X  Z(W_k,z)\|_{\Gamma} \leq C\|W_k\|.
\end{equation}
Now,  By Assumption~\ref{assumption:data}, $\epsilon_{k+1}$ is $\eta$-subgaussian. Since $\eta$ is a constant, $\|\sqrt{K_i}\|_{\operatorname{op}} \leq C$, by \eqref{eq:bound:grad:varphi} and we are working with the stopped processes, then for large enough $C$ it follows
\begin{equation}\label{eq:grad:subgaussian:bounds}
    \|(\widehat{r}_{k,i} , \epsilon_{k+1})\|_{\psi_2} \leq C,\quad\text{and}\quad \|\ip{\sqrt{K_i}\nabla \varphi (X_k),  v_{k+1} }\|_{\psi_2} \leq C.
\end{equation}
Here we omit the dependence of the constant $C$ on $i$ as we could always redefine it as $C = \max\{C_1, C_2\}$. We will proceed with the same logic throughout the rest of the proof.

Hence, using \eqref{eq:bounds:for:tails:grad} and for the stopped processes $\delta\|X_k\| \leq \delta\|W_k\| \leq CM$, $\|\hat{X}_k^\top \mu_i \| \leq \|W_k\| \|\mu_i\| \leq CM$. We then may apply the following argument to bound the tails of the projected quantities for $\beta > CM$
\begin{equation}
    \begin{aligned}
    \Proba(\|(g_{k,i}, \epsilon_{k+1})\|> \beta) 
    &
    \leq \Proba(\| (r_k, \epsilon_{k+1})\| > \beta - CM)\leq 2\exp\left(-\frac{(\beta - CM)^2}{2C}\right).
    \end{aligned}
\end{equation}
 Since we set $\beta = d^{\zeta}$ for $\zeta >0$, the condition $\beta > CM$ holds for sufficiently large $d$ hence it follows that 
 \begin{equation}
    \begin{gathered}\label{eq:tail:bounds:proj:grad}
    \Proba(\|(g_{k,i}, \epsilon_{k+1})\|> \beta) \leq  \exp(-\Omega(\beta^2)),\\
    %\Proba(|\langle \nabla \varphi (X_k),  a_{k+1, i} + \delta s_{k,i} \tfrac{X_k}{\|X_k\|} \rangle|> \beta) \leq \exp(-\Omega(\beta^2)).
    \end{gathered}
 \end{equation}
 This tail bounds will be useful when we bound the difference between $\mathcal{M}_k^{\grad}$ and $\mathcal{M}_k^{\grad, \beta}$.

 From the definition of the projection and the bounds on the absolute value of $f'_i$ in Lemma~\ref{lem:lipschitz:proj} and Lemma~\ref{lem:growth_grad_f}, we have
\begin{equation}\label{eq:grad:bounds:for:secterm}
    \begin{gathered}
        |f_i'(\text{Proj}_\beta(g_{k,i}, \epsilon_{k+1}))| \leq CL(f)\beta^{\max(1,\alpha)},\\
        |\text{Proj}_{\beta} (\ip{\nabla \varphi (X_k),  \sqrt{K_i} v_{k+1} }  )| \leq \beta.
    \end{gathered}
\end{equation}
It then follows from Cauchy--Schwarz
\begin{equation}\label{eq:grad:bound:F_k^beta}
    |f_i'(\text{Proj}_\beta(g_{k,i}, \epsilon_{k+1})) \cdot \text{Proj}_{\beta} (\ip{\nabla \varphi (X_k),  \sqrt{K_i} v_{k+1} }  )| \leq CL(f)\beta^{2 + \alpha}.
\end{equation}
Since we are working with the stopped process and $|s_{k,i}|\leq 1$, we also have from Cauchy--Schwarz and \eqref{eq:bound:grad:varphi} 
\begin{equation}\label{eq:grad:bound:F_k^beta:1}
    \begin{gathered}
    |f_i'(\text{Proj}_\beta(g_{k,i}, \epsilon_{k+1})) \cdot \ip{\nabla \varphi (\hat{X}_k), \mu_i }   | \leq CML(f)\beta^{1+\alpha}, \\
    |f_i'(\text{Proj}_\beta(g_{k,i}, \epsilon_{k+1})) \cdot \delta s_{k, {I_{k+1}}} \ip{\nabla \varphi (\hat{X}_k), \tfrac{X_k}{\|X_k\|} }   | \leq CML(f)\beta^{1+\alpha}.
    \end{gathered}
\end{equation}
Since this is an almost sure bound, it also holds under expectation from which it follows that 
\begin{equation}\label{eq:grad:bound:F_k^beta:good}
   | \Delta \mathcal{M}_k^{\grad, \beta} | \leq \frac{2\bar{\gamma}}{d}CL(f)\beta^{2 + \alpha},\quad a.s.
\end{equation}
By Azuma's inequality then for $0 \leq k \leq Td$
\begin{equation}
   \Proba\left( | \mathcal{M}_{k}^{\grad, \beta} | > t\right)\leq 2 \exp\left(\frac{-t^2}{2\sum_{i=1}^k Cd^{-2}\beta^{2(2 + \alpha)}}\right)\leq 2 \exp\left(\frac{-t^2}{TCd^{-1}\beta^{2(2 + \alpha)}}\right)
\end{equation}
Thus, with overwhelming probability, for every $0 \leq k \leq Td$
\begin{equation}
    |\mathcal{M}_{k}^{\grad, \beta}| < d^{-1/2}\beta^{3+\alpha}.
\end{equation}
Taking a union bound from $0$ to $Td$ does not modify the overwhelming probability statement. Hence, with overwhelming probability
\begin{equation}\label{eq:grad:azuma:union}
    \sup_{0 \leq k \leq Td}|\mathcal{M}_{k}^{\grad, \beta}| < d^{-1/2}\beta^{3+\alpha}.
\end{equation}
It remains to bound the difference between $ \mathcal{M}_k^{\grad}$ and $\mathcal{M}_k^{\grad, \beta}$. For notational convenience, we write 
\begin{equation}
    \begin{aligned}
        F_{k,i}&\defas\frac{\gamma_k}{d} f_{i}'(g_{k, i} )\ip{\nabla \varphi (\hat{X}_k),  \sqrt{K_i}v_{k+1} + \mu_i+ \delta s_{k, {i}} \tfrac{X_k}{\|X_k\|} }  ,\\
        F_{k,i}^\beta &\defas\frac{\gamma_k}{d}  f_{i}'(\text{Proj}_{\beta}(g_{k, i} )) \cdot \bigg(\text{Proj}_{\beta} (\ip{\nabla \varphi (\hat{X}_k),  \sqrt{K_i}v_{k+1}}) \\
        &+
        \ip{\nabla \varphi (\hat{X}_k), \mu_i } + \delta s_{k, {i}} \ip{\nabla \varphi (\hat{X}_k), \tfrac{X_k}{\|X_k\|} } \bigg)  .
    \end{aligned}
\end{equation}
Our goal is to show the quantity 
\begin{equation}
    |(F_{k,i} - \E[F_{k,i}]) - (F^{\beta}_{k,i} - \E[F^{\beta}_{k,i}])| \leq |F_{k,i} -F^{\beta}_{k,i}| + |\E [F_{k,i} -F^{\beta}_{k,i}]|,
\end{equation}
vanishes as $d\to\infty$. Notice that 
\begin{equation*}
    \begin{aligned}
        &|F_{k,i} -F^{\beta}_{k,i}| 
        \\
        &
        \leq \frac{\gamma_k}{d} |f_{i}'(g_{k, i} )\ip{\nabla \varphi (\hat{X}_k),  \sqrt{K}v_{k+1 }} - f_{i}'(\text{Proj}_{\beta}(g_{k, i} )) \cdot \text{Proj}_{\beta} (\ip{\nabla \varphi (\hat{X}_k),  \sqrt{K}v_{k+1 }})| 
        \\
        &
        +
        \frac{\gamma_k}{d} |\ip{\nabla \varphi (\hat{X}_k),\mu_i}|\cdot | f_{i}'(g_{k, i} )  -  f_{i}'(\text{Proj}_{\beta}(g_{k, i} )) |
        \\
        &
        +
        \frac{\gamma_k}{d}\delta  |\ip{\nabla \varphi (\hat{X}_k),\tfrac{X_k}{\|X_k\|}}|\cdot | f_{i}'(g_{k, i} )  -  f_{i}'(\text{Proj}_{\beta}(g_{k, i} )) |
    \end{aligned}
\end{equation*}
Thus, using the subgaussian tail bounds we previously derived in \eqref{eq:tail:bounds:proj:grad} and \eqref{eq:tail:bounds:proj:grad}, we obtain 
\begin{equation}\label{eq:grad:diff:F_k}
    \begin{aligned}
    \Proba(|F_{k,{I_{k+1}}} - F^{\beta}_{k,I_{k+1}}| > 0) 
    &
    = 
    \Proba(\left\{\|(g_{k,I_{k+1}}, \epsilon_{k+1})\| > \beta\right\} \cup \{|\ip{\nabla \varphi (\hat{X}_k),  \sqrt{K_{I_{k+1}}}v_{k+1 }}| > \beta\})
    \\
    &
    \leq \sum_{i=1}^2 p_i( \Proba(\|(g_{k,i}, \epsilon_{k+1})\| > \beta) + \Proba( |\ip{\nabla \varphi (\hat{X}_k),  \sqrt{K_i}v_{k+1 }}| > \beta))
    \\
    &\leq 4 \exp\left(-\Omega(\beta^2)\right).
    \end{aligned}
\end{equation}
Since we defined $\beta = d^{\zeta}$ for some $\zeta>0$, then \eqref{eq:grad:diff:F_k} implies $F_k - F^{\beta}_k = 0$ with overwhelming probability. It then follows from a union bound that 
\begin{equation}
    \sum_{k=1}^{Td} |F_{k,{I_{k+1}}} - F^{\beta}_{k,I_{k+1}}| = 0,
\end{equation}
with overwhelming probability. We must now bound the difference in expectation. In the following expectations, we omit the conditional expectation on $\mathcal{F}_k$ for notational simplicity. We have by Cauchy--Schwarz
\begin{equation}
    \begin{aligned}
       | \E[F_{k,i} - F_{k,i}^{\beta}]| &
        = 
        |\E[(F_{k,i} - F_{k,i}^{\beta})\mathds{1}_{|F_{k,i} -F_{k,i}^{\beta}|>0}]|
        \\
        &
        \leq |\E[F_{k,i} \mathds{1}_{|F_{k,i} -F_{k,i}^{\beta}|>0}] | + |\E[F_{k,i}^{\beta}\mathds{1}_{|F_{k,i} -F_{k,i}^{\beta}|>0}]|.
    \end{aligned}
\end{equation}
Using a similar argument as in (A.78)-(A.82) in \citep{collinswoodfin2023hitting}, we will show both of the terms in the upper bound are exponentially small with high-probability. First, using \eqref{eq:grad:bound:F_k^beta} and \eqref{eq:grad:bound:F_k^beta:1} we have
\begin{equation}
    \begin{aligned}
        |\E[F_{k,i}^{\beta}\mathds{1}_{|F_{k,i} -F_{k,i}^{\beta}|>0}]| 
        &
        \leq \frac{\bar{\gamma}}{d}CL(f)\beta^{2 + \alpha} \Proba\left(|F_{k,i} -F_{k,i}^{\beta}|>0\right)
        \\
        &
         \leq  \frac{\bar{\gamma}}{d}CL(f)\beta^{2 + \alpha} \exp\left(-\Omega(\beta^2)\right).
    \end{aligned}
\end{equation}
Then for the second term we have 
\begin{equation}
\begin{split}
|\EE [F_{k,i}\cdot \mathds{1}_{|F_{k,i} -F_{k,i}^{\beta}|>0}]|&\leq
\EE [ |F_{k,i}| \cdot \mathds{1}_{E_1} ]+\EE [ |F_{k,i} |\cdot \mathds{1}_{E_2} ]+\EE [ |F_{k,i}| \cdot \mathds{1}_{E_3} ],\\
\text{where }\;E_1\defas & \{\|(g_{k,i}, \epsilon_{k+1})\| \leq\beta \}\cap\{|\ip{\nabla \varphi (\hat{X}_k),  \sqrt{K_i}v_{k+1 }}|>\beta \},\\
E_2\defas &\{\|(g_{k,i}, \epsilon_{k+1})\| >\beta \}\cap\{|\ip{\nabla \varphi (\hat{X}_k),  \sqrt{K_i}v_{k+1 }}|\leq\beta \},\\
E_3\defas &\{\|(g_{k,i}, \epsilon_{k+1})\| >\beta \}\cap\{|\ip{\nabla \varphi (\hat{X}_k),  \sqrt{K_i}v_{k+1 }}|>\beta \}.\\
\end{split}
\end{equation}
From the bounds in \eqref{eq:grad:bounds:for:secterm} and \eqref{eq:grad:bound:F_k^beta:1}, we obtain for the term involving $E_1$ 
\begin{equation}
    \begin{aligned}
        \EE [ |F_{k,i}| \cdot \mathds{1}_{E_1} ] &\leq \frac{\bar{\gamma}}{d} CML(f)\beta^{1+\alpha} \Proba(|\ip{\nabla \varphi (\hat{X}_k),  \sqrt{K_i}v_{k+1 }}|>\beta)
        \\
        &
        + \frac{\bar{\gamma}}{d} CL(f)\beta^{1+\alpha}\E\left[|\ip{\nabla \varphi (\hat{X}_k),  \sqrt{K_i}v_{k+1 }}|\mathds{1}_{|\ip{\nabla \varphi (\hat{X}_k),  \sqrt{K_i}v_{k+1 }}|>\beta }\right].
    \end{aligned}
\end{equation}
We will show $\E\left[|\ip{\nabla \varphi (\hat{X}_k),  \sqrt{K_i}v_{k+1 }}|\mathds{1}_{|\ip{\nabla \varphi (\hat{X}_k),  \sqrt{K_i}v_{k+1 }}|>\beta }\right]$ is exponentially small. Using \eqref{eq:tail:bounds:proj:grad}, we have for some constant $c>0$
\begin{equation}
    \begin{aligned}
        \E\bigg[&\ip{\nabla \varphi (\hat{X}_k),  \sqrt{K_i}v_{k+1 }} |\mathds{1}_{|\ip{\nabla \varphi (\hat{X}_k),  \sqrt{K_i}v_{k+1 }}|>\beta }\bigg] 
        \\
        &
        = \beta \Proba(|\ip{\nabla \varphi (\hat{X}_k),  \sqrt{K_i}v_{k+1 }}| > \beta)+ \int_{\beta}^\infty \Proba( |\ip{\nabla \varphi (\hat{X}_k),  \sqrt{K_i}v_{k+1 }}| > t) \dif t 
        \\
        &
        \leq 2\beta \exp\left(-\Omega(\beta^2)\right)+
        \int_{\beta}^\infty 2\exp\left(- \frac{t^2}{2C}\right) \dif t 
        \\
        &= 2\beta \exp\left(-\Omega(\beta^2)\right)+ 2\sqrt{2\pi C}  \;\Proba\left(\mathcal{N}(0,1) > \frac{\beta}{\sqrt{C}}\right)
        \\
        &\leq C\beta \exp\left(-\Omega(\beta^2)\right),
    \end{aligned}
\end{equation}
as Gaussian random variables are subgaussian. Since $\beta = d^{\zeta}$ for some $\zeta >0$ then the bound is exponentially small in $d$. Thus, we obtain
\begin{equation}\label{eq:grad:E_1:bound}
    \begin{aligned}
        \EE [ |F_{k,i}| \cdot \mathds{1}_{E_1} ] \leq \frac{\bar{\gamma}}{d}L(f) C\beta^{2+\alpha}\exp\left(-\Omega(\beta^2)\right).
    \end{aligned}
\end{equation}
Since we are working with the stopped processes, from Cauchy--Schwarz, \eqref{eq:tail:bounds:proj:grad}, \eqref{eq:grad:bounds:for:secterm} and the moments bound in Lemma~\ref{lem:growth_grad_f}, we have for the term involving $E_2$
\begin{equation}
    \begin{aligned}
        \EE [ &|F_{k,i}| \cdot \mathds{1}_{E_2} ] 
        \leq \frac{\bar{\gamma}}{d} CM\beta\E\left[|f'_i(g_{k,i}, \epsilon_{k+1})|\mathds{1}_{\|(g_{k,i}, \epsilon_{k+1})\|>\beta }\right]
        \\
        &
        \leq \frac{\bar{\gamma}}{d}CM\beta \left(\E\left[|f'_i(g_{k,i}, \epsilon_{k+1})|^2\right]\right)^{1/2} \Proba(\|(g_{k,i}, \epsilon_{k+1})\|>\beta )^{1/2}
        \\
        &
        \leq \frac{\bar{\gamma}}{d} CML(f)\beta \left(1+(\|K_i\|_{\operatorname{op}}^{1/2}+\|\mu_i\|)\|\hat{X}_k\| +  \delta \|X_k\| + \eta\right)^{1+\alpha} \exp(-\Omega(\beta^2))
        \\
    \end{aligned}
\end{equation}
Since we are working with the stopped processes and from Assumption~\ref{assumption:data}, $\|\hat{X}_k\|$, $\|W_k\|$, $\|K\|_{\operatorname{op}}$, $\eta^2$, $\delta$ are bounded by constants (independent of $d$) so we obtain
\begin{equation}
    \begin{aligned}
        \EE [ |F_{k,i}| \cdot \mathds{1}_{E_2} ] 
        &
        \leq \frac{\bar{\gamma}}{d}CL(f)\beta \exp(-\Omega(\beta^2))
        \\
    \end{aligned}
\end{equation}
For the term involving $E_3$, we use Cauchy--Schwarz to obtain
\begin{equation}
    \begin{aligned}
    \EE [ &|F_{k,i}| \cdot \mathds{1}_{E_3} ] \\
    &
    \leq \frac{\bar{\gamma}}{d} (\E\left[|\ip{\nabla \varphi (\hat{X}_k),  \sqrt{K_i}v_{k+1 }}|^2\mathds{1}_{\|(g_{k,i}, \epsilon_{k+1})\|>\beta }\right])^{1/2}(\E\left[|f'_i(g_{k,i}, \epsilon_{k+1})|^2\mathds{1}_{\|(g_{k,i}, \epsilon_{k+1})\|>\beta }\right])^{1/2}
    \\
    &
    +
    \frac{\bar{\gamma}}{d} CM \E\left[|f'_i(g_{k,i}, \epsilon_{k+1})|\mathds{1}_{\|(g_{k,i}, \epsilon_{k+1})\|>\beta }\right].
    \end{aligned}
\end{equation}
Using similar steps as for $E_1$ and $E_2$, this is also exponentially small. Putting this altogether, we obtain
\begin{equation}
    |\EE [F_{k,i}\cdot \mathds{1}_{|F_{k,i} -F_{k,i}^{\beta}|>0}]| \leq Cd^{-1} \text{poly}(\beta)\exp(-\Omega(\beta^2))
\end{equation}
For simplicity of notation, let $\mathcal{E} = \{\sum_{k=1}^{Td} |F_{k,I_{k+1}} -F^{\beta}_{k,I_{k+1}}| =0\}$ then with $\beta = d^{\zeta}$ we have
\begin{equation}
    \begin{aligned}
    &
    \Proba( \sup_{0\leq k \leq Td} |\mathcal{M}_{k}^{\grad}| > t)
     \\
     &
     \leq  \Proba( \sup_{0\leq k \leq Td} |\mathcal{M}_{k}^{\grad, \beta}| > t -  \sum_{k=1}^{Td}| F_{k, I_{k+1}} -F^{\beta}_{k, I_{k+1}}| -  \sum_{k=1}^{Td}\sum_{i=1}^2 p_i|\E [F_{k,i} -F^{\beta}_{k,i}]|) 
    \\
    &
    \leq  
    \Proba( \sup_{0\leq k \leq Td} |\mathcal{M}_{k}^{\grad, \beta}| > t -   \sum_{k=1}^{Td}\sum_{i=1}^2 p_i|\E [F_{k,i} -F^{\beta}_{k,i}]|, \;\mathcal{E}) + \Proba(\mathcal{E}^c)
    \\
    &
    \leq  
    \Proba( \sup_{0\leq k \leq Td} |\mathcal{M}_{k}^{\grad, \beta}| > t -  \text{poly}(d^{\zeta})\exp(-\Omega(d^{2\zeta}))) + \exp(-\Omega(d^{2\zeta})).
    \end{aligned}
\end{equation}
Setting $t \defas d^{-\frac{1}{2}+(3+\alpha)\zeta }$ then for sufficiently large $d$, $t- \text{poly}(d^{\zeta }) \exp(-\Omega(d^{2\zeta}))) \geq \frac{t}{2}$. Thus, it follows from Azuma's inequality and a union bound as done in \eqref{eq:grad:azuma:union} that, with overwhelming probability
\begin{equation}
    \sup_{0 \leq k \leq Td}|\mathcal{M}_{k}^{\grad}| < d^{-\frac{1}{2}+(3+\alpha)\zeta}.
\end{equation}
Hence for any arbitrarily small value of $\zeta < \frac{1}{2(3+\alpha)}$, the result follows.
\end{proof}
In the following proof of Proposition~\ref{prop:hessmartingale}, we build upon Proposition A.3 of \citep{collinswoodfin2023hitting} and Lemma 10 in \citep{collinswoodfin2025Exact}.
\begin{proposition}[Hessian martingale]\label{prop:hessmartingale}
    For $i=1,2$, suppose $f_i \, : \,\R^{3} \to \mathbb{R}$ is $\alpha$-pseudo-Lipschitz function with constant $L(f)$ (see Assumption~\ref{assumption:pseudo_lipschitz}). Let the statistic $Z \, : \, \R^d \to \R^{4 \times 4}$ be defined as in \eqref{eq:def:stats:S:M:Z}. Then, for any $\zeta>0$ and $T > 0$, and with overwhelming probability,
    \beq
    \sup_{0 \leq t\leq T} \|\cM_{\lfloor(t \wedge \vartheta)d\rfloor}^{\hess}(Z(\cdot, z))\|<d^{-\frac{1}{2}+\zeta}.
    \eeq
    \end{proposition}
\begin{proof}
    Similarly to in the proof of Proposition~\ref{prop:gradmartingale}, we denote $\varphi (\hat{X}_k)\defas Z_{j\ell}(W,z)$ as the $j\ell$-coordinate of $Z$ for a fixed $z \in \Gamma^2$. We will also work with the stopped process $\cM_{(t \wedge \vartheta)d}^{\hess}(\varphi )$ and, for notational simplicity, omit the dependence on the stopping time $\vartheta$ and on $\varphi$ such that $\cM_{k}^{\hess} \defas \cM_{(k/d \wedge \vartheta)d}^{\hess}(\varphi )$. We also omit the dependence on $X$ in $\nabla^2 \varphi = \nabla_X^2 \varphi$.
    
    For some $\zeta >0$, set $\beta = d^{\zeta}$. We will later optimize the value of $\zeta$ similarly to the proof of Proposition~\ref{prop:gradmartingale}. Recall the hessian martingale increment
\begin{equation}
    \begin{aligned}
    \Delta \mathcal{M}_k^{\text{Hess}}
    &
    = 
    \frac{\gamma_k^2}{2d^2} \bigg ( \ip{\nabla^2 \varphi(X_k), \Delta_k^{\otimes 2}} - \EE[ \ip{\nabla^2 \varphi(X_k), \Delta_k^{\otimes 2}} \, | \, \mathcal{F}_k] \bigg ),
    \end{aligned}
\end{equation}
where $\Delta_k = f_{I_{k+1}}'(g_{k, I_{k+1}} )(a_{k+1, I_{k+1}} + \delta s_{k, {I_{k+1}}} \frac{X_k}{\|X_k\|} ) +  \lambda X_k$. We first split this term into nine separate increments
\begin{equation}
    \Delta \mathcal{M}_k^{\text{Hess}} = \sum_{q=1}^9\Delta \mathcal{M}_k^{\text{Hess},q},
\end{equation}
where
\begin{equation}
    \begin{aligned}
        \Delta \mathcal{M}_k^{\text{Hess},1} 
        &= \frac{\gamma_k^2}{2d^2} \bigg (   (f_{I_{k+1}}'(g_{k, I_{k+1}} ))^2 \ip{\nabla^2 \varphi(X_k), (\sqrt{K_{I_{k+1}}}v_{k+1})^{\otimes 2}  }
        \\
        &
        \qquad 
        - \sum_{i=1}^2 p_i \EE[ (f_{i}'(g_{k, i} ))^2 \ip{\nabla^2 \varphi(X_k), (\sqrt{K_{i}}v_{k+1})^{\otimes 2}  } \, | \, \mathcal{F}_k] \bigg )
        \\
        \Delta \mathcal{M}_k^{\text{Hess},2} 
        &= \frac{\gamma_k^2}{d^2} \bigg (   (f_{I_{k+1}}'(g_{k, I_{k+1}} ))^2 \ip{\nabla^2 \varphi(X_k),   \sqrt{K_{I_{k+1}}}v_{k+1}\otimes \mu_{I_{k+1}} }
        \\
        &
        \qquad 
        - \sum_{i=1}^2 p_i \EE[ (f_{i}'(g_{k, i} ))^2 \ip{\nabla^2 \varphi(X_k),\sqrt{K_{I_{k+1}}}v_{k+1}\otimes \mu_{I_{k+1}}} \, | \, \mathcal{F}_k] \bigg )
        \\
        \Delta \mathcal{M}_k^{\text{Hess},3} 
        &= \frac{\gamma_k^2}{2d^2} \bigg (   (f_{I_{k+1}}'(g_{k, I_{k+1}} ))^2 \ip{\nabla^2 \varphi(X_k),  \mu_{I_{k+1}}^{\otimes 2}}
        \\
        &   
        \qquad 
        - \sum_{i=1}^2 p_i \EE[ (f_{i}'(g_{k, i} ))^2 \ip{\nabla^2 \varphi(X_k), \mu_{i}^{\otimes 2}} \, | \, \mathcal{F}_k] \bigg )
        \\
        \Delta \mathcal{M}_k^{\text{Hess},4} 
        &
        = \frac{\gamma_k^2}{d^2} \bigg ( \delta s_{k, {I_{k+1}}}(f_{I_{k+1}}'(g_{k, I_{k+1}} ))^2 \ip{\nabla^2 \varphi(X_k),  \mu_{I_{k+1}}\otimes \frac{X_k}{\|X_k\|} }      
\\
&
\qquad- \sum_{i=1}^2 p_i\EE[\delta s_{k, i}(f_i'(g_{k, i} ))^2 \ip{\nabla^2 \varphi(X_k),  \mu_{i}\otimes \frac{X_k}{\|X_k\|} }    \, | \, \mathcal{F}_k] \bigg )
\\
\Delta \mathcal{M}_k^{\text{Hess},5} 
        &
        = \frac{\gamma_k^2}{d^2} \bigg ( \delta s_{k, {I_{k+1}}}(f_{I_{k+1}}'(g_{k, I_{k+1}} ))^2 \ip{\nabla^2 \varphi(X_k),   \sqrt{K_{I_{k+1}}}v_{k+1}\otimes \frac{X_k}{\|X_k\|} }      
\\
&
\qquad- \sum_{i=1}^2 p_i\EE[\delta s_{k, i}(f_i'(g_{k, i} ))^2 \ip{\nabla^2 \varphi(X_k),   \sqrt{K_i}v_{k+1}\otimes \frac{X_k}{\|X_k\|} }    \, | \, \mathcal{F}_k] \bigg )
\\
\Delta \mathcal{M}_k^{\text{Hess},6} &= \frac{\gamma_k^2}{2d^2} \delta^2 \ip{\nabla^2 \varphi(X_k),  \frac{X_k^{\otimes 2}}{\|X_k\|^2}}  \bigg (  (s_{k, {I_{k+1}}})^2 (f_{I_{k+1}}'(g_{k, I_{k+1}} ))^2 
\\
&
\qquad 
- \sum_{i=1}^2 p_i\EE[  (s_{k, i}^2 (f_{i}'(g_{k, i} ))^2  \, | \, \mathcal{F}_k] \bigg )
\\
\Delta \mathcal{M}_k^{\text{Hess},7}
    &
    = 
    \frac{\gamma_k^2}{d^2} \bigg ( f_{I_{k+1}}'(g_{k, I_{k+1}} )\ip{\nabla^2 \varphi(X_k),  \sqrt{K_{I_{k+1}}}v_{k+1}  \otimes \lambda X_k}
    \\
    &
    \qquad
     - \sum_{i=1}^2 p_i \EE[ f_{i}'(g_{k, i} )\ip{\nabla^2 \varphi(X_k),  \sqrt{K_{i}}v_{k+1}  \otimes \lambda X_k}\, | \, \mathcal{F}_k] \bigg )
\\
\Delta \mathcal{M}_k^{\text{Hess},8}
    &
    = 
    \frac{\gamma_k^2}{d^2} \bigg ( f_{I_{k+1}}'(g_{k, I_{k+1}} )\ip{\nabla^2 \varphi(X_k),  \mu_{I_{k+1}}  \otimes \lambda X_k}
    \\
    &
    \qquad
     - \sum_{i=1}^2 p_i \EE[ f_{i}'(g_{k, i} )\ip{\nabla^2 \varphi(X_k),  \mu_i  \otimes \lambda X_k}\, | \, \mathcal{F}_k] \bigg )
\\
\Delta \mathcal{M}_k^{\text{Hess},9}
    &
    = 
    \frac{\gamma_k^2}{d^2} \delta \lambda \ip{\nabla^2 \varphi(X_k),  \frac{X_k^{\otimes 2}}{\|X_k\|}} \bigg ( s_{k, {I_{k+1}}} f_{I_{k+1}}'(g_{k, I_{k+1}} )
    \\
    &
    \qquad
     - \sum_{i=1}^2 p_i \EE[ s_{k, i} f_{i}'(g_{k,i} )\, | \, \mathcal{F}_k] \bigg )
    \end{aligned}
\end{equation}
As done in the proof of Proposition~\ref{prop:gradmartingale}, we introduce the projected increments
\begin{equation}
    \begin{aligned}
        \Delta \mathcal{M}_k^{\text{Hess},1, \beta} 
        &= \frac{\gamma_k^2}{2d^2} \bigg (   (f_{I_{k+1}}'(\operatorname{Proj}_{\beta}(g_{k, I_{k+1}} )))^2 \operatorname{Proj}_{d\beta}(\ip{\nabla^2 \varphi(X_k), (\sqrt{K_{I_{k+1}}}v_{k+1})^{\otimes 2}  })
        \\
        &
        \qquad 
        - \sum_{i=1}^2 p_i \EE[ (f_{i}'(\operatorname{Proj}_{\beta}(g_{k, i} )))^2 \operatorname{Proj}_{d\beta}(\ip{\nabla^2 \varphi(X_k), (\sqrt{K_{i}}v_{k+1})^{\otimes 2}  }) \, | \, \mathcal{F}_k] \bigg )
        \\
        \Delta \mathcal{M}_k^{\text{Hess},2, \beta} 
        &= \frac{\gamma_k^2}{d^2} \bigg (   (f_{I_{k+1}}'(\operatorname{Proj}_{\beta}(g_{k, I_{k+1}} )))^2 \operatorname{Proj}_{\beta}(\ip{\nabla^2 \varphi(X_k),   \sqrt{K_{I_{k+1}}}v_{k+1} \otimes \mu_{I_{k+1}} })
        \\
        &
        \qquad 
        - \sum_{i=1}^2 p_i \EE[ (f_{i}'(\operatorname{Proj}_{\beta}(g_{k, i} )))^2 \operatorname{Proj}_{\beta}(\ip{\nabla^2 \varphi(X_k), \sqrt{K_{I_{k+1}}}v_{k+1} \otimes \mu_{I_{k+1}} }) \, | \, \mathcal{F}_k] \bigg )
        \\
        \Delta \mathcal{M}_k^{\text{Hess},3, \beta} 
        &= \frac{\gamma_k^2}{2d^2} \bigg (   (f_{I_{k+1}}'(\operatorname{Proj}_{\beta}(g_{k, I_{k+1}} )))^2 \ip{\nabla^2 \varphi(X_k),  \mu_{I_{k+1}}^{\otimes 2}}
        \\
        &   
        \qquad 
        - \sum_{i=1}^2 p_i \EE[ (f_{i}'(\operatorname{Proj}_{\beta}(g_{k, i} )))^2 \ip{\nabla^2 \varphi(X_k), \mu_{i}^{\otimes 2}} \, | \, \mathcal{F}_k] \bigg )
        \\
        \Delta \mathcal{M}_k^{\text{Hess},4,\beta} 
        &
        = \frac{\gamma_k^2}{d^2} \bigg ( \delta s_{k, {I_{k+1}}}(f_{I_{k+1}}'(\operatorname{Proj}_{\beta}(g_{k, I_{k+1}} )))^2 \ip{\nabla^2 \varphi(X_k),  \mu_{I_{k+1}}\otimes \frac{X_k}{\|X_k\|} }      
\\
&
\qquad- \sum_{i=1}^2 p_i\EE[\delta s_{k, i}(f_{i}'(\operatorname{Proj}_{\beta}(g_{k, i} )))^2 \ip{\nabla^2 \varphi(X_k),  \mu_{i}\otimes \frac{X_k}{\|X_k\|} }    \, | \, \mathcal{F}_k] \bigg )
\\
\Delta \mathcal{M}_k^{\text{Hess},5, \beta} 
        &
        = \frac{\gamma_k^2}{d^2} \bigg ( \delta s_{k, {I_{k+1}}}(f_{I_{k+1}}'(\operatorname{Proj}_{\beta}(g_{k, I_{k+1}} )))^2 \operatorname{Proj}_{\beta}( \ip{\nabla^2 \varphi(X_k),   \sqrt{K_{I_{k+1}}}v_{k+1}\otimes \frac{X_k}{\|X_k\|} } )     
\\
&
\qquad- \sum_{i=1}^2 p_i\EE[\delta s_{k, i}(f_{i}'(\operatorname{Proj}_{\beta}(g_{k, i} )))^2 \operatorname{Proj}_{\beta}(\ip{\nabla^2 \varphi(X_k),   \sqrt{K_i}v_{k+1}\otimes \frac{X_k}{\|X_k\|} } )   \, | \, \mathcal{F}_k] \bigg )
\\
\Delta \mathcal{M}_k^{\text{Hess},6,\beta} &= \frac{\gamma_k^2}{2d^2} \delta^2 \ip{\nabla^2 \varphi(X_k),  \frac{X_k^{\otimes 2}}{\|X_k\|^2}}  \bigg (  s_{k, {I_{k+1}}}^2 (f_{I_{k+1}}'(\operatorname{Proj}_{\beta}(g_{k, I_{k+1}} )))^2 
\\
&
\qquad 
- \sum_{i=1}^2 p_i\EE[  (s_{k, i}^2 (f_{i}'(\operatorname{Proj}_{\beta}(g_{k, i} )))^2  \, | \, \mathcal{F}_k] \bigg )
\\
\Delta \mathcal{M}_k^{\text{Hess},7, \beta}
    &
    = 
    \frac{\gamma_k^2}{d^2} \bigg ( f_{I_{k+1}}'(\operatorname{Proj}_{\beta}(g_{k, I_{k+1}} ))\operatorname{Proj}_{\beta}(\ip{\nabla^2 \varphi(X_k),  \sqrt{K_{I_{k+1}}}v_{k+1}  \otimes \lambda X_k})
    \\
    &
    \qquad
     - \sum_{i=1}^2 p_i \EE[ f_{i}'(\operatorname{Proj}_{\beta}(g_{k, i} ))\operatorname{Proj}_{\beta}(\ip{\nabla^2 \varphi(X_k),  \sqrt{K_{i}}v_{k+1}  \otimes \lambda X_k})\, | \, \mathcal{F}_k] \bigg )
\\
\Delta \mathcal{M}_k^{\text{Hess},8, \beta}
    &
    = 
    \frac{\gamma_k^2}{d^2}\lambda \bigg ( f_{I_{k+1}}'(\operatorname{Proj}_{\beta}(g_{k, I_{k+1}} ))\ip{\nabla^2 \varphi(X_k),  \mu_{I_{k+1}}  \otimes  X_k}
    \\
    &
    \qquad
     - \sum_{i=1}^2 p_i \EE[ f_{i}'(\operatorname{Proj}_{\beta}(g_{k, i} ))\ip{\nabla^2 \varphi(X_k),  \mu_i  \otimes  X_k}\, | \, \mathcal{F}_k] \bigg )
\\
\Delta \mathcal{M}_k^{\text{Hess},9 , \beta}
    &
    = 
    \frac{\gamma_k^2}{d^2} \delta \lambda \ip{\nabla^2 \varphi(X_k),  \frac{X_k^{\otimes 2}}{\|X_k\|}} \bigg ( s_{k, {I_{k+1}}} f_{I_{k+1}}'(\operatorname{Proj}_{\beta}(g_{k, I_{k+1}} ))
    \\
    &
    \qquad
     - \sum_{i=1}^2 p_i \EE[ s_{k, i} f_{i}'(\operatorname{Proj}_{\beta}(g_{k, i} ))\, | \, \mathcal{F}_k] \bigg )
    \end{aligned}
\end{equation}
Note that the projection radius of $\Delta \mathcal{M}_k^{\text{Hess},1 , \beta}$ differs from the projection radius $\beta$ due to the order of the fluctuations of $\langle \nabla^2 \varphi(X_k), (\sqrt{K_{i}}v_{k+1})^{\otimes 2}\rangle$. We begin by bounding the terms not involving $v_{k+1}$ which are $\mathcal{M}_k^{\text{Hess},3}$, $\mathcal{M}_k^{\text{Hess},4}$, $\mathcal{M}_k^{\text{Hess},6}$, $\mathcal{M}_k^{\text{Hess},8}$, and $\mathcal{M}_k^{\text{Hess},9}$. From Lemma~\ref{lem:S_derivative_bounds} and a similar reasoning as \eqref{eq:bound:grad:varphi}, we have
\begin{equation}
    \begin{gathered}
 |\ip{\nabla^2 \varphi(X_k),  \mu_i^{\otimes 2}}| \leq C \|\mu_i\|^2,\quad
|\ip{\nabla^2 \varphi(X_k),  \mu_i \otimes \frac{X_k}{\|X_k\|}}| \leq C\|\mu_i\|,
\\
|\ip{\nabla^2 \varphi(X_k),  \mu_i \otimes X_k}| \leq C\|\mu_i\|\|X_k\|,
\quad |\ip{\nabla^2 \varphi(X_k),  \frac{X_k^{\otimes 2}}{\|X_k\|}}| \leq C \|X_k\|,
\\
\text{and}\quad |\ip{\nabla^2 \varphi(X_k),  \frac{X_k^{\otimes 2}}{\|X_k\|^2}}| \leq C.
    \end{gathered}
\end{equation}
Here we omit the dependence of the constant $C$ on $i$ as we could always redefine it as $C = \max\{C_1, C_2\}$. We will proceed with the same logic throughout the rest of the proof.

Since we are working with the stopped processes and $\|\mu_i\|^2 \leq C$, it follows
\begin{equation}
    \begin{gathered}
 |\ip{\nabla^2 \varphi(X_k),  \mu_i^{\otimes 2}}| \leq C ,\quad
|\ip{\nabla^2 \varphi(X_k),  \mu_i \otimes \frac{X_k}{\|X_k\|}}| \leq C,
\\
|\ip{\nabla^2 \varphi(X_k),  \mu_i \otimes X_k}| \leq CM,
\quad |\ip{\nabla^2 \varphi(X_k),  \frac{X_k^{\otimes 2}}{\|X_k\|}}| \leq C M,
\\
\text{and}\quad |\ip{\nabla^2 \varphi(X_k),  \frac{X_k^{\otimes 2}}{\|X_k\|^2}}| \leq C
    \end{gathered}
\end{equation}
Using the bounds on the absolute value of $f'$ in Lemma~\ref{lem:lipschitz:proj} and Lemma~\ref{lem:growth_grad_f}, we have
\begin{equation}\label{eq:hess:bounds:for:secterm}
    \begin{aligned}
        |f_{i}'(\operatorname{Proj}_{\beta}(g_{k, i} ))| 
        &
        \leq CL(f)\beta^{1+\alpha},
        \\
         |(f_{i}'(\operatorname{Proj}_{\beta}(g_{k, i} )))^2|
          &
          \leq (CL(f)\beta^{1+\alpha})^2.
    \end{aligned}
\end{equation}
Since $|s_{k}|\leq 1$, it immediately follows
\begin{equation}\label{eq:hess:bounds:for:secterm:sk}
    \begin{aligned}
        |s_kf_{i}'(\operatorname{Proj}_{\beta}(g_{k, i} ))| 
        &
        \leq CL(f)\beta^{\max(1+\alpha)},
        \\
         |s_k^2 (f_{i}'(\operatorname{Proj}_{\beta}(g_{k, i} )))^2|
          &
          \leq (CL(f)\beta^{1+\alpha})^2.
    \end{aligned}
\end{equation}
These bounds holding almost surely imply they also hold under expectation from which we obtain
\begin{equation}
    \begin{gathered}
        |\Delta \mathcal{M}_k^{\text{Hess},3,\beta} |
        \leq 
        C \frac{\bar{\gamma}^2}{d^2}L(f)^2\beta^{2+2\alpha},
        \quad
         |\Delta \mathcal{M}_k^{\text{Hess},4,\beta} |
         \leq 
         C \frac{\bar{\gamma}^2\delta}{d^2}L(f)^2\beta^{2+2\alpha},
         \\
         |\Delta \mathcal{M}_k^{\text{Hess},6,\beta} |
        \leq 
        C\frac{\bar{\gamma}^2\delta^2}{d^2}L(f)^2\beta^{2+2\alpha},
        \quad
         |\Delta \mathcal{M}_k^{\text{Hess},8,\beta} |
         \leq 
         2 C\frac{\bar{\gamma}^2 \lambda}{d^2}\lambda L(f)\beta^{ 1+ \alpha},
         \\
         \text{and}
         \quad
         |\Delta \mathcal{M}_k^{\text{Hess},9,\beta} |
        \leq C\frac{\bar{\gamma}^2}{d^2}\lambda \delta L(f)\beta^{1+\alpha}.
    \end{gathered}
\end{equation}
Using a similar argument involving Azuma's inequality as in the proof of Proposition~\ref{prop:gradmartingale}, we obtain with overwhelming probability
\begin{equation}\label{eq:hess:sup:boundM:35}
    \begin{gathered}
        \sup_{0 \leq k \leq Td}|\mathcal{M}_{k}^{\hess,3, \beta}| 
        < 
        d^{-3/2} \beta^{3+2\alpha},
        \quad
         \sup_{0 \leq k \leq Td}|\mathcal{M}_{k}^{\hess,4, \beta}| 
         < 
         d^{-3/2}\beta^{3+2\alpha},
        \\
        \sup_{0 \leq k \leq Td}|\mathcal{M}_{k}^{\hess,6, \beta}| 
        < 
        d^{-3/2} \beta^{3+2\alpha},
        \quad
         \sup_{0 \leq k \leq Td}|\mathcal{M}_{k}^{\hess,8, \beta}| 
         < 
         d^{-3/2}\beta^{2+\alpha},
        \\
        \text{and}
        \quad 
        \sup_{0 \leq k \leq Td}|\mathcal{M}_{k}^{\hess,9, \beta}| 
        < 
        d^{-3/2} \beta^{2+\alpha}.
    \end{gathered}
\end{equation}
Using similar steps as in the proof of Proposition~\ref{prop:gradmartingale}, it then follows that the differences $|\mathcal{M}_k^{\text{Hess},q} - \mathcal{M}_k^{\text{Hess},q,\beta}|$ for $q\in \{3,4,6,8,9\}$ are exponentially small with overwhelming probability. 

We now continue by bounding the terms linear in $v_{k+1}$ which are $\mathcal{M}_k^{\text{Hess},2,\beta}$, $\mathcal{M}_k^{\text{Hess},5,\beta}$ and $\mathcal{M}_k^{\text{Hess},7,\beta}$. The argument is once again extremely similar to the proof of Proposition~\ref{prop:gradmartingale}. Setting $\widehat{r}_{k,i} \defas \hat{X}_k \sqrt{K_i}v_{k+1}$ and using a similar argument as \eqref{eq:grad:subgaussian:bounds} and from Lemma~\ref{lem:S_derivative_bounds}, for large enough constant $C$ we have 
\begin{equation}\label{eq:hess:subgaussian:bounds}
    \begin{gathered}
    \|(\widehat{r}_{k,i}, \epsilon_{k+1})\|_{\psi_2} \leq C,\quad \| \ip{\sqrt{K_i}\nabla^2 \varphi(X_k)\mu_i ,  v_{k+1}}\|_{\psi_2} \leq C
    \\
    \| \ip{\sqrt{K_i}\nabla^2 \varphi(X_k)\frac{X_k}{\|X_k\|} ,  v_{k+1}}\|_{\psi_2} \leq C,\quad
   \text{and}\quad \| \lambda\ip{\sqrt{K_i}\nabla^2 \varphi(X_k) X_k,  v_{k+1}}\|_{\psi_2} \leq C.
    \end{gathered}
\end{equation}
By \eqref{eq:hess:bounds:for:secterm}, \eqref{eq:hess:bounds:for:secterm:sk}, the definition of the projection operator and using a similar argument as in \eqref{eq:grad:bound:F_k^beta} in the proof of Proposition~\ref{prop:gradmartingale}, we have
\begin{equation}
    \begin{gathered}
       |\Delta \mathcal{M}_k^{\text{Hess},2,\beta}| \leq \frac{\bar{\gamma}^2}{d^2} C(L(f)) ^2 \beta^{3 + 2\alpha}\quad,| \Delta \mathcal{M}_k^{\text{Hess},5,\beta}| \leq \delta\frac{\bar{\gamma}^2}{d^2}  C L(f) \beta^{3 + 2\alpha},
       \\
       \text{and}\quad \Delta \mathcal{M}_k^{\text{Hess},7,\beta}| \leq \frac{\bar{\gamma}^2}{d^2}  C L(f) \beta^{2 + \alpha}.
    \end{gathered}
\end{equation}
from which it follows
\begin{equation}\label{eq:hess:sup:boundM:24}
    \begin{gathered}
        \sup_{0 \leq k \leq Td}|\mathcal{M}_{k}^{\hess,2, \beta}| < d^{-3/2} \beta^{4+2\alpha},\quad \sup_{0 \leq k \leq Td}|\mathcal{M}_{k}^{\hess,5, \beta}| < d^{-3/2}\beta^{4+2\alpha},
        \\
        \text{and}\quad \sup_{0 \leq k \leq Td}|\mathcal{M}_{k}^{\hess,7, \beta}| < d^{-3/2}\beta^{3+\alpha}.
    \end{gathered}
\end{equation}
From the subgaussian norm bounds \eqref{eq:hess:subgaussian:bounds} and bounds on $f_i'$ from Lemma~\ref{lem:growth_grad_f}, we can use an almost identical argument as in Proposition~\ref{prop:gradmartingale} to show the differences $|\mathcal{M}_k^{\text{Hess},q} - \mathcal{M}_k^{\text{Hess},q,\beta}|$ are exponentially small with overwhelming probability for $q=2,5,7$. 

Finally, it remains to bound the term $\mathcal{M}_{k}^{\hess,1, \beta}$ which is quadratic in $v_{k+1}$. The argument will once again be similar to the proof of Proposition~\ref{prop:gradmartingale} and the steps above. We have
\begin{equation}
    \langle \nabla^2 \varphi(X_k) , (\sqrt{K_i}v_{k+1})^{\otimes 2}\rangle = \langle \sqrt{K_i} \nabla^2 \varphi(X_k)\sqrt{K_i}, v_{k+1}^{\otimes 2}\rangle.
\end{equation}
It then follows from Hanson-Wright's inequality
\begin{equation}
    \begin{aligned}
    \Proba&\left(|\langle \sqrt{K_i} \nabla^2 \varphi(X_k)\sqrt{K_i}, v_{k+1}^{\otimes 2}\rangle - \Tr( \nabla^2 \varphi(X_k)K_i)| > t\right) 
    \\
    &\leq 2\exp\left(-C\min\left\{\frac{t^2}{\| \sqrt{K_i} \nabla^2 \varphi(X_k)\sqrt{K_i}\|^2}, \frac{t}{\| \sqrt{K_i} \nabla^2 \varphi(X_k)\sqrt{K_i}\|_{\operatorname{op}}}\right\}\right).
    \end{aligned}
\end{equation}
From Lemma~\ref{lem:S_derivative_bounds} and since $\|K_i\|_{\operatorname{op}} \leq C$ then
\begin{equation}
    \begin{aligned}
        \| \sqrt{K_i} \nabla^2 \varphi(X_k)\sqrt{K_i}\|_{\operatorname{op}} \leq \|\sqrt{K_i}\|_{\operatorname{op}}^2 \|\nabla^2 \varphi(X_k)\|_{\operatorname{op}}  \leq C.
    \end{aligned}
\end{equation}
Combining this with $\|\cdot\| \leq \sqrt{d} \|\cdot\|_{\operatorname{op}}$ then it follows
\begin{equation}
    \begin{aligned}
    \Proba&\left(|\langle \sqrt{K_i} \nabla^2 \varphi(X_k)\sqrt{K_i}, v_{k+1}^{\otimes 2}\rangle - \Tr( \nabla^2 \varphi(X_k)K_i)| > t\right) 
    \leq 2\exp\left(-\frac{\min(t^2 d^{-1}, t)}{C}\right).
    \end{aligned}
\end{equation}
From Cauchy--Schwarz and Lemma~\ref{lem:S_derivative_bounds}, $|\Tr( \nabla^2 \varphi(X_k)K)| \leq d\|\nabla^2 \varphi(X_k)\|_{\opt} \|K\|_{\opt} \leq C d$. Hence, we obtain 
\begin{equation}\label{}
    \Proba\left(\left|\langle \sqrt{K_i} \nabla^2 \varphi(X_k)\sqrt{K_i}, v_{k+1}^{\otimes 2}\rangle \right| > t\right) 
    \leq 2\exp\left(-\frac{\min((t-Cd)^2 d^{-1}, t-Cd)}{C}\right)
\end{equation}
Hence, for any $\zeta >0$, setting $t = d^{ 1+ \zeta}$ then with overwhelming probability
\begin{equation}\label{eq:hess:outerprod:a^2:overwhelming}
    \left|\langle \sqrt{K_i} \nabla^2 \varphi(X_k)\sqrt{K_i}, v_{k+1}^{\otimes 2}\rangle \right| < d^{1 + \zeta}.
\end{equation}
Now, from \eqref{eq:hess:bounds:for:secterm} and $\text{Proj}_{d \beta}( \langle \nabla^2 \varphi(X_k) , (\sqrt{K_i}v_{k+1})^{\otimes 2}\rangle) \leq d \beta$, we obtain the almost sure bound
\begin{equation}
    \begin{aligned}
       |  (f_i'(\text{Proj}_{\beta}(g_{k,i})))^2 \cdot\text{Proj}_{d\beta}(\ip{\nabla^2 \varphi(X_k), (\sqrt{K_i}v_{k+1})^{\otimes 2}   } )| \leq  dC(L(f))^2\beta^{3+2\alpha}
    \end{aligned}
\end{equation}
Since an almost sure bound holds under expectation, it follows
\begin{equation}
    |\mathcal{M}_k^{\text{Hess},1,\beta }| \leq \frac{\bar{\gamma}^2}{d}C(L(f))^2\beta^{3+2\alpha}.
\end{equation} 
By Azuma's inequality then for $0 \leq k \leq Td$
\begin{equation}
   \Proba\left( | \mathcal{M}_{k}^{\hess, 1,\beta } | > t\right)\leq 2 \exp\left(\frac{-t^2}{2\sum_{i=1}^k Cd^{-2}\beta^{2(3+2\alpha)}}\right)\leq 2 \exp\left(\frac{-t^2}{T  Cd^{-1}\beta^{2(3+2\alpha)}}\right)
\end{equation}
Thus, with overwhelming probability, for every $0 \leq k \leq Td$
\begin{equation}
    |\mathcal{M}_{k}^{\hess, 1,\beta }| < d^{-1/2}\beta^{4+2\alpha}.
\end{equation}
Since taking a union bound from $0$ to $Td$ does not modify the overwhelming probability statement, then with overwhelming probability
\begin{equation}\label{eq:hess:sup:boundM:1}
    \sup_{0 \leq k \leq Td}|\mathcal{M}_{k}^{\hess, 1,\beta }| < d^{-1/2}\beta^{4+2\alpha}.
\end{equation}
It remains to bound the difference between $ \mathcal{M}_{k}^{\hess, 1, }$ and $\mathcal{M}_{k}^{\hess, 1,\beta }$. Once again, the argument is extremely similar to the one used in the proof of Proposition~\ref{prop:gradmartingale}. For notational convenience, we write 
\begin{equation}
    \begin{aligned}
        F_{k,i}&\defas\frac{\gamma_k^2}{2d^2}  (f_{i}'(g_{k, i} ))^2\ip{\nabla^2 \varphi(X_k), (\sqrt{K_i}v_{k+1})^{\otimes 2}   }  ,\\
        F_{k,i}^\beta &\defas\frac{\gamma_k^2}{2d^2}  (f_i'(\text{Proj}_{\beta}(g_{k, i} )))^2 \cdot\text{Proj}_{d\beta}(\ip{\nabla^2 \varphi(X_k), (\sqrt{K_i}v_{k+1})^{\otimes 2}   } ).
    \end{aligned}
\end{equation}
We will show the following expression is 
\begin{equation}
    |(F_{k,i} - \E[F_{k,i}]) - (F^{\beta}_{k,i} - \E[F^{\beta}_{k,i}])| \leq |F_{k,i} -F^{\beta}_{k,i}| + |\E [F_{k,i} -F^{\beta}_{k,i}]|,
\end{equation}
exponentially small in $d$. Using the subgaussian tail bound \eqref{eq:tail:bounds:proj:grad} we previously derived in the proof of Proposition~\ref{prop:gradmartingale} and \eqref{eq:hess:outerprod:a^2:overwhelming}, then 
\begin{equation}
    \begin{aligned}
    \Proba(|F_{k,i} - F^{\beta}_{k,i}| > 0) 
    &
    \leq \sum_{i=1}p_i (\Proba(\|(g_{k,i}, \epsilon_{k+1})\| > \beta) + \Proba( |\ip{\nabla^2 \varphi(X_k), (\sqrt{K_{i}}v_{k+1} )^{\otimes 2}  }|  > d \beta))
    \\
    &\leq 2\exp\left(-\Omega(\beta^2)\right) + 2\exp\left(-\frac{d\min((\beta-C)^2 , \beta-C)}{C}\right)
    \end{aligned}
\end{equation}
Since we defined $\beta = d^{\zeta}$ for some $\zeta>0$, then for sufficiently large $d$, we have 
$\beta - C \geq \frac{\beta}{2}$. Hence, for sufficiently large $d$ we have
\begin{equation}\label{eq:hess:diff:F_k}
    \begin{aligned}
    \Proba(|F_{k,i} - F^{\beta}_{k,i}| > 0) 
    &\leq C \exp\left(-\Omega(d\min(\beta^2 , \beta))\right).
    \end{aligned}
\end{equation}
which implies $F_{k,i} - F^{\beta}_{k,i} = 0$ with overwhelming probability. It then follows from a union bound that 
\begin{equation}
    \sum_{k=1}^{Td} |F_{k,I_{k+1}} - F^{\beta}_{k,I_{k+1}}| = 0,
\end{equation}
with overwhelming probability.  We omit the proof that $|\E [F_{k,i} -F^{\beta}_{k,i}]|$ is exponentially small in $d$, since the steps are almost identical to the ones in the proof of Proposition~\ref{prop:gradmartingale} using the growth bound on $f'$, \eqref{eq:tail:bounds:proj:grad} and \eqref{eq:hess:outerprod:a^2:overwhelming}. Thus, the difference $|\mathcal{M}_k^{\text{Hess},1} - \mathcal{M}_k^{\text{Hess},1,\beta}|$ is exponentially small with overwhelming probability.

Finally, combining \eqref{eq:hess:sup:boundM:35}, \eqref{eq:hess:sup:boundM:24} and \eqref{eq:hess:sup:boundM:1} with the bounds we derived for the differences $|\mathcal{M}_k^{\text{Hess},q} - \mathcal{M}_k^{\text{Hess},q,\beta}|$ which hold with overwhelming probability for $q=1,...,9$, with a similar argument as the end of the proof of Proposition~\ref{prop:gradmartingale}, we obtain
\begin{equation}
    \sup_{0 \leq k \leq Td}|\mathcal{M}_{k}^{\hess}| < d^{-\frac{1}{2}+(4+2\alpha)\zeta}.
\end{equation}
Hence for any arbitrarily small value of $\zeta < \frac{1}{2(4+2\alpha)}$, the result follows.
\end{proof}
\subsubsection{Bounds on the Error Term $\mathcal{E}_k^{\hess}$ of the Hessian Martingale Increments \label{sec:lower_term_Hessian}}
In this section, we will bound the lower order error terms arising from the Hessian martingale increments defined in \eqref{eq:error_terms}. We build upon Proposition A.4 in \citep{collinswoodfin2023hitting} and Lemma~11 in \citep{collinswoodfin2025Exact}. Once again, we provide the proof for the soft label setting as the proof for hard labels follows the same argument with $X^{\star} = 0$.

\begin{proposition}[Hessian error term]\label{prop:Hessian-error} For $i=1,2$, suppose $f_i \, : \,\R^{3} \to \mathbb{R}$ is $\alpha$-pseudo-Lipschitz function with constant $L(f)$ (see Assumption~\ref{assumption:pseudo_lipschitz}). Let the statistic $Z \, : \, \R^d \to \R^{4 \times 4}$ be defined as in \eqref{eq:def:stats:S:M:Z}. Then, for any $T > 0$, with overwhelming probability
    \begin{equation}
      \sup_{0 \le t \le T} \sum_{k=0}^{\lfloor(t \wedge \vartheta) d\rfloor-1} | \EE[ \mathcal{E}_k^{\text{\rm Hess}}(Z(\cdot, z)) \, | \, \mathcal{F}_{k}] | \le C (L(f))^2  d^{-1/2}.
    \end{equation}
    \end{proposition}
\begin{proof}
    Similarly to in the proofs of Propositions~\ref{prop:gradmartingale} and \ref{prop:hessmartingale}, we denote $\varphi \defas Z_{j\ell}(W,z)$ as the $j\ell$-coordinate of $Z$ for a fixed $z \in \Gamma^2$. Throughout this proof, we will be working with the stopped versions of the stochastic processes involved. For notational simplicity,  we omit the dependence on the stopping time $\vartheta$ and on $\varphi$ such that $\mathcal{E}_{k}^{\hess} \defas \mathcal{E}_{(k/d \wedge \vartheta)d}^{\hess}(\varphi )$. As in the proofs of Propositions~\ref{prop:gradmartingale} and \ref{prop:hessmartingale}, for constant $C>0$, we omit the dependence on $i$ as we could redefine $C = \max\{C_1, C_2\}$. We also write $\nabla \varphi(X_k) \defas \nabla_X \varphi(X_k)$.
    
    Recall we defined $\Pi_{k,i} = Q_{k,i} Q_{k,i}^\top$ where $Q_{k,i} \in \R^{d\times \ell_i}$ a matrix with orthonormal columns with $\ell_i = \text{rank}(\sqrt{K_i}\hat{X}_k) \leq 2$. It is then easy to see that $\|\Pi_{k,i}\|^2 \leq 2$. Before proceeding with the proof, we introduce the definition of the nuclear norm
    \begin{equation}
        \|A \|_{*} \defas \sup_{\|B\|_{\operatorname{op}}=1}\langle A, B\rangle.
    \end{equation}
    Since $Q_{k,i}$ is a matrix with orthonormal columns, note that $\|\Pi_{k,i}\|_{*} \leq 2$. Recall from \eqref{eq:error_terms} the definition of $\mathcal{E}_k^{\hess}$
    \begin{equation}
        \begin{aligned}
            \mathcal{E}_k^{\text{Hess}}(\varphi)
    &
    = \frac{\gamma_k^2}{d^2}\sum_{i=1}^2 p_i (f_{i}'(g_{k, i} ))^{ 2} \cdot \bigg(-\ip{\nabla^2 \varphi(X_k), \sqrt{K_i} \Pi_{k,i}\sqrt{K_i}  } 
    \\
    &
    +
       \ip{ \nabla^2 \varphi (\hat{X}_k), \big ( \sqrt{K_i} \Pi_{k,i} v_{k+1} \big )^{\otimes 2}  }
       +
       2\ip{ \nabla^2 \varphi (\hat{X}_k), \sqrt{K_i} \Pi_{k,i} v_{k+1}\mu_i^\top }
       \bigg)
       \\
       &
       +
       \frac{\gamma_k^2}{d^2} \sum_{i=1}^2 p_i\delta s_{k, i} ( f_{i}'(g_{k, i} ))^{2}\ip{ \nabla^2 \varphi (\hat{X}_k),     (\sqrt{K_i} \Pi_{k,i} v_{k+1} + \mu_i) \otimes \tfrac{X_k}{\|X_k\|}    } 
       \\
        &
        +
    \frac{\gamma_k^2}{2d^2}\sum_{i=1}^2 p_i \delta^2 s_{k, i}^2 ( f_{i}'(g_{k, i} ))^{2} \ip{ \nabla^2 \varphi (\hat{X}_k),   \tfrac{X_k^{\otimes 2}}{\|X_k\|^2}} 
    \\
    & 
    +
    \frac{\gamma_k^2}{2d^2} \ip{ \nabla^2 \varphi(X_k), (\lambda X_k)^{\otimes 2}}
    \\
    &
    +
    \frac{\gamma_k^2}{d^2} \sum_{i=1}^2 p_i f_{i}'(g_{k, i} )\ip{ \nabla^2 \varphi(X_k), (\sqrt{K_i} \Pi_{k,i} v_{k+1} + \mu_i + \delta s_{k, i} \tfrac{X_k}{\|X_k\|}) \otimes \lambda X_k   }
    .
\end{aligned}
\end{equation}
Starting with the first term, we have
\begin{equation}
    \begin{aligned}
      | ( f_{i}'(g_{k, i} ))^{ 2}\cdot\ip{\nabla^2 \varphi(X_k), \sqrt{K_i} \Pi_k \sqrt{K_i}  } | 
      &
      \leq 
      | f_{i}'(g_{k, i} )|^{2 }\|\nabla^2 \varphi(X_k)\| \|\sqrt{K_i} \Pi_k \sqrt{K_i} \|
      \\
      &
      \leq
      | f_{i}'(g_{k, i} )|^{2 }\|\nabla^2 \varphi(X_k)\| \|K_i\|_{\operatorname{op}} \|\Pi_k\|_{*}
      \\
      &
      \leq
      2| f_{i}'(g_{k, i} )|^{2 }\|\nabla^2 \varphi(X_k)\| \|K_i\|_{\operatorname{op}}.
    \end{aligned}
\end{equation}
From Lemma~\ref{lem:S_derivative_bounds}, we have
\begin{equation}\label{eq:errort:Sgrad:bound}
    \|\nabla^2 \varphi(X_k)\| = \|\nabla^2 Z_{j\ell}(W_k,z)\| \leq \sup_{z \in \Gamma}\|\nabla^2 Z_{j\ell}(W_k,z)\| \leq \|\nabla^2 Z(W_k,z)\|_{\Gamma} \leq C \sqrt{d}.
\end{equation}
By Assumption~\ref{assumption:data}, $\|K_i\|_{\operatorname{op}} \leq C$ and from Lemma~\ref{lem:growth_grad_f}, we have $\E[|f_{i}'(g_{k, i} )|^2\mid \mathcal{F}_k]\leq CL(f)^2(1 +( \|K_i\|_{\operatorname{op}}^{1/2} + \|\mu_i\|)\|\hat{X}_k\|+\delta\|X_k\|+\eta)^{\max(1, 2\alpha)}$. Since we are working with the stopped processes, it then follows
\begin{equation}
    \begin{aligned}
    \E[|  f_{i}'(g_{k, i} )^{2 }&\cdot\ip{\nabla^2 \varphi(X_k), \sqrt{K_i} \Pi_k \sqrt{K_i}  } | \mid \mathcal{F}_k]  
    \\
    &
    \leq
     C\sqrt{d}L(f)^2(1 + ( \|K_i\|_{\operatorname{op}}^{1/2} + \|\mu_i\|)\|\hat{X}_k\|+\delta\|X_k\|+\eta)^{\max(1, 2\alpha)} = \mathcal{O}(d^{1/2}).
    \end{aligned}
\end{equation}
For the second term, we obtain
\begin{equation}\label{eq:hessError:upper:bound}
    \begin{aligned}
       | f_{i}'(g_{k, i} )^{2 } \ip{ \nabla^2 \varphi (X_k), \big ( \sqrt{K_i} \Pi_{k,i} v_{k+1} \big )^{\otimes 2}  }|
       & 
       \leq
       |f_{i}'(g_{k, i} )|^{2 } \|\nabla^2 \varphi (X_k)\| \| \sqrt{K_i} \Pi_{k,i} v_{k+1} \big \|^2
       \\
       &
       \leq
       |f_{i}'(g_{k, i} )|^{2 } \|\nabla^2 \varphi (X_k)\| \|K_i\|_{\operatorname{op}} \|\Pi_{k,i} v_{k+1} \big \|^2
    \end{aligned}
\end{equation}
Once again, from Lemma~\ref{lem:growth_grad_f}, we have $\E[|f_{i}'(g_{k, i} )|^4\mid \mathcal{F}_k]\leq CL(f)^4(1 + (\|K\|_{\operatorname{op}}^{1/2} + \|\mu_i\|)\|\hat{X}_k\|+\delta\|X_k\|+\eta)^{\max(1, 4\alpha)}$. Now, since $Q_{k,i}$ has $\ell \leq 2$ orthonormal columns then $\Pi_{k,i}$ is a projection matrix and since $Q_{k,i} v_{k+1} \sim \mathcal{N}(0, \Id_{\ell}) $ it follows $\E[\|\Pi_{k,i} v_{k+1}\|^4] = \E[(\|Q_kv_{k+1}\|^2)^2] = \E[(\chi^2_{\ell})^2] \leq 8$. Hence, by \eqref{eq:errort:Sgrad:bound} and \eqref{eq:hessError:upper:bound} 
\begin{equation}
    \begin{aligned}
       \E[| f_{i}'(g_{k, i} )^{2 } &\ip{ \nabla^2 \varphi (X_k), \big ( \sqrt{K_i} \Pi_{k,i}v_{k+1} \big )^{\otimes 2}  }|\mid \mathcal{F}_k]
        \\
       &
       \leq
      C\sqrt{d}L(f)^2(1 + (\|K_i\|_{\operatorname{op}}^{1/2} + \|\mu_i\|)\|\hat{X}_k\|+\delta\|X_k\|+\eta)^{\max(\frac{1}{2}, 2\alpha)} = \mathcal{O}(d^{1/2}).
    \end{aligned}
\end{equation}
Using a similar argument, we obtain for the third term
\begin{equation}
    \begin{aligned}
        \E[|f_{i}'(g_{k, i} )^{2}\cdot \ip{ \nabla^2 \varphi (\hat{X}_k), \sqrt{K_i} \Pi_{k,i} v_{k+1}\mu_i^\top }|\mid \mathcal{F}_k]= \mathcal{O}(d^{1/2}).
    \end{aligned}
\end{equation}
Similarly, since $|s_k|\leq 1$ for the fourth term 
\begin{equation}\label{eq:hessErr:bound:2}
    \begin{aligned}
        |\delta s_{k, i} ( f_{i}'(g_{k, i} ))^{2}&\ip{ \nabla^2 \varphi (\hat{X}_k),     (\sqrt{K_i} \Pi_{k,i} v_{k+1} + \mu_i) \otimes \tfrac{X_k}{\|X_k\|}}| 
        \\
        &
        \leq 
        |f_{i}'(g_{k, i} )|^{2} \|\nabla^2 \varphi(X_k)\|(\|K\|^{1/2}_{\operatorname{op}} \|\Pi_k v_{k+1}\| + \|\mu_i\|)
    \end{aligned}
\end{equation}
Using a similar argument to what we just did, from Cauchy--Schwarz, Lemma~\ref{lem:growth_grad_f}, \eqref{eq:errort:Sgrad:bound}, \eqref{eq:hessErr:bound:2} and since $\E[\|\Pi_k v_{k+1}\|^2] = \E[\chi_\ell^2] \leq 2$ we obtain
\begin{equation}
    \begin{aligned}
    \E[|s_k^2 f_{i}'(g_{k, i} )^{2}  &\cdot  \ip{\nabla^2 \varphi(X_k), \sqrt{K_i} \Pi_k v_{k+1}  \otimes \tfrac{X_k}{\|X_k\|}}| \mid \mathcal{F}_k] 
    \\
    &
    \leq C\sqrt{d}L(f)^2(1 + (\|K_i\|_{\operatorname{op}}^{1/2} + \|\mu_i\|)\|\hat{X}_k\|+\delta\|X_k\|+\eta)^{\max(\frac{1}{2}, 2\alpha)} = \mathcal{O}(d^{1/2}).
    \end{aligned}
\end{equation}
For the fifth term, since $|s_{k,i}|\leq 1$ we have
\begin{equation}
    \begin{aligned}
        | s_{k, i}^2 ( f_{i}'(g_{k, i} ))^{2} \ip{ \nabla^2 \varphi (\hat{X}_k),   \tfrac{X_k^{\otimes 2}}{\|X_k\|^2}} | \leq |f_{i}'(g_{k, i} )|^2 \|\nabla^2 \varphi (X_k)\|.
    \end{aligned}
\end{equation}
Hence, it follows from Lemma~\ref{lem:growth_grad_f} and \eqref{eq:errort:Sgrad:bound}
\begin{equation}
    \begin{aligned}
    \E[  | s_{k, i}^2 &( f_{i}'(g_{k, i} ))^{2} \ip{ \nabla^2 \varphi (\hat{X}_k),   \tfrac{X_k^{\otimes 2}}{\|X_k\|^2}} | \mid \mathcal{F}_k]
    \\
    &
    \leq 
    C\sqrt{d}L(f)^2(1 + (\|K_i\|_{\operatorname{op}}^{1/2} + \|\mu_i\|)\|\hat{X}_k\|+\delta\|X_k\|+\eta)^{\max(1, 2\alpha)}=\mathcal{O}(d^{1/2}).
    \end{aligned}
\end{equation}
For the firth term, from \eqref{eq:errort:Sgrad:bound} and recalling that we are working with the stopped processes we have 
\begin{equation}
    \E[|\ip{ \nabla^2 \varphi(X_k), (\lambda X_k)^{\otimes 2}}|\mid \mathcal{F}_k] \leq C\sqrt{d}.
\end{equation}
Finally, for the last term, since $|s_k|\leq 1$ we have
\begin{equation}
    \begin{aligned}
        | f_{i}'(g_{k, i} )&\ip{ \nabla^2 \varphi(X_k), (\sqrt{K_i} \Pi_{k,i} v_{k+1} + \mu_i + \delta s_{k, i} \tfrac{X_k}{\|X_k\|}) \otimes \lambda X_k   }|
        \\
        &\leq \lambda |f_{i}'(g_{k, i} )|\|\nabla^2 \varphi(X_k)\| \|X_k\|(\|K_i\|_{\operatorname{op}}^{1/2} \|\Pi_{k,i} v_{k+1}\| + \|\mu_i\|+ \delta)
    \end{aligned}
\end{equation}
Thus, using a similar argument as above, from Lemma~\ref{lem:growth_grad_f}, \eqref{eq:errort:Sgrad:bound} and Cauchy--Schwarz we have 
\begin{equation}
    \begin{aligned}
        \E[ |f_{i}'(g_{k, i} )\cdot &\ip{\nabla^2 \varphi(X_k), (\sqrt{K_i} \Pi_k v_{k+1} +  \delta s_k \tfrac{X_k}{\|X_k\|})\otimes\lambda X_k}|\mid \mathcal{F}_k]
        \\
        &
        \leq
        C\sqrt{d}L(f)(1 + (\|K_i\|_{\operatorname{op}}^{1/2} + \|\mu_i\|)\|\hat{X}_k\|+\delta\|X_k\|+\eta)^{\max(\frac{1}{2}, \alpha)} = \mathcal{O}(d^{1/2}).
    \end{aligned}
\end{equation}
Since $\tfrac{\gamma_k^2}{d^2} \leq \tfrac{\bar{\gamma}^2}{d^2}$ for all $k\geq 0$ and these bounds hold uniformly in $i$, the result follows from summing over $i$, from Jensen's inequality, the triangle inequality and noting that the summation contributes a factor $Td$.
\end{proof}
\section{Integro-Differential Equation Analysis}\label{sec:IntegroDiffAnalysis}
We provide alternative characterization for the solution to the integro-differential equation \eqref{eq:ODE_resolvent_2} which we recall below.
\begin{mdframed}[style=exampledefault]
    \textbf{Integro-Differential Equation for $\mathcal{Z}(t, z)$.} For any contour $\Gamma \subset \mathbb{C}$ enclosing the eigenvalues of $K_1$ and $K_2$, we have an expression for the derivative of $\mathcal{Z}$ for $z\in \Gamma^2$:
    \begin{equation}\label{eq:ODE_resolvent_2:app}
        \dif \mathcal{Z}(t, \cdot) = \mathscr{F}(z, \mathcal{Z}(t, \cdot)) \, \dif t,\quad \text{with initial condition}\quad  \mathcal{Z}(0, z) = Z(W_0, z),
        \end{equation}
        where
        \begin{equation}
            \begin{gathered}
            \mathcal{Z}(t, z) \defas \begin{bmatrix}
                \mathcal{S}(t, z) &  \mathcal{M}(t, z) \\
                \mathcal{M}(t, z)^\top & \mu^\top \mathscr{R}(z;\mathcal{K})\mu
            \end{bmatrix},
            \\
            \text{and}\quad
            \mathscr{F}(z, \mathcal{Z}(t, \cdot)) \defas  \begin{bmatrix}
                \mathscr{F}_\mathcal{S}(z, \mathcal{Z}(t, \cdot)) & \mathscr{F}_\mathcal{M}(z,\mathcal{Z}(t, \cdot)) \\
                \mathscr{F}_\mathcal{M}(z, \mathcal{Z}(t, \cdot))^\top & 0_{2\times 2}
            \end{bmatrix},
        \end{gathered}
        \end{equation}
The maps $\mathscr{F}_\mathcal{S}$ and $\mathscr{F}_\mathcal{M}$ are defined as follows
\begin{equation}\label{eq:ODE_resolvent_2:M:app}
    \begin{aligned}
        &\dif \mathcal{M}(t,z) 
        =
        -2\gamma(t) \sum_{i=1}^2 p_i\bigg(H_{1,i}^\top( \widehat{\mathrsfs{B}}_i(t))  \left(z_i\mathcal{M}(t,z)- \frac{1}{2\pi\operatorname{i}} \oint_{\Gamma} \mathcal{M}(t, z) \dif z_i\right)
        \\
        &
        + H_{2,i} (\widehat{\mathrsfs{B}}_i(t)) \otimes \mu^\top \mathscr{R}(z; \mathcal{K})\mu_i
        + \left(\frac{\lambda}{2} + \partial_{44} h_i (\widehat{\mathrsfs{B}}_i(t) )\right)   \operatorname{D}  \mathcal{M}(t,z)\bigg)\dif t,
    \end{aligned}
    \end{equation}
    and
    \begin{equation}\label{eq:ODE_resolvent_2:S:app}
        \begin{aligned}
    &\dif \mathcal{S}(t,z) 
    =
    -2\gamma(t) \sum_{i=1}^2 p_i\bigg(\left(z_i\mathcal{S}(t,z)- \frac{1}{2\pi\operatorname{i}} \oint_{\Gamma} \mathcal{S}(t, z) \dif z_i\right) H_{1,i}(\widehat{\mathrsfs{B}}_i(t)) 
    \\
    &
     + H_{1,i}^\top( \widehat{\mathrsfs{B}}_i(t))  \left(z_i\mathcal{S}(t,z)- \frac{1}{2\pi\operatorname{i}} \oint_{\Gamma} \mathcal{S}(t, z) \dif z_i\right)
    \\
    &
    + \mathcal{M}_i(t, z) \otimes H_{2,i} (\widehat{\mathrsfs{B}}_i(t)) 
    + H_{2,i} (\widehat{\mathrsfs{B}}_i(t)) \otimes \mathcal{M}_i(t, z)
    \\
    &
    + \left(\frac{\lambda}{2} + \partial_{44} h_i (\widehat{\mathrsfs{B}}_i(t) )\right) (\mathcal{S}(t,z) \operatorname{D}  +  \operatorname{D}  \mathcal{S}(t,z))\bigg)\dif t
    \\
    &
    + \frac{\gamma(t)^2}{d} \sum_{i=1}^2 p_i \mathcal{I}_i(\widehat{\mathrsfs{B}}_i(t))\Tr(( K_i + \mu_i \mu_i^\top ) \mathscr{R}(z;\mathcal{K}))  \operatorname{D} \, \dif t,
\end{aligned}
\end{equation}
\begin{gather}\label{eq:def:B(t):app}
\text{where} \quad
H_{1,i}(\widehat{\mathrsfs{B}}_i)=\left [ \begin{array}{c|c} 
   \partial_{11} h_i(\widehat{\mathrsfs{B}}_i)  & 0 \\
    \hline
   \partial_{21} h_i(\widehat{\mathrsfs{B}}_i) & 0 
    \end{array} \right ],\quad
    H_{2,i}(B)= \left [ \begin{array}{c} 
        \partial_{13} h_i(\widehat{\mathrsfs{B}}_i)   \\
         \hline
        0 
         \end{array} \right ],\quad
 \operatorname{D} = \left [ \begin{array}{c|c} 
 1 & 0\\
 \hline
 0 & 0
 \end{array} \right ]
 \nonumber
 \\
\widehat{\mathrsfs{B}}_i(t) = \begin{bmatrix}
    \mathrsfs{B}_i(t) & \mathfrak{m}_i(t) & 0  \\
    \mathfrak{m}_i(t)^\top &  \oint \mu_i^\top  \mathscr{R}(z;\mathcal{K}) \mu_i \Dif z & 0 \\
    0 & 0 & \oint \mathcal{S}_{11}(t,z) \Dif z
 \end{bmatrix}.
\end{gather}
Here $\mathrsfs{B}_i(t) = \oint z_i \mathcal{S}(t,z) \Dif z \in \R^{2\times 2}$, $\mathfrak{m}_i(t) = \oint \mathcal{M}_i(t) \Dif z$ where $\mathcal{M}_i(t) \in \R^{2}$ denotes the $i$th column of $\mathcal{M}$, associated with $\mu_i$. 
\end{mdframed}
We proceed using the spectral decomposition of the resolvent $\mathscr{R}(z;\mathcal{K})$ similarly to Lemma~2 in \citep{collinswoodfin2025Exact} as follows:
\begin{equation}
    \mathscr{R}(z;\mathcal{K}) = \sum_{j=1}^d \frac{1}{\lambda^{(1)}_j - z_1} \omega_j\omega_j^\top \sum_{\ell=1}^d \frac{1}{\lambda^{(2)}_\ell - z_2} \omega_\ell\omega_\ell^\top = \sum_{j=1}^d \frac{1}{\prod_{\ell=1}^2(\lambda^{(\ell)}_j - z_\ell)} \omega_j \omega_j^\top
\end{equation} 
We will now assume the following decompositions of $\mathcal{S}(t,z)$ and $\mathcal{M}(t,z)$ for some functions $\mathrsfs{V}_j$ and $\mathfrak{m}_{j,\ell}$ which we will determine later.
\begin{equation}\label{eq:app:decomps:S:M}
    \mathcal{S}(t,z) = \sum_{j=1}^d \frac{\mathrsfs{V}_j(t)}{\prod_{\ell=1}^2(\lambda^{(\ell)}_j - z_\ell)}\quad \text{and}\quad \mathcal{M}_i(t,z) = \sum_{j=1}^d \frac{\mathfrak{m}_{j,i}(t)}{\prod_{\ell=1}^2(\lambda^{(\ell)}_j - z_\ell)}
\end{equation}
Using \eqref{eq:ODE_resolvent_2:app}, we determine $\mathrsfs{V}_j$ and $\mathfrak{m}_{j,\ell}$ are defined as follows 
\begin{equation}\label{eq:coupled:ODES:app}
    \begin{aligned}
        &\frac{\dif \mathrsfs{V}_j(t)}{\dif t} = -2\gamma(t) \sum_{i=1}^2 p_i \bigg(\lambda_j^{(i)}(\mathrsfs{V}_j(t) H_{1,i}(\widehat{\mathrsfs{B}}_i(t)) + H_{1,i}(\widehat{\mathrsfs{B}}_i(t))^\top \mathrsfs{V}_j(t) ) 
        + \mathfrak{m}_{j,i}(t) \otimes H_{2,i}(\widehat{\mathrsfs{B}}_i(t)) 
        \\
        &
        \quad
        + H_{2,i}(\widehat{\mathrsfs{B}}_i(t)) \otimes \mathfrak{m}_{j,i}(t) 
        +\left(\frac{\lambda}{2} + \partial_{44} h_i (\widehat{\mathrsfs{B}}_i(t) )\right) (\mathrsfs{V}_j(t)\operatorname{D} + \operatorname{D} \mathrsfs{V}_j(t))\bigg)
        \\
        &
        \quad
        +\frac{\gamma(t)^2}{d} \left(\sum_{i=1}^2 p_i \mathcal{I}_i(\widehat{\mathrsfs{B}}_i(t))(\lambda_j^{(i)} + (\omega_j^\top \mu_i)^2)\right) \operatorname{D},
        \\
        &
        \frac{\dif \mathfrak{m}_{j,\ell}(t)}{\dif t}  = -2\gamma(t) \sum_{i=1}^2 p_i \bigg( \lambda_j^{(i)}H_{1,i}(\widehat{\mathrsfs{B}}_i(t))^\top \mathfrak{m}_{j,\ell}(t)  + \langle \omega_j, \mu_\ell\rangle \langle \omega_j, \mu_i\rangle \cdot H_{2,i}(\widehat{\mathrsfs{B}}_i(t))
        \\
         &\quad+ \left(\frac{\lambda}{2} + \partial_{44} h_i (\widehat{\mathrsfs{B}}_i(t)) \right)  \operatorname{D}\mathfrak{m}_{j,\ell}(t)\bigg),
    \end{aligned}
\end{equation}
with $\mathrsfs{V}_j(0) = \hat{X}_0^\top \omega_j \omega_j^\top \hat{X}_0$ and $\mathfrak{m}_{j,\ell}(0) = \hat{X}_0^\top \omega_j \omega_j^\top\mu_{\ell}$. It immediately follows from Lemma 2 in \citep{collinswoodfin2025Exact} that $\mathcal{S}(t,z)$ and $\mathcal{M}_{\ell}(t,z)$ defined in \eqref{eq:app:decomps:S:M} with $\mathrsfs{V}_j$ and $\mathfrak{m}_{j,\ell}$ satisfying \eqref{eq:coupled:ODES:app} are solutions to \eqref{eq:ODE_resolvent_2:app}. We omit the proof as the argument is identical. We will now use a similar argument to Appendix A in \citep{collinswoodfin2023hitting} to determine the explicit dynamics of $\widehat{\mathrsfs{B}}_i(t))$ for $i=1,2$. Let $\Phi$ be the fundamental matrix of the ODE:
\begin{equation}\label{eq:app:base:ODE}
   \dot{\Phi}_j =2\gamma(t)\sum_{i=1}^2 p_i \bigg( \lambda_j^{(i)} H_{1,i}(\widehat{\mathrsfs{B}}_i(t))  + \left(\frac{\lambda}{2} + \partial_{44} h_i (\widehat{\mathrsfs{B}}_i(t)) \right) \operatorname{D}\bigg)\Phi_j,\quad \Phi_j(0)= \Id_{2}. 
\end{equation}
Define the following terms:
\[
\begin{aligned}
   U_{0,\mathrsfs{V}_j}(t) &
   \defas -2\gamma(t) \sum_{i=1}^2 p_i \bigg( \mathfrak{m}_{j,i}(t) \otimes H_{2,i}(\widehat{\mathrsfs{B}}_i(t)) + H_{2,i}(\widehat{\mathrsfs{B}}_i(t)) \otimes \mathfrak{m}_{j,i}(t) \bigg)
   \\
   &
   \quad
   +\frac{\gamma(t)^2}{d} \left(\sum_{i=1}^2 p_i \mathcal{I}_i(\widehat{\mathrsfs{B}}_i(t))(\lambda_j^{(i)} + (\omega_j^\top \mu_i)^2)\right) \operatorname{D},
\end{aligned}
\]
and
\[
\begin{aligned}
   U_{0,\mathfrak{m}_{j,\ell}}(t) =
   -2\gamma(t) \sum_{i=1}^2 p_i  \langle \omega_j, \mu_\ell\rangle \langle \omega_j, \mu_i\rangle \cdot H_{2,i}(\widehat{\mathrsfs{B}}_i(t)).
\end{aligned}
\]
Then it follows that 
\[
\begin{aligned}
    \dot{(\Phi_j^\top \mathfrak{m}_{j,\ell})} 
    &
    =
    \dot{\Phi_j}^\top \mathfrak{m}_{j,\ell}+ \Phi_j^\top \dot{\mathfrak{m}_{j,\ell}}
    \\
    &=
    \Phi_j^\top\bigg(2\gamma(t)\sum_{i=1}^2 p_i \bigg( \lambda_j^{(i)} H_{1,i}(\widehat{\mathrsfs{B}}_i(t))^\top  + \left(\frac{\lambda}{2} + \partial_{44} h_i (\widehat{\mathrsfs{B}}_i(t)) \right) \operatorname{D}\bigg)\bigg) \mathfrak{m}_{j,\ell}
    \\
    &
    + \Phi_j^\top \bigg(U_{0,\mathfrak{m}_{j,\ell}}-2\gamma(t)\sum_{i=1}^2 p_i \bigg( \lambda_j^{(i)} H_{1,i}(\widehat{\mathrsfs{B}}_i(t))^\top  + \left(\frac{\lambda}{2} + \partial_{44} h_i (\widehat{\mathrsfs{B}}_i(t)) \right) \operatorname{D}\bigg)\mathfrak{m}_{j,\ell}\bigg)
    \\
    &= \Phi_j^\top U_{0,\mathfrak{m}_{j,\ell}}.
\end{aligned}
\]
Similarly, we obtain
\begin{equation}\label{eq:app:similar:arg:ODE}
\begin{aligned}
    \dot{(\Phi_j^\top \mathrsfs{V}_j\Phi_j)} 
    &
    =
    \dot{\Phi_j}^\top \mathrsfs{V}_j\Phi_j + \Phi_j^\top \dot{\mathrsfs{V}_j}\Phi_j + \Phi_j^\top \mathrsfs{V}_j\dot{\Phi_j}
    \\
    &=
    \Phi_j^\top \bigg(2\gamma(t)\sum_{i=1}^2 p_i \bigg( \lambda_j^{(i)} H_{1,i}(\widehat{\mathrsfs{B}}_i(t))^\top  + \left(\frac{\lambda}{2} + \partial_{44} h_i (\widehat{\mathrsfs{B}}_i(t)) \right) \operatorname{D}\bigg)\bigg)  \mathrsfs{V}_j\Phi_j 
    \\
    &
    + \Phi_j^\top \bigg(U_{0,\mathrsfs{V}_j}- \mathrsfs{V}_j\bigg(2\gamma(t)\sum_{i=1}^2 p_i \bigg( \lambda_j^{(i)} H_{1,i}(\widehat{\mathrsfs{B}}_i(t))  + \left(\frac{\lambda}{2} + \partial_{44} h_i (\widehat{\mathrsfs{B}}_i(t)) \right) \operatorname{D}\bigg)\bigg)
    \\
    & - \bigg(2\gamma(t)\sum_{i=1}^2 p_i \bigg( \lambda_j^{(i)} H_{1,i}(\widehat{\mathrsfs{B}}_i(t)) ^\top + \left(\frac{\lambda}{2} + \partial_{44} h_i (\widehat{\mathrsfs{B}}_i(t)) \right) \operatorname{D}\bigg)\bigg)\mathrsfs{V}_j \bigg)\Phi_j
    \\
    &
    + \Phi_j^\top \mathrsfs{V}_j\bigg(2\gamma(t)\sum_{i=1}^2 p_i \bigg( \lambda_j^{(i)}  H_{1,i}(\widehat{\mathrsfs{B}}_i(t))  + \left(\frac{\lambda}{2} + \partial_{44} h_i (\widehat{\mathrsfs{B}}_i(t)) \right) \operatorname{D}\bigg)\bigg)\Phi_j
    \\
    &= \Phi_j^\top U_{0, \mathrsfs{V}_j}\Phi_j.
\end{aligned}
\end{equation}
These ODEs are solvable.
\begin{mdframed}[style=exampledefault]
    \textbf{Resolvent formula.}
    We have
    \begin{equation}
        \label{eq:resolvent_form:M}
        \mathfrak{m}_{j,\ell}(t) = \Phi_j^{-\top}(t) \mathfrak{m}_{j,\ell}(0) + \int_0^t \Phi_j^{-\top}(t) \Phi_j^\top(s) U_{0,\mathfrak{m}_{j,\ell}}(s) \, \dif s
        \end{equation}
    and
    \begin{equation}
    \label{eq:resolvent_form:S}
    \mathrsfs{V}_j(t) = \Phi_j^{-\top}(t) \mathrsfs{V}_j(0) \Phi_j^{-1}(t) + \int_0^t \Phi_j^{-\top}(t) \Phi_j^\top(s) U_{0,\mathrsfs{V}_j}(s) \Phi_j(s) \Phi_j^{-1}(t) \, \dif s
    \end{equation}
    
    \begin{gather*}
    \text{where} \quad
    \mathrsfs{V}_j(0) = \hat{X}_0^\top \omega_j \omega_j^\top \hat{X}_0,\quad \mathfrak{m}_{j,\ell}(0) = \hat{X}_0^\top \omega_j \omega_j^\top\mu_{\ell},
     \\
     \text{$\Phi_j(t)$ is the solution to $\dot{\Phi}_j =2\gamma(t)\sum_{i=1}^2 p_i \bigg( \lambda_j^{(i)} H_{1,i}(\widehat{\mathrsfs{B}}_i(t))  + \left(\frac{\lambda}{2} + \partial_{44} h_i (\widehat{\mathrsfs{B}}_i(t)) \right) \operatorname{D}\bigg)\Phi_j$}
     \\
      \text{with \, $\Phi(0, z) = I_{2}$,}  \quad 
   U_{0,\mathrsfs{V}_j}(t) 
   = -2\gamma(t) \sum_{i=1}^2 p_i \bigg( \mathfrak{m}_{j,i}(t) \otimes H_{2,i}(\widehat{\mathrsfs{B}}_i(t)) 
   \\
   + H_{2,i}(\widehat{\mathrsfs{B}}_i(t)) \otimes \mathfrak{m}_{j,i}(t)\bigg) 
   +\frac{\gamma(t)^2}{d} \left(\sum_{i=1}^2 p_i \mathcal{I}_i(\widehat{\mathrsfs{B}}_i(t))(\lambda_j^{(i)} + (\omega_j^\top \mu_i)^2)\right) \operatorname{D}
      \\
\text{and}\quad
U_{0,\mathfrak{m}_{j,\ell}}(t) 
   =
   -2\gamma(t) \sum_{i=1}^2 p_i  \langle \omega_j, \mu_\ell\rangle \langle \omega_j, \mu_i\rangle \cdot H_{2,i}(\widehat{\mathrsfs{B}}_i(t)).
     \end{gather*}
    \end{mdframed}
    This alternative characterization requires to solve the ODE, 
    \begin{equation}\label{eq:Phi:ODE}
        \dot{\Phi}_j =2\gamma(t)\sum_{i=1}^2 p_i \bigg( \lambda_j^{(i)} H_{1,i}(\widehat{\mathrsfs{B}}_i(t))  + \left(\frac{\lambda}{2} + \partial_{44} h_i (\widehat{\mathrsfs{B}}_i(t)) \right) \operatorname{D}\bigg)\Phi_j,\quad \Phi_j(0)= \Id_{2}.
    \end{equation}
    We must now solve for $\mathfrak{m}_{j,\ell}(t)$ and $\mathrsfs{V}_j(t)$ in \eqref{eq:resolvent_form:M} and \eqref{eq:resolvent_form:S} respectively. Most of of the work towards this is to solve the matrix ODE \eqref{eq:Phi:ODE} which we express in matrix form:
\begin{equation}
    \dot{\Phi}_j =2\gamma(t) \sum_{i=1}^2 p_i \begin{bmatrix}
        \lambda_j^{(i)}\partial_{11} h_i(\widehat{\mathrsfs{B}}_i(t)) + \tfrac{\lambda}{2} + \partial_{44}h_i(\widehat{\mathrsfs{B}}_i(t))   & 0 \\
        \lambda_j^{(i)}\partial_{21} h_i(\widehat{\mathrsfs{B}}_i(t))&0
   \end{bmatrix}\Phi,\quad \Phi(0) = \Id_2.
\end{equation}
Hence, this reduces to the following system of first-order linear ODEs:
\begin{equation}
    \begin{aligned}
        \dot{\Phi}_{j,11} 
        & 
        = 2\gamma(t)\sum_{i=1}^2 p_i ( \lambda_j^{(i)}\partial_{11} h_i(\widehat{\mathrsfs{B}}_i(t)) + \tfrac{\lambda}{2} + \partial_{44}h_i(\widehat{\mathrsfs{B}}_i(t))) \Phi_{j,11} \quad \Phi_{j,11}(0) = 1
        \\
        \dot{\Phi}_{j,21}
        &
        = 2\gamma(t) (\sum_{i=1}^2p_i \lambda_j^{(i)}\partial_{21} h_i(\widehat{\mathrsfs{B}}_i(t)))\Phi_{j,11}, \quad \Phi_{j,21}(0) = 0,
        \\
        \dot{\Phi}_{j,12} 
        & 
        = 2\gamma(t) \sum_{i=1}^2 p_i (\lambda_j^{(i)}\partial_{11} h_i(\widehat{\mathrsfs{B}}_i(t)) + \tfrac{\lambda}{2} + \partial_{44}h_i(\widehat{\mathrsfs{B}}_i(t)))  \Phi_{j,12}, \quad \Phi_{j,12}(0) = 0,
        \\
        \dot{\Phi}_{j,22}
        &
        = 2\gamma(t) (\sum_{i=1}^2p_i \lambda_j^{(i)}\partial_{21} h_i(\widehat{\mathrsfs{B}}_i(t)))\Phi_{j,12}, \quad \Phi_{j,22}(0) = 1.
    \end{aligned}
\end{equation}
Using integration factors to solve the ODEs of $\Phi_{j,11}$ and $\Phi_{j,12}$, we obtain
\begin{equation}\label{eq:sol:Phi(t)}
    \begin{gathered}
    \Phi_{j,11}(t) = \exp\left(\int_0^t 2\gamma(s) \sum_{i=1}^2 p_i \bigg(\lambda_j^{(i)}\partial_{11} h_i(\widehat{\mathrsfs{B}}_i(s)) + \tfrac{\lambda}{2} + \partial_{44}h_i(\widehat{\mathrsfs{B}}_i(s))\bigg)\dif s\right),\\
     \Phi_{j,21}(t) = \int_0^t  2\gamma(s) \bigg(\sum_{i=1}^2 p_i \lambda_j^{(i)}\partial_{21} h_i(\widehat{\mathrsfs{B}}_i(s))\bigg)\Phi_{j,11} (s)\dif s,\quad
    \Phi_{j,22}(t) = 1,\quad \text{and}\quad \Phi_{j,12}(t) = 0.
    \end{gathered}
\end{equation}
Given that $\Phi_j \in \R^{2\times 2}$ then its inverse is given by
\begin{equation}
    \Phi_j^{-1}(t) = \frac{1}{\Phi_{j,11}(t)} \begin{bmatrix}
        1 & 0 \\
        - \Phi_{j,21}(t) & \Phi_{j,11}(t)
    \end{bmatrix}.
\end{equation}
It remains to compute the expressions in \eqref{eq:resolvent_form:M} and \eqref{eq:resolvent_form:S}. We have
\begin{equation}
    \begin{aligned}
        &\Phi_j^{-\top}(t) \mathrsfs{V}_j(0) \Phi_j^{-1}(t) \\
        &
        =
        \frac{1}{\Phi_{j,11}^2(t)} \begin{bmatrix}
            \substack{X_0^\top \omega_j\omega_j^\top X_0\\ - 2 \Phi_{j,21}(t) X_0^\top \omega_j\omega_j^\top X^{\star} \\+ \Phi_{j,21}^2(t) (X^{\star})^\top \omega_j\omega_j^\top X^{\star}}
            & \substack{\Phi_{j,11}(t) X_0^\top \omega_j\omega_j^\top X^{\star}\\ - \Phi_{j,11}(t) \Phi_{j,21}(t) (X^{\star})^\top\omega_j\omega_j^\top X^{\star}   }
            \\
            \\
            \substack{\Phi_{j,11}(t) X_0^\top \omega_j\omega_j^\top X^{\star} \\- \Phi_{j,21}(t) \Phi_{j,11}(t) (X^{\star})^\top \omega_j\omega_j^\top X^{\star}}& \Phi_{11}^2(t) (X^{\star})^\top \omega_j\omega_j^\top X^{\star}
        \end{bmatrix}.
    \end{aligned}
\end{equation}
 Now, it is easy to see that 
\begin{equation}
    \begin{aligned}
        \Phi_j(s) \Phi_j^{-1}(t) 
        &
        = 
        \begin{bmatrix}
            \frac{\Phi_{j,11}(s)}{\Phi_{j,11}(t)} & 0\\\frac{\Phi_{j,21}(s) - \Phi_{j,21}(t)}{\Phi_{j,11}(t)} & 1
        \end{bmatrix},
    \end{aligned}
\end{equation}
from which it follows that
\begin{equation*}
    \begin{aligned}
        & \frac{\gamma(t)^2}{d} \left(\sum_{i=1}^2 p_i \mathcal{I}_i(\widehat{\mathrsfs{B}}_i(t))(\lambda_j^{(i)} + (\omega_j^\top \mu_i)^2)\right)\Phi_j^{-\top}(t) \Phi_j^\top(s) 
        D \Phi_j(s) \Phi_j^{-1}(t) 
        \\
        &
        =
        \frac{\gamma(t)^2}{d} \frac{\Phi_{j,11}^2(s)}{\Phi_{j,11}^2(t)}\left(\sum_{i=1}^2 p_i \mathcal{I}_i(\widehat{\mathrsfs{B}}_i(t))(\lambda_j^{(i)} + (\omega_j^\top \mu_i)^2)\right) D
       ,
    \end{aligned}
\end{equation*}
and
\begin{equation}
    \begin{aligned}
        &\Phi_j^{-\top}(t) \Phi_j^\top(s) 
       ( \mathfrak{m}_{j,i}(s) \otimes H_{2,i}(\widehat{\mathrsfs{B}}_i(s))  )
        \Phi_j(s) \Phi_j^{-1}(t) 
        \\
        &=\partial_{13} h_i(\widehat{\mathrsfs{B}}_i(s)) \frac{\Phi_{j,11}(s)}{\Phi_{j,11}(t)} 
       \begin{bmatrix}
        \frac{\Phi_{j,11}(s)}{\Phi_{j,11}(t)}\mathfrak{m}_{j,i, 1}(s)  + \frac{\Phi_{j,21}(s) - \Phi_{j,21}(t)}{\Phi_{j,11}(t)} \mathfrak{m}_{j,i, 2}(s) & 0 \\
        \mathfrak{m}_{j,i,2}(s)&0
       \end{bmatrix},
    \end{aligned}
\end{equation}
where $\mathfrak{m}_{j,i,1}(t)$ and $\mathfrak{m}_{j,i,2}(t)$ denote the two elements composing $\mathfrak{m}_{j,i}(t) \in \R^2$. Similarly, we obtain
 \[
\begin{aligned}
    &
    \Phi_j^{-\top}(t) \Phi_j^\top(s) 
        (H_{2,i} (\widehat{\mathrsfs{B}}_i(s)) \otimes \mathfrak{m}_{j,i}(s))
        \Phi_j(s) \Phi_j^{-1}(t) 
        \\
        &=\partial_{13} h_i(\widehat{\mathrsfs{B}}_i(s)) \frac{\Phi_{j,11}(s)}{\Phi_{j,11}(t)} 
       \begin{bmatrix}
        \frac{\Phi_{j,11}(s)}{\Phi_{j,11}(t)}\mathfrak{m}_{j,i, 1}(s)  + \frac{\Phi_{j,21}(s) - \Phi_{j,21}(t)}{\Phi_{j,11}(t)} \mathfrak{m}_{j,i, 2}(s) & \mathfrak{m}_{j,i,2}(s) \\
        0&0
       \end{bmatrix}
\end{aligned}
\]
We obtain
\begin{equation} \label{eq:deep_appendix:S}
    \begin{aligned}
        &\mathrsfs{V}_j(t) 
        =
        \frac{1}{\Phi_{j,11}^2(t)} \begin{bmatrix}
            \substack{X_0^\top \omega_j\omega_j^\top X_0\\ - 2 \Phi_{j,21}(t) X_0^\top \omega_j\omega_j^\top X^{\star} \\+ \Phi_{j,21}^2(t) (X^{\star})^\top \omega_j\omega_j^\top X^{\star}}
            & \substack{\Phi_{j,11}(t) X_0^\top \omega_j\omega_j^\top X^{\star}\\ - \Phi_{j,11}(t) \Phi_{j,21}(t) (X^{\star})^\top\omega_j\omega_j^\top X^{\star}   }
            \\
            \substack{\Phi_{j,11}(t) X_0^\top \omega_j\omega_j^\top X^{\star} \\- \Phi_{j,21}(t) \Phi_{j,11}(t) (X^{\star})^\top \omega_j\omega_j^\top X^{\star}}& \Phi_{11}^2(t) (X^{\star})^\top \omega_j\omega_j^\top X^{\star}
        \end{bmatrix}
        \\
        &
        -
        \int_0^t  2\gamma(s) \sum_{i=1}^2 p_i\partial_{13} h_i(\widehat{\mathrsfs{B}}_i(s)) \frac{\Phi_{j,11}(s)}{\Phi_{j,11}(t)} 
        \\
        &
        \qquad 
        \times
        \begin{bmatrix}
         2\left(\frac{\Phi_{j,11}(s)}{\Phi_{j,11}(t)}\mathfrak{m}_{j,i, 1}(s)  + \frac{\Phi_{j,21}(s) - \Phi_{j,21}(t)}{\Phi_{j,11}(t)} \mathfrak{m}_{j,i, 2}(s)\right) & \mathfrak{m}_{j,i, 2}(s) \\
         \mathfrak{m}_{j,i,2}(s)&0
        \end{bmatrix}\l \dif s
        \\
        & 
        +
         \int_0^t   \frac{\gamma(s)^2}{d} \frac{\Phi_{j,11}^2(s)}{\Phi_{j,11}^2(t)}\left(\sum_{i=1}^2 p_i \mathcal{I}_i(\widehat{\mathrsfs{B}}_i(s))(\lambda_j^{(i)} + (\omega_j^\top \mu_i)^2)\right) \operatorname{D} \, \dif s. \nonumber
        \end{aligned}
\end{equation}
We also have
\begin{equation}
    \begin{aligned}
        &\Phi_j^{-\top}(t,z) \mathfrak{m}_{j,\ell}(0) \\
        &
        =
        \frac{1}{\Phi_{j,11}(t)} \begin{bmatrix}
            \substack{X_0^\top \omega_j \omega_j^\top \mu_{\ell}\\ - \Phi_{j,21}(t) (X^{\star})^{\top}  \omega_j \omega_j^\top \mu_{\ell} }
            \\
           \Phi_{j,11}(t) (X^{\star})^{\top} \omega_j \omega_j^\top\mu_{\ell}
        \end{bmatrix}.
    \end{aligned}
\end{equation}
We also have
\begin{equation}
    \begin{aligned}
        &\Phi_j^{-\top}(t) \Phi_j^\top(s)\bigg(  -2\gamma(s) \sum_{i=1}^2 p_i  \langle \omega_j, \mu_\ell\rangle \langle \omega_j, \mu_i\rangle \cdot H_{2,i}(\widehat{\mathrsfs{B}}_i(s)) \bigg)
        \\
        &
        =-2\gamma(s)\sum_{i=1}^2 p_i \frac{\Phi_{j,11}(s)}{\Phi_{j,11}(t)}  \langle \omega_j, \mu_\ell\rangle \langle \omega_j, \mu_i\rangle \cdot H_{2,i}(\widehat{\mathrsfs{B}}_i(s)).
    \end{aligned}
\end{equation}
Hence, we obtain
\begin{equation} \label{eq:deep_appendix:M}
    \begin{aligned}
        \mathfrak{m}_{j,\ell}(t) &=
        \frac{1}{\Phi_{j,11}(t)} \begin{bmatrix}
            \substack{X_0^\top \omega_j \omega_j^\top \mu_{\ell}\\ - \Phi_{j,21}(t) (X^{\star})^{\top}  \omega_j \omega_j^\top \mu_{\ell} }
            \\
            \\
           \Phi_{j,11}(t) (X^{\star})^{\top} \omega_j \omega_j^\top\mu_{\ell}
        \end{bmatrix}
        \\
        &
         -\int_0^t 2\gamma(s)\sum_{i=1}^2 p_i \frac{\Phi_{j,11}(s)}{\Phi_{j,11}(t)}  \langle \omega_j, \mu_\ell\rangle \langle \omega_j, \mu_i\rangle \cdot H_{2,i}(\widehat{\mathrsfs{B}}_i(s)) \dif s.
    \end{aligned}
\end{equation}
The result \eqref{eq:C_scalar} now immediately follows from summing over $j$ from $1$ to $d$. Indeed, we have for $\ell=1,2$
\begin{equation}
    \mathrsfs{B}_\ell(t) = \sum_{j=1}^d \lambda_{j}^{(\ell)} \mathrsfs{V}_j(t),\quad \mathfrak{m}_{\ell} = \sum_{j=1}^d \mathfrak{m}_{j,\ell} (t)\quad\text{and}\quad \widehat{\mathrsfs{B}}_{\ell,44}(t)  =  \sum_{j=1}^d \mathrsfs{V}_{j,11}(t).
\end{equation}
\begin{mdframed}[style=exampledefault]
    \textbf{Evolution of $\widehat{\mathrsfs{B}}_\ell(t)$ for $\ell=1,2$.}
    \begin{equation} \label{eq:C_scalar}
        \begin{aligned}
            \widehat{\mathrsfs{B}}_\ell(t) 
            &
            = \begin{bmatrix}
                \mathrsfs{B}_\ell(t) & \mathfrak{m}_\ell(t) & 0
                \\
                \mathfrak{m}_\ell(t)^\top & \|\mu_\ell\|^2 & 0\\
                0 & 0 &\widehat{\mathrsfs{B}}_{\ell,44}(t)
            \end{bmatrix},\quad  \mathrsfs{B}_\ell(t)=\begin{bmatrix}
            \mathrsfs{B}_{\ell,11}(t) & \mathrsfs{B}_{\ell,12}(t)
            \\
            \mathrsfs{B}_{\ell,21}(t) & \mathrsfs{B}_{\ell,22}(t)
        \end{bmatrix}
            \\
            \text{where} \quad 
            \mathrsfs{B}_{\ell,11}(t) 
            &
            = 
            X_0^\top \tfrac{K_\ell}{\Phi_{11}^{2}(t, \mathcal{K})}X_0 - 2 X_0^\top \tfrac{K_\ell\Phi_{21}(t,  \mathcal{K})}{\Phi_{11}^{2}(t,  \mathcal{K})}  X^{\star} + (X^{\star})^\top \tfrac{K_\ell \Phi_{21}^2(t,  \mathcal{K})}{\Phi_{11}^{2}(t, \mathcal{K})}  X^{\star}  
            \\
            &
            -
        4\int_0^t \gamma(s) \sum_{i=1}^2 p_i\partial_{13} h_i(\widehat{\mathrsfs{B}}_i(s)) \sum_{j=1}^d\lambda_j^{(\ell)}  \frac{\Phi_{j,11}(s)}{\Phi_{j,11}(t)} 
        \bigg(\frac{\Phi_{j,11}(s)}{\Phi_{j,11}(t)}\mathfrak{m}_{j,i, 1}(s)  
        \\
        &
        \qquad
        + \frac{\Phi_{j,21}(s) - \Phi_{j,21}(t)}{\Phi_{j,11}(t)} \mathfrak{m}_{j,i, 2}(s)\bigg) \dif s
            \\
            &
            + \int_0^t \frac{\gamma(s)^2}{d} \sum_{i=1}^2 p_i \mathcal{I}_i(\widehat{\mathrsfs{B}}_i(s)) 
               \Tr(( K_i + \mu_i \mu_i^\top )K_{\ell}  \tfrac{\Phi_{11}^2(s, \mathcal{K})}{\Phi_{11}^2(t, \mathcal{K})} ) \;\dif s,
            \\
            \mathrsfs{B}_{\ell, 12}(t)
            &
            =
            X_0^\top \tfrac{K_{\ell}}{\Phi_{11}(t,\mathcal{K})} X^{\star} - (X^{\star})^\top \tfrac{K_{\ell}\Phi_{21}(t,\mathcal{K})}{\Phi_{11}(t,\mathcal{K})} X^{\star} 
            \\
            &
            \qquad
            -
            \int_0^t  2\gamma(s) \sum_{i=1}^2 p_i\partial_{13} h_i(\widehat{\mathrsfs{B}}_i(s)) \sum_{j=1}^d \lambda_{j}^{(\ell)} \frac{\Phi_{j,11}(s)}{\Phi_{j,11}(t)} 
             \mathfrak{m}_{j,i, 2}(s) \dif s,
            \\
            \mathfrak{m}_{\ell,1}(t) &= X_0^\top \Phi_{11}^{-1}(t, \mathcal{K})\mu_{\ell} - (X^{\star})^\top \tfrac{\Phi_{21}(t, \mathcal{K})}{\Phi_{11}(t, \mathcal{K})}\mu_{\ell} 
            \\
            &
            - 2 \int_0^t \gamma(s) \sum_{i=1}^2 p_i \partial_{13} h_i(\widehat{\mathrsfs{B}}_i(s))\cdot \mu_i^\top \tfrac{\Phi_{11}(s,\mathcal{K})}{\Phi_{11}(t,\mathcal{K})}\mu_{\ell} \;\dif s,
            \\
            \mathrsfs{B}_{\ell, 44}(t) 
            &
            = 
            X_0^\top \Phi_{11}^{-2}(t, \mathcal{K})X_0 - 2 X_0^\top \tfrac{\Phi_{21}(t,  \mathcal{K})}{\Phi_{11}^{2}(t,  \mathcal{K})}  X^{\star} + (X^{\star})^\top \tfrac{\Phi_{21}^2(t,  \mathcal{K})}{\Phi_{11}^{2}(t, \mathcal{K})}  X^{\star}  
            \\
            &
            -
        4\int_0^t \gamma(s) \sum_{i=1}^2 p_i\partial_{13} h_i(\widehat{\mathrsfs{B}}_i(s)) \sum_{j=1}^d \frac{\Phi_{j,11}(s)}{\Phi_{j,11}(t)} 
        \bigg(\frac{\Phi_{j,11}(s)}{\Phi_{j,11}(t)}\mathfrak{m}_{j,i, 1}(s)  
        \\
        &
        \qquad
        + \frac{\Phi_{j,21}(s) - \Phi_{j,21}(t)}{\Phi_{j,11}(t)} \mathfrak{m}_{j,i, 2}(s)\bigg) \dif s
            \\
            &
            + \int_0^t \frac{\gamma(s)^2}{d} \sum_{i=1}^2 p_i \mathcal{I}_i(\widehat{\mathrsfs{B}}_i(s)) 
               \Tr(( K_i + \mu_i \mu_i^\top ) \tfrac{\Phi_{11}^2(s, \mathcal{K})}{\Phi_{11}^2(t, \mathcal{K})} ) \;\dif s,
            \\
            \mathrsfs{B}_{\ell,21}(t) 
            &
            =
            \mathrsfs{B}_{\ell,12}(t), \quad
            \mathrsfs{B}_{\ell, 22}(t) = (X^{\star})^\top K_{\ell} X^{\star}\quad \text{and} \quad \mathfrak{m}_{\ell,2} = (X^{\star})^\top \mu_{\ell}
        \end{aligned}
    \end{equation}
    The function $\Phi_{11}(t, \mathcal{K})$ and $\Phi_{21}(t,\mathcal{K})$, by solving a differential equation, are given by 
    \begin{equation}\label{eq:Phi:def:solved}
    \begin{gathered}
        \Phi_{11}(t,\mathcal{K}) = \exp\left(\int_0^t 2\gamma(s)\sum_{i=1}^2 p_i \bigg( K_i\partial_{11} h_i(\widehat{\mathrsfs{B}}_i(s)) + \frac{\lambda}{2} + \partial_{44}h_i(\widehat{\mathrsfs{B}}_i(s))\bigg)\dif s\right),\\
        \text{and}\quad\Phi_{21}(t,\mathcal{K}) = \int_0^t  2\gamma(s) \bigg(\sum_{i=1}^2p_i K_i\partial_{21} h_i(\widehat{\mathrsfs{B}}_i(s))\bigg)\Phi_{11}(s,\mathcal{K}) \dif s. 
        \end{gathered}
    \end{equation}
    \end{mdframed}
    Let $  \mathrsfs{R}^{{\operatorname{adv}}}(t)$ denote the deterministic equivalent of $\mathcal{R}^{{\operatorname{adv}}}$. Given that $\mathcal{R}^{{\operatorname{adv}}}(X) = \sum_{i=1}^2 p_i\cdot h_i(\widehat{B}_i(W))$ then we obtain
    \[
        \mathrsfs{R}^{{\operatorname{adv}}}(t) = \sum_{i=1}^2 p_i \cdot h_i( \widehat{\mathrsfs{B}}_i(t)).
     \]
Analogously to the spectral decomposition in \eqref{eq:app:decomps:S:M}, we may use a similar argument to decompose $\mathcal{Z}(t,z)$ as follows
\begin{equation}\label{eq:sum:Z:app}
    \mathcal{Z}(t,z) = \sum_{j=1}^d \frac{\mathscr{Z}_j(t)}{\prod_{\ell=1}^2(\lambda^{(\ell)}_j - z_\ell)},
\end{equation}
where 
\begin{equation}
    \mathscr{Z}_j(t) \defas \begin{bmatrix}
        \mathrsfs{V}_j(t) & \mathfrak{M}_j(t) \\
        \mathfrak{M}_j(t)^\top & \mu^\top \omega_j \omega_j^\top \mu
    \end{bmatrix}\quad \text{with}\quad\mathfrak{M}_j(t) \defas \begin{bmatrix}
        \mathfrak{m}_{j,1} &  \mathfrak{m}_{j,2}
    \end{bmatrix}.  
\end{equation}
Using \eqref{eq:ODE_resolvent_2:app}, we obtain the differential equation
\begin{equation}
    \frac{\dif\mathscr{Z}_j(t) }{\dif t}= - \mathscr{Z}_j(t) \mathscr{B}_j(t) -  \mathscr{B}_j(t)^\top  \mathscr{Z}_j(t)+ \mathscr{G}_j(t),
\end{equation}
where we define for $\operatorname{e}_i $ the standard basis vectors of $\R^2$
\begin{equation}
    \begin{gathered}
    \mathscr{B}_j(t) \defas \begin{bmatrix}
        \mathscr{A}_j(t) & 0 \\
        \mathscr{C}(t) & 0 
    \end{bmatrix},\quad  \mathscr{C}(t) = 2\gamma(t)\sum_{i=1}^2 p_i \operatorname{e}_i \otimes H_{2,i}(\widehat{\mathrsfs{B}}_i(t)),
    \\
    \mathscr{A}_j(t) \defas 2\gamma(t)\sum_{i=1}^2 p_i \bigg( \lambda_j^{(i)} H_{1,i}(\widehat{\mathrsfs{B}}_i(t))  + \left(\frac{\lambda}{2} + \partial_{44} h_i (\widehat{\mathrsfs{B}}_i(t)) \right) \operatorname{D}\bigg),\\
    \text{and}\quad \mathscr{G}_j(t) \defas \begin{bmatrix}
        \frac{\gamma(t)^2}{d} \left(\sum_{i=1}^2 p_i \mathcal{I}_i(\widehat{\mathrsfs{B}}_i(t))(\lambda_j^{(i)} + (\omega_j^\top \mu_i)^2)\right)\operatorname{D} & 0 \\
        0 & 0 
    \end{bmatrix}.
\end{gathered}
\end{equation}
Analogously to \eqref{eq:app:base:ODE}, we define the fundamental ODE
\begin{equation}
    \dot{\Xi}_j = \mathscr{B}_j(t)\Xi_j,\quad \Xi_j(0)= \Id_{4}. 
 \end{equation}
 Given that $\mathscr{B}_j(t)$ is a block lower triangular matrix then we obtain from the initial conditions
 \begin{equation}
    \Xi_j (t) = \begin{bmatrix}
        \Phi_j(t) & 0 \\
        \Upsilon_j(t) & \Id_2
    \end{bmatrix},\quad\text{with }\quad \Upsilon_j(t) = \int_0^t \mathscr{C}_j(s) \Phi_j(s) \dif s,
 \end{equation}
 and where $\Phi_j$ is defined in \eqref{eq:sol:Phi(t)}. Using a similar argument as \eqref{eq:app:similar:arg:ODE}, we obtain
 \begin{equation}
    \mathscr{Z}_j(t) = \Xi_j(t)^{-\top}\mathscr{Z}_j(0)\Xi_j(t)^{-1} + \int_0^t \Xi_j(t)^{-\top}\Xi_j(s)^{\top}\mathscr{G}_j(s)\Xi_j(s)\Xi_j(t)^{-1} \dif s.
 \end{equation}
 Now, setting $x_j \defas \hat{X}_0^\top \omega_j \in \R^2$ and $q_j \defas \mu^\top \omega_j \in \R^2$ then we have
 \begin{equation}
    \mathrsfs{V}_j(0) = x_jx_j^\top, \mathscr{M}_j(t) = x_jq_j^\top\quad\text{and}\quad \mu^\top \omega_j \omega_j^\top \mu = q_j q_j^\top.
 \end{equation}
 It thus follows that 
 \begin{equation}
    \mathscr{Z}_j(0) = \begin{bmatrix}
        x_j \\
        q_j
    \end{bmatrix}\begin{bmatrix}
        x_j \\
        q_j
    \end{bmatrix}^\top \succeq 0.
 \end{equation}
 By Assumption~\ref{assumption:fisher}, it is clear that for $i=1,2$, $I_i(\mathrsfs{B}(s)) \geq 0$ for $0\leq s\leq t$ and every $t\geq 0$. Since $\mathscr{Z}_j(0)$ are positive semidefinite (psd) and congruence preserves positive semidefiniteness, it follows that $\mathscr{Z}_j(t)$ is psd for every $1\leq j \leq d$ and $t\geq 0$.
 
 Since the diagonal block matrices of a psd matrix are psd, it follows that $\mathrsfs{V}_j(t)$ is psd for all $1\leq j \leq d$ and $t\geq 0$ and thus $\mathrsfs{B}_{\ell}(t)$ is psd for all $t\geq 0$ and $\ell = 1,2$. Since the diagonal entries of a psd matrix are nonnegative, it follows that $\mathrsfs{V}_{j,11}(t) \geq 0$ for all $1\leq j \leq d$ and $t\geq 0$ and thus $\mathrsfs{B}_{\ell, 44}(t) \geq 0$ for every $t \geq 0$ and $\ell=1,2$.

 Now, define $ \mathscr{Z}(t) = \oint  \mathcal{Z}(t) \Dif z$ then by \eqref{eq:sum:Z:app} we have $\mathscr{Z}(t)  = \sum_{j=1}^d \mathscr{Z}_j(t)$ which implies $\mathscr{Z}(t) $ is psd for all $t\geq 0$. Every principal submatrix of $\mathscr{Z}(t) $ is then psd which implies the following two matrices obtained by removing the third column and row to obtain $\mathscr{Z}_2(t)$ and by removing the fourth column and row to obtain $\mathscr{Z}_1(t)$ are psd:
 \begin{equation}
    \mathscr{Z}_i(t) \defas \begin{bmatrix}
        \sum_{j=1}^d \mathrsfs{V}_j(t) & \mathfrak{m}_i \\
        \mathfrak{m}_i^\top & \|\mu_i\|^2
    \end{bmatrix}\quad \text{for}\quad i=1,2.
 \end{equation}
 When $\mu_1\not = 0$ and $\mu_2 \not =0$ then from the Schur complement, we obtain the following inequality for $i=1,2$
\begin{equation}\label{eq:upper bound:norm:m(t)}
    \|\mathfrak{m}_i\|^2 \leq \|\mu_i\|^2 \Tr(\sum_{j=1}^d \mathrsfs{V}_j(t)) \leq \|\mu_i\|^2 \mathrsfs{N}(t).
\end{equation}
 We close off this section by providing the hard label version of the system of ODEs \eqref{eq:ODE:V_i(t):def} which simplifies to
 \begin{equation}\label{eq:odes:hard:label}
     \begin{aligned}
         &\frac{\dif \mathrsfs{V}_{j,11}(t)}{\dif t} = -4\gamma(t) \sum_{i=1}^2 p_i \bigg( 
         \left( \lambda_j^{(i)}  \partial_{11}h_i(\widehat{\mathrsfs{B}}_i(t))+ \frac{ \lambda}{2} + \partial_{44} h_i (\widehat{\mathrsfs{B}}_i(t) )\right) \mathrsfs{V}_{j,11}(t)  
         \\
         &
         \quad
         +  \partial_{13}h_i(\widehat{\mathrsfs{B}}_i(t)) \mathfrak{m}_{j,i, 1}(t) \bigg)
         +\frac{\gamma(t)^2}{d} \sum_{i=1}^2 p_i \mathcal{I}_i(\widehat{\mathrsfs{B}}_i(t))(\lambda_j^{(i)} + (\omega_j^\top \mu_i)^2),
         \\
         &
         \frac{\dif \mathfrak{m}_{j,\ell, 1}(t)}{\dif t}  = -2\gamma(t) \sum_{i=1}^2 p_i \bigg( \left( \lambda_j^{(i)}\partial_{11}h_i(\widehat{\mathrsfs{B}}_i(t)) + \frac{\lambda}{2} + \partial_{44} h_i (\widehat{\mathrsfs{B}}_i(t)) \right)   \mathfrak{m}_{j,\ell, 1}(t) 
         \\
         &
         \quad +\partial_{13}h_i(\widehat{\mathrsfs{B}}_i(t))\langle \omega_j, \mu_\ell\rangle \langle \omega_j, \mu_i\rangle \bigg),
     \end{aligned}
 \end{equation}
 with initial condition $\mathrsfs{V}_{j,11}(0) = X_0^\top \omega_j \omega_j^\top X_0$ and $\mathfrak{m}_{j,\ell, 1}(0) = X_0^\top \omega_j \omega_j^\top\mu_{\ell}$.

\section{Polyak Stepsize}\label{sect:Polyakstepsize}
In this section, we determine the learning rate threshold for which the functions $\mathscr{D}^2(X)\defas\|X-X^\star\|^2$ and $\mathscr{D}^{2, \Adv}(X)\defas\|X- X^{\star,\Adv}\|^2$ are decreasing when $d\to \infty$ using their deterministic equivalents. We will work with the regularization parameter $\lambda = 0$ for simplicity. Our analysis will build upon the system of ODEs \eqref{eq:ODE:V_i(t):def}. The deterministic equivalent for the distance to least squares optimality $\|X-X^\star\|^2$ is given by
\begin{equation}
    \mathrsfs{D}^2(t)  = \sum_{i=1}^d V_{j,11}(t) - 2  V_{j,12}(t) +  V_{j,22}(t).
\end{equation}
Using a similar argument as in the proof of Proposition~\ref{prop:nonexplosiveness}, we can show that $\mathrsfs{D}^2(t)$ solves the following ODE.
\begin{corollary}
    $\mathrsfs{D}^2(t)$ solves the differential equation
    \begin{equation}\label{eq:D^2(t):ode}
        \begin{aligned}
    \frac{\dif}{\dif t}\mathrsfs{D}^2(t)
        &=
        -2\gamma(t)\sum_{i=1}^2p_i  Q (\widehat{\mathrsfs{B}}_i(t))+ \frac{\gamma(t)^2}{d}\sum_{i=1}^2 p_i\Tr(K_i + \mu_i\mu_i^\top )  \E_{v,\epsilon}[f_i'(q_{t,i}, \epsilon)^2] ,
        \end{aligned}
    \end{equation}
    where 
    \begin{equation}
        \begin{aligned}
           Q (\widehat{\mathrsfs{B}}_i(t)) =\E_{v,\epsilon}[f_i'(q_{t,i}, \epsilon) (x_{t,i} - x_i^\star)] 
  - \delta \frac{\mathrsfs{A}(t)}{\sqrt{\widehat{\mathrsfs{B}}_{i,44}(t)}}\E_{v,\epsilon}[ u_{t,i}f_i'(q_{t,i}, \epsilon) ],
        \end{aligned}
    \end{equation}
    for
    \begin{equation}
        \begin{aligned}
        q_{t,i} &=x_{t,i}\oplus x_i^{\star}= \sqrt{\mathrsfs{B}_i(t)} v  +\mathfrak{m}_i(t) + \begin{bmatrix}
            \delta u_{t,i} \sqrt{\widehat{\mathrsfs{B}}_{i,44}(t)} & 0
        \end{bmatrix}^\top,
        \\
         u_{t,i} &= \operatorname{argmax}_{|u|\leq 1} f_i(\sqrt{\mathrsfs{B}_i(t)} v  +\mathfrak{m}_i(t) + \begin{bmatrix}
            \delta u \sqrt{\widehat{\mathrsfs{B}}_{i,44}(t)} & 0
        \end{bmatrix}^\top, \epsilon),
    \end{aligned}
    \end{equation}
    for $v\sim \mathcal{N}(0, \Id_2)$. Here $\mathrsfs{A}(t)$ is the deterministic equivalent of $X_{\lfloor td \rfloor}^\top X^\star$ and solves the following ODE
    \begin{equation}
        \begin{aligned}
            \frac{\dif}{\dif t}\mathrsfs{A}(t)
        &
        =
        -\gamma(t)\sum_{i=1}^2p_i \bigg(\E_{v,\epsilon}[f_i'(q_{t,i}, \epsilon) x_i^\star] 
  + \delta\frac{\mathrsfs{A}(t)}{\sqrt{\widehat{\mathrsfs{B}}_{i,44}(t)}}\E_{v,\epsilon}[ u_{t,i}f_i'(q_{t,i}, \epsilon) ]\bigg).
        \end{aligned}
    \end{equation}
\end{corollary}
Given \eqref{eq:D^2(t):ode}, the exact threshold for monotone decrease of $\|X- X^\star\|^2$ for SGD (i.e. $\dif \mathrsfs{D}^2(t) < 0$) which we call $\gamma^{\operatorname{Stable}, \star}$ is given by 
\begin{equation}\label{eq:ls_gr:stab:det}
    \gamma(t) < \gamma_{\lfloor td \rfloor}^{\operatorname{Stable}, \star} = \frac{\sum_{i=1}^2p_i  Q (\widehat{\mathrsfs{B}}_i(t))}{\sum_{i=1}^2 p_i\frac{\Tr(K_i + \mu_i\mu_i^\top )}{2d}  \E_{v,\epsilon}[f_i'(q_{t,i}, \epsilon)^2] }.
\end{equation}
Given that adversarial training is an adaptive algorithm which depends on the sampled data point at each step, it follows that $\gamma_{\lfloor td \rfloor}^{\operatorname{Stable}, \star}$ is adaptive. Translating this quantity back to SGD, we obtain the threshold for stability with respect to descent optimality towards $X^\star$
\begin{equation}\label{eq:ls_gr:stab:sgd}
    \gamma_{k}^{\operatorname{Stable}, \star} = \frac{\langle  X_k - X^{\star}, \nabla \mathcal{R}^{\Adv}(X_k) \rangle  }{\sum_{i=1}^2 p_i\frac{\Tr(K_i + \mu_i\mu_i^\top )}{2d} \E_{v,\epsilon}[f_i'(g_{k,i}, \epsilon_{k+1})^2 ]}.
\end{equation}
This stability threshold corresponds to the adversarial version of the results presented in Section~3 and Appendix E of \citep{collinswoodfin2024highline} and should be compared to the greedy learning rate $\gamma^{\operatorname{Polyak}, \star}(t)$ that maximizes the decrease of $\mathrsfs{D}^2(t)$ at each iteration: $\gamma^{\operatorname{Polyak}, \star}(t) \in \operatorname{argmin}_{\gamma} \dif \mathrsfs{D}^2(t)$. Analogously to \citep{collinswoodfin2024highline}, we call this learning rate the \emph{Polyak stepsize}. Solving for $\gamma_k$ and $\gamma(t)$ respectively in \eqref{eq:ls_gr:stab:det} and \eqref{eq:ls_gr:stab:sgd}, we obtain the following closed forms for the Polyak learning rate
\begin{equation}
    \gamma_{k}^{\operatorname{Polyak}, *} = \tfrac{1}{2}\gamma_{k}^{\operatorname{Stable}, *}\quad \text{with deterministic equivalent} \quad \gamma^{\operatorname{Polyak}, *}(t) = \tfrac{1}{2}\gamma^{\operatorname{Stable}, *}(t).
\end{equation}
Setting $\delta  =0$ recovers the results in \citep{collinswoodfin2024highline}. Hence, $\ell_2$-adversarial training preserves the ratio of the Polyak learning rate to stable learning rate with respect to $X^\star$. We also obtain a threshold for descent of $\|X-X^{\star, \Adv}\|^2$ whose deterministic equivalent is given by
\begin{equation}
    \mathrsfs{D}^{2,\Adv}(t)  = \sum_{i=1}^d V_{j,11}(t) - 2   \mathrsfs{U}_j(t)+ \|X^{\star,\operatorname{Adv}}\|^2,
\end{equation}
where $ \mathrsfs{U}_j(t)$ is the deterministic equivalent of $X_{\lfloor td\rfloor}^\top \omega_j \omega_j^\top X^{\star,\operatorname{Adv}}$. Using a similar argument as in the proof of Proposition~\ref{prop:nonexplosiveness}, we can show that $  \mathrsfs{D}^{2,\Adv}(t) $ solves the following ODE.
\begin{corollary}\label{cor:D^2_adv(t):ode}
    $  \mathrsfs{D}^{2,\Adv}(t) $ solves the differential equation
    \begin{equation}\label{eq:D^2_adv(t):ode}
        \begin{aligned}
    \frac{\dif}{\dif t}  \mathrsfs{D}^{2,\Adv}(t) 
        &=
        -2\gamma(t)\sum_{i=1}^2p_i  Q^{\Adv} (\widehat{\mathrsfs{B}}_i(t))+ \frac{\gamma(t)^2}{d}\sum_{i=1}^2 p_i\Tr(K_i + \mu_i\mu_i^\top )  \E_{v,\epsilon}[f_i'(q_{t,i}, \epsilon)^2] ,
        \end{aligned}
    \end{equation}
    where
    \begin{equation}
        \begin{aligned}
           Q^{\Adv} (\widehat{\mathrsfs{B}}_i^\Adv(t)) &=\E_{w,\epsilon}[f_i'(q_{t,i}, \epsilon) (x_{t,i} - x_i^{\star,\Adv})]   - \delta \frac{\mathrsfs{A}^{\Adv}(t)}{\sqrt{\widehat{\mathrsfs{B}}_{i,44}(t)}}\E_{v,\epsilon}[ u_{t,i}f_i'(q_{t,i}, \epsilon) ],
        \end{aligned}
    \end{equation}
    for
    \begin{equation}
        \begin{aligned}
        q_{t,i} &= \sqrt{\mathrsfs{B}_i(t)} v  +\mathfrak{m}_i(t) + \begin{bmatrix}
            \delta u_{t,i} \sqrt{\widehat{\mathrsfs{B}}_{i,44}(t)} & 0
        \end{bmatrix}^\top,
        \\
         u_{t,i} &= \operatorname{argmax}_{|u|\leq 1} f_i(\sqrt{\mathrsfs{B}_i(t)} v  +\mathfrak{m}_i(t) + \begin{bmatrix}
            \delta u \sqrt{\widehat{\mathrsfs{B}}_{i,44}(t)} & 0
        \end{bmatrix}^\top, \epsilon),
    \end{aligned}
    \end{equation}
    for $w = v \oplus y \sim \mathcal{N}(0, \Id_3)$ for $v\in \R^2$ and $y\in \R$, where 
    \begin{equation}
        x_{t,i}\oplus x_i^{\star} \oplus x_i^{\star, \Adv} = \sqrt{ \mathrsfs{B}^{\Adv}_i(t)}w + \mathfrak{m}_i(t)\oplus( \mu_i^\top  X^{\star, \Adv} ) + \begin{bmatrix}
            \delta u_{t,i} \sqrt{\widehat{\mathrsfs{B}}_{i,44}(t)} & 0 & 0 
        \end{bmatrix}^\top.
    \end{equation} 
    Here $ \mathrsfs{B}^{\Adv}_i(t)$ is defined in \eqref{def:det:equiv:covariance:adv} and $\mathrsfs{A}^{\Adv}(t)$ is the deterministic equivalent of $X_{\lfloor td \rfloor}^\top X^{\star, \Adv}$ which solves the following ODE
    \begin{equation}\label{eq:AlignAdvGR:ode}
        \begin{aligned}
            \frac{\dif}{\dif t}\mathrsfs{A}^{\Adv}(t)
        &
        =
        -\gamma(t)\sum_{i=1}^2p_i \bigg(\E_{w,\epsilon}[f_i'(q_{t,i}, \epsilon) x^{\star, \Adv}] 
  + \delta\frac{\mathrsfs{A}^{\Adv}(t)}{\sqrt{\widehat{\mathrsfs{B}}_{i,44}(t)}}\E_{v,\epsilon}[ u_{t,i}f_i'(q_{t,i}, \epsilon) ]\bigg).
        \end{aligned}
    \end{equation}
\end{corollary}
Given \eqref{eq:D^2_adv(t):ode}, the exact threshold for monotone decrease
of $\|X- X^{\star, \Adv}\|^2$ for SGD (i.e. $\dif \mathrsfs{D}^{2, \Adv}(t) < 0$) is given by 
\begin{equation}\label{eq:adv_gr:stab:det}
    \gamma(t) < \gamma_{\lfloor td \rfloor}^{\operatorname{Stable}, \Adv} = \frac{\sum_{i=1}^2p_i  Q^{\Adv} (\widehat{\mathrsfs{B}}_i^\Adv(t))}{\sum_{i=1}^2 p_i\frac{\Tr(K_i + \mu_i\mu_i^\top )}{2d}  \E_{v,\epsilon}[f_i'(q_{t,i}, \epsilon)^2]}.
\end{equation}
 Translating this quantity back to SGD, we obtain the threshold for stability for descent towards $X^{\star, \Adv}$:
\begin{equation}\label{eq:adv_gr:stab:sgd}
    \begin{aligned}
    \gamma_{k}^{\operatorname{Stable}, \Adv} 
    &=
    \frac{\langle  X_k - X^{\star, \Adv}, \nabla \mathcal{R}^\Adv(X_k) \rangle  }{\sum_{i=1}^2 p_i\frac{\Tr(K_i + \mu_i\mu_i^\top )}{2d} \E_{v,\epsilon}[f_i'(g_{k,i}, \epsilon_{k+1})^2 ]}
    \\
    &
    =
    \frac{\sum_{i=1}^2p_i  \E_{a,\epsilon}\left[f_i'(g_{k,i}, \epsilon)(X_k - X^{\star})^\top a_{k+1, i}\right] }{\sum_{i=1}^2 p_i\frac{\Tr(K_i + \mu_i\mu_i^\top )}{2d} \E_{v,\epsilon}[f_i'(g_{k,i}, \epsilon_{k+1})^2 ]}
    \\
    &
    \quad
    +
    \frac{\sum_{i=1}^2p_i \E_{a,\epsilon}\left[f_i'(g_{k,i}, \epsilon) (X^{\star} - X^{\star, \Adv})^\top a_{k+1, i} \right] }{\sum_{i=1}^2 p_i\frac{\Tr(K_i + \mu_i\mu_i^\top )}{2d} \E_{v,\epsilon}[f_i'(g_{k,i}, \epsilon_{k+1})^2 ]}
    \\
    &
    \quad
    +
    \frac{\delta\left(\| X_k \| -  \langle  X^{\star, \Adv}, \tfrac{X_k}{\|X_k\|}\rangle\right)\sum_{i=1}^2p_i  \E_{a,\epsilon}\left[f_i'(g_{k,i}, \epsilon)  s_{k,i} \right] }{\sum_{i=1}^2 p_i\frac{\Tr(K_i + \mu_i\mu_i^\top )}{2d} \E_{v,\epsilon}[f_i'(g_{k,i}, \epsilon_{k+1})^2 ]}.
    \end{aligned}
\end{equation}
This stability threshold admits a natural decomposition into three contributions. The first term is the standard descent contribution towards $X^{\star}$. The second term is the correction for the shift in ground truth between $X^\star$ and $X^{\star, \Adv}$. The third term is the contribution from the adversarial perturbation of the input. 

A natural question is how the stability threshold differs for SGD on the population risk $\mathcal{R}(X_k)$ versus the adversarial risk $\mathcal{R}^{\operatorname{Adv}}(X_k)$. We focus on the case of SGD ran on the least squares loss in the soft label setting with anisotropic centered Gaussian data $a\sim \mathcal{N}(0, K)$ which corresponds to the following minimization problem
\begin{equation} \label{eq:loss:least squares}
    \min_{X \in \mathbb{R}^{d}} \Big \{ \mathcal{R}(X) = \tfrac{1}{2} \EE_{a,\epsilon} [((X- X^\star)^\top a - \epsilon)^2] = \tfrac{1}{2}\eta^2 + \tfrac{1}{2}(X- X^\star)^\top K (X- X^\star)\Big\}.
    \end{equation}
Translating the relevant deterministic equivalents to their non-deterministic counterparts, the Polyak stepsize for \eqref{eq:loss:least squares} which we call $ \mathfrak{g}^{\operatorname{Polyak}}$ corresponds to \citep{collinswoodfin2024highline}
\begin{equation}
    \mathfrak{g}_{k}^{\operatorname{Polyak}} = \frac{2 \mathcal{R}(X_k) - \eta^2}{\frac{2\Tr(K)}{d}\mathcal{R}(X_k)} \text{ and on noiseless least squares,}\quad  \mathfrak{g}_{k}^{\operatorname{Polyak}}=\frac{d}{\Tr(K)}.
\end{equation}
 In the case of $\ell_2$-adversarial least squares, we define $\gamma_k^{\operatorname{eff}}$ and $\lambda_k^{\operatorname{eff}}$ as the effective learning rate and effective regularization
\begin{equation}
    \gamma_k^{\operatorname{eff}} = \gamma_k\bigg(1 + \delta \sqrt{\frac{2}{\pi}}\frac{\|X_k\|}{\sqrt{2\mathcal{R}(X_k)}}\bigg)\quad\text{and}\quad \lambda_k^{\operatorname{eff}} = \frac{\delta^2 + \delta \sqrt{\frac{2}{\pi}}\frac{\sqrt{2\mathcal{R}(X_k)}}{\|X_k\|}}{1 + \delta \sqrt{\frac{2}{\pi}}\frac{\|X_k\|}{\sqrt{2\mathcal{R}(X_k)}}},
\end{equation} 
with $\gamma_k$ corresponding to the learning rate defined in Assumption~\ref{assumption:lr}. See Section~\ref{sect:AdvHSGD:vs:HSGD} for more details on $\gamma_k^{\operatorname{eff}}$ and $\lambda_k^{\operatorname{eff}}$. Let $\tilde{\gamma}^{\operatorname{eff}}(t)$ and $\tilde{\lambda}^{\operatorname{eff}}(t)$ denote the deterministic equivalents of $\gamma_k^{\operatorname{eff}}$ and $\lambda_k^{\operatorname{eff}}$ written in terms of $\mathrsfs{R}(t)$ and $\widehat{\mathrsfs{B}}_{44}(t)$. Then the descent to least squares optimality for $\ell_2$-adversarial least squares solves the following ODE:
\begin{equation}\label{eq:ls:descent:optim:ode}
    \begin{aligned}
        \frac{\dif}{\dif t}\mathrsfs{D}^2(t)
        &=
        -2\tilde{\gamma}^{\operatorname{eff}}(t)\left(2\mathrsfs{R}(t) - \eta^2 + \tilde{\lambda}^{\operatorname{eff}}(t)(\widehat{\mathrsfs{B}}_{44}(t) - \mathrsfs{A}(t))\right) +\frac{2\gamma^2(t)}{d}\Tr(K)\mathrsfs{R}^{\operatorname{Adv}}(t).
    \end{aligned}
\end{equation} 
where $\mathrsfs{R}^{\operatorname{Adv}}(t)$ is the deterministic equivalent of $\mathcal{R}^{\operatorname{Adv}}(X_k)$. Noting the relation
\begin{equation}
    \mathrsfs{R}^{{\operatorname{adv}}}(t) = \frac{\tilde{\gamma}^{\operatorname{eff}}(t)}{\gamma(t)} \left(\mathrsfs{R}(t) + \frac{\tilde{\lambda}^{\operatorname{eff}}(t)}{2}\widehat{\mathrsfs{B}}_{44}(t)\right),
\end{equation}
then translating these quantities back to SGD, the Polyak stepsize with respect to $X^\star$ is naturally adaptive and given by
\begin{equation}\label{eq:Polyak:ls}
    \gamma_{k}^{\operatorname{Polyak, ls}} = \frac{2 \mathcal{R}(X_k) - \eta^2 + \lambda_k^{\operatorname{eff}} (\|X_k\|^2 - \langle X_k , X^\star\rangle )}{\frac{\Tr(K)}{d}(2\mathcal{R}(X_k) +\lambda_k^{\operatorname{eff}} \|X_k\|^2)}.
\end{equation}
In the noiseless setting, this simplifies to 
    \begin{equation}\label{eq:Polyak:ls:adv}
        \begin{aligned}
        \gamma_{k}^{\operatorname{Polyak, ls}}
         &
         = 
         \frac{d}{\Tr(K)}\left(1 - \lambda_k^{\operatorname{eff}}  \frac{ \langle X_k , X^\star\rangle }{2\mathcal{R}(X_k) + \lambda_k^{\operatorname{eff}} \|X_k\|^2}\right)
         \\
         &
         = 
         \frac{d}{\Tr(K)}\left(1 -\left(\delta^2 + \delta \sqrt{\frac{2}{\pi}}\frac{\sqrt{2\mathcal{R}(X_k)}}{\|X_k\|}\right)  \frac{ \langle X_k , X^\star\rangle }{2\mathcal{R}^{\operatorname{Adv}}(X_k)}\right).
        \end{aligned}
    \end{equation}
In the case of $\ell_2$-adversarial least squares, the descent to adversarial optimality solves the following ODE:
\begin{equation}\label{eq:ls:descent:optim:ode:adv}
    \begin{aligned}
        \frac{\dif}{\dif t}\mathrsfs{D}^{2,\Adv}(t)
        &=
        -2\gamma^{\operatorname{eff}}(t)\bigg(2\mathrsfs{R}(t) - \eta^2 + \lambda^{\operatorname{eff}}(t)(\widehat{\mathrsfs{B}}_{44}(t) - \mathrsfs{A}^\Adv(t))
        \\
        &+ \sum_{j=1}^d \lambda_j \left(\mathrsfs{V}_{j,21}(t) - \mathrsfs{V}_{j,22}(t) - \mathrsfs{U}_j(t) + \langle X^\star, \omega_j\rangle \langle X^{\star,\Adv}, \omega_j\rangle \right) \bigg)
        \\
        & +\frac{2\gamma^2(t)}{d}\Tr(K)\mathrsfs{R}^{\operatorname{Adv}}(t).
    \end{aligned}
\end{equation} 
Whenever $X^{\star, \Adv} = X^\star$, we recover \eqref{eq:ls:descent:optim:ode}. Using a similar argument as \eqref{eq:Polyak:ls}, the Polyak stepsize $\gamma_{k}^{\operatorname{Polyak, adv}}$ with respect to $X^{\star, \Adv}$ is naturally adaptive and given by
\begin{equation}\label{eq:Polyak:adv}
    \begin{aligned}
    \gamma_{k}^{\operatorname{Polyak, adv}}
     &= 
      \frac{2\mathcal{R}(X_k) - \eta^2 - (X_k- X^\star)^\top K(X^{\star, \Adv}- X^\star)+ \tilde{\lambda}_k^{\operatorname{eff}}(\|X_k\|^2  -\langle X_k , X^{\star, \Adv}\rangle) }{\frac{\Tr(K) }{d}\left(2\mathcal{R}(X_k) + \tilde{\lambda}_k^{\operatorname{eff}} \|X_k\|^2\right)}.
    \end{aligned}
\end{equation}

    Adversarial training is designed to improve robust performance by minimizing the worst-case loss over admissible perturbations, but this can come at the expense of standard risk performance. This standard risk versus robust risk trade-off has been extensively studied empirically and theoretically. See Section~\ref{sect:related:work} for more details.
    
    This brings us to the question: can this trade-off be controlled through the learning rate scheduling of adversarial training itself. In the context of our framework, this can be phrased as: given values of $\delta$, $X^\star$ and $K$, does there exist a learning rate schedule which minimizes the distance to optimality towards $X^{\star}$? The thresholds \eqref{eq:ls_gr:stab:sgd} and \eqref{eq:adv_gr:stab:sgd}
are pointwise stability thresholds along the SGD trajectory. In particular,
descent towards $X^\star$ at iteration $k$ is possible precisely when
\begin{equation}\label{eq:descentrequirement:ls}
    \langle X_k-X^\star,\nabla\mathcal R^\Adv(X_k)\rangle>0.
\end{equation}
This condition may hold along a particular path even when
$X^\star$ is not the minimizer of $\mathcal R^\Adv$ which implies it is possible
for adversarial training to move towards $X^\star$ over a region of state
space despite the adversarial optimum satisfying
$X^{\star,\Adv}\neq X^\star$. In order for \eqref{eq:descentrequirement:ls} to hold globally on $\mathcal{U}$ (see Assumption~\ref{assumption:risk}), this requires $\nabla \mathcal{R}^\Adv(X^{\star})=0$ as shown in the following Lemma.
\begin{lemma}\label{lemma:sufficient:cond:X^star}
    Suppose $\mathcal{R}^\Adv$ is continuously differentiable at $X^\star$ then
    \[
        \langle X_k - X^{\star}, \nabla \mathcal{R}^\Adv(X_k) \rangle >0 \; \forall \; X\not = X^{\star}\implies \nabla \mathcal{R}^\Adv(X^{\star})=0.
    \]
\end{lemma}
\begin{proof}
    We proceed by proof by contradiction. Suppose $\nabla \mathcal{R}^\Adv(X^{\star})\not= 0$ and for some $\varepsilon >0$, define $X_{\varepsilon} = X^\star + \varepsilon u$ for $u = -\frac{\nabla \mathcal{R}^\Adv(X^{\star})}{\|\nabla \mathcal{R}^\Adv(X^{\star})\|}$. Then we have $ \|X_{\varepsilon} - X^\star\| = \varepsilon$.

    Since $\mathcal{R}^\Adv$ is continuously differentiable at $X^\star$ then $ u^\top \nabla \mathcal{R}^\Adv(X_\varepsilon) \to u^\top \nabla \mathcal{R}^\Adv(X^\star)$ as $\varepsilon \to 0$. Now, $u^\top \nabla \mathcal{R}^\Adv(X^\star) = - \|\nabla \mathcal{R}^\Adv(X^\star)\| <0$. By continuity, for sufficiently small $\varepsilon$ we have
    \[
        \|\mathcal{R}^\Adv(X_\varepsilon) - \mathcal{R}^\Adv(X^\star)\| < \frac{1}{2}\|\mathcal{R}^\Adv(X_\star)\|.
    \]
    Hence, we obtain
    \[
    \begin{aligned}
        u^\top \nabla \mathcal{R}^\Adv(X_\varepsilon) &
        = u^\top \mathcal{R}^\Adv(X^\star) + u^\top (\nabla \mathcal{R}^\Adv(X_\varepsilon) - \mathcal{R}^\Adv(X^ \star))
        \\
        &
        \leq -\frac{1}{2} \|\nabla \mathcal{R}^\Adv(X^\star)\| < 0.
    \end{aligned}
    \]
    Thus, it follows that $(X_{\varepsilon} - X^\star)^\top\nabla \mathcal{R}^\Adv(X_\varepsilon) = \varepsilon u^\top \mathcal{R}^\Adv(X_\varepsilon)  < 0 $ which is a contradiction.
\end{proof}
Sufficient conditions for $\langle X_k - X^{\star, \Adv}, \nabla \mathcal{R}^\Adv(X_k) \rangle >0$ for all $X\not = X^{\star, \Adv}$ are either $ \mathcal{R}^\Adv$ is convex and has a unique minimizer or $ \mathcal{R}^\Adv$ is $\upsilon$-strongly convex for some $\upsilon > 0$. 

A consequence of Lemma~\ref{lemma:sufficient:cond:X^star} is that if $ X^{\star, \Adv}$ is the unique minimizer fo $\mathcal{R}^\Adv$ then $\langle X_k - X^{\star}, \nabla \mathcal{R}^\Adv(X_k) \rangle >0$ holding globally requires $ X^{\star, \Adv} = X^\star$. As proven in \citep{xing2021same_minimizer}, in the case of $\ell_2$-adversarial least squares with $K \succ 0$, $X^{\star,\Adv} = X^\star$ whenever
\begin{equation}
    \delta < \sqrt{\frac{2}{\pi}}\frac{\|X^\star\|}{(X^\star)^\top K^{-1} X^\star}.
\end{equation}
 By Lemma~\ref{lem:strong:con:Radv}, for least squares $ \mathcal{R}^{\operatorname{Adv}}$ is $\upsilon$-strongly convex with unique minimizer $X^{\star, \Adv}$ which entails that for any iterate $X_k\not = X^{\star, \Adv}$ of SGD, there exists a learning rate for which descent towards $X^{\star,\Adv}$ is possible.

\section{Line Search}\label{sect:linesearch}
\subsection{General Line Search}
Similarly to the Polyak stepsize \ref{sect:Polyakstepsize}, we determine the learning rate threshold for which the adversarial risk $\mathcal{R}^\Adv(X)$ is decreasing at each iteration when $d\to \infty$ using its deterministic equivalent. We will work with the regularization parameter $\lambda = 0$ for simplicity. Given that the expressions in this section are very involved for mixture of Gaussians, we work with a single covariance $K$, a single mean $\mu = 0$ and regularization parameter $\lambda =0$. We thus drop the subscripts $i$ on $\widehat{B}_i(W)$ for example. Similarly to Equation (90) in \citep{collinswoodfin2024highline}, from Lemma~\ref{lem:derivative_adv_risk} and Assumptions~\ref{assumption:risk}~and~\ref{assumption:fisher}, there exists functions $m$ and $v$ such that
\begin{equation}
    \begin{gathered}
    \|\nabla \mathcal{R}^{\Adv}(X)\|^2 =   m\left( \widehat B(W), \hat{X}^\top K^2 \hat{X}  \right),
    \\
    \text{and}\quad \Tr\left(\nabla^2 \mathcal{R}^{\Adv}(X) K \right)  = v
    \left(
        \widehat B(W), \hat{X}^\top K^3 \hat{X}, \hat{X}^\top K^2\hat{X}, K
    \right).
    \end{gathered}
\end{equation}
From Theorem~\ref{thm:concentration_statistic:adv} and since $m$ and $v$ satisfy Assumption~\ref{assumption:smooth:stats} then define $\mathrsfs{M}(t)$ and $\upsilon(t)$ as their respective deterministic equivalents. It follows using a similar argument as in the proof of Proposition~\ref{prop:nonexplosiveness} that $\mathrsfs{R}^\Adv$ satisfies the following ODE
 \begin{equation}\label{eq:ode:Radv}
    \dif \mathrsfs{R}^{\Adv}(t) = -\gamma(t) \mathrsfs{M}(t)\dif t + \frac{\gamma^2(t)}{d} \upsilon(t)\mathcal{I}(\widehat{\mathrsfs{B}}(t)) \dif t.
 \end{equation}
Similarly to the section on the Polyak stepsize, we get a maximal learning rate threshold which guarantees monotone decrease of $\mathrsfs{R}^{\Adv}(t)$ for a given time $t$ (i.e. $ \dif \mathrsfs{R}^{\Adv}(t) < 0$). Translating relevant deterministic equivalents to SGD, we obtain 
\begin{equation}
    \begin{gathered}
        \gamma_k^{\operatorname{Stable}, \mathrsfs{R}^\Adv}= \frac{ \|\nabla \mathcal{R}^{\Adv}(X_k)\|^2}{\frac{\Tr(\nabla^2 \mathcal{R}^{\Adv}(X_k) K)}{d}\mathcal{I}(\widehat{B}(W_k))}\quad
        \text{and det. equiv.} \quad \gamma^{\operatorname{Stable}, \mathrsfs{R}^\Adv}(t) = \frac{ \mathrsfs{M}(t)}{\frac{1}{d}\upsilon(t)\mathcal{I}(\widehat{\mathrsfs{B}}(t))}.
    \end{gathered}
\end{equation}
Analogously to \citep{collinswoodfin2024highline}, optimizing over $\gamma$ in the ODE \eqref{eq:ode:Radv} and translating relevant deterministic equivalents back to SGD, we obtain the learning rate schedule of $\emph{exact line search}$
\begin{equation}
    \gamma_k^{\operatorname{line}} = \frac{ \gamma_k^{\operatorname{Stable}, \mathrsfs{R}^\Adv}}{2} \quad
    \text{and deterministic equivalent}\quad \gamma^{\operatorname{line}}(t) = \frac{ \gamma^{\operatorname{Stable}, \mathrsfs{R}^\Adv}}{2}.
\end{equation}

\subsection{Line search for $\ell_2$-Adversarial Least Squares}
We study the learning rate schedule produced by exact line search on $\ell_2$-adversarial least squares. See \eqref{eq:loss:least squares} in Section~\ref{sect:Polyakstepsize} for the problem setup. In this case, the deterministic equivalent $\mathrsfs{R}^{\Adv}$ satisfies the following ODE
\begin{equation}
    \begin{aligned}
        \frac{\dif }{\dif t}\mathrsfs R^{\mathrm{Adv}}(t)
&=
-\gamma
\left(1+\delta\sqrt{\frac{2}{\pi}}q(t)\right)^2
\sum_{j=1}^d
\Bigl(
\bigl(\lambda_j+\tilde{\lambda}^{\mathrm{eff}}(t)\bigr)^2
\mathrsfs V_{j,11}(t)
\\
&
-2\lambda_j
\bigl(\lambda_j+\tilde{\lambda}^{\mathrm{eff}}(t)\bigr)
\mathrsfs V_{j,12}(t)
+
\lambda_j^2
\mathrsfs V_{j,22}(t)
\Bigr)
\\
&
+
\frac{\gamma^2}{d}
\left(1+\delta\sqrt{\frac{2}{\pi}}\,q(t)\right)
\mathrsfs R^{\mathrm{Adv}}(t)
\left(
\operatorname{Tr}(K^2)
+
\tilde{\lambda}^{\mathrm{eff}}(t)\operatorname{Tr}(K)
\right).
\end{aligned}
\end{equation}
Hence, for $\ell_2$-adversarial least squares we obtain
\begin{equation}\label{eq:line search:ls}
    \begin{aligned}
    &\gamma^{\operatorname{line}}(t)
    = \frac{
    \sum_{j=1}^d
    \Bigl(
    \bigl(\lambda_j+\tilde{\lambda}^{\mathrm{eff}}(t)\bigr)^2
    \mathrsfs V_{j,11}(t)
    -2\lambda_j
    \bigl(\lambda_j+\tilde{\lambda}^{\mathrm{eff}}(t)\bigr)
    \mathrsfs V_{j,12}(t)
    +
    \lambda_j^2
    \mathrsfs V_{j,22}(t)
    \Bigr) 
    }{
        \frac{2}{d}
        \left(
            \mathrsfs R(t)
            +
            \frac{\tilde{\lambda}^{\operatorname{eff}}(t)}{2}
            \widehat{\mathrsfs B}_{44}(t)
            \right)
    \left(
    \operatorname{Tr}(K^2)
    +
    \tilde{\lambda}^{\mathrm{eff}}(t)\operatorname{Tr}(K)
    \right)}.
    \end{aligned}
\end{equation}
We now provide the proof of Proposition~\ref{prop:linesearch:explicit}.
\begin{proof}[Proof of Proposition~\ref{prop:linesearch:explicit}]
For this given covariance setup, we have
\[
\frac1d
\left(
\operatorname{Tr}(K^2)
+
\tilde{\lambda}^{\operatorname{eff}}(t)\operatorname{Tr}(K)
\right)
=
\sum_{j=1}^n \frac{\#\lambda_j}{d}\lambda_j (\lambda_j + \tilde{\lambda}^{\operatorname{eff}}(t)).
\]
Given the relation of $\ell_2$-adversarial least squares to ridge regression presented in Section~\ref{sect:AdvHSGD:vs:HSGD}, we define the ridge objective function $ \Phi_{\tilde{\lambda}^{\operatorname{eff}}(t)}(t)
=
\mathrsfs R(t)
+
\frac{\tilde{\lambda}^{\operatorname{eff}}(t)}{2}
\widehat{\mathrsfs B}_{44}(t)$.
For a fixed $\tilde{\lambda}^{\operatorname{eff}}(t)$ and $J_+ \defas \{j : \lambda_j >0\}$, the minimum value this objective is 
\begin{equation}\label{eq:Phimin:ridge:defas}
\Phi_{\min}(\tilde{\lambda}^{\operatorname{eff}}(t))
=
\frac12
\sum_{j \in J_+} \lambda_j
\frac{\tilde{\lambda}^{\operatorname{eff}}(t)}
{\lambda_j+\tilde{\lambda}^{\operatorname{eff}}(t)}
\mathrsfs V_{j,22}(t) + \frac{\eta^2}{2}
\end{equation}
Define the excess risk $\mathcal{E}(t)
\defas
\Phi_{\tilde{\lambda}^{\operatorname{eff}}(t)}(t)
-
\Phi_{\min}(\tilde{\lambda}^{\operatorname{eff}}(t))$ then we obtain
\[
\mathcal{E}(t)
=
\frac{1}{2}\sum_{j \in J_+} \left(
\bigl(\lambda_j+\tilde{\lambda}^{\operatorname{eff}}(t)\bigr)\mathrsfs V_{j,11}(t)
-
2\lambda_j\mathrsfs V_{j,12}(t)
+
\frac{\lambda_j^2}{\lambda_j+\tilde{\lambda}^{\operatorname{eff}}(t)}
\mathrsfs V_{j,22}(t)
\right).
\]
It is thus easy to see that $\mathrsfs{E}(t) = \sum_{j \in J_+ } \mathrsfs{E}_j(t) $. Denote the numerator of \eqref{eq:line search:ls} as $\mathrsfs H(t)$ then it satisfies
\[
\mathrsfs H(t)
=
2\sum_{j=1}^n \left(\lambda_j+\tilde{\lambda}^{\operatorname{eff}}(t)\right)\mathcal{E}_j(t).
\]
It thus follows that
\[
\gamma^{\operatorname{line}}(t)
=
\frac{
    \sum_{j=1}^n \left(\lambda_j+\tilde{\lambda}^{\operatorname{eff}}(t)\right)\mathcal{E}_j(t)
}{
\left(
\mathrsfs R(t)
+
\frac{\tilde{\lambda}^{\operatorname{eff}}(t)}{2}
\widehat{\mathrsfs B}_{44}(t)
\right)
\sum_{j=1}^n \frac{\#\lambda_j}{d}\lambda_j (\lambda_j + \tilde{\lambda}^{\operatorname{eff}}(t))
}.
\]
Whenever $\mathcal{E}(t)>0$, define
\[
\alpha_j(t)
\defas
\frac{\mathcal{E}_j(t)}
{\mathcal{E}(t)}\quad \text{then}\quad 0\le \alpha(t)\le 1,
\]
and
\[
    \sum_{j=1}^n \left(\lambda_j+\tilde{\lambda}^{\operatorname{eff}}(t)\right)\mathcal{E}_j(t) =
    \mathcal{E}(t) \sum_{j=1}^n \alpha_j \left(\lambda_j+\tilde{\lambda}^{\operatorname{eff}}(t)\right),
\]
from which it follows that
\[
\gamma^{\operatorname{line}}(t)
=
\frac{
    \sum_{j=1}^n \alpha_j \left(\lambda_j+\tilde{\lambda}^{\operatorname{eff}}(t)\right)
}{
    \sum_{j=1}^n \frac{\#\lambda_j}{d}\lambda_j (\lambda_j + \tilde{\lambda}^{\operatorname{eff}}(t))
}
\frac{
\mathcal{E}(t)
}{
\mathrsfs R(t)+\frac{\tilde{\lambda}^{\operatorname{eff}}(t)}{2}\widehat{\mathrsfs B}_{44}(t)
}.
\]
The result follows from the definition of $\Phi_{\tilde{\lambda}^{\operatorname{eff}}(t)}(t)$ which implies
\[
\frac{
\mathcal{E}(t)
}{
\mathrsfs R(t)+\frac{\tilde{\lambda}^{\operatorname{eff}}(t)}{2}\widehat{\mathrsfs B}_{44}(t)
}
=
1-
\frac{
\Phi_{\min}(\tilde{\lambda}^{\operatorname{eff}}(t))
}{
 \Phi_{\tilde{\lambda}^{\operatorname{eff}}(t)}(t)
}.
\]
The bounds follow from the fact that $0\le \alpha(t)\le 1$ and the isotropic result follows from simply $\lambda_1 =  ... =\lambda_n$.
\end{proof}
We also provide the proof of Proposition~\ref{lem:Polyak:ridge:general}.
\begin{proof}[Proof of Proposition~\ref{lem:Polyak:ridge:general}]
    From \citep[Proposition 1]{xing2021same_minimizer} and Proposition~\ref{prop:fixed:point:gamma:to:0:main}, the minimizer of $\mathcal{R}^{\Adv}$ is given by $X^{\star, \Adv} = \sum_{j \in J_+} \frac{\lambda_j}{\lambda_j + \tilde{\lambda}_\star} \langle X^\star, \omega_j\rangle \omega_j$ for $J_+ \defas \{j : \lambda_j >0\}$ where $\tilde{\lambda}_\star$ is defined in Proposition~\ref{prop:fixed:point:gamma:to:0:main}. Hence we have for $1\leq j \leq J_+$
    \[
        \langle X^{\star,\Adv} ,\omega_j \rangle = \frac{\lambda_j}{\lambda_j +\tilde{\lambda}_{\star} }\langle X^\star,\omega_j \rangle.
    \]
    such that $X^{\star,\Adv} =0$ when $\tilde{\lambda}_\star = 0$. Through this representation, we can express the deterministic equivalent of \eqref{eq:Polyak:adv} as follows
    \[
\gamma^{\operatorname{Polyak},\Adv}(t)
=
\frac{
    \mathcal{F}^{\operatorname{Polyak}}(t)
}{
\frac{2\operatorname{Tr}(K)}{d}
\left(
\Phi_{\tilde{\lambda}^{\mathrm{eff}}(t)}(t)
\right)
},
\]
where 
\[
\begin{aligned}
\mathcal{F}^{\operatorname{Polyak}}(t)
&
\defas
\sum_{j=1}^d
\bigg(
\left(\lambda_j+\tilde{\lambda}^{\mathrm{eff}}(t)\right)
\mathrsfs V_{j,11}(t)
-
\left(
\lambda_j
+
s_j^\star
\left(\lambda_j+\tilde{\lambda}^{\mathrm{eff}}(t)\right)
\right)
\mathrsfs V_{j,12}(t)
+
\lambda_j s_j^\star
\mathrsfs V_{j,22}(t)
\bigg).
\end{aligned}
\]
From a simple computation, \eqref{eq:Phimin:ridge:defas} and \eqref{eq:convention:limit:lambdastar}, we obtain
\[
    \sum_{j=1}^d\lambda_j s_j^\star
\mathrsfs V_{j,22}(t)  =  \eta^2 +\sum_{j=1}^d \lambda_j \mathrsfs V_{j,22}(t) - 2\Phi_{\min}(\tilde{\lambda}_{\star}). 
\]
which leads to $\mathcal{F}^{\operatorname{Polyak}}(t) = S(t) - 2\Phi_{\min}(\tilde{\lambda}_{\star})$ where we define
 \[ 
 S(t) \defas \eta^2 + \sum_{j=1}^d \left( \left(\lambda_j+\widetilde{\lambda}^{\mathrm{eff}}(t)\right) \mathrsfs V_{j,11}(t) - \left(
    \lambda_j
    +
    s_j^\star
    \left(\lambda_j+\tilde{\lambda}^{\mathrm{eff}}(t)\right)
    \right)\mathrsfs V_{j,12}(t) + \lambda_j \mathrsfs V_{j,22}(t) \right).
  \] 
  Our goal is to rewrite $S(t)$ in terms of $\Phi_{\tilde{\lambda}^{\operatorname{Polyak}}(t)}(t)$ for some explicit $\tilde{\lambda}^{\operatorname{Polyak}}(t)$. It is easy to see from the definition of $\mathrsfs{R}(t)$ and the identity
  \[
    \lambda_j
    -
    s_j^\star
    \left(\lambda_j+\tilde{\lambda}^{\mathrm{eff}}(t)\right)
= \vartheta_{j}^\star(t),
  \] 
  that 
   \[ 
   S(t)-2\mathrsfs R(t) = \widetilde{\lambda}^{\mathrm{eff}}(t) \widehat{\mathrsfs B}_{44}(t)  +\sum_{j=1}^d \vartheta_{j}^\star(t)\mathrsfs V_{j,12}(t). 
   \] 
Assuming $\widehat{\mathrsfs B}_{44}(t) > 0$, define
 \[ 
    \tilde{\lambda}^{\operatorname{Polyak}}(t) \defas \widetilde{\lambda}^{\mathrm{eff}}(t) + \left( \tilde{\lambda}_{\star}-\widetilde{\lambda}^{\mathrm{eff}}(t) \right)   \sum_{j=1}^d s_j^\star \frac{ \mathrsfs V_{j,12}(t) } { \widehat{\mathrsfs B}_{44}(t)},
    \] 
    then it is clear that $ S(t) = 2\Phi_{ \tilde{\lambda}^{\operatorname{Polyak}}(t)}(t)$ from which the results follows.
\end{proof}

\section{Analysis of Examples} \label{sec:analysis_examples}

In this section, we derive the functions $h_i$, $\nabla h_i$, $f_i$, $f_i'$, $\EE_{a, \epsilon}[f_i'^2]$ for two examples fitting our framework: least squares with soft labels and binary logistic regression with hard labels. See Section 2 and Appendix D in \cite{collinswoodfin2023hitting} for more examples.
\subsection{Least Squares} \label{sec:lsq_analysis}
We study the dynamics of the $\ell_2$-adversarial least squares with regularization parameter $\lambda=0$ under the assumption that the data $a\sim \mathcal{N}(0, K)$. Given that this is equivalent to $p_1=1$ and $p_2=0$ with $\mu_1=0$, we drop the subscript $i$ in this section. The classical least squares setting consists of minimizing $X \in \R^{d}$ over the risk function
\[
\begin{aligned}
    \mathcal{R}(X) &= \frac{1}{2} \EE_{a, \epsilon}\left[\left( X^\top a - ((X^{\star})^\top a+ \epsilon)\right)^2\right] =\frac{1}{2}(X- X^{\star})^\top K(X- X^{\star})   + \frac{\eta^2}{2},
\end{aligned}
\]
with respect to targets $y =  X^{\star\top} a+ \epsilon$. It is easy to identify the functions $f$ as 
\begin{equation}\label{eq:identif:f:linreg}
    f( X^\top a  , \epsilon) = \frac{1}{2}\left( X^\top a - ((X^{\star})^\top a+ \epsilon)\right)^2.
\end{equation}
The adversarial formulation consists of the following risk
\begin{equation}\label{eq:Def:Radv:lr}
    \begin{aligned}
    &\mathcal{R}^{{\operatorname{adv}}}(X) 
    \defas 
    \frac{1}{2} \EE_{a, \epsilon}\left[\max_{|s|\leq 1}\left((X^\top a - ( (X^{\star})^\top a + \epsilon)) + \delta s \|X\| \right)^2\right]
    \\
    &
    =
    \frac{1}{2} \EE_{a, \epsilon}\left[\max_{s \in \{-1,1\}}\left(((X - X^{\star})^\top a - \epsilon) + \delta s \|X\| \right)^2\right]
    \\
    &
    =
    \frac{1}{2} \EE_{a,  \epsilon}\left[\left(\left|(X - X^{\star})^\top a - \epsilon\right| +\delta  \|X\|\right)^2\right]
    \\
    &
    =
    \frac{1}{2} \bigg(\EE_{a, \epsilon}\left[\left((X - X^{\star})^\top a - \epsilon\right)^2  \right] + 2\delta\|X\|\;\EE_{a, \epsilon}\left[\left|(X - X^{\star})^\top a - \epsilon\right|\right] +\frac{\delta^2}{2}  \|X\|^2\bigg)
    \\
    &
    =
    \mathcal{R}(X)+ \delta\sqrt{\frac{2}{\pi}}\|X\| \sqrt{2\mathcal{R}(X)} +\frac{\delta^2}{2}  \|X\|^2
    \\
    \end{aligned}
    \end{equation}
since $u \mapsto u^2$ is convex. The function $h: \R^{4\times 4} \to \R$ must satisfy $h(\widehat{B}(W)) = \mathcal{R}^{\operatorname{adv}}(X)$. Similarly to Section B.1 in \citep{collinswoodfin2023hitting}, we identify the following terms:
\begin{gather*}
    z_{11} = X^\top K X ,\quad z_{12} = X^\top K X^{\star}\quad z_{21} = X^{\star\top}KX, \\
     z_{22} =  X^{\star\top} K X^{\star},\quad  z_{44} = \|X\|^2
\end{gather*}
Hence, we have the identifications
\[
    \mathcal{R}(X)=  \frac{1}{2} \left( z_{11} - z_{12} - z_{21} + z_{22} \right)+ \frac{\eta^2}{2}.
\]
which leads to
\[
\begin{aligned}
    h\left(\begin{bmatrix}
        z_{11} & z_{12} & 0 & 0 \\
        z_{21} & z_{22} & 0 & 0 \\
        0 & 0 & 0 & 0\\
        0 & 0 & 0  & z_{44}
    \end{bmatrix}\right) &=  \frac{1}{2} \left( z_{11} - z_{12} - z_{21} + z_{22} \right) \\
    &\qquad+ \frac{\eta^2}{2}+ \delta\sqrt{\frac{2}{\pi}} \sqrt{z_{44}\left(\left( z_{11} - z_{12} - z_{21} + z_{22} \right) + \eta^2\right)} + \frac{\delta^2}{2} z_{44}.
\end{aligned}
\]
Consequently, it is clear that 
\begin{equation}\label{eq:R(t):detequiv}
    \mathrsfs{R}(t) = \frac{1}{2} \left(\mathrsfs{B}_{11}(t) - \mathrsfs{B}_{12}(t) - \mathrsfs{B}_{21}(t) + \mathrsfs{B}_{22}(t)\right)+ \frac{\eta^2}{2},
\end{equation}
where we denote $\mathrsfs{R}(t)$ as the deterministic equivalent of $\mathcal{R}(X)$. Similarly, from the definition of $z_{44}$, it is clear that $\widehat{\mathrsfs{B}}_{44}(t)$ is the deterministic equivalent of $\|X\|^2$. From \eqref{eq:Def:Radv:lr} and $\mathrsfs{R}^{{\operatorname{adv}}}(t) = h \circ (\widehat{\mathrsfs{B}}(t))$, we obtain
\begin{equation}\label{eq:Radv(t):detequiv}
    \mathrsfs{R}^{{\operatorname{adv}}}(t) =\mathrsfs{R}(t) + \delta\sqrt{\frac{2}{\pi}}\sqrt{2\widehat{\mathrsfs{B}}_{44}(t)\mathrsfs{R}(t)} +\frac{\delta^2}{2} \widehat{ \mathrsfs{B}}_{44}(t).
\end{equation}
Continuing, we then see that
\[
    \nabla h(B) = \left[ \begin{array}{c|c|c|c} 
       \partial_{11} h(B) &\partial_{12} h (B) & 0 & 0 
        \\
        \hline
       \partial_{21} h(B) &\partial_{22} h (B) & 0 & 0
       \\
       \hline
       0 & 0 & 0 &0
       \\
        \hline
        0&0 & 0 &\partial_{44} h(B)
        \end{array} \right],
\]
with
\[
\begin{aligned}
    \partial_{ij} h(B) &= \frac{1}{2} + \frac{\delta}{2} \sqrt{\frac{2}{\pi}} \sqrt{\frac{z_{44}}{z_{11} - z_{12} - z_{21} + z_{22}  + \eta^2}},\text{when}\quad (i,j)\in\{(1,1), (2,2)\},\\
    \partial_{ij} h(B) &= -\frac{1}{2} - \frac{\delta}{2} \sqrt{\frac{2}{\pi}} \sqrt{\frac{z_{44}}{z_{11} - z_{12} - z_{21} + z_{22}  + \eta^2}},\text{when}\quad (i,j)\in\{(1,2), (2,1)\}\\
    \text{and}\quad\partial_{44} h(B) &= \frac{\delta^2}{2} + \frac{\delta}{2} \sqrt{\frac{2}{\pi}} \sqrt{\frac{z_{11} - z_{12} - z_{21} + z_{22}  + \eta^2}{z_{44}}}.
\end{aligned}
\]
It immediately follows that
\[
\begin{gathered}
        H_{1}(\widehat{\mathrsfs{B}}(t)) = \frac{1}{2}\begin{bmatrix}
            1+ \delta \sqrt{\frac{2}{\pi}} \sqrt{\frac{\widehat{\mathrsfs{B}}_{44}(t)}{2\mathrsfs{R}(t)}} & 0 \\
            -1 -\delta \sqrt{\frac{2}{\pi}} \sqrt{\frac{\widehat{\mathrsfs{B}}_{44}(t)}{2\mathrsfs{R}(t)}} & 0
        \end{bmatrix},\quad H_{2}(\widehat{\mathrsfs{B}}(t)) = \frac{1}{2}\begin{bmatrix}
            0 \\
            0
        \end{bmatrix},
        \\
        \text{and}\quad \partial_{44} h(\widehat{\mathrsfs{B}}(t)) =\frac{\delta^2}{2} + \frac{\delta}{2} \sqrt{\frac{2}{\pi}} \sqrt{\frac{2 \mathrsfs{R}(t)}{\widehat{\mathrsfs{B}}_{44}(t)}}
    \end{gathered}.
\]
It is clear that we lose differentiability at $\mathrsfs{R}(t)= \frac{1}{2} \left(\mathrsfs{B}_{11}(t) - \mathrsfs{B}_{12} - \mathrsfs{B}_{21} + \mathrsfs{B}_{22}\right)+ \frac{\eta^2}{2} = 0$ and $\widehat{\mathrsfs{B}}_{44}(t) = 0$. . To apply the theory in this paper, we work on the set $\mathcal{U}$ defined as follows
\begin{equation}
    \mathcal{U} \defas \left\{ \widehat{B}: \widehat{B}_{11} - \widehat{B}_{12} - \widehat{B}_{21} + \widehat{B}_{22} + \eta > 0, \quad \widehat{B}_{44} >0\right\}
\end{equation}
Recall $\bar{s} = \operatorname{argmax}_{s\in  \{-1,1\}}\left( X^\top a  + \delta s \|X\|- ((X^{\star})^\top a+ \epsilon)\right)^2$. From \eqref{eq:identif:f:linreg}, considering $f$ as a map $r\mapsto \tfrac{1}{2}(r - (X^{\star\top} a + \epsilon))^2$ we see that the almost everywhere derivative of $f$ with respect to $X^\top a$ is given by 
\[
        f'(x) = x - (X^{\star\top} a - \epsilon)
\]
So it follows from our initial computation of $\mathcal{R}^{{\operatorname{adv}}}$ that
\[
\begin{aligned}
        \E_{a,\epsilon}[f'((X^\top a + \bar{s}\delta\|X\| , \epsilon))^2] 
        &
        = \E_{a,\epsilon}[(X^\top a+ \bar{s}\delta\|X\|-(X^{\star\top}a - \epsilon))^2]
        \\
        &
        = \E_{a,\epsilon}[\max_{s\in \{-1,1\}}(( (X- X^{\star})^\top a - \epsilon) + s\delta\|X\|)^2]
        \\
        &
        =2\mathcal{R}^{{\operatorname{adv}}}(X),
\end{aligned}
\]
which leads to
\[
    \mathcal{I}(\widehat{\mathrsfs{B}}(t)) =2\mathrsfs{R}(t) + 2\delta\sqrt{\frac{2}{\pi}}\sqrt{2\widehat{\mathrsfs{B}}_{44}(t)\mathrsfs{R}(t)} +\delta^2  \widehat{\mathrsfs{B}}_{44}(t) =2 \mathrsfs{R}^{{\operatorname{adv}}}(t).
\]
We now have all the components to explicitly compute $\Phi$ as defined in \eqref{eq:Phi:def:solved} and $\widehat{\mathrsfs{B}}(t)$ as defined in \eqref{eq:C_scalar}. We have
\begin{equation}\label{eq:Phi:def:solved:lr}
    \begin{gathered}
         \Phi_{11}(t,K) 
         = \exp\bigg(\int_0^t \gamma(s) \bigg(\left(1+ \delta \sqrt{\frac{2}{\pi}} \sqrt{\frac{\widehat{\mathrsfs{B}}_{44}(s)}{2\mathrsfs{R}(s)}}\right)K 
         \\
         \qquad+ \left( \delta^2 + \delta \sqrt{\frac{2}{\pi}} \sqrt{\frac{2\mathrsfs{R}(s)}{\widehat{\mathrsfs{B}}_{44}(s)}}\right)\Id_d\bigg) \dif s\bigg), \, \, 
         \\
         \text{and} \, \,
        \Phi_{21}(t, K) 
        =-\int_0^t  \gamma(s)  \left(1 + \delta \sqrt{\frac{2}{\pi}} \sqrt{\frac{\widehat{\mathrsfs{B}}_{44}(s)}{2\mathrsfs{R}(s)}}\right)K\Phi_{11}(s,K) \dif s. 
        \end{gathered}
    \end{equation}
We define the following adaptive learning rate and regularization parameter respectively
\begin{equation}\label{eq:adapt:gamma:lambda}
    \tilde{\gamma}^{\operatorname{eff}}(t) = \gamma(t) \left(1 + \delta \sqrt{\frac{2}{\pi}} \sqrt{\frac{\widehat{\mathrsfs{B}}_{44}(t)}{2\mathrsfs{R}(t)}}\right)
    \quad \text{and}\quad 
    \tilde{\lambda}^{\operatorname{eff}}(t) = \frac{ \delta^2 + \delta \sqrt{\frac{2}{\pi}} \sqrt{\frac{2\mathrsfs{R}(t)}{\widehat{\mathrsfs{B}}_{44}(t)}}}{1 + \delta \sqrt{\frac{2}{\pi}} \sqrt{\frac{\widehat{\mathrsfs{B}}_{44}(t)}{2\mathrsfs{R}(t)}}}.
\end{equation}
Define the following integral terms
\begin{equation}\label{eq:adapt:gamma:lambda:int}
    \tilde{\Gamma}(t) = \int_0^t \tilde{\gamma}^{\operatorname{eff}}(s) \dif s\quad \text{and}\quad  \Lambda(t) = \int_0^t \tilde{\gamma}^{\operatorname{eff}}(s)\tilde{\lambda}^{\operatorname{eff}}(s) \dif s.
\end{equation}
Then plugging these terms into \eqref{eq:Phi:def:solved:lr}, we obtain
\begin{equation}
    \begin{aligned}
         \Phi_{11}(t,K) 
         &
          = e^{K \tilde{\Gamma}(t) + \Lambda(t)}, \quad
         \text{and} \quad
        \Phi_{21}(t, K) 
        &
        =-\int_0^t  K\tilde{\gamma}^{\operatorname{eff}}(s) e^{K\tilde{\Gamma}(s) + \Lambda(s)}\dif s. 
        \end{aligned}
    \end{equation}
We can also view $\Phi_{21}(t, K) $ the following way
\[
\begin{aligned}
    \Phi_{21}(t, K) 
    &
    =1 - \Phi_{11}(t,K) +\int_0^t \tilde{\gamma}^{\operatorname{eff}}(s)\tilde{\lambda}^{\operatorname{eff}}(s) \Phi_{11}(s,K)\dif s.
\end{aligned}
\]
 From \eqref{eq:C_scalar}, the first three terms of $\mathrsfs{B}_{11}(t)$ and $\widehat{\mathrsfs{B}}_{44}(t)$ form a perfect square from which it follows
\begin{equation}\label{eq:def:l2norm2:det:equiv:lr:withPhi}
    \begin{aligned}
        \mathrsfs{B}_{11}(t) 
        &= 
        \langle K , \left(\Phi_{11}^{-1}(t,K)\left(X_0 - \Phi_{21}(t,K) X^\star \right)\right)^{\otimes 2}\rangle
        \\
        &
        \qquad + \frac{2}{d} \int_0^t \gamma^2(s)  \mathrsfs{R}^{{\operatorname{adv}}}(s) \tr \big ( K^2 e^{-2K(\tilde{\Gamma}(t) - \tilde{\Gamma}(s)) -2(\Lambda(t) -\Lambda(s))} \big ) \, \dif s,\\
        \widehat{\mathrsfs{B}}_{44}(t) 
        &= \|\Phi_{11}^{-1}(t,K) (X_0 - \Phi_{21}(t,K) X^\star)\|^2
        \\
        &
        \qquad+ \frac{2}{d} \int_0^t \gamma^2(s)  \mathrsfs{R}^{{\operatorname{adv}}}(s) \tr \big ( K e^{-2K(\tilde{\Gamma}(t) - \tilde{\Gamma}(s)) -2(\Lambda(t) -\Lambda(s))} \big ) \, \dif s,
    \end{aligned}
\end{equation}
Since 
\[
\begin{aligned}
    \mathrsfs{B}_{12}(t)  = \mathrsfs{B}_{21} (t) &= 
    X_0^\top \tfrac{K}{\Phi_{11}(t,K)} X^{\star} - (X^{\star})^\top \tfrac{K\Phi_{21}(t,K)}{\Phi_{11}(t,K)} X^{\star} 
    \\
    &
    = \left(X_0 - \Phi_{21}(t,K)X^{\star}\right)^\top\Phi_{11}^{-1}(t,K)K X^{\star}
\end{aligned}
\]
From this, it follows that
\begin{equation}
    \begin{aligned}
        &\mathrsfs{B}_{11}(t) - \mathrsfs{B}_{12}(t) - \mathrsfs{B}_{21}(t) + \mathrsfs{B}_{22}(t)\\
        &
        =\langle K , (e^{-K\tilde{\Gamma}(t)  - \Lambda(t)}\left(X_0 + \left(\int_0^t  K\tilde{\gamma}^{\operatorname{eff}}(s) e^{K \tilde{\Gamma}(s) + \Lambda(s)}\dif s\right) X^\star\right)- X^{\star})^{\otimes 2}\rangle
        \\
        &
        \qquad + \frac{2}{d} \int_0^t \gamma^2(s)  \mathrsfs{R}^{{\operatorname{adv}}}(s) \tr \big ( K^2 e^{-2K(\tilde{\Gamma}(t) - \tilde{\Gamma}(s)) -2(\Lambda(t) -\Lambda(s))} \big ) \, \dif s.
    \end{aligned}
\end{equation}
Hence, we obtain
\begin{equation}\label{eq:R(t):equiv:final:lr}
\begin{aligned}
    \mathrsfs{R}(t)&= \frac{1}{2}\langle K , (e^{-K\tilde{\Gamma}(t)  - \Lambda(t)}\left(X_0 + \left(\int_0^t  K\tilde{\gamma}^{\operatorname{eff}}(s) e^{K \tilde{\Gamma}(s) + \Lambda(s)}\dif s\right) X^\star \right)- X^{\star})^{\otimes 2}\rangle + \frac{\eta^2}{2}
    \\
    &
    \qquad + \frac{1}{d} \int_0^t \gamma^2(s)  \mathrsfs{R}^{{\operatorname{adv}}}(s) \tr \big ( K^2 e^{-2K(\tilde{\Gamma}(t) - \tilde{\Gamma}(s)) -2(\Lambda(t) -\Lambda(s))} \big ) \, \dif s.
\end{aligned}
\end{equation}
Plugging the values of $\Phi_{11}$ and $\Phi_{21}$ in \eqref{eq:def:l2norm2:det:equiv:lr:withPhi}, we obtain
\begin{equation}\label{eq:def:l2norm2:det:equiv:lr}
\begin{aligned}
    \widehat{\mathrsfs{B}}_{44}(t) 
    &= \left\| e^{-K\tilde{\Gamma}(t)  -\Lambda(t)}\left(X_0 +\left(\int_0^t  K\tilde{\gamma}^{\operatorname{eff}}(s) e^{K \tilde{\Gamma}(s) + \Lambda(s)}\dif s\right) X^\star\right)\right\|^2
    \\
    &
    \qquad+ \frac{2}{d} \int_0^t \gamma^2(s)  \mathrsfs{R}^{{\operatorname{adv}}}(s) \tr \big ( K e^{-2K(\tilde{\Gamma}(t) - \tilde{\Gamma}(s)) -2(\Lambda(t) -\Lambda(s))} \big ) \, \dif s.
\end{aligned}
\end{equation}
Finally, combining \eqref{eq:R(t):equiv:final:lr} and \eqref{eq:def:l2norm2:det:equiv:lr}, we recover $ \mathrsfs{R}^{{\operatorname{adv}}}(t)$ using \eqref{eq:Radv(t):detequiv}. Note that we may express $\mathrsfs{R}^{{\operatorname{adv}}}(t)$ in terms of $\tilde{\gamma}^{\operatorname{eff}}(t)$ and $\tilde{\lambda}^{\operatorname{eff}}(t)$ as follows
\begin{equation}\label{eq:Radv:adapt}
    \mathrsfs{R}^{{\operatorname{adv}}}(t) = \frac{\tilde{\gamma}^{\operatorname{eff}}(t)}{\gamma(t)} \left(\mathrsfs{R}(t) + \frac{\tilde{\lambda}^{\operatorname{eff}}(t)}{2}\widehat{\mathrsfs{B}}_{44}(t)\right).
\end{equation}
Using this reformulation, we obtain
\begin{equation}\label{eq:def:R(t):det:equiv:lr:reform}
    \begin{aligned}
        \mathrsfs{R}(t)&= \frac{1}{2}\langle K , (e^{-K\tilde{\Gamma}(t)  - \Lambda(t)}\left(X_0 + \left(\int_0^t  K\tilde{\gamma}^{\operatorname{eff}}(s) e^{K \tilde{\Gamma}(s) + \Lambda(s)}\dif s\right) X^\star \right)- X^{\star})^{\otimes 2}\rangle + \frac{\eta^2}{2}
        \\
        &
        \qquad + \frac{1}{d} \int_0^t \gamma(s) \tilde{\gamma}^{\operatorname{eff}}(s)  \left(\mathrsfs{R}(s) + \frac{\tilde{\lambda}^{\operatorname{eff}}(s)}{2}\widehat{\mathrsfs{B}}_{44}(s)\right) \tr \big ( K^2 e^{-2K(\tilde{\Gamma}(t) - \tilde{\Gamma}(s)) -2(\Lambda(t) -\Lambda(s))} \big ) \, \dif s.
    \end{aligned}
\end{equation}
and 
\begin{equation}\label{eq:def:l2norm2:det:equiv:lr:reform}
    \begin{aligned}
        \widehat{\mathrsfs{B}}_{44}(t) 
        &= \left\| e^{-K\tilde{\Gamma}(t)  -\Lambda(t)}\left(X_0 +\left(\int_0^t  K\tilde{\gamma}^{\operatorname{eff}}(s) e^{K \tilde{\Gamma}(s) + \Lambda(s)}\dif s\right) X^\star\right)\right\|^2
        \\
        &
        \qquad+ \frac{2}{d} \int_0^t \gamma(s) \tilde{\gamma}^{\operatorname{eff}}(s)  \left(\mathrsfs{R}(s) + \frac{\tilde{\lambda}^{\operatorname{eff}}(s)}{2}\widehat{\mathrsfs{B}}_{44}(s)\right) \tr \big ( K e^{-2K(\tilde{\Gamma}(t) - \tilde{\Gamma}(s)) -2(\Lambda(t) -\Lambda(s))} \big ) \, \dif s.
    \end{aligned}
    \end{equation}
     The system of ODEs \eqref{eq:ODE:V_i(t):def} in this setting simplifies to
\begin{equation}\label{eq:odes:linrreg}
   \begin{aligned}
       \frac{\dif \mathrsfs{V}_{j,11}(t)}{\dif t} 
       &
       = 
       -2 \gamma(t)\bigg(\lambda_j \bigg(1 + \delta\sqrt{\frac{2}{\pi}} \sqrt{\frac{\widehat{\mathrsfs{B}}_{44}(t)}{2\mathrsfs{R}(t)}}\bigg) + \delta^2 + \delta\sqrt{\frac{2}{\pi}} \sqrt{\frac{2\mathrsfs{R}(t)}{\widehat{\mathrsfs{B}}_{44}(t)}}\bigg) \mathrsfs{V}_{j,11}(t) 
       \\
       &
       + 2\gamma(t) \lambda_j \bigg(1+ \delta\sqrt{\frac{2}{\pi}} \sqrt{\frac{\widehat{\mathrsfs{B}}_{44}(t)}{2\mathrsfs{R}(t)}}\bigg)\mathrsfs{V}_{j,21}(t) + 2\lambda_j\frac{\gamma^2(t)}{d} \mathrsfs{R}^{\operatorname{adv}}(t)
        \\
        \frac{\dif \mathrsfs{V}_{j,21}(t)}{\dif t} 
        &
        = 
        - \gamma(t)\bigg(\lambda_j \bigg(1 + \delta\sqrt{\frac{2}{\pi}} \sqrt{\frac{\widehat{\mathrsfs{B}}_{44}(t)}{2\mathrsfs{R}(t)}}\bigg) + \delta^2 + \delta\sqrt{\frac{2}{\pi}} \sqrt{\frac{2\mathrsfs{R}(t)}{\widehat{\mathrsfs{B}}_{44}(t)}} \bigg)\mathrsfs{V}_{j,21}(t) 
        \\
        &
        +\gamma(t) \lambda_j \bigg(1+ \delta\sqrt{\frac{2}{\pi}} \sqrt{\frac{\widehat{\mathrsfs{B}}_{44}(t)}{2\mathrsfs{R}(t)}}\bigg) \mathrsfs{V}_{j,22}(t),
   \end{aligned}
\end{equation}
where we recall $\mathrsfs{V}_{j,22}(t) = X^{\star\top} \omega_j\omega_j^\top X^\star$. Hence, from \eqref{eq:adapt:gamma:lambda} and \eqref{eq:Radv:adapt}, we obtain the simplified representation
\begin{equation}
    \begin{aligned}
        \frac{\dif \mathrsfs{V}_{j,11}(t)}{\dif t} 
        &
        = 
        -2 \tilde{\gamma}^{\operatorname{eff}}(t)\left(\lambda_j  + \tilde{\lambda}^{\operatorname{eff}}(t)\right) \mathrsfs{V}_{j,11}(t) 
        \\
        &
        + 2\tilde{\gamma}^{\operatorname{eff}}(t) \lambda_j \mathrsfs{V}_{j,21}(t) + 2\lambda_j\frac{\gamma(t) \tilde{\gamma}^{\operatorname{eff}}(t)}{d}    \left(\mathrsfs{R}(t) + \frac{\tilde{\lambda}^{\operatorname{eff}}(t)}{2}\widehat{\mathrsfs{B}}_{44}(t)\right) 
         \\
         \frac{\dif \mathrsfs{V}_{j,21}(t)}{\dif t} 
         &
         = 
         - \tilde{\gamma}^{\operatorname{eff}}(t)\left(\lambda_j  + \tilde{\lambda}^{\operatorname{eff}}(t)\right)\mathrsfs{V}_{j,21}(t) 
         +\tilde{\gamma}^{\operatorname{eff}}(t) \lambda_j  \mathrsfs{V}_{j,22}(t).
    \end{aligned}
 \end{equation}
 \subsubsection{Proof of Proposition~\ref{prop:adv:reg:volterra:stability}}\label{sect:compareADVHSGD:to:HSGD}
Recall from \eqref{eq:main:AdvHSGD:def}, Adversarial homogenized SGD for $\ell_2$-adversarial least squares is given by
 \begin{equation}\label{eq:AdvHSGD:linreg:adv}
    \dif \WHSGD_t = -\gamma(t) \nabla \mathcal{R}^{\operatorname{adv}}(\WHSGD_t) \dif t + \gamma(t) \sqrt{\frac{2}{d} \mathcal{R}^{\operatorname{adv}} K} \dif B_t.
    \end{equation}
From a simple computation we see that
\begin{equation}
    \nabla \mathcal{R}^{\operatorname{adv}}(\WHSGD_t) = \left(1 + \delta \sqrt{\frac{2}{\pi}}\frac{\|\WHSGD_t\|}{\sqrt{2\mathcal{R}(\WHSGD_t)}}\right)K (\WHSGD_t - X^{\star}) + \left(\delta^2 + \delta \sqrt{\frac{2}{\pi}}\frac{\sqrt{2\mathcal{R}(\WHSGD_t)}}{\|\WHSGD_t\|}\right)X.
\end{equation}
Since $ \nabla \mathcal{R}(\WHSGD_t) =K (\WHSGD_t - X^{\star})$, with \eqref{eq:adapt:gamma:lambda} and \eqref{eq:Radv:adapt} in mind, we obtain the reformulation
\begin{equation}\label{eq:AdvHSGD:app}
    \dif \WHSGD_t = -\gamma^{\operatorname{eff}}(t)\left(\nabla \mathcal{R}(\WHSGD_t) + \lambda^{\operatorname{eff}}(t)\WHSGD_t\right) \dif t + \sqrt{\frac{2\gamma(t)\gamma^{\operatorname{eff}}(t)}{d} \left(\mathcal{R}(\WHSGD_t) + \frac{\lambda^{\operatorname{eff}}(t)}{2}\|\WHSGD_t\|^2\right)K} \dif B_t,
\end{equation}
where we define analogously
\begin{equation}
    \gamma^{\operatorname{eff}}(t) = \gamma(t)\bigg(1 + \delta \sqrt{\frac{2}{\pi}}\frac{\|\WHSGD_t\|}{\sqrt{2\mathcal{R}(\WHSGD_t)}}\bigg)\quad\text{and}\quad \lambda^{\operatorname{eff}}(t) = \frac{\delta^2 + \delta \sqrt{\frac{2}{\pi}}\frac{\sqrt{2\mathcal{R}(\WHSGD_t)}}{\|\WHSGD_t\|}}{1 + \delta \sqrt{\frac{2}{\pi}}\frac{\|\WHSGD_t\|}{\sqrt{2\mathcal{R}(\WHSGD_t)}}}.
\end{equation}
Now, $ \gamma^{\operatorname{eff}}(t)$ and $\lambda^{\operatorname{eff}}(t)$ correspond to the effective learning rate and regularization of $\ell_2$-adversarial least squares. As originally shown in \citep{CollinsWoodfinPaquette01}, Homogenized SGD (HSGD) is a similar SDE to AdvHSGD for $\ell_2$-regularized least squares with parameter $\lambda  >0$ and is the solution to the following SDE
\begin{equation}\label{eq:HSGD:reg}
    \dif \WHSGD_t^{\operatorname{Reg}} = -\gamma(t) \left(\nabla \mathcal{R}(\WHSGD_t^{\operatorname{Reg}}) + \lambda \WHSGD_t^{\operatorname{Reg}} \right)\dif t + \gamma(t) \sqrt{\frac{2}{d} \mathcal{R}(\WHSGD_t^{\operatorname{Reg}}) K} \dif B_t.
    \end{equation}
Now, introduce the adaptive learning rate and regularization parameters
\begin{equation}
    \gamma^{\operatorname{Reg}}(t) = \gamma(t)\bigg(1 + \delta \sqrt{\frac{2}{\pi}}\frac{\|\WHSGD_t^{\operatorname{Reg}}\|}{\sqrt{2\mathcal{R}(\WHSGD_t^{\operatorname{Reg}})}}\bigg)\quad\text{and}\quad \lambda^{\operatorname{Reg}}(t) = \frac{\delta^2 + \delta \sqrt{\frac{2}{\pi}}\frac{\sqrt{2\mathcal{R}(\WHSGD_t^{\operatorname{Reg}})}}{\|\WHSGD_t^{\operatorname{Reg}}\|}}{1 + \delta \sqrt{\frac{2}{\pi}}\frac{\|\WHSGD_t^{\operatorname{Reg}}\|}{\sqrt{2\mathcal{R}(\WHSGD_t^{\operatorname{Reg}})}}}.
\end{equation}
Then the SDE for $\ell_2$-regularized least squares amounts to
\begin{equation}\label{eq:LRHSGD:reg:adaptive}
    \begin{aligned}
    \dif \WHSGD_t^{\operatorname{Reg}} &= -\gamma^{\operatorname{Reg}}(t) \left(\nabla \mathcal{R}(\WHSGD_t^{\operatorname{Reg}}) + \lambda^{\operatorname{Reg}}(t) \WHSGD_t^{\operatorname{Reg}} \right)\dif t 
    \\
    &
    \qquad + \gamma(t) \sqrt{\frac{2}{d}\bigg( \mathcal{R}(\WHSGD_t^{\operatorname{Reg}})+ \delta \sqrt{\frac{2}{\pi}} \|\WHSGD_t^{\operatorname{Reg}}\| \sqrt{2 \mathcal{R}(\WHSGD_t^{\operatorname{Reg}})}  + \frac{\delta^2}{\pi}\|\WHSGD_t^{\operatorname{Reg}}\|^2\bigg)K} \dif B_t.
    \end{aligned}
\end{equation}
 Let $\mathscr{R}(t)$ and $\Phi(t)$ be the deterministic equivalents of $\mathcal{R}(\WHSGD_t^{\operatorname{Reg}})$ and $\|\WHSGD_t^{\operatorname{Reg}}\|^2$ respectively and define the deterministic equivalents of $\gamma^{\operatorname{Reg}}(t)$ and $\lambda^{\operatorname{Reg}}(t) $ as follows
\begin{equation}\label{eq:adapt:gamma:lambda:reg}
   \mathfrak{g}(t) = \gamma(t) \left(1 + \delta \sqrt{\frac{2}{\pi}} \sqrt{\frac{\Phi(t)}{2\mathscr{R}(t)}}\right)\quad \text{and}\quad \mathfrak{l}(t) = \frac{ \delta^2 + \delta \sqrt{\frac{2}{\pi}} \sqrt{\frac{2\mathscr{R}(t)}{\Phi(t)}}}{1 + \delta \sqrt{\frac{2}{\pi}} \sqrt{\frac{\Phi(t)}{2\mathscr{R}(t)}}}.
\end{equation}
We also introduce the integral terms
\begin{equation}\label{eq:Gamma:Lambda:reg}
    \tilde{\Gamma}^{\operatorname{Reg}}(t)
    =
    \int_0^t
    \mathfrak{g}(s)
    \dif s\quad \text{and}
    \quad
    \Lambda^{\operatorname{Reg}}(t)
    =
    \int_0^t
    \mathfrak{g}(s)
    \mathfrak{l}(s)
    \dif s.
\end{equation}
Using \eqref{eq:adapt:gamma:lambda:reg}, we see that
\begin{equation}\label{eq:reg:diffusion:forcing}
    \left(\mathfrak{g}(t)\right)^2
    \mathscr{R}(t)
    =
    \gamma^2(t)
    \bigg(
        \mathscr{R}(t)
    +
    \delta \sqrt{\frac{2}{\pi}}
    \sqrt{
    2\mathscr{R}(t)
   \Phi(t)
    }
    +
    \frac{\delta^2}{\pi}
    \Phi(t)
    \bigg).
\end{equation}
Using similar methods to \citep{collinswoodfin2024highline, CollinsWoodfinPaquette01,collinswoodfin2023hitting}, we infer the deterministic equivalents of the least squares risk $\mathcal{R}(\WHSGD_t^{\operatorname{Reg}})$ and $\|\WHSGD_t^{\operatorname{Reg}}\|^2$ are given by 
\begin{equation}\label{eq:R:reg:det:equiv:volterra}
\begin{aligned}
    \mathscr{R}(t)
&=
    \frac{1}{2}
    \bigg\langle
    K,
    \bigg(
    e^{-K\tilde{\Gamma}^{\operatorname{Reg}}(t)-\Lambda^{\operatorname{Reg}}(t)}
    \left(
    X_0
    +
    \left(
    \int_0^t
    K\mathfrak{g}(u)
    e^{K\tilde{\Gamma}^{\operatorname{Reg}}(u)+\Lambda^{\operatorname{Reg}}(u)}
    \dif u
    \right)
    X^\star
    \right)
    -
    X^\star
    \bigg)^{\otimes 2}
    \bigg\rangle
    \\
    &\qquad
    +
    \frac{\eta^2}{2}
    +
    \frac{1}{d}
    \int_0^t
    \gamma^2(s)
    \bigg(
   \mathscr{R}(s)
    +
    \delta \sqrt{\frac{2}{\pi}}
    \sqrt{
    2\mathscr{R}(s)
    \Phi(s)
    }
    +
    \frac{\delta^2}{\pi}
    \Phi(s)
    \bigg)
    \\
    &\qquad\qquad
    \times
    \tr\left(
    K^2
    e^{-2K\left(\tilde{\Gamma}^{\operatorname{Reg}}(t)
    -
    \tilde{\Gamma}^{\operatorname{Reg}}(s)\right)
    -2\left(\Lambda^{\operatorname{Reg}}(t)-\Lambda^{\operatorname{Reg}}(s)\right)}
    \right)
    \dif s.
\end{aligned}
\end{equation}
and
\begin{equation}\label{eq:l2norm2:reg:det:equiv:volterra}
\begin{aligned}
   \Phi(t)
    &=
    \left\|
    e^{-K\tilde{\Gamma}^{\operatorname{Reg}}(t)-\Lambda^{\operatorname{Reg}}(t)}
    \left(
    X_0
    +
    \left(
    \int_0^t
    K\mathfrak{g}(u)
    e^{K\tilde{\Gamma}^{\operatorname{Reg}}(u)+\Lambda^{\operatorname{Reg}}(u)}
    \dif u
    \right)
    X^\star
    \right)
    \right\|^2
    \\
    &\qquad
    +
    \frac{2}{d}
    \int_0^t
    \gamma^2(s)
    \bigg(
   \mathscr{R}(s)
    +
    \delta \sqrt{\frac{2}{\pi}}
    \sqrt{
    2\mathscr{R}(s)
    \Phi(s)
    }
    +
    \frac{\delta^2}{\pi}
    \Phi(s)
    \bigg)
    \\
    &\qquad\qquad
    \times
    \tr\left(
    K
    e^{-2K\left(\tilde{\Gamma}^{\operatorname{Reg}}(t)
    -
    \tilde{\Gamma}^{\operatorname{Reg}}(s)\right)
    -2\left(\Lambda^{\operatorname{Reg}}(t)-\Lambda^{\operatorname{Reg}}(s)\right)}
    \right)
    \dif s.
\end{aligned}
\end{equation}
We are now ready to prove Proposition~\ref{prop:adv:reg:volterra:stability}.
    \begin{proof}[Proof of Proposition~\ref{prop:adv:reg:volterra:stability}]
    For notational simplicity, we write $D(t) \defas \left|\mathrsfs{R}(t)-\mathscr{R}(t)\right|
    +
    \left|
    \widehat{\mathrsfs{B}}_{44}(t)-\Phi(t)
    \right|$ and define
    \[
        q_{\Adv}(t)
        \defas
        \sqrt{
        \frac{\widehat{\mathrsfs{B}}_{44}(t)}
        {2\mathrsfs{R}(t)}} \quad \text{and}\quad q_{\operatorname{Reg}}(t)
        \defas
        \sqrt{
        \frac{\Phi(t)}
        {2\mathscr{R}(t)}
        }
        .
    \]
    By the assumptions on $\mathrsfs{R},\mathscr{R},
    \widehat{\mathrsfs{B}}_{44}$ and $\Phi$, there exist constants
    $0<q_-<q_+<\infty$ such that
    \[
        q_- \leq q_{\Adv}(t),q_{\operatorname{Reg}}(t)\leq q_+
    \]
    for all $t\in[0,T]$. We also define for notational simplicity 
    \[
    \begin{gathered}
        A_{\Adv}(t)
        \defas
        1+
        \delta\sqrt{\frac{2}{\pi}}q_{\Adv}(t),\quad B_{\Adv}(t)
        \defas
        \delta^2+
        \delta\sqrt{\frac{2}{\pi}}\frac{1}{q_{\Adv}(t)}
        \\
        A_{\operatorname{Reg}}(t)
        \defas
        1+
        \delta\sqrt{\frac{2}{\pi}}q_{\operatorname{Reg}}(t)\quad \text{and}\quad B_{\operatorname{Reg}}(t)
        \defas
        \delta^2+
        \delta\sqrt{\frac{2}{\pi}}\frac{1}{q_{\operatorname{Reg}}(t)}.
    \end{gathered}
    \]
    We first show that $\tilde{\Gamma}$, $\tilde{\Gamma}^{\operatorname{Reg}}$, $\Lambda$ and $\Lambda^{\operatorname{Reg}}$ are Lipschitz. Since the map $ (R,N)\mapsto \sqrt{\frac{N}{2R}}$ is smooth on compact subsets of $(0,\infty)\times (0, \infty)$, the assumptions imply
    \[
        |q_{\Adv}(t)-q_{\operatorname{Reg}}(t)|
        \leq
        C
        \left(
        |\mathrsfs{R}(t)-\mathscr{R}(t)|
        +
        |\widehat{\mathrsfs{B}}_{44}(t)-\Phi(t)|
        \right)
        =
        CD(t).
    \]
    Since $q_{\Adv},q_{\operatorname{Reg}}\in[q_-,q_+]$, the map $q\mapsto q^{-1}$ is also
    Lipschitz on the relevant interval. Hence
    \begin{equation}
    \label{eq:A:B:lipschitz}
        |A_{\Adv}(t)-A_{\operatorname{Reg}}(t)|
        +
        |B_{\Adv}(t)-B_{\operatorname{Reg}}(t)|
        \leq
        CD(t).
    \end{equation}
    Therefore
    \begin{equation}
    \label{eq:Gamma:lipschitz}
        \left|
            \tilde{\Gamma}(t)
        -
        \tilde{\Gamma}^{\operatorname{Reg}}(t)
        \right|
        \leq
        C(\bar{\gamma})
        \int_0^t
        D(s)\dif s \quad \text{and}\quad \left|
            \Lambda(t)
            -
            \Lambda^{\operatorname{Reg}}(t)
            \right|
            \leq
            C(\bar{\gamma})
            \int_0^t
            D(s)\dif s.
    \end{equation}
    More generally, for $0\leq s\leq t$,
    \begin{equation}
    \label{eq:Gamma:Lambda:increment:lipschitz}
    \begin{gathered}
    \left|
    \left(
        \tilde{\Gamma}(t)-\tilde{\Gamma}(s)
    \right)
    -
    \left(
    \tilde{\Gamma}^{\operatorname{Reg}}(t)-\tilde{\Gamma}^{\operatorname{Reg}}(s)
    \right)
    \right|+
    \left|
    \left(
    \Lambda(t)-\Lambda(s)
    \right)
    -
    \left(
    \Lambda^{\operatorname{Reg}}(t)-\Lambda^{\operatorname{Reg}}(s)
    \right)
    \right|
    \\
    \leq
    C(\bar{\gamma})
    \int_s^t
   D(u)\dif u.
    \end{gathered}
    \end{equation}
    For $p=1,2$, define the kernels
    \[
        \mathcal{K}_{p}^{\Adv}(t,s)
        \defas
        \tr\left(
        K^p
        e^{-2K( \tilde{\Gamma}(t)-\tilde{\Gamma}(s))
        -2(\Lambda(t)-\Lambda(s))}
        \right),
    \]
    and
    \[
        \mathcal{K}_{p}^{\operatorname{Reg}}(t,s)
        \defas
        \tr\left(
        K^p
        e^{-2K(\tilde{\Gamma}^{\operatorname{Reg}}(t)-\tilde{\Gamma}^{\operatorname{Reg}}(s))
        -2(\Lambda^{\operatorname{Reg}}(t)-\Lambda^{\operatorname{Reg}}(s))}
        \right).
    \]    
    From the eigendecomposition of $K$, we obtain the bounds
    \[
        \mathcal{K}_{p}^{\Adv}(t,s)
        \leq
        \Tr(K^p),
        \qquad
        \mathcal{K}_{p}^{\operatorname{Reg}}(t,s)
        \leq
        \Tr(K^p),
    \]
    which imply
    \begin{equation}
    \label{eq:kernel:bounded}
        \frac{1}{d}\mathcal{K}_{p}^{\Adv}(t,s)
        +
        \frac{1}{d}\mathcal{K}_{p}^{\operatorname{Reg}}(t,s)
        \leq C(\|K\|_{\opt}).
    \end{equation}
    Moreover, since $ |e^{-x}-e^{-y}|\leq |x-y|$ for all $x,y\geq 0$ and using \eqref{eq:Gamma:Lambda:increment:lipschitz}, we obtain
    \begin{equation}
    \label{eq:kernel:lipschitz}
        \frac{1}{d}
        \left|
        \mathcal{K}_{p}^{\Adv}(t,s)
        -
        \mathcal{K}_{p}^{\operatorname{Reg}}(t,s)
        \right|
        \leq
        C(\|K\|_{\opt},\bar{\gamma})
        \int_s^t
       D(u)\dif u,
        \quad \text{for}\quad p=1,2.
    \end{equation}
    Now, let $\mathcal{R}_{\operatorname{Reg}}(t)\defas    \mathscr{R}(t)
    +
    \delta\sqrt{\frac{2}{\pi}}
    \sqrt{
    2\mathscr{R}(t)\Phi(t)
    }
    +
    \frac{\delta^2}{\pi}
    \Phi(t)$. Then we obtain adding and subtracting $\frac{\delta^2}{\pi}\widehat{\mathrsfs{B}}_{44}(t)$
    \[
    \begin{aligned}
        \mathrsfs{R}^{\Adv}(t)-\mathcal{R}_{\operatorname{Reg}}(t)
        &=
        \mathrsfs{R}(t)-\mathscr{R}(t)
        +
        \delta\sqrt{\frac{2}{\pi}}
        \left(
        \sqrt{
        2\mathrsfs{R}(t)\widehat{\mathrsfs{B}}_{44}(t)
        }
        -
        \sqrt{
        2\mathscr{R}(t)\Phi(t)
        }
        \right)
        \\
        &
+\frac{\delta^2}{\pi}
\left(
\widehat{\mathrsfs{B}}_{44}(t)-\Phi(t)
\right)
+
\delta^2
\left(
\frac12-\frac1\pi
\right)
\widehat{\mathrsfs{B}}_{44}(t)
    \end{aligned}
    \]
    Since the map $(R,N)\mapsto \sqrt{2RN}$ is smooth on compact subsets of
    $(0,\infty) \times (0, \infty)$, the boundedness assumption implies
    \[
        \left|
        \sqrt{
        2\mathrsfs{R}(t)\widehat{\mathrsfs{B}}_{44}(t)
        }
        -
        \sqrt{
        2\mathscr{R}(t)\Phi(t)
        }
        \right|
        \leq
        CD(t).
    \]
    Therefore
    \begin{equation}
    \label{eq:forcing:difference}
        \left|
            \mathrsfs{R}^{\Adv}(t)-\mathcal{R}_{\operatorname{Reg}}(t)
        \right|
        \leq
        CD(t)
        +
        \delta^2
        \left(
        \frac12-\frac1\pi
        \right)
        \widehat{\mathrsfs{B}}_{44}(t).
    \end{equation}
    For simplicity of notation, we define the following terms which appear in the first term of the deterministic equivalents:
\[
\begin{gathered}
    m_{\Adv}(t)
    \defas
    e^{-K\tilde{\Gamma}(t)-\Lambda(t)}
    \left(
    X_0+
    \left(
    \int_0^t
    K\gamma(u)A_{\Adv}(u)
    e^{K\tilde{\Gamma}_Z(u)+\Lambda_Z(u)}
    \dif u
    \right)X^\star
    \right),
    \\
    m_{\operatorname{Reg}}(t)
    \defas
    e^{-K\tilde{\Gamma}^{\operatorname{Reg}}(t)-\Lambda^{\operatorname{Reg}}(t)}
    \left(
    X_0+
    \left(
    \int_0^t
    K\gamma(u)A_{\operatorname{Reg}}(u)
    e^{K\tilde{\Gamma}^{\operatorname{Reg}}(u)+\Lambda^{\operatorname{Reg}}(u)}
    \dif u
    \right)X^\star
    \right).
\end{gathered}
\]
The corresponding deterministic first terms are
\[
\begin{gathered}
    G_{\mathrsfs{R}}(t)
    \defas
    \frac12
    \left\langle
    K,
    (m_{\Adv}(t)-X^\star)^{\otimes 2}
    \right\rangle
    +
    \frac{\eta^2}{2}, \quad G_{\mathrsfs{B}}(t)
    \defas
    \|m_{\Adv}(t)\|^2,
    \\
    G_{\mathscr{R}}(t)
    \defas
    \frac12
    \left\langle
    K,
    (m_{\operatorname{Reg}}(t)-X^\star)^{\otimes 2}
    \right\rangle
    +
    \frac{\eta^2}{2}, \quad G_{\Phi }(t)
    \defas
    \|m_{\operatorname{Reg}}(t)\|^2.
\end{gathered}
\]
We now compare $m_{\Adv}$ and $m_{\operatorname{Reg}}$ using their explicit 
representations. For conciseness, let $P_{\Adv}(t)
\defas
e^{-K\tilde{\Gamma}(t)-\Lambda(t)}.$ and $P_{\operatorname{Reg}}(t)
\defas
e^{-K\tilde{\Gamma}^{\operatorname{Reg}}(t)-\Lambda^{\operatorname{Reg}}(t)}$ then we have
\[
\begin{aligned}
    m_{\Adv}(t)-&m_{\operatorname{Reg}}(t)
    =
    \left(
    P_{\Adv}(t)-P_{\operatorname{Reg}}(t)
    \right)
    \left(
    X_0+
    \int_0^t
    K\gamma(u)A_{\Adv}(u)
    (P_{\Adv}(u))^{-1}X^\star
    \dif u
    \right)
    \\
    &
    +
    P_{\operatorname{Reg}}(t)
    \int_0^t
    K\gamma(u)
    \left(
    A_{\Adv}(u)(P_{\Adv}(u))^{-1}
    -
    A_{\operatorname{Reg}}(u)(P_{\operatorname{Reg}}(u))^{-1}
    \right)
    X^\star
    \dif u.
\end{aligned}
\]
Once again from $ |e^{-x}-e^{-y}|\leq |x-y|$ for all $x,y\geq 0$, we obtain
\[
    \|P_{\Adv}(t)-P_{\operatorname{Reg}}(t)\|_{\opt}
    \leq
    C(\|K\|_{\opt}, \bar{\gamma})
    \int_0^t
    D(s)\dif s.
\]
Similarly, we obtain $
    \|P_{\Adv}(u)^{-1}-P_{\operatorname{Reg}}(u)^{-1}\|_{\opt}
    \leq
    C(\|K\|_{\opt}, \bar{\gamma})
    \int_0^u
   D(s)\dif s$. 
 Since $\|(P_{\Adv}(t))^{-1}\|_\opt$, $\|(P_{\operatorname{Reg}}(t))^{-1}\|_\opt \leq C(T, \|K\|_{\opt}, \bar{\gamma})$ and from \eqref{eq:A:B:lipschitz}, it follows that
\[
    \|m_{\Adv}(t)-m_{\operatorname{Reg}}(t)\|
    \leq
    C(T, \|K\|_{\opt}, \bar{\gamma})
    \int_0^t
    D(s)\dif s.
\]
Thus, we obtain the following bounds on the first terms
\begin{equation}\label{eq:G:R:difference}
\begin{gathered}
\left|
    G_{\mathrsfs{R}}(t)
-
G_{\mathscr{R}}(t)
\right|
\leq
C(T, \|K\|_{\opt}, \bar{\gamma})
\int_0^t
D(s)\dif s
\\
\text{and}\quad
\left|
    G_{\mathrsfs{B}}(t)
-
G_{\Phi}(t)
\right|
\leq
C(T, \|K\|_{\opt}, \bar{\gamma})
\int_0^t
D(s)\dif s.
\end{gathered}
\end{equation}
 For the difference in risks, we have
    \[
    \begin{aligned}
        \mathrsfs{R}(t)-\mathscr{R}(t)
        &=
        G_{\mathrsfs{R}}(t)
        -
        G_{\mathscr R}(t)
        +
        \frac1d
        \int_0^t
        \gamma^2(s)
        \left[
        \mathcal{F}_{\Adv}(s)\mathcal{K}_{2}^{\Adv}(t,s)
        -
        \mathcal{R}_{\operatorname{Reg}}(s)\mathcal{K}_{2}^{\operatorname{Reg}}(t,s)
        \right]
        \dif s.
    \end{aligned}
    \]
    Thus, since the boundedness assumption implies $ |\mathcal{R}_{\operatorname{Reg}}(t)|\leq C$ and by \eqref{eq:kernel:bounded},
    \eqref{eq:kernel:lipschitz}, \eqref{eq:forcing:difference} and \eqref{eq:G:R:difference} we obtain
    \begin{equation} \label{eq:R:difference:bound}
    \begin{aligned}
        \left|\mathrsfs{R}(t)-\mathscr{R}(t)\right|
        &\leq
        C(T, \|K\|_{\opt}, \bar{\gamma})
\int_0^t
D(s)\dif s
        +
        C(\|K\|_{\opt}, \bar{\gamma})
        \delta^2
        \left(
        \frac{1}{2}-\frac{1}{\pi}
        \right)
        \int_0^t
        \widehat{\mathrsfs{B}}_{44}(s)\dif s.
    \end{aligned}
\end{equation}
    Using a similar argument, we obtain
    \begin{equation}
    \label{eq:B:difference:bound}
    \begin{aligned}
        \left|
        \widehat{\mathrsfs{B}}_{44}(t)-\Phi(t)
        \right|
        &\leq
        C(T, \|K\|_{\opt}, \bar{\gamma})
        \int_0^t
        D(s)\dif s
                +
                C(\|K\|_{\opt}, \bar{\gamma})
                \delta^2
                \left(
                \frac{1}{2}-\frac{1}{\pi}
                \right)
                \int_0^t
                \widehat{\mathrsfs{B}}_{44}(t)\dif s.
    \end{aligned}
    \end{equation}
    Adding \eqref{eq:R:difference:bound} and \eqref{eq:B:difference:bound}, we obtain
    \begin{equation}
    \label{eq:D:Grönwall:pre}
        D(t)
        \leq
        C(T, \|K\|_{\opt}, \bar{\gamma})
        \int_0^t
        D(s)\dif s
                +
                C(\|K\|_{\opt}, \bar{\gamma}, T)
                \delta^2
                \left(
                \frac{1}{2}-\frac{1}{\pi}
                \right)
                \int_0^t
                \widehat{\mathrsfs{B}}_{44}(t)\dif s.
    \end{equation}
    An application of Grönwall's inequality and taking the supremum over $[0,T]$ completes the proof.
    \end{proof}
\subsubsection{Analysis of Learning Dynamics}\label{proof:prop:q(t):bounded}
In this section, we provide the proofs of the results presented in Section~\ref{sect:lsadv:analysis:main}. We start by showing $q(t)$ is bounded above and below using a barrier argument conditionally on $\mathrsfs{R}(t)$ remaining positive.
\begin{proposition}\label{prop:q(t):bound:cond}
    For any $T>0$, if $\mathrsfs{R}(t) >0$ for all $t \in [0,T]$ then for all $t \in [0, T]$
    \begin{equation}
        \min(q(0), q_L) \leq q(t) \leq \max(q(0), q_U),
    \end{equation}
    where $q_L$ and $q_U$ are defined as follows
    \begin{equation}
        \begin{gathered}
            q_L = \min\left(\frac{1}{\delta}, \sqrt{\frac{\Tr(K)}{2 \Tr(K^2)}}, \frac{\gamma \Tr(K)}{16 d \left( \|K\|_{\opt}^{1/2} + \delta \sqrt{\frac{2}{\pi}} \right)}, \sqrt{\frac{\gamma \Tr(K)}{16 d \left(2 \|K\|_{\opt} + \delta^2 \right)}}\right)
            \\
            q_U = \max\left(1, \frac{1}{\delta},\sqrt{\frac{2\Tr(K)}{ \Tr(K^2)}}, \frac{8d \left( \|K\|_{\opt} \left( 1+ \delta \sqrt{\frac{2}{\pi}}\right) + \delta^2 (1 + \|K\|^{1/2}_{\opt})\right)}{\gamma\delta^2 \Tr(K^2)} \right).
        \end{gathered}
    \end{equation}
Suppose further that there exists a $c >0$ such that $\frac{\Tr(K^2)}{d} \geq c$ uniformly in $d$ and $0 < c_0 \leq q(0) \leq C_0 < \infty$ where $c_0, C_0 >0$ are independent of $d$. Then we obtain the following bounds uniform in $d$:
\begin{equation}
    \begin{gathered}
        q_L = \min\left(\frac{1}{\delta}, \frac{1}{\sqrt{2 \|K\|_{\opt}}}, \frac{c\gamma}{16 \|K\|_{\opt} \left( \|K\|_{\opt}^{1/2} + \delta \sqrt{\frac{2}{\pi}} \right)}, \sqrt{\frac{c\gamma}{16  \|K\|_{\opt} \left(2 \|K\|_{\opt} + \delta^2 \right)}}\right),
        \\
        q_U = \max\left(1, \frac{1}{\delta},\sqrt{ \frac{2\|K\|_{\opt}}{c}}, \frac{8  \left( \|K\|_{\opt} \left( 1+ \delta \sqrt{\frac{2}{\pi}}\right) + \delta^2 (1 + \|K\|^{1/2}_{\opt})\right)}{\gamma \delta^2 c} \right).
    \end{gathered}
\end{equation}
\end{proposition}
\begin{proof}
    For simplicity of notation, we define the following notation
    \begin{equation}
        \begin{gathered}
        \operatorname{term}_1 \defas \frac{\gamma}{\mathrsfs{R}(t)} (q(t))^2\left( 1 + \delta\sqrt{\frac{2}{\pi}}q(t) \right) \sum_{j=1}^d \lambda_j^2 \left(\mathrsfs{V}_{j,11}(t) - 2 \mathrsfs{V}_{j,12}(t) + \mathrsfs{V}_{j,22}(t)\right) ,
        \\
        \operatorname{term}_2 \defas \frac{\gamma}{\mathrsfs{R}(t)} (\delta^2 (q(t))^2 - 1)  \sum_{j=1}^d \lambda_j\left(\mathrsfs{V}_{j,11}(t) - \mathrsfs{V}_{j,12}(t) \right),
        \\
        \operatorname{term}_3 \defas -2 \gamma q(t)\left(\delta^2 q(t) + \delta \sqrt{\frac{2}{\pi}}\right),
        \\
        \operatorname{term}_4 \defas  \frac{\gamma^2}{d} \left(1+ 2 \delta \sqrt{\frac{2}{\pi}} q(t) + \delta^2 (q(t))^2\right)(\Tr(K) - \Tr(K^2) (q(t))^2).
        \end{gathered}
    \end{equation}
    We will refer to each term using ordinals. Suppose $0 <q(t) \leq q_L$. Since $\mathrsfs{V}_j(t)$ is psd then $ \mathrsfs{V}_{j,11}(t) - 2 \mathrsfs{V}_{j,12}(t) + \mathrsfs{V}_{j,22}(t) \geq 0$ which implies the first term (i.e. $\operatorname{term}_1$) is nonnegative. Since $q(t) \leq q_L$ then 
    \[
        q(t) \leq \sqrt{\frac{\Tr(K)}{2 \Tr(K^2)}},
    \]
    and thus
    \begin{equation}
        \Tr(K) - (q(t))^2 \Tr(K^2) \geq \frac{1}{2}\Tr(K).
    \end{equation}
    Hence, we obtain 
    \begin{equation}
        \operatorname{term}_4 \geq \frac{\gamma^2}{2d}\Tr(K).
    \end{equation}
    Since $q(t) \leq q_L$ then $q(t) \leq \frac{1}{\delta}$ which implies $|\delta^2 (q(t))^2 - 1| \leq 1$. From the second bound in Lemma~\ref{lem:helper:bounds:diffforQ} then we obtain
    \begin{equation}
        \operatorname{term}_2 \geq -\gamma \left(4 \|K\|_{\opt}  (q(t))^2+ 2 \|K\|_{\opt}^{1/2}q(t)\right).
    \end{equation}
 Hence, for all $0 \leq q(t) \leq q_L$ we have
    \begin{equation}
        \begin{aligned}
        \frac{\dif}{\dif t} (q(t))^2 \geq  \frac{\gamma^2}{2d}\Tr(K) -2\gamma\left(2 \|K\|_{\opt} +\delta^2  \right)(q(t))^2 - 2\gamma \left( \|K\|_{\opt}^{1/2} +   \delta \sqrt{\frac{2}{\pi}}\right)q(t).
        \end{aligned}
    \end{equation}
    Since $q(t) \leq q_L$ then we have
    \[
    \begin{gathered}
        2\gamma\left(2 \|K\|_{\opt} +\delta^2  \right)(q(t))^2 \leq  \frac{\gamma^2}{8d}\Tr(K),
        \\
        \text{and}\quad 2\gamma \left( \|K\|_{\opt}^{1/2} +   \delta \sqrt{\frac{2}{\pi}}\right)q(t) \leq  \frac{\gamma^2}{8d}\Tr(K).
    \end{gathered}
    \]
    Thus, whenever $0 \leq q(t) \leq q_L$ then 
    \begin{equation}
        \frac{\dif}{\dif t} (q(t))^2 \geq  \frac{\gamma^2}{4d}\Tr(K) \geq 0.
    \end{equation}
    Now, let $q(t) \geq q_U$ then $q(t) \geq 1$, $q(t) \geq \frac{1}{\delta}$ and 
    \[
        q(t) \geq \sqrt{\frac{2\Tr(K)}{\Tr(K^2)}},
    \]
    from which it follows that 
    \begin{equation}
        \Tr(K) - (q(t))^2 \Tr(K^2) \leq -\frac{1}{2}(q(t))^2 \Tr(K^2).
        \end{equation}
    Thus, we obtain the bound for the fourth term
    \begin{equation}
        \operatorname{term}_4 \leq- \frac{\gamma^2 \delta^2 }{2d}(q(t))^4  \Tr(K^2).
    \end{equation}
    Since $\lambda_j^2 \leq \|K\|_{\opt} \lambda_j$ then 
    \[
        \sum_{j=1}^d \lambda_j^2 \left(\mathrsfs{V}_{j,11}(t) - 2 \mathrsfs{V}_{j,12}(t) + \mathrsfs{V}_{j,22}(t)\right) \leq 2\|K\|_{\opt} \mathrsfs{R}(t).
    \]
    Using this and since $q(t) \geq 1$ then for the first term we obtain
    \begin{equation}
        \operatorname{term}_1 \leq 2\gamma \|K\|_{\opt}\left( 1 + \delta\sqrt{\frac{2}{\pi}} \right)  (q(t))^3.
    \end{equation}
    Since $q(t) \geq \frac{1}{\delta}$ then $\delta^2 (q(t))^2 - 1 \geq 0$. Using this, $q(t) \geq 1$ and the first bound in Lemma~\ref{lem:helper:bounds:diffforQ}, we obtain
    \begin{equation}
        \operatorname{term}_2 \leq 2\gamma \delta^2  \left(\|K\|_{\opt}^{1/2}  +1\right)(q(t))^3.
    \end{equation}
    Finally, the third term is nonpositive so we drop it to obtain the upper bound
    \begin{equation}
        \begin{aligned}
        \frac{\dif}{\dif t} (q(t))^2 \leq  2\gamma\left( \|K\|_{\opt}\left( 1 + \delta\sqrt{\frac{2}{\pi}} \right)  + \delta^2  \left(\|K\|_{\opt}^{1/2}  +1\right)\right) (q(t))^3 -\frac{\gamma^2 \delta^2 }{2d}(q(t))^4  \Tr(K^2).
        \end{aligned}
    \end{equation}
    Since $q(t) \geq q_U$ then we have
    \[
        2\gamma\left( \|K\|_{\opt}\left( 1 + \delta\sqrt{\frac{2}{\pi}} \right)  + \delta^2  \left(\|K\|_{\opt}^{1/2}  +1\right)\right) (q(t))^3 \leq \frac{\gamma^2 \delta^2 }{4d}(q(t))^4  \Tr(K^2),
    \]
    from which it follows that for all $q(t) \geq q_U$ 
    \begin{equation}
        \frac{\dif}{\dif t} (q(t))^2 \leq -\frac{\gamma^2 \delta^2 }{4d}(q(t))^4  \Tr(K^2) < 0.
    \end{equation}
    Hence, by a first exit argument, we obtain the first bound. To obtain the bound uniform in $d$, we note that given there exists a $c>0$ such that $\frac{1}{d}\Tr(K^2) \geq c >0$ then 
    \[
       \Tr(K) = \sum_{j=1}^d \lambda_j \geq \frac{1}{\|K\|_\opt} \sum_{j=1}^d \lambda_j^2 = \frac{\Tr(K^2 )}{\|K\|_\opt}.
    \]
    Hence, $\frac{1}{d}\Tr(K) \geq \frac{c}{\|K\|_{\opt}}$ from which it follows that 
    \begin{equation}
        \frac{\Tr(K)}{\Tr(K^2)} \geq \frac{1}{\|K\|_{\opt}}\quad \text{and}\quad  \frac{\Tr(K)}{\Tr(K^2)} \leq \frac{\|K\|_{\opt}}{c}.
    \end{equation}
\end{proof}
Conditionally on $\mathrsfs{R}(0) >0$, we can extend Proposition~\ref{prop:q(t):bound:cond} to all time. Note that the argument is not necessary if $\eta>0$ as then $\mathrsfs{R}(t) \geq \frac{\eta^2}{2}$ for all $t\geq 0$.
\begin{proof}[Proof of Proposition~\ref{prop:q(t):bounded}]
    Without loss of generality, we assume $q_L \leq q(0) \leq q_U$. We will proceed by contradiction. Given the definition of $\mathrsfs{R}^\Adv$, if $\mathrsfs{R}(t) >0$ then $\mathrsfs{R}^\Adv(t) >0$. Suppose there exists a $\tau >0$ such that $\mathrsfs{R}(\tau)= 0$. Formally, $\tau$ is defined
    \begin{equation}
        \tau = \inf \{ t>0: \mathrsfs{R}(t) = 0\} < \infty.
    \end{equation}
    By assumption, we have $\mathrsfs{R}(0) >0$ and since $\mathrsfs{R}(t)$ is continuous then $\mathrsfs{R}(t) >0 $ for every $t\in [0,\tau)$. By Proposition~\ref{prop:q(t):bound:cond} then $q(t)$ is well-defined on $[0,\tau)$ such that $0< q_L \leq q(t) \leq q_U <\infty$. On $[0,\tau)$, we have
    \[
        \tilde{\Gamma}'(t) = \gamma \left(1 + \delta \sqrt{\frac{2}{\pi}}q(t)\right) \leq \gamma \left(1 + \delta \sqrt{\frac{2}{\pi}}q_U\right) .
    \]
    Thus, we obtain
    \[
        \tilde{\Gamma}(t) \leq \gamma \left(1 + \delta \sqrt{\frac{2}{\pi}}q_U\right) \tau,
    \]
    and given $\tilde{\Gamma}(t)$ is nondecreasing then $\tilde{\Gamma}(\tau^-) = \lim_{t\uparrow \tau} \tilde{\Gamma}(t) < \infty$. Similarly, we have
    \[
        \Lambda'(t) = \gamma\left(\delta^2 + \delta \sqrt{\frac{2}{\pi}}\frac{1}{q(t)}\right)  \leq  \gamma\left(\delta^2 + \delta \sqrt{\frac{2}{\pi}} \frac{1}{q_L}\right),
    \]
    from which it follows as $\Lambda(t)$ is nondecreasing that $\Lambda(\tau^-) = \lim_{t\uparrow \tau} \Lambda(t) < \infty$. Set $0 <\varepsilon <\tau$. Then for every $s \in [0,\varepsilon]$, we have $\mathrsfs{R}(s) >0$. Since $\mathrsfs{R}$ is continuous then it attains a minimum on $[0,\varepsilon]$ which we denote $m_{\varepsilon} = \min_{0 \leq s\leq \varepsilon} \mathrsfs{R}(s)$. Now, for any $\varepsilon \leq t <\tau$ we have
    \[
    \begin{gathered}
    0 \leq  \tilde{\Gamma}(t) - \tilde{\Gamma}(s) \leq \tilde{\Gamma}(\tau^-) - \tilde{\Gamma}(0)  = \tilde{\Gamma}(\tau^-) ,
    \\
    0 \leq  \Lambda(t) - \Lambda(s) \leq \Lambda(\tau^-) -\Lambda(0)  = \Lambda(\tau^-) ,
    \end{gathered}
    \]
    from which it follows that 
    \[
        e^{-2\lambda_j (\tilde{\Gamma}(t) - \tilde{\Gamma}(s)) - 2(\Lambda(t) - \Lambda(s))} \geq e^{-2 (\|K\|_{\opt} \tilde{\Gamma}(\tau^-) +\Lambda(\tau^-) )}.
    \]
    Thus, we have $ \tr \big ( K^2 e^{-2K(\tilde{\Gamma}(t) - \tilde{\Gamma}(s)) -2(\Lambda(t) -\Lambda(s))} \big ) \geq e^{-2 (\|K\|_{\opt} \tilde{\Gamma}(\tau^-) +\Lambda(\tau^-) )} \Tr(K^2) > 0$ which implies from \eqref{eq:R(t):equiv:final:lr}
    \begin{equation}
        \begin{aligned}
        \mathrsfs{R}(t) 
        &
        \geq 
         \frac{\gamma^2}{d} \int_0^t \mathrsfs{R}(s)\tr \big ( K^2 e^{-2K(\tilde{\Gamma}(t) - \tilde{\Gamma}(s)) -2(\Lambda(t) -\Lambda(s))} \big ) \, \dif s 
         \\
         &
         \geq
         \frac{\gamma^2}{d} \int_0^\varepsilon m_{\varepsilon} e^{-2 (\|K\|_{\opt} \tilde{\Gamma}(\tau^-) +\Lambda(\tau^-) )} \Tr(K^2) \, \dif s
         \\
         &
         \geq   \frac{\gamma^2}{d} \Tr(K^2)m_{\varepsilon}e^{-2 (\|K\|_{\opt} \tilde{\Gamma}(\tau^-) +\Lambda(\tau^-) )} \varepsilon >0 .
        \end{aligned}
    \end{equation}
    Hence, $\mathrsfs{R}(t) >0 $ for all $t \in [\varepsilon, \tau)$ which implies since $\mathrsfs{R}(t)$ is continuous that $\mathrsfs{R}(\tau) = \lim_{t\uparrow \tau} \mathrsfs{R}(t) \geq \frac{\gamma^2}{d} \Tr(K^2)m_{\varepsilon}e^{-2 (\|K\|_{\opt} \tilde{\Gamma}(\tau^-) +\Lambda(\tau^-) )} \varepsilon > 0$ which is a contradiction. Therefore, $\mathrsfs{R}(t) >0$ for every $t >0$. Hence, Proposition~\ref{prop:q(t):bound:cond} is valid on $[0,\infty)$.
\end{proof}
Before the proof of Proposition~\ref{prop:fixed:point:main}, we introduce the following notation to simplify the derivations:
\begin{equation}\label{eq:def:rho_j(t)}
    \begin{gathered}
        \rho_j(t) \defas  e^{-\lambda_j \tilde{\Gamma}(t) - \Lambda(t)} \left( x_{0,j} + \lambda_j x_{j}^\star \int_0^t \tilde{\gamma}^{\operatorname{eff}}(s) e^{\lambda_j \tilde{\Gamma}(s) + \Lambda(s)} \dif s \right),
        \\
        \text{where}\quad x_{0,j} \defas X_0^\top \omega_j \quad \text{and}\quad x_{j}^\star \defas (X^\star)^\top \omega_j.
    \end{gathered}
\end{equation}
Define $\ell(u) \defas \Lambda(\tilde{\Gamma}^{-1}(u))$ then $\ell'(u) = \tilde{\lambda}^{\operatorname{eff}}(\tilde{\Gamma}^{-1}(u))$. We now express $\rho_j(t)$ as a function of $\tilde{\Gamma}$ using the time change $u= \tilde{\Gamma}(s)$:
\begin{equation}
    \begin{gathered}
        \rho_j(u ) \defas  e^{-\lambda_j u - \ell(u) } \left( x_{0,j} + \lambda_j x_{j}^\star \int_0^{u} e^{\lambda_j s + \ell(s)} \dif s \right).
    \end{gathered}
\end{equation}
We will work with each definition of $\rho_j$ depending on the context. Given that $\frac{\dif}{\dif s} e^{\lambda_j s + \ell(s)} = (\lambda_j  + \ell'(s))e^{\lambda_j s + \ell(s)}$ then
\[
    \int_0^u \lambda_j e^{\lambda_j s + \ell(s)}\dif s = e^{\lambda_j u + \ell(u)} - 1 - \int_0^u\ell'(s)e^{\lambda_j s + \ell(s)} \dif s.
\]
From this identity, we also obtain the representation
\begin{equation}
    \rho_j(u ) \defas x^\star_j + e^{-\lambda_j u - \ell(u) } \left( x_{0,j} - x^\star_j\right) - x^\star_j \int_0^{u} \ell'(s) e^{-\lambda_j (u-s) - (\ell(u)-\ell(s))} \dif s .
\end{equation}
Using these identities, we define the following functions
\begin{equation}\label{eq:def:F:lr}
    \begin{aligned}
        F(u) 
        &
        \defas  
        \frac{1}{2}\sum_{j=1}^d \lambda_j (\rho_j(u) - x^\star_j)^2
        \\
        &
=
        \frac{1}{2}\sum_{j=1}^d \lambda_j (e^{-\lambda_j u - \ell(u) } \left( x_{0,j} - x^\star_j\right) - x^\star_j \int_0^{u} \ell'(s) e^{-\lambda_j (u-s) - (\ell(u)-\ell(s))} \dif s)^2,
    \end{aligned}
\end{equation}
and
\begin{equation}\label{eq:def:L:lr}
    \begin{aligned}
        L(u) 
        &
        \defas  
        \sum_{j=1}^d (\rho_j(u))^2=
        \sum_{j=1}^d  e^{-2(\lambda_j u + \ell(u)) } \left( x_{0,j} + \lambda_j x_{j}^\star \int_0^{u} e^{\lambda_j s + \ell(s)} \dif s \right)^2.
    \end{aligned}
\end{equation}
We also introduce the following notation
\[
\begin{gathered}
    M_j(t) \defas \gamma^2 \int_0^t\mathrsfs{R}^{{\operatorname{adv}}}(s)  \lambda_j^2  e^{-2\lambda_j(\tilde{\Gamma}(t) - \tilde{\Gamma}(s)) -2(\Lambda(t) -\Lambda(s))}  \, \dif s,
    \\
    \text{and}\quad N_j(t)\defas  2\gamma^2 \int_0^t\mathrsfs{R}^{{\operatorname{adv}}}(s)  \lambda_j  e^{-2\lambda_j(\tilde{\Gamma}(t) - \tilde{\Gamma}(s)) -2(\Lambda(t) -\Lambda(s))}  \, \dif s.
\end{gathered}
\]
This leads to
        \[
        \begin{gathered}
            \mathrsfs{R}(t)= F(\tilde{\Gamma}(t)) 
            + \frac{1}{d}\sum_{j=1}^d M_j(t),
            \quad
            \text{and}\quad  \widehat{\mathrsfs{B}}_{44}(t)= L(\tilde{\Gamma}(t)) 
            + \frac{1}{d}\sum_{j=1}^d N_j(t).
        \end{gathered}
        \]
We now prove a more general version of Proposition~\ref{prop:fixed:point:main}.
    \begin{proposition}\label{prop:fixed:point}
            Suppose $\mathrsfs{R}(0), \widehat{\mathrsfs{B}}_{44}(0) >0$, $(X^\star)^\top K X^\star> 0$, $\sup_{t\geq 0} \mathrsfs{R}^\Adv(t)< \infty$ and $q(t) \to q_{\infty} \in (0,\infty)$, $\delta >0$. If, for all $q \in [\min(q(0), q_L ), \max(q(0), q_U) ]$, the following condition is satisfied
        \begin{equation}\label{eq:unique:fp:cond}
            \frac{B_d(q)}{A_d(q)}
            \geq
            \max\left(
                \frac{P_\delta(q)}{Q_\delta(q)},0
                \right) \defas \psi_{\delta}(q),
        \end{equation}
        where 
        \[
        \begin{gathered}
            A_d(q)
            \defas
            \frac{1}{d}\sum_{j=1}^d
            \frac{
            \lambda_j
            }{
           \delta^2 \left(
            \delta q+\sqrt{\frac{2}{\pi}}
            +
            r_j\delta q(1+\sqrt{\frac{2}{\pi}}\delta q)
            \right)^2
            },
            \\
             B_d(q)
            \defas
            \frac{1}{d}\sum_{j=1}^d
            \frac{
            \lambda_j^2 
            }{
            \delta^4\left(
            \delta q+\sqrt{\frac{2}{\pi}}
            +
            r_j\delta q(1+\sqrt{\frac{2}{\pi}}\delta q)
            \right)^2
            }
            \\
            P_\delta(q)
            \defas
            \left(2-\frac{6}{\pi}\right)(\delta q)^4
            +
            \left(
            3\sqrt{\frac{2}{\pi}}
            -
            7\left(\frac{2}{\pi}\right)^{3/2}
            \right)(\delta q)^3  \\
            -
            \left(
            4\left(\frac{2}{\pi}\right)^2
            +
           \frac{4}{\pi}
            \right)(\delta q)^2
            -
            \left(
            3\left(\frac{2}{\pi}\right)^{3/2}
            +
            \sqrt{\frac{2}{\pi}}
            \right)(\delta q)
            -
            \frac{2}{\pi},
            \\
        \text{and}
        \quad
            Q_\delta(q)
            \defas
            \frac{2}{\pi}(\delta q)^6
            +
            \left(
            7\left(\frac{2}{\pi}\right)^{3/2}
            -
            \sqrt{\frac{2}{\pi}}
            \right)(\delta q)^5 
            +
            \left(
            8\left(\frac{2}{\pi}\right)^2
            +
           \frac{16}{\pi}
            -
            2
            \right)(\delta q)^4  \\
            +
            \left(
            15\left(\frac{2}{\pi}\right)^{3/2}
            +
            \sqrt{\frac{2}{\pi}}
            \right)(\delta q)^3
            +
            \frac{18}{\pi}(\delta q)^2
            +
            2\sqrt{\frac{2}{\pi}}\delta q.
        \end{gathered}
        \]
        Then $q_\infty$ is the unique fixed point satisfying the fixed point equation
        \begin{equation}\label{eq:fixed:point}
            q^2_\infty = (1 - G(q_\infty) C_{\mathrsfs{R}}(q_\infty)) \frac{L_\infty(q_\infty)}{2 F_\infty(q_\infty)} + \frac{1}{2}G(q_\infty)C_\mathrsfs{B}(q_\infty)
        \end{equation}
        and the admissibility condition $G(q_\infty)C_{\mathrsfs{R}}(q_\infty) < 1$ for $G(q) = 1+ 2 \delta \sqrt{\frac{2}{\pi}} q + \delta^2 q^2$ and $C_\mathrsfs{R}(q) = \frac{\gamma}{2d}\sum_{j=1}^d \frac{\lambda_j^2}{D_j(q)}$ with $D_j(q) = \lambda_j(1+ \delta \sqrt{\frac{2}{\pi}}q)  + \delta^2 + \delta \sqrt{\frac{2}{\pi}}\frac{1}{q}$. The limiting values of $\mathrsfs{R}$, $\widehat{\mathrsfs{B}}_{44}$ and $\mathrsfs{R}^\Adv$ are given by
        \begin{equation*}
            \begin{gathered}
                \mathrsfs{R}(\infty) = \frac{F_{\infty}(q_\infty)}{1- G(q_\infty)C_{\mathrsfs{R}}(q_\infty)},\quad 
                \widehat{\mathrsfs{B}}_{44}(\infty) = 2 q^2_\infty  \mathrsfs{R}(\infty),\quad
                \text{and}\quad \mathrsfs{R}^\Adv(\infty) = G(q_\infty)\mathrsfs{R}(\infty),
            \end{gathered}
        \end{equation*}
    where $C_{\mathrsfs{B}}(q) = \frac{\gamma}{d}\sum_{j=1}^d \frac{\lambda_j}{D_j(q)}$, $F_\infty(q) = \frac{1}{2}\sum_{j=1}^d \lambda_j \left(\frac{\tilde{\lambda}^{\operatorname{eff}}(q)}{\lambda_j +\tilde{\lambda}^{\operatorname{eff}}(q)}\right)^2(x_j^\star)^2$, $ L(\infty) = \sum_{j=1}^d  \left(\frac{\lambda_j}{\lambda_j + \tilde{\lambda}^{\operatorname{eff}}(\infty)}\right)^2( x_j^\star)^2$, and $\tilde{\lambda}^{\operatorname{eff}}(q) = \frac{\delta^2 + \delta \sqrt{\frac{2}{\pi}} \frac{1}{q}}{1 + \delta \sqrt{\frac{2}{\pi}}q }$. In the case of isotropic covariance, the condition in \eqref{eq:unique:fp:cond} simplifies to
    \[
            \delta^2
            \geq
            \max\left(
                \frac{P_\delta(q)}{Q_\delta(q)},0
                \right).
    \]
    Furthermore, let $ r_\star
        \defas
        \sup_{q>0}\psi_\delta(q)$ with $r_\star\approx 8.35\times 10^{-5}$. Then a less precise sufficient condition than \eqref{eq:unique:fp:cond} is 
        \[
            \delta \leq \sqrt{\frac{\lambda_{\operatorname{min}}^+}{ r_\star}},
        \] 
        which leads to $\delta \leq r_\star^{-1/2}$ for isotropic covariance.
    \end{proposition}
    \begin{proof}
        Recall we have
        \begin{equation}
            \tilde{\Gamma}'(t) = \tilde{\gamma}^{\operatorname{eff}}(t) = \gamma\left(1+ \delta \sqrt{\frac{2}{\pi}}q(t)\right) \quad \text{and}\quad \Lambda'(t) = \tilde{\gamma}^{\operatorname{eff}}(t) \tilde{\lambda}^{\operatorname{eff}}(t) = \gamma\left(\delta^2+ \delta \sqrt{\frac{2}{\pi}}\frac{1}{q(t)}\right).
        \end{equation}
        Recall $\ell(u) = \Lambda(\tilde{\Gamma}^{-1}(u))$ then $\ell'(u) = \tilde{\lambda}^{\operatorname{eff}}(\tilde{\Gamma}^{-1}(u))$. Given that $\tilde{\Gamma}(t) \geq \gamma t$ we have $\tilde{\Gamma}^{-1}(u) \to \infty$ when $u \to \infty$ which implies $q(\tilde{\Gamma}^{-1}(u)) \to q_\infty$ and 
        \[
            \ell'(u) \to \tilde{\lambda}^{\operatorname{eff}}(\infty) = \frac{\delta^2 + \delta\sqrt{\frac{2}{\pi}} \frac{1}{q_\infty}}{1 + \delta \sqrt{\frac{2}{\pi}}q_\infty}.
        \]
        We now compute the limits of $F(u)$ and $L(u)$. From the definition of $ \rho_j(u ) $, it satisfies the following ODE
        \[
            \frac{\dif}{\dif u}\rho_j(u ) = - (\lambda_j + \ell'(u)) \rho_j(u ) + \lambda_j  x_{j}^\star.
        \]
        We will show that $\rho_j(u )  \to \rho_j(\infty) \defas \frac{\lambda_j}{\lambda_j + \tilde{\lambda}^{\operatorname{eff}}(\infty)}x_{j}^\star$. Let $y_j(u) \defas \rho_j(u ) - \rho_j(\infty)$ then $y'_j(u) = \rho'_j(u )$ from which we obtain after substituting $\rho_j(u ) = y_j(u) + \rho_j(\infty)$ and $\lambda_j x_j^\star = (\lambda_j + \tilde{\lambda}^{\operatorname{eff}}(\infty))\rho_j(\infty)$
        \begin{equation}
            y'_j(u) = - (\lambda_j + \ell'(u)) y_j(u)  + (\tilde{\lambda}^{\operatorname{eff}}(\infty) - \ell'(u))\rho_j(\infty)
        \end{equation}
        This ODE takes the form $y'(u) = - A(u) y(u) + B(u)$ with $A_j(u) = (\lambda_j + \ell'(u)) $ and $B_j(u) = (\tilde{\lambda}^{\operatorname{eff}}(\infty) - \ell'(u))\rho_j(\infty)$. Since $\ell'(u) \to \tilde{\lambda}^{\operatorname{eff}}(\infty)$ and $\lambda_j + \tilde{\lambda}^{\operatorname{eff}}(\infty) >0$, for large enough $u$, $A_j(u)$ is bounded below by a positive constant while $B_j(u) \to 0$ as $u\to \infty$. Thus, from variation of constants, we obtain $y_j(u) \to 0$ hence $\rho_j(u) \to \frac{\lambda_j}{\lambda_j + \tilde{\lambda}^{\operatorname{eff}}(\infty)}x_j^\star$. 

        From the definitions of $F$ and $L$, we obtain 
        \[
        \begin{gathered}
            F(\infty) =   \frac{1}{2}\sum_{j=1}^d  \lambda_j \left(\frac{\tilde{\lambda}^{\operatorname{eff}}(\infty)}{\lambda_j + \tilde{\lambda}^{\operatorname{eff}}(\infty)}\right)^2( x_j^\star)^2\quad \text{and}\quad
            L(\infty) = \sum_{j=1}^d  \left(\frac{\lambda_j}{\lambda_j + \tilde{\lambda}^{\operatorname{eff}}(\infty)}\right)^2( x_j^\star)^2.
        \end{gathered}
        \]
        Now, the adversarial risk can be written as $  \mathrsfs{R}^\Adv(t) = G(q(t)) \mathrsfs{R}(t)$ for $G(q) = 1+ 2\delta \sqrt{\frac{2}{\pi}}q+ \delta^2 q^2$. Since $\sup_{t\geq 0} \mathrsfs{R}^\Adv(t) < C< \infty$ for some constant $C>0$ then we have the uniform bounds
      \begin{equation}
        \begin{gathered}
        0 \leq M_j(t) \leq \gamma^2 \lambda_j^2 C\int_0^\infty e^{-2\gamma(\lambda_j + \delta^2)u} \dif u  =  \frac{C\gamma\lambda_j^2}{2 \left(\lambda_j + \delta^2\right)},
        \\
        0 \leq N_j(t) \leq  2\gamma^2\lambda_j C\int_0^\infty e^{-2\gamma(\lambda_j + \delta^2)u} \dif u  =  \frac{C\gamma\lambda_j}{\lambda_j + \delta^2} .
        \end{gathered}
    \end{equation}
        Proposition~\ref{prop:q(t):bounded} also implies there exists $q_-, q_+ >0$ such that $q_- \leq q(t) \leq q_+$. Differentiating $M_j(t)$ and $N_j(t)$, we obtain
        \[
            \begin{gathered}
            \frac{\dif}{\dif t}M_j(t) =-2 \gamma D_j(q(t)) M_j(t) + \gamma^2 \mathrsfs{R}^\Adv(t) \lambda_j^2  ,
            \\
            \frac{\dif}{\dif t}N_j(t) =-2 \gamma D_j(q(t)) N_j(t) + 2\gamma^2  \mathrsfs{R}^\Adv(t) \lambda_j, 
            \end{gathered}
        \]
        where we define $D_j(q) \defas \lambda_j \left(1+ \delta \sqrt{\frac{2}{\pi}} q\right) + \delta^2 + \delta \sqrt{\frac{2}{\pi}}\frac{1}{q}$. Since 
        \[
            \mathrsfs{R}^\Adv(t) = G(q(t)) \mathrsfs{R}(t)= G(q(t))\left(F(\tilde{\Gamma}(t)) 
            + \frac{1}{d}\sum_{j=1}^d M_j(t)\right)
        \]
         then 
        \begin{equation}\label{eq:M:fixedpoint}
            \frac{\dif}{\dif t}M_j(t) =-2 \gamma D_j(q(t)) M_j(t) + \gamma^2 \lambda_j^2 G(q(t))\left(F(\tilde{\Gamma}(t)) 
            + \frac{1}{d}\sum_{k=1}^d M_k(t)\right)
        \end{equation}
        Now, suppose $M_j(u) \to M_j(\infty)$ then the stationary equation would be
        \begin{equation}\label{eq:stationary:M}
            0=-2 \gamma D_j(q_\infty) M_j(\infty) + \gamma^2 \lambda_j^2 G(q_\infty)\left(F(\infty) 
            + \frac{1}{d}\sum_{j=1}^d M_j(\infty)\right).
        \end{equation}
        Let $S(t) = \frac{1}{d}\sum_{j=1}^d M_j(t)$ then $M_j(\infty) = \frac{\gamma \lambda_j^2}{2 D_j(q_\infty) } G(q_\infty)\left(F(\infty)  +  S(\infty)\right)$. Summing over $j$ and multiplying by $\tfrac{1}{d}$, we obtain
        \begin{equation}
            S(\infty) = G(q_\infty)C_{\mathrsfs{R}}(q_\infty)(F(\infty) + S(\infty)).
        \end{equation}
        Thus, if $G(q_\infty)C_{\mathrsfs{R}}(q_\infty) < 1$ then 
        \begin{equation}
            S(\infty) = \frac{G(q_\infty)C_{\mathrsfs{R}}(q_\infty)}{1- G(q_\infty)C_{\mathrsfs{R}}(q_\infty)}F(\infty).
        \end{equation}
        Given that $\mathrsfs{R}(t) = F(\tilde{\Gamma}(t)) + S(t)$ then the candidate limit for $\mathrsfs{R}$ is 
        \begin{equation}\label{eq:limit:R(inf)}
            \mathrsfs{R}(\infty) = F(\infty) + S(\infty) = \frac{F(\infty)}{1- G(q_\infty)C_{\mathrsfs{R}}(q_\infty)}.
        \end{equation}
        We will now show $M_j(t) \to M_j(\infty)$ as $t \to \infty$. Let $y_j(t) = M_j(t) - M_j(\infty)$ . Subtracting \eqref{eq:stationary:M} from \eqref{eq:M:fixedpoint}, we obtain
        \[
        y'_j(t) = -2 \gamma D_j(q_\infty)y_j(t)+ \gamma^2 \lambda_j^2 G(q_\infty)\sum_{k=1}^d y_k(t) + \operatorname{err}_j(t),
        \]
        where
        \[
        \begin{aligned}
            \operatorname{err}_j(t) 
            &
            = \gamma^2 \lambda_j^2 G(q(\infty))(F(\tilde{\Gamma}(t)) - F(\infty)) + \gamma^2\lambda_j^2\left(G(q(t)) -G(q(\infty))\right) \mathrsfs{R}(t)
            \\
            &
            -2\gamma\left(\lambda_j \delta \sqrt{\frac{2}{\pi}}(q(t)-q_\infty) + \delta \sqrt{\frac{2}{\pi}}\left(\frac{1}{q(t)} - \frac{1}{q_\infty}\right) \right) M_j(t).
        \end{aligned}
        \] 
        Since $\mathrsfs{R}(t) \leq \sup_{t\geq 0} \mathrsfs{R}^\Adv(t)< \infty$, $M_j(t)$ are uniformly bounded, $F(\tilde{\Gamma}(t)) \to F(\infty)$ and $q(t) \to q_\infty$ then the error term $ \operatorname{err}_j(t) \to 0 $ as $t\to \infty$. Now, the vector $y(t)$ composed of $y_j(t)$ satisfies the following ODE
        \[
            y'(t) = A_\infty y(t) + \operatorname{err}(t),
        \]
        where $A_\infty$ is defined such that 
        \[
            [A_\infty]_{jk} = -2 \gamma D_j(q_\infty) \mathds{1}_{jk} + \gamma^2 \lambda_j^2 G(q_\infty),
        \]
        where $\mathds{1}_{jk} = 1$ when $j = k$ and 0 if not. Hence, we have $A_\infty = - Q + u \mathbf{1}^\top$ where $\mathbf{1}$ is a vector of $1$, $D$ is diagonal such that 
        \[
            Q_{jj} = 2 \gamma D_j(q_\infty)\quad \text{and}\quad u_j = \gamma^2 \lambda_j^2 G(q_\infty).
        \]
        The condition $G(q_\infty)C_\mathrsfs{R}(q_\infty) < 1$ is equivalent to $\mathbf{1}^\top Q^{-1} u < 1$ as $\mathbf{1}^\top Q^{-1} u = G(q_\infty)C_\mathrsfs{R}(q_\infty)$. We will now show that $A_\infty$ is a Hurwitz-stable matrix. Let $\upsilon$ be an eigenvalue of $A_\infty$ then $\det(A_\infty - \upsilon \Id_d) = 0$. 
        
        We proceed by contradiction. Suppose there is an eigenvalue $\upsilon$ with $\Re[\upsilon] \geq 0$. Then $-Q - \upsilon \Id_d$ is invertible which implies $\det(A_\infty - \upsilon \Id_d) = \det(-Q - \upsilon \Id_d) (1- 1^\top (Q + \upsilon \Id_d)^{-1} u )$. 
        
        Now, every diagonal entry of $-Q - \upsilon \Id_d$ has negative real part hence $\det(-Q - \upsilon \Id_d) \not= 0$. Since $\upsilon$ is an eigenvalue then this enforces that $1= 1^\top (Q + \upsilon \Id_d)^{-1} u$. Since $\Re[\upsilon] \geq 0$ and $Q_{jj} > 0$ then taking absolute values, we obtain
        \[
            1 \leq \sum_{j=1}^d \frac{u_j}{|\upsilon + Q_{jj}|} \leq \sum_{j=1}^d \frac{u_j}{Q_{jj}} = \mathbf{1}^\top Q^{-1} u < 1,
        \]
        which is a contradiction. Hence $A_\infty$ is Hurwitz and it follows that there exists $C>0$ and $\kappa >0$ such that $\|e^{A_\infty t}\| \leq C e^{-\kappa t}$. Since $\operatorname{err}(t) \to 0$ then by variation of constants then $y(t) \to 0$ which implies $M_j(t) \to M_j(\infty)$.

        Thus, it follows that $\mathrsfs{R}(t) \to \mathrsfs{R}(\infty)$ and since $q(t) \to q_\infty$ then $\mathrsfs{R}^\Adv(t) \to\mathrsfs{R}^\Adv(\infty)   = \mathrsfs{R}(\infty)G(q_\infty)$ and $\widehat{\mathrsfs{B}}_{44}(t) \to \widehat{\mathrsfs{B}}_{44}(\infty) = 2 q_\infty^2 \mathrsfs{R}(\infty)$. We could also use a similar argument to show that $\widehat{\mathrsfs{B}}_{44}(t) \to \widehat{\mathrsfs{B}}_{44}(\infty)$ as $t\to\infty$ with
        \begin{equation}\label{eq:limit:B_44(inf)}
            \widehat{\mathrsfs{B}}_{44}(\infty) = L(\infty) + G(q_\infty )C_{\mathrsfs{B}}(q_\infty) \mathrsfs{R}(\infty)\quad \text{where}\quad C_{\mathrsfs{B}}(q) = \frac{\gamma}{d}\sum_{j=1}^d \frac{\lambda_j}{D_j(q)}.
        \end{equation}
        Finally, we have $q^2_\infty = \frac{\widehat{\mathrsfs{B}}_{44}(\infty)}{2\mathrsfs{R}(\infty)}$. Plugging \eqref{eq:limit:R(inf)} and \eqref{eq:limit:B_44(inf)} into this expression, we obtain
        \begin{equation}
            q^2_\infty = (1 - G(q_\infty) C_{\mathrsfs{R}}(q_\infty)) \frac{L_\infty(q_\infty)}{2 F_\infty(q_\infty)} + \frac{1}{2}G(q_\infty)C_\mathrsfs{B}(q_\infty),
        \end{equation}
        where we emphasize the dependence on $q_\infty$ by writing $L_\infty(q_\infty) \defas L(\infty)$ and $F_\infty(q_\infty) \defas F(\infty)$. We now prove uniqueness of the admissible fixed point. We have
        \[
            \frac{\dif}{\dif q}\big(q\tilde{\lambda}^{\operatorname{eff}}(q)\big)
            =
            \frac{\delta^2-\delta^2 \frac{2}{\pi }}{(1+\delta\sqrt{\frac{2}{\pi}}q)^2}
            =
            \frac{\delta^2(1-\frac{2}{\pi})}{\left(1+\delta\sqrt{\frac{2}{\pi}}q\right)^2}>0.
        \]
        Thus $q\mapsto (q\tilde{\lambda}^{\operatorname{eff}}(q))^2$ is strictly increasing. Multiplying the fixed point equation by $\tilde{\lambda}^{\operatorname{eff}}(q)^2$, it is sufficient to show that
        \[
            \Theta(q)
            \defas
            \tilde{\lambda}^{\operatorname{eff}}(q)^2\left[
            (1-G(q)C_{\mathrsfs R}(q))
            \frac{L_\infty(q)}{2F_\infty(q)}
            +
            \frac12G(q)C_{\mathrsfs B}(q)
            \right]
        \]
        is strictly decreasing on $ \mathcal I=\{q>0:G(q)C_{\mathrsfs R}(q)<1\}$ as the fixed point equation is equivalent to $(q\tilde{\lambda}^{\operatorname{eff}}(q))^2=\Theta(q)$. First, we have
        \[
            \frac{\dif}{\dif q} (\tilde{\lambda}^{\operatorname{eff}}(q))
            =-\frac{\delta \sqrt{\frac{2}{\pi}}\left( \frac{1}{q^2} + 2\delta \sqrt{\frac{2}{\pi}} \frac{1}{q} + \delta^2\right)}{\left(1+\delta \sqrt{\frac{2}{\pi}} q\right)^2} < 0.
        \]
        Next, define the functions
        \[
            A(r) = \sum_{j=1}^d \lambda_j w_j(r), \quad \text{with} \quad w_j(r) = \frac{
                \frac{\lambda_j(x_j^\star)^2}{(\lambda_j+r)^2}
                }{
                \sum_{k=1}^d
                \frac{\lambda_k(x_k^\star)^2}{(\lambda_k+r)^2},
                }
        \]
        such that $ \tilde{\lambda}^{\operatorname{eff}}(q)^2\frac{L_\infty(q)}{2F_\infty(q)} = A(\tilde{\lambda}^{\operatorname{eff}}(q))$. Differentiating we obtain
        \[
\frac{\dif}{\dif r} A(r) = 2 \sum_{1 \leq i < j \leq d}w_i(r) w_j(r)  \frac{(\lambda_i - \lambda_j) ^2}{(\lambda_i+r)(\lambda_j+r)}.
        \] 
        Apart from some degenerate cases such as $X^\star = c \omega_j$ with $\lambda_j = 0$ or $K=0$ then $\frac{\dif}{\dif r} A(r)  \geq 0$. From the chain rule then we obtain
        \[
            \frac{\dif}{\dif q} A(\tilde{\lambda}^{\operatorname{eff}}(q)) =   A'(\tilde{\lambda}^{\operatorname{eff}}(q))   \frac{\dif}{\dif q} (\tilde{\lambda}^{\operatorname{eff}}(q)) \leq 0.
        \] 
         Finally, $q\mapsto G(q)C_{\mathrsfs R}(q)$
        is strictly increasing, because for each $j$,
        \[
        \begin{gathered}
            \frac{\dif}{\dif q}\left(\frac{G(q)}{D_j(q)}\right)
            =
            \frac{\delta}
            {q^2D_j(q)^2}\bigg(2 \delta ^3 q^3 + \delta^2\lambda_j q^4 \sqrt{\frac{2}{\pi}} + 5 \delta^2 q^2 \sqrt{\frac{2}{\pi}} \\
            + 2\delta \lambda_j q^3 + \delta q \frac{8}{\pi} +\lambda_j q^2 \sqrt{\frac{2}{\pi}} + \sqrt{\frac{2}{\pi}}\bigg)
            >0.
        \end{gathered}
        \]
        Therefore $1-G(q)C_{\mathrsfs R}(q)$ is strictly decreasing and positive on $\mathcal I$ from which it follows that the first term in $\Theta(q)$ is nonincreasing on $\mathcal I$. For the second term, we have
        \[
        \begin{aligned}
            \frac{(\tilde{\lambda}^{\operatorname{eff}}(q))^2}{2}G(q) C_{\mathrsfs{B}} =  \frac{\gamma}{2d} \sum_{j=1}^d T_{\lambda_j, \delta}(q),
        \end{aligned}
        \]
        where we define
        \[
            T_{\lambda, \delta}(q)=\lambda \frac{\left(\delta^2 + \delta \sqrt{\frac{2}{\pi}} \frac{1}{q}\right)^2\left(1 + 2\delta \sqrt{\frac{2}{\pi}} q + \delta^2 q^2\right)}{\lambda\left(1+ \delta \sqrt{\frac{2}{\pi}} q \right)^3 + \left(\delta^2 + \delta \sqrt{\frac{2}{\pi}}\frac{1}{q}\right)\left(1 + \delta \sqrt{\frac{2}{\pi}}q \right)^2}.
        \]
        Our goal is to find a sufficient condition on $\delta$ for which $\frac{d}{dq}\sum_{j=1}^d T_{\lambda_j,\delta}(q)\le 0$. For a single eigenvalue, define $ r \defas \frac{\lambda}{\delta^2}.$ and note that
        \[
            \delta^2+ \delta \sqrt{\frac{2}{\pi}}\frac{1}{q}
                    =
            \delta^2\frac{\delta q+\sqrt{\frac{2}{\pi}}}{\delta q}.
        \]
        Substituting this and $\lambda=\delta^2r$ into $T_{\lambda,\delta}$, we obtain after a simple computation
        \[
        \begin{aligned}
            T_{\lambda,\delta}(q)
            &=\delta^3 
            \frac{
            r
         \left(\delta q+\sqrt{\frac{2}{\pi}}\right)^2
            \left(1+2\sqrt{\frac{2}{\pi}}\delta q+\delta^2q^2\right)
            }{q\left(1+\delta \sqrt{\frac{2}{\pi}}q\right)^2\left(\delta q r \left(1+\delta \sqrt{\frac{2}{\pi}}q\right)+ \delta q + \sqrt{\frac{2}{\pi}}\right)}.
        \end{aligned}
        \]

        Differentiating with respect to $q$, we obtain
        \[
        \begin{aligned}
            \frac{\dif }{\dif q}T_{\lambda,\delta}(q)
            &=
            \delta^3
            \frac{
            r\left(\delta q+\sqrt{\frac{2}{\pi}}\right)
            \left\{
            P_\delta(q)-rQ_\delta(q)
            \right\}
            }{q^2\left(1+\delta \sqrt{\frac{2}{\pi}}q\right)^3\left(\delta q r \left(1+\delta \sqrt{\frac{2}{\pi}}q\right)+ \delta q + \sqrt{\frac{2}{\pi}}\right)^2}.
        \end{aligned}
        \]
        Since all coefficients of $Q_\delta(q)$ are positive then we have $Q_\delta(q)>0$ for all $q >0$. For each eigenvalue, define  $r_j\defas \frac{\lambda_j}{\delta^2}$ then
        \[
        \begin{aligned}
           \frac{\gamma}{2d} \frac{\dif}{\dif q}\sum_{j=1}^d T_{\lambda_j,\delta}(q)
            &=
            \frac{\delta^3\gamma}{2d}
            \frac{
                \delta q+\sqrt{\frac{2}{\pi}}
            }{q^2\left(1+\delta \sqrt{\frac{2}{\pi}}q\right)^3} 
            \sum_{j=1}^d
            \frac{
            r_j\left\{
            P_\delta(q)-r_jQ_\delta(q)
            \right\}
            }{\left(\delta q r_j \left(1+\delta \sqrt{\frac{2}{\pi}}q\right)+ \delta q + \sqrt{\frac{2}{\pi}}\right)^2}.
        \end{aligned}
        \]
        The sign of the summed derivative is determined by the sign of the difference $P_\delta(q)A_d(q)-
        Q_\delta(q) B_d(q)$ where we defined the terms
        \[
        \begin{gathered}
            A_d(q)
            \defas
            \frac{1}{d}\sum_{j=1}^d
            \frac{
            r_j
            }{
            \left(
            \delta q+\sqrt{\frac{2}{\pi}}
            +
            r_j\delta q(1+\sqrt{\frac{2}{\pi}}\delta q)
            \right)^2
            },
            \\
            \text{and}\quad B_d(q)
            \defas
            \frac{1}{d}\sum_{j=1}^d
            \frac{
            r_j^2 
            }{
            \left(
            \delta q+\sqrt{\frac{2}{\pi}}
            +
            r_j\delta q(1+\sqrt{\frac{2}{\pi}}\delta q)
            \right)^2
            }.
        \end{gathered}
        \]
        Since $Q_\delta(q)>0$ and ignoring the degenerate case $K=0$, the condition $P_\delta(q)A_d(q)-Q_\delta(q)B_d(q)\leq 0$ is equivalent to
        \[
            \frac{B_d(q)}{A_d(q)}
            \ge
            \frac{P_\delta(q)}{Q_\delta(q)},
        \]
        whenever $P_\delta(q)>0$ and whenever $P_\delta(q)\le0$ then $P_\delta(q)A_d(q)-Q_\delta(q)B_d(q)<0$. Thus, a sufficient condition for the second term in $\Theta(q)$ to be nonincreasing is for all $q \in [q_- , q_+ ]$
        \[
            \frac{B_d(q)}{A_d(q)}
            \geq
            \max\left(
                \frac{P_\delta(q)}{Q_\delta(q)},0
                \right) \defas \psi_\delta(q).
        \]
        We now derive a sufficient conditions on $\delta$ for which this holds depending on the assumptions on the eigenvalues of $K$. From the definitions of $P_\delta$ and $ Q_\delta$, it is clear that they only depend on $q$ through $\delta q$. Define the quantity $ r_\star
        \defas
        \sup_{q>0}\psi_\delta(q)$. Then $r_\star$ is independent of the value of $\delta$. Indeed, we have
        \[
            r_\star
            \defas
            \sup_{q>0}\psi_\delta(q) =  \sup_{q>0} \max\left(
                \frac{P_\delta(q)}{Q_\delta(q)},0
                \right) =  \sup_{q>0} \max\left(
                    \frac{\widehat{P}(\delta q)}{\widehat{Q}(\delta q)},0
                    \right) = \sup_{s>0} \max\left(
                        \frac{\widehat{P}(s)}{\widehat{Q}(s)},0
                        \right)
        \]
        where $\widehat{P}(\delta q) = P_\delta(q)$, $\widehat{Q}(\delta q) = Q_\delta(q)$, $s = \delta q$ and since for $\delta>0$ then $q > 0 \iff s >0$.  
        Hence a simpler sufficient condition which guarantees monotonicity is for all $q \in [q_- , q_+ ]$
        \[
            \frac{B_d(q)}{A_d(q)}\ge r_\star.
        \]
        Numerically, we obtain $r_\star\approx 8.35\times 10^{-5}$. We rigorously show $0<r_\star < \infty$. Define $\widehat{\psi}$ such that $\widehat{\psi}(\delta q) = \psi_\delta(q)$ We have $\widehat{P}(0) = 0$. The leading coefficient of $\widehat{P}$ is greater than zero which implies that $\widehat{P}(s) \to \infty$ as $s\to \infty$ and thus there exists $s_0$ such that $\widehat{P}(s_0) >0$. It follows that $r_\star >0$. Now, as $s\to \infty$ we have $\widehat{P}(s) \sim \left(2 - \frac{6}{\pi}\right)s^4$ while $\widehat{Q}\sim \frac{2}{\pi}s^6$ from which it follows that 
        \[
        \frac{\widehat{P}(s)}{\widehat{Q}(s)} \sim \frac{2- \frac{6}{\pi}}{\frac{2}{\pi}} \frac{1}{s^2} \to 0 \quad \text{as}\quad s\to \infty.
        \]
        Since $\widehat{P}(0) <0$ then by continuity $\widehat{P}(s) < 0$ for sufficiently small $s>0$ which implies $\widehat{\psi}(s)= 0$ near $0$. Since $\widehat{\psi}(s) \to 0$ as $s\to \infty$ and $\widehat{\psi}$ is continuous on $(0, \infty)$, then $r_\star <\infty$.

        From the definitions of $B_d(q)$ and $A_d(q)$, it is easy to see that if $\lambda_j = 0$ then its contribution to $B_d(q)$ and $A_d(q)$ is zero. Let $J_+ \defas \{j: \lambda_j >0 \}$ and ignoring the degenerate case $K=0$ then we have
        \[
            \frac{B_d(q)}{A_d(q)} =\sum_{j \in J_+} r_j w_j,
        \]
        where we define the weights $\sum_{j \in J_+} w_j = 1$ as follows
        \[
            w_j \defas \frac{\frac{r_j}{\left(
                \delta q+\sqrt{\frac{2}{\pi}}
                +
                r_j\delta q\left(1+\sqrt{\frac{2}{\pi}}\delta q\right)
                \right)^2}}{\sum_{k \in J_+}^d
                \frac{
                r_k 
                }{
                \left(
                \delta q+\sqrt{\frac{2}{\pi}}
                +
                r_k\delta q\left(1+\sqrt{\frac{2}{\pi}}\delta q\right)
                \right)^2
                }}.
        \]
        Consequently, $B_d(q)/ A_d(q)$ is a weighted average which implies
        \[
            \frac{\lambda_{\operatorname{min}}^+}{\delta^2}= \underset{j \in J_+} {\operatorname{min}}r_j \leq \frac{B_d(q)}{A_d(q)},
        \]
        where $\lambda_{\operatorname{min}}^+$ is the smallest non-zero eigenvalue. Thus, a sufficient condition is 
        \[
            \delta \leq \sqrt{\frac{\lambda_{\operatorname{min}}^+}{ r_\star}}.
        \]
        The condition for isotropic covariance then follows immediately. 
    \end{proof}
This leads us to the proof of Proposition~\ref{prop:fixed:point:gamma:to:0:main}.
    \begin{proof}[Proof of Proposition~\ref{prop:fixed:point:gamma:to:0:main}]
        We omit the proof of $\delta_- < \delta_+$ as it follows from two simple applications of Cauchy--Schwarz. Since $K$ is symmetric, we have $\R = \operatorname{Range}(K) \oplus \ker(K)$. Hence, we can decompose $X = X_+ +X^0$ for $X_+ \in \operatorname{Range}(K)$ and $X^0 \in \ker(K)$. Similarly, we can decompose $X^\star = X^\star_+ + X^{\star, 0}$ and $X^{\star, \Adv} =X^{\star, \Adv}_+ + X^{\star, \Adv, 0}$. 
        
        We now show by contradiction that $ X^{\star, \Adv, 0}=0$. Assume that $ X^{\star, \Adv, 0} \not=0$. From the above decompositions, we have that$(X-X^\star)^\top K(X-X^\star) = (X_+ - X^\star_+)^\top K(X_+ - X^\star_+)$ and thus $\mathcal{R}(X ) = \mathcal{R}(X_+)$. 
        
        Since $ \operatorname{Range}(K) \perp \ker(K)$ then $\|X\|^2 = \|X_+\|^2 + \|X^0\|^2 \geq \|X_+\|^2$. Therefore, we have
        \[
            \mathcal{R}^{\Adv}(X) - \mathcal{R}^{\Adv}(X_+) = \delta \sqrt{\frac{2}{\pi}} \sqrt{2\mathcal{R}(X_+)}\left(\|X\| - \|X_+\|\right) + \frac{\delta^2}{2}\|X^0\|^2.
        \]
        Since $X^{\star, \Adv, 0} \not = 0$ then $\|X^{\star, \Adv, 0}\|>0$ and $\|X^{\star, \Adv}\|^2 > \|X^{\star, \Adv}_+\|^2$ which implies $ \mathcal{R}^{\Adv}(X^{\star, \Adv} ) - \mathcal{R}^{\Adv}(X^{\star, \Adv}_+) > 0$. Thus, we have reached a contradiction as $X^{\star, \Adv}$ is a minimizer of $\mathcal{R}^{\Adv}$ from which it follows that $X^{\star, \Adv} \in \operatorname{Range}(K)$.

        On $\operatorname{Range}(K)$, the operator $\left. K \right|_{\operatorname{Range}(K)}$ is positive definite. The preceding argument shows that minimizing $\mathcal R^{\Adv}$ over $\R^d$ is equivalent to minimizing it over $\operatorname{Range}(K)$:
        \[
            \argmin_{X\in\mathbb R^d}\mathcal R^{\Adv}(X) = \argmin_{X\in\operatorname{Range}(K)}\mathcal R^{\Adv}(X)
        \]
        Plugging the definition of $q_\star$ into $\tilde{\lambda}_{\star}$ in \eqref{eq:fixed:point:gamma:0} recovers Equation 4 in \citep{xing2021same_minimizer} restricted to $\operatorname{Range}(K)$. Thus, we apply Proposition 1 in \citep{xing2021same_minimizer} to the positive definite covariance matrix $\left. K \right|_{\operatorname{Range}(K)}$ and obtain the minimizer
        \[
            X^{\star,\Adv} =\sum_{j\in J_+}
            \frac{\lambda_j}{\lambda_j+\tilde{\lambda}_{\star}}
            x_j^\star\omega_j,
        \]
        from which we obtain the limiting values as $\gamma \to 0$ in \eqref{eq:result:q:limits}. It remains to prove that $q_{\gamma}(\infty) \to q_{\star}$. We start with the regime $\delta_- < \delta < \delta_+$. From Proposition~\ref{prop:fixed:point} we have
        \[
        D_j(q) = \lambda_j(1+ \delta \sqrt{\frac{2}{\pi}}q)  + \delta^2 + \delta \sqrt{\frac{2}{\pi}}\frac{1}{q} \geq \delta^2 > 0.
        \]
        For all $q \in [a,b] \subset (0, \infty)$, this implies that $C_\mathrsfs{R}(q)  = \mathcal{O}(\gamma)$ and $C_{\mathrsfs{B}}(q) =\mathcal{O}(\gamma)$. The functions $G(q)$, $F_{\infty}(q)$, $L_{\infty}(q)$ are continuous on $[a,b]$ and do not depend on $\gamma$. Moreover, $F_{\infty}(q) > 0$ on $[a,b]$ from which it follows that
        \[
            \sup_{q \in [a,b]} \left| (1 - G(q) C_{\mathrsfs{R}}(q)) \frac{L_\infty(q)}{2 F_\infty(q)} + \frac{1}{2}G(q)C_\mathrsfs{B}(q) - \frac{L_{\infty}(q)}{2F_{\infty}(q)}\right| \underset{\gamma \to 0}{\longrightarrow} 0.
        \]
        Hence, for a subsequence $\gamma_n \downarrow 0$ such that $q_{\gamma_n}\to q_0\in(0,\infty)$, then 
        \begin{equation}\label{eq:fixed:point:gamma:0}
        q_0^2=\frac{L_{\infty}(q_0)}{2F_{\infty}(q_0)}.
        \end{equation}
         Thus every finite positive accumulation point of $q_\gamma$ must be a solution of the fixed point equation \eqref{eq:fixed:point:gamma:0}.  In the regime $\delta_- < \delta < \delta_+$, the solution to \eqref{eq:fixed:point:gamma:0} is equivalent to the solution to Equation 4 in \citep{xing2021same_minimizer}. Moreover, \citep{xing2021same_minimizer} show $\tilde{\lambda}_{\star}$ is the unique stationnary fixed point satisfying Equation 4. They also prove that Equation 4 has no solution in $(0, \infty)$ whenever $\delta \leq  \delta_-$ or $\delta_+ \leq \delta$. Therefore, outside the middle regime, $q_\gamma(\infty)$ cannot remain in any $\gamma$-independent compact subset of $(0,\infty)$.

         With this in mind, we prove the limits in the boundary regimes. First suppose, towards contradiction, that along some sequence
         $\gamma_n\downarrow0$, $ q_{\gamma_n}(\infty)\to0$. For $0<q\le1$, $D_j(q)
         \ge
         \delta \sqrt{\frac{2}{\pi}} \frac{1}{q}$. Therefore
         \[
             C_{\mathrsfs{R}}(q)
             \le
             \sqrt{\frac{\pi}{2}}\frac{\gamma q}{2d\delta }
             \sum_{j=1}^d\lambda_j^2  \quad \text{and}\quad C_{\mathrsfs{B}}(q)
             \le
             \sqrt{\frac{\pi}{2}}\frac{\gamma q}{d\delta }
             \sum_{j=1}^d\lambda_j.
         \]
         Since $G(q)\to1$ as $q\downarrow0$, this implies $G(q_{\gamma_n})C_{\mathrsfs{R}}(q_{\gamma_n})\to0$ and $G(q_{\gamma_n})C_{\mathrsfs{B}}(q_{\gamma_n})\to0$. Using the fixed point equation and the fact that $C_{\mathrsfs{B}}(q)\ge0$, we get
         \[
             (q_{\gamma_n}(\infty))^2
             \ge
             \left(1-G(q_{\gamma_n}(\infty))C_{\mathrsfs{R}}(q_{\gamma_n}(\infty))\right)
             \frac{L_{\infty}(q_{\gamma_n}(\infty) )}{2 F_{\infty}(q_{\gamma_n}(\infty))}.
         \]
         Dividing by $ (q_{\gamma_n}(\infty))^2$ gives
         \[
             1
             \ge
             \left(1-G(q_{\gamma_n}(\infty))C_{\mathrsfs{R}}(q_{\gamma_n}(\infty))\right)
             \frac{L_{\infty}(q_{\gamma_n}(\infty) )}{2 F_{\infty}(q_{\gamma_n}(\infty))} \frac{1}{q^2_{\gamma_n}(\infty)}.
         \]
         Since $\tilde{\lambda}(q) = \sim \delta \sqrt{\frac{2}{\pi}}\frac{1}{q}$ as $q\downarrow0$ then
         \[ F_\infty(q) \to \frac12\sum_{j\in J_+}\lambda_j(x_j^\star)^2\quad \text{and}\quad L_\infty(q) \sim \frac{\pi}{2\delta^2}q^2 \sum_{j\in J_+}\lambda_j^2(x_j^\star)^2,
         \]
         from which it follows that
         \[
            \frac{L_{\infty}(q_{\gamma_n}(\infty) )}{2 F_{\infty}(q_{\gamma_n}(\infty))} \frac{1}{q^2_{\gamma_n}(\infty)}
             \longrightarrow
             \left(\frac{\delta_+}{\delta}\right)^2.
         \]
         Thus, we obtain $1\ge \left(\frac{\delta_+}{\delta}\right)^2$ which implies whenever $\delta<\delta_+$ and then $q_\gamma$ cannot converge to $0$. A similar argument shows that along some subsequence
         $\gamma_n\downarrow 0$, when $\delta > \delta_+$, $q_{\gamma_n}$ cannot converge to $\infty$.
         
         If $\delta\le\delta_-$, then $\delta<\delta_+$, so no subsequence of $q_\gamma(\infty)$ can converge to $0$ nor converge to any point in $(0,\infty)$ as the fixed point equation \eqref{eq:fixed:point:gamma:0} does not admit a solution in this regime. This implies $ q_\gamma\to\infty$ when $\gamma \to 0$ and thus $\tilde{\lambda}(q_\gamma)\to0$. If $\delta_-<\delta<\delta_+$, then $q_\gamma$ cannot converge to $0$ nor $\infty$. Hence every sequence $\gamma_n\downarrow0$ has
         a subsequence along which $q_{\gamma_n}\to q_0\in(0,\infty)$ which solves the fixed point equation \eqref{eq:fixed:point:gamma:0}. Since $q_\star$ is the unique solution to \eqref{eq:fixed:point:gamma:0}, then every subsequential limit equals $q_\star$ from which it follows that $ q_\gamma\to q_\star$ and by continuity then $ \tilde{\lambda}^{\operatorname{eff}}(q_\gamma)\to\tilde{\lambda}^{\operatorname{eff}}(q_\star)$.
         
         Finally, if $\delta\ge\delta_+$, then $\delta>\delta_-$, so no subsequence of $q_\gamma(\infty)$ can converge to $\infty$ nor to a point in $(0, \infty)$ as the fixed point equation \eqref{eq:fixed:point:gamma:0} does not admit a solution in this regime. This implies $ q_\gamma\to0$ when $\gamma \to 0$ and thus $\tilde{\lambda}(q_\gamma)\to \infty$. 
         
        The limits for the isotropic setting follow from a simple computation.
    \end{proof}
\subsubsection{Analysis of Learning Dynamics when $d\to \infty$}\label{app:power:law}
In this section, we provide the proofs for the results in Section~\ref{sect:plaw}. For each finite $d\geq 1$, we reexpress \eqref{eq:def:R(t):det:equiv:lr:reform} and \eqref{eq:def:l2norm2:det:equiv:lr:reform} as follows
\begin{equation}
    \begin{aligned}
        \mathrsfs{R}(t)&= F(\tilde{\Gamma}(t)) 
+ \frac{ \gamma^2}{d} \int_0^t\mathrsfs{R}^{{\operatorname{adv}}}(s) \tr \big ( K^2 e^{-2K(\tilde{\Gamma}(t) - \tilde{\Gamma}(s)) -2(\Lambda(t) -\Lambda(s))} \big ) \, \dif s.
    \end{aligned}
\end{equation}
and 
\begin{equation}
    \begin{aligned}
        \widehat{\mathrsfs{B}}_{44}(t) 
        &= L(\tilde{\Gamma}(t))
+
\frac{2\gamma^2}{d} \int_0^t   \mathrsfs{R}^{{\operatorname{adv}}}(s) \tr \big ( K e^{-2K(\tilde{\Gamma}(t) - \tilde{\Gamma}(s)) -2(\Lambda(t) -\Lambda(s))} \big ) \, \dif s
    \end{aligned}
    \end{equation}
From Lemma~\ref{lem:adv:risk:L2_ls}, we have $\mathrsfs{R}^\Adv(t) \asymp \mathrsfs{R}(t) +\frac{\delta^2}{2}\widehat{\mathrsfs{B}}_{44}(t)$. Hence, we will study the dynamics of the latter term which are simpler. We define its analogous term
\begin{equation}
    \begin{aligned}
        F^\Adv(u) &
        \defas F(u) + \frac{\delta^2}{2} L(u)
        =
        \frac{1}{2}\sum_{j=1}^d (\lambda_j + \delta^2 )\bigg(e^{-\lambda_j u - \ell(u) } \left( x_{0,j} - x^\star_j\right) 
        \\
        &
        + x^\star_j \left(\frac{\delta^2}{\lambda_j + \delta^2 } - \left(\int_0^{u} \ell'(s) e^{-\lambda_j (u-s) - (\ell(u)-\ell(s))} \dif s\right)\right) \bigg)^2
        \\
        &
        + \frac{1}{2} \sum_{j=1}^d \frac{\delta^2 \lambda_j}{ \lambda_j + \delta^2} (x^\star_j)^2.
    \end{aligned}
\end{equation}  
We also define the kernel term
\begin{equation}
    \mathscr{K}(x,y) \defas \frac{ \gamma^2}{d}\sum_{j=1}^d (\lambda_j^2 + \delta^2 \lambda_j) e^{-2 (\lambda_j x + y)}.
\end{equation}
It then follows that 
\begin{equation}
    \begin{aligned}
    \mathrsfs{R}(t) +\frac{\delta^2}{2}\widehat{\mathrsfs{B}}_{44}(t) 
    &
    =  F^\Adv(\tilde{\Gamma}(t)) +\int_0^t \mathscr{K}(\tilde{\Gamma}(t) - \tilde{\Gamma}(s), \Lambda(t) - \Lambda(s))\mathrsfs{R}^\Adv(s) \dif s.
    \end{aligned}
\end{equation}
        \begin{proof}[Proof of Proposition~\ref{prop:isotropic_Rasymp}]
            It is immediate that $ F^\Adv(\tilde{\Gamma}(t)) \leq \mathrsfs{R}(t) +\frac{\delta^2}{2}\widehat{\mathrsfs{B}}_{44}(t) \leq \mathrsfs{R}^\Adv(t)$. From Lemma~\ref{lem:R_up_low_bound}, we have
            \begin{align*}
                \mathrsfs{R}(t) +\frac{\delta^2}{2}\widehat{\mathrsfs{B}}_{44}(t) \leq  F^\Adv(\tilde{\Gamma}(t)) +2\gamma^2(1+\delta^2)\int_0^t e^{-2\gamma(1+\delta^2)(1-\gamma)(t-s)} F^\Adv(\tilde{\Gamma}(s)) \dif s.
            \end{align*}
            We require an upper and lower bound on $F^\Adv(\tilde{\Gamma}(t))$. From the definition of $F^\Adv(\tilde{\Gamma}(s))$, it follows immediately that 
            \begin{equation}\label{eq:lower:bound:Fadv}
                C_1 = \frac{1}{2} \sum_{j=1}^d \frac{\delta^2 \lambda_j}{ \lambda_j + \delta^2} (x^\star_j)^2 \leq F^\Adv(\tilde{\Gamma}(t)).
            \end{equation}
            Now, for the upper bound, noting that 
            \begin{equation}
                \left|\frac{\delta^2}{\lambda_j + \delta^2 } - \int_0^{u} \ell'(s) e^{-\lambda_j (u-s) - (\ell(u)-\ell(s))} \dif s \right|\leq 1
            \end{equation}
            then 
            \[
            \begin{gathered}
                \bigg(e^{-\lambda_j u - \ell(u) } \left( x_{0,j} - x^\star_j\right) 
+ x^\star_j \left(\frac{\delta^2}{\lambda_j + \delta^2 } - \left(\int_0^{u} \ell'(s) e^{-\lambda_j (u-s) - (\ell(u)-\ell(s))} \dif s\right)\right) \bigg)^2
\\
                 \leq 2e^{-2(\lambda_j u +\ell(u) )} \left( x_{0,j} - x^\star_j\right) ^2 + 2(x^\star_j)^2.
            \end{gathered}
            \]
            Thus, we obtain
            \begin{equation}\label{eq:upper:bound:Fadv}
            \begin{aligned}
                F^\Adv(u) &
                \leq 
                \sum_{j=1}^d (\lambda_j + \delta^2 )\left(e^{-2(\lambda_j u + \ell(u) )} \left( x_{0,j} - x^\star_j\right) ^2 + (x^\star_j)^2\right)
                + \frac{1}{2} \sum_{j=1}^d \frac{\delta^2 \lambda_j}{ \lambda_j + \delta^2} (x^\star_j)^2
                \\
                &
                \leq 
                2\sum_{j=1}^d (\lambda_j + \delta^2 ) x_{0,j}^2 +\sum_{j=1}^d \left(3(\lambda_j + \delta^2 ) +  \frac{\delta^2 \lambda_j}{ 2(\lambda_j + \delta^2)}\right)(x^\star_j)^2.
            \end{aligned}
        \end{equation}
            Setting $\lambda_j=1$ and summing over $j$, we obtain $ C_1 \leq F^\Adv(\tilde{\Gamma}(t)) \leq  C_{F}$ after noting that 
            \[
                \|X^\star\|^2 = \sum_{j=1}^d (x_j^\star)^2 = \int_0^1 y^{-\alpha \beta} \dif y + \mathcal{O}(d^{-1}) = 1+ \mathcal{O}(d^{-1}),
            \]
          and $\|X_0\|\leq C$ by Assumption~\ref{assumption:scaling}. It thus follows that 
            \begin{equation}
                4\gamma^2(1+\delta^2)\int_0^t e^{-2\gamma(1+\delta^2)(1-\gamma)(t-s)} F^\Adv(\tilde{\Gamma}(s)) \dif s \leq C_F \int_0^\infty e^{-2\gamma(1+\delta^2)(1-\gamma)(t-s)} \dif s,
            \end{equation}
            from which the upper bound follows.
        \end{proof}
        We omit the proof of Proposition~\ref{prop:powerlaw_Rasymp} as it follows a similar argument to Proposition~\ref{prop:isotropic_Rasymp}. We now present the proofs of Corollaries~\ref{cor:powerlaw_R_B_44_bounded:isotropic}~and~\ref{cor:powerlaw_R_B_44_bounded:plaw}.

        \begin{proof}[Proof of Corollary~\ref{cor:powerlaw_R_B_44_bounded:isotropic}]
                From Corollaries~\ref{prop:q(t):bounded}~and~\ref{prop:isotropic_Rasymp}, since $\mathrsfs{R}(t)$ can be expressed in terms of $\mathrsfs{R}^\Adv(t)$ and $q(t)$ as follows
                \[
                    \mathrsfs{R}(t) = \frac{\mathrsfs{R}^\Adv(t)}{1+ 2\delta \sqrt{\frac{2}{\pi}} q(t) + \delta^2 (q(t))^2},
                \]
                then there exists constants $C_1, C_2>0$ such that $C_1 \leq \mathrsfs{R}(t) \leq C_2$. Similarly, since $\widehat{\mathrsfs{B}}_{44}(t) = 2\mathrsfs{R}(t) (q(t))^2$ then there exists constants $C_3, C_4>0$ such that $C_3 \leq \widehat{\mathrsfs{B}}_{44}(t) \leq C_4$. The lower bounds $F(\tilde{\Gamma}(t)) \leq \mathrsfs{R}(t)$ and $L(\tilde{\Gamma}(t)) \leq \widehat{\mathrsfs{B}}_{44}(t)$ are immediate from the definitions of $\mathrsfs{R}(t)$ and $\widehat{\mathrsfs{B}}_{44}(t)$. From Lemma~\ref{lem:adv:risk:L2_ls}, we have the upper bound
                \begin{equation}
                    \begin{aligned}
                        \mathrsfs{R}(t)&\leq F(\tilde{\Gamma}(t)) 
                        + 2\gamma^2 \int_0^t  e^{-2\gamma(1 + \delta^2)(t-s)}\left(\mathrsfs{R}(s) + \frac{\delta^2}{2}\widehat{\mathrsfs{B}}_{44}(s)\right)\, \dif s.
                    \end{aligned}
                \end{equation}
                It then follows from Lemma~\ref{lem:R_up_low_bound} that 
                \begin{equation}
                    \begin{aligned}
                        \mathrsfs{R}(t)&\leq  F(\tilde{\Gamma}(t)) 
                        + 2\gamma^2\int_0^t   e^{-2\gamma(1 + \delta^2)(t-s)} F^\Adv(\tilde{\Gamma}(s)) \dif s 
                        \\
                        &
                        +  4\gamma^4(1+\delta^2)\int_0^t  e^{-2\gamma(1+\delta^2)(2-\gamma)(t-s)} F^\Adv(\tilde{\Gamma}(s))\dif s .
                    \end{aligned}
                \end{equation}
                The bound for $\widehat{\mathrsfs{B}}_{44}(t)$ follows a similar argument.
            \end{proof}

            \begin{proof}[Proof of Corollary~\ref{cor:powerlaw_R_B_44_bounded:plaw}]
                                The upper and lower bounds by constants follow a similar argument as in the proof of Corollary~\ref{cor:powerlaw_R_B_44_bounded:isotropic} using Corollaries~\ref{prop:q(t):bounded}~and~\ref{prop:powerlaw_Rasymp}. The lower bound $F(\tilde{\Gamma}(t)) \leq \mathrsfs{R}(t)$ is immediate from the definition of $\mathrsfs{R}(t)$ and similarly for $\widehat{B}_{44}(t)$. From Lemma~\ref{lem:adv:risk:L2_ls}, we have the upper bound
                                \begin{equation}
                                    \begin{aligned}
                                        \mathrsfs{R}(t)&\leq F(\tilde{\Gamma}(t)) 
                                        +  \int_0^t \overline{\mathscr{K}}_{\mathrsfs{R}}(t-s ) \left(\mathrsfs{R}(s) + \frac{\delta^2}{2}\widehat{\mathrsfs{B}}_{44}(s)\right)\, \dif s.
                                    \end{aligned}
                                \end{equation}
                                It then follows from Lemma~\ref{lem:R_up_low_bound} that 
                                \begin{equation}
                                    \begin{aligned}
                                        \mathrsfs{R}(t)&\leq  F(\tilde{\Gamma}(t)) 
                                        + \int_0^t \overline{\mathscr{K}}_{\mathrsfs{R}}(t-s ) F^\Adv(\tilde{\Gamma}(s)) \dif s +C \int_0^t  (\overline{\mathscr{K}}_{\mathrsfs{R}} *\overline{\mathscr{K}})(t-s) F^\Adv(\tilde{\Gamma}(s)) \dif s .
                                    \end{aligned}
                                \end{equation}
                    where we switched the order of integration and used a change of variables to obtain the bound and where 
                    \[
                    C = \left(\frac{\overline{\mathscr{K}}(0)}{\overline{\mathscr{K}}(T)(2\varepsilon +1)}+2\right)\frac{1}{1-2\|\overline{\mathscr{L}}\|_1(1+\varepsilon)
                    }.
                    \]
                                The bounds for $\widehat{\mathrsfs{B}}_{44}(t)$ follow a similar argument.
                                        \end{proof}

                        We close off this section of the appendix with the proof of Proposition~\ref{prop:fixed:point:plaw}.
                        \begin{proof}[Proof of Proposition~\ref{prop:fixed:point:plaw}]
                            Under the power law assumption with $\alpha \geq0$, we have $\|K\|_{\opt} =1$ and as 
                            \[
                                \frac{1}{d}\Tr(K^2)  = \frac{1}{d} \sum_{j=1}^d \left( \frac{j}{d}\right)^{2\alpha} \to \int_0^1 y^{2\alpha} \dif y = \frac{1}{2\alpha + 1} > 0.
                            \]
                            Thus, by Proposition~\ref{prop:q(t):bounded}, we have that there exists $0< q_- < q_+ < \infty$ independent on $d$ such that $q_- \leq q_d(t) \leq q_+$. Taking the limit as $t\to \infty$, it thus follows that $q_- \leq q_{d,\infty} \leq q_+$ which implies $\{ q_{d,\infty} \}_{d\geq 1}$ is precompact in $(0,\infty)$. For each $d\geq 1$, we define 
                            \[
                                H_d(q) = q^2- \left((1 - G(q) C_{\mathrsfs{R}, d}(q)) \frac{L_{\infty, d}(q)}{2 F_{\infty,d}(q)} + \frac{1}{2}G(q)C_{\mathrsfs{B}, d}(q)\right),
                            \]
                            where $G$, $C_{\mathrsfs{R}, d}$, $L_{\infty, d}$, $F_{\infty,d}$ and $C_{\mathrsfs{B}, d}$ are defined in Proposition~\ref{prop:fixed:point} and $H_d(q_{d,\infty}) = 0$. We add the subscript $d$ to emphasize their dependence on dimension. Analogously, we define 
                            \[
                                H(q) = q^2- \left((1 - G(q) C_{\mathrsfs{R}}(q)) \frac{L_{\infty}(q)}{2 F_{\infty}(q)} + \frac{1}{2}G(q)C_{\mathrsfs{B}}(q)\right).
                            \]
                            We will show that
                            \[
                                \sup_{q \in [q_-, q_+]}\left|H_d(q) - H(q)\right| \underset{d \to \infty}{\longrightarrow} 0.
                            \]
                            First, the integrands defining $C_{\mathrsfs{R}}$ and $C_{\mathrsfs{B}}$ are continuous on $[q_-, q_+] \times [0,1]$ hence uniformly continuous. It is then easy to see that the Riemann sums $C_{\mathrsfs{R},d}$ and $C_{\mathrsfs{B},d}$ respectively converge uniformly on $[q_-, q_+]$ towards $C_{\mathrsfs{R}}$ and $C_{\mathrsfs{B}}$ as $d\to \infty$. We have
                            \[
                            \begin{gathered}
                            F_{\infty, d}(q) = \frac{1}{2d} \sum_{j=1}^d \left(\frac{j}{d}\right)^{\alpha(1-\beta)}\left(\frac{\tilde{\lambda}^{\operatorname{eff}}(q) }{ \left(\frac{j}{d}\right)^{\alpha} + \tilde{\lambda}^{\operatorname{eff}}(q)}\right)^{2},
                            \\
                            \text{and}\quad L_{\infty, d}(q) = \frac{1}{d}\sum_{j=1}^d \left(\frac{j}{d}\right)^{-\alpha \beta} \left(\frac{\left(\frac{j}{d}\right)^{\alpha } }{\left(\frac{j}{d}\right)^{\alpha } + \tilde{\lambda}^{\operatorname{eff}}(q)}\right)^2.
                            \end{gathered}
                            \]
                            Fix $\varepsilon \in (0, 1]$. On $[q_-, q_+]\times [\varepsilon, 1]$, the integrands of $F_{\infty}$ and $L_{\infty}$ are uniformly continuous. Hence, with the same argument as above, the Riemann sums converge uniformly in $q \in [q_-, q_+]$ for $[\varepsilon, 1]$. It remains to bound the contribution from $(0, \varepsilon]$. Uniformly over $q\in [q_-, q_+]$, we have
                        \[
                        0\le \frac12 y^{\alpha(1-\beta)}
                        \left(
                        \frac{\widetilde\lambda(q)}
                        {y^\alpha+\widetilde\lambda(q)}
                        \right)^2
                        \le
                        \frac12 y^{\alpha(1-\beta)}
                        \le
                        \frac12 y^{-\alpha \beta},
                        \qquad 0<y\le1.
                        \]
                        Hence
                        \[
                        \sup_{q\in [q_-, q_+]}\int_0^\varepsilon  \frac12 y^{\alpha(1-\beta)}
                        \left(
                        \frac{\widetilde\lambda(q)}
                        {y^\alpha+\widetilde\lambda(q)}
                        \right)^2\,dy
                        \le
                        \frac12\int_0^\varepsilon y^{-\alpha\beta}\,dy
                        =
                        \frac{1}{2(1-\alpha\beta)}\varepsilon^{1-\alpha\beta}.
                        \]
                        Similarly, we have
                        \[
                        \sup_{q\in [q_-, q_+]}
                        \frac{1}{2d}\sum_{j/d\le\varepsilon}
                        \left(\frac{j}{d}\right)^{\alpha(1-\beta)}\left(\frac{\tilde{\lambda}^{\operatorname{eff}}(q) }{ \left(\frac{j}{d}\right)^{\alpha} + \tilde{\lambda}^{\operatorname{eff}}(q)}\right)^{2}
                        \le
                        \frac{1}{2d}
                        \sum_{j/d\le\varepsilon}
                        \left(\frac jd\right)^{-\alpha\beta}.
                        \]
                        Since for sufficiently large $d$
                        \[
                        \frac{1}{d}\sum_{j/d\le\varepsilon}
                        \left(\frac jd\right)^{-\alpha\beta}
                        =
                        d^{-(1-\alpha\beta)}\sum_{j\le d\varepsilon}j^{-\alpha\beta}
                        \le
                        d^{-(1-\alpha\beta)}
                        +
                        \frac{\varepsilon^{1-\alpha\beta}}{1-\alpha\beta},
                        \]
                        then it follows that
                        \[
                        \limsup_{d\to\infty}
                        \sup_{q\in [q_-, q_+]}
                        \frac{1}{2d}\sum_{j/d\le\varepsilon}
                        \left(\frac{j}{d}\right)^{\alpha(1-\beta)}\left(\frac{\tilde{\lambda}^{\operatorname{eff}}(q) }{ \left(\frac{j}{d}\right)^{\alpha} + \tilde{\lambda}^{\operatorname{eff}}(q)}\right)^{2}
                        \le
                        \frac{1}{2(1-\alpha\beta)}\varepsilon^{1-\alpha\beta}.
                        \]
                        Therefore, taking $d\to \infty$, the contribution of $F_{\infty, d}- F_{\infty}$ from $(0,\varepsilon]$ is bounded by $\frac{1}{1-\alpha\beta}\varepsilon^{1-\alpha\beta}$. Since $\alpha\beta<1$, taking $\varepsilon \to 0$, we obtain uniform convergence of $ F_{\infty, d}$ to $ F_{\infty}$ as $d\to \infty$. A similar argument shows that  $ L_{\infty, d}$ converges uniformly to $ L_{\infty}$ on $[q_-, q_+]$ as $d\to \infty$.
                        
                        Now, $F_{\infty}(q) > 0$ for all $q\in [q_-, q_+]$ and is continuous on this compact interval so its minimum is greater than $0$. Since $ F_{\infty, d}$ to $ F_{\infty}$ converges uniformly on $[q_-, q_+]$ as $d\to \infty$, for sufficiently large $d$ it follows that $ F_{\infty, d} >0$ for all $[q_-, q_+]$. Hence, the ratio $\frac{L_{\infty,d}}{F_{\infty,d}}$ converges uniformly on $[q_-, q_+]$ to $\frac{L_{\infty}}{F_{\infty}}$ as $d\to \infty$. Since the function $G(q)$ is continuous on $[q_-, q_+]$, combining all these results we conclude that
                        \begin{equation}\label{eq:unif:conv:H}
                            \sup_{q \in [q_-, q_+]}\left|H_d(q) - H(q)\right| \underset{d \to \infty}{\longrightarrow} 0.
                        \end{equation}
                        Now, take any subsequence of $\{q_{d, \infty}\}_{d\geq 1}$. Then by compactness, let $q_{d_k,\infty}$ be a further subsequence which converges $q_{d_k, \infty} \to q_{\star,\infty} \in [q_-, q_+]$. Since $H_{d_k} (q_{d_k, \infty})= 0$, $H$ is continuous on $[q_-, q_+]$ and by \eqref{eq:unif:conv:H}, it thus follows
                        \begin{equation}
                            |H(q_{\star, \infty}) | \leq |H(q_{\star, \infty}) - H(q_{d_k, \infty})| + |H(q_{d_k, \infty})  - H_{d_k}(q_{d_k, \infty})| \underset{k \to \infty}{\longrightarrow} 0.
                        \end{equation}
                        Thus, all limit points of $\{q_{d,\infty}\}_{d\geq1}$ must satisfy the fixed point equation at $d\to \infty$. Finally, if $q_{\infty}$ is the unique admissible solution to the fixed point equation $H(q)$ then every convergent subsequence converges to $q_\infty$ from which it follows that $q_{d,\infty} \to  q_{\infty}$. 
                        
                        We now derive sufficient conditions for which the fixed point satisfying $H(q)$ and the admissibility condition is unique. Recall we previously defined $J_+ \defas \{j: \lambda_j >0 \}$. The condition for $\alpha =0$ is shown in Proposition~\ref{prop:fixed:point}. We now derive sufficient conditions for $\alpha >0$. Noting that
                            \[
                                \left(
                                \delta q+\sqrt{\frac{2}{\pi}}
                                +
                                r_j\delta q(1+\sqrt{\frac{2}{\pi}}\delta q)
                                \right)^2 =  \left( \delta q+\sqrt{\frac{2}{\pi}}\right)^2\left(1 + \frac{r_j\delta q\left(1+\sqrt{\frac{2}{\pi}}\delta q\right)}{ \delta q+\sqrt{\frac{2}{\pi}}}\right)^2,
                            \]
                            then we obtain from a simple computation
                            \[
                                \frac{B_d(q)}{A_d(q)} = \frac{\sum_{j \in J_+}
                                \frac{
                                1
                                }{\left(\frac{1}{r_j} + \frac{\delta q\left(1+\sqrt{\frac{2}{\pi}}\delta q\right)}{ \delta q+\sqrt{\frac{2}{\pi}}}\right)^2}}{\sum_{k \in J_+}
                                \frac{
                                \frac{1}{r_k}
                                }{\left(\frac{1}{r_k} + \frac{\delta q\left(1+\sqrt{\frac{2}{\pi}}\delta q\right)}{ \delta q+\sqrt{\frac{2}{\pi}}}\right)^2
                                }}.
                            \]
                            We can express this as follows 
                            \[
                                \frac{B_d(q)}{A_d(q)} = \frac{1}{\sum_{j \in J_+} p_j(q) \frac{1}{r_j}}, \quad \text{with} \quad
                                p_j(q) \defas \frac{\left(\frac{1}{r_j} + \frac{\delta q\left(1+\sqrt{\frac{2}{\pi}}\delta q\right)}{ \delta q+\sqrt{\frac{2}{\pi}}}\right)^{-2}}{\sum_{k \in J_+}
                                \left(\frac{1}{r_k} + \frac{\delta q\left(1+\sqrt{\frac{2}{\pi}}\delta q\right)}{ \delta q+\sqrt{\frac{2}{\pi}}}\right)^{-2}
                                },
                            \]
                            such that $\sum_{j\in J_+} p_j(q) = 1$. The map $z \mapsto \left(z + \frac{\delta q\left(1+\sqrt{\frac{2}{\pi}}\delta q\right)}{ \delta q+\sqrt{\frac{2}{\pi}}}\right)^{-2}$ is decreasing in $z$ so $p_j(q)$ from $1$ to $|J_+|$ are in descending order if $\lambda_j$ are organized in decreasing order. Thus, by Chebyshev's sum inequality we obtain
                            \[
                                \sum_{j \in J_+} p_j(q) \frac{1}{r_j} \leq \frac{1}{|J_+|}\sum_{j \in J_+} \frac{1}{r_j} = \frac{1}{|J_+|}\sum_{j \in J_+} \frac{\delta^2}{\lambda_j}.
                            \]
                            Thus, a sufficient condition is then
                            \begin{equation*}
                                \left(\frac{1}{|J_+|}\sum_{j \in J_+} \frac{\delta^2}{\lambda_j}\right)^{-1} \geq r_\star,
                            \end{equation*}
                            which is equivalent to
                            \begin{equation}
                              \delta^2  \leq  \frac{1}{r_\star } \left(\frac{1}{|J_+|}\sum_{j \in J_+} \lambda_j^{-1}\right)^{-1}.
                            \end{equation}
                            In the case of power law covariance with $\alpha >0$, for finite $d$ this condition is equivalent to
                            \[
                             \delta \leq \sqrt{\frac{d^{1-\alpha}}{r_\star \sum_{j=1}^d j^{-\alpha}}}.
                            \]
                            Hence, the condition for $0 < \alpha < 1$ follows by the asymptotic $\sum_{j=1}^d j^{-\alpha} \sim \frac{d^{1-\alpha}}{1 -\alpha}$ as $d\to \infty$ which implies as $d\to \infty$
                            \[
                             \delta^2 < \frac{1-\alpha}{r_\star}.
                            \]
                            In the case of $\alpha \geq 1$, the previous argument does not provide a sufficient condition under which $\delta >0$ as $d\to \infty$. The following argument provides a sufficient condition on $\delta$ depending on $\alpha$ for $\alpha > 0$. We have $\lambda_j = \left(\frac{j}{d}\right)^{\alpha}$ such that $r_j = \frac{(j/d)^{\alpha}}{\delta^2}$. Using a similar argument as above, we obtain
                            \[
                                \frac{B_d(q)}{A_d(q)} = \frac{1}{\delta^2}\frac{\frac{1}{d}\sum_{j=1}^d
                                \frac{
                                \left(\frac{j}{d}\right)^{2\alpha}
                                }{   \left(1+ \frac{\delta q\left(1+\sqrt{\frac{2}{\pi}}\delta q\right)}{\delta^2\left( \delta q+\sqrt{\frac{2}{\pi}}\right)}(\frac{j}{d})^{\alpha}\right)^2    }}{\frac{1}{d}\sum_{k=1 }^d
                                \frac{
                                    \left(\frac{k}{d}\right)^{\alpha}
                                }{\left(1+ \frac{\delta q\left(1+\sqrt{\frac{2}{\pi}}\delta q\right)}{\delta^2\left( \delta q+\sqrt{\frac{2}{\pi}}\right)}(\frac{k}{d})^{\alpha}\right)^2
                                }}.
                            \]
                            Define the function 
                            \[
                                U_{\alpha, d}(z) = \frac{\frac{1}{d}\sum_{j=1}^d
                                \frac{
                                \left(\frac{j}{d}\right)^{2\alpha}
                                }{   \left(1+ z (\frac{j}{d})^{\alpha}\right)^2    }}{\frac{1}{d}\sum_{k=1 }^d
                                \frac{
                                    \left(\frac{k}{d}\right)^{\alpha}
                                }{ \left(1+ z (\frac{k}{d})^{\alpha}\right)^2
                                }},
                            \]
                            then we have $ \frac{B_d(q)}{A_d(q)} = \frac{1}{\delta^2}  U_{\alpha, d}\left(\frac{ q\left(1+\sqrt{\frac{2}{\pi}}\delta q\right)}{\delta\left( \delta q+\sqrt{\frac{2}{\pi}}\right)}\right)$. Our goal is to bound $ U_{\alpha, d}$ from below uniformly in $z \geq 0$. Let $M_{\alpha, d} \defas \underset{1\leq j \leq d}{\operatorname{max}} \frac{d- j +1}{d}\left(\frac{j}{d}\right)^{2\alpha}$. Then for any $1\leq \ell \leq d$, we have
                            \[
                                \frac{1}{d}\sum_{j=1}^d \frac{
                                    \left(\frac{j}{d}\right)^{2\alpha}
                                    }{   \left(1+ z (\frac{j}{d})^{\alpha}\right)^2    } \geq \frac{1}{d}\sum_{j=\ell}^d \frac{
                                        \left(\frac{j}{d}\right)^{2\alpha}
                                        }{   \left(1+ z (\frac{j}{d})^{\alpha}\right)^2    } \geq \frac{d-\ell+1}{d}\frac{(\ell/d)^{2\alpha} }{(1+z)^2}.
                            \]
                            Taking the maximum over $1\leq \ell \leq d$, we obtain
                            \[
                                \frac{1}{d}\sum_{j=1}^d \frac{
                                    \left(\frac{j}{d}\right)^{2\alpha}
                                    }{   \left(1+ z (\frac{j}{d})^{\alpha}\right)^2    }  \geq \frac{M_{\alpha, d} }{(1+z)^2}.
                            \]
                            For the denominator, we have
                            \[
                                \frac{1}{d}\sum_{j=1}^d \frac{
                                    \left(\frac{j}{d}\right)^{\alpha}
                                    }{   \left(1+ z (\frac{j}{d})^{\alpha}\right)^2    }  \leq  \frac{1}{d}\sum_{j=1}^d \left(\frac{j}{d}\right)^{\alpha},
                            \]
                            and for $z>0$ 
                            \[
                                \frac{
                                    \left(\frac{j}{d}\right)^{\alpha}
                                    }{   \left(1+ z (\frac{j}{d})^{\alpha}\right)^2    }\leq \frac{1}{z},
                            \]
                            from which we obtain
                            \[
                                \frac{1}{d}\sum_{j=1}^d \frac{
                                    \left(\frac{j}{d}\right)^{\alpha}
                                    }{   \left(1+ z (\frac{j}{d})^{\alpha}\right)^2    } \leq \min\left( \frac{1}{d}\sum_{j=1}^d \left(\frac{j}{d}\right)^{\alpha}, \frac{1}{z}\right).
                            \]
                            Thus, combining these results leads to for $z>0$
                            \[
                                U_{\alpha, d}(z)  \geq \frac{M_{\alpha, d}}{(1+z)^2}   \max\left( \left(\frac{1}{d}\sum_{j=1}^d \left(\frac{j}{d}\right)^{\alpha}\right)^{-1}, z\right).
                            \]
                            We have $\max\left( \left(\frac{1}{d}\sum_{j=1}^d \left(\frac{j}{d}\right)^{\alpha}\right)^{-1}, z\right) \geq \left(1+z\right)\left(1+\frac{1}{d}\sum_{j=1}^d \left(\frac{j}{d}\right)^{\alpha}\right)^{-1}$ from which we obtain
                            \[
                                U_{\alpha, d}(z)  \geq \frac{\kappa_{\alpha,d}}{1+z}\quad \text{with}\quad  \kappa_{\alpha,d}\defas \frac{M_{\alpha, d}}{ 1+\frac{1}{d}\sum_{j=1}^d \left(\frac{j}{d}\right)^{\alpha}}.
                            \]
                            Define $s_{\star} = \sup_{q >0} \frac{\delta q\left(1+\sqrt{\frac{2}{\pi}} \delta q\right)}{\delta q + \sqrt{\frac{2}{\pi}}} \psi_\delta(q)$ where $\psi_\delta(q)$ is defined in Proposition~\ref{prop:fixed:point}. Using a similar argument as in the proof of Proposition~\ref{prop:fixed:point}, $s_{\star}$ is independent of $\delta$ and $0<s_\star < \infty$. Numerically, $s_{\star}\approx1.77 \times 10^{-3}$. Since $\psi_{\delta}(q) \leq r_{\star}$ then 
                            \[
                               \left(\delta^2+ \frac{\delta q\left(1+\sqrt{\frac{2}{\pi}} \delta q\right)}{\delta q + \sqrt{\frac{2}{\pi}}} \right)\psi_\delta(q) \leq \delta^2 r_\star +s_\star.
                            \]
                            Thus, it follows that if $\delta^2   \leq\frac{\kappa_{\alpha,d }- s_\star}{ r_\star} $ then 
                            \[
                                \frac{B_d(q)}{A_d(q)} \geq \psi_\delta(q) \quad \text{for all}\quad q>0.
                            \]
                            Now, as $d\to \infty$ we have
                            \[
                                \frac{1}{d}\sum_{j=1}^d \left(\frac{j}{d}\right)^{\alpha} \to \int_0^1 y^{\alpha}\dif y = \frac{1}{\alpha + 1},
                            \]
                            and $M_{\alpha, d} \to \max_{y\in [0,1]}(1-y)y^{2\alpha}$ with maximum attained at $y = \frac{2\alpha}{2\alpha + 1}$. Thus, we obtain for $d\to \infty$
                            \[
                                \kappa_{\alpha, d}\to \kappa_{\alpha} \defas \frac{\alpha+ 1}{\alpha+2} \frac{1}{2\alpha + 1}\left(\frac{2\alpha}{2\alpha + 1}\right)^{2\alpha}.
                            \]
                            Thus, for $d\to \infty$, in order to have $\delta >0$ we require $\kappa_{\alpha} > s_\star$. Let $\alpha_0$ denote the unique positive solution of $\kappa_{\alpha} = s_\star$ then numerically $\alpha_0 \approx 102.87$. The uniqueness follows from the intermediate value theorem and showing $\frac{\dif}{\dif \alpha }\kappa_{\alpha} <0$ for all $\alpha >0$ after noting that $\lim_{\alpha \to 0}\kappa_{\alpha} = \frac{1}{2}$ and $\lim_{\alpha \to \infty}\kappa_{\alpha} = 0$. 

                            Under any of the stated alternatives, the sufficient uniqueness condition from Proposition~\ref{prop:fixed:point} holds for sufficiently large $d$ and thus in the limit.
                        \end{proof}
    
                        \subsubsection{Supporting Results}
\begin{lemma}\label{lem:adv:risk:L2_ls}
    The $\ell_2$-adversarial least squares risk satisfies 
    \[
    \mathrsfs{R}^\Adv(t) \asymp\mathrsfs{R}(t) +\frac{\delta^2}{2}\widehat{\mathrsfs{B}}_{44}(t).
    \]
    For the upper bound, the inequality is given by $\mathrsfs{R}^\Adv(t) \leq 2 \left(\mathrsfs{R}(t) +\frac{\delta^2}{2}\widehat{\mathrsfs{B}}_{44}(t)\right)$.
\end{lemma}
\begin{proof}
    The lower bound follows immediately. The upper bound follows from an application of Young's inequality. For $a,b\in \R$, we have $ab \leq \tfrac{a^2}{2} + \frac{ b^2}{2}$. Setting $a = \sqrt{2 \mathrsfs{R}(t)}$ and $b = \delta \sqrt{\frac{2}{\pi} \widehat{\mathrsfs{B}}_{44}(t)}$ then it follows that 
    \[
        \mathrsfs{R}^\Adv(t) \leq 2\mathrsfs{R}(t) + \delta^2\left(\frac{1}{2}+ \frac{1}{\pi}\right)\widehat{\mathrsfs{B}}_{44}(t) \leq 2\left(\mathrsfs{R}(t) +\frac{\delta^2}{2}\widehat{\mathrsfs{B}}_{44}(t)\right).
    \]
\end{proof}
\begin{lemma}\label{lem:strong:con:Radv}
    As per Definition~\ref{def:adversarial:ground:truth}, recall $X^{\star, \Adv} = \underset{X}{\operatorname{argmin}} \;\mathcal{R}^\Adv(X)$. The $\ell_2$-adversarial risk $ \mathcal{R}^\Adv(X)$ is strongly convex with parameter $\delta^2(1- \tfrac{2}{\pi})$ and $X^{\star, \Adv}$ is its unique minimizer. 
\end{lemma}
\begin{proof}
    Recall 
    \begin{equation}
        \begin{aligned}
        &\mathcal{R}^{{\operatorname{adv}}}(X) =
        \mathcal{R}(X)+ \delta\sqrt{\frac{2}{\pi}}\|X\| \sqrt{2\mathcal{R}(X)} +\frac{\delta^2}{2}  \|X\|^2.
        \end{aligned}
        \end{equation}
        Completing the square, we can rewrite this as 
        \begin{equation}
            \mathcal{R}^{{\operatorname{adv}}}(X) =
        \frac{1}{2}\left(\sqrt{2\mathcal{R}(X)}+ \delta\sqrt{\frac{2}{\pi}}\|X\| \right)^2 +\frac{\delta^2}{2} \left(1-\frac{2}{\pi} \right)\|X\|^2.
        \end{equation}
        The map $X\mapsto \frac{1}{2}\left(\sqrt{2\mathcal{R}(X)}+ \delta\sqrt{\frac{2}{\pi}}\|X\| \right)^2$ is convex. For all $\delta >0$, the second term is  $\delta^2(1- \tfrac{2}{\pi})$-strongly convex. Hence, the result follows.
\end{proof}
\begin{lemma}\label{lem:helper:bounds:diffforQ}
    Recall $q(t) = \sqrt{\frac{\widehat{\mathrsfs{B}}_{44}(t)}{2\mathrsfs{R}(t)}}$. For any $T>0$, suppose $\mathrsfs{R}(t) >0$ for all $ t\in [0,T]$ then for all $ t\in [0,T]$
    \begin{equation}
        \begin{gathered}
            \frac{\left|\sum_{j=1}^d \lambda_j \left(\mathrsfs{V}_{j,11}(t) - \mathrsfs{V}_{j,12}(t) \right)\right|}{\mathrsfs{R}(t)} \leq 2\|K\|_{\opt}^{1/2} q(t) +2,
        \end{gathered}
    \end{equation}
    and
    \begin{equation}
        \begin{gathered}
            \frac{\left|\sum_{j=1}^d \lambda_j \left(\mathrsfs{V}_{j,11}(t) - \mathrsfs{V}_{j,12}(t) \right)\right|}{\mathrsfs{R}(t)} \leq 4 \|K\|_{\opt}  (q(t))^2+ 2 \|K\|_{\opt}^{1/2}q(t).
        \end{gathered}
    \end{equation}
\end{lemma}
\begin{proof}
    For simplicity of notation, define $\mathcal{E}_j(t) = \mathrsfs{V}_{j,11}(t) -  2\mathrsfs{V}_{j,12}(t) +  \mathrsfs{V}_{j,22}(t)$ then $2\mathrsfs{R}(t) = \sum_{j=1}^d \lambda_j \mathcal{E}_j(t) + \eta^2$. Since $\mathrsfs{V}_j(t)$ is psd then its determinant is nonnegative which implies $(\mathrsfs{V}_{j,12}(t))^2 \leq \mathrsfs{V}_{j,11}(t)\mathrsfs{V}_{j,22}(t)$. Using this inequality, it follows that 
    \begin{equation}\label{eq:lower:bound:E_j}
        \mathcal{E}_j(t) \geq \left(\sqrt{\mathrsfs{V}_{j,11}(t)} - \sqrt{\mathrsfs{V}_{j,22}(t)}\right)^2.
    \end{equation}
    We also obtain the following bounds using Cauchy--Schwarz
    \begin{equation}\label{eq:lower:bound:V_22_j}
        \begin{gathered}
        2\mathrsfs{R}(t) \geq \left(\sqrt{\sum_{j=1}^d \lambda_j \mathrsfs{V}_{j,11}(t) } - \sqrt{\sum_{j=1}^d \lambda_j \mathrsfs{V}_{j,22}(t) }\right)^2,
        \\
        \text{and}\quad\left|\sum_{j=1}^d \lambda_j \mathrsfs{V}_{j,12}(t) \right| \leq \sqrt{\sum_{j=1}^d \lambda_j \mathrsfs{V}_{j,11}(t)}\sqrt{\sum_{j=1}^d \lambda_j \mathrsfs{V}_{j,22}(t)}.
        \end{gathered}
    \end{equation}
    We start by deriving the first bound. We may rewrite the difference $\mathrsfs{V}_{j,11}(t) - \mathrsfs{V}_{j,12}(t)$ as follows
    \begin{equation}
        \begin{aligned}
            \mathrsfs{V}_{j,11}(t) - \mathrsfs{V}_{j,12}(t) = \frac{1}{2}  \mathcal{E}_j(t) + \frac{1}{2}\left( \mathrsfs{V}_{j,11}(t) - \mathrsfs{V}_{j,22}(t)\right).
        \end{aligned}
    \end{equation}
Summing over $\lambda_j$, it follows
\begin{equation}
    \sum_{j=1}^d \lambda_j \left( \mathrsfs{V}_{j,11}(t) - \mathrsfs{V}_{j,12}(t) \right) \leq \mathrsfs{R}(t) + \frac{1}{2}\sum_{j=1}^d \lambda_j  \left( \mathrsfs{V}_{j,11}(t) - \mathrsfs{V}_{j,22}(t)\right).
\end{equation}
Since $\mathrsfs{V}_{j,11}(t), \mathrsfs{V}_{j,22}(t) \geq 0$ and using \eqref{eq:lower:bound:E_j} then 
\[
\begin{aligned}
    \left|\mathrsfs{V}_{j,11}(t) - \mathrsfs{V}_{j,22}(t) \right|
     &
     =
     \left|\sqrt{\mathrsfs{V}_{j,11}(t)} - \sqrt{ \mathrsfs{V}_{j,22}(t)}\right| \left(\sqrt{\mathrsfs{V}_{j,11}(t)} + \sqrt{\mathrsfs{V}_{j,22}(t)}\right)
     \\
     &
     \leq
    \sqrt{  \mathcal{E}_j(t) }\left(\sqrt{\mathrsfs{V}_{j,11}(t)} + \sqrt{\mathrsfs{V}_{j,22}(t)}\right).
\end{aligned} 
\]
Summing over $\lambda_j$ and from Cauchy--Schwarz, it follows
\begin{equation}
    \begin{aligned}
    \left|\sum_{j=1}^d \lambda_j \left( \mathrsfs{V}_{j,11}(t) - \mathrsfs{V}_{j,22}(t)\right) \right| \leq \sqrt{2 \mathrsfs{R}(t)}\left(\sqrt{\sum_{j=1}^d \lambda_j \mathrsfs{V}_{j,11}(t) } + \sqrt{\sum_{j=1}^d \lambda_j \mathrsfs{V}_{j,22}(t) }\right).
    \end{aligned}
\end{equation}
Now, using the definition of $q(t)$ we obtain the simple bound
\begin{equation}\label{eq:bound:lambda_jV_j_11:sum}
    \sum_{j=1}^d \lambda_j \mathrsfs{V}_{j,11}(t) \leq \|K\|_{\opt} \widehat{\mathrsfs{B}}_{44}(t) = 2 \|K\|_{\opt} \mathrsfs{R}(t) (q(t))^2. 
\end{equation}
Using this bound and \eqref{eq:lower:bound:V_22_j}, it follows that
\begin{equation}
     \sqrt{\sum_{j=1}^d \lambda_j \mathrsfs{V}_{j,22}(t) } \leq  \sqrt{\sum_{j=1}^d \lambda_j \mathrsfs{V}_{j,11}(t) }+ \sqrt{2\mathrsfs{R}(t)} \leq \sqrt{2\mathrsfs{R}(t)} \left(1 +  \|K\|_{\opt}^{1/2} q(t) \right).
\end{equation}
Thus, from these previous bounds, we obtain
\begin{equation}
    \begin{aligned}
    \frac{1}{2}\left|\sum_{j=1}^d \lambda_j \left( \mathrsfs{V}_{j,11}(t) - \mathrsfs{V}_{j,22}(t)\right) \right| \leq  \mathrsfs{R}(t)\left(2\|K\|_{\opt}^{1/2} q(t) +1 \right).
    \end{aligned}
\end{equation}
Dividing by $\mathrsfs{R}(t)$ obtains the first bound. For the second bound, from \eqref{eq:lower:bound:V_22_j}, \eqref{eq:bound:lambda_jV_j_11:sum} and \eqref{eq:lower:bound:V_22_j} we have
\begin{equation*}
    \begin{aligned}
    \frac{\left|\sum_{j=1}^d \lambda_j \left( \mathrsfs{V}_{j,11}(t) - \mathrsfs{V}_{j,12}(t)\right) \right|}{\mathrsfs{R}(t)} 
    &
    \leq \frac{\left|\sum_{j=1}^d \lambda_j \mathrsfs{V}_{j,11}(t)\right| +  \sqrt{\sum_{j=1}^d \lambda_j \mathrsfs{V}_{j,11}(t)}\sqrt{\sum_{j=1}^d \lambda_j \mathrsfs{V}_{j,22}(t)}}{\mathrsfs{R}(t)}
    \\
    &
    \leq 4 \|K\|_{\opt}  (q(t))^2+ 2 \|K\|_{\opt}^{1/2}q(t)    .
    \end{aligned}
\end{equation*}
\end{proof}
\begin{lemma}\label{lem:R_up_low_bound}
    Let
    \[
        \mathrsfs{R}(t) +\frac{\delta^2}{2}\widehat{\mathrsfs{B}}_{44}(t) 
        =  F^\Adv(\tilde{\Gamma}(t)) + \int_0^t \mathscr{K}(\tilde{\Gamma}(t) - \tilde{\Gamma}(s), \Lambda(t) - \Lambda(s))\mathrsfs{R}^\Adv(s) \dif s.
        \]
    Then for all $t$, 
    \begin{align*}
        \mathrsfs{R}(t) +\frac{\delta^2}{2}\widehat{\mathrsfs{B}}_{44}(t)&\geq F^\Adv(\tilde{\Gamma}(t)) +\int_0^t\mathscr{K}(\tilde{\Gamma}(t)-\tilde{\Gamma}(s), \Lambda(t)-\Lambda(s))F^\Adv(\tilde{\Gamma}(s)) \dif s
    \end{align*}
    and
    \[
        \mathrsfs{R}(t) +\frac{\delta^2}{2}\widehat{\mathrsfs{B}}_{44}(t) 
        \leq  F^\Adv(\tilde{\Gamma}(t)) 
+ \int_0^{t} \overline{\mathscr{K}}(t - s)\left(\mathrsfs{R}(s) +\frac{\delta^2}{2}\widehat{\mathrsfs{B}}_{44}(s) \right)\dif s,
    \]
    where $\overline{\mathscr{K}}(x) \defas \frac{2\gamma^2}{d}\sum_{j=1}^d  (\lambda_j^2 + \delta^2 \lambda_j) e^{-2  \gamma(\lambda_j + \delta^2)x} $. Define $\overline{\mathscr{L}}(x) = e^{2\gamma \delta^2 x}\overline{\mathscr{K}}(x)$. If $\lambda_j =1 $ for all $1\leq j\leq d$. Then for $\gamma < 1$ we have
    \begin{align*}
        \mathrsfs{R}(t) +\frac{\delta^2}{2}\widehat{\mathrsfs{B}}_{44}(t) \leq   F^\Adv(\tilde{\Gamma}(t))  +  2\gamma^2(1+\delta^2)\int_0^t e^{-2\gamma(1+\delta^2)(1-\gamma)(t-s)} F^\Adv(\tilde{\Gamma}(s)) \dif s.
    \end{align*}
    \end{lemma}
    \begin{proof}
        The lower bound holds immediately from Lemma~\ref{lem:adv:risk:L2_ls} and $\mathrsfs{R}(s) +\frac{\delta^2}{2}\widehat{\mathrsfs{B}}_{44}(s)\geq F^\Adv(\tilde{\Gamma}(s))$. Notice that for $x,y \geq 0$
        \begin{equation*}
           \mathscr{K}(x,y) \leq \sum_{j=1}^d (\lambda_j^2 + \delta^2 \lambda_j) e^{-2 \lambda_j x} .
        \end{equation*}
        Given that $\tilde{\Gamma}(t) \geq \gamma t$ and $\Lambda(t) \geq  \gamma \delta^2 t$ for all $t\geq 0$ then we obtain for $0\leq s\leq t$:
        \begin{equation*}
            \mathscr{K}(\tilde{\Gamma}(t) - \tilde{\Gamma}(s), \Lambda(t) - \Lambda(s)) \leq \frac{\gamma^2}{d}\sum_{j=1}^d (\lambda_j^2 + \delta^2 \lambda_j) e^{-2  \gamma(\lambda_j + \delta^2)(t-s)} .
         \end{equation*}
         It then follows from Lemma~\ref{lem:adv:risk:L2_ls} that
    \[
    \begin{aligned}
        \mathrsfs{R}(t) +&\frac{\delta^2}{2}\widehat{\mathrsfs{B}}_{44}(t) 
        \leq  F^\Adv(\tilde{\Gamma}(t)) 
+ \int_0^{t} \overline{\mathscr{K}}(t - s)\left(\mathrsfs{R}(s) +\frac{\delta^2}{2}\widehat{\mathrsfs{B}}_{44}(s) \right)\dif s,
    \end{aligned}
        \]
        which is the upper bound. We define the convolution map 
     \[
     \mathcal{G}(f)(t) = (\overline{\mathscr{K}}*f)(t) =\int_0^{t} \overline{\mathscr{K}}(t - s) f(s) \dif s.
     \]
    We now show this map is contracting
    \begin{align*}\label{eq:G_G_f}  
    \mathcal{G}^2(f)(t) = \mathcal{G}(\mathcal{G}(f))(t) &=  \int_0^t \overline{\mathscr{K}}(t - s) \mathcal{G}(f)(s) \dif s
    \\ &= \int_0^t \overline{\mathscr{K}}(t - s)  \int_0^s \overline{\mathscr{K}}(s - u) f(u) \dif u \dif s
    \\ &= 
    \int_0^t \left(\int_{u}^t \overline{\mathscr{K}}(t - s)\overline{\mathscr{K}}(s- u)\dif s \right)f(u) \dif u 
    \\
    &= 
    \int_0^t \overline{\mathscr{K}}^{*2}(t-u) f(u) \dif u,
    \end{align*}
    where the last equality is obtained through a change of variables. Hence, repeatedly applying the convolution map yields by induction that
    \[
        \mathcal{G}^j(f)(t) =  \int_0^t \overline{\mathscr{K}}^{*j}(t-u) f(u) \dif u.
    \]
    Under isotropic covariance, we have $ \overline{\mathscr{K}}(x) \defas 2\gamma^2  (1+ \delta^2 ) e^{-2  \gamma(1+ \delta^2)x} $ from which it follows that $\|\overline{\mathscr{K}}\|_1 = \gamma$. Thus, for $\gamma < 1$ we have
 , we obtain
    \[
    \begin{aligned}
        \mathrsfs{R}(t) &+\frac{\delta^2}{2}\widehat{\mathrsfs{B}}_{44}(t) 
        \leq 
         F^\Adv(\tilde{\Gamma}(t)) + \sum_{j=1}^\infty  \mathcal{G}^j(F^\Adv\circ \tilde{\Gamma})(t)
        \\
        &
        \leq 
         F^\Adv(\tilde{\Gamma}(t)) + \left(\left( \sum_{j=1}^\infty\overline{\mathscr{K}}^{*j}\right) * (F^\Adv \circ  \tilde{\Gamma})\right)(t) 
        \\
        & 
        \leq
        F^\Adv(\tilde{\Gamma}(t))  +  2\gamma^2(1+\delta^2)\int_0^t e^{-2\gamma(1+\delta^2)(1-\gamma)(t-s)} F^\Adv(\tilde{\Gamma}(s)) \dif s,
    \end{aligned}
    \]
    where we used $\sum_{j=1}^\infty\overline{\mathscr{K}}^{*j}(t) = \sum_{j=1}^\infty \frac{(2\gamma^2(1+\delta^2))^j t^{j-1}}{(j-1)!} e^{-2\gamma(1+\delta^2)t} = 2\gamma^2(1+\delta^2)e^{-2\gamma(1+\delta^2)(1-\gamma)t}$.
\end{proof}
The following lemma follows a similar argument as Lemma D.3 in \citep{collinswoodfin2024highline}.
\begin{lemma}\label{lem:R_up_low_bound:2}
    Let $\overline{\mathscr{K}}_\infty \defas\lim_{d\to \infty}\overline{\mathscr{K}}$ and $\overline{\mathscr{L}}_\infty \defas\lim_{d\to \infty}\overline{\mathscr{L}}$ where $\overline{\mathscr{K}}$ and $\overline{\mathscr{L}}$ are defined in Lemma~\ref{lem:R_up_low_bound}. If $\overline{\mathscr{K}}_\infty $ and $\overline{\mathscr{L}}_\infty$ are monotonically decreasing, with $\|\overline{\mathscr{K}}_\infty\|_1, \|\overline{\mathscr{L}}_\infty\|_1<\infty$ and if there exist $\varepsilon >0$ and $T>0$ such that, for all $t>T$,
    \begin{equation}\label{eq:power_law_assmp_K}   
        \int_0^t\overline{\mathscr{L}}_\infty(s)\overline{\mathscr{L}}_\infty(t-s)\dif s\leq 2(1+\varepsilon)\|\overline{\mathscr{L}}_\infty\|_1\overline{\mathscr{L}}_\infty(t)\quad\text{and}\quad 2\|\overline{\mathscr{L}}_\infty\|_1(1+\varepsilon)
    <1
    \end{equation}
    then for all $t\geq 1$
    \begin{align*}
        \mathrsfs{R}_\infty(t) +\frac{\delta^2}{2}\widehat{\mathrsfs{B}}_{\infty,44}(t) \leq  F_\infty^\Adv(\tilde{\Gamma}(t)) + C\int_0^t\overline{\mathscr{K}}_\infty(t-s) F_{\infty}^\Adv(\tilde{\Gamma}(s)) \dif s,
    \end{align*}
    where $\mathrsfs{R}_\infty(t) \defas \lim_{d\to \infty}\mathrsfs{R}(t)$, $\widehat{\mathrsfs{B}}_{_\infty,44}(t)\defas \lim_{d\to \infty}\widehat{\mathrsfs{B}}_{44}(t)$, $F_\infty^\Adv(\tilde{\Gamma}(t)) \defas \lim_{d\to \infty} F^\Adv(\tilde{\Gamma}(t))$ and for
    \[
    C = \left(\frac{\overline{\mathscr{K}}_\infty(0)}{\overline{\mathscr{K}}_\infty(T)(2\varepsilon +1)}+2\right)\frac{1}{1-2\|\overline{\mathscr{L}}_\infty\|_1(1+\varepsilon)
    }.
    \]
    \end{lemma}
    \begin{proof}
        We now take $d\to \infty$ and omit the $\infty$ subscript for notational simplicity. We define the convolution map 
     \[
     \mathcal{G}(f)(t) = (\overline{\mathscr{K}}*f)(t) =\int_0^{t} \overline{\mathscr{K}}(t - s) f(s) \dif s.
     \]
   Similarly to Lemma~\ref{lem:R_up_low_bound}, we obtain by induction
    \[
        \mathcal{G}^j(f)(t) =  \int_0^t \overline{\mathscr{K}}^{*j}(t-u) f(u) \dif u.
    \]
    Now from the Assumption of the Lemma, we have
    \begin{equation*}
        \overline{\mathscr{K}}^{*2} = \int_0^t  e^{-2\gamma \delta^2 s}\overline{\mathscr{L}}(s)e^{-2\gamma \delta^2 (t-s)}\overline{\mathscr{L}}(t-s) \dif s = e^{-2\gamma \delta^2 t} \int_0^t \overline{\mathscr{L}}(s)\overline{\mathscr{L}}(t-s) \leq  2(1+\varepsilon)\|\overline{\mathscr{L}}\|_1\overline{\mathscr{K}}(t).
    \end{equation*}
    Thus, from Lemma~\ref{lem:Kesten} and noting that $\|\overline{\mathscr{K}}\|_1 \leq \|\overline{\mathscr{L}}\|_1$ we obtain the upper bound
    \[
        \sup_{t\ge 0} \frac{\overline{\mathscr{K}}^{*n}(t)}{\overline{\mathscr{K}}(t)} 
        \leq 
         (2\|\overline{\mathscr{L}}\|_1(1+\varepsilon))^{n-1}C_1.
    \]
    with $C_1 =\frac{\overline{\mathscr{K}}(0)}{\overline{\mathscr{K}}(T)(2\varepsilon +1)}+1$. Expanding the upper bound and using this identity with $q = 2(1+\varepsilon) \|\overline{\mathscr{L}}\|_1 $ such that $q<1$, we obtain
    \[
    \begin{aligned}
        \mathrsfs{R}(t) &+\frac{\delta^2}{2}\widehat{\mathrsfs{B}}_{44}(t) 
        \leq 
         F^\Adv(\tilde{\Gamma}(t)) + \sum_{j=1}^\infty  \mathcal{G}^j(F^\Adv\circ \tilde{\Gamma})(t)
        \\
        &
        \leq 
         F^\Adv(\tilde{\Gamma}(t)) + \sum_{j=1}^\infty \int_0^t \overline{\mathscr{K}}^{*j}(t-u) F^\Adv(\tilde{\Gamma}(u)) \dif u
        \\
        &
        \leq 
         F^\Adv(\tilde{\Gamma}(t)) + \left( 1+C_1\left(\sum_{j=0}^\infty  (2\|\overline{\mathscr{L}}\|_1(1+\varepsilon))^{j} -1\right)\right)\int_0^t\overline{\mathscr{K}}(t-u) F^\Adv(\tilde{\Gamma}(u)) \dif u
        \\
        &
        \leq 
         F^\Adv(\tilde{\Gamma}(t)) + \left(1+C_1\frac{q}{1-q}\right)\int_0^t\overline{\mathscr{K}}(t-u) F^\Adv(\tilde{\Gamma}(u)) \dif u.
    \end{aligned}
    \]
\end{proof}
The proof of the following lemma can be found in Lemma D.4 of Appendix D in \citep{collinswoodfin2024highline}.
\begin{lemma}[Lemma IV.4.7 in \cite{athreya2004branching}, Lemma D.4 in \citep{collinswoodfin2024highline}]\label{lem:Kesten}
    Suppose $\mathrsfs{K}$ is  monotonically decreasing, with $\|\mathrsfs{K}\|_1<\infty$, and that there exists $T>0$ such that $\forall t\geq T$, and $\varepsilon \ge 0$, 
    \begin{align} 
    \int_0^t\mathrsfs{K}(s)\mathrsfs{K}(t-s)\dif s\leq2(1+\varepsilon)\|\mathrsfs{K}\|_1\mathrsfs{K}(t).
    \end{align}
    Then, 
    \begin{align*}
    \sup_{t\ge 0} \frac{\mathrsfs{K}^{*n}(t)}{\mathrsfs{K}(t)} 
    \leq 
     (2\|\mathrsfs{K}\|_1(1+\varepsilon))^{n-1}\left(\frac{\mathrsfs{K}(0)}{\mathrsfs{K}(T)(2\varepsilon +1)}+1\right).
    \end{align*}
    \end{lemma}
    This Lemma is analogous to Lemma D.5 in \citep{collinswoodfin2024highline}.
    \begin{lemma}\label{lem:powerlaw:Linearreg}
        Suppose $\alpha  > 0$ and $\lambda_j = \left(\frac{j}{d}\right)^\alpha$. Recall $\overline{\mathscr{K}}(x) \defas \frac{2\gamma^2}{d}\sum_{j=1}^d  (\lambda_j^2 + \delta^2 \lambda_j) e^{-2  \gamma(\lambda_j + \delta^2)x} $ and $\overline{\mathscr{L}}(x) = e^{2\gamma \delta^2 x}\overline{\mathscr{K}}(x)$. Taking $d\to \infty$, then for $\delta > 0$ and $\kappa = 1 + \frac{1}{\alpha}$, $\overline{\mathscr{K}}(t) \asymp_{\alpha, \gamma,\delta} e^{-2  \gamma\delta^2 t} t^{-\kappa}$. Note that when $\delta = 0$, $\kappa = 2 + \frac{1}{\alpha}$, $\overline{\mathscr{K}}(t) \asymp_{\alpha, \gamma} t^{-\kappa}$. The condition \eqref{eq:power_law_assmp_K} is satisfied for all $\alpha >0$ with 
        \[
            \gamma < \frac{1}{2\left(\delta^2 + \frac{1}{\alpha + 1}\right)}.
        \]
    \end{lemma}
    \begin{proof}
        From the definition of $\overline{\mathscr{K}}$ in Lemma~\ref{lem:R_up_low_bound}, we have
         \begin{equation}\label{eq:asymp:Kernel}
            \begin{aligned}
            \overline{\mathscr{K}}(x) &= \frac{2\gamma^2}{d}\sum_{j=1}^d (\lambda_j^2 + \delta^2 \lambda_j) e^{-2  \gamma(\lambda_j +\delta^2 )x} 
            \\
            &
            =
            e^{-2  \gamma\delta^2 x}   \overline{\mathscr{L}}(x).
            \\
            &
            =
            2\gamma^2 e^{-2  \gamma\delta^2 x} \left(\int_0^1 y^{2\alpha}  e^{-2  \gamma  x y^\alpha } \dif y+ \delta^2 \int_0^1y^\alpha e^{-2  \gamma  x y^\alpha } \dif y\right)+ \mathcal{O}(d^{-1})
            \end{aligned}
         \end{equation}
         Using a change of variable
         \begin{equation*}
            \begin{gathered}
            I_p(x) = \int_0^1 y^{p}  e^{-2  \gamma  x y^\alpha } \dif y = \frac{1}{\alpha (2\gamma)^{\frac{p+1}{\alpha}}} x^{-\frac{p+1}{\alpha}} \boldsymbol{\gamma}\left(\frac{p+1}{\alpha},2\gamma x\right).
            \end{gathered}
         \end{equation*}
         where $\boldsymbol{\gamma}(s,z) \defas \int_0^z t^{s-1}e^{-t} \dif t$ is the incomplete gamma function. We have $\lim_{z\to \infty} \boldsymbol{\gamma}(s,z)  = \Gamma(s)$, the complete gamma function. Thus, we obtain  
         \begin{equation*}
            I_\alpha(x) \sim C(\gamma,\alpha,\delta) x^{-(1+ \frac{1}{\alpha})} \quad \text{and} \quad I_{2\alpha}(x)  \sim C(\gamma,\alpha,\delta) x^{-(2+ \frac{1}{\alpha})},
         \end{equation*}
         from which it follows that
\begin{equation}\label{eq:powerlaw:scaling}
    \overline{\mathscr{L}}(x)\sim \begin{cases} 
        C(\gamma,\alpha,\delta)x^{-(1+ \frac{1}{\alpha})} \quad \text{for}\quad \delta >0,\\
        C(\gamma,\alpha)x^{-(2+ \frac{1}{\alpha})}  \quad \text{for}\quad \delta =0.
    \end{cases} 
\end{equation}
Taking $d\to \infty$, from a simple computation we obtain for 
\begin{align*}
  \begin{gathered}
  \|\overline{\mathscr{L}}\|_1 =\gamma\left( \frac{\Tr(K)}{d} +  \delta^2\right) \to \gamma\left(\frac{1}{\alpha + 1}+ \delta^2\right) < \infty,\quad 
  \|\overline{\mathscr{K}}\|_1 =\gamma \frac{\Tr(K)}{d} \to \frac{\gamma}{\alpha + 1}  < \infty.
  \end{gathered}
\end{align*}
We now show for $d\to \infty$ that $ \overline{\mathscr{L}}(t) \asymp_{\gamma,\alpha,\delta} t^{-\kappa}$ satisfies \eqref{eq:power_law_assmp_K} for $\kappa > 1$. To do so, we define $\overline{\mathscr{L}}_{\infty} = \lim_{d\to \infty} \overline{\mathscr{L}}$ we show that for $t\to \infty$
\[
 \frac{(\overline{\mathscr{L}}_{\infty} * \overline{\mathscr{L}}_{\infty})(t)}{\overline{\mathscr{L}}_{\infty}(t)} \to 2 \|\overline{\mathscr{L}}_{\infty}\|_{1}.
\]
Taking $d\to \infty$ in \eqref{eq:asymp:Kernel} and , we have 
\begin{equation*}
        \overline{\mathscr{L}}_{\infty}(x) =2\gamma^2  \left(\int_0^1 y^{2\alpha}  e^{-2  \gamma  x y^\alpha } \dif y+ \delta^2 \int_0^1y^\alpha e^{-2  \gamma  x y^\alpha } \dif y\right).
 \end{equation*}
 It follows that $\overline{\mathscr{L}}_{\infty} $ is nonincreasing. We split $\overline{\mathscr{L}_{\infty}} * \overline{\mathscr{L}_{\infty}}$ as follows
 \begin{equation}
    \begin{aligned}
    (\overline{\mathscr{L}}_{\infty} * \overline{\mathscr{L}}_{\infty})(t) &= \int_0^{\sqrt{t}} \overline{\mathscr{L}}_{\infty}(s)\overline{\mathscr{L}}_{\infty}(t-s) \dif s 
    \\
    &+ \int_{\sqrt{t}}^{t-\sqrt{t}} \overline{\mathscr{L}}_{\infty}(s)\overline{\mathscr{L}}_{\infty}(t-s) \dif s + \int_{t - \sqrt{t}}^t\overline{\mathscr{L}}_{\infty}(s)\overline{\mathscr{L}}_{\infty}(t-s) \dif s .
    \end{aligned}
 \end{equation}
 Since $\overline{\mathscr{L}}_{\infty}(t) \sim C(\gamma, \alpha, \delta) t^{-\kappa}$ then it follows that $\overline{\mathscr{L}}$ is a regularly varying function with index $-\kappa$. Hence, we can write $\overline{\mathscr{L}}(t)= t^{-\kappa} L(t)$ for some slowly varying function $L(t) \defas t^{\kappa} \overline{\mathscr{L}}(t)$. From $\overline{\mathscr{L}}_{\infty}(t)  \sim C(\gamma, \alpha, \delta)  t^{-\kappa}$, it follows immediately that $L(t) \to C>0$ as $t\to \infty$ for some constant $C(\gamma, \alpha, \delta)>0$. For $0 \leq s \leq \sqrt{t}$, we have $\frac{t-s}{t} \in [1 - \frac{1}{\sqrt{t}}, 1]$ from which it follows that $\left(\frac{t-s}{t}\right)^{-\kappa} \to 1$ uniformly in $s \in [0, \sqrt{t}]$. Similarly, since $\sqrt{t} = o(t)$ and $L(t) \to C$ as $t\to \infty$ then $\frac{L(t-s)}{L(t)} \to 1$ uniformly in $[0, \sqrt{t}]$. Combining these results, we obtain $\frac{\overline{\mathscr{L}}_{\infty}(t-s)}{\overline{\mathscr{L}}_{\infty}(t)}$ also converges uniformly to $1$ on $0 \leq s\leq \sqrt{t}$ from which it follows
\[
    \frac{1}{ \overline{\mathscr{L}}_{\infty}(t)}\int_0^{\sqrt{t}} \overline{\mathscr{L}}_{\infty}(s)\overline{\mathscr{L}}_{\infty}(t-s) \dif s  = (1+ o(1))\int_0^{\sqrt{t}} \overline{\mathscr{L}}_{\infty}(s)\dif s \underset{t \to \infty}{\longrightarrow} \|\overline{\mathscr{L}}_{\infty}\|_1.
\]
By symmetry, we also obtain $\frac{1}{ \overline{\mathscr{L}}_{\infty}(t)}\int_{t-\sqrt{t}}^{t} \overline{\mathscr{L}}_{\infty}(s)\overline{\mathscr{L}}_{\infty}(t-s) \dif s\underset{t \to \infty}{\longrightarrow} \|\overline{\mathscr{L}}_{\infty}\|_1$. For the middle integral, for large enough $t$ we have $\sqrt{t} \leq\frac{t}{2}$ so by symmetry and since $\overline{\mathscr{L}}_{\infty}$ is nonincreasing we obtain
\[
\begin{aligned}
    \int_{\sqrt{t}}^{t-\sqrt{t}} \overline{\mathscr{L}}_{\infty}(s)\overline{\mathscr{L}}_{\infty}(t-s) \dif s = 2  \int_{\sqrt{t}}^{\frac{t}{2}}\overline{\mathscr{L}}_{\infty}(s)\overline{\mathscr{L}}_{\infty}(t-s) \dif s \leq 2\overline{\mathscr{L}}_{\infty}(t/2) \int_{\sqrt{t}}^{\infty} \overline{\mathscr{L}}_{\infty}(s)\dif s .
\end{aligned}
\]
Since $L(t)$ is a slowly varying function, we have $\frac{\overline{\mathscr{L}}_{\infty}(t/2) }{\overline{\mathscr{L}}_{\infty}(t) } \to 2^{\kappa}$ as $t\to \infty$ from which it follows that $\frac{1}{\overline{\mathscr{L}}_{\infty}(t)} \int_{\sqrt{t}}^{t-\sqrt{t}} \overline{\mathscr{L}}_{\infty}(s)\overline{\mathscr{L}}_{\infty}(t-s) \dif s \to 0$ as $t\to \infty$. Putting this altogether, we obtain
\[
 \frac{(\overline{\mathscr{L}}_{\infty} * \overline{\mathscr{L}}_{\infty})(t)}{\overline{\mathscr{L}}_{\infty}(t)} \underset{t\to \infty}{\longrightarrow} 2 \|\overline{\mathscr{L}}_{\infty}\|_{1}.
\]
Thus, for every $\varepsilon >0$ there exists a $T \geq 1$ such that for all $t \geq T$, $(\overline{\mathscr{L}}_{\infty} * \overline{\mathscr{L}}_{\infty})(t) \leq 2(1+\varepsilon) \|\overline{\mathscr{L}}_{\infty} \|_1 \overline{\mathscr{L}}_{\infty}(t)$.
    \end{proof}
\subsection{Binary Logistic Regression}\label{ex:binary:logreg}
We consider a binary logistic regression problem with $\epsilon=0$ where we are trying to classify two classes. Setting the hard labels $y\in\{-1,1\}$, the standard logistic regression objective function (cross-entropy) can be rewritten as
\begin{align}
    f_I(X^\top a_I) = \log(1 + e^{-y_I\cdot X^\top a_I}).
\end{align}
Since the map $t\mapsto \log(1+ e^{-t})$ is convex and strictly decreasing, the inner-maximization problem in the adversarial risk simplifies to \citep{bach2023adversarial}
\begin{equation}\label{eq:def:LOG:REG:risk}
\begin{aligned}
    \mathcal{R}^{{\operatorname{adv}}}(X) 
   &
    =  
    \mathbb{E}_{a, I}\left[\max_{|s|\leq 1}\log(1 + \exp(-y_I\cdot( X^\top a_I + \delta s \|X\|)))\right] 
    \\
    &
    = 
    \mathbb{E}_{a, I}\left[\log(1 + \exp(-y_I\cdot(X^\top a_I) + \delta \|X\|))\right]
    \\
    &
    =
    \sum_{i=1}^2p_i \mathbb{E}_{v}\left[\log(1 + \exp((-1)^{i}\cdot(\sqrt{X^\top K_i X} v + X^\top \mu_i) +\delta \|X\|))\right]
\end{aligned}
\end{equation}
where $v \sim \mathcal{N}(0,1)$. To simplify the derivation, we take $\mu_1 = \mu$ and $\mu_2 = - \mu$. We note that this is equivalent to any other choices of means for binary logistic regression by recentering as follows $\tilde{a} = a - c$ with $c = \tfrac{\mu_1 + \mu_2}{2}$ such that
\[
    \mu_1 - c = \frac{\mu_1-\mu_2}{2} \defas \mu \quad \text{and}\quad  \mu_2 - c = -\frac{\mu_1-\mu_2}{2} \defas -\mu.
\]
We are studying the hard label setting, hence consider $X^{\star} = 0$. We identify the following terms 
\begin{gather*}
    z_{11} = X^\top K_i X ,\quad z_{13} = X^\top \mu_i , \quad z_{31} =  z_{13} \quad \text{and}\quad  z_{44} = \|X\|^2,
\end{gather*}
which leads to
\[
\begin{aligned}
    h_i\left(\begin{bmatrix}
        z_{11} & 0 &  z_{13} & 0 \\
       0 &0 & 0& 0 \\
     z_{31} &  0& 0 & 0\\
        0 & 0 & 0  & z_{44}
    \end{bmatrix}\right) &=  \mathbb{E}_{v}\left[\log\left(1 + \exp\left((-1)^{i}\left(\sqrt{z_{11}} v  + \frac{z_{13} + z_{31}}{2}\right)+ \delta \sqrt{z_{44}}\right)\right)\right]
\end{aligned}
\]
Here we include both $z_{13}$ and $z_{31}$ to stay consistent with the gradient derived in Lemma~\ref{lem:derivative_adv_risk}.
It follows from \eqref{eq:def:LOG:REG:risk} that 
\begin{equation}
    \mathcal{R}^{{\operatorname{adv}}}(X) = \sum_{i=1}^2 p_i h_i(\widehat{B}_i(W)),
\end{equation}
with  
\begin{equation*}
\begin{gathered}
    h_i(\widehat{B}_i(W)) \defas  \mathbb{E}_{v}\left[\log\left(1 + \exp\left((-1)^{i}\left(\sqrt{B_{i,11}(W)} v  + \widehat{B}_{i,13}(W)\right)+ \delta \sqrt{\widehat{B}_{1,44}(W)}\right)\right)\right].
\end{gathered}
\end{equation*}
since $\mu_2 = -\mu_1$ and $\widehat{B}_{1,44}(W) = \widehat{B}_{2,44}(W)$. Hence, we obtain the following deterministic equivalents for $\mathcal{R}$ and $\mathcal{R}^{\operatorname{adv}}$
\begin{equation}\label{eq:det:equiv:LR:risks}
    \begin{gathered}
        \mathrsfs{R}(t) = \sum_{i=1}^2 p_i  \mathbb{E}_{v}\left[\log\left(1 + \exp\left((-1)^{i}\sqrt{\mathrsfs{B}_{i,11}(t)} v - \mathfrak{m}(t)\right)\right)\right],
        \\
        \mathrsfs{R}^{\operatorname{adv}}(t) = \sum_{i=1}^2 p_i  \mathbb{E}_{v}\left[\log\left(1 + \exp\left((-1)^{i}\sqrt{\mathrsfs{B}_{i,11}(t)} v -\mathfrak{m}(t) + \delta \sqrt{\widehat{\mathrsfs{B}}_{1,44}(t)}\right)\right)\right].
    \end{gathered}
\end{equation}
Here $\mathfrak{m}_{1} = - \mathfrak{m}_2$ so we drop the subscript in $\mathfrak{m}$ as we only need to track one of the means. Analogously to \citep{collinswoodfin2025Exact}, we define the variables $w_{ij}(k)$ for $i=1,2$ and $j=1,2$ to denote the robust class probabilities with $r_i(k) =\sqrt{X_{k}^\top K_iX_k} v + (-1)^{i-1} X_k^\top \mu$ at iteration $k$ of SGD
\begin{equation*}
    \begin{gathered}
        w_{11}(k) = \frac{e^{r_1(k) - \delta \|X\|}}{1+ e^{r_1(k) - \delta \|X\|}}, \quad w_{12}(k) = \frac{1}{1+ e^{r_1(k) - \delta \|X\|}}, 
        \\
         w_{21}(k) = \frac{e^{r_2(k) - \delta \|X\|}}{1+ e^{r_2(k) - \delta \|X\|}}\quad \text{and}\quad w_{22}(k) = \frac{1}{1+ e^{r_2(k) - \delta \|X\|}}.
    \end{gathered}
\end{equation*}
Hence, their deterministic equivalents are given by
\begin{equation}\label{eq:def:logits:detequiv}
    \begin{gathered}
        \mathfrak{u}_{12}(t) \defas \left(1+ e^{\sqrt{\mathrsfs{B}_{1,11}(t)} v + \mathfrak{m}(t) -\delta \sqrt{\widehat{\mathrsfs{B}}_{1,44}}(t)}\right)^{-1},\quad \mathfrak{u}_{11}(t) \defas 1 - \mathfrak{u}_{12} (t)
        \\
        \mathfrak{u}_{22}(t) \defas \left(1+ e^{\sqrt{\mathrsfs{B}_{2,11}(t)} v -\mathfrak{m}(t) -\delta \sqrt{\widehat{\mathrsfs{B}}_{1,44}} (t)}\right)^{-1},\quad \mathfrak{u}_{21}(t) \defas 1 - \mathfrak{u}_{22}(t)
    \end{gathered}
\end{equation}
We now compute $\partial_{11} h_i$ (i.e. $H_{1,t}$), $\partial_{13} h_i$ (i.e. $H_{2,t}$) and $\partial_{44} h_i$. From the chain rule, we have
\begin{equation}
    \begin{aligned}
        \partial_{11} h_i(\widehat{B}_i)
         &
         = 
         \frac{(-1)^{i}}{2\sqrt{B_{i,11}}}\mathbb{E}_{v}\left[v \cdot \sigma\left((-1)^{i}\sqrt{B_{i,11}} v - \widehat{B}_{1,13} +\delta \sqrt{\widehat{B}_{1,44}}\right)\right],
         \\
         \partial_{13} h_i(\widehat{B}_i)
         &= 
         \frac{(-1)^{i}}{2}\mathbb{E}_{v}\left[\sigma\left((-1)^{i}\sqrt{B_{i,11}} v - \widehat{B}_{1,13} +\delta \sqrt{\widehat{B}_{1,44}}\right)\right],
         \\
         \partial_{44} h_i(\widehat{B}_i)
         &
         = 
         \frac{\delta}{2\sqrt{\widehat{B}_{1,44}}}\mathbb{E}_{v}\left[\sigma\left((-1)^{i}\sqrt{B_{i,11}} v - \widehat{B}_{1,13} +\delta \sqrt{\widehat{B}_{1,44}}\right)\right].
         \\
    \end{aligned}
\end{equation}
where $\sigma(x) \defas \tfrac{1}{1+e^{-x}}$ denotes the sigmoid function and $\partial_{13} h_i(\widehat{B}_i) = \partial_{31} h_i(\widehat{B}_i)$. From Stein's Lemma and since $\sigma'(x) = \sigma(x)(1-\sigma(x))$, we obtain
\[
\begin{aligned}
    \partial_{11} h_i(\widehat{B}_i)
    &= 
    \frac{1}{2}\mathbb{E}_{v}\bigg[\sigma\left((-1)^{i}\sqrt{B_{i,11}} v - \widehat{B}_{1,13} +\delta \sqrt{\widehat{B}_{1,44}}\right) \\
    &\qquad \times (1-\sigma\left((-1)^{i}\sqrt{B_{i,11}} v - \widehat{B}_{1,13} +\delta \sqrt{\widehat{B}_{1,44}}\right))\bigg].
\end{aligned}
\]
Using the notation in \eqref{eq:def:logits:detequiv}, we obtain the relations
\begin{equation}
    \begin{gathered}
        \partial_{11} h_i(\widehat{\mathrsfs{B}}_i(t))
         = 
         \frac{1}{2}\mathbb{E}_{v}\left[  \mathfrak{u}_{i1}(t)   \mathfrak{u}_{i2}(t) \right],\quad
         \partial_{13} h_1(\widehat{\mathrsfs{B}}_1(t))
         = 
         -\frac{1}{2}\mathbb{E}_{v}\left[  \mathfrak{u}_{12}(t) \right],
         \\
         \partial_{13} h_2(\widehat{\mathrsfs{B}}_2(t))
         = 
         \frac{1}{2}\mathbb{E}_{v}\left[  \mathfrak{u}_{21}(t) \right],
         \quad
         \partial_{44} h_1(\widehat{\mathrsfs{B}}_1(t))
         = 
         \tfrac{\delta}{2\sqrt{\widehat{\mathrsfs{B}}_{1,44}(t)}}\mathbb{E}_{v}\left[ \mathfrak{u}_{12}(t) \right]
         \\
         \text{and}\quad \partial_{44} h_2(\widehat{\mathrsfs{B}}_2(t))
         = 
         \tfrac{\delta}{2\sqrt{\widehat{\mathrsfs{B}}_{1,44}(t)}}\mathbb{E}_{v}\left[ \mathfrak{u}_{21}(t) \right].
    \end{gathered}
\end{equation}

As a function of the input $X^\top a$, we have 
\begin{equation}
    f_i'(x) = (-1)^{i}\sigma((-1)^{i}\cdot x).
\end{equation}
Hence, we obtain 
\begin{equation}
    \mathcal{I}_i(\widehat{B}_i) = \E_v\bigg[ \sigma^2\bigg((-1)^{i}\sqrt{X^\top K_i X} v -  X^\top \mu +  \delta \|X\|)\bigg],
\end{equation}
and its deterministic equivalents
\begin{equation}
    \mathcal{I}_1(\widehat{\mathrsfs{B}}_1(t)) = \E_v[( \mathfrak{u}_{12}(t))^2] \quad \text{and}\quad \mathcal{I}_2(\widehat{\mathrsfs{B}}_2(t)) = \E_v[( \mathfrak{u}_{21}(t))^2].
\end{equation}
Using this notation, the system of ODEs \eqref{eq:ODE:V_i(t):def} is given by
\begin{equation}\label{eq:odes:hard:label:logreg}
    \begin{aligned}
        &\frac{\dif \mathrsfs{V}_{j,11}(t)}{\dif t} = -2\gamma(t)\bigg( \sum_{i=1}^2 p_i  \lambda_j^{(i)} \mathbb{E}_{v}\left[  \mathfrak{u}_{i1}(t)   \mathfrak{u}_{i2}(t) \right] 
        \\
        &
        \qquad
 +   \tfrac{\delta}{\sqrt{\widehat{\mathrsfs{B}}_{1,44}(t)}}(p_1\mathbb{E}_{v}\left[ \mathfrak{u}_{12}(t) \right] +p_2\mathbb{E}_{v}\left[ \mathfrak{u}_{21}(t)\right])  \bigg)\mathrsfs{V}_{j,11}(t)
        \\
        &
        \qquad
        +2\gamma(t)(p_1\mathbb{E}_{v}\left[ \mathfrak{u}_{12}(t) \right]+p_2\mathbb{E}_{v}\left[ \mathfrak{u}_{21}(t)\right])\mathfrak{m}_{j,1, 1}(t)
        \\
        &
        \qquad
        +\frac{\gamma(t)^2}{d} \sum_{i=1}^2 p_i \E_v[( \mathfrak{u}_{i(\sim i)}(t))^2](\lambda_j^{(i)} + (\omega_j^\top \mu)^2),
        \\
        &
        \frac{\dif \mathfrak{m}_{j,1, 1}}{\dif t}  = -\gamma(t)\bigg(  \sum_{i=1}^2 p_i   \lambda_j^{(i)} \mathbb{E}_{v}\left[  \mathfrak{u}_{i1}(t)   \mathfrak{u}_{i2}(t) \right] 
        \\
        &   
        \qquad 
        +\tfrac{\delta}{\sqrt{\widehat{\mathrsfs{B}}_{1,44}(t)}}(p_1\mathbb{E}_{v}\left[ \mathfrak{u}_{12}(t) \right] +p_2\mathbb{E}_{v}\left[ \mathfrak{u}_{21}(t)\right])  \bigg)\mathfrak{m}_{j,1, 1}(t)
        \\
        &
        \qquad + \gamma(t)\left(p_1\mathbb{E}_{v}\left[  \mathfrak{u}_{12}(t)\right]+p_2  \E\left[\mathfrak{u}_{21}(t)\right] \right)(\mu^\top \omega_j)^2,
    \end{aligned}
\end{equation}
with initial condition $\mathrsfs{V}_{j,11}(0) = X_0^\top \omega_j \omega_j^\top X_0$, $\mathfrak{m}_{j,1, 1}(0) = X_0^\top \omega_j \omega_j^\top\mu$ and where $\sim i$ denotes not $i$ such that $\mathfrak{u}_{i(\sim i)} = \mathfrak{u}_{12}$ for $i=1$ and $\mathfrak{u}_{i(\sim i) } = \mathfrak{u}_{21}$ for $i=2$. Note that if $K_1 = K_2$ then we have $\mathfrak{u}_{12}(t) \overset{d}{=} \mathfrak{u}_{21}(t)$ and $\mathfrak{u}_{11}(t) \overset{d}{=} \mathfrak{u}_{22}(t)$.

\section{Numerical Simulation Details}\label{app:captions}

\paragraph{Figure \ref{fig:concentration}:} \textbf{Concentration of $\ell_2$-adversarial risk} on noiseless $\ell_2$-adversarial least squares with a single class $a\sim N(0,K)$ (left) and noiseless binary logistic regression with hard labels on a mixture of Gaussians (right) with different means and same covariance. For logistic regression, see Appendix~\ref{ex:binary:logreg}. For least squares, 30 runs of SGD with constant learning rate $\gamma = 0.3$ and $\delta = 0.3$ under power law covariance (See Assumption~\ref{ass:powerlaw:main}) with $\alpha = 1.25$, $\beta = 0.75$, $\mu = 0$, $p_1=1$ and $X_0 \sim \mathcal{N}(0, \Id_d/d)$. For logistic regression, 30 runs of SGD with constant learning rate $\gamma = 0.3$ and $\delta = 0.1$ under power law covariance (See Assumption~\ref{ass:powerlaw:main}) with $\alpha = 1.25$ and $K_1=K_2$. In the hard label setup, we impose a power law assumption on the mean \citep{collinswoodfin2025Exact} instead of the ground truth such that $(\langle \mu, \omega_j\rangle)^2 = \frac{1}{d}\left(\frac{j}{d}\right)^{\beta}$ with $\beta = 0.75$. We also set $\mu_1 = - \mu_2$, $p_1=p_2= 0.5$ and $X_0 \sim \mathcal{N}(0, \Id_d/d)$. The shaded regions represent the max and min values at each step across the 30 SGD runs. As dimension $d$ increases, in both plots the adversarial risk concentrates around the deterministic limit \textcolor{red}{(red)} described by the system of ODEs \eqref{eq:ODE:V_i(t):def} as predicted by Theorem~\ref{thm:concentration_statistic:adv}. See Section~\ref{app:captions} for simulation details.

\paragraph{Figure~\ref{fig:HSGD}:} \textbf{Comparison of $\ell_2$-adversarial least squares and $\ell_2$-regularized least squares with adaptive learning rate and regularization.} The left plot compares the paths of the deterministic equivalents of $\mathcal{R}^{\Adv}$ for AdvHSGD and HSGD and confirms Proposition~\ref{prop:adv:reg:volterra:stability}. The right plot compares the path of $\mathcal{R}^{\Adv}(X_k)$ for SGD with adaptive learning rate $\gamma^{\operatorname{Reg}}(t)$ and regularization $\lambda^{\operatorname{Reg}}(t)$ versus the deterministic equivalent computed from \eqref{eq:R:reg:det:equiv:volterra} and \eqref{eq:l2norm2:reg:det:equiv:volterra} and shows close agreement between them for a variety of constant $\gamma(t) \equiv \gamma$. The $\ell_2$-adversarial risk is computed using $\mathcal{R}^{\Adv}(X_k) = \mathcal{R}(X_k) + \delta \sqrt{\frac{2}{\pi}} \sqrt{2\mathcal{R}(X_k)}\|X_k\| + \frac{\delta^2}{2}\|X_k\|^2$. For the left plot, we fixed the learning rate $\gamma = 1$, the noise variance $\eta^2 =0$, $K = \Id_d$, $\mu=0$,  $p_1=1$, $X_0 \sim \mathcal{N}(0, 4\Id_d/d)$ and $X^\star \sim \mathcal{N}(0, \Id_d/d)$ across all runs and only vary $\delta$. For the right plot, we fix $\delta = 0.3$, $\eta^2 =0$, power law covariance (See Assumption~\ref{ass:powerlaw:main}) with $\alpha = 1.25$, $\beta = 0.75$, $\mu = 0$, $p_1=1$ and $X_0 \sim \mathcal{N}(0, 4\Id_d/d)$. All experiments are done for $d = 800$.

\paragraph{Figure~\ref{fig:polyak:linesearch:evidence}:} \textbf{SGD with exact line search $\gamma_{k}^{\operatorname{line}}$ or Polyak stepsize $\gamma_{k}^{\operatorname{Polyak, adv}}$ } matches closely the path of our system of ODEs \eqref{eq:ODE:V_i(t):def} with deterministic learning rates schedules $\gamma(t)=\gamma^{\operatorname{line}}(t)$ and $\gamma(t)=\gamma^{\operatorname{Polyak, adv}}(t)$ for $\mathcal{R}^\Adv(X_k)$ on noiseless $\ell_2$-adversarial least squares. For the plot with Polyak stepsize (left), we set $d=800$, $\eta = 0$, $X_0 \sim \mathcal{N}(0, 4\Id_d/d)$, $X^\star\sim \mathcal{N}(0, \Id_d/d)$, $\mu=0$, $p_1=1$, $K_1 = K_2 = \Id_d$. For the plot with exact line search, we set $d=800$, $\eta = 0$, $X_0 \sim \mathcal{N}(0, 4\Id_d/d)$, $\mu=0$, $p_1=1$, $K$ and $X^\star$ satisfying a power law setup (See Assumption~\ref{ass:powerlaw:main}) with $\alpha = 1.25$ and $\beta = 0.75$.

\paragraph{Figure~\ref{fig:polyak:anisotropy}:} \textbf{Comparison between Exact Line Search and Polyak Stepsize under weak anisotropy} on noiseless $\ell_2$-adversarial least squares for different values of $\delta$ in the three regimes of $X^{\star,\Adv}$ (See Proposition~\ref{prop:fixed:point:gamma:to:0:main}). The three plots illustrate the convergence of the $\ell_2$-adversarial risk and that, under weak anisotropy, exact line search and the Polyak stepsize perform similarly. See Figure~\ref{fig:polyak:delta} for results in the case of strong anisotropy. Across all the plots, we fix $X_0 \sim \mathcal{N}(0, \frac{1}{d}\Id_d)$, $X^\star = \mathbf{1}_d /\sqrt{d}$ where $\mathbf{1}$ is a vector of ones of size $d$, $d=800$, $\eta=0$, $\mu=0$, $p_1=1$ and $K$ satisfying the covariance setup presented in \eqref{eq:s_values:eigs} with $s =7.5$ which produces $\frac{1}{d}\Tr(K^2) \approx 1.12$. See~\eqref{eq:ratio:linesearch:Polyak} for more details.

\paragraph{Figure~\ref{fig:polyak:delta}:} \textbf{Comparison for Exact Line Search and Polyak Stepsize under strong anisotropy} on noiseless $\ell_2$-adversarial least squares for different values of $\delta$ in the three regimes of $X^{\star,\Adv}$ (See Proposition~\ref{prop:fixed:point:gamma:to:0:main}). The three plots illustrate the convergence of the $\ell_2$-adversarial risk and how $\delta$ and $\tilde{\lambda}^{\operatorname{eff}}(t)$ mitigate the impact of strong anisotropy on the discrepancy between the Polyak stepsize and exact line search. Across all the plots, we fix $X_0 \sim \mathcal{N}(0, \frac{1}{d}\Id_d)$, $X^\star = \mathbf{1}_d /\sqrt{d}$ where $\mathbf{1}$ is a vector of ones of size $d$, $d=800$, $\eta=0$, $\mu=0$, $p_1=1$ and $K$ satisfying the covariance setup presented in \eqref{eq:s_values:eigs} with $s =2.8$ which produces $\frac{1}{d}\Tr(K^2) \approx 5.25$. See~\eqref{eq:ratio:linesearch:Polyak} for more details.

\paragraph{Figure~\ref{fig:q}:} \textbf{Numerical evidence that $q(t) \defas \sqrt{\frac{\widehat{\mathrsfs{B}}_{44}(t)}{2\mathrsfs{R}(t)}}$ converges} for noiseless $\ell_2$-adversarial least squares. The first and second plots provide evidence that $q(t)$ converges for a variety of constant learning rates $\gamma$ (left) and $\delta$ (middle). In either plot, we fix one of the parameters and vary the other. Here we fix the parameters $d=800$, $\eta = 0$, $X_0 \sim \mathcal{N}(0, 4\Id_d/d)$, $K$ and $X^\star$ satisfy a power law relationship (See Assumption~\ref{ass:powerlaw:main}) with $\alpha = 1.25$ and $\beta = 0.75$. We set $\delta= 0.3$ when we vary $\gamma$ and fix $\gamma = 0.3$ when we vary $\delta$. The third plot (right) provides evidence that $q(t)$ converges for a variety of $X^\star$ and covariances $K$ under different power law setups (See Assumption~\ref{ass:powerlaw:main}) and validates that the fixed point equation in Proposition~\ref{prop:fixed:point} predicts the limiting value $q_\infty$. For this plot, we use the same setup with $\delta=  \gamma = 0.3$ fixed and vary $\alpha$ and $\beta$ according to the values given in the legend.

\end{document}

%% file: macros.tex
\def\N{\mathbb{N}}
\def\R{\mathbb{R}}

\def\Adv{{\operatorname{Adv}}}
\def\opt{{\operatorname{op}}}

\def\dif{\mathop{}\!\mathrm{d}}
\def\Proba{\mathbb{P}}

\def\EE{{\mathbb E}\,}

\def\defas{\stackrel{\operatorname{def}}{=}}

\DeclareDocumentCommand{\Prto} {o} {
  \IfNoValueTF {#1}
  {\overset{\Pr}{\longrightarrow}}
  { \xrightarrow[ #1 \to \infty]{\Pr }}
}

\DeclareDocumentCommand{\Asto} {o} {
  \IfNoValueTF {#1}
  {\overset{\text{\rm a.s.}}{\longrightarrow}}
  { \xrightarrow[ #1 \to \infty]{\text{\rm a.s.} }}
}

\DeclareDocumentCommand{\law} {o} {
  \IfNoValueTF {#1}
  {\overset{\text{law}}{=}}
  { \xrightarrow[ #1 \to \infty]{\Pr }}
}

\DeclareMathOperator{\E}{\mathbb{E}}

\newcommand{\Id}{I}

\newcommand{\cM}{\mathcal{M}}

\DeclareMathOperator{\grad}{Grad}
\DeclareMathOperator{\hess}{Hess}
\newcommand{\beq}{ \begin{equation} }
\newcommand{\eeq}{ \end{equation} }

\newcommand{\ip}[1]{\langle {#1} \rangle }

\newcommand{\tr}{\text{Tr}}

\newcommand{\argmin}{\text{arg\,min}}

\newcommand{\vertiii}[1]{{\left\vert\kern-0.25ex\left\vert\kern-0.25ex\left\vert #1 
    \right\vert\kern-0.25ex\right\vert\kern-0.25ex\right\vert}}

\newcommand{\WHSGD}{ \mathscr{X} }

\newcommand{\HSGD}{ \mathscr{W} }

\def\Dif{\mathop{}\!\mathrm{D}}